\documentclass[reqno,12pt]{amsart}
\usepackage{a4}
\usepackage[latin1]{inputenc}
\usepackage[T1]{fontenc}

\usepackage[normalem]{ulem} 
\usepackage{color}

\numberwithin{equation}{section}

\usepackage{amssymb}
\usepackage{amsmath}
\usepackage{amsfonts}
\usepackage{amsthm}
\usepackage{latexsym}
\usepackage{amscd}
\usepackage{verbatim}
\usepackage{url}
\usepackage{stmaryrd}
\usepackage{enumerate}
\usepackage{colonequals}
\usepackage{dsfont}
\usepackage{appendix}
\usepackage{mathtools}
\usepackage{stackrel}
\usepackage[all,cmtip,2cell]{xy}
\UseAllTwocells
\usepackage{enumitem}
\usepackage{fullpage}
\usepackage{tikz-cd}

\usepackage{txfonts}
\usepackage[mathcal]{euscript}

\usepackage[refpage]{nomencl}
\usepackage{multicol}
\makenomenclature
\renewcommand{\nompreamble}{\begin{multicols}{2}}
\renewcommand{\nompostamble}{\end{multicols}}

\newcommand{\sym}[1]{#1\nomenclature{#1}{}%
}
\newcommand{\symn}[1]{\nomenclature{#1}{}%
}
\newcommand{\shriek}{!}

\usepackage{imakeidx}
\makeindex
\newcommand{\nc}[1]{#1\index{#1}%
}
\newcommand{\ncn}[1]{\index{#1}%
}

\usepackage{hyperref}
\usepackage{cleveref}
\let\cref\Cref

{\bf}{\it}
\newtheorem*{thm*}{Theorem}
\newtheorem*{sch*}{Scholium}
\newtheorem{thm}{Theorem}[subsection]{\bf}{\it}
\newtheorem{prop}[thm]{Proposition}
\newtheorem{lemma}[thm]{Lemma}
\newtheorem{cor}[thm]{Corollary}

\theoremstyle{definition}
\newtheorem{dfn}[thm]{Definition}
\theoremstyle{remark}
\newtheorem{rmk}[thm]{Remark}
\theoremstyle{remark}
\newtheorem{exm}[thm]{Example}
\newtheorem{assu}[thm]{Assumption}
\newtheorem{cons}[thm]{Construction}
\newtheorem{nota}[thm]{Notation}

\renewcommand{\P}{\mathbb{P}}

\newcommand{\A}{\mathbb{A}}
\newcommand{\B}{\mathbb{B}}
\newcommand{\F}{\mathbb{F}}
\newcommand{\G}{\mathbb{G}}
\newcommand{\N}{\mathbb{N}}
\newcommand{\Q}{\mathbb{Q}}
\newcommand{\Z}{\mathbb{Z}}
\newcommand{\CAlg}{\mathrm{CAlg}}

\newcommand{\Cat}{\mathrm{Cat}}
\newcommand{\CAT}{\mathrm{CAT}}

\newcommand{\one}{\mathbf{1}}
\newcommand{\Prlmon}{{\rm Pr}^{{\rm L},\,\otimes}}

\newcommand{\Prl}{{\rm Pr}^{\rm L}}
\newcommand{\Prr}{{\rm Pr}^{\rm R}}
\newcommand{\yon}{\mathrm{y}}

\newcommand{\Adic}{\mathrm{Adic}}
\newcommand{\an}{\mathrm{an}}
\newcommand{\An}{\mathrm{An}}
\newcommand{\af}{\mathrm{af}}
\newcommand{\colim}{\mathrm{colim}}
\newcommand{\cst}{\mathrm{cst}}

\newcommand{\Corr}{\mathrm{Corr}}
\newcommand{\diag}{\mathrm{diag}}
\newcommand{\eff}{\mathrm{eff}}
\newcommand{\et}{\acute{\rm e}{\rm t}}
\newcommand{\fet}{{\rm f}\acute{\rm e}{\rm t}}

\newcommand{\Et}{\acute{\rm E}{\rm t}}
\newcommand{\Etgr}{\acute{\rm E}{\rm t}^{\rm gr}}
\newcommand{\FDA}{\mathbf{FDA}}
\newcommand{\Fun}{\mathrm{Fun}}

\newcommand{\hocolim}{\mathrm{hocolim}}
\newcommand{\Hom}{\mathrm{Hom}}
\newcommand{\id}{\mathrm{id}}
\newcommand{\inc}{\mathrm{inc}}
\newcommand{\Ind}{\mathrm{Ind}}
\newcommand{\Pro}{\mathrm{Pro}}
\newcommand{\Map}{\mathrm{Map}}
\newcommand{\Nis}{\mathrm{nis}}
\newcommand{\op}{\mathrm{op}}
\newcommand{\qcqs}{\mathrm{qcqs}}
\newcommand{\red}{\mathrm{red}}
\newcommand{\rig}{\mathrm{rig}}
\newcommand{\riget}{{\rm rig}\acute{\rm e}{\rm t}}
\newcommand{\rigNis}{\mathrm{rignis}}

\newcommand{\Rig}{\mathrm{Rig}}
\newcommand{\RigSpc}{\mathrm{RigSpc}}
\newcommand{\RigSch}{\mathrm{RigSch}}
\newcommand{\WComp}{\mathrm{WComp}}
\newcommand{\Comp}{\mathrm{Comp}}
\newcommand{\FSpc}{\mathrm{FSpc}}
\newcommand{\RigSm}{\mathrm{RigSm}}
\newcommand{\FRigSm}{\mathrm{FRigSm}}
\newcommand{\FRigEt}{{\rm FRig}\acute{\rm E}{\rm t}}
\newcommand{\Sm}{\mathrm{Sm}}
\newcommand{\Spa}{\mathrm{Spa}}
\newcommand{\Spec}{\mathrm{Spec}}
\newcommand{\Spect}{\mathrm{Spt}}
\newcommand{\Sp}{\mathcal{S}\mathrm{p}}
\newcommand{\Spf}{\mathrm{Spf}}
\newcommand{\st}{\mathrm{st}}
\newcommand{\Sus}{\mathrm{Sus}}
\newcommand{\Fin}{\mathrm{Fin}_*}
\newcommand{\M}{\mathrm{M}}
\newcommand{\Th}{\mathrm{Th}}

\newcommand{\DA}{\mathbf{DA}}
\newcommand{\FSch}{\mathrm{FSch}}

\newcommand{\FSm}{\mathrm{FSm}}
\newcommand{\FSH}{\mathbf{FSH}}
\newcommand{\Mdl}{\mathrm{Mdl}}
\newcommand{\Mod}{\mathrm{Mod}}

\newcommand{\RigDA}{\mathbf{RigDA}}
\newcommand{\RigDM}{\mathbf{RigDM}}
\newcommand{\RigSH}{\mathbf{RigSH}}
\newcommand{\PSh}{\mathrm{PSh}}

\newcommand{\Sch}{\mathrm{Sch}}
\newcommand{\Set}{\mathrm{Set}}
\newcommand{\Shv}{\mathrm{Shv}}
\newcommand{\SH}{\mathbf{SH}}

\newcommand{\uSH}{\mathbf{uSH}}
\newcommand{\Zar}{\mathrm{zar}}

\newcommand{\Rder}{\mathrm{R}}
\newcommand{\Lder}{\mathrm{L}}
\newcommand{\sat}{\mathrm{sat}}
\newcommand{\mot}{\mathrm{mot}}
\newcommand{\Tate}{\mathrm{T}}
\newcommand{\U}{\mathbb{U}}
\renewcommand{\H}{\mathrm{H}}
\newcommand{\hyp}{\wedge}
\newcommand{\pr}{\mathrm{pr}}
\newcommand{\Op}{\mathrm{Op}}
\newcommand{\ellnil}{\ell\text{-}\mathrm{nil}}
\newcommand{\ellcpl}{\ell\text{-}\mathrm{cpl}}

\newcommand{\fiber}{\mathrm{fib}}

\newcommand{\lft}{\mathrm{lft}}

\newcommand{\sft}{\mathrm{sft}}
\newcommand{\wc}{\mathrm{wc}}
\newcommand{\cp}{\mathrm{cp}}

\newcommand{\form}{\mathrm{for}}
\newcommand{\proper}{\mathrm{prop}}

\newcommand{\pvcd}{\mathrm{pvcd}}

\begin{document}
\title{The six-functor formalism for rigid analytic motives}

\author[Ayoub]{Joseph Ayoub}
\address{University of Zurich / 
LAGA - Universit\'e Sorbonne Paris Nord}
\email{joseph.ayoub@math.uzh.ch}
\urladdr{user.math.uzh.ch/ayoub/}
\author[Gallauer]{Martin Gallauer}
\address{Max-Planck-Institut f\"ur Mathematik}
\email{gallauer@mpim-bonn.mpg.de}
\urladdr{people.maths.ox.ac.uk/gallauer/}
\author[Vezzani]{Alberto Vezzani}
\address{Universit\`a degli Studi di Milano}
\email{alberto.vezzani@unimi.it}
\urladdr{users.mat.unimi.it/users/vezzani/}

\thanks{The first author is partially supported by the 
{\it Swiss National Science Foundation} (SNF), 
project 200020\_178729. The second author is supported by a Titchmarsh Fellowship of the University of Oxford.
The third author is partially 
supported by the Italian {\it Ministero dell'Universit\`a e della Ricerca} (MUR),  
project  PRIN 2022B24AY.
}

\newcommand{\kw}{Motives (algebraic, formal and rigid analytic), 
six-functor formalism, proper base change theorem}
\keywords{\kw}

\hypersetup{
    pdfcreator    = {\LaTeX{}},
%    psdextra, %% I had to remove this option as the document wouldn't compile; is it important?
    pdfauthor     = {Joseph Ayoub and Martin Gallauer and 
                     Alberto Vezzani},
    pdftitle      = {The six-functor formalism for rigid 
                     analytic motives},
    pdfkeywords   = {\kw}}

\begin{abstract}
We offer a systematic study of rigid analytic motives over 
general rigid analytic spaces, and we develop their six-functor 
formalism. A key ingredient is an extended proper base 
change theorem that we are able to justify by reducing to the case of 
algebraic motives. In fact, more generally, 
we develop a powerful technique for 
reducing questions about rigid analytic 
motives to questions about algebraic motives, which is likely to be useful in other contexts as well. 
We pay special attention to establishing our results 
without noetherianity assumptions on rigid analytic spaces. 
This is indeed possible using Raynaud's approach to rigid analytic geometry. 
\end{abstract}

\maketitle

\tableofcontents

\section*{Introduction}

In this paper, we study rigid analytic motives over general 
rigid analytic spaces and we develop a 
six-functor formalism for them. We have tried to free
our treatment from unnecessary hypotheses, 
and many of our main results hold in great 
generality, with the notable exception of 
Theorems \ref{thm:main-thm-}(2) and \ref{thm:compute-chi}
where we impose \'etale descent.
(This is necessary for the former but might be superfluous 
for the latter.) 
In this introduction, we 
restrict to \'etale rigid analytic motives with rational 
coefficients, for which our results are the most 
complete.\footnote{In fact, everything we say here 
holds more generally with coefficients in an arbitrary ring
when the class of rigid analytic spaces is accordingly 
restricted. For instance, if one is only considering rigid analytic 
spaces over $\Q_p$ of finite \'etale cohomological dimension, 
the results discussed in the introduction are valid 
with $\Z[p^{-1}]$-coefficients.}

\subsection*{The six-functor formalism} 

$\empty$

\smallskip

Rigid analytic motives were introduced in 
\cite{ayoub-rig} as a natural extension of the 
notion of a motive associated 
to a scheme. Given a rigid analytic space $S$, we denote by
$\RigDA_{\et}(S;\Q)$ the $\infty$-category of \'etale rigid analytic 
motives over $S$ with rational coefficients. 
By construction, it is naturally equipped with the structure 
of a symmetric monoidal $\infty$-category (see Definition 
\ref{dfn:rigsh-stable}).
 
Given a morphism of rigid analytic spaces $f:T \to S$, 
the functoriality of the construction yields an adjunction 
\begin{equation}
\label{eqn:intro-f-star}
f^*:\RigDA_{\et}(S;\Q) \rightleftarrows \RigDA_{\et}(T;\Q):f_*.
\end{equation}
When $f$ is locally of finite type, 
we construct in this paper another adjunction 
(see Definition \ref{dfn:f-!-up-down})
\begin{equation}
\label{eqn:intro-f-!}
f_!:\RigDA_{\et}(T;\Q) \rightleftarrows \RigDA_{\et}(S;\Q):f^!,
\end{equation}
i.e., we define the ``exceptional direct image'' and the ``exceptional inverse image'' functors associated to $f$. Our main goal in this paper is to show the following result.

\begin{sch*}
\label{scholium:six-functor}
The functors $f^*$, $f_*$, $f_!$, $f^!$, $\otimes$
and $\underline{\Hom}$ satisfy the usual properties of 
a six-functor formalism. These include: 
\begin{itemize}

\item the compatibility with composition of morphisms
(see Proposition
\ref{prop:basic-functorial-rigsh} and Corollary 
\ref{cor:ext-f-!-all-ft-mor});

\item the localization formula (see Proposition \ref{prop:loc1}(2));

\item the base change theorems (see Proposition 
\ref{prop:6f1}(3) for the smooth base change,
Theorem \ref{thm:general-base-change-thm}
for the quasi-compact base change, Theorem
\ref{thm:prop-base}
for the extended proper base change, and Proposition
\ref{prop:exchange-base-change-in-general-1}
for the exchange between the ``ordinary inverse image'' and the
``exceptional direct image'' functors);

\item the canonical equivalences $f_!\simeq f_*$, when $f$ is proper 
(see Example \ref{exm:f!-2-cases}), 
and the equivalence $f^!\simeq f^*(d)[2d]$,\footnote{Strictly speaking,
Theorem \ref{thm:f-!-for-f-smooth} only gives 
that $f^!$ is equivalent to the twist of $f^*$ 
by the Thom space associated to $\Omega_f$. However, in the case of 
$\RigDA_{\et}(-;\Q)$, 
Thom spaces are globally given by Tate twists. 
Indeed, arguing as in \cite[Remarque 11.3]{ayoub-etale},
this would follow from the property that the mapping space 
$\Map_{\RigDA_{\et}(S;\Q)}(\Q,\Q)\simeq 
\Q^{\pi_0(S)}$ is coconnective. The latter property can be established 
using Theorem \ref{thm:compute-chi}
and \cite[Proposition 11.1]{ayoub-etale}. We leave the 
details for the interested reader.} 
when $f$ smooth of pure relative dimension $d$ 
(see Theorem \ref{thm:f-!-for-f-smooth});

\item the compatibility with tensor product, the
projection formula and duality (see Proposition
\ref{prop:basic-functorial-rigsh}, and Corollaries 
\ref{cor:strong-dual-}, \ref{cor:proj-form-f-!} and 
\ref{cor:weak--dual-}).
 
\end{itemize}
\end{sch*}

Of course, our six-functor formalism matches 
the one developed by Huber \cite{huber} 
for the \'etale cohomology of adic spaces.
(Similar formalisms for \'etale cohomology were
also developed by Berkovich 
\cite{berk-etale-coh} and de Jong--van der Put \cite{dJ-vdP}.)

A partial six-functor formalism for rigid analytic motives 
was obtained in \cite[\S 1.4]{ayoub-rig} at a minimal cost
as an application of the theory developed 
in \cite[Chapitre 1]{ayoub-th1}. Given a
non-Archimedean field $K$ and a classical affinoid $K$-algebra $A$,
the assignment sending a finite type $A$-scheme $X$ to the 
$\infty$-category $\RigDA_{\et}(X^{\an};\Q)$ gives rise to a stable homotopical functor in the sense of 
\cite[D\'efinition 1.4.1]{ayoub-th1}. (Here $X^{\an}$ is the 
analytification of $X$.) Applying 
\cite[Scholie 1.4.2]{ayoub-th1}, 
we have in particular an adjunction as in 
\eqref{eqn:intro-f-!} 
under the assumption that $f$ is algebraizable, i.e., that $f$ 
is the analytification of a morphism 
between finite type $A$-schemes, 
for some unspecified classical affinoid $K$-algebra $A$. 
Clearly, it is unnatural and unsatisfactory to 
restrict to algebraizable morphisms, and it is our objective in this
paper to remove this restriction. The key ingredient for doing so 
is Theorem \ref{thm:prop-base} which we may consider as an 
extended proper base change theorem for commuting 
direct images along proper morphisms and extension by zero 
along open immersions.
It is worth noting that in the algebraic setting, the extended 
proper base change theorem is essentially a reformulation 
of the usual one, but this is far from true in the 
rigid analytic setting.
In fact, the usual proper base change theorem 
in rigid analytic geometry is a 
particular case of the so-called quasi-compact base change theorem
(see Theorem \ref{thm:general-base-change-thm})
which is an easier property.

The extended proper base change theorem is already known 
if one restricts to projective morphisms in which case it can
be deduced from the partial six-functor formalism developed in 
\cite[\S 1.4]{ayoub-rig}. However, in the rigid analytic setting, 
it is not possible to deduce the general case of proper
morphisms from the special case of projective morphisms. 
Indeed, the classical argument 
used in \cite[Expos\'e XII]{SGAIV3}
for reducing the proper case to the projective case 
relies on Chow's lemma 
for which there is no analogue in rigid analytic geometry.
(For instance, there are proper rigid analytic tori which are 
not algebraizable \cite[\S 6.6]{fvdp}.)
Therefore, a new approach was necessary for proving 
Theorem \ref{thm:prop-base} in general.

\subsection*{Rigid motives as modules in algebraic motives}

$\empty$

\smallskip

Our approach is based on another 
contact point with algebraic geometry: instead of using 
the analytification functor from schemes to rigid analytic spaces, 
we go backward 
and associate to a rigid analytic space $X$ the pro-scheme consisting 
of the special fibers of the different formal models of $X$. 
We now sketch the main idea of our construction, which is detailed in Section \ref{subsect:const-tilde-xi}.

Let $\mathcal{S}$ be a formal scheme. We may associate to it the 
$\infty$-category of formal motives $\FDA_{\et}(\mathcal{S};\Q)$ which is canonically equivalent to the 
$\infty$-category of (algebraic) motives 
$\DA_{\et}(\mathcal{S}_{\sigma};\Q)$ over the special fiber 
$\mathcal{S}_{\sigma}$ (see 
Theorem \ref{thm:formal-mot-alg-mot}). The ``generic fiber'' functor induces a functor
$$\xi_{\mathcal{S}}:
\DA_{\et}(\mathcal{S}_{\sigma};\Q)\simeq \FDA_{\et}(\mathcal{S};\Q)
\to \RigDA_{\et}
(\mathcal{S}^{\rig};\Q)$$
where $\mathcal{S}^{\rig}$ is the rigid analytic space associated to 
$\mathcal{S}$. It is immediate to see that $\xi_{\mathcal{S}}$ is monoidal and has a right adjoint $\chi_{\mathcal{S}}$
sending the unit object of $\RigDA_{\et}(\mathcal{S}^{\rig};\Q)$
to a commutative algebra object $\chi_{\mathcal{S}}\Q$ of 
$\DA_{\et}(\mathcal{S}_{\sigma};\Q)$. 
Moreover, the functor $\chi_{\mathcal{S}}$ admits a factorization 
$$\RigDA_{\et}(\mathcal{S}^{\rig};\Q)
\xrightarrow{\widetilde{\chi}_{\mathcal{S}}}
\DA_{\et}(\mathcal{S}_{\sigma};\chi_{\mathcal{S}} \Q) 
\xrightarrow{\rm ff} \DA_{\et}(\mathcal{S}_{\sigma};\Q)$$
where $\DA_{\et}(\mathcal{S}_{\sigma};\chi_{\mathcal{S}} \Q)$ 
denotes the $\infty$-category of 
$\chi_{\mathcal{S}}\Q$-modules in 
$\DA_{\et}(\mathcal{S}_{\sigma};\Q)$
and where ${\rm ff}$ is the forgetful functor. Also, the functor 
$\widetilde{\chi}_{\mathcal{S}}$ admits a left adjoint 
$$\begin{array}{rcl}
\widetilde{\xi}_{\mathcal{S}}:\DA_{\et}(\mathcal{S}_{\sigma};\chi_{\mathcal{S}} \Q) & \to & \RigDA_{\et}(\mathcal{S}^{\rig};\Q)\\
& \vspace{-.3cm} &\\
M &\mapsto & \xi_{\mathcal{S}}(M)\otimes_{\xi_{\mathcal{S}}
\chi_{\mathcal{S}}\Q}\Q.
\end{array}$$

Now, if $S$ is a quasi-compact and quasi-separated 
rigid analytic space, we may consider the cofiltered category 
$\Mdl(S)$ of formal models of $S$
(see Notation \ref{nota:Mdl}).
The above construction yields a functor 
$$\widetilde{\xi}_S:\underset{\mathcal{S}\in 
\Mdl(S)}{\colim}\, \DA_{\et}(\mathcal{S}_{\sigma};\chi_{\mathcal{S}}\Q)
\to \RigDA_{\et}(S;\Q).$$
One of our main results is the following (see Theorem 
\ref{thm:main-thm-} and Remark \ref{rmk:main-on-rig}).

\begin{thm*}
Restrict to rigid analytic spaces which are quasi-compact, 
quasi-separated and having finite Krull dimension. 
The natural transformation $\widetilde{\xi}$ exhibits the functor 
$S\mapsto \RigDA_{\et}(S;\Q)$ as 
the \'etale sheafification of the functor
$S\mapsto \colim_{\mathcal{S}\in \Mdl(S)}\, \DA_{\et}(\mathcal{S}_{\sigma};\chi_{\mathcal{S}}\Q)$
viewed as a presheaf valued in presentable $\infty$-categories. 
\end{thm*} 

We use the above description of the $\infty$-categories
$\RigDA_{\et}(S;\Q)$  
to deduce the extended proper base change theorem in 
rigid analytic geometry
(i.e., Theorem \ref{thm:prop-base}) 
from its algebraic analogue. In fact, it turns 
out that we only need a formal consequence of this description which 
happens to be also a key ingredient in its proof, namely Theorem 
\ref{thm:proj-form-chi-xi} (see also Theorem \ref{thm:chijsharp}). 
Once Theorem \ref{thm:prop-base} is proven, it is easy to construct 
the adjunction \eqref{eqn:intro-f-!}.

We also point out that the commutative algebras 
$\chi_{\mathcal{S}}\Q$ admit a concrete description.
For precise statements, see Theorem 
\ref{thm:compute-chi} and Remark 
\ref{rmk:chi-B-computable}. Moreover, the special case
of the above theorem where we take for $S=\Spf(k[[\pi]])^{\rig}$,
with $k$ a field of characteristic zero, is 
tightly connected to the main result of \cite[Chapitre 1]{ayoub-rig}.
This will be explained in Remark \ref{rmk:k((pi))}.

\subsection*{Further results and applications}

$\empty$

\smallskip

Besides the six-functor formalism, the paper contains several 
foundational results on motives of rigid analytic spaces which are of independent interest.
In particular, we study the descent property of the 
$\infty$-categories $\RigDA_{\et}(S;\Q)$ for the 
\'etate topology; see Theorem 
\ref{thm:hyperdesc}.

Another notable result is Theorem 
\ref{thm:anstC} which, roughly speaking, asserts that 
$\RigDA_{\et}(-;\Q)$ transforms limits of certain rigid analytic 
pro-spaces into colimits of presentable $\infty$-categories.
A similar property is also true for $\DA_{\et}(-;\Q)$ but 
the proof in the rigid analytic setting is much more involved
and relies on approximation techniques as those used in 
the proof of \cite[Proposition 4.5]{vezz-fw}. 
We also like to mention that this continuity property for 
$\RigDA_{\et}(-;\Q)$ plays a crucial role (along with 
many of the results described above) in the 
recent paper \cite{lbv} where a new relative cohomology theory for 
rigid analytic varieties over a positive characteristic 
perfectoid space $P$ 
is defined and studied. Interestingly, this relative cohomology 
theory takes values in solid quasi-coherent sheaves over 
the relative Fargues--Fontaine curve associated to $P$.

\subsubsection*{Acknowledgments}
We thank Tony Yue Yu for asking if there was
a good duality theory for the motives of smooth and proper
rigid analytic spaces. Answering this question was
one of our motivations for this work.
We thank Kazuhiro Fujiwara and Fumiharu Kato
for clarifying some points in their book 
``Foundations of Rigid Geometry, I''. 
We also thank Denis Nardin and Marco Robalo 
for helpful discussions about Proposition 
\ref{prop:Prlomega-comp-gen}, and the referees for their constructive comments and recommendations.

\subsection*{Notation and conventions} 

\subsubsection*{\texorpdfstring{$\infty$}{oo}-Categories}
We use the language of $\infty$-categories
following Lurie's books \cite{lurie} and \cite{lurie:higher-algebra}.
The reader familiar with the content of 
these books will have no problem understanding 
our notation pertaining to higher category theory and higher algebra 
which are often very close to those in loc.~cit. 
Nevertheless, we list below some of these notational conventions
which we use frequently.
\ncn{categories@$\infty$-categories}

As usual, we employ the device of Grothendieck universes, and we 
denote by \sym{$\Cat_{\infty}$} the $\infty$-category of small 
$\infty$-categories and \sym{$\CAT_{\infty}$} the $\infty$-category 
of locally small, but possibly large $\infty$-categories. We denote by 
\sym{$\CAT_{\infty}^{\Lder}$} (resp. \sym{$\CAT_{\infty}^{\Rder}$}) the 
wide sub-$\infty$-category of $\CAT_{\infty}$ spanned by 
functors which are left (resp. right) adjoints.
Similarly, we denote by \sym{$\Prl$} (resp. \sym{$\Prr$}) the 
$\infty$-category of presentable $\infty$-categories and 
left adjoint (resp. right adjoint) functors. 
We denote by $\Prl_{\omega}\subset \Prl$ 
(resp. $\Prr_{\omega}\subset \Prr$) the 
sub-$\infty$-category of compactly generated $\infty$-categories 
and compact-preserving functors (resp. functors commuting with 
filtered colimits). 
\ncn{categories@$\infty$-categories!presentable}
\symn{$\Prl_{\omega}$}
\symn{$\Prr_{\omega}$}

$1$-Categories are typically referred to as just ``categories'' and viewed as $\infty$-categories via the nerve construction. For emphasis, we sometimes call them ``ordinary categories''.
We denote by \sym{$\mathcal{S}$} the $\infty$-category of small \nc{spaces}, 
by \sym{$\Sp$} the $\infty$-category of small \nc{spectra} and by 
$\Sp_{\geq 0}\subset \Sp$ its full sub-$\infty$-category of 
connective spectra.
\symn{$\Sp_{\geq 0}$}

Given an $\infty$-category $\mathcal{C}$, we denote by 
$\Map_{\mathcal{C}}(x,y)$ the mapping space between two objects 
$x$ and $y$ in $\mathcal{C}$. Given another $\infty$-category 
$\mathcal{D}$, we denote by $\Fun(\mathcal{C},\mathcal{D})$
the $\infty$-category of functors from $\mathcal{C}$ to $\mathcal{D}$.
If $\mathcal{C}$ is small, we denote by 
$\mathcal{P}(\mathcal{C})=\Fun(\mathcal{C}^{\op},\mathcal{S})$ the 
$\infty$-category of 
presheaves on $\mathcal{C}$ and by $\yon:\mathcal{C}
\to \mathcal{P}(\mathcal{C})$ the Yoneda embedding. 
\symn{$\Map$}
\symn{$\Fun$}
\symn{$\mathcal{P}$}
\symn{$\yon$}

\subsubsection*{Monoidal structures}
By ``\nc{monoidal} $\infty$-category'' we always mean ``symmetric
monoidal $\infty$-category'', i.e., a 
coCartesian fibration $\mathcal{C}^{\otimes} 
\to \Fin$ such that the induced functor 
$(\rho^i_!)_i:\mathcal{C}_{\langle n\rangle}
\to \prod_{1\leq i \leq n} \mathcal{C}_{\langle 1\rangle}$
is an equivalence for all $n\geq 0$. 
(Recall that $\Fin$ is the category of finite pointed sets,
$\langle n\rangle=\{1,\ldots, n\}\cup\{*\}$ and 
$\rho^i:\langle n\rangle \to \langle 1\rangle$ 
is the unique map such that $(\rho^i)^{-1}(1)=\{i\}$.)
If $\mathcal{C}^{\otimes}$ is a monoidal $\infty$-category, 
we denote by $\CAlg(\mathcal{C})$ the $\infty$-category of 
commutative algebras in $\mathcal{C}$. If $A\in \CAlg(\mathcal{C})$,
we denote by $\Mod_A(\mathcal{C})$ the $\infty$-category of $A$-modules.
Using Lurie's straightening construction, a monoidal category 
can be considered as an object of $\CAlg(\CAT_{\infty})$, i.e.,
as a commutative algebra in $\CAT_{\infty}^{\times}$.
\symn{$\Fin$}
\symn{$\CAlg$}
\symn{$\Mod$}

The $\infty$-categories
$\Prl$ and $\Prl_{\omega}$ underly monoidal $\infty$-categories 
\sym{$\Prlmon$} and \sym{$\Prlmon_{\omega}$}.
A monoidal $\infty$-category is said to be 
presentable (resp. compactly generated) if it belongs
to $\CAlg(\Prl)$ (resp. $\CAlg(\Prl_{\omega})$).

\subsubsection*{Sites and topoi}
If $\mathcal{C}$ is an $\infty$-category endowed with 
a topology $\tau$, we denote by $(\mathcal{C},\tau)$ 
the corresponding \nc{site}. 
We denote by $\Shv_{\tau}(\mathcal{C})
\subset \mathcal{P}(\mathcal{C})$ the 
full sub-$\infty$-category of $\tau$-sheaves 
and by $\Shv^{\hyp}_{\tau}(\mathcal{C})\subset 
\Shv_{\tau}(\mathcal{C})$ its full sub-$\infty$-category
of $\tau$-hypersheaves. We use 
\sym{$\Lder_{\tau}$} to denote the $\tau$-sheafification functor 
as well as the $\tau$-hypersheafification functor. 
(The context will make it clear which of the two we mean.)
Morphisms of sites $(\mathcal{C},\tau)\to (\mathcal{C}',\tau')$ 
are underlain by functors in the opposite direction 
$\mathcal{C}'\to\mathcal{C}$.
In particular, a cofiltered inverse system of sites
$(\mathcal{C}_{\alpha},\tau_{\alpha})_{\alpha}$ 
is underlain by a filtered direct system of $\infty$-categories,
and we write $\lim_{\alpha}(\mathcal{C}_{\alpha},\tau_{\alpha})$
for the site $(\colim_{\alpha}\mathcal{C}_{\alpha},\tau)$
where $\tau$ is the topology generated by the $\tau_{\alpha}$'s 
in the obvious way.
\ncn{site!morphism}
\ncn{site!limit}
\symn{$\Shv$}
\symn{$\Shv^{\hyp}$}

In the cases of most interest to us, 
the sites are underlain by ordinary categories.
In these cases, we follow the classical terminology and 
say that a morphism of sites is an equivalence of sites 
if it induces an equivalence on the associated ordinary topoi.
(This will be repeated each time, to avoid any possible confusion.)
\ncn{site!equivalence}

\subsubsection*{Formal and rigid analytic geometries}
We use Raynaud's approach to rigid analytic geometry \cite{Raynaud}
which is systematically developed in the books of Abbes \cite{abbes}
and Fujiwara--Kato \cite{fujiwara-kato}. In fact, we mainly
use \cite{fujiwara-kato} where rigid analytic spaces are 
introduced without noetherianness assumptions.

We denote formal schemes by calligraphic letters $\mathcal{X}$,
$\mathcal{Y}$, etc., and rigid analytic spaces by roman letters
$X$, $Y$, etc. Formal schemes are always assumed adic of 
finite ideal type in the sense of \cite[Chapter I, Definitions 1.1.14 \& 1.1.16]{fujiwara-kato}. Morphisms of formal schemes are always assumed
adic in the sense of \cite[Chapter I, Definition 1.3.1]{fujiwara-kato}.
Given a formal scheme $\mathcal{X}$, we denote by 
$\mathcal{X}^{\rig}$ its associated rigid analytic space
which we call the Raynaud generic fiber 
(or simply the generic fiber) of $\mathcal{X}$. 
Recall that $\mathcal{X}^{\rig}$ is simply $\mathcal{X}$ considered
in the localisation of the category of formal schemes with respect 
to admissible blowups. A general rigid analytic space is 
locally isomorphic to the generic fiber of a formal scheme.
As we show in Corollary \ref{cor:from-Huber-to-FK}, 
the category of stably uniform adic spaces (see \cite{buzz-ver}) 
embeds fully faithfully in the category of rigid analytic spaces.
\ncn{formal schemes}
\ncn{schemes!formal|see{formal schemes}}
\ncn{formal schemes!generic fiber}
\ncn{rigid analytic spaces}
\symn{$(-)^{\rig}$}

Given a rigid analytic space $X$, we denote by $|X|$ 
the associated topological space {(see Notation \ref{not:visualisation-})}. This is constructed in 
\cite[Chapter II, \S 3.1]{fujiwara-kato} where it is called the 
Zariski--Riemann space of $X$. The space $|X|$ is 
endowed with a sheaf of rings $\mathcal{O}_X$, called the 
structure sheaf, and a subsheaf of rings
$\mathcal{O}_X^+\subset \mathcal{O}_X$, 
called the integral structure sheaf.
(In \cite[Chapter II, \S 3.2]{fujiwara-kato}, the integral structure sheaf
is denoted by $\mathcal{O}^{\rm int}_X$, but we prefer to 
follow Huber's notation in \cite{huber}.)
\symn{$\lvert(-)\rvert$}
\symn{$\mathcal{O}$}
\symn{$\mathcal{O}^+$}

\subsubsection*{Motives (algebraic, formal and rigid analytic)}
We fix a commutative ring spectrum $\Lambda$, i.e., an 
object of $\CAlg(\mathcal{S}{\rm p})$ 
which we assume to be connective for simplicity.

Given a scheme $S$, we denote by $\SH_{\tau}(S;\Lambda)$ 
the Morel--Voevodsky $\infty$-category of $\tau$-motives on $S$ 
with coefficients in 
$\Lambda$ {(see, for example, \cite{jardine-mot})}. Here $\tau\in \{\Nis,\et\}$ 
is either the Nisnevich or the \'etale topology. 
When $\tau$ is the Nisnevich topology, we sometimes
omit the subscript ``$\Nis$'' and speak simply of 
motives over $S$. If $\Lambda$ is the Eilenberg--Mac Lane spectrum 
associated to a commutative dg-ring (also denoted by $\Lambda$), we
usually write $\DA_{\tau}(S;\Lambda)$ instead of 
$\SH_{\tau}(S;\Lambda)$.
\symn{$\SH$}
\symn{$\DA$}

Given a formal scheme $\mathcal{S}$, we denote by 
$\FSH_{\tau}(\mathcal{S};\Lambda)$ the $\infty$-category of 
formal $\tau$-motives on $\mathcal{S}$ with coefficients in 
$\Lambda$ {(see Definition \ref{def:DAeff-form})}. 
Similarly, given a rigid analytic space $S$, 
we denote by $\RigSH_{\tau}(S;\Lambda)$
the $\infty$-category of rigid analytic $\tau$-motives on 
$S$ with coefficients in $\Lambda$ {(see 
Definition \ref{def:DAeff})}.
Here again, $\tau\in \{\Nis,\et\}$ is 
either the Nisnevich or the \'etale topology, and
when $\tau$ is the Nisnevich topology we sometimes omit the subscript 
``$\Nis$''. If $\Lambda$ is the Eilenberg--Mac Lane spectrum 
associated to a commutative dg-ring (also denoted by $\Lambda$), 
we usually write 
$\FDA_{\tau}(\mathcal{S};\Lambda)$ and $\RigDA_{\tau}(S;\Lambda)$
instead of $\FSH_{\tau}(\mathcal{S};\Lambda)$ and 
$\RigSH_{\tau}(S;\Lambda)$.
\symn{$\FSH$}
\symn{$\FDA$}
\symn{$\RigDA$}
\symn{$\RigSH$}

We also consider the unstable (aka., effective) and/or hypercomplete 
variants of these motivic $\infty$-categories, 
which we refer to using superscripts ``$\eff$'' and/or ``$\hyp$''.
For example, $\SH_{\tau}^{\hyp}(S;\Lambda)$ is the 
Morel--Voevodsky $\infty$-category of hypercomplete $\tau$-motives
and $\SH_{\tau}^{\eff,\,\hyp}(S;\Lambda)$ is its effective version.
If a statement is equally valid 
for the $\Tate$-stable
and the effective motivic $\infty$-categories, we use 
the superscript ``$(\eff)$''. For example, the sentence 
``the $\infty$-category $\RigDA^{(\eff)}_{\tau}(S;\Lambda)$ is 
presentable'' means that both $\infty$-categories
$\RigDA^{\eff}_{\tau}(S;\Lambda)$ and $\RigDA_{\tau}(S;\Lambda)$
are presentable. 
We use the superscripts ``$(\hyp)$'', ``$(\eff,\,\hyp)$''
in a similar way. For example, the sentence 
``$S\mapsto \SH_{\tau}^{(\eff,\,\hyp)}(S;\Lambda)$
is a $\Prl$-valued $\tau$-(hyper)sheaf''
means that we have two $\tau$-sheaves, namely
$\SH_{\tau}^{\eff}(-;\Lambda)$ and 
$\SH_{\tau}(-;\Lambda)$, and two $\tau$-hypersheaves,
namely $\SH_{\tau}^{\eff,\,\hyp}(-;\Lambda)$ and 
$\SH_{\tau}^{\hyp}(-;\Lambda)$.
\symn{$\RigDA^{(\eff)}$}

\section{Formal and rigid analytic geometry}

\label{sect:rigid-analytic-spaces}

In this section, we gather a few results in rigid 
analytic geometry which we need later in the paper.
We use Raynaud's approach \cite{Raynaud} which can be 
summarised roughly as follows: the category of rigid analytic 
spaces is the localisation of the category of formal schemes 
with respect to admissible blowups. This is correct up to 
imposing the right conditions on formal schemes and slightly enlarging 
the localised category to allow gluing along open immersions. 
Raynaud's approach has been systematically developed 
by Abbes \cite{abbes} and 
Fujiwara--Kato \cite{fujiwara-kato}.
We will mainly follow the book \cite{fujiwara-kato}
where rigid analytic spaces are 
introduced without noetherianness assumptions. 
Indeed, one of the aims of the paper is to show that there 
are reasonable $\infty$-categories of rigid analytic motives over 
general rigid analytic spaces. 
We warn the readers that many results in \cite{fujiwara-kato}
require noetherianness assumptions, especially when it comes 
to the study of quasi-coherent sheaves. 
However, the theory of quasi-coherent sheaves is 
largely irrelevant for what we
do in this paper.

The reader who is only interested in motives of classical rigid 
analytic varieties in the sense of Tate 
and who is accustomed with Raynaud's notion of  
formal models, may skip this section and refer back 
to it when needed.

\subsection{Recollections}

$\empty$

\smallskip

\label{subsubsect:recoll-rig}

Unless otherwise stated, \nc{adic rings} are always assumed to
be complete of finite ideal type in the sense of 
\cite[Chapter I, Definitions 1.1.3 \& 1.1.6]{fujiwara-kato}.
(This is also the convention of \cite[D\'efinition 1.8.4]{abbes}
and \cite[Section 1]{huber1}.)
Thus, an adic ring $A$ is a complete linearly topologized ring 
whose topology is $I$-adic for some ideal $I\subset A$ of finite type. 
Morphisms between adic rings are always assumed to be adic in the sense
of \cite[Chapter I, Definition 1.1.15]{fujiwara-kato}. Thus, 
a morphism of adic rings
$A \to B$ is a ring homomorphism such that
$IB$ is an ideal of definition of $B$ 
for one (and hence every) ideal of definition $I$ of $A$.
\ncn{adic rings!morphism}

A useful basic fact when dealing with adic rings is the existence of 
$I$-adic completions in the sense of  
\cite[Chapter 0, Definition 7.2.6]{fujiwara-kato}.

\begin{lemma}
\label{lem:I-adic-compl}
Let $A$ be a ring, $I\subset A$ a finitely generated ideal
and $M$ an $A$-module. The Hausdorff completion 
$\widehat{M}=\lim_{n\in\N}M/I^nM$ of the $A$-module
$M$ endowed with the 
$I$-adic topology is itself an $I$-adic topological $A$-module.
More precisely, for $m\geq 0$ we have: 
\begin{itemize}

\item $I^m\widehat{M}$ is closed in $\widehat{M}$ 
and coincides with $\widehat{I^mM}=\lim_{n\in\N}I^mM/I^{m+n}M$,
which is the Hausdorff completion of $I^mM$;

\item $M/I^mM\to \widehat{M}/I^m\widehat{M}$ is an isomorphism.

\end{itemize}
\symn{$\widehat{(-)}$}
\end{lemma}

\begin{proof}
This follows from 
\cite[Chapter III, \S2 n$^{\circ}$ 11, Proposition 14 \& Corollary 1]{bourbaki-algcomm}
when $M$ is finitely generated. See 
\cite[Chapter 0, 
Corollary 7.2.9 \& Propositions 7.2.15 \& 7.2.16]{fujiwara-kato} 
for general $M$.
\end{proof}

\begin{nota}
\label{not:on-adic-quotient}
If $A$ is an adic ring and $T=(T_i)_i$ is a family of 
indeterminates, we denote by $A\langle T \rangle$
the algebra of restricted power series in $T$ with 
coefficients in $A$, i.e., 
the $I$-adic completion of $A[T]$ for 
an ideal of definition $I\subset A$. 
Unless otherwise stated, given an ideal $J\subset A\langle T\rangle$, 
we denote by 
$A\langle T\rangle/J$ the $I$-adically complete quotient, i.e.,
the quotient of $A\langle T\rangle$ by the closure of the ideal $J$.
\symn{$\langle T\rangle$}
\end{nota}

Unless otherwise stated, formal schemes are always assumed to be adic 
of finite ideal type in the sense of 
\cite[Chapter I, Definitions 1.1.14 \& 1.1.16]{fujiwara-kato}.
Thus, a formal scheme 
$\mathcal{X}=(|\mathcal{X}|,\mathcal{O}_{\mathcal{X}})$ 
is a ringed space
with is locally isomorphic to $\Spf(A)$, where $A$ is an adic 
ring (of finite ideal type, as always). Morphisms of formal schemes 
are assumed to be adic, i.e., are locally of the form 
$\Spf(B) \to \Spf(A)$, with $A \to B$ an adic morphism.
\symn{$\Spf$}

Let $\mathcal{X}$ be a formal scheme. An ideal
$\mathcal{I}\subset \mathcal{O}_{\mathcal{X}}$ is said to 
be an ideal of definition if locally it is of the form
$I\mathcal{O}_{\Spf(A)}$ where $A$ is an adic ring and 
$I\subset A$ an ideal of definition. 
In this case, the ringed space 
$(|\mathcal{X}|,\mathcal{O}_{\mathcal{X}}/\mathcal{I})$
is an ordinary scheme which we simply denote by 
$\mathcal{X}/\mathcal{I}$. By
\cite[Chapter I, Corollary 3.7.12]{fujiwara-kato},
every quasi-compact and quasi-separated formal scheme
admits an ideal of definition which we may assume 
to be finitely generated.
\ncn{formal schemes!ideal of definition}
\symn{$(-)/(-)$}

\begin{dfn}
\label{dfn:monogenic-principal}
Let $A$ be an adic ring. We say that $A$ is of 
principal ideal type if it admits an ideal of definition which 
is principal (i.e., generated by a nonzero divisor). 
We will say that $A$ is of monogenic ideal type if it admits 
an ideal of definition which is monogenic (i.e., generated by 
one element). Similarly, we say that a formal scheme is 
of principal ideal type (resp. of monogenic ideal type) if 
it admits an ideal of definition which is principal (resp. monogenic).
There are also obvious local versions of these notions where we
only require that an ideal of definition of 
a specific type exists locally.
\ncn{adic rings!principal ideal type}
\ncn{adic rings!monogenic ideal type}
\ncn{formal schemes!monogenic ideal type}
\ncn{formal schemes!principal ideal type}
\end{dfn}

\begin{rmk}
\label{rmk:monogenic-vs-principal}
Let $A$ be an adic ring of monogenic ideal type and 
$\pi\in A$ a generator of an ideal of definition of $A$. 
Then $A$ is of principal ideal type if and only if 
$A$ is $\pi$-torsion-free.
\end{rmk}

\begin{nota}
\label{not:categ-formal-scheme}
We denote by \sym{$\FSch$} the category of formal schemes and by 
\sym{$\FSch^{\qcqs}$} its full subcategory spanned by quasi-compact and 
quasi-separated formal schemes (in the sense of
\cite[Chapter I, Definitions 1.6.1 \& 1.6.5]{fujiwara-kato}). 
Note that the category \sym{$\Sch$} (resp. \sym{$\Sch^{\qcqs}$}) of 
schemes (resp. of quasi-compact and quasi-separated schemes)
can be identified with the full subcategory of $\FSch$
(resp. $\FSch^{\qcqs}$) spanned by those formal schemes 
for which $(0)$ is an ideal of definition.
\end{nota}

\begin{nota}
\label{not:special-fiber}
The inclusion of the category of reduced schemes into $\FSch$
admits a right adjoint which we denote by 
$\mathcal{X}\mapsto \mathcal{X}_{\sigma}$.
It commutes with gluing along open immersions 
and satisfies $\mathcal{X}_{\sigma}=(\mathcal{X}/\mathcal{I})_{\red}$
whenever $\mathcal{X}$ admits an ideal of definition
$\mathcal{I}\subset \mathcal{O}_{\mathcal{X}}$. The scheme 
$\mathcal{X}_{\sigma}$ is called the special fiber
of $\mathcal{X}$.
\symn{$(-)_{\sigma}$}
\end{nota}

The following notions agree with the ones introduced in
\cite[Chapter I, Definitions 4.2.2 \& 4.3.4 \& 4.7.1 \& 4.8.12
\& 5.3.10 \& 5.3.16]{fujiwara-kato}.

\begin{dfn}
\label{dfn:smooth-etale-proper-finite-}
Let $f:\mathcal{Y} \to \mathcal{X}$ be a morphism of 
formal schemes.
\begin{enumerate}

\item[(1)] We say that $f$ is a closed immersion (resp. finite, proper)
if locally on $\mathcal{X}$ there is an ideal of definition 
$\mathcal{I}\subset \mathcal{O}_X$ such that the induced morphism
of schemes $\mathcal{Y}/\mathcal{I}\to \mathcal{X}/\mathcal{I}$
is a closed immersion (resp. finite, proper).
\ncn{formal schemes!closed immersion}
\ncn{formal schemes!proper morphism}
\ncn{formal schemes!finite morphism}

\item[(2)] We say that $f$ is an open immersion (resp. adically flat, 
\'etale, smooth)
if locally on $\mathcal{X}$ there is an ideal of definition 
$\mathcal{I}\subset \mathcal{O}_X$ such that the induced morphism
of schemes $\mathcal{Y}/\mathcal{I}^n\to \mathcal{X}/\mathcal{I}^n$
is an open immersion (resp. flat, \'etale, 
smooth) for every $n\in \N$.
\ncn{formal schemes!open immersion}
\ncn{formal schemes!adically flat morphism}
\ncn{formal schemes!smooth morphism}
\ncn{formal schemes!\'etale morphism}

\end{enumerate}
\end{dfn}

Let $\mathcal{X}$ be a formal scheme. 
An ideal $\mathcal{J}\subset \mathcal{O}_{\mathcal{X}}$ is 
said to be admissible if, locally on $\mathcal{X}$, it is 
finitely generated and contains an ideal of definition.
An admissible blowup of $\mathcal{X}$ is the blowup of 
an admissible ideal. For more details, see
\cite[Chapter II, \S1.1]{fujiwara-kato}. 
We recall here that the composition 
$\mathcal{X}''\to \mathcal{X}$
of two admissible blowups
$\mathcal{X}''\to \mathcal{X}'$ and $\mathcal{X}'\to \mathcal{X}$
is itself an admissible blowup if $\mathcal{X}$ is 
quasi-compact and quasi-separated.
(This is \cite[Chapter II, Proposition 1.1.10]{fujiwara-kato}.)
We denote by $\mathfrak{B}(\mathcal{X})$ the 
category of admissible blowups and morphisms of formal
$\mathcal{X}$-schemes.
If $\mathcal{X}$ is quasi-compact and quasi-separated, then 
$\mathfrak{B}(\mathcal{X})$ 
is cofiltered 
(by \cite[Chapter II, Proposition 1.3.1]{fujiwara-kato})
and if $\mathcal{U}\to \mathcal{X}$ 
is a quasi-compact open immersion, then the obvious functor
$\mathfrak{B}(\mathcal{X})\to \mathfrak{B}(\mathcal{U})$ 
is surjective (by 
\cite[Chapter II, Proposition 1.1.9]{fujiwara-kato}).
\ncn{formal schemes!admissible blowup}
\symn{$\mathfrak{B}$}

\begin{nota}
\label{not:rigid-analytic-spaces}
(See \cite[Chapter II, \S2]{fujiwara-kato})
We denote by \sym{$\RigSpc^{\qcqs}$}
the {$1$-categorical} localisation of the category $\FSch^{\qcqs}$ with respect to 
admissible blowups. More concretely, there is a functor 
$(-)^{\rig}:\FSch^{\qcqs}\to \RigSpc^{\qcqs}$
which is a bijection on objects and, given two 
quasi-compact and quasi-separated formal schemes 
$\mathcal{X}$ and $\mathcal{Y}$, we have
\begin{equation}
\label{eqn-not:rigid-analytic-space}
\Hom_{\RigSpc^{\qcqs}}(\mathcal{Y}^{\rig},\mathcal{X}^{\rig})=
\underset{\mathcal{Y}'\to \mathcal{Y} 
\,\in\, \mathfrak{B}(\mathcal{Y})}{\colim}\,
\Hom_{\FSch^{\qcqs}}(\mathcal{Y}',\mathcal{X}).
\end{equation}
The objects of $\RigSpc^{\qcqs}$ are the quasi-compact and 
quasi-separated \nc{rigid analytic spaces} (according to 
\cite[Chapter II, Definitions 2.1.1 \& 2.1.2]{fujiwara-kato}).
If $\mathcal{X}$ is a quasi-compact and quasi-separated formal scheme,
$\mathcal{X}^{\rig}$ is called the Raynaud generic fiber
(or simply the generic fiber) of $\mathcal{X}$. For this reason, 
we sometimes write ``$\mathcal{X}_{\eta}$'' instead of 
``$\mathcal{X}^{\rig}$''.
A map in $\RigSpc^{\qcqs}$ is an open immersion if it is 
isomorphic to the generic fiber of an open immersion in
$\FSch^{\qcqs}$. General rigid analytic spaces are 
obtained by gluing along open immersions from objects in 
$\RigSpc^{\qcqs}$ as in
\cite[Chapter II, \S 2.2.(c)]{fujiwara-kato}.
The resulting category is denoted by \sym{$\RigSpc$} and its 
objects are the rigid analytic spaces. There is also a 
generic fiber functor
$(-)^{\rig}:\FSch\to \RigSpc$
extending the one on quasi-compact and quasi-separated formal schemes.
\symn{$(-)_\eta$}
\end{nota}

{
\begin{nota}
\label{nota:Mdl}
Let $X$ be a rigid analytic space. A formal model for 
$X$ is a formal scheme $\mathcal{X}$ endowed with an 
isomorphism $X\simeq \mathcal{X}^{\rig}$
(see \cite[Chapter II, Definition 2.1.7]{fujiwara-kato}).
Formal models of $X$ form a category which we denote by $\Mdl(X)$. 
When $X$ is quasi-compact and quasi-separated, 
$\Mdl(X)$ is cofiltered by \cite[Chapter II,
Proposition 2.1.10]{fujiwara-kato}.
Similarly, given a morphism $f:Y \to X$ of rigid analytic spaces, 
we have a category $\Mdl(f)$ of formal models of $f$
whose objects are morphisms of formal schemes 
$\phi:\mathcal{Y} \to \mathcal{X}$ 
together with an isomorphism $f\simeq \phi^{\rig}$
in $\RigSpc^{\Delta^1}$. When $X$ and $Y$ are quasi-compact and 
quasi-separated, the category $\Mdl(f)$ is cofiltered.
\symn{$\Mdl$}
\end{nota}
}
\ncn{rigid analytic spaces!formal model}

\begin{rmk}
\label{rmk:essential-form-sch}
If $\mathcal{X}$ is a formal scheme and $\mathcal{I}$ is an ideal 
of definition of $\mathcal{X}$, then the admissible blowup 
of $\mathcal{I}$ is locally of principal ideal type
(in the sense of Definition
\ref{dfn:monogenic-principal}). 
Therefore, every quasi-compact and quasi-separated 
rigid analytic space $X$ admits formal models which are locally 
of principal ideal type and 
these form a cofinal subcategory of $\Mdl(X)$ 
which we denote by $\Mdl'(X)$.
\symn{$\Mdl'$}
\end{rmk}

\begin{nota}
\label{not:visualisation-}
Let $X$ be a quasi-compact and quasi-separated rigid analytic space.
We define a locally ringed space $(|X|,\mathcal{O}^+_X)$ by
$$(|X|,\mathcal{O}_X^+)=
\lim_{\mathcal{X}\in \Mdl(X)}
(|\mathcal{X}|,\mathcal{O}_{\mathcal{X}}).$$
If $\mathcal{X}_0$ is a formal model of $X$ and $\mathcal{I}\subset 
\mathcal{O}_{\mathcal{X}_0}$ is an ideal of definition, 
then $\mathcal{I}\mathcal{O}^+_X$ is an invertible ideal in 
$\mathcal{O}^+_X$. We set $\mathcal{O}_X=\bigcup_{n\geq 0}(\mathcal{I}\mathcal{O}^+_X)^{-n}$. Then $\mathcal{O}_X$ is a sheaf of rings 
which does not depend on $\mathcal{I}$ and which contains
$\mathcal{O}_X^+$. By gluing along open immersions, the assignment 
$X\mapsto (|X|,\mathcal{O}_X,\mathcal{O}^+_X)$ 
can be extended to any rigid analytic space $X$. For more details, 
we refer the reader to \cite[Chapter II, \S 3]{fujiwara-kato}.
We say that $|X|$ is the topological space associated to 
$X$, that $\mathcal{O}_X$ is the structure sheaf of $X$, and that 
$\mathcal{O}_X^+$ is the integral structure sheaf of $X$. 
\end{nota}

\begin{rmk}
\label{rmk:on-top-of-X}
Let $X$ be a rigid analytic space. The topological space 
$|X|$ is valuative, in the sense of 
\cite[Chapter 0, Definition 2.3.1]{fujiwara-kato},
and spectral if $X$ is quasi-compact and quasi-separated.
The Krull dimension (or simply the dimension) of $X$ is defined to be 
the Krull dimension of $|X|$, i.e., the supremum of the 
lengths of chains of irreducible closed subsets of $|X|$.
\end{rmk}

\begin{nota}
\label{not:point-residue-field-}
Let $X$ be a rigid analytic space and $x\in |X|$ a point.
By \cite[Chapter II, Proposition 3.2.6]{fujiwara-kato}, 
the local ring $\mathcal{O}_{X,x}^+$ is a prevaluative 
ring. (Here we use the terminology of 
\cite[D\'efinition 1.9.1]{abbes}.)
More precisely, there is a nonzero divisor $a\in \mathcal{O}_{X,x}^+$
with the following properties:
\begin{itemize}

\item every finitely generated ideal of $\mathcal{O}_{X,x}^+$
containing a power of $a$ is principal;

\item $\mathcal{O}_{X,x}^+[a^{-1}]=\mathcal{O}_{X,x}$;

\item $\mathfrak{m}_{X,x}=\bigcap_{n\in \N}
a^n\mathcal{O}_{X,x}^+$ where $\mathfrak{m}_{X,x}$ is the maximal ideal
of $\mathcal{O}_{X,x}$;

\item $\mathcal{O}_{X,x}^+/\mathfrak{m}_x$ is a 
valuation ring of the residue field 
$\mathcal{O}_{X,x}/\mathfrak{m}_x$. We denote by 
$\Gamma_x$ its value group (denoted multiplicatively).

\end{itemize}
We let $\kappa^+(x)$ be the $a$-adic completion of 
$\mathcal{O}_{X,x}^+$, $\kappa(x)$ its fraction 
field and $\widetilde{\kappa}(x)$ the residue field of $\kappa^+(x)$.
We also let $\kappa^{\circ}(x)\subset \kappa(x)$ be 
the subring of power bounded elements. Then 
$\kappa^{\circ}(x)$ is the unique height $1$
valuation ring containing $\kappa^+(x)$. 
Moreover, $\kappa(x)$ is a non-Archimedean complete field 
for the norm induced by $\kappa^{\circ}(x)$.
\symn{$\kappa^+$}
\symn{$\kappa$}
\symn{$\widetilde{\kappa}$}
\symn{$\kappa^{\circ}$}
\end{nota}

\begin{dfn}
\label{dfn:smooth-etale-proper-finite-rig}
Let $f:Y \to X$ be a morphism of rigid analytic spaces.
\begin{enumerate}

\item[(1)] We say that $f$ is a closed immersion (resp. finite, proper)
if, locally on $X$, $f$ admits a formal model which is a
closed immersion (resp. finite, proper). 
\ncn{rigid analytic spaces!closed immersion}
\ncn{rigid analytic spaces!finite morphism}
\ncn{rigid analytic spaces!proper morphism}

\item[(2)] We say that 
$f$ is a locally closed immersion if it can be written
as the composition of a closed immersion $Y\to U$ followed by an 
open immersion $U\to X$.
\ncn{rigid analytic spaces!locally closed immersion}

\item[(3)] We say that $f$ is \'etale (resp. smooth)
with good reduction if, locally on $X$, 
$f$ admits a formal model which is \'etale 
(resp. smooth).
\ncn{rigid analytic spaces!\'etale morphism with good reduction}
\ncn{rigid analytic spaces!smooth morphism with good reduction}

\end{enumerate}
\end{dfn}

We next discuss the analytification functor 
following \cite[Chapter II, \S 9.1]{fujiwara-kato}.

\begin{cons}
\label{cons:analytification-}
Let $A$ be an adic ring, $I\subset A$ an ideal of definition, 
$U=\Spec(A)\smallsetminus \Spec(A/I)$ and
$S=\Spf(A)^{\rig}$.
There exists an \nc{analytification functor} 
$$(-)^{\an}:\Sch^{\lft}/U\to \RigSpc/S,$$
where $\Sch^{\lft}/U$ is the category
of $U$-schemes locally of finite type. 
This functor is uniquely determined by the following two properties.
\symn{$(-)^{\an}$}
\symn{$\Sch^{\lft}$}
\begin{enumerate}

\item[(1)] It is compatible with gluing along open immersions.

\item[(2)] For a separated finite type $U$-scheme $X$ with an open 
immersion $X\to \overline{X}$ into a proper $A$-scheme, and 
complement $Y=\overline{X}\smallsetminus X$, we have 
\begin{equation}
\label{eq-cons:analytification-}
X^{\an}=(\widehat{\overline{X}})^{\rig}\smallsetminus 
(\widehat{Y})^{\rig}
\end{equation}
where, for an $A$-scheme $W$, $\widehat{W}=\colim_n\,W\otimes_A A/I^n$ 
is the $I$-adic completion of $W$.

\end{enumerate}
In the second property, one may replace 
$Y$ with the closure in $\overline{X}$ 
of $\overline{X}\times_A U\smallsetminus X$.
That \eqref{eq-cons:analytification-} is independent of the 
choice of the compactification, follows from
\cite[Chapter II, Propositions 9.1.5 \& 
9.1.9]{fujiwara-kato}.\footnote{Proposition 
9.1.9 of loc.~cit. is stated under the assumption that $A$ is 
topologically universally rigid-noetherian, but this assumption is 
unnecessary.}
\end{cons}

\subsection{Relation with adic spaces}

$\empty$

\smallskip

\label{subsect:rel-adic-sp}

Recall from \cite[page 37]{huber}
that a \nc{Tate ring} is a topological ring $A$ admitting 
a topologically nilpotent unit and 
an open subring $A_0\subset A$ which is adic.
(Here, by convention, Tate rings are assumed complete.)
The ring $A_0$ is called a ring of definition. If $\pi\in A$ is 
a topologically nilpotent unit contained in $A_0$, then 
the topology of $A_0$ is $\pi$-adic, i.e., 
the $\pi^nA_0$ form a fundamental system of open neighbourhoods of $0$.
A morphism of Tate rings $f:A \to B$ is a continuous morphism of rings 
for which there exists rings of definitions 
$A_0\subset A$ and $B_0\subset B$ with $f(A_0)\subset B_0$.

\begin{nota}
\label{not:A-circ-circ}
Given a Tate ring $A$, we denote by $A^{\circ}\subset A$ 
the subring of 
power bounded elements and $A^{\circ\circ}\subset A^{\circ}$ 
the ideal of topologically nilpotent elements. We say that $A$ is 
uniform if $A^{\circ}$ is bounded (which is equivalent to 
ask that $A^{\circ}$ is a ring of definition).
\symn{$(-)^{\circ}$}
\symn{$(-)^{\circ\circ}$}
\end{nota}

A \nc{Tate affinoid ring} $A$
is a pair $(A^{\pm},A^+)$ where $A^{\pm}$ is a Tate ring and
$A^+$ is an integrally closed open subring of $A^{\pm}$ contained in 
$(A^{\pm})^{\circ}$.

\begin{cons}
\label{cons:huber-pair-and-back}
$\empty$

\begin{enumerate}

\item[(1)] Let $A$ be an adic ring of principal ideal type and 
$\pi\in A$ a generator of an ideal of definition. We associate to 
$A$ a Tate 
affinoid ring $A^{\natural}=(A^{\natural\pm},A^{\natural+})$ where 
$A^{\natural\pm}=A[\pi^{-1}]$ and $A^{\natural +}$ is the integral
closure of $A$ in $A[\pi^{-1}]$.
\symn{$(-)^{\natural(\pm)}$}

\item[(2)] The functor $A \mapsto A^{\natural}$, from adic 
rings of principal ideal type to Tate affinoid rings,
admits a ind-right adjoint. The latter associates to a Tate affinoid
ring $R=(R^{\pm},R^+)$ the ind-adic ring $R_{\natural}$ consisting 
of those rings of definition of $R^{\pm}$ contained in $R^+$.
\symn{$(-)_{\natural}$}
\end{enumerate}
\end{cons}

\begin{rmk}
\label{rmk:natural-no-ind}
When the Tate affinoid ring $R$ is uniform, 
then the associated ind-adic ring $R_{\natural}$ is isomorphic to an 
adic ring. In fact, we have $R_{\natural}=R^+$.
\end{rmk}

\begin{lemma}
\label{lem:natural-full-faith}
The functor $R\mapsto R_{\natural}$, from the category of Tate 
affinoid rings to the category of ind-adic rings of 
principal ideal type, is fully faithful.
\end{lemma}

\begin{proof}
Indeed, let $R$ and $R'$ be two Tate affinoid rings
and $f:R_{\natural} \to R'_{\natural}$ a morphism of 
ind-adic rings. There exists rings of definition 
$R_0\subset R^{\pm}$ and $R'_0\subset R'^{\pm}$ 
contained in $R^+$ and $R'^+$ such that 
$f$ restricts to a morphism of adic rings
$f_0:R_0\to R'_0$. Then $f_0$ induces a morphism of Tate rings 
$f^{\pm}:R^{\pm} \to R'^{\pm}$. Since $f_0$ is the restriction of $f$, 
for every ring of definition 
$R_1\subset R^{\pm}$ containing $R_0$ and contained in $R^+$, 
there exists a ring of definition 
$R'_1\subset R'^{\pm}$ contained in $R'^+$ and 
a morphism $f_1:R_1\to R'_1$ extending $f_0$.
This shows that $f^{\pm}$ maps $R^+$ into $R'^+$ as needed.
\end{proof}

Given a Tate affinoid ring $A=(A^{\pm},A^+)$, we denote by 
$\Spa(A)=(|\Spa(A)|,\mathcal{O}_{\Spa(A)},\mathcal{O}_{\Spa(A)}^+)$ 
the preadic space associated to $A$ as in 
\cite[pages 38--39]{huber}.
In general, $\mathcal{O}^+\subset \mathcal{O}$ are
presheaves of rings on the topological space $|\Spa(A)|$
which might fail to be sheaves.
\symn{$\Spa$}

\begin{prop}
\label{prop:compar-Spa-Spf}
$\empty$

\begin{enumerate}

\item[(1)] Let $A$ be an adic ring of principal ideal type. There is a homeomorphism 
$|\Spf(A)^{\rig}|\simeq |\Spa(A^{\natural})|$ modulo which 
$\mathcal{O}^+_{\Spf(A)^{\rig}}$ 
(resp. $\mathcal{O}_{\Spf(A)^{\rig}}$)
is isomorphic to the sheafification of 
$\mathcal{O}^+_{\Spa(A^{\natural})}$
(resp. $\mathcal{O}_{\Spa(A^{\natural})}$).

\item[(2)] Let $R$ be an affinoid ring. There exists a homeomorphism
$|\Spa(R)|\simeq \lim |\Spf(R_{\natural})^{\rig}|$ modulo which
$\mathcal{O}^+_{\lim \Spf(R_{\natural})^{\rig}}$
(resp. $\mathcal{O}_{\lim \Spf(R_{\natural})^{\rig}}$)
is isomorphic to the sheafification of 
$\mathcal{O}^+_{\Spa(R)}$ (resp. $\mathcal{O}_{\Spa(R)}$).

\end{enumerate}
\end{prop}

\begin{proof}
A point $x\in |\Spf(A)^{\rig}|$ determines a morphism of 
adic rings $A \to \kappa^+(x)$, and hence a continuous 
valuation $v_x:A \to \Gamma_x\cup\{0\}$
landing in $\Gamma^+_x\cup\{0\}$. (Here $\Gamma^+_x\subset \Gamma$ 
denotes the submonoid defined by the inequality $\leq 1$.)
Since the image of $\pi$ in $\kappa^+(x)$ is nonzero, 
$v_x$ extends uniquely to a continuous valuation $v_x:A^{\natural\pm}
\to \Gamma_x\cup\{0\}$. Moreover, $v_x$ maps $A^{\natural+}$
into $\Gamma^+_x\cup\{0\}$ since $A^{\natural+}$ is integral over $A$.
Therefore, $v_x$ belongs to $\Spa(A^{\natural})$. It is easy to 
see that $x\mapsto v_x$ is a bijection, which is continuous and open.
More precisely, given elements $a_0,\ldots, a_n$ in $A$ 
generating an admissible ideal of $A$, 
the open subset 
$|\Spf(A\langle \frac{a_1}{a_0},\ldots, \frac{a_n}{a_0} \rangle)^{\rig}| \subset |\Spf(A)^{\rig}|$
is mapped bijectively to the rational subset 
$|\Spa(A^{\natural}\langle \frac{a_1}{a_0},\ldots, 
\frac{a_n}{a_0}\rangle)|\subset |\Spa(A^{\natural})|$.
This also shows that $\mathcal{O}_{\Spf(A)^{\rig}}$
is the sheafification of $\mathcal{O}_{\Spa(A^{\natural})}$.

Assertion (2) can be deduced from assertion (1) and the fact that 
the counit map $(R_{\natural})^{\natural} \to R$ identifies 
the Tate affinoid ring $R$ with the colimit of the ind-Tate affinoid
ring $(R_{\natural})^{\natural}$.
\end{proof}

\begin{dfn}
\label{dfn:uniform-adic-space}
A uniform adic space is a triple 
$X=(|X|,\mathcal{O}_X,\mathcal{O}^+_X)$, 
consisting of a topological space $|X|$ and sheaves of rings
$\mathcal{O}_X^+\subset \mathcal{O}_X$, which is locally 
isomorphic to $\Spa(A)$, where $A$ is a stably uniform 
Tate affinoid ring
in the sense of \cite[pages 30--31]{buzz-ver}.
(This is reasonable since by \cite[Theorem 7]{buzz-ver}
every stably uniform Tate affinoid ring is sheafy.)
\ncn{adic spaces!uniform}
\end{dfn}

\begin{cor}
\label{cor:from-Huber-to-FK}
Let \sym{$\Adic$} be the category of uniform
adic spaces. Then there exists a fully faithful embedding 
$\Adic \to \RigSpc$ which is compatible with gluing along open 
immersions and which sends $\Spa(R)$ to $\Spf(R^+)^{\rig}$. 
\end{cor}

\begin{proof}
It suffices to treat the affinoid case; the general 
case follows then by gluing along open immersions.
Given two stably uniform Tate affinoid rings $A$ and $B$, 
the fact that $A$ is sheafy implies that there is a bijection 
$\Hom(A,B) \simeq \Hom(\Spa(B),\Spa(A))$.
It follows from Remark \ref{rmk:natural-no-ind} that 
there is a functor 
$\Spa(A) \mapsto \Spf(A^+)^{\rig}$, from affinoid uniform 
adic spaces to rigid analytic spaces, and it remains to show 
that the map 
$$\Hom(A^+,B^+)\to \Hom(\Spf(B^+)^{\rig},\Spf(A^+)^{\rig}),$$
with $A$ and $B$ as above, is a bijection.
An element of the right-hand side can be represented by 
a morphism $\mathcal{Y}\to \Spf(A^+)$, where 
$\mathcal{Y}\to \Spf(B^+)$ is an admissible blowup. 
We may assume that $\mathcal{O}_{\mathcal{Y}}$ is $\pi$-torsion-free, with 
$\pi$ a generator of an ideal of definition in $B^+$.
We claim that $\mathcal{O}(\mathcal{Y})=B^+$ which 
implies that $\mathcal{Y}\to \Spf(A^+)$ factors uniquely through 
$\Spf(B^+)$, finishing the proof.

Let $(\mathcal{Y}_i)_i$ be an affine open covering of 
$\mathcal{Y}$ and set $\mathcal{Y}_{ij}=\mathcal{Y}_i\cap \mathcal{Y}_j$.
Let $B_i$ and $B_{ij}$ be the 
Tate affinoid rings associated to the adic rings 
$\mathcal{O}(\mathcal{Y}_i)$ and $\mathcal{O}(\mathcal{Y}_{ij})$
respectively. Then $(\Spa(B_i))_i$ is an open covering of 
$\Spa(B)$, and $\Spa(B_{ij})=\Spa(B_i) \cap \Spa(B_j)$. 
Since $B$ is sheafy, we deduce that 
$B^+$ is the equaliser of the usual pair of arrows
$\prod_i B_i^+ \rightrightarrows
\prod_{ij}B_{ij}^+$. Since 
$\mathcal{O}_{\mathcal{Y}}$ is $\pi$-torsion-free, 
we have inclusions $\mathcal{O}(\mathcal{Y}_i)\subset 
B_i^+$ and $\mathcal{O}(\mathcal{Y}_{ij})\subset B^+_{ij}$. 
This proves that $\mathcal{O}(\mathcal{Y})$, which is the 
equaliser of $\prod_i \mathcal{O}(\mathcal{Y}_i) \rightrightarrows
\prod_{ij}\mathcal{O}(\mathcal{Y}_{ij})$, 
is contained in $B^+$ as needed.
\end{proof}

\subsection{\'Etale and smooth morphisms}

$\empty$

\smallskip

\label{subsect:etale-smooth-rig}

In Definition \ref{dfn:smooth-etale-proper-finite-rig}
we introduced the classes of \'etale and smooth morphisms
with good reduction. These classes are too small, and we  
need to enlarge them to get the correct notions of 
\'etaleness and smoothness in rigid analytic geometry.
First, we introduce a notation.

\begin{nota}
\label{not:satura-adic}
Let $A$ be an adic ring and $J\subset A$ an ideal. We denote by
$J^{\sat}$ the ideal of $A$ consisting of those elements $a\in A$ 
for which there exists an ideal of definition $I\subset A$ such that
$aI\subset J$. The ideal $J^{\sat}$ is called the 
saturation of $J$.
\ncn{adic rings!saturated ideal}
\symn{$(-)^{\sat}$}
\end{nota}

We say that $J$ is saturated if $J=J^{\sat}$. The saturation of 
an ideal is a saturated ideal.

\begin{rmk}
\label{rmk:on-sat-ideal}
If $A$ is an adic ring of principal ideal type and $J\subset A$ a
saturated ideal, then $J$ is closed and $A/J$ is also of 
principal ideal type. Moreover, for a closed ideal $J\subset A$, the 
quotient $A/J$ is of principal ideal type if and only if $J$ is 
saturated.
\end{rmk}

Our definition of \'etaleness uses the Jacobian matrix. Compare with 
\cite[Definition 1.3.1]{fujiwara-tubular}.

\begin{dfn}
\label{dfn:new-etale-rig}
$\empty$

\begin{enumerate}

\item[(1)] Let $A$ be an adic ring and $B$ an adic $A$-algebra.
We say that $B$ is rig-\'etale over $A$ if there exists a 
presentation $B\simeq A\langle t_1,\ldots, t_n\rangle/J$
and elements $f_1,\ldots, f_n\in J$ such that 
$(f_1,\ldots, f_n)^{\sat}=J^{\sat}$ 
and the determinant of the Jacobian matrix 
$\det(\partial f_i/\partial t_j)$ generates an open ideal in $B$.
\ncn{adic rings!rig-\'etale over}

\item[(2)] A morphism $\mathcal{Y}\to \mathcal{X}$ of formal 
schemes is said to be rig-\'etale if, locally for the 
rig topology on $\mathcal{X}$ and $\mathcal{Y}$ (see
Definition \ref{dfn:rig-Zar-Nis} below), it is isomorphic 
to $\Spf(B) \to \Spf(A)$ with $B$ rig-\'etale over $A$.
(When $\mathcal{X}$ and $\mathcal{Y}$ are quasi-compact, 
this simply means that after replacing $\mathcal{X}$ and 
$\mathcal{Y}$ by admissible blowups, the resulting morphism
is locally isomorphic to 
$\Spf(B) \to \Spf(A)$ with $B$ rig-\'etale over $A$.)
\ncn{formal schemes!rig-\'etale morphism}

\item[(3)] A morphism of rigid analytic spaces $Y\to X$ is said 
to be \'etale if, locally on $X$ and $Y$, it admits formal
models which are rig-\'etale. 
\ncn{rigid analytic spaces!\'etale morphism}

\end{enumerate}
\end{dfn}

\begin{rmk}
\label{rmk:etale-abstract-defn}
If the rigid analytic space $X$ is assumed to be universally 
noetherian (in the sense of 
\cite[Chapter II, Definition 2.2.23]{fujiwara-kato}), then 
a morphism $f:Y\to X$ is \'etale if and only if it is flat and 
neat (i.e., $\Omega_f=0$). This follows from 
\cite[Propositions 1.7.1 and 1.7.5]{huber}
together with
\cite[Chapter II, Theorem A.5.2]{fujiwara-kato}.
See also
\cite[Proposition 5.1.6]{fujiwara-tubular}
which is proven under more restrictive assumptions.
\end{rmk}

\begin{rmk}
\label{rmk:rig-et-sat}
Let $A$ be an adic ring and $B$ a rig-\'etale adic 
$A$-algebra given by $A\langle t_1,\ldots, t_n\rangle/J$
with $J$ containing $f_1,\ldots, f_n$ as in Definition 
\ref{dfn:new-etale-rig}. Consider the adic $A$-algebras
$$B'=A\langle t_1,\ldots, t_n\rangle/(f_1,\ldots, f_n)
\quad \text{and} \quad
B''=A\langle t_1,\ldots, t_n\rangle/(f_1,\ldots, f_n)^{\sat}.$$
We have surjective maps
$B'\to B\to B''$ inducing isomorphisms
$\Spf(B'')^{\rig}\simeq \Spf(B)^{\rig}\simeq \Spf(B')^{\rig}$.
Moreover, $B'$ and $B''$ are rig-\'etale over $A$.
The case of $B''$ is clear. For $B'$, we need to prove the following statement. Let $C$ be an adic ring and $c\in C$ an element.
Then $c$ generates an open ideal in $C$ if and only if
it generates an open ideal in $C/(0)^{\sat}$. Indeed, let 
$I$ be an ideal of definition and assume that 
$I\subset (c)+(0)^{\sat}$. We need to show that a power of $I$ is contained in $(c)$. Since $I$ is finitely generated, we may 
find elements $v_1,\ldots, v_m$ in $(0)^{\sat}$ such that 
$I\subset (c)+(v_1,\ldots, v_m)$.
Let $r$ be an integer such that $v_iI^r=0$ for all
$1\leq i \leq m$. Then clearly $I^{r+1}
\subset cI^r \subset (c)$ as needed.
\end{rmk}

\begin{lemma}
\label{lem:rigid-of-rig-etal}
Let $A$ be an adic ring of monogenic ideal type and 
$\pi\in A$ a generator of an ideal of definition of $A$.
Let $B$ be a rig-\'etale $A$-algebra. Then there exists an 
integer $N\in \N$ such that for 
every $\pi$-torsion-free adic $A$-algebra $C$, the map 
$\Hom_A(B,C) \to \Hom_{A/\pi^N}(B/\pi^N,C/\pi^N)$
is injective.
\end{lemma}

\begin{proof}
The proof of \cite[Proposition 2.1.1]{fujiwara-tubular}
can be easily adapted to the situation considered in the 
statement. For the reader's convenience we recall the argument.

For $m\in\N$, we set 
$A_m=A/\pi^m$, $B_m=B/\pi^m$ and $C_m=C/\pi^m$.
Since $B$ is rig-\'etale over $A$, there exists an integer 
$c$ such that $\Omega^1_{B_m/A_m}$ is annihilated by 
$\pi^c$ independently of $m$. (Indeed, if $B$ is given as 
in Definition \ref{dfn:new-etale-rig}, 
it suffices to take $c$ so that $\pi^c$ belongs to the ideal
generated by $\det(\partial f_i/\partial t_j)$.)
Now let $f,f':B \to C$ be two morphisms of $A$-algebras inducing 
the same morphism $f_m:B_m \to C_m$ for some $m\geq c+1$. We will
show that $f_{m+1}=f'_{m+1}$, which suffices to conclude using 
induction.

We may consider $f_{2m}$ and 
$f'_{2m}$ as deformations of $f_{m}$. The difference between 
these deformations is classified
by an element $\epsilon\in\Hom(C_m\otimes_{B_m} \Omega^1_{B_m/A_m}, 
\pi^mC/\pi^{2m}C)$. 
Since $\pi$ is a nonzero divisor of $C$ and 
$\Omega^1_{B_m/A_m}$ annihilated by 
$\pi^c$, the image of any $C$-linear morphism
$C_m\otimes_{B_m} \Omega^1_{B_m/A_m} \to 
\pi^mC/\pi^{2m}C$ is contained in $\pi^{2m-c}C/\pi^{2m}C$.
In particular, the map 
$$\Hom(C_m\otimes_{B_m} \Omega^1_{B_m/A_m}, 
\pi^mC/\pi^{2m}C)\to \Hom(C_m\otimes_{B_m} \Omega^1_{B_m/A_m}, 
\pi^mC/\pi^{m+1}C)$$
is identically zero. Since the image of $\epsilon$ by this map
classifies the difference between 
$f_{m+1}$ and $f'_{m+1}$, we get the equality 
$f_{m+1}=f'_{m+1}$.
\end{proof}

\begin{prop}
\label{prop:approx-solution}
Let $A$ be an adic ring of monogenic ideal type and 
$\pi\in A$ a generator of an ideal of definition of $A$.
Let $t=(t_1,\ldots, t_n)$ be
a system of coordinates and $f=(f_1,\ldots, f_n)$ 
an $n$-tuple in $A\langle t \rangle$. 
Let $J\subset A\langle t\rangle$ be an ideal such that 
$(f)\subset J\subset (f)^{\sat}$ and set $B=A\langle t \rangle/J$.
Assume that $\det(\partial f_i /\partial t_j)$ 
generates an open ideal in $B$, so that $B$ is a rig-\'etale 
adic $A$-algebra.
Then, there exists a positive integer $N$
such that for every 
$\pi$-torsion-free adic $A$-algebra $C$ and every integer
$e\geq N$, the map 
\begin{equation}
\label{eq-prop:approx-solution}
\Hom_A(B,C) \to {\rm im}\left\{\Hom_{A/\pi^{2e}}(B/\pi^{2e},C/\pi^{2e})
\to \Hom_{A/\pi^e}(B/\pi^e,C/\pi^e)\right\}
\end{equation}
is bijective.
Moreover, the integer $N$ depends continuously on $f$, i.e., 
we may find one which works for every
$n$-tuple $f'=(f'_1,\ldots, f'_n)$ in $A\langle t\rangle$
which is $\pi$-adically sufficiently close to $f$.
\end{prop}

\begin{proof}
For $N$ sufficiently large, the injectivity of 
\eqref{eq-prop:approx-solution}
follows from Lemma 
\ref{lem:rigid-of-rig-etal}.
The fact that there is an $N$ which works for all $f'$ 
close enough to $f$ follows from the proof of Lemma 
\ref{lem:rigid-of-rig-etal}. (Indeed, the $N$ depends only on the ideal 
generated by $\det(\partial f_i/\partial t_j)$.)

For the surjectivity of \eqref{eq-prop:approx-solution}, 
it is enough to solve the following problem: given an $n$-tuple 
$c_0=(c_{0,1},\ldots, c_{0,n})$ in $C$ such that 
the components of $f(c_0)$ belong to $\pi^{2e}C$, 
find an $n$-tuple $c=(c_1,\ldots, c_n)$ in $C$ such that 
$f(c)=0$ and the components of $c-c_0$ belong to $\pi^eC$. 
(Indeed, since $C$ is $\pi$-torsion-free an $n$-tuple $c$ such that 
$f(c)=0$ determines an $A$-morphism $B\to C$.)

This problem can be solved using the Newton method as 
in the first step of the proof of \cite[Lemme 1.1.52]{ayoub-rig}.
In fact, one can also remark that the argument in loc.~cit. 
is valid more generally for non-Archimedean Banach rings, i.e., complete normed rings with a non-Archimedean norm. In particular, it applies
with ``$A$'', ``$C$'' and ``$R$'' in loc.~cit. replaced with 
$A[\pi^{-1}]$, $B[\pi^{-1}]$ and $C[\pi^{-1}]$
endowed with the natural norms for which $A/(0)^{\sat}$, 
$B/(0)^{\sat}$ and $C=C/(0)^{\sat}$ 
are the unit balls. (More precisely, for $a\in A[\pi^{-1}]$, we set
$\|a\|=e^{-v(a)}$ where $v(a)$ is the largest integer such that 
$a\in \pi^{v(a)}(A/(0)^{\sat})$, and similarly for $B$ and $C$.) 
Since $\pi$ is a nonzero divisor of $C$, 
a solution $c=(c_1,\ldots, c_n)$ in $(C[\pi^{-1}])^n$
of the system of equations $f=0$, close enough to $c_0$, 
determines a solution in $C^n$. 
We may take for $N$ an integer which is larger than 
$\ln(2M^2)$ with $M$ as in 
\cite[page 46]{ayoub-rig}.\footnote{Here and below,
the page references to \cite{ayoub-rig} correspond to the 
published version.} It is clear that $N$ depends 
$\pi$-adically continuously on $f$. 
\end{proof}

\begin{prop}
\label{prop:on-etale-presentation}
Let $A$ be an adic ring of monogenic ideal type and 
$\pi\in A$ a generator of an ideal of definition of $A$. 
Let $B$ be a rig-\'etale adic $A$-algebra admitting a presentation
$B=A\langle t\rangle/(f)^{\sat}$, with $t=(t_1,\ldots, t_n)$
a system of coordinates and $f=(f_1,\ldots, f_n)$ an $n$-tuple
in $A\langle t\rangle$
such that $\det(\partial f_i /\partial t_j)$ generates
an open ideal in $B$.
Then there exists an integer $N$ such that the following holds. 
For every $n$-tuple 
$f'=(f'_1,\ldots, f'_n)$ in $A\langle t\rangle$
such that $f'-f$ belongs $(\pi^NA\langle t\rangle)^n$, 
the adic $A$-algebra $B'=A\langle t\rangle/(f')^{\sat}$
is isomorphic to $B$. Moreover, there is 
an isomorphism $B\simeq B'$ induced by $n$-tuple 
$g=(g_1,\ldots, g_n)$ in $A\langle t\rangle$ such that 
$g-t$ belongs to $(\pi A\langle t\rangle)^n$. 
\end{prop}

\begin{proof}
This follows by applying Proposition
\ref{prop:approx-solution} to the rig-\'etale adic $A$-algebras
$B$ and $B'$.
\end{proof}

\begin{nota}
\label{not:affine-rig-etale-}
Let $A$ be an adic ring. We denote by 
$\mathcal{E}_A$ the category of rig-\'etale $A$-algebras 
and $\mathcal{E}'_A$ its full subcategory 
spanned by those adic $A$-algebras whose zero ideal is saturated.
(Thus, every object $B\in \mathcal{E}'_A$ admits a presentation
$B\simeq A\langle t \rangle/(f)^{\sat}$ with 
$t=(t_1,\ldots, t_n)$ and $f=(f_1,\ldots, f_n)$ such that 
$\det(\partial f_i /\partial t_j)$ generates
an open ideal in $B$.)
The inclusion $\mathcal{E}_A'\to \mathcal{E}_A$
admits a left adjoint given by $B\mapsto B/(0)^{\sat}$.
Given a morphism of adic rings $A_1 \to A_2$,
there is are induced functors 
$\mathcal{E}_{A_1}\to \mathcal{E}_{A_2}$ and
$\mathcal{E}'_{A_1} \to \mathcal{E}'_{A_2}$ given by 
$B\mapsto A_2\,\widehat{\otimes}_{A_1}\,B$ and 
$B\mapsto (A_2\,\widehat{\otimes}_{A_1}\,B)/(0)^{\sat}$
respectively.
\symn{$\mathcal{E}$}
\symn{$\mathcal{E}'$}
\end{nota}

\begin{cor}
\label{cor:proj-limi-cal-E-A}
Let $(A_{\alpha})_{\alpha}$ be a filtered inductive system 
of adic rings of monogenic ideal type with colimit $A$
(in the category of adic rings).
Then the obvious functor 
\begin{equation}
\label{eq-cor:proj-limi-cal-E-A}
\underset{\alpha}{\colim}\,\mathcal{E}'_{A_{\alpha}}
\to \mathcal{E}'_A
\end{equation}
is an equivalence of categories. 
\end{cor}

\begin{proof}
Let $R$ be the colimit of $(A_{\alpha})_{\alpha}$ taken in the category
of discrete rings. We may assume that there is a smallest index 
$o$ and we fix $\pi\in A_o$ generating an ideal of definition of $A_o$.
Then $A=\lim_{n\in \N} R/\pi^nR$, and  
there is a map of rings $R\to A$ with kernel $J=\bigcap_n \pi^nR$ 
and with dense image $\widetilde{R}\subset A$. We split the proof 
into two steps.

\paragraph*{Step 1}
\noindent
First, we prove that 
\eqref{eq-cor:proj-limi-cal-E-A}
is essentially surjective. 
By Proposition 
\ref{prop:on-etale-presentation}, 
an object $B\in \mathcal{E}'_A$ admits a presentation of the form
$B=A\langle t\rangle/(\widetilde{f})^{\sat}$ where 
$t=(t_1,\ldots, t_n)$ is
a system of coordinates and $\widetilde{f}=
(\widetilde{f}_1,\ldots, \widetilde{f}_n)$ an $n$-tuple in
$\widetilde{R}[t]$ such that $\widetilde{g}=
\det(\partial \widetilde{f}_i/\partial t_j)$ generates
an open ideal in $B$.
Using Remark \ref{rmk:rig-et-sat}, 
we can find an integer $N$ and an element 
$\widetilde{h}\in A\langle t\rangle$ such that 
$\pi^N-\widetilde{h}\widetilde{g}$ belongs to the closure of 
the ideal $(\widetilde{f})\subset \widetilde{R}[t]$ in 
$A\langle t\rangle$.
In particular, we may write
$$\pi^N-\widetilde{h}\widetilde{g}=\sum_{i=1}^n 
\widetilde{a}_i \widetilde{f}_i+\widetilde{v}\pi^{N+1}$$
with $\widetilde{v}\in A\langle t \rangle$ and 
$\widetilde{a}_1,\ldots, \widetilde{a}_n\in \widetilde{R}[t]$.
Write $\widetilde{h}=\widetilde{h}_0+\widetilde{h}_1\pi^{N+1}$ 
with $\widetilde{h}_0\in \widetilde{R}[t]$ and $\widetilde{h}_1
\in A\langle t\rangle$. 
Replacing $\widetilde{h}$ by $\widetilde{h}_0$ and $\widetilde{v}$ 
by $\widetilde{v}+\widetilde{h}_1$, we may assume that 
$\widetilde{h}$ belongs to $\widetilde{R}[t]$. 
It follows from Lemma
\ref{lem:I-adic-compl}
that the expression 
$\pi^N-\widetilde{h}\widetilde{g}-\sum_{i=1}^n \widetilde{a}_i
\widetilde{f}_i\in \widetilde{R}[t]$
belongs to $\pi^{N+1}\widetilde{R}[t]$. 
Said differently, we may also assume that 
$\widetilde{v}\in\widetilde{R}[t]$. We now choose a lift
$f=(f_1,\ldots, f_n)$ of 
$\widetilde{f}$ to an $n$-tuple in $R[t]$ 
and set $g=\det(\partial f_i/\partial t_j)$.
We also choose lifts $h, a_1,\ldots, a_n\in R[t]$ of 
$\widetilde{h},\widetilde{a}_1,\ldots, \widetilde{a}_n$.
Since the elements of $J$ are divisible by any power of 
$\pi$, we may also find a lift $v\in R[t]$ of $\widetilde{v}$ 
such that 
$$\pi^N-hg=\sum_{i=1}^n a_if_i+v\pi^{N+1}.$$
For $\alpha$ sufficiently big, the previous equality can be lifted
to an equality
$$\pi^N-h_{\alpha}g_{\alpha}=\sum_{i=1}^n a_{\alpha,i}f_{\alpha,i}+v_{\alpha}\pi^{N+1}$$
in $A_{\alpha}[t]$ with the property that 
$g_{\alpha}=\det(\partial f_{\alpha,i}/\partial t_j)$.
Since $1-v_{\alpha}\pi$ is invertible in $A_{\alpha}\langle t\rangle$,
it follows that 
$B_{\alpha}=A_{\alpha}\langle t\rangle/(f_{\alpha})^{\sat}$
is a rig-\'etale $A_{\alpha}$-algebra. Clearly, 
the functor 
\eqref{eq-cor:proj-limi-cal-E-A}
sends $B_{\alpha}$ to $B$.

\paragraph*{Step 2}
\noindent
We now prove that 
\eqref{eq-cor:proj-limi-cal-E-A} 
is fully faithful. We fix two objects 
$B_o, C_o\in \mathcal{E}'_{A_o}$.
For an index $\alpha$, we set 
$B_{\alpha}=(B_o\widehat{\otimes}_{A_o}
A_{\alpha})/(0)^{\sat}$ and define $C_{\alpha}$ similarly.
We also set $B=(B_o\widehat{\otimes}_{A_o}A)/(0)^{\sat}$ 
and define $C$ similarly.
We want to show that 
$$\underset{\alpha}{\colim}\,
\Hom_{A_{\alpha}}(B_{\alpha},C_{\alpha})\to \Hom_A(B,C)$$
is a bijection. (This is enough since we are free to
change the smallest index $o$.
We also used that the colimit in \eqref{eq-cor:proj-limi-cal-E-A} is filtered in order to describe the hom-set in the domain.) 
The above map can be rewritten as
$$\underset{\alpha}{\colim}\,
\Hom_{A_o}(B_o,C_{\alpha})\to \Hom_{A_o}(B_o,C).$$
Since $C$ and the $C_{\alpha}$'s are $\pi$-torsion-free, 
we may replace $B_o$ by any rig-\'etale $A_o$-algebra
$B'_o$ such that $B_o\simeq B'_o/(0)^{\sat}$. 
By Remark 
\ref{rmk:rig-et-sat},
we may choose $B'_o$ topologically finitely 
presented. We now apply
Proposition \ref{prop:approx-solution}: there exists an 
integer $N$ such that the maps
$$\Hom_{A_o}(B'_o,C_{\alpha}) \to
{\rm im}\left\{\Hom_{A_o/\pi^{2N}}(B'_o/\pi^{2N},C_{\alpha}/\pi^{2N})
\to \Hom_{A_o/\pi^{N}}(B'_o/\pi^{N},C_{\alpha}/\pi^{N})\right\}$$
are bijections and similarly for $C$ (instead of $C_{\alpha}$).
Since filtered colimits commute with taking images, we are left to
show that 
$$\underset{\alpha}{\colim}\,
\Hom_{A_o/\pi^e}(B'_o/\pi^e,C_{\alpha}/\pi^e)\to 
\Hom_{A_o/\pi^e}(B'_o/\pi^e,C/\pi^e)$$
is a bijection for any positive integer $e$.
This is clear since $B'_o/\pi^e$ is a finitely presented 
$A_o/\pi^e$-algebra and $C/\pi^e$ is the colimit of the 
filtered system $(C_{\alpha}/\pi^e)_{\alpha}$.
\end{proof}

For later use, we record the following two results.

\begin{lemma}
\label{lem:sect-open-immer}
Let $e:X'\to X$ be an \'etale morphism of rigid analytic spaces, 
and let $s:X \to X'$ be a section of $e$. Then $s$ is an open 
immersion. 
\end{lemma}

\begin{proof}
The question is local on $X$ and around $s(X)$. Thus, we may assume 
that $X=\Spf(A)^{\rig}$ with $A$ an adic ring of principal 
ideal type, that $X'=\Spf(A')^{\rig}$ with $A'$ a rig-\'etale 
adic $A$-algebra, and that $s$ is induced by a morphism
of $A$-algebras $h:A' \to A$. Fix a generator $\pi$ of 
an ideal of definition of $A$.
By Proposition 
\ref{prop:on-etale-presentation}, 
we may assume that 
$A'=A\langle t\rangle/(f)^{\sat}$ with $t=(t_1,\ldots, t_n)$
a system of coordinates 
and $f=(f_1,\ldots, f_n)$ an $n$-tuple in 
$A[t]$ such that $\det(\partial f_i/\partial t_j)$
generates an open ideal in $A'$. Consider the $A$-algebra
$C=A[t]/(f)$. Then, $C[\pi^{-1}]$ is \'etale over $A[\pi^{-1}]$
and $h$ induces a morphism of $A[\pi^{-1}]$-algebras
$C[\pi^{-1}]\to A[\pi^{-1}]$. From standard properties of
ordinary \'etale algebras, we deduce that 
$\Spec(A[\pi^{-1}])\to \Spec(C[\pi^{-1}])$ 
is a clopen immersion. 
Passing to the analytification over $A$ in the sense of 
Construction \ref{cons:analytification-}, 
we deduce a clopen immersion $\Spf(A)^{\rig}
\to \Spec(C[\pi^{-1}])^{\an}$. But the latter factors as 
follows:
$$\Spf(A)^{\rig} \xrightarrow{s} \Spf(A')^{\rig} \to 
\Spec(C[\pi^{-1}])^{\an},$$
where the second map is an open immersion.
This finishes the proof.
\end{proof}

\begin{prop}
\label{prop:locally-around-section}
Let $i:Z \to X$ be a closed immersion of rigid analytic spaces.
Let $X'$ be an \'etale rigid analytic $X$-space and 
$s:Z \to X'$ a partial section. Then, locally on $X$,
$s$ extends to a section 
$\tilde{s}:U \to X'$ defined on an open neighbourhood $U$ of $Z$.
Moreover, $\tilde{s}$ is an open immersion.
\end{prop}

\begin{proof}
The question being local on $X$, 
we may assume that $X=\Spf(A)^{\rig}$ with 
$A$ an adic ring of principal ideal type, and $Z=\Spf(B)^{\rig}$
with $B$ a quotient of $A$ by a closed ideal $I\subset A$.
We may also assume that $X'=\Spf(A')^{\rig}$ with 
$A'$ a rig-\'etale $A$-algebra, and that the section $s$ 
is induced by a morphism of $A$-adic rings
$h:A'\to B$. Let $\pi\in A$ be a generator of
an ideal of definition. Without loss of generality, 
we may assume that $B$ and $A'$ are $\pi$-torsion-free.

For $N\in \N$ and $J\subset I$ a finitely generated ideal,
consider the adic $A$-algebra $C_{J,\,N}=A\langle J/\pi^N\rangle$
given as the $\pi$-adic completion of the 
sub-$A$-algebra $A[J/\pi^N]\subset A[\pi^{-1}]$
generated by fractions $a/\pi^N$ with $a\in J$.
Then $B$ is the filtered colimit in the category of adic rings
of the $C_{J,\,N}$'s when $N$ and $J$ vary.
Applying Corollary \ref{cor:proj-limi-cal-E-A}
to this inductive system, we can find $J$ and $N$ such that the image of
$\Hom_A(A',C_{J,\,N}) \to \Hom_A(A',B)$
contains $h$. This means that the section $s$ extends to an
$X$-morphism $\Spf(C_{J,\,N})^{\rig} \to \Spf(A')^{\rig}$.
Since $\Spf(C_{J,\,N})^{\rig}$ 
is an open subspace of $X$, this proves the existence of
$\tilde{s}$ as in the proposition. 
That $\tilde{s}$ is an open immersion follows 
from Lemma \ref{lem:sect-open-immer}.
\end{proof}

\begin{dfn}
\label{dfn:rig-smooth}
$\empty$

\begin{enumerate}

\item[(1)] Let $A$ be an adic ring and $B$ an adic $A$-algebra. 
We say that $B$ is rig-smooth over $A$ if, locally on $B$, there
exists a rig-\'etale morphism of adic $A$-algebras
$A\langle t_1,\ldots, t_m\rangle
\to B$.
\ncn{adic rings!rig-smooth over}

\item[(2)] A morphism $\mathcal{Y}\to \mathcal{X}$ of formal schemes is said to be rig-smooth if, locally for the rig topology 
on $\mathcal{X}$ and $\mathcal{Y}$ (see
Definition \ref{dfn:rig-Zar-Nis} below), it is 
isomorphic to $\Spf(B) \to \Spf(A)$ with $B$ rig-smooth over $A$.
\ncn{formal schemes!rig-smooth morphism}

\item[(3)] A morphism of rigid analytic spaces $Y\to X$ is said to be 
smooth if, locally on $X$ and $Y$, it admits a formal model which is 
rig-smooth.
\ncn{rigid analytic spaces!smooth morphism}

\end{enumerate}
\end{dfn}

\begin{rmk}
\label{rmk:compare-adic-sm-et} 
By \cite[Corollary 1.6.10 \& Proposition 1.7.1]{huber}, 
we see that, via the embedding of Corollary 
\ref{cor:from-Huber-to-FK}, a map of uniform adic spaces 
is smooth (resp. \'etale) if and only if the associated map of 
rigid analytic spaces is.
\end{rmk}

The next proposition is similar to
\cite[page 582, Th\'eor\`eme 7]{elkik},
but we do not assume the adic ring $A$ to be noetherian.

\begin{prop}
\label{prop:formal-compl-smooth}
Let $A$ be an adic ring of monogenic ideal type and 
$\pi\in A$ a generator of an ideal of definition of $A$. 
Let $B$ be a rig-\'etale (resp. rig-smooth) 
adic $A$-algebra, and assume that 
$B$ is $\pi$-torsion-free.
Then, locally on $B$, there exists a finitely generated
$\pi$-torsion-free
$A$-algebra $P$ such that $P[\pi^{-1}]$ is \'etale (resp. smooth)
over $A[\pi^{-1}]$
and its $\pi$-adic completion 
$\widehat{P}=\lim_{n\in \N} P/\pi^n$ is isomorphic to $B$.
\end{prop}

\begin{proof}
According to 
\cite[pages 588--589]{elkik}, the proof of 
\cite[page 582, Th\'eor\`eme 7]{elkik}
can be adapted to cover the above statement. 
Alternatively, one can use Proposition 
\ref{prop:on-etale-presentation} as follows. 
By this proposition,
we may assume that the adic $A$-algebra $B$ is of the form
$$B=A\langle t_1,\ldots, t_m,s_1,\ldots, s_n\rangle/
(f_1,\ldots, f_n)^{\sat},$$
with $f_1,\ldots, f_n\in A[t_1,\ldots, t_m,s_1,\ldots, s_n]$, and  
such that $\det(\partial f_i/\partial s_j)$ generates an open 
ideal in $B$. (The rig-\'etale case corresponds to $m=0$.) 
Consider the $A$-algebra
$$P'=A[t_1,\ldots, t_m,s_1,\ldots, s_n]/(f_1,\ldots, f_n)^{\sat}$$
whose $\pi$-adic completion is $B$. Let 
$e\in P'$ be the image of 
$\det(\partial f_i/\partial s_j)$ in $P'$.
By assumption, a power of $\pi$ is a multiple of
$e$ in the $\pi$-adic completion of $P'$. Thus, there are elements
$b,c\in B$ and an integer $N$ such that 
$\pi^N=e\cdot b+c\pi^{N+1}$. The $A$-algebra 
$P=P'[(1-c\pi)^{-1}]$ satisfies the properties
required in the statement.
\end{proof}

The following is a variant of Proposition
\ref{prop:locally-around-section}
for smooth morphisms. It will play a crucial role in the 
proof of the localization property for rigid analytic motives 
(see Proposition \ref{prop:loc1}).

\begin{prop}
\label{prop:local-struct-smooth-sect}
Let $Z \to X$ be a closed immersion of rigid analytic spaces.
Let $X'$ be a smooth rigid analytic $X$-space and 
$s:Z \to X'$ a partial section. Then, locally on $X$, 
we may find an open neighbourhood $U\subset X$ of $Z$, 
an open neighbourhood $U'\subset X'$ of $s(Z)$ and 
an isomorphism $U'\simeq \B^m_U$, for some integer $m\geq 0$,
modulo which $s:Z \to U'$ is the zero section over $Z$. 
\end{prop}

\begin{proof}
The problem being local on $X$ and around $s(Z)$, 
we may assume that $X'$ is \'etale over 
$\B^m_X$ and, by change of coordinates, that the composition 
$$Z \xrightarrow{s} X' \to \B^m_X$$
is the zero section over $Z$. Applying 
Proposition \ref{prop:locally-around-section}
to the \'etale morphism $X' \to \B^m_X$ and the 
closed immersion $Z \to \B^m_X$ given by the zero section over $Z$, 
we find locally an open neighbourhood $U'\subset X'$ of $s(Z)$ 
such that $U'\to \B^m_X$ is also an open immersion. 
Letting $U$ be the inverse
image of $U'$ by the zero section $X \to \B^m_X$ and replacing 
$U'$ by $U'\times_X U$, we may assume that 
$U'$ is an open neighbourhood of the zero section of 
$\B^n_U$. Since the zero section of 
$\B^n_U$ admits a system of fundamental neighbourhoods which are 
$m$-dimensional relative balls, we may 
also assume that $U'$ is isomorphic to $\B^m_U$ as needed.
\end{proof}

We end this subsection with the following result.

\begin{prop}
\label{prop:open-image-for-rig-smooth-morphism}
Let $f:Y \to X$ be a smooth morphism of rigid analytic spaces.
Then the induced map $|f|:|Y|\to |X|$ is open.
\end{prop}

\begin{proof}
It is enough to show that $f(|Y|)$ is open in $|X|$.
The question is local on $X$ and $Y$. 
By Proposition \ref{prop:formal-compl-smooth}
we may assume that $X=\Spf(A)^{\rig}$, with $A$ an adic ring 
of principal ideal type, and $Y=\Spf(B)^{\rig}$, 
with $B=\widehat{P}$ the $\pi$-adic completion of a finitely 
presented $A$-algebra $P$ such that $P[\pi^{-1}]$ 
is smooth over $A[\pi^{-1}]$. (As usual, $\pi$ is a generator
of an ideal of definition of $A$. Also, 
note that finite presentation in 
Proposition \ref{prop:formal-compl-smooth}
can be assumed if we don't insist on $\pi$-torsion-freeness.)
By the Raynaud--Gruson platification theorem 
\cite[Theorem 5.2.2]{platification}, and working locally 
over $X$, we may further assume that 
$P$ is flat over $A$. 
By \cite[Chapitre IV, Th\'eor\`eme 2.4.6]{EGAIV2},
the morphism $\Spec(P) \to \Spec(A)$ is then open, and we denote 
by $U\subset \Spec(A)$ its image.
Let $(a_i)_i$ be a family in $A$ generating the ideal defining 
the complement of $U$ in $\Spec(A)$. Let $A_i$ be the 
$\pi$-adic completion of $A[a_i^{-1}]$
and $B_i$ the $\pi$-adic completion of $P_i=P\otimes_A A_i$.
Set $X_i=\Spf(A_i)^{\rig}$ and $Y_i=\Spf(B_i)^{\rig}$.
By construction, $(Y_i)_i$ is an open covering of $Y$
and it is enough to show that $f(Y_i)$ is open in $X$. 
We will show more precisely that $f(Y_i)=X_i$, i.e., that
$Y_i \to X_i$ is surjective.

Replacing $X$ and $Y$ by $X_i$ and 
$Y_i$, we are reduced to showing that 
$f:Y \to X$ is surjective, for $X=\Spf(A)^{\rig}$ and 
$Y=\Spf(B)^{\rig}\simeq \Spf(\widehat{P})^{\rig}$ as above, 
assuming furthermore that the $A$-algebra $P$ 
is faithfully flat. 
To do so, it will be enough to show the following assertion. 
If $\mathcal{X}' \to \Spf(A)$ is an admissible blowup
and $\mathcal{Y}'=(\mathcal{X}'\otimes_AB)/(0)^{\sat}$,
the induced map $\mathcal{Y}'_{\sigma}
\to \mathcal{X}'_{\sigma}$ is surjective. 
(Indeed, by 
\cite[Chapter III, Proposition 3.1.5]{fujiwara-kato},
the obvious map 
$|Y|\to |\mathcal{Y}'_{\sigma}|$
is surjective.)
Since $P$ is flat over $A$, the formal scheme 
$\mathcal{X}'\times_{\Spf(A)}\Spf(B)$ is already saturated
and we have an isomorphism 
$\mathcal{Y}'/\pi\simeq \mathcal{X}'/\pi\otimes_AP$.
In particular, we see 
that the map $\mathcal{Y}'/\pi \to \mathcal{X}'/\pi$ is 
faithfully flat, and hence surjective as needed.
\end{proof}

\subsection{Topologies}

$\empty$

\smallskip

\label{subsect:topol-rig}

Open covers define the Zariski topologies on schemes and 
formal schemes, and the analytic topology on 
rigid analytic spaces.
In this subsection, we introduce 
various finer Grothendieck topologies which we use 
when discussing motives.
On schemes, we mainly consider the \'etale and 
Nisnevich topologies. These topologies extend naturally to formal schemes: 
a family $(\mathcal{Y}_i \to \mathcal{X})_i$ consisting of 
\'etale morphisms is an \'etale (resp. a Nisnevich) cover 
if $(\mathcal{Y}_{i,\sigma}\to \mathcal{X}_{\sigma})_i$ 
is an \'etale (resp. a Nisnevich) cover.

\begin{nota}
\label{not:small-sites}
Given a scheme $S$, we denote by $\Et/S$ the category of 
\'etale $S$-schemes. Similarly, given a formal scheme $\mathcal{S}$,
we denote by $\Et/\mathcal{S}$ the category of \'etale
formal $\mathcal{S}$-schemes.
\symn{$\Et$}
\end{nota}

\begin{lemma}
\label{lem:iso-site}
Let $\mathcal{S}$ be a formal scheme.
The functor $\mathcal{X}\mapsto \mathcal{X}_{\sigma}$ induces 
an \'equivalence of categories
$\Et/\mathcal{S} \to \Et/\mathcal{S}_{\sigma}$
respecting the \'etale and Nisnevich topologies.
\end{lemma}

\begin{proof}
This follows immediately from 
\cite[Chapitre IV, Th\'eor\`eme 18.1.2]{EGAIV4}.
\end{proof}

\begin{nota}
\label{not:small-sites-rigid}
Given a rigid analytic space $S$, we denote by $\Et/S$ 
the category of \'etale rigid analytic $S$-spaces
(in the sense of Definition \ref{dfn:new-etale-rig}). 
We denote by $\Etgr/S$ the full subcategory of $\Et/S$ spanned by 
those \'etale rigid analytic $S$-spaces with good reduction
(in the sense of Definition
\ref{dfn:smooth-etale-proper-finite-rig}).
\symn{$\Etgr$}
\end{nota}

\begin{dfn}
\label{dfn:rig-top-fsch}
Let $(Y_i \to X)_i$ be a family of \'etale morphisms of rigid analytic 
spaces. We say that this family is a Nisnevich cover if, 
locally on $X$ and after refinement, it admits a formal model 
$(\mathcal{Y}_i \to \mathcal{X})_i$ which is a Nisnevich cover. 
Nisnevich covers generate a topology on 
rigid analytic spaces which we call the Nisnevich topology.
\ncn{rigid analytic spaces!Nisnevich cover}
\end{dfn}

\begin{dfn}
\label{dfn:etale-top-rigspc}
Let $(f_i:Y_i \to X)_i$ be a family of \'etale morphisms of rigid analytic 
spaces. We say that this family is an \'etale cover if it is 
jointly surjective, i.e., $|X|=\bigcup_i f_i(|Y_i|)$. 
\'Etale covers generate the \'etale topology on rigid analytic spaces.
\ncn{rigid analytic spaces!\'etale cover}
\end{dfn}

\begin{rmk}
\label{rmk:compare-etale-cover-adic}
By means of Proposition \ref{prop:compar-Spa-Spf} and
Remark \ref{rmk:compare-adic-sm-et}, we see that the above 
definition of \'etale covers agrees with the one for uniform 
adic spaces in \cite[Section 2.1]{huber}.
Also, note that if $X$ is quasi-compact, then every 
\'etale cover of $X$ can be refined by a finite subfamily. 
This follows from Proposition
\ref{prop:open-image-for-rig-smooth-morphism}.
\end{rmk}

\begin{nota}
\label{nota:topol-et-gr-nis}
The \'etale topology is generally denoted by ``\sym{$\et$}'' and the 
Nisnevich topology is denoted by ``\sym{$\Nis$}''.
Also, the Zariski topology is generally 
denoted by ``\sym{$\Zar$}'' and the analytic 
topology is denoted by ``\sym{$\an$}''.
\end{nota}

\begin{rmk}
\label{rmk:small-sites}
If $S$ is a scheme 
and $\tau\in \{\Nis,\et\}$, we call $(\Et/S,\tau)$ 
the small $\tau$-site of $S$, and similarly for a formal scheme. 
If $S$ is a rigid analytic space, we call 
$(\Etgr/S,\Nis)$ the small Nisnevich site of $S$ and 
$(\Et/S,\et)$ the small \'etale site of $S$.
\end{rmk}

The big smooth sites introduced below are used for 
constructing the categories of motives.

\begin{nota}
\label{not:big-sites}
$\empty$

\begin{enumerate}

\item[(1)] If $S$ is a scheme, we denote by $\Sch/S$ the overcategory
of $S$-schemes and $\Sm/S$ its full subcategory consisting of smooth
objects. For $\tau\in \{\Nis,\et\}$, we call $(\Sm/S,\tau)$
the big smooth site of $S$.
\symn{$\Sm$}

\item[(2)] If $\mathcal{S}$ is a formal scheme, we denote by 
$\FSch/\mathcal{S}$ the overcategory of formal $\mathcal{S}$-schemes
and $\FSm/\mathcal{S}$ its full subcategory consisting of smooth
objects. For $\tau\in \{\Nis,\et\}$, we call $(\FSm/S,\tau)$ 
the big smooth site of $\mathcal{S}$.
\symn{$\FSm$}

\item[(3)] If $S$ is a rigid analytic space, we denote by 
$\RigSpc/S$ the overcategory of rigid analytic $S$-spaces 
and $\RigSm/S$ its full subcategory consisting of smooth 
objects (in the sense of Definition
\ref{dfn:rig-smooth}). For $\tau\in \{\Nis,\et\}$, we call 
$(\RigSm/S,\tau)$ the big smooth site of $S$.
\symn{$\RigSm$}

\end{enumerate}
\end{nota}

We next discuss the class of rig topologies 
on formal schemes.

\begin{dfn}
\label{dfn:rig-Zar-Nis}
Let $(\mathcal{Y}_i \to \mathcal{X})_i$ be a family of morphisms
of formal schemes. We say that this family is a rig cover if 
the induced family $(\mathcal{Y}^{\rig}_i \to \mathcal{X}^{\rig})_i$
is an open cover. The topology generated by rig covers is called 
the rig topology and it is denoted by ``\sym{$\rig$}''.
\ncn{formal schemes!rig cover}
\end{dfn}

\begin{rmk}
\label{rmk:descri-rig-covers}
Let $\mathcal{X}$ be a quasi-compact and quasi-separated formal
scheme. Then every rig cover of $\mathcal{X}$ can be refined 
by the composition of an admissible blowup
$\mathcal{X}'\to \mathcal{X}$ and a Zariski cover of $\mathcal{X}'$.
\end{rmk}

By ``equivalence of sites'' we mean a continuous functor inducing 
an equivalence between the associated ordinary topoi.

\begin{lemma}
\label{lem:equiv-top-rig-an}
Consider full subcategories $\underline{\mathcal{V}}\subset 
\FSch$ (resp. $\underline{\mathcal{V}}\subset 
\FSch/\mathcal{S}$ for a formal scheme $\mathcal{S}$) 
and $\underline{V}\subset \RigSpc$ (resp. 
$\underline{V}\subset \RigSpc/S$ with $S=\mathcal{S}^{\rig}$) 
such that:
\begin{itemize}

\item $\underline{\mathcal{V}}$ is 
stable by admissible blowups and quasi-compact open formal subschemes;

\item $\underline{V}$ contains $\mathcal{X}^{\rig}$ for every 
$\mathcal{X}\in \underline{\mathcal{V}}$, and every object of 
$\underline{V}$ is locally of this form.

\end{itemize}
Then the functor $(-)^{\rig}:\underline{\mathcal{V}}\to
\underline{V}$ defines an equivalence of sites 
$(\underline{V},\an) \xrightarrow{\sim}
(\underline{\mathcal{V}},\rig)$. In particular, we have an 
equivalence of sites 
$(\RigSpc,\an)\xrightarrow{\sim}(\FSch,\rig)$
(resp. $(\RigSpc/S,\an)\xrightarrow{\sim}(\FSch/\mathcal{S},\rig)$).
\end{lemma}

\begin{proof}
The statement would have been a particular case of
\cite[Corollary A.4]{huber}, except that we don't know a priori 
that the continuous functor $(-)^{\rig}$ defines a morphism of sites  
and that we do not assume that our categories have finite limits. 
(In fact, we are particularly interested in the case where 
$\underline{\mathcal{V}}$ is the category of rig-smooth 
formal $\mathcal{S}$-schemes, which does not admit
finite limits.) Instead of trying to modify the proof of 
\cite[Corollary A.4]{huber}, we present an independent argument. 
We only treat the absolute case since the relative case is similar. 

By \cite[Expos\'e III, Th\'eor\`eme 4.1]{SGAIV1}, we may assume that 
$\underline{\mathcal{V}}\subset \FSch^{\qcqs}$
and that $\underline{V}$ is the full subcategory of 
$\RigSpc^{\qcqs}$ spanned by objects of the form
$\mathcal{X}^{\rig}$ for $\mathcal{X}\in \underline{\mathcal{V}}$.
The rig topology on $\underline{\mathcal{V}}$ is not subcanonical
(except for very special choices of $\underline{\mathcal{V}}$). 
We denote by $\underline{\mathcal{V}}'$ the full subcategory of
the category of sheaves of sets on $(\underline{\mathcal{V}},\rig)$ 
spanned by sheafifications of representable presheaves. 
The obvious functor 
$a:\underline{\mathcal{V}}\to\underline{\mathcal{V}}'$,
sending a formal scheme $\mathcal{X}$ to the sheaf associated of
the presheaf represented by $\mathcal{X}$, 
induces an equivalence of sites
$(\underline{\mathcal{V}}',\rig)\simeq (\underline{\mathcal{V}},\rig)$,
where the topology of $(\underline{\mathcal{V}}',\rig)$ 
is the one induced from the canonical topology on the topos
of sheaves on $(\underline{\mathcal{V}},\rig)$. (This is a 
well-known fact which follows, for example, from \cite[Expos\'e IV, 
Corollaire 1.2.1]{SGAIV1}; see also 
\cite[Corollaire 4.4.52]{ayoub-th2}.)
To prove the lemma, we remark that there is an equivalence 
of categories $\underline{\mathcal{V}}'\simeq \underline{V}$
which identifies the rig topology on $\underline{\mathcal{V}}'$
with the analytic topology on $\underline{V}$.
Indeed, for an admissible blowup
$\mathcal{Y}'\to\mathcal{Y}$ in $\underline{\mathcal{V}}$, 
the diagonal map 
$\mathcal{Y}'\to \mathcal{Y}'\times_{\mathcal{Y}}\mathcal{Y}'$
is a rig cover, which implies that 
$a\mathcal{Y}'\to a\mathcal{Y}$ is an isomorphism.
Using that the Zariski topology is subcanonical 
on $\underline{\mathcal{V}}$, we deduce that 
$$\Hom_{\underline{\mathcal{V}}'}(a\mathcal{Y},a\mathcal{X})=
\underset{\mathcal{Y}'\to \mathcal{Y}\,\in\, 
\mathfrak{B}(\mathcal{Y})}{\colim}\,
\Hom_{\underline{\mathcal{V}}}(\mathcal{Y}',\mathcal{X})$$
for any $\mathcal{X},\mathcal{Y}\in \underline{\mathcal{V}}$.
The result follows then by comparison with
\eqref{eqn-not:rigid-analytic-space}.
\end{proof}

\begin{cor}
\label{cor:rig-topol-formal-schemes}
Let $\tau\in \{\Nis,\et\}$ be one of the topologies introduced
above on rigid analytic spaces.
Consider full subcategories $\underline{\mathcal{V}}\subset 
\FSch$ (resp. $\underline{\mathcal{V}}\subset 
\FSch/\mathcal{S}$ for a formal scheme $\mathcal{S}$) 
and $\underline{V}\subset \RigSpc$ (resp. 
$\underline{V}\subset \RigSpc/S$ with $S=\mathcal{S}^{\rig}$) 
satisfying the following conditions.
\begin{itemize}

\item If $\tau=\Nis$, then 
$\underline{\mathcal{V}}$ is 
stable by admissible blowups and every \'etale morphism whose target is in 
$\underline{\mathcal{V}}$ lies entirely in $\underline{\mathcal{V}}$.

\item If $\tau=\et$, then every rig-\'etale morphism whose target is in 
$\underline{\mathcal{V}}$ lies entirely in $\underline{\mathcal{V}}$.

\item $\underline{V}$ contains $\mathcal{X}^{\rig}$ for every 
$\mathcal{X}\in \underline{\mathcal{V}}$, and every object of 
$\underline{V}$ is locally of this form.

\end{itemize}
Then there exists a unique topology
$\rig\text{-}\tau$ on $\underline{\mathcal{V}}$ such that the
functor $(-)^{\rig}:\underline{\mathcal{V}}\to \underline{V}$
defines an equivalence of sites
$(\underline{V},\tau) \xrightarrow{\sim}
(\underline{\mathcal{V}},\rig\text{-}\tau)$.
In particular, we have an 
equivalence of sites 
$(\RigSpc,\tau)\xrightarrow{\sim}(\FSch,\rig\text{-}\tau)$
(resp. 
$(\RigSpc/S,\tau)\xrightarrow{\sim}(\FSch/\mathcal{S},\rig\text{-}\tau)$).
\end{cor}

\begin{rmk}
\label{rmk:rig-topol-}
Corollary \ref{cor:rig-topol-formal-schemes}
gives us two more topologies on formal schemes: the 
rig-Nisnevich topology (denoted by ``\sym{$\rigNis$}'') and
the rig-\'etale topology (denoted by ``\sym{$\riget$}'').
These topologies can be described more directly 
by their corresponding notions of covers.
A family $(\mathcal{Y}_i \to \mathcal{X})_i$
of morphisms of formal schemes is a rig-Nisnevich cover 
if the induced family $(\mathcal{Y}^{\rig}_i \to \mathcal{X}^{\rig})_i$
is a Nisnevich cover. In particular, if
$\mathcal{X}$ is a quasi-compact and quasi-separated formal scheme,
then every rig-Nisnevich cover of $\mathcal{X}$ can be refined by the
composition of an admissible blowup $\mathcal{X}'\to \mathcal{X}$
and a Nisnevich cover of $\mathcal{X}'$. Proposition 
\ref{prop:descr-rig-etale} below
gives an analogous result for rig-\'etale covers.
\ncn{formal schemes!rig-\'etale cover}
\ncn{formal schemes!rig-Nisnevich cover}
\end{rmk}

\begin{rmk}
\label{rmk:summary-morph-sites}
Summarizing, we have a diagram of morphisms of sites:
$$\xymatrix{
(\FSch,\et) \ar[d] & (\FSch,\riget) \ar[l] \ar[d] &
(\RigSpc,\et) \ar[l]^{\sim} \ar[d] \\
(\FSch,\Nis) \ar[d] & (\FSch,\rigNis) \ar[l] \ar[d] &
(\RigSpc,\Nis) \ar[l]^{\sim} \ar[d] \\
(\FSch,\Zar) & (\FSch,\rig) \ar[l] & (\RigSpc,\an). \ar[l]^{\sim}}$$
\end{rmk}

\begin{dfn}
\label{dfn:finite-rig-etale}
$\empty$

\begin{enumerate}

\item[(1)] Let $A$ be an adic ring and $B$ an adic $A$-algebra. 
We say that $B$ is finite rig-\'etale if $B$ is finite over 
$A$ and \'etale over $\Spec(A)\smallsetminus \Spec(A/I)$ 
for an ideal of definition $I$ of $A$.
\ncn{adic rings!finite rig-\'etale over}

\item[(2)] A morphism of formal schemes $\mathcal{Y}\to \mathcal{X}$
is said to be finite rig-\'etale if it is affine and, locally 
over $\mathcal{X}$, isomorphic to $\Spf(B) \to \Spf(A)$ 
with $B$ a finite rig-\'etale adic $A$-algebra.
\ncn{formal schemes!finite rig-\'etale morphism}

\item[(3)] A morphism of formal schemes $\mathcal{Y}\to \mathcal{X}$
is said to be a finite rig-\'etale covering if it is finite 
rig-\'etale and the induced morphism $|\mathcal{Y}^{\rig}|
\to |\mathcal{X}^{\rig}|$ is surjective.
\ncn{formal schemes!finite rig-\'etale cover}

\end{enumerate}
\end{dfn}

\begin{lemma}
\label{lem:surj-finite-rig-etale}
Let $A$ be an adic ring and $B$ a finite adic $A$-algebra. 
Then $\Spf(B) \to \Spf(A)$ is 
a finite rig-\'etale covering if and only if 
\begin{equation}
\label{eq-lem:surj-finite-rig-etale}
\Spec(B)\smallsetminus \Spec(B/IB) \to 
\Spec(A)\smallsetminus \Spec(A/I)
\end{equation}
is a finite \'etale covering, when $I$ is an ideal of definition of $A$.
\end{lemma}

\begin{proof}
The morphism \eqref{eq-lem:surj-finite-rig-etale}
is finite \'etale if and only if 
$\Spf(B) \to \Spf(A)$ is finite rig-\'etale.
So we need to show that 
\eqref{eq-lem:surj-finite-rig-etale} 
is surjective if and only if 
$|\Spf(B)^{\rig}|\to |\Spf(A)^{\rig}|$ is surjective.
This follows easily from the description of 
$|\Spf(A)^{\rig}|$ in terms of valuation rings 
of residue fields of points of $\Spec(A)\smallsetminus \Spec(A/I)$
and \cite[Chapter VI, \S8, n$^{\circ}$ 6, 
Proposition 6]{bourbaki-algcomm}.
\end{proof}

\begin{rmk}
\label{rmk:finite-etale-adic-space}
Using the embedding of Corollary \ref{cor:from-Huber-to-FK}, 
it follows from Lemma \ref{lem:surj-finite-rig-etale}
that a map of uniform adic spaces is finite \'etale 
(as in \cite[Example 1.6.6.(ii)]{huber}) if and only if it 
has a finite rig-\'etale formal model.
\end{rmk}

\begin{prop}
\label{prop:descr-rig-etale}
Let $\mathcal{X}$ be a quasi-compact and quasi-separated 
formal scheme. 
Then every rig-\'etale cover of $\mathcal{X}$ can be 
refined by the composition of an admissible blowup
$\mathcal{X}'\to \mathcal{X}$, a Nisnevich cover
$(\mathcal{Y}'_i \to \mathcal{X}')_i$, 
and finite rig-\'etale coverings
$\mathcal{Z}'_i \to \mathcal{Y}'_i$.
\end{prop} 

\begin{proof}
Let $(\mathcal{U}_j\to \mathcal{X})_{j\in J}$ 
be a rig-\'etale cover. We may assume that $J$ is finite 
(see Remark \ref{rmk:compare-etale-cover-adic}) and that 
$\mathcal{X}=\Spf(A)$ is affine with $A$ an adic ring of principal 
ideal type. We fix a generator $\pi\in A$ of an ideal of definition 
of $A$. By Proposition
\ref{prop:formal-compl-smooth}, we may refine the 
rig-\'etale cover and assume that each $\mathcal{U}_j$ is the 
adic completion of a finite presentation $A$-scheme $U_j$
which is \'etale over $A[\pi^{-1}]$. 
(Note that finite presentation in 
Proposition \ref{prop:formal-compl-smooth}
can be assumed if we don't insist on $\pi$-torsion-freeness.)
By the Raynaud--Gruson platification theorem 
\cite[Theorem 5.2.2]{platification}, there exists an admissible 
blowup $X'\to X=\Spec(A)$ such that the strict transform
$U'_j \to X'$ of $U_j \to X$ is flat for every $j$.
In particular, the morphism $U'_j \to X'$ is also quasi-finite.

Let $\mathcal{U}'_j$ and $\mathcal{X}'$ be the adic completions
of $U'_j$ and $X'$. By construction, we have 
$\mathcal{X}'^{\rig}\simeq \mathcal{X}^{\rig}$ and 
$\mathcal{U}'^{\rig}_j\simeq \mathcal{U}_j^{\rig}$.
Thus, $(\mathcal{U}'_j \to \mathcal{X}')_j$ is also a rig-\'etale
cover. Since $\mathcal{O}_{\mathcal{X}'}$ and the $\mathcal{O}_{\mathcal{U}_j'}$'s are $\pi$-torsion-free,
we deduce that the family $(\mathcal{U}'_j \to \mathcal{X}')_j$
is jointly surjective. Equivalently, the family of quasi-finite
morphisms $(U'_j/\pi \to X'/\pi)_j$ is jointly surjective.
Using standard properties of the Nisnevich 
topology, we can find a family of \'etale morphisms
$(Y'_i \to X')_i$
such that:
\begin{enumerate}

\item[(1)] $(Y'_i/\pi \to X'/\pi)_i$ is a Nisnevich cover of $X'/\pi$;

\item[(2)] for every index $i$ there is a index $j$ and a clopen immersion 
$Z'_i \to U'_j\times_{X'} Y'_i$ such that $Z'_i \to Y'_i$
is finite and $Z'_i/\pi \to Y'_i/\pi$ is surjective.

\end{enumerate}
In addition to being finite, the morphism
$Z'_i \to Y'_i$ is flat and \'etale over 
$Y'_i[\pi^{-1}]$. Since $Z'_i/\pi \to Y'_i/\pi$
is surjective, we may replace $Y'_i$ by 
an open neighbourhood of $Y'_i/\pi$ and assume that 
$Z'_i \to Y'_i$ is also surjective. 
In particular, we see that $Z'_i[\pi^{-1}]\to Y'_i[\pi^{-1}]$ 
is a finite \'etale covering. If $\mathcal{Y}'_i$ and
$\mathcal{Z}'_i$ denote the adic completions of $Y'_i$ and $Z_i$,
Lemma \ref{lem:surj-finite-rig-etale}
implies that the morphisms $\mathcal{Z}'_i \to \mathcal{Y}'_i$ are
finite rig-\'etale coverings.
Moreover, the family $(\mathcal{Y}_i'\to \mathcal{X}')_i$ is a Nisnevich 
cover by point (1) above. Finally, the family 
$(\mathcal{Z}'_i\to \mathcal{X})_i$ refines the initial 
rig-\'etale cover as needed.
\end{proof}

\begin{cor}
\label{cor:limit-etale-topi-rig-an}
Let $(\mathcal{S}_{\alpha})_{\alpha}$ be a cofiltered inverse system 
of quasi-compact and quasi-separated formal schemes 
with affine transition maps, and let $\mathcal{S}=\lim_{\alpha}
\mathcal{S}_{\alpha}$ be the limit of this system. 
We set $S_{\alpha}=\mathcal{S}^{\rig}_{\alpha}$ and 
$S=\mathcal{S}^{\rig}$. Then, there is an equivalence of sites
$(\Et/S,\et) \simeq \lim_{\alpha}(\Et/S_{\alpha},\et)$.
\end{cor}

\begin{proof}
Without loss of generality, we may assume that the indexing category 
of the inverse system
$(\mathcal{S}_{\alpha})_{\alpha}$ admits a final object $o$.
We may replace $\mathcal{S}_o$ by the blowup of a finitely generated 
ideal of definition and each $\mathcal{S}_{\alpha}$ by its 
strict transform, and assume that the $\mathcal{S}_{\alpha}$'s are
locally of principal ideal type. 
The question being local for the Zariski topology on $\mathcal{S}_o$, 
we may assume that the formal schemes
$\mathcal{S}_{\alpha}$'s are affine of principal ideal type. 
We set $A_{\alpha}=\mathcal{O}(\mathcal{S}_{\alpha})$
and $A=\mathcal{O}(S)$, and we employ Notation
\ref{not:affine-rig-etale-}.
Using Corollary 
\ref{cor:rig-topol-formal-schemes}, 
it is enough to show that the morphism of sites
$$(\mathcal{E}'_A,\riget) \to \lim_{\alpha}(\mathcal{E}'_{A_{\alpha}},\riget)$$
is an equivalence. 
Corollary 
\ref{cor:proj-limi-cal-E-A}
gives an equivalence on the underlying categories
and it remains to show that the topologies match. For this, 
we need to show that every rig-\'etale cover in 
$\mathcal{E}'_A$ can be refined by the image of a rig-\'etale 
cover in $\mathcal{E}'_{A_{\alpha}}$ for $\alpha$ small enough.
This follows readily from Proposition
\ref{prop:descr-rig-etale}.
\end{proof}

\begin{rmk}
\label{rmk:analog-Nis-limit-etale-topi-rig-an}
Keeping the notation of Corollary 
\ref{cor:limit-etale-topi-rig-an}, we similarly have an 
equivalence of sites $(\Etgr/S,\Nis)\simeq \lim_{\alpha}
(\Etgr/S_{\alpha},\Nis)$. This is easier to prove: one reduces to 
the analogous statement for the small Nisnevich sites of 
formal schemes, and then further to the analogous statement for 
the small Nisnevich sites of ordinary schemes using 
Lemma \ref{lem:iso-site}. 
\end{rmk}

We end this subsection with a short discussion 
of points in the rigid analytic setting.

\begin{dfn}
\label{dfn:point-like-rigid-analytic-space}
A rigid point $s$ is a rigid analytic space of the form
$\Spf(V)^{\rig}$ where $V$ is an adic valuation ring 
of principal ideal type; compare with 
\cite[Chapter II, Definition 8.2.1]{fujiwara-kato}.
We also write $s$ for the unique closed point of $|s|$.
Using Notation \ref{not:point-residue-field-},
we then have $V=\kappa^+(s)$. Also, $\kappa(s)$ is the fraction
field of $V$, $\widetilde{\kappa}(s)$ is the residue field of $V$
and $\kappa^{\circ}(s)$ is the localisation of $V$ at its 
height $1$ prime ideal. A morphism of rigid points 
$s'\to s$ is a morphism of rigid analytic spaces 
sending the closed point of $|s'|$ to the closed point of $|s|$.
Said differently, the induced morphism $\kappa^+(s) \to \kappa^+(s')$ 
is local. 
\ncn{rigid points}
\ncn{rigid analytic spaces!rigid points|see {rigid points}}
\end{dfn}

\begin{rmk}
\label{rmk:algebraic-morphism-etale-etc}
A morphism of rigid points 
$\overline{s} \to s$ is said to be algebraic if the complete field
$\kappa(\overline{s})$ contains a dense separable extension 
of $\kappa(s)$. Algebraic rigid points over 
$s$ are all obtained by the following recipe. Start with 
a separable extension $L/\kappa(s)$ and choose a 
valuation ring $V\subset L$ such that $V\cap \kappa(s)=\kappa^+(s)$. 
(By \cite[Chapter VI, \S8, n$^{\circ}$ 6, 
Proposition 6 \& Corollary 1]{bourbaki-algcomm}
such valuation rings exist, and they 
are conjugate under the automorphism group of the extension
$L/\kappa(s)$ if the latter is Galois.)
Then define a rigid point $\overline{s}$ 
by taking $\kappa^+(\overline{s})$ to be the adic completion of 
$V$ (considered as a $\kappa^+(s)$-algebra). 
By \cite[Proposition 3.4.1/6]{BGR}, 
if $L$ is a separable closure of $\kappa(s)$, 
then $\kappa(\overline{s})$ is algebraically closed
(and not only separably closed).
\ncn{rigid points!algebraic morphism}
\end{rmk}

\begin{dfn}
\label{dfn:geometric-points-tau}
Let $\overline{s}$ be a rigid point. 
\begin{enumerate}

\item[(1)] We say that $\overline{s}$ is $\Nis$-geometric if the 
valuation ring $\kappa^+(\overline{s})$ is Henselian. 
\ncn{rigid points!nis-geometric}

\item[(2)] We say that $\overline{s}$ is $\et$-geometric 
(or, simply, geometric) if the field $\kappa(\overline{s})$ 
is algebraically closed.
\ncn{rigid points!\'et-geometric}
\end{enumerate}
\end{dfn}

\begin{rmk}
\label{rmk:tau-geom-pt-of-S}
Let $S$ be a rigid analytic space.

\begin{enumerate}

\item[(1)] A point $s\in S$ determines a rigid point, which we denote 
again by $s$, given by $\Spf(\kappa^+(s))^{\rig}$. Moreover, we have 
an obvious morphism of rigid analytic spaces $s\to S$
sending the closed point of $|s|$ to $s\in |S|$.

\item[(2)] A morphism of rigid analytic spaces
$\overline{s} \to S$ from a rigid point $\overline{s}$
is called a rigid point of $S$. It factors uniquely as 
$\overline{s} \to s \to S$, 
where $s\in |S|$ is the image of the closed
point of $|\overline{s}|$. By abuse of language, we 
say that ``$s$ is the image of $\overline{s} \to S$''
or that ``$\overline{s}$ is over $s$''. 
We say that a rigid point $\overline{s}\to S$ of $S$ is
algebraic if the morphism of rigid points
$\overline{s} \to s$ is algebraic.
(See Remark \ref{rmk:algebraic-morphism-etale-etc}.)

\end{enumerate}
\end{rmk}

\begin{lemma}
\label{lem:limit-of-formal-schemes}
Let $\mathcal{S}$ be a formal scheme and set $S=\mathcal{S}^{\rig}$.
\begin{enumerate}

\item[(1)] Given a point $s\in S$, there is a canonical isomorphism
$$\Spf(\kappa^+(s))\simeq \lim_{\Spf(\kappa^+(s)) \to 
\mathcal{U} \to \mathcal{S}}\mathcal{U},$$
where the limit is over factorizations of 
$\Spf(\kappa^+(s))\to \mathcal{S}$ with 
$\mathcal{U}$ affine and such that 
$\mathcal{U}^{\rig}$ is an open neighbourhood of $s$ in $S$.

\item[(2)] Given an algebraic rigid point $\overline{s} \to S$,
there is a canonical isomorphism
$$\Spf(\kappa^+(\overline{s}))\simeq 
\lim_{\Spf(\kappa^+(\overline{s})) \to 
\mathcal{U} \to \mathcal{S}}\mathcal{U},$$
where the limit is over factorizations of 
$\Spf(\kappa^+(\overline{s}))\to \mathcal{S}$ with 
$\mathcal{U}$ affine and rig-\'etale over $\mathcal{S}$.

\end{enumerate}
\end{lemma}

\begin{proof}
Assertion (1) follows immediately from 
\cite[Chapter II, Proposition 3.2.6]{fujiwara-kato}
and the definition of $\kappa^+(s)$; see
Notation \ref{not:point-residue-field-}.
To prove assertion (2), we may assume that $\mathcal{S}=\Spf(A)$
is affine and prove that the $A$-algebra $\kappa^+(\overline{s})$
is a filtered colimit of rig-\'etale adic $A$-algebras
in the category of adic rings.
Let $s\in S$ be the image of $\overline{s}$. 
Using assertion (1), we may write 
$$\kappa^+(s)=\underset{\alpha}{\colim}\, A_{\alpha},$$
in the category of adic rings, where $A_{\alpha}$ are adic 
$A$-algebras such that the $\Spf(A_{\alpha})^{\rig}$ are open 
neighbourhoods of $s$ in $S=\Spf(A)^{\rig}$. 
Applying Corollary \ref{cor:proj-limi-cal-E-A}
to the inductive system $(A_{\alpha})_{\alpha}$, 
we see that every rig-\'etale $\kappa^+(s)$-algebra 
whose zero ideal is saturated is a filtered colimit in the category of 
adic rings of rig-\'etale adic $A$-algebras.
Thus, it is enough to show that 
$\kappa^+(\overline{s})$ is a filtered colimit of 
adic rig-\'etale $\kappa^+(s)$-algebras.
This follows immediately from Remark
\ref{rmk:algebraic-morphism-etale-etc}
and the following fact. If $L/\kappa(s)$ 
is a finite separable extension and $R\subset L$ 
is a sub-$\kappa^+(s)$-algebra of finite type with fraction field 
$L$, then $R$ is a rig-\'etale $\kappa^+(s)$-algebra.
(We leave it to the reader to find a presentation of 
$R$ as in Definition \ref{dfn:new-etale-rig}(1).)
\end{proof}

\begin{cons}
\label{cons:point-an-nis-et}
Let $\tau\in \{\Nis,\et\}$.
Let $S$ be a rigid analytic space and let $s\in S$ be a point. 
We may construct an algebraic $\tau$-geometric rigid point 
$\overline{s} \to S$ over $s$ as follows. 
\begin{enumerate}

\item[(1)] (The case $\tau=\Nis$) 
Let $\widetilde{\kappa}(\overline{s})/\widetilde{\kappa}(s)$
be a separable extension 
and denote by $\overline{\kappa}{}^+(s)$
the Henselisation of $\kappa^+(s)$ at the point
$\Spec(\widetilde{\kappa}(\overline{s}))\to \Spec(\kappa^+(s))$.
Then $\overline{\kappa}{}^+(s)$ is again a valuation ring.
(This follows from \cite[Chapter VI, \S8, n$^{\circ}$ 6, 
Proposition 6]{bourbaki-algcomm}.)
We denote by $\kappa^+(\overline{s})$ the adic completion
of $\overline{\kappa}{}^+(s)$ and set 
$\overline{s}=\Spf(\kappa^+(\overline{s}))^{\rig}$.
We have an obvious map $\overline{s} \to S$, 
which factors through $s\to S$. The map 
$\overline{s} \to S$ is a $\Nis$-geometric rigid point of $S$.

\item[(2)] (The case $\tau=\et$) Let $\overline{\kappa}(s)$ be a 
separably closed algebraic extension of $\kappa(s)$. 
(We do not require this extension to be separable.)
Let $\overline{\kappa}{}^+(s) \subset \overline{\kappa}(s)$
be a valuation ring which extends $\kappa^+(s)\subset 
\kappa(s)$. We denote by 
$\kappa^+(\overline{s})$ the adic completion of 
$\overline{\kappa}{}^+(s)$ and set 
$\overline{s}=\Spf(\kappa^+(\overline{s}))^{\rig}$.
(As mentioned above, by \cite[Proposition 3.4.1/6]{BGR},
the fraction field $\kappa(\overline{s})$ 
of $\kappa^+(\overline{s})$ is always algebraically closed.)
We have an obvious map 
$\overline{s}\to S$ which factors through $s\to S$.
The map $\overline{s}\to S$ is an \'etale geometric rigid point of $S$. 

\end{enumerate}
In the situation of (1) (resp. (2)), 
given a presheaf $\mathcal{F}$ on $\Etgr/S$ (resp. $\Et/S$)
with values in an $\infty$-category admitting filtered colimits,
we set: 
$$\mathcal{F}_{\overline{s}}=\underset{\overline{s}\to U \to S}{\colim}
\mathcal{F}(U),$$
where the colimit is over the \'etale neighbourhoods 
with good reduction
(resp. \'etale neighbourhoods) 
of $\overline{s}$ in $S$. The object 
$\mathcal{F}_{\overline{s}}$ is called the \nc{stalk} of $\mathcal{F}$ 
at $\overline{s}$.
\end{cons}

\begin{rmk}
\label{rmk:stalk-at-rig-point}
The functors 
$\mathcal{F}\mapsto \mathcal{F}_{\overline{s}}$
introduced in Construction 
\ref{cons:point-an-nis-et}
admit a more basic version for the analytic topology,
given by $\mathcal{F}\mapsto 
\mathcal{F}_s=\colim_{s\in U  \subset X}
\,\mathcal{F}(U)$, where the colimit is over the 
open neighbourhoods of $s$ in $S$.
\end{rmk}

\begin{prop}
\label{prop:enough-points-rig-an}
Let $S$ be a rigid analytic space.

\begin{enumerate}

\item[(1)] The site $(\Etgr/S,\Nis)$ 
admits a conservative family of points
given by $\mathcal{F}\mapsto \mathcal{F}_{\overline{s}}$,
where $\overline{s}\to S$ run over the $\Nis$-geometric rigid
points as in Construction \ref{cons:point-an-nis-et}(1).

\item[(2)] The site $(\Et/S,\et)$ admits a conservative family of points 
given by $\mathcal{F}\mapsto \mathcal{F}_{\overline{s}}$,
where $\overline{s}\to S$ run over the geometric rigid points 
as in Construction \ref{cons:point-an-nis-et}(2).

\end{enumerate}
\end{prop}

\begin{proof}
We only treat the second part.
By a standard argument, one 
reduces to prove the following two assertions.
\begin{enumerate}

\item[(1)] Every \'etale cover of a geometric rigid point 
$\overline{s}$ splits.

\item[(2)] A family $(Y_i \to X)_i$ in $\Et/S$ is an \'etale cover
if, for every geometric rigid point $\overline{s} \to S$ 
and every $S$-morphism $\overline{s} \to X$, there exists $i$ and 
an $X$-morphism $\overline{s}\to Y_i$.

\end{enumerate}
The first assertion follows from Proposition
\ref{prop:descr-rig-etale} (and Corollary
\ref{cor:rig-topol-formal-schemes}).
The second assertion follows from Definition 
\ref{dfn:etale-top-rigspc}. 
\end{proof}

\begin{cor}
\label{cor:refining-etale-cover-open}
Let $S$ be a rigid analytic space and $U\subset S$ 
a nonempty open subspace. Assume that $U$ and $S$ are
quasi-compact. Then, every \'etale cover
of $U$ can be refined by the base change of an \'etale cover of $S$.
\end{cor}

\begin{proof}
Fix an \'etale cover $(U_i \to U)_i$ of $U$ with $U_i$ quasi-compact
and quasi-separated.
Given an algebraic geometric rigid point 
$\overline{s} \to S$, we consider 
$\overline{u}=\overline{s}\times_SU$.
This is a quasi-compact open rigid analytic subspace of $\overline{s}$. 
Thus, $\overline{u}$ is either empty or 
$\overline{u}\to U$ is an algebraic geometric rigid point of 
$U$. In both cases, the morphism 
$\overline{u}\to U$ factors through $U_i$ 
for some $i$. Using Corollary 
\ref{cor:limit-etale-topi-rig-an} and Lemma
\ref{lem:limit-of-formal-schemes}, there exists an 
\'etale neighbourhood 
$V_{\overline{s}}\to S$ of $\overline{s}$ such 
that $V_{\overline{s}}\times_S U$
factors through $U_i$. This shows that the base change
of the \'etale cover
$(V_{\overline{s}}\to S)_{\overline{s}}$
refines $(U_i \to U)_i$ as needed.
\end{proof}

\section{Rigid analytic motives}

\label{sect:rigid-motives}

In this section, we recall the construction of rigid analytic motives 
following \cite{ayoub-rig} and prove some of their basic properties. 
In particular, we prove in Subsection
\ref{subsect:descent-rigsh} that the functor 
$\RigSH_{\tau}(-;\Lambda)$, sending a rigid analytic 
space $S$ to the $\infty$-category of rigid analytic motives 
over $S$, is a $\tau$-sheaf with values in $\Prl$.
An important result obtained in 
this section is Theorem \ref{thm:anstC}
asserting that this sheaf
transforms certain limits of rigid analytic spaces into colimits 
of presentable $\infty$-categories.
This result plays an important role at several places in the paper,  
notably for constructing direct images with compact support in 
Subsection \ref{subsect:exceptional-functors}.
In Subsection \ref{sec:stalks}, we use this result  
for computing the stalks of 
$\RigSH_{\tau}(-;\Lambda)$.

\subsection{The construction}

$\empty$

\smallskip

\label{subsect:dfn-rigda}

From now on,
we fix a connective commutative ring spectrum 
$\Lambda \in \CAlg(\Sp_{\geq 0})$ and denote by 
$\Mod_{\Lambda}$ the $\infty$-category of $\Lambda$-modules. 
Connectivity of $\Lambda$ is assumed here for convenience. 
It implies that $\Mod_{\Lambda}$ admits a $t$-structure whose
heart is the ordinary category of $\pi_0\Lambda$-modules. 
Examples of $\Lambda$ include
localisations of the sphere spectrum at various primes 
and Eilenberg--Mac Lane spectra of ordinary rings such as
$\Z$, $\Z/n$, $\Q$, etc.
\symn{$\Mod$}

\begin{nota}
\label{not:sheaves-hom-type}
Given an $\infty$-category $\mathcal{C}$, 
we denote by $\mathcal{P}(\mathcal{C})$ the 
$\infty$-category 
of presheaves on $\mathcal{C}$ with values in the $\infty$-category 
$\mathcal{S}$ of Kan complexes. 
If $\mathcal{C}$ is endowed with a Grothendieck topology $\tau$, 
we denote by $\Shv^{(\hyp)}_{\tau}(\mathcal{C})$ the full 
sub-$\infty$-category of $\mathcal{P}(\mathcal{C})$ 
spanned by the $\tau$-(hyper)sheaves. Thus, 
$\Shv_{\tau}(\mathcal{C})$ is the $\infty$-topos associated 
to the site $(\mathcal{C},\tau)$ as in 
\cite[Definition 6.2.2.6]{lurie} and 
$\Shv_{\tau}^{\hyp}(\mathcal{C})$
is its hypercompletion in the sense of \cite[\S 6.5.2]{lurie}. 
\symn{$\mathcal{P}$}
\symn{$\Shv^{(\hyp)}$}
\end{nota}

\begin{nota}
\label{not:pre-sheaves-general}
Given an $\infty$-category $\mathcal{C}$, we denote by 
$\PSh(\mathcal{C};\Lambda)$ the $\infty$-category of 
presheaves of $\Lambda$-modules on $\mathcal{C}$, i.e., 
contravariant functors from $\mathcal{C}$ to $\Mod_{\Lambda}$.
If $\mathcal{C}$ is endowed with a Grothendieck topology $\tau$, 
we denote by $\Shv^{(\hyp)}_{\tau}(\mathcal{C};\Lambda)$ the full 
sub-$\infty$-category of $\PSh(\mathcal{C};\Lambda)$ 
spanned by the $\tau$-(hyper)sheaves.
(For the precise meaning, see Definition
\ref{hypersheaves} below.) We denote by 
\begin{equation}
\label{eq-not:pre-sheaves-general}
\Lder_{\tau}:\PSh(\mathcal{C};\Lambda)
\to \Shv^{(\hyp)}_{\tau}(\mathcal{C};\Lambda)
\end{equation}
the left adjoint to the 
obvious inclusion. This functor 
is called $\tau$-(hyper)sheafification. 
We also denote by 
\begin{equation}
\label{eq-not:pre-sheaves-general-2}
(-)^{\hyp}:\Shv_{\tau}(\mathcal{C};\Lambda)
\to \Shv^{\hyp}_{\tau}(\mathcal{C};\Lambda)
\end{equation}
the left adjoint to the obvious inclusion. This functor 
is called \nc{hypercompletion}. 
\symn{$\PSh$}
\symn{$\Lder_\tau$}
\symn{$(-)^{\hyp}$}
\end{nota}

\begin{rmk}
\label{rmk:t-structure-connect-sheaves}
The $\infty$-category 
$\Shv^{(\hyp)}_{\tau}(\mathcal{C};\Lambda)$ is stable 
and admits a $t$-structure whose truncation functors 
are denoted by $\tau_{\geq m}$ and $\tau_{\leq n}$,
and whose heart is the category of ordinary sheaves of 
$\pi_0\Lambda$-modules on the homotopy category of 
$\mathcal{C}$. An object $\mathcal{F}\in 
\Shv^{(\hyp)}_{\tau}(\mathcal{C};\Lambda)$ 
is said to be $m$-connective (resp. $n$-coconnective)
if the natural map $\tau_{\geq m}\mathcal{F} \to \mathcal{F}$
(resp. $\mathcal{F}\to \tau_{\leq n}\mathcal{F}$) 
is an equivalence. As usual, when $m=0$ (resp. $n=0$)
we say that $\mathcal{F}$ is connective (resp. coconnective).
\symn{$\tau_\geq$}
\symn{$\tau_\leq$}
\end{rmk}

We record the following lemma which 
we will use at several occasions.
Similar results can be found in~\cite[Lemma~C.3]{hoyois-quadratic} and~\cite[Proposition~2.22]{MR3545934}.

\begin{lemma}
\label{lem:equi-of-sites-infty-topoi}
Consider two sites $(\mathcal{C},\tau)$ and $(\mathcal{C}',\tau')$,
where $\mathcal{C}$ and $\mathcal{C}'$ are
ordinary categories, and let $F:\mathcal{C} \to \mathcal{C}'$
be a functor. Assume the following conditions.
\begin{enumerate}

\item \label{old-assumption-1} 
The topologies $\tau$ and $\tau'$ 
are induced by pretopologies ${\rm Cov}_{\tau}$ and 
${\rm Cov}_{\tau'}$ in the sense of \cite[Expos\'e II, 
D\'efinition 1.3]{SGAIV1}. 

\item \label{old-assumption-2}
For $X\in \mathcal{C}$, 
$F$ takes a family in ${\rm Cov}_{\tau}(X)$
to a family in ${\rm Cov}_{\tau'}(F(X))$.
Moreover, if $a:U \to X$ is an arrow which is a member 
of a family belonging to ${\rm Cov}_{\tau}(X)$ and $b:V \to X$ 
a second arrow in $\mathcal{C}$, we have 
$F(U\times_X V)\simeq F(U)\times_{F(X)}F(V)$.

\end{enumerate}
Then, the inverse image functors on presheaves 
induce by sheafification the following functors:
\begin{equation}
\label{eq-lem:equi-of-sites-infty-topoi}
F^*:\Shv_{\tau}^{(\hyp)}(\mathcal{C})\to 
\Shv_{\tau'}^{(\hyp)}(\mathcal{C}')
\qquad \text{and} \qquad 
F^*:\Shv_{\tau}^{(\hyp)}(\mathcal{C};\Lambda)\to 
\Shv_{\tau'}^{(\hyp)}(\mathcal{C}';\Lambda).
\end{equation}
Assume now, in addition, the following conditions.
\begin{enumerate}
\setcounter{enumi}{2}

\item \label{new-assumption1}
For $X\in \mathcal{C}$, any family in ${\rm Cov}_{\tau'}(F(X))$
can be refined by the image by $F$ of a family in ${\rm Cov}_{\tau}(X)$.

\renewcommand{\thefootnote}{\fnsymbol{footnote}}
\item \label{new-assumption2}
Every object $Y\in \mathcal{C}'$ admits a   
$\tau'$-hypercover by objects lying in the essential image of $F$.
Moreover, in the nonhypercomplete case, this 
hypercover can be chosen to be truncated. 
\footnote[1]{A previous version of this article did not include this condition but instead asked that the functor in \eqref{new-assumption3} be an equivalence. This is not enough (in the nonhypercomplete case): cfr. Simon Henry's answer in \url{https://mathoverflow.net/q/317085}.}
\renewcommand{\thefootnote}{\arabic{footnote}}
\setcounter{footnote}{4}
\item \label{new-assumption3}
The functor $F$ induces a fully faithful embedding
$F^*:\Shv_{\tau}(\mathcal{C})_{\leq 0}
\to \Shv_{\tau'}(\mathcal{C}')_{\leq 0}$
between the associated ordinary topoi.

\end{enumerate}
Then the functors \eqref{eq-lem:equi-of-sites-infty-topoi} are equivalences of $\infty$-categories.
\end{lemma}

\begin{proof}
The case of (hyper)sheaves of $\Lambda$-modules
follows from the case of (hyper)sheaves of spaces
using \Cref{recollection-hypersheaves}(2). Consider 
the functors 
$$\xymatrix{\mathcal{P}(\mathcal{C}) \ar@<.6pc>[r]^-{F^*} 
\ar@<-.6pc>[r]_-{F^!} &
\mathcal{P}(\mathcal{C}') \ar[l]|-{F_*}}$$
where $F^*$ is given by left Kan extension along $F$,
$F^!$ by right Kan extension along $F$ and 
$F_*$ by composition with $F$. Recall that 
$F_*$ is right adjoint to $F^*$ and 
$F^!$ is right adjoint to $F_*$.
We will prove the following assertions.
\begin{enumerate}

\item[(A)] Under the assumptions 
\eqref{old-assumption-1} and 
\eqref{old-assumption-2}, $F^*$
takes $\tau$-local equivalences to $\tau'$-local equivalences
(in both the hypercomplete and nonhypercomplete cases).
Equivalently, $F_*$ takes $\tau'$-(hyper)sheaves to 
$\tau$-(hyper)sheaves.

\item[(B)] Under the assumptions \eqref{old-assumption-1} and
\eqref{new-assumption1},
$F_*$ takes $\tau'$-local equivalences to $\tau$-local 
equivalences 
(in both the hypercomplete and nonhypercomplete cases). 
Equivalently, $F^!$ takes $\tau$-(hyper)sheaves
to $\tau'$-(hyper)sheaves.

\end{enumerate}
Assertion (A) is clear. Indeed, the assumptions (1) and (2) 
imply that $F^*$ takes $\tau$-(hyper)covers to $\tau'$-(hyper)covers.
The argument for (B) is standard (see 
\cite[Lemmas~2.13, 2.14, 2.18, 2.19]{MR3545934}),
but we give a proof for convenience. 
We first treat the nonhypercomplete case. 
In this case, the class of
$\tau'$-local equivalences is the smallest strongly saturated class
(in the sense of \cite[Definition 5.5.4.5]{lurie})
containing the inclusions of 
$\tau'$-covering sieves $R\hookrightarrow \yon(Y)$, 
for $Y\in \mathcal{C}'$.
Since $F_*$ is colimit-preserving, we only need to show that 
$F_*(R)\to F_*(\yon(Y))$ is a $\tau$-local equivalence. By the universality 
of colimits \cite[Proposition 6.1.3.10]{lurie}, 
it is enough to show that for every $X\in \mathcal{C}$
and every morphism $u:\yon(X) \to F_*(\yon(Y))$ in 
$\mathcal{P}(\mathcal{C})$,
the monomorphism 
$$P=\yon(X)\times_{u,\,F_*(\yon(Y))}F_*(R) \to \yon(X)$$ 
is a $\tau$-covering sieve. 
The morphism $u$ corresponds by adjunction to a morphism 
$v:F(X) \to Y$ in $\mathcal{C}'$, and it is easy to see that 
$P\hookrightarrow \yon(X)$ is the inclusion of the sieve of $X$ 
consisting of 
those morphisms $X'\to X$ in $\mathcal{C}$
such that $F(X')\to F(X)$ belongs to
$\yon(F(X))\times_{v,\,Y}R$, 
which is a $\tau'$-covering sieve of $F(X)$. Assumption 
\eqref{new-assumption1}
implies that $P$ is a $\tau$-covering sieve of $X$ 
as needed. 
Next, we explain how to deduce (B) in the hypercomplete
case from the nonhypercomplete case. Let $f$ be a
$\tau'$-local equivalence in $\mathcal{P}(\mathcal{C}')$. 
Then, for every $n\in \N$, $\tau_{\leq n}(f)$ is a 
truncated $\tau'$-local equivalence, and 
$F_*\tau_{\leq n}(f)=\tau_{\leq n}F_*(f)$
is a truncated $\tau$-local equivalence by (B) in
the nonhypercomplete case. Since, in the hypercomplete 
case, $\tau$-local equivalences are detected on the
truncations, the result follows.

It is now easy to conclude. Assertion (A) implies the 
existence of the first functor in 
\eqref{eq-lem:equi-of-sites-infty-topoi}.
To prove that this functor is fully faithful, 
we verify that the unit morphism $\id \to F_*F^*$
is an equivalence. (We stress that $F^*$ denotes the first functor in 
\eqref{eq-lem:equi-of-sites-infty-topoi} and
$F_*$ is the direct image functor on (hyper)sheaves.) 
Assertion (B) implies that $F_*$ has a right adjoint, 
which is the restriction of $F^!$ to (hyper)sheaves.
In particular, $F_*$ is colimit-preserving, and it is enough to 
check that the unit morphism $\id \to F_*F^*$
is an equivalence on objects of the form $\Lder_{\tau}\yon(X)$,
for $X\in \mathcal{C}$. Now $\Lder_{\tau}\yon(X)$ is $0$-truncated
and the same is true for $F^*\Lder_{\tau}\yon(X)\simeq
\Lder_{\tau}\yon(F(X))$.
Since $F_*$ preserves $0$-truncated objects (being exact), 
it follows that the unit morphism 
$\Lder_{\tau}\yon(X)\to F_*F^*\Lder_{\tau}\yon(X)$
identifies with the unit of the adjunction $(F^*,F_*)$
on the ordinary topoi associated to $(\mathcal{C},\tau)$
and $(\mathcal{C}',\tau')$. Thus we can now conclude using the
assumption \eqref{new-assumption3}.

Finally, to finish the proof, 
it remains to see that $\Shv_{\tau'}(\mathcal{C}')$ 
can be generated under colimits by the image of the functor $F^*$ in 
\eqref{eq-lem:equi-of-sites-infty-topoi}.
This follows immediately from the assumption
\eqref{new-assumption2}.
\end{proof}

Below and elsewhere in this paper, 
``monoidal'' always means ``symmetric monoidal''.

\begin{rmk}
\label{rmk:monoid-str-shv-general}
Recall that $\Mod_{\Lambda}$ underlies a monoidal 
$\infty$-category $\Mod_{\Lambda}^{\otimes}$. 
Applying 
\cite[Proposition 3.1.2.1]{lurie} to the coCartesian fibration 
$\Mod_{\Lambda}^{\otimes}\to \Fin$, we deduce that 
$$\Fun(\mathcal{C}^{\op},\Mod_{\Lambda}^{\otimes})\times_{\Fun(\mathcal{C}^{\op},\,\Fin)}\Fin \to \Fin$$ 
defines a monoidal $\infty$-category 
$\PSh(\mathcal{C};\Lambda)^{\otimes}$ whose underlying 
$\infty$-category is $\PSh(\mathcal{C};\Lambda)$.
By \cite[Proposition 2.2.1.9]{lurie:higher-algebra},
$\Shv^{(\hyp)}_{\tau}(\mathcal{C};\Lambda)$
underlies a unique monoidal $\infty$-category 
$\Shv_{\tau}^{(\hyp)}(\mathcal{C};\Lambda)^{\otimes}$
such that 
\eqref{eq-not:pre-sheaves-general}
lifts to a monoidal functor.
\symn{$\Mod^{\otimes}$}
\symn{$\PSh^{\otimes}$}
\symn{$\Shv^{\otimes}$}
\end{rmk}

\begin{rmk}
\label{rmk:Lambda-tau-X}
There is a monoidal functor 
$\Lambda \otimes-:\mathcal{S}^{\times} \to 
\Mod^{\otimes}_{\Lambda}$
sending a Kan complex to the 
associated free $\Lambda$-module.
(More precisely, this is the composition of the 
infinite suspension functor $\Sigma^\infty:\mathcal{S}^{\times}
\to \Sp^{\otimes}$
with the change of algebra functor 
$\Lambda \otimes-:\Sp^{\otimes} \to 
\Mod^{\otimes}_{\Lambda}$ provided by 
\cite[Theorem 4.5.3.1]{lurie:higher-algebra}.)
It induces monoidal functors
$$\mathcal{P}(\mathcal{C})^{\times} 
\to \PSh(\mathcal{C};\Lambda)^{\otimes}
\qquad \text{and} \qquad
\Shv_{\tau}^{(\hyp)}(\mathcal{C})^{\times}
\to \Shv_{\tau}^{(\hyp)}(\mathcal{C};\Lambda)^{\otimes}.$$
Composing with the Yoneda functors $\yon:\mathcal{C}\to \mathcal{P}(\mathcal{C})$ and $\Lder_{\tau}\circ \yon:
\mathcal{C}\to \Shv_{\tau}^{(\hyp)}(\mathcal{C})$, we get functors
$$\Lambda(-):\mathcal{C} \to \PSh(\mathcal{C};\Lambda)
\qquad \text{and} \qquad \Lambda_{\tau}(-):\mathcal{C}
\to \Shv_{\tau}^{(\hyp)}(\mathcal{C};\Lambda).$$
If $\mathcal{C}$ has finite direct products, 
the above functors lift to monoidal functors 
from $\mathcal{C}^{\times}$ to 
$\PSh(\mathcal{C};\Lambda)^{\otimes}$ and 
$\Shv^{(\hyp)}_{\tau}(\mathcal{C};\Lambda)^{\otimes}$.
In particular, 
the monoidal structure on $\PSh(\mathcal{C};\Lambda)$
described in Remark \ref{rmk:monoid-str-shv-general} 
coincides with the one given by Day convolution according to
\cite[Corollary 4.8.1.12 \& Remark 4.8.1.13]{lurie:higher-algebra}.
\symn{$\Lambda(-)$}
\end{rmk}

\begin{dfn}
\label{defn:prl-prr-omeg-etc}
$\empty$

\begin{enumerate}

\item[(1)] We denote by \sym{$\Prl$} (resp. \sym{$\Prr$}) the 
$\infty$-category of presentable $\infty$-categories and left 
adjoint (resp. right adjoint) functors; 
see \cite[Definition 5.5.3.1]{lurie}.
There is an equivalence $\Prr\simeq (\Prl)^{\op}$
(see \cite[Corollary 5.5.3.4]{lurie}), and 
both $\Prl$ and $\Prr$ are sub-$\infty$-categories of 
\sym{$\CAT_\infty$}, the $\infty$-category of (possibly large) 
$\infty$-categories.
The $\infty$-category $\Prl$ underlies a monoidal 
$\infty$-category \sym{$\Prlmon$} by
\cite[Proposition 4.8.1.15]{lurie:higher-algebra}.
\ncn{categories@$\infty$-categories!presentable}
\ncn{categories@$\infty$-categories!monoidal presentable}
\ncn{categories@$\infty$-categories!compactly generated}

\item[(2)] 
We also denote by \sym{$\Prl_\omega$} the $\infty$-category of compactly
generated $\infty$-categories and left adjoint compact-preserving
functors. It is opposite to \sym{$\Prr_\omega$}, the
$\infty$-category of compactly generated $\infty$-categories and
right adjoint functors which commute with filtered colimits.
See \cite[Definition 5.5.7.1, \& Notations 5.5.7.5 \& 5.5.7.7]{lurie}.
By \cite[Lemma 5.3.2.11]{lurie:higher-algebra}, $\Prl_{\omega}$ 
underlies a monoidal $\infty$-category \sym{$\Prlmon_{\omega}$} and
the inclusion $\Prl_{\omega}\to \Prl$ lifts to a monoidal functor
$\Prlmon_{\omega}\to \Prlmon$.

\item[(3)] 
A monoidal $\infty$-category $\mathcal{M}^{\otimes}$ 
is said to be presentable (resp. compactly generated) if the 
underlying $\infty$-category $\mathcal{M}$ is presentable 
(resp. compactly generated) and the endofunctor 
$A\otimes -$ is a left adjoint functor for all $A\in \mathcal{M}$
(resp. is a left adjoint compact-preserving functor 
for all compact $A\in \mathcal{M}$). 
This is equivalent to say that $\mathcal{M}^{\otimes}$
belongs to $\CAlg(\Prl)$ (resp. $\CAlg(\Prl_{\omega})$).

\end{enumerate}
\end{dfn}

\begin{rmk}
\label{rmk:presentable-shv-general}
The $\infty$-categories $\PSh(\mathcal{C};\Lambda)$
and $\Shv_{\tau}^{(\hyp)}(\mathcal{C};\Lambda)$ are presentable
(by \cite[Proposition 5.5.3.6 \& Remark 5.5.1.6]{lurie})
and they are respectively generated under colimits
by the objects $\Lambda(X)$ and $\Lambda_{\tau}(X)$, 
for $X\in \mathcal{C}$. 
In fact, the objects $\Lambda(X)$
are compact, so that $\PSh(\mathcal{C};\Lambda)$
is compactly generated. More is true: 
the monoidal $\infty$-categories 
$\PSh(\mathcal{C};\Lambda)^{\otimes}$ and 
$\Shv_{\tau}^{(\hyp)}(\mathcal{C};\Lambda)^{\otimes}$
are presentable, and, if  
$\mathcal{C}$ has finite direct products, 
$\PSh(\mathcal{C};\Lambda)^{\otimes}$ 
is even compactly generated. 
\end{rmk}

To define the $\infty$-category of rigid analytic motives over
a rigid analytic space $S$, we consider the case where
$(\mathcal{C},\tau)$ is the big smooth site 
$(\RigSm/S,\tau)$ with $\tau\in \{\Nis,\et\}$.
(See Notation \ref{not:big-sites}(3).) 
Before proceeding to the definition, we make a
remark concerning these sites.

\begin{rmk}
\label{rmk:big-sites-set-theoretic}
The category $\RigSm/S$ is not small, and some care is
needed when speaking about presheaves and $\tau$-(hyper)sheaves
on it. In fact, the only problem that one needs to keep in 
mind is that the $\infty$-category 
$\PSh(\RigSm/S;\Lambda)$ is not locally small.
However, this problem disappears when passing to the 
sub-$\infty$-category $\Shv_{\tau}^{(\hyp)}(\RigSm/S;\Lambda)$.
Indeed, it is easy to see that this $\infty$-category is equivalent
to $\Shv_{\tau}^{(\hyp)}((\RigSm/S)^{<\alpha};\Lambda)$,
where $\alpha$ is an infinite cardinal and 
$(\RigSm/S)^{<\alpha} \subset \RigSm/S$ is the full subcategory
spanned by those rigid analytic $S$-spaces that can
be covered by $<\alpha$ opens which are quasi-compact and quasi-separated.
(This uses Lemma
\ref{lem:equi-of-sites-infty-topoi}.) 
Clearly, $(\RigSm/S)^{<\alpha}$ is essentially small and thus 
$\Shv_{\tau}^{(\hyp)}(\RigSm/S;\Lambda)$ 
is a presentable $\infty$-category.
The same remark applies to other sites such as
$(\Et/S,\tau)$, etc.
Below, whenever we need to speak about 
general presheaves on $\RigSm/S$, $\Et/S$, etc., we implicitly 
fix an infinite cardinal $\alpha$ and replace these categories by 
$(\RigSm/S)^{<\alpha}$, $(\Et/S)^{<\alpha}$, etc.

\end{rmk}

We will use the
following notation.

\begin{nota}
\label{not:relative-ball}
$\empty$

\begin{enumerate}

\item[(1)] Let $\mathcal{X}$ be a formal scheme. We denote by 
$\A^n_{\mathcal{X}}$ the relative $n$-dimensional \nc{affine space}
given by $\Spf(\mathcal{O}_{\mathcal{X}}
\langle t_1,\ldots, t_n\rangle)$.
By abuse of notation, we also write 
``$\mathcal{X}\times \A^n$'' instead of 
``$\A^n_{\mathcal{X}}$'' although $\FSch$ has no direct products 
(nor a final object).
\symn{$\A^n$}

\item[(2)] Let $X$ be a rigid analytic space. If $X$ admits a 
formal model $\mathcal{X}$, we set 
$\B^n_X=(\A^n_{\mathcal{X}})^{\rig}$. This is independent of 
the choice of $\mathcal{X}$ and, in general, we may define 
$\B^n_X$ by gluing along open immersions. 
The rigid analytic $X$-space $\B^n_X$ is called the 
relative $n$-dimensional \nc{ball}.
By abuse of notation, we also write 
``$X\times \B^n$'' instead of 
``$\B^n_X$'' although $\RigSpc$ has no direct products 
(nor a final object).
\symn{$\B^n$}

\item[(3)] If $X$ is a rigid analytic space, we denote by 
$\U^1_X\subset \B^1_X$ the open rigid analytic subspace of 
$\B^1_X$ which is locally given by 
$\Spf(\mathcal{O}_{\mathcal{X}}\langle t,t^{-1} \rangle)
\subset \Spf(\mathcal{O}_{\mathcal{X}}\langle t \rangle)$. 
The rigid analytic $X$-space $\U^1_X$ is called the 
relative \nc{unit circle}.\footnote{In \cite{ayoub-rig}, the relative unit 
circle is denoted by $\partial \B^1_X$ and, in other places in 
the literature, it is denoted by $\mathbb{T}^1_X$.}
\symn{$\U^1$}
\end{enumerate}
\end{nota}

We fix a rigid analytic space $S$ and $\tau\in \{\Nis,\et\}$.

\begin{dfn}
\label{def:DAeff}
Let $\RigSH^{\eff,\,(\hyp)}_{\tau}(S;\Lambda)$
be the full sub-$\infty$-category 
of $\Shv_{\tau}^{(\hyp)}(\RigSm/S;\Lambda)$
spanned by those objects which are local with respect to 
the collection of maps of the form $\Lambda_{\tau}(\B^1_X)\to
\Lambda_{\tau}(X)$, for $X\in \RigSm/S$, and their desuspensions. Let
\begin{equation}
\label{eq-def:DAeff-1}
\Lder_{\B^1}:\Shv^{(\hyp)}_{\tau}(\RigSm/S;\Lambda)
\to \RigSH^{\eff,\,(\hyp)}_{\tau}(S;\Lambda)
\end{equation}
be the left adjoint to the obvious inclusion. 
This is called the $\B^1$-localisation functor.
We also set $\Lder_{\B^1,\,\tau}=\Lder_{\B^1}\circ \Lder_{\tau}$
with $\Lder_{\tau}$ the $\tau$-(hyper)sheafification functor, 
see \eqref{eq-not:pre-sheaves-general}. The functor 
$\Lder_{\B^1,\,\tau}$ is called the $(\B^1,\tau)$-localisation
functor.
Given a smooth rigid analytic $S$-space $X$, 
we set $\M^{\eff}(X)=
\Lder_{\B^1}(\Lambda_{\tau}(X))$. This is
the effective motive of $X$.
\symn{$\RigSH^{\eff,\,(\hyp)}$}
\symn{$\Lder_{\B^1}$}
\symn{$\M^{\eff}$}
\end{dfn}

\begin{rmk}
\label{rmk:on-B-1-local}
The defining condition for a $\tau$-(hyper)sheaf of $\Lambda$-modules
$\mathcal{F}$ to belong to the sub-$\infty$-category
$\RigSH^{\eff,\,(\hyp)}_{\tau}(S;\Lambda)$ is equivalent to the 
condition that $\mathcal{F}$ is $\B^1$-invariant in the following sense:
for every $X\in \RigSm/S$, 
the map of $\Lambda$-modules $\mathcal{F}(X)\to \mathcal{F}(\B^1_X)$
is an equivalence. Since $\mathcal{F}$ is a $\tau$-(hyper)sheaf, 
it is enough to ask this condition for $X$ varying in a 
subcategory $\mathcal{C}\subset \RigSm/S$ 
such that every object of $\RigSm/S$ admits a $\tau$-hypercover
by objects in $\mathcal{C}$ which is moreover truncated in the 
non-hypercomplete case.
\end{rmk}

\begin{rmk}
\label{rmk:mon-str-rigsh-eff}
The $\infty$-category $\RigSH^{\eff,\,(\hyp)}_{\tau}(S;\Lambda)$
is stable and, 
by \cite[Proposition 2.2.1.9]{lurie:higher-algebra}, 
it underlies a unique monoidal $\infty$-category
$\RigSH^{\eff,\,(\hyp)}_{\tau}(S;\Lambda)^{\otimes}$ such that  
$\Lder_{\B^1}$ lifts to a monoidal functor.
Moreover, this monoidal $\infty$-category is 
presentable, i.e., belongs to $\CAlg(\Prl)$, since we localise
with respect to a small set of morphisms.
\end{rmk}

\begin{rmk}
\label{rmk:another-site-for-rigsh}
There is another site that one can use for constructing
$\RigSH^{\eff,\,(\hyp)}_{\tau}(S;\Lambda)$,
at least when $S$ admits a formal model $\mathcal{S}$
(e.g., $S$ quasi-compact and quasi-separated).
Indeed, by Corollary 
\ref{cor:rig-topol-formal-schemes},
the site $(\RigSm/S;\tau)$ is equivalent to the site 
$(\FRigSm/\mathcal{S};\rig\text{-}\tau)$
where $\FRigSm/\mathcal{S}$ denotes the full subcategory of 
$\FSch/\mathcal{S}$ whose objects are the rig-smooth 
formal $\mathcal{S}$-schemes. (See Definition
\ref{dfn:rig-smooth} and
Remark \ref{rmk:rig-topol-}).
Using Lemma \ref{lem:equi-of-sites-infty-topoi},
we deduce an equivalence of $\infty$-categories
$$\Shv^{(\hyp)}_{\rig\text{-}\tau}(\FRigSm/\mathcal{S};\Lambda)
\simeq \Shv^{(\hyp)}_{\tau}(\RigSm/S;\Lambda)$$
and $\RigSH^{\eff,\,(\hyp)}_{\tau}(S;\Lambda)$ is equivalent to the 
sub-$\infty$-category of 
$\Shv^{(\hyp)}_{\rig\text{-}\tau}(\FRigSm/\mathcal{S};\Lambda)$ 
spanned by those objects which are local with respect to the collection
of maps $\Lambda_{\rig\text{-}\tau}(\A^1_{\mathcal{X}}) \to
\Lambda_{\rig\text{-}\tau}(\mathcal{X})$, with 
$\mathcal{X}\in \FRigSm/\mathcal{S}$, and their desuspensions.
\symn{$\FRigSm$}
\end{rmk}

\begin{dfn}
\label{dfn:rigsh-stable}
Let $\Tate_S$ (or simply \sym{$\Tate$} if $S$ is clear from
the context) be the image by $\Lder_{\B^1}$ of the 
cofiber of the split inclusion
$\Lambda_{\tau}(S)\to\Lambda_{\tau}(\U^1_S)$
induced by the unit section. With the notation of 
\cite[Definition 2.6]{robalo:k-theory-bridge}, we set 
\begin{equation}
\label{eq-def:DAeff-2}
\RigSH_{\tau}^{(\hyp)}(S;\Lambda)^{\otimes}
=\RigSH^{\eff,\,(\hyp)}_{\tau}(S;\Lambda)^{\otimes}[\Tate_S^{-1}].
\end{equation}
More precisely, there is a morphism
$\Sigma^{\infty}_{\Tate}:
\RigSH_{\tau}^{\eff,\,(\hyp)}(S;\Lambda)^{\otimes}
\to \RigSH^{(\hyp)}_{\tau}(S;\Lambda)^{\otimes}$
in $\CAlg(\Prl)$, sending 
$\Tate_S$ to a $\otimes$-invertible object, and which is 
initial for this property. We denote by 
$\Omega^{\infty}_{\Tate}:\RigSH^{(\hyp)}_{\tau}(S;\Lambda)
\to \RigSH^{\eff,\,(\hyp)}_{\tau}(S;\Lambda)$
the right adjoint to $\Sigma_{\Tate}^{\infty}$.
Given a smooth rigid analytic $S$-space $X$, 
we set $\M(X)=\Sigma^{\infty}_{\Tate}\M^{\eff}(X)$. This is  
the motive of $X$.
\symn{$\RigSH^{(\hyp)}$}
\symn{$\Sigma^{\infty}_{\Tate}$}
\symn{$\Omega^{\infty}_{\Tate}$}
\symn{$\M$}
\end{dfn}

\begin{dfn}
\label{dfn:rig-an-mot}
Objects of $\RigSH^{(\hyp)}_{\tau}(S;\Lambda)$ 
are called \nc{rigid analytic motives}
over $S$. We will denote by $\Lambda$ (or $\Lambda_S$ if we need
to be more precise) 
the monoidal unit of $\RigSH^{(\hyp)}_{\tau}(S;\Lambda)$. 
For any $n\in\N$, 
we denote by \sym{$\Lambda(n)$} the image of 
$\Tate_S^{\otimes n}[-n]$ by $\Sigma^{\infty}_{\Tate}$, and by 
$\Lambda(-n)$ the $\otimes$-inverse of $\Lambda(n)$.
For $n\in \Z$, we denote by 
$M\mapsto M(n)$ the \nc{Tate twist} given by tensoring with $\Lambda(n)$.
\end{dfn}

\begin{rmk}
\label{rmk:symmetric-obj-} 
The object $\Tate_S$ is symmetric in the sense of 
\cite[Definition 2.16]{robalo:k-theory-bridge}.
(See, for example, \cite[Lemma 3.13]{jardine-mot} whose proof extends
immediately to the rigid analytic setting.)
By \cite[Corollary 2.22]{robalo:k-theory-bridge}, it follows that 
the $\infty$-category $\RigSH^{(\hyp)}_{\tau}(S;\Lambda)$ underlying 
\eqref{eq-def:DAeff-2} is equivalent to 
the colimit in $\Prl$ 
of the $\N$-diagram whose transition maps are given by 
tensoring with $\Tate_S$.
Also, by \cite[Corollary 2.23]{robalo:k-theory-bridge},
the monoidal $\infty$-category 
\eqref{eq-def:DAeff-2}
is stable.
\end{rmk}

\begin{rmk}
\label{rmk:RigSH-DA}
When $\Lambda$ is the Eilenberg--Mac Lane spectrum associated
to an ordinary ring, also denoted by $\Lambda$, the $\infty$-category 
$\RigSH_{\tau}^{(\eff,\,\hyp)}(S;\Lambda)$ is more commonly denoted 
by $\RigDA_{\tau}^{(\eff,\,\hyp)}(S;\Lambda)$. Also, when 
$\tau$ is the Nisnevich topology, we sometimes 
drop the subscript ``$\Nis$''.
\symn{$\RigDA_{\tau}^{(\eff,\,\hyp)}$}
\end{rmk}

\begin{rmk}
\label{rmk:rigda-model-cat}
There is a more traditional description of the $\infty$-category
$\RigSH^{(\eff,\,\hyp)}_{\tau}(S;\Lambda)$ using the language of model 
categories. This is the approach taken in \cite[\S 1.4.2]{ayoub-rig}.

Assume that $\Lambda$ is given as a symmetric $S^1$-spectrum, and 
denote by $\Mod_{\Delta}(\Lambda)$ the simplicial category of 
$\Lambda$-modules which we endow with the model structure
described in \cite[Corollary 5.4.2]{symmetric-spectra-HSS}. Note that the 
$\infty$-category $\Mod_{\Lambda}$ is equivalent to the 
simplicial nerve of the full subcategory of $\Mod_{\Delta}(\Lambda)$
consisting of cofibrant-fibrant objects. 
Let $\PSh_{\Delta}(\RigSm/S;\Lambda)$ be the simplicial category 
whose objects are the presheaves on $\RigSm/S$ with values in 
$\Mod_{\Delta}(\Lambda)$, which we endow with 
its projective global model structure. 
The projective $(\B^1,\tau)$-local structure
on $\PSh_{\Delta}(\RigSm/S;\Lambda)$,
also known as the motivic model structure, 
is obtained from the latter via 
the Bousfield localization with respect to the union 
of the following classes of maps:
\symn{$\Mod_{\Delta}$}
\symn{$\PSh_{\Delta}$}
\begin{enumerate}

\item[(1)] morphisms of presheaves 
inducing isomorphisms on the $\tau$-sheaves associated to 
their homotopy presheaves;

\item[(2)] morphisms of the form $\Lambda(\B^1_X)[n]\to\Lambda(X)[n]$ 
induced by the canonical projection, 
for $X \in \RigSm/S$ and $n \in \Z$.

\end{enumerate}
The $\infty$-category 
$\RigSH^{\eff,\,\hyp}_{\tau}(S;\Lambda)$ is equivalent to
the simplicial nerve of the full simplicial subcategory of 
$\PSh_{\Delta}(\RigSm/S;\Lambda)$ consisting of 
motivically cofibrant-fibrant objects.
This follows from \cite[Propositions 4.2.4.4 \& A.3.7.8]{lurie}.

To obtain the $\Tate$-stable version, we form the category 
$\Spect_T(\PSh_{\Delta}(\RigSm/S;\Lambda))$
of $T$-spectra of presheaves of $\Lambda$-modules on $\RigSm/S$.
(Here $T$ is any cofibrant replacement
of $\Lambda(\U^1_S)/\Lambda(S)$.) 
The $(\B^1,\tau)$-local model structure induces the 
stable $(\B^1,\tau$)-local model structure on $T$-spectra,
which is also known as the motivic model structure. 
The $\infty$-category $\RigSH^{\hyp}_{\tau}(S;\Lambda)$ is equivalent to
the simplicial nerve of the full simplicial subcategory of 
$\Spect_T(\PSh_{\Delta}(\RigSm/S;\Lambda))$ 
consisting of motivically cofibrant-fibrant objects.
This follows from 
\cite[Theorem 2.26]{robalo:k-theory-bridge}.
\symn{$\Spect_T$}

The above discussion can be adapted to the non-hypercomplete case.
One only needs to replace the class of maps in (1) above 
by a smaller one, namely the class of maps of the form 
$\hocolim_{[n]\,\in\, \mathbf{\Delta}}\,
\Lambda(Y_n) \to \Lambda(Y_{-1})$ 
where $Y_{\bullet}$ is a truncated $\tau$-hypercover of $Y_{-1}
\in \RigSm/S$.
In both cases, the weak equivalences of the 
(stable) $(\B^1,\tau)$-local model structure are called the 
(stable) $(\B^1,\tau)$-local equivalences.
\end{rmk}

\begin{lemma}
\label{lem:generation-rigsh}
The monoidal $\infty$-category 
$\RigSH^{(\eff,\,\hyp)}_{\tau}(S;\Lambda)^{\otimes}$
is presentable and its underlying $\infty$-category is 
generated under colimits, and up to desuspension and negative 
Tate twists when applicable, by the motives $\M^{(\eff)}(X)$ with
$X\in \RigSm/S$ quasi-compact and quasi-separated.
\end{lemma}

\begin{proof}
That the monoidal $\infty$-category of the statement is 
presentable was mentioned above. The claim about the generators
follows from Remark
\ref{rmk:presentable-shv-general}
in the effective case. In the $\Tate$-stable 
case, we then use the 
universal property of $\otimes$-inversion given by 
\cite[Proposition 2.9]{robalo:k-theory-bridge}.
\end{proof}

\begin{prop}
\label{prop:basic-functorial-rigsh}
The assignment $S\mapsto \RigSH_{\tau}^{(\eff,\,\hyp)}
(S;\Lambda)^{\otimes}$
extends naturally into a functor
\begin{equation}
\label{eq-cor:compactness-rigda}
\RigSH_{\tau}^{(\eff,\,\hyp)}(-;\Lambda)^{\otimes}:
\RigSpc^{\op}\to \CAlg(\Prl).
\end{equation}
\end{prop}

\begin{proof}
We refer to \cite[\S 9.1]{robalo} for the construction of 
an analogous functor in the algebraic setting.
\end{proof}

\begin{nota}
\label{nota:f-stars-for-rig-an}
Let $f:Y\to X$ be a morphism of rigid analytic spaces. 
The image of $f$ by 
\eqref{eq-cor:compactness-rigda}
is the inverse image functor 
$$f^*:\RigSH^{(\eff,\,\hyp)}_{\tau}(X;\Lambda)
\to 
\RigSH^{(\eff,\,\hyp)}_{\tau}(Y;\Lambda)$$
which has the structure of a monoidal functor.
Its right adjoint $f_*$ is the direct image functor.
It has the structure of a right-lax monoidal functor.
(See Lemma \ref{lem:adj-monoidal-cat} below.)
\end{nota}

\subsection{Previously available functoriality}

$\empty$

\smallskip

\label{subsect:functoriality-rigsh}

We gather here part of what is known about 
the functor $S\mapsto \RigSH_{\tau}^{(\eff,\,\hyp)}(S;\Lambda)$
introduced in Subsection \ref{subsect:dfn-rigda}.
The results that we discuss here 
were obtained in \cite[\S 1.4]{ayoub-rig}
under the assumption that $S$ 
is of finite type over a non-Archimedean field.
However, the proofs apply also to the general case with very little 
modification.

\begin{prop}
\label{prop:6f1}
Let $f:Y\to X$ be a smooth morphism of rigid analytic spaces. 
\begin{enumerate}

\item[(1)] The functor $f^*$, as in Notation
\ref{nota:f-stars-for-rig-an}, admits a left adjoint 
$$f_{\sharp}:\RigSH^{(\eff,\,\hyp)}_{\tau}(Y;\Lambda)
\to \RigSH_{\tau}^{(\eff,\,\hyp)}(X;\Lambda)$$ 
sending the motive of 
a smooth rigid analytic $Y$-space $V$ to the motive of $V$ considered 
as a smooth rigid analytic $X$-space in the obvious way.

\item[(2)] \emph{(Smooth projection formula)}
The canonical map
$$f_{\sharp}(f^*M\otimes N)\to M\otimes f_{\sharp} N$$
is an equivalence for all 
$M \in \RigSH_{\tau}^{(\eff,\,\hyp)}(X;\Lambda)$ and 
$N\in\RigSH_{\tau}^{(\eff,\,\hyp)}(Y;\Lambda)$.

\item[(3)] \emph{(Smooth base change)} 
Let $g:X'\to X$ be a morphism of rigid analytic spaces and form
a Cartesian square 
$$\xymatrix{Y'\ar[d]^{f'}\ar[r]^{g'} & Y\ar[d]^{f}\\
X' \ar[r]^{g} & X.}$$
The natural transformations 
$f'_\sharp\circ g'^*\to g^*\circ f_\sharp$ and 
$f^*\circ g_*\to g'_*\circ f'^*$,
between functors 
from $\RigSH_{\tau}^{(\eff,\,\hyp)}(Y;\Lambda)$ to 
$\RigSH_{\tau}^{(\eff,\,\hyp)}(X';\Lambda)$ and
back, are equivalences. 

\end{enumerate}
\end{prop}

\begin{proof}
The functor $f^*:\RigSm/X \to \RigSm/Y$ 
admits a left adjoint $f_{\sharp}$
sending a smooth rigid analytic $Y$-space $V$ to 
$V$ considered as a smooth rigid analytic $X$-space.
The adjunction $(f_{\sharp},f^*)$ induces an adjunction 
between categories of motives. This is discussed 
in \cite[Th\'eor\`emes 1.4.13 \& 1.4.16]{ayoub-rig} 
using the language of model categories.
For the second assertion, we refer to the proof of
\cite[Proposition 4.5.31]{ayoub-th2}.
For the third assertion, we refer to the proof of
\cite[Lemme 1.4.32]{ayoub-rig}.
Both proofs are formal and extend 
readily to the context we are considering. 
\end{proof}

\begin{cor}
\label{cor:j-sharp-star-fully-faithful}
Let $j:U \to X$ be an open immersion of rigid analytic spaces.
Then the functors 
$$j_{\sharp},\,j_*:\RigSH^{(\eff,\,\hyp)}_{\tau}(U;\Lambda)
\to \RigSH^{(\eff,\,\hyp)}_{\tau}(X;\Lambda)$$
are fully faithful.
\end{cor}

\begin{proof}
This follows from Proposition 
\ref{prop:6f1}(3) with $f$ and $g$ equal to $j$.
\end{proof}

\begin{prop}
\label{prop:loc1}
Let $i: Z\to X$ be a closed immersion of rigid analytic spaces
(as in Definition 
\ref{dfn:smooth-etale-proper-finite-rig})
and $j:U\to X$ the complementary open immersion
(i.e., such that $|U|=|X|\smallsetminus |Z|$).
\begin{enumerate}

\item[(1)] The functor $i_*:\RigSH_{\tau}^{(\eff,\,\hyp)}(Z;\Lambda)
\to \RigSH_{\tau}^{(\eff,\,\hyp)}(X;\Lambda)$ is fully faithful.

\item[(2)] \emph{(Localization)} The counit of the adjunction 
$(j_{\sharp},j^*)$
and the unit of the adjunction $(i^*,i_*)$
form a cofiber sequence 
\begin{equation}
\label{eq-prop:loc1-1}
j_{\sharp}j^{*}\to\id\to i_{*}i^{*}
\end{equation}
of endofunctors of $\RigSH_{\tau}^{(\eff,\,\hyp)}(X;\Lambda)$.
In particular, the pair $(i^*,j^*)$ is conservative.
\ncn{localization}

\item[(3)] \emph{(Closed projection formula)} 
The canonical map 
\begin{equation}
\label{eq-prop:loc1-2}
M\otimes i_*N\to i_*(i^*M\otimes N)
\end{equation}
is an equivalence for all
$M \in \RigSH_{\tau}^{(\eff,\,\hyp)}(X;\Lambda)$ 
and $N \in 
\RigSH_{\tau}^{(\eff,\,\hyp)}(Z;\Lambda)$.
\ncn{closed projection formula}

\item[(4)] \emph{(Closed base change)}
Let $g:X'\to X$ be a morphism of rigid analytic spaces and form
a Cartesian square 
$$\xymatrix{Z'\ar[d]^{i'}\ar[r]^{g'} & Z\ar[d]^{i}\\
X' \ar[r]^{g} & X.}$$
The natural transformation
$g^*\circ i_*\to i'_*\circ g'^*$,
between functors 
from $\RigSH_{\tau}^{(\eff,\,\hyp)}(Z;\Lambda)$ to 
$\RigSH_{\tau}^{(\eff,\,\hyp)}(X';\Lambda)$,
is an equivalence. If moreover $g$ is smooth, then the natural 
transformation 
$g_{\sharp}\circ i'_* \to i_* \circ g'_{\sharp}$, from 
$\RigSH_{\tau}^{(\eff,\,\hyp)}(Z';\Lambda)$ to 
$\RigSH_{\tau}^{(\eff,\,\hyp)}(X;\Lambda)$,
is an equivalence.
\ncn{closed base change}

\end{enumerate}
\end{prop}

\begin{proof}
Assertion (2) implies all the others. 
Indeed, applying 
$i^*$ to the cofiber sequence 
\eqref{eq-prop:loc1-1} and using that $i^*j_{\sharp}\simeq 0$
(which follows from Proposition 
\ref{prop:6f1}(3)), we deduce that 
$i^*i_*i^* \to i^*$ is an equivalence. Assertion (1)
follows then from Lemma
\ref{lem:i-star-colimits-1} below. 
We may check that 
\eqref{eq-prop:loc1-2}
is an equivalence after applying $i^*$ and $j^*$. 
Assertion (3) follows then by using that $j^*i_*\simeq 0$
(by Proposition \ref{prop:6f1}(2))
and $i^*i_*\simeq \id$ (by assertion (1)). Similarly, to 
prove assertion (4) we use that the pairs $(i^*,j^*)$
and $(i'^*,j'^*)$ are conservative
(with $j':U'\to X'$ the base change of $j$), and the equivalences
$j^*i_*\simeq 0$, $j'^*i'_*\simeq 0$, 
$i^*i_*\simeq \id$ and $i'^*i'_*\simeq \id$, 
and smooth base change as in
Proposition \ref{prop:6f1}(3) for the second natural transformation.

We now discuss the proof of assertion (2).
When $X$ is of finite type over a non-Archimedean 
field, assertion (2) can be found in \cite[\S 1.4.3]{ayoub-rig}.
(See \cite[Th\'eor\`eme 1.4.20]{ayoub-rig} for the effective case
and the proof of \cite[Corollaire 1.4.28]{ayoub-rig} 
for the $\Tate$-stable case.)
We claim that the proofs of loc.~cit. extend
to general rigid analytic spaces.

The key step is to show that 
\cite[Th\'eor\`eme 1.4.20]{ayoub-rig}
is still valid for general rigid analytic spaces, i.e., 
that assertion (2) holds true in the effective case.
This is the statement that for any $\mathcal{F}$
in $\RigSH^{\eff,\,(\hyp)}_{\tau}(X;\Lambda)$, the square
\begin{equation}
\label{eqn:Morel-Voev-square}
\begin{split}
\xymatrix{j_{\sharp} j^*\mathcal{F} \ar[r] \ar[d]
& \mathcal{F}\ar[d]\\
0 \ar[r] & i_*i^*\mathcal{F}}
\end{split}
\end{equation}
is coCartesian in 
$\RigSH^{\eff,\,(\hyp)}_{\tau}(X;\Lambda)$.
Using Lemma \ref{lem:generation-rigsh} and 
Lemma \ref{lem:i-star-colimits-1}
below, we may assume that 
$\mathcal{F}=\Lder_{\B^1,\,\tau}\Lambda(X')$
with $X'\in \RigSm/X$. (See Definition
\ref{def:DAeff}.)
Using Lemma \ref{lem:i-lower-star-sheaf} below, 
we have an equivalence
$$i_*i^*\Lder_{\B^1,\,\tau}\Lambda(X')
\simeq \Lder_{\B^1,\,\tau}i_* \Lambda_{t_{\emptyset}}(X'_Z)$$
where $X'_Z=X'\times_X Z$ and 
$t_{\emptyset}$ the 
topology on $\RigSpc$ generated by one family, namely the 
empty family considered as a cover of the empty rigid analytic space.
Thus, it is enough to show that 
$\Lder_{\B^1,\,\tau}$ transforms the square 
$$\xymatrix{\Lambda_{t_{\emptyset}}(X'_U) \ar[r] \ar[d]
& \Lambda_{t_{\emptyset}}(X')\ar[d]\\
0 \ar[r] & i_*\Lambda_{t_{\emptyset}}(X'_Z)}$$
into a coCartesian one.
Using the analogues of \cite[Corollaire 4.5.40 \&
Lemme 4.5.41]{ayoub-th2}, we reduce to show that 
\cite[Proposition 1.4.21]{ayoub-rig} is valid for 
general rigid analytic spaces. More precisely, 
given a partial section $s:Z \to X'$ defined 
over $Z$, we need to show that the morphism
$T_{X',\,s}\otimes \Lambda \to \{*\}\otimes \Lambda$
is a $(\B^1,\tau)$-equivalence (i.e., becomes an equivalence
after applying $\Lder_{\B^1,\,\tau}$).
Here $T_{X',\,s}$ is the presheaf of sets on 
$\RigSm/X$ given by 
$$T_{X',\,s}(P)=\left\{
\begin{array}{ccc}
\Hom_X(P,X')\times_{\Hom_Z(P\times_X Z,\,X')}\{*\}
& \text{if} & P\times_X Z\neq \emptyset,\\
\{*\} & \text{if} & P\times_X Z=\emptyset.
\end{array}\right.$$
Arguing as in the first and second steps of the proof of 
\cite[Proposition 1.4.21]{ayoub-rig}
one proves that the problem is local on $X$ and 
around $s(Z)$ for the analytic topology. 
(In loc.~cit., we only consider hypersheaves, but the reader 
can easily check that hypercompletion is not used in this reduction.)
Using Proposition
\ref{prop:local-struct-smooth-sect}, it is thus enough to treat the
case $X'=\B^m_X$ and $s$ the zero section restricted to $Z$.
In this case, we may use an explicit homotopy to conclude as in 
the third step of the proof of \cite[Proposition 4.5.42]{ayoub-th2}.

Now that assertion (2) is proven in the effective case, 
we explain how it extends to the $\Tate$-stable case.  
Since assertion (2) in the effective case implies assertion (3) 
in the effective case, the functor 
$$i_*:\RigSH^{\eff,\,(\hyp)}_{\tau}(Z;\Lambda) 
\to \RigSH^{\eff,\,(\hyp)}_{\tau}(X;\Lambda)$$
commutes with tensoring with $\Tate$, i.e., there is an equivalence
of functors $\Tate_X\otimes i_*(-)\simeq i_*(\Tate_Z\otimes -)$. 
(See Definition \ref{dfn:rigsh-stable}.)
Using Remark \ref{rmk:symmetric-obj-}
and the fact that $i_*$ belongs to $\Prl$ (by 
Lemma \ref{lem:i-star-colimits-1} below),
we deduce that $i_*$ commutes with $\Sigma^{\infty}_{\Tate}$, i.e.,
there is an equivalence $\Sigma_{\Tate}^{\infty}\circ i_*
\simeq i_*\circ  \Sigma_{\Tate}^{\infty}$.
Therefore, applying 
$\Sigma^{\infty}_{\Tate}$ to the coCartesian squares
\eqref{eqn:Morel-Voev-square}, we deduce that 
$$\xymatrix{j_{\sharp} j^*M \ar[r] 
\ar[d]
& M\ar[d]\\
0 \ar[r] & i_*i^*M}$$
is coCartesian for any $M$ in the image of 
$\Sigma_{\Tate}^{\infty}(-)$ up to a twist. 
Using Lemma \ref{lem:generation-rigsh} and 
Lemma \ref{lem:i-star-colimits-1}, we deduce that the above square 
is coCartesian for any $M\in \RigSH_{\tau}^{(\hyp)}(X;\Lambda)$.
\end{proof}

\begin{lemma}
\label{lem:i-lower-star-sheaf}
Let $i:Z \to X$ be a closed immersion of rigid analytic spaces.
The functor 
\begin{equation}
\label{eq-lem:i-lower-star-sheaf}
i_*:\Shv_{t_{\emptyset}}(\RigSm/Z;\Lambda)
\to \Shv_{t_{\emptyset}}(\RigSm/X;\Lambda)
\end{equation}
commutes with $\tau$-(hyper)sheafification and the 
$(\B^1,\tau)$-localisation functor.
\end{lemma}

\begin{proof}
This is a generalisation of \cite[Lemma 1.4.18]{ayoub-rig}.
For the proof of loc.~cit. to extend to our context, we need 
to show the following property. Given a smooth rigid analytic 
$X$-space $X'$ such that $X'_Z=X'\times_X Z$ is nonempty, 
every $\tau$-cover of $X'_Z$ can be refined by 
the inverse image of a $\tau$-cover of $X'$. To prove this, we may 
assume that $X'=X$. The question is local on $X$. 
Thus, we may assume that $X=\Spf(A)^{\rig}$, with $A$ an adic 
ring of principal ideal type, and $Z=\Spf(B)^{\rig}$
with $B$ a quotient of $A$ by a saturated closed ideal $I$. 
Let $\pi$ be a generator of an ideal of definition of $A$. 
Then $B$ is the filtered colimit in the category of adic rings
of $C_{J,\,N}=A\langle J/\pi^N\rangle$ where $N\in \N$ and 
$J\subset I$ is a finitely generated ideal. 
Set $Y_{J,\,N}=\Spf(C_{J,\,N})^{\rig}$.

By Corollary \ref{cor:limit-etale-topi-rig-an}
and Remark \ref{rmk:analog-Nis-limit-etale-topi-rig-an},
every $\tau$-cover $(V_i \to Z)_i$ 
can be refined by the restriction to $Z$ of 
a $\tau$-cover $(U_j\to Y_{J,\,N})_j$
for well chosen $J$ and $N$. 
We get a 
$\tau$-cover of $X$ with the required
property by adding to the family $(U_j\to X)_j$
the open inclusion $X\smallsetminus Z \to X$.
\end{proof}

\begin{lemma}
\label{lem:i-star-colimits-1}
Let $i:Z \to X$ be a closed immersion of rigid analytic spaces.
\begin{enumerate}

\item[(1)] The functor $i_*:\RigSH^{(\eff,\,\hyp)}_{\tau}(Z;\Lambda)
\to \RigSH^{(\eff,\,\hyp)}_{\tau}(X;\Lambda)$
commutes with colimits. Thus, it admits a right 
adjoint which we denote by $i^!$.

\item[(2)] The image of the functor $i^*:
\RigSH^{(\eff,\,\hyp)}_{\tau}(X;\Lambda)
\to \RigSH^{(\eff,\,\hyp)}_{\tau}(Z;\Lambda)$
generates the $\infty$-category $\RigSH^{(\eff,\,\hyp)}_{\tau}(Z;\Lambda)$
by colimits.

\end{enumerate}
\end{lemma}

\begin{proof}
In the effective case, assertion (1) follows from Lemma
\ref{lem:i-lower-star-sheaf}. Indeed, for a rigid 
analytic space $S$, the colimit of a diagram in 
$\RigSH^{(\eff,\,\hyp)}_{\tau}(S;\Lambda)$
is computed by applying $\Lder_{\B^1,\,\tau}$
to the colimit of the same diagram in 
$\Shv_{t_{\emptyset}}(\RigSm/S;\Lambda)$. So, 
it is enough to show that \eqref{eq-lem:i-lower-star-sheaf}
commutes with colimits, which is obvious. The passage from 
the effective case to the $\Tate$-stable
case follows from Remark
\ref{rmk:symmetric-obj-} and the commutation 
$\Tate_X\otimes i_*(-)\simeq i_*(\Tate_Z\otimes -)$. 
(This relies on assertion (2) of Proposition 
\ref{prop:loc1}, but only in the effective case, so there is no
vicious circle.)

We now prove assertion (2). 
By Lemma \ref{lem:generation-rigsh}, 
it is enough to show that the motive $\M^{(\eff)}(V)$ of
a smooth rigid analytic $Z$-space $V$ is a colimit of 
objects in the image of $i^*$. The problem is local on $X$ and $V$,
so we may assume that $X=\Spf(A)^{\rig}$, $Z=\Spf(B)^{\rig}$
and $V=\Spf(F)^{\rig}$ where 
$A$ is an adic ring of principal ideal type, $B$ a quotient 
of $A$ by a saturated closed ideal and 
$F\in \mathcal{E}'_{B\langle s\rangle}$
with $s=(s_1,\ldots, s_m)$ a system of coordinates.
(For the definition of the category 
$\mathcal{E}'_{B\langle s\rangle}$, see
Notation \ref{not:affine-rig-etale-}.)
Writing $B$ as the colimit of 
$C_{J,\,N}$ as in the proof of
Lemma \ref{lem:i-lower-star-sheaf}, we may apply Corollary 
\ref{cor:proj-limi-cal-E-A} to find 
$E\in \mathcal{E}'_{C_{J,\,N}\langle s\rangle}$, for some $J$ and $N$,
such that $E\,\widehat{\otimes}_{C_{J,\,N}}\,B/(0)^{\sat}\simeq F$.
Thus, $U=\Spf(E)^{\rig}$ is a smooth rigid analytic 
$X$-space such $U\times_X Z\simeq V$, and we have 
$i^*\M^{(\eff)}(U)\simeq \M^{(\eff)}(V)$ as needed.
\end{proof}

One of the aims of this paper is to define 
the full six-functor formalism  
for rigid analytic motives. We have seen above that 
the functors $f^*$, $f_*$, $f_\sharp$, $\otimes$ and 
$\underline{\Hom}$ can be defined with little effort. 
We now state what was known so far concerning 
the exceptional functors $f_!$ and $f^!$ following
\cite[\S 1.4.4]{ayoub-rig} (see also 
\cite[Theorem 2.9]{bamb-vezz}). 

\begin{rmk}
\label{rmk:analytific-funct}
Let $A$ be an adic ring, $I\subset A$ an ideal of 
definition, and $U=\Spec(A)\smallsetminus \Spec(A/I)$.
Recall from Construction
\ref{cons:analytification-} that there exists an 
analytification functor 
\begin{equation}
\label{eqn:analytific-funct}
(-)^{\an}:\Sch^{\lft}/U \to \RigSpc/U^{\an}
\end{equation}
from the category $\Sch^{\lft}/U$, of $U$-schemes which
are locally of finite type, to the category of rigid analytic 
$U^{\an}$-spaces. (Note that $U^{\an}=\Spf(A)^{\rig}$.)
This functor preserves \'etale and smooth morphisms, closed immersions 
and complementary open immersions, as well as proper morphisms.
\end{rmk}

The following result follows immediately from 
Propositions 
\ref{prop:6f1} and 
\ref{prop:loc1}, and the construction.

\begin{prop}
\label{prop:homotopic-functor}
Keep the notation as in Remark \ref{rmk:analytific-funct}.
The contravariant functor 
$$X \mapsto \RigSH^{(\hyp)}_{\tau}(X^{\an};\Lambda),
\quad f\mapsto f^{\an,\,*}$$ 
from $\Sch^{\lft}/U$ to $\Prl$
is a stable homotopical functor in the sense that it satisfies
the $\infty$-categorical versions of the
properties (1)--(6) listed in \cite[\S 1.4.1]{ayoub-th1}.
\end{prop}

\begin{rmk}
\label{rmk:homotopical-functor}
The $\infty$-categorical versions of the
properties (1)--(6) listed in \cite[\S 1.4.1]{ayoub-th1}
can be checked after passing to the homotopy categories. 
Thus, we may as well reformulate Proposition
\ref{prop:homotopic-functor}
by saying that the functor from $\Sch^{\lft}/U$ 
to the $2$-category of triangulated categories, sending 
$X$ to the homotopy category associated to $\RigSH^{(\hyp)}_{\tau}(X^{\an};\Lambda)$, is a stable homotopical functor in the 
sense of \cite[D\'efinition 1.4.1]{ayoub-th1}.
\end{rmk}

Proposition 
\ref{prop:homotopic-functor}
gives access to the results developed in 
\cite[Chapitres 1--3]{ayoub-th1,ayoub-th2}
yielding a limited six-functor formalism
for rigid analytic motives. 
We will not list explicitly all the properties that 
form this formalism since a full six-functor formalism 
will be obtained later in Section \ref{sec:6ff}.
We content ourselves with the following preliminary statement 
which we actually need in establishing 
the full six-functor formalism for rigid analytic motives.

\begin{cor}
\label{cor:6ffalg}
Keep the notation as in Remark \ref{rmk:analytific-funct}.
Given a morphism $f:Y\to X$ between quasi-projective 
$U$-schemes, there is an adjunction 
$$f^{\an}_!:\RigSH^{(\hyp)}_{\tau}(Y^{\an};\Lambda) 
\rightleftarrows \RigSH^{(\hyp)}_{\tau}(X^{\an};\Lambda):f^{\an,\,!}$$
Moreover, the following properties are satisfied.
\begin{enumerate}

\item[(1)] The assignments $f\mapsto f^{\an}_!$ and 
$f\mapsto f^{\an,\,!}$
are compatible with composition.\footnote{Here, we only 
claim the compatibility with composition up to non-coherent 
homotopies. A more structured version of this will be obtained later
in a more general situation; see Theorem
\ref{thm:exist-functo-exceptional-image}.}

\item[(2)] Given a Cartesian square of quasi-projective $U$-schemes
$$\xymatrix{Y' \ar[r]^-{g'} \ar[d]^-{f'} & Y \ar[d]^-f\\
X' \ar[r]^-g & X,\!}$$
there is an equivalence 
$g^{\an,\,*}\circ f^{\an}_!\simeq f'^{\an}_!\circ g'^{\an,\,*}$.

\item[(3)] There is a natural transformation 
$f_!^{\an} \to f_*^{\an}$ which is an equivalence 
if $f$ is projective.

\item[(4)] If $f$ is smooth, there are equivalences 
$f^{\an,\,!}\simeq \Th(\Omega_f)\circ f^{\an,\,*}$
and $f^{\an}_!\simeq f^{\an}_{\sharp}\circ \Th^{-1}(\Omega_f)$ 
where $\Th(\Omega_f)$ and $\Th^{-1}(\Omega_f)$
are the Thom equivalences associated 
to $\Omega_f$ as in \cite[\S 1.5.3]{ayoub-th1}.

\end{enumerate}
\end{cor}

\begin{proof}
This follows from Proposition
\ref{prop:homotopic-functor} and 
\cite[Scholie 1.4.2]{ayoub-th1}.
\end{proof}

\begin{rmk}
\label{rmk:thom-equivalence}
\symn{$\Th^{-1}$}
\nc{Thom equivalences} can be defined for any 
$\mathcal{O}_X$-module $\mathcal{M}$ 
which is locally free of finite rank
on a rigid analytic space $X$. Indeed, $\mathcal{M}$
determines a vector bundle $p:M\to X$ whose fiber at a point 
$x\in X$ is given by $\Spec(\kappa(x)[\mathcal{M}_x])^{\an}$.
We set $\Th(\mathcal{M})=p_{\sharp}s_*$ and 
$\Th^{-1}(\mathcal{M})=s^!p^*$,
where $s:X\to M$ is the zero section. If $\mathcal{M}$ is free 
of rank $m$, then $\Th(\mathcal{M}) \simeq (-)(m)[2m]$
and $\Th^{-1}(\mathcal{M}) \simeq (-)(-m)[-2m]$. 
That said, we may write ``$\Th(\Omega_{f^{\an}})$''
instead of ``$\Th(\Omega_f)$'' in Corollary 
\ref{cor:6ffalg}(4).
(If $h$ is a smooth morphism of rigid analytic spaces, 
there is an associated $\mathcal{O}$-module $\Omega_h$ 
which is locally free of finite rank. It can be defined
locally as the cokernel of the Jacobian matrix.)
\symn{$\Th$}
\end{rmk}

\begin{dfn}
\label{dfn:projective-morphism}
$\empty$

\begin{enumerate}

\item[(1)] If $S$ is a rigid analytic space, we denote 
by $\P^n_S$ the relative $n$-dimensional 
\nc{projective space} over $S$. If $S=\Spf(A)^{\rig}$, for an adic ring 
$A$, then $\P^n_S=(\P^n_{\Spf(A)})^{\rig}$, and for general $S$, 
$\P^n_S$ is defined by gluing. If $A$ and $U$ are as in 
Remark \ref{rmk:analytific-funct},
we also have $\P^n_{U^{\an}}\simeq (\P^n_U)^{\an}$.
\symn{$\P^n$}

\item[(2)] Let $f:Y \to X$ be a morphism of rigid analytic spaces. 
We say that $f$ is locally projective if, locally on $X$, 
$f$ can be factored as a closed immersion followed by a projection 
of the form $\P^n_X \to X$.
\ncn{rigid analytic spaces!locally projective morphism}

\end{enumerate}
\end{dfn}

For later use, we also record the following statement.

\begin{prop}
\label{prop:base-change-finite-and-projection}
Let $f:Y \to X$ be a locally projective morphism 
of rigid analytic spaces.
\begin{enumerate}

\item[(1)] \emph{(Projective projection formula)} 
The canonical map $M\otimes f_*N\to f_*(f^*M\otimes N)$
is an equivalence for all $M\in \RigSH^{(\hyp)}_{\tau}(X;\Lambda)$
and $N\in \RigSH^{(\hyp)}_{\tau}(Y;\Lambda)$.
\ncn{projective projection formula}

\item[(2)] \emph{(Projective extended base change)} 
Let $g:X'\to X$ be a morphism of rigid analytic spaces and form
a Cartesian square
$$\xymatrix{Y' \ar[r]^-{g'} \ar[d]^-{f'} & Y \ar[d]^-f\\
X' \ar[r]^-g & X.}$$
The natural transformation
$g^*\circ f_*\to f'_*\circ g'^*$,
between functors from $\RigSH_{\tau}^{(\hyp)}(Y;\Lambda)$ to 
$\RigSH_{\tau}^{(\hyp)}(X';\Lambda)$,
is an equivalence. If moreover $g$ is smooth, then the natural 
transformation 
$g_{\sharp}\circ f'_* \to f_* \circ g'_{\sharp}$, from 
$\RigSH_{\tau}^{(\hyp)}(Y';\Lambda)$ to 
$\RigSH_{\tau}^{(\hyp)}(X;\Lambda)$,
is an equivalence.
\ncn{projective extended base change}

\end{enumerate}
\end{prop}

\begin{proof}
If $f=f_1\circ f_2$, then the assertions for $f$ follow from their
analogues for $f_1$ and $f_2$. Also, the assertions 
can be checked locally on $X$.
Thus, it is enough to treat the case
of a closed immersion $i:Z \to X$ and the case of 
$p:\P^n_X \to X$.
The case of a closed immersion follows from Proposition 
\ref{prop:loc1}. For $p:\P^n_X\to X$, we use 
Corollary \ref{cor:6ffalg} which provides us with a canonical 
equivalence $p_*\simeq p_!=p_{\sharp}\circ \Th^{-1}(\Omega_p)$. The 
result follows then from Proposition
\ref{prop:6f1}.
\end{proof}

We now go back to the notation introduced in 
Remark \ref{rmk:analytific-funct}.
Given a $U$-scheme $X$ which is locally of finite type,
the analytification functor 
\eqref{eqn:analytific-funct}
induces a premorphism of sites
\begin{equation}
\label{eqn:analytific-funct-ran}
\An:(\RigSm/X^{\an},\tau) \to (\Sm/X,\tau).
\end{equation}
(Indeed, the analytification of an \'etale cover is an 
\'etale cover, and the analytification of a Nisnevich cover
can be refined by an open cover; see \cite[Th\'eor\`eme
1.2.39]{ayoub-rig} whose proof can be adapted to our context.)
By the functoriality of the construction of the 
$\infty$-categories of motives, \eqref{eqn:analytific-funct-ran}
induces a functor
\begin{equation}
\label{eqn:analytific-funct-ran-star}
\An^*:\SH^{(\eff,\,\hyp)}_{\tau}(X;\Lambda)
\to \RigSH^{(\eff,\,\hyp)}_{\tau}(X^{\an};\Lambda).
\end{equation}
In \cite{ayoub-rig}, this functor is denoted by $\Rig^*$.
\symn{$\An$}
\symn{$\An^*$}

\begin{prop}
\label{prop:analytif-presheaf}
The functors \eqref{eqn:analytific-funct-ran-star} 
are part of a morphism of $\CAlg(\Prl)$-valued presheaves
\begin{equation}
\label{eq-prop:analytif-presheaf-1}
\SH_{\tau}^{(\eff,\,\hyp)}(-;\Lambda)^{\otimes}
\to \RigSH_{\tau}^{(\eff,\,\hyp)}((-)^{\an};\Lambda)^{\otimes}
\end{equation}
on $\Sch^{\lft}/U$.
In particular, the functors $\An^*$ are monoidal and 
commute with the inverse image functors. Moreover, if 
$f$ is a smooth morphism in $\Sch^{\lft}/U$, the natural transformation 
$$f^{\an}_{\sharp}\circ \An^* \to \An^*\circ f_{\sharp}$$
is an equivalence.
\end{prop}

\begin{proof}
One argues as in \cite[\S 9.1]{robalo} for  
the first assertion.
The second assertion is clear.
\end{proof}

\begin{prop}
\label{prop:ran-star-com-f-lower-star-prelim}
Let $f:Y\to X$ be a proper morphism in $\Sch^{\lft}/U$.
Then, the natural transformation 
$$\An^*\circ f_*\to f^{\an}_*\circ \An^*,$$ 
between functors from $\SH^{(\hyp)}_{\tau}(Y;\Lambda)$ to
$\RigSH^{(\hyp)}_{\tau}(X^{\an};\Lambda)$, is invertible.
\end{prop}

\begin{proof}
We split the proof into two steps.

\paragraph*{Step 1}
\noindent 
Here we assume that $f$ is projective.
It is enough to prove the claim when $f$ is a closed immersion 
and when $f$ is the projection $\P^n_X\to X$. 
In the first case, one uses Proposition 
\ref{prop:loc1} and its algebraic analogue.
In the second case, one uses Corollary \ref{cor:6ffalg}
and its algebraic analogue to reduce to show that 
$f_{\sharp}^{\an}\circ \An^*\simeq \An^* \circ f_{\sharp}$
which holds by Proposition
\ref{prop:analytif-presheaf}.

\paragraph*{Step 2}
\noindent
Here we deal with the general case. 
We may assume that $X$ is quasi-compact and quasi-separated. 
Using Proposition \ref{prop:analytif-presheaf},
we reduce easily to show that
$$\An^*\circ f_*\circ j_{\sharp}\to f^{\an}_*\circ j^{\an}_{\sharp}\circ 
\An^*$$
is an equivalence for every open immersion 
$j:V \to Y$, with $V$ affine.
By the refined version of Chow's lemma given in
\cite[Corollary 2.6]{nagata-deligne},
there is a blowup $e:Y'\to Y$, with centre
disjoint from $V$, such that $f':Y'\to X$ is projective.
Let $j':V \to Y'$ be the obvious inclusion.
by Proposition
\ref{prop:base-change-finite-and-projection}(2) and its algebraic 
version, we have equivalences 
$e_*\circ j'_{\sharp}\simeq j_{\sharp}$ and $e^{\an}_*\circ j'^{\an}_{\sharp}
\simeq j^{\an}_{\sharp}$. Thus, it is enough to prove the 
proposition for $f'=f\circ e$. Since this morphism is projective, 
we may conclude by the first step. 
\end{proof}

\begin{rmk}
\label{rmk:proper-vs-projective-base-change-alg}
The method used in the second step of the proof of 
Proposition \ref{prop:ran-star-com-f-lower-star-prelim}
will be used again in the second part of the proof of Proposition
\ref{prop:prop-base-proformal} below 
to deduce the proper base change theorem 
for $\SH^{(\hyp)}_{\tau}(-;\Lambda)$ from its special 
case for projective morphisms which is covered by 
\cite[Corollaire 1.7.18]{ayoub-th1}.
(For a slightly different method using the usual version of 
Chow's lemma but requiring the schemes to be noetherian, see
the proof of \cite[Proposition 2.3.11(2)]{cd}.)
Similarly, this method can be used to generalise 
Proposition \ref{prop:base-change-finite-and-projection}
to the case where $f$ is locally the analytification 
of a proper morphism of schemes. However, our aim is to prove
a more substantial generalisation of that proposition
which cannot be reached using this method. This will 
be achieved in Theorem \ref{thm:prop-base} below.
\end{rmk}

\subsection{Descent}

$\empty$

\smallskip

\label{subsect:descent-rigsh}

In this subsection, we prove that the functor 
$S\mapsto \RigSH_{\tau}^{(\eff,\,\hyp)}(S;\Lambda)$,
$f\mapsto f^*$, whose existence is claimed in Proposition
\ref{prop:basic-functorial-rigsh}, 
admits (hyper)descent for the topology $\tau$.  
This can be considered as a folklore theorem, but we reproduce 
the proof here for completeness. For a comparable result 
in the algebraic setting, see 
\cite[Proposition 4.8]{hoyois-6op}.

For later use, we recall the precise definition
of a (hyper)sheaf valued in a general $\infty$-category.
(Compare with \cite[Definition 2.1]{drew:motivic-hodge}.)

\begin{dfn}
\label{hypersheaves}
Let $(\mathcal{C},\tau)$ be a site and let $\mathcal{V}$ 
be an $\infty$-category admitting all limits.
A functor $F:\mathcal{C}^{\op}\to\mathcal{V}$ is called a 
$\tau$-(hyper)sheaf (or is said to satisfy $\tau$-(hyper)descent) 
if its right Kan extension 
$\overline{F}:\mathcal{P}(\mathcal{C})^{\op}\to\mathcal{V}$,
along the Yoneda embedding,
factors through the opposite of the localisation functor 
$\Lder_{\tau}:\mathcal{P}(\mathcal{C}) \to
\Shv_{\tau}^{(\hyp)}(\mathcal{C})$. This is equivalent to 
the condition that $\overline{F}$ induces an equivalence 
\begin{equation}
\label{eq-hypersheaves}
\overline{F}(X_{-1})\xrightarrow{\sim}
\lim_{[n]\in\Delta}
\overline{F}(X_n)
\end{equation}
for every effective $\tau$-hypercover $X_{\bullet}$. 
(An effective $\tau$-hypercover $X_{\bullet}$ 
is an augmented simplicial object of 
$\mathcal{P}(\mathcal{C})$ such that $\Lder_{\tau}(X_{\bullet})$
is an effective hypercovering of the $\infty$-topos
$\Shv^{(\hyp)}_{\tau}(\mathcal{C})_{/\Lder_{\tau} X_{-1}}$,
in the sense of \cite[Definition 6.5.3.2]{lurie}.)
We denote by $\Shv^{(\hyp)}_{\tau}(\mathcal{C};\mathcal{V})$
the full sub-$\infty$-category of  
$\PSh(\mathcal{C};\mathcal{V})=\Fun(\mathcal{C}^{\op},\mathcal{V})$
spanned by $\tau$-(hyper)sheaves. 
When $\mathcal{V}$ is the 
$\infty$-category $\mathcal{S}$ of spaces, we get back the
$\infty$-topos $\Shv_{\tau}^{(\hyp)}(\mathcal{C})$.
\ncn{hypersheaves@(hyper)sheaves}
\ncn{hypercovering!effective}
\symn{$\Shv^{(\hyp)}$}
\end{dfn}

We gather a few facts about (hyper)sheaves with values in general 
$\infty$-categories. 
We refer the reader to \cite[\S 2]{drew:motivic-hodge}
for proofs and more details.

\begin{rmk}
\label{rmk:effective-hypercover}
Keep the notation as in Definition
\ref{hypersheaves}. 
\begin{enumerate}

\item[(1)] In the hypercomplete case, every 
$\tau$-hypercover is effective. Therefore, for $F$ to be a 
$\tau$-hypersheaf, the equivalence
\eqref{eq-hypersheaves} needs to hold for 
every $\tau$-hypercover, but see Remark 
\ref{recollection-hypersheaves}(3) below.

\item[(2)] In the non-hypercomplete case, 
for $F$ to be a $\tau$-sheaf, it is enough that the equivalence
\eqref{eq-hypersheaves}
holds for $X_{\bullet}$ the {\v C}ech nerve 
associated to a $\tau$-cover in $\mathcal{C}$.
It then holds for any truncated 
$\tau$-hypercover.
See \cite[Definition 6.2.2.6 \& Lemma 6.5.3.9]{lurie}.

\end{enumerate}
\end{rmk}

\begin{rmk}
\label{recollection-hypersheaves}
Keep the notation as in Definition
\ref{hypersheaves}. 
\begin{enumerate}

\item[(1)] Let $\phi:\mathcal{V} \to \mathcal{V}'$ be a limit-preserving 
functor between $\infty$-categories admitting all limits. 
Then the induced functor $\Phi:\PSh(\mathcal{C};\mathcal{V})
\to \PSh(\mathcal{C};\mathcal{V}')$ preserves $\tau$-(hyper)sheaves.
If moreover $\phi$ detects limits, then $\Phi$ detects 
$\tau$-(hyper)sheaves.

\item[(2)] Assume that $\mathcal{V}$ is presentable. 
Then the $\infty$-category
$\Shv_{\tau}^{(\hyp)}(\mathcal{C};\mathcal{V})$ 
is an accessible left-exact 
localization of $\PSh(\mathcal{C};\mathcal{V})$. In particular, 
it is also presentable. 
We denote by 
$$\Lder_{\tau}:\PSh(\mathcal{C};\mathcal{V})
\to \Shv_{\tau}^{(\hyp)}(\mathcal{C};\mathcal{V})$$
the $\tau$-(hyper)sheafification functor 
defined as the left adjoint to the obvious inclusion.
(This was introduced in Notation 
\ref{not:pre-sheaves-general}
for $\mathcal{V}=\Mod_{\Lambda}$.)
With respect to the monoidal structure on $\Prl$ of 
\cite[\S 4.8.1]{lurie:higher-algebra},
we have $\Shv^{(\hyp)}_{\tau}(\mathcal{C};\mathcal{V})
\simeq \Shv^{(\hyp)}_{\tau}(\mathcal{C})\otimes\mathcal{V}$;
see \cite[Proposition 2.4(1)]{drew:motivic-hodge} whose proof
is also valid in the non-hypercomplete case.
\symn{$\Lder_\tau$}

\item[(3)] If $(\mathcal{C},\tau)$ is a Verdier site 
(in the sense of \cite[Definition 9.1]{DHI})
satisfying the assumptions (1-3) of \cite[\S 10]{DHI}, 
the condition of $F$ being a $\tau$-(hyper)sheaf can be expressed without recourse to its right Kan extension $\overline{F}$.
More precisely, $F$ is a $\tau$-(hyper)sheaf if 
$F$ transforms representable 
coproducts in $\mathcal{C}$ into products in $\mathcal{V}$
and if for every internal $\tau$-hypercover $X_{\bullet}$
(in the sense of \cite[Definition 10.1]{DHI}) which is effective, 
$F$ induces an equivalence
$$F(X_{-1})\xrightarrow{\sim} 
\lim_{[n]\in\Delta} F(X_n).$$
(As explained in Remark 
\ref{rmk:effective-hypercover}, in the hypercomplete case,
effectivity is an empty condition and,
in the non-hypercomplete case, we may replace it with the 
condition that $X_{\bullet}$ is truncated or better with the condition 
that $X_{\bullet}$
is the {\v C}ech nerve of a basal morphism $X_0 \to X_{-1}$
which is a $\tau$-cover.) This is proven in 
\cite[Proposition 2.7]{drew:motivic-hodge} in the hypercomplete
case and is clear in the non-hypercomplete case. It applies to the
sites we consider in this paper, such as the big smooth sites
of Notation \ref{not:big-sites}.

\end{enumerate}
\end{rmk}

The main result of this subsection is the following.

\begin{thm}
\label{thm:hyperdesc}
Let $\tau\in \{\Nis,\et\}$ be a topology on rigid analytic spaces. 
The contravariant functor 
$$S\mapsto \RigSH_{\tau}^{(\eff,\,\hyp)}(S;\Lambda), \quad 
f\mapsto f^*$$ 
defines a $\tau$-(hyper)sheaf on $\RigSpc$ 
valued in $\Prl$.
\end{thm}

\begin{rmk}
\label{rmk:sheaf-monoidal-structure-rigsh}
The forgetful functor 
$\CAlg(\Prl) \to \Prl$ being limit-preserving and conservative 
(by \cite[Corollary 3.2.2.5 \& Lemma 3.2.2.6]{lurie:higher-algebra}),
Theorem \ref{thm:hyperdesc} and Remark 
\ref{recollection-hypersheaves}(1) 
imply that 
$\RigSH^{(\eff,\,\hyp)}_{\tau}(-;\Lambda)^{\otimes}$ 
is also a $\tau$-(hyper)sheaf valued in $\CAlg(\Prl)$.
\end{rmk}

Before we can give the proof of Theorem \ref{thm:hyperdesc}
we need a digression about general (hyper)sheaves on general sites.
Let $\mathcal{C}$ be a small $\infty$-category and $X$ 
an object of $\mathcal{C}$. Composition with the obvious projection 
${\rm j}_X:\mathcal{C}_{/X}\to\mathcal{C}$ 
induces a functor ${\rm j}_X^*:\mathcal{P}(\mathcal{C})
\to \mathcal{P}(\mathcal{C}_{/X})$ which preserves limits and 
colimits. We denote by  
${\rm j}_{X,\,!}$ the left adjoint of ${\rm j}_X^*$ and 
${\rm j}_{X,\,*}$ its right adjoint. A topology $\tau$ on 
$\mathcal{C}$ induces a topology on $\mathcal{C}_{/X}$ which 
we also denote by $\tau$. It is easy to see that ${\rm j}_X^*$ and 
${\rm j}_{X,\,*}$ preserve $\tau$-(hyper)sheaves. 
(For ${\rm j}_{X,\,*}$, note that 
modulo the equivalence $\mathcal{P}(\mathcal{C}_{/X})\simeq 
\mathcal{P}(\mathcal{C})_{/\yon(X)}$, the functor ${\rm j}_{X,\,*}$ 
takes a presheaf $F$ on $\mathcal{C}_{/X}$ to the presheaf 
$U\mapsto \Map_{\mathcal{P}(\mathcal{C})_{/\yon(X)}}
(\yon(U)\times \yon(X),F)$.)
We get in this way an adjunction 
$${\rm j}_X^*:\Shv^{(\hyp)}_{\tau}(\mathcal{C})
\rightleftarrows
\Shv_{\tau}^{(\hyp)}(\mathcal{C}_{/X}):{\rm j}_{X,\,*}$$
where ${\rm j}_X^*$ commutes with all limits and colimits. 
In particular, ${\rm j}_X^*$ admits a left adjoint (on the level of 
(hyper)sheaves) which we denote by ${\rm j}_{X,\,!}^{\tau}$. 
It is related to ${\rm j}_{X,\,!}$ by an equivalence 
${\rm j}_{X,\,!}^{\tau}\circ \Lder_{\tau}\simeq 
\Lder_{\tau}\circ ({\rm j}_X)_!$.
The following lemma is well-known. We include a proof for 
completeness.

\begin{lemma}
\label{slice-topoi}
Let $(\mathcal{C},\tau)$ be a site and 
$X\in \mathcal{C}$.
The functor ${\rm j}_{X,\,!}^{\tau}$ factors through an equivalence
$${\rm e}_X:\Shv^{(\hyp)}_{\tau}(\mathcal{C}_{/X})
\xrightarrow{\sim}
\Shv^{(\hyp)}_{\tau}(\mathcal{C})_{/\Lder_{\tau}\yon(X)}.$$
\end{lemma}

\begin{proof}
The functor ${\rm j}_{X,\,!}^{\tau}:
\Shv^{(\hyp)}_{\tau}(\mathcal{C}_{/X})\to
\Shv^{(\hyp)}_{\tau}(\mathcal{C})$ sends the final 
object $\Lder_{\tau}\yon(\id_X)$ of 
$\Shv^{(\hyp)}_{\tau}(\mathcal{C}_{/X})$
to $\Lder_{\tau}\yon(X)$.
This gives the functor ${\rm e}_X$. By construction, we have a commutative square
$$\xymatrix{\mathcal{P}(\mathcal{C}/X) \ar[r]^-{{\rm e}'_X} 
\ar[d]^-{\Lder_{\tau}} & 
\mathcal{P}(\mathcal{C})_{/\yon(X)} \ar[d]^-{\Lder'_{\tau}}\\
\Shv^{(\hyp)}_{\tau}(\mathcal{C}/X)\ar[r]^-{{\rm e}_X} & 
\Shv^{(\hyp)}_{\tau}(\mathcal{C})_{/\Lder_{\tau}\yon(X)}.\!}$$
By \cite[Corollary~5.1.6.12]{lurie}, ${\rm e}'_X$ is an equivalence.
Note that $\Lder'_{\tau}$ is essentially surjective
on objects. Indeed, given a morphism of $\tau$-(hyper)sheaves 
$F\to \Lder_{\tau}\yon(X)$, there is an equivalence 
$F\simeq \Lder_{\tau}(F\times_{\Lder_{\tau}\yon(X)}\yon(X))$
since $\Lder_{\tau}$ is exact and idempotent.
To finish the proof, it will suffice 
to show that ${\rm e}_X$ is fully faithful. Let ${\rm f}_X$ 
be a right adjoint to ${\rm e}_X$ and 
${\rm f}'_X$ 
a right adjoint to ${\rm e}'_X$.
We know that the unit $\id \to {\rm f}'_X\circ {\rm e}'_X$
is an equivalence, and  
we need to prove that the unit $\id \to 
{\rm f}_X\circ {\rm e}_X$ 
is an equivalence. By \cite[Proposition 5.2.5.1]{lurie}, 
${\rm f}_X$ sends a map
$F\to \Lder_{\tau}\yon(X)$ to the fiber product 
${\rm j}_X^*F\times_{{\rm j}_X^*\Lder_{\tau}\yon(X)}\{*\}$
and ${\rm f}'_X$ sends a map
$F'\to \yon(X)$ to the fiber product 
${\rm j}_X^*F'\times_{{\rm j}_X^*\yon(X)}\{*\}$.
Since (hyper)sheafification is exact, 
we deduce that the natural transformation
$\Lder_{\tau}\circ {\rm f}'_X \to {\rm f}_X\circ \Lder'_{\tau}$
is an equivalence. Using the commutative square
$$\xymatrix{\Lder_{\tau} \ar[r] \ar[d]^-{\sim} & {\rm f}_X\circ {\rm e}_X\circ \Lder_{\tau} \ar[d]^-{\sim} \\
\Lder_{\tau} \circ {\rm f}_X'\circ {\rm e}'_X \ar[r]^-{\sim} & {\rm f}_X\circ \Lder_{\tau}' \circ {\rm e}'_X,}$$
if follows that the natural transformation
$\Lder_{\tau} \to 
{\rm f}_X \circ {\rm e}_X\circ \Lder_{\tau}$
is an equivalence, which is enough to conclude 
since $\Lder_{\tau}$ is essentially surjective.
\end{proof}

We denote by ${\rm Top}^{\Lder}$
the $\infty$-category of $\infty$-topoi and 
exact left adjoint functors, 
as defined in \cite[Definition 6.3.1.5]{lurie}.

\begin{prop}
\label{prop:sheaf-values-topos}
Let $(\mathcal{C},\tau)$ be a site. The functor
$\Shv^{(\hyp)}_{\tau}(\mathcal{C}_{/(-)}):\mathcal{C}^{\op}
\to {\rm Top}^{\Lder}$,
taking an object $X$ of $\mathcal{C}$ to the $\infty$-topos 
$\Shv_{\tau}^{(\hyp)}(\mathcal{C}_{/X})$ and a morphism 
$f$ in $\mathcal{C}$
to the functor ${\rm j}_f^*$,
is a $\tau$-(hyper)sheaf. 
\end{prop}

\begin{proof}
Every $\infty$-topos $\mathcal{X}$ determines a 
$\mathrm{Top}^{\rm L}$-valued sheaf on itself:
by \cite[Proposition~6.3.5.14]{lurie},
the functor $\chi:\mathcal{X}^{\op}
\to {\rm Top}^{\rm L}$, sending $X\in \mathcal{X}$ to 
$\mathcal{X}_{/X}$, 
preserves limits.
Take $\mathcal{X}=\Shv^{(\hyp)}_{\tau}(\mathcal{C})$.
Since $\Lder_{\tau}:\mathcal{P}(\mathcal{C})\to \mathcal{X}$
preserves colimits, we deduce that 
$\chi\circ \Lder_{\tau}:\mathcal{P}(\mathcal{C})^{\op}
\to {\rm Top}^{\rm L}$
preserves limits. It follows that the functor  
$\chi\circ \Lder_{\tau}$ is a right Kan extension of 
$\chi\circ \Lder_{\tau}\circ \yon:\mathcal{C}^{\op}\to {\rm Top}^{\rm L}$.
Since $\chi\circ \Lder_{\tau}$ clearly factors 
through $\Shv_{\tau}^{(\hyp)}(\mathcal{C})$, the functor 
$\chi\circ \Lder_{\tau}\circ \yon$ is a $\tau$-(hyper)sheaf. Now, by 
Lemma \ref{slice-topoi}, the functor $\chi\circ \Lder_{\tau}\circ \yon$
is equivalent to the one sending $X\in \mathcal{C}$ to 
$\Shv^{(\hyp)}_{\tau}(\mathcal{C}_{/X})$.
\end{proof}

\begin{cor}
\label{cor:sheafy-sheaves}
Let $(\mathcal{C},\tau)$ be a site and $\mathcal{V}$ 
a presentable $\infty$-category. The functor
$\Shv^{(\hyp)}_{\tau}(\mathcal{C}_{/(-)};\mathcal{V}):
\mathcal{C}^{\op}\to \Prl$,
taking an object $X$ of $\mathcal{C}$ to the $\infty$-category 
$\Shv^{(\hyp)}_{\tau}(\mathcal{C}_{/X};\mathcal{V})$ and a morphism 
$f$ in $\mathcal{C}$ to the functor ${\rm j}_f^*$, 
is a $\tau$-(hyper)sheaf. 
\end{cor}

\begin{proof}
By Proposition \ref{prop:sheaf-values-topos}, the result holds 
when $\mathcal{V}$ is the $\infty$-category of spaces $\mathcal{S}$,
and we want to reduce to this case. We denote by 
$\mathcal{X}(-;\mathcal{V}):\mathcal{C}^{\op}\to \Prl$
the functor sending $X\in \mathcal{C}$ to 
$\Shv^{(\hyp)}_{\tau}(\mathcal{C}_{/X};\mathcal{V})$. 
By Remark \ref{recollection-hypersheaves},
we have an equivalence of functors 
$\mathcal{X}(-;\mathcal{S})\otimes \mathcal{V}
\xrightarrow{\sim} \mathcal{X}(-;\mathcal{V})$,
where the tensor product is taken in $\Prl$
(see \cite[\S4.8.1]{lurie:higher-algebra}).
Moreover, for any $f:Y\to X$ in $\mathcal{C}$, 
the functor ${\rm j}_f^*:\mathcal{X}(X;\mathcal{S})
\to \mathcal{X}(Y;\mathcal{S})$ commutes with all limits.
It follows from Lemma
\ref{lem:tensor-prod-pres-cat} 
below that there is an 
equivalence of functors
$$\mathcal{X}(-;\mathcal{S})\otimes \mathcal{V}
\simeq \Fun^{\lim}(\mathcal{V}^{\op},\mathcal{X}(-;\mathcal{S})).$$
Thus, it is enough to show that 
$\Fun^{\lim}(\mathcal{V}^{\op},\mathcal{X}(-;\mathcal{S})):
\mathcal{C}^{\op}\to \CAT_{\infty}$
is a $\tau$-(hyper)sheaf. 
This follows from Proposition
\ref{prop:sheaf-values-topos}
since the endofunctor
$\Fun^{\lim}(\mathcal{V}^{\op},-)$ of $\CAT_{\infty}$ 
preserves limits.
\end{proof}

\begin{lemma}
\label{lem:tensor-prod-pres-cat}
Let ${\rm Pr^{LR}}$ be the wide sub-$\infty$-category of 
$\Prl$ where morphisms are the limit-preserving left adjoints.
Let $\mathcal{D}$ be a presentable $\infty$-category. 
Then the functor
$\mathcal{D}\otimes-:{\rm Pr^{LR}}\to \CAT_{\infty}$, 
obtained by restriction from the tensor product of $\Prl$,
is equivalent to the functor 
$$\Fun^{\lim}(\mathcal{D}^{\op};-):{\rm Pr^{LR}}\to \CAT_{\infty},$$
where $\Fun^{\lim}(-,-)\subset \Fun(-,-)$
indicates the full sub-$\infty$-category of limit-preserving functors.
\end{lemma}

\begin{proof}
The endofunctor $\mathcal{D}\otimes-$ of $\Prl$
induces an endofunctor of $\Prr$ given by the composition of
$$\Prr \xrightarrow{\sim} (\Prl)^{\op}
\xrightarrow{\mathcal{D}\otimes-}
(\Prl)^{\op}
\xrightarrow{\sim} \Prr.$$
By \cite[Proposition 4.8.1.17]{lurie:higher-algebra}, 
this coincides with the endofunctor 
$\Fun^{\lim}(\mathcal{D}^{\op},-)$ of $\Prr$.
It follows that the endofunctor $\mathcal{D}\otimes-$ 
of $\Prl$ is given by the composition of 
$$\Prl \xrightarrow{\sim} (\Prr)^{\op}
\xrightarrow{\Fun^{\lim}(\mathcal{D}^{\op},-)}
(\Prr)^{\op} \xrightarrow{\sim} \Prl.$$
It remains to show that the composition of 
$${\rm Pr^{LR}}
\to \Prl
\xrightarrow{\sim} (\Prr)^{\op}
\xrightarrow{\Fun^{\lim}(\mathcal{D}^{\op},-)}
(\Prr)^{\op}
\xrightarrow{\sim} 
\Prl \to \CAT_{\infty}$$
is also given by $\Fun^{\lim}(\mathcal{D}^{\op},-)$.
On objects, this is clear. On morphisms, this is also true by
the following observation: if  
$F:\mathcal{E}\to \mathcal{E}'$ is in ${\rm Pr^{LR}}$ with 
right adjoint $G$, then $\Fun^{\lim}(\mathcal{D}^{\op},F)$
is left adjoint to $\Fun^{\lim}(\mathcal{D}^{\op},G)$. To address 
higher coherences, we employ the formalism of Cartesian fibrations.

Let $S$ be a simplicial set and $p:\mathcal{M} \to S$ 
a coCartesian fibration classified by a map $l:S\to {\rm Pr^{LR}}$. 
Then $p$ is also a Cartesian fibration which is classified by a map
$r:S\to (\Prr)^{\op}$ equivalent to the composition of 
$$S \xrightarrow{l} \Prl \xrightarrow{\sim} (\Prr)^{\op}.$$
Moreover, $p$-Cartesian 
and $p$-coCartesian edges of $\mathcal{M}$ 
are preserved by small limits in the following sense. 
Let $a:s\to s'$ be an edge 
in $S$, $\overline{e}:K^{\lhd}\to \mathcal{M}_s$ and 
$\overline{e}':K^{\lhd}\to \mathcal{M}_{s'}$ limit diagrams, and 
$f:\overline{e}\to \overline{e}'$ an edge in 
$\Fun(K^{\lhd},\mathcal{M})$ over $a$. If $f(k)$ is 
$p$-coCartesian (resp. $p$-Cartesian) for every $k\in K$, 
then the same is true for $f(\infty)$, 
where $\infty\in K^{\rhd}$ is the cone point.
This is simply a reformulation of the fact that $l$ 
(resp. $r$) takes an edge of $S$ to a limit-preserving functor.
Consider the simplicial set
$\mathcal{N}=\mathcal{M}^{\mathcal{D}^{\op}}
\times_{S^{\mathcal{D}^{\op}}} S$
whose $n$-simplices correspond to pairs consisting of 
an $n$-simplex $[n]\to S$ and an $S$-morphism 
$[n]\times\mathcal{D}^{\op}\to \mathcal{M}$.
Let $\mathcal{N}'\subset \mathcal{N}$ be the largest simplicial 
subset whose vertices correspond to 
limit-preserving functors 
$\mathcal{D}^{\op}\to \mathcal{M}_s$, for some $s\in S$.
Let $q:\mathcal{N}\to S$ and $q':\mathcal{N}'\to S$
be the obvious projections.
By \cite[Proposition 2.4.2.3(2) \& Proposition 3.1.2.1]{lurie}, 
$q$ is again a coCartesian fibration, classified by 
$\Fun(\mathcal{D}^{\op},-)\circ l:S\to \CAT_{\infty}$, and 
a Cartesian fibration classified by 
$\Fun(\mathcal{D}^{\op},-)\circ r:S \to (\CAT_{\infty})^{\op}$.
Since $p$-coCartesian (resp. $p$-Cartesian) edges are preserved by 
small limits, it follows readily that a $q$-coCartesian 
(resp. $q$-Cartesian) edge whose domain (resp. target) 
belongs to $\mathcal{N}'$ lies 
entirely in $\mathcal{N}'$. This shows that
$q'$ is a coCartesian 
fibration, classified by 
$l'=\Fun^{\lim}(\mathcal{D}^{\op},-)\circ l:S\to \CAT_{\infty}$, 
and a Cartesian fibration classified by 
$r'=\Fun^{\lim}(\mathcal{D}^{\op},-)\circ r:S \to (\CAT_{\infty})^{\op}$.
It follows that $l'$ factors through $\CAT_{\infty}^{\rm L}$, 
$r'$ factors through $\CAT_{\infty}^{\rm R}$, and 
$l'$ coincides with the composition of 
$$\mathcal{S}^{\op}\xrightarrow{r'}  (\CAT_{\infty}^{\rm R})^{\op}
\xrightarrow{\sim} \CAT_{\infty}^{\rm L}.$$ 
Unravelling the definitions, this gives what we want.
\end{proof}

\begin{proof}[Proof of Theorem \ref{thm:hyperdesc}]
It suffices to prove that for every rigid analytic space $S$, the functor 
$$\RigSH_{\tau}^{(\eff,\,\hyp)}(-;\Lambda):(\Et/S)^{\op}\to\Prl,$$
is a $\tau$-(hyper)sheaf. (When $\tau=\Nis$, one
can restrict further to $(\Etgr/S)^{\op}$, but this does not change
the argument.)
This functor transforms coproducts in $\Et/S$ into
products in $\Prl$. Thus, it suffices to show that it admits descent 
with respect to internal hypercovers of $(\Et/S,\tau)$,
which are truncated in the non-hypercomplete case.

For $U\in \Et/S$, we have $(\RigSm/S)/U\simeq \RigSm/U$.
Corollary \ref{cor:sheafy-sheaves} implies that the functor 
$$\Shv^{(\hyp)}_{\tau}(\RigSm/-;\Lambda):
(\Et/S)^{\op}\to \Prl$$
is a $\tau$-(hyper)sheaf. Let $U_{\bullet}$ 
be an internal hypercover of $(\Et/S,\tau)$ which we assume 
to be truncated in the non-hypercomplete case.
For all $n\geq -1$, $\RigSH^{\eff,\,(\hyp)}_{\tau}(U_n;\Lambda)$ is a full 
sub-$\infty$-category of $\Shv^{(\hyp)}_{\tau}(\RigSm/U_n;\Lambda)$. 
Since limits in $\CAT_{\infty}$ preserve fully faithful embeddings, 
we deduce that 
$\lim_{[n]\in \Delta} \RigSH_{\tau}^{\eff,\,(\hyp)}(U_n;\Lambda)$
can be naturally identified with the sub-$\infty$-category of 
$$\Shv^{(\hyp)}_{\tau}(\RigSm/U_{-1};\Lambda)\simeq 
\lim_{[n]\in \Delta}
\Shv^{(\hyp)}_{\tau}(\RigSm/U_n;\Lambda)$$
spanned by the objects $\mathcal{F}\in 
\Shv^{(\hyp)}_{\tau}(\RigSm/U_{-1};\Lambda)$
such that $f^*\mathcal{F}$ belongs to 
$\RigSH_{\tau}^{\eff,\,(\hyp)}(U_0;\Lambda)$, with $f:U_0\to U_{-1}$. 
Thus, to prove that 
$\RigSH_{\tau}^{\eff,\,(\hyp)}(-;\Lambda)$ has descent for the 
$\tau$-hypercover $U_{\bullet}$, we need to check the following property:
if $\mathcal{F}$ is a $\tau$-(hyper)sheaf on $\RigSm/S$ such 
that $f^*\mathcal{F}$ is $\B^1$-invariant, then so is $\mathcal{F}$.
This follows immediately from the equivalence 
$\underline{\Hom}(\B^1_{U_0},f^*\mathcal{F})\simeq 
f^*\underline{\Hom}(\B^1_{U_{-1}},\mathcal{F})$ 
and the fact that $f^*$ is conservative.

We now explain how to deduce the $\Tate$-stable
case from the
effective case. We temporarily denote by 
$\RigSH_{\tau}^{(\eff,\,\hyp)}(-;\Lambda)^*$
(resp. $\RigSH_{\tau}^{(\eff,\,\hyp)}(-;\Lambda)_*$)
the presheaf (resp. copresheaf) given informally by $U\mapsto 
\RigSH_{\tau}^{(\eff,\,\hyp)}(-;\Lambda)$ and $f\mapsto f^*$
(resp. $f\mapsto f_*$). Recall from Remark
\ref{rmk:symmetric-obj-} 
that the presheaf 
$\RigSH^{(\hyp)}_{\tau}(-;\Lambda)^*$ 
can be defined as the colimit in $\PSh(\Et/S;\Prl)$
of the $\N$-diagram of presheaves: 
$$\RigSH_{\tau}^{\eff,\,(\hyp)}(-;\Lambda)^* \xrightarrow{\Tate\otimes-}
\RigSH_{\tau}^{\eff,\,(\hyp)}(-;\Lambda)^* 
\xrightarrow{\Tate\otimes-}\cdots.$$ 
It follows from \cite[Corollary 5.5.3.4 \&Theorem 5.5.3.18]{lurie} that the copresheaf $\RigSH^{(\hyp)}_{\tau}(U;\Lambda)_*$, can be computed 
as the limit in $\Fun(\Et/S,\CAT_{\infty})$ of the $\N^{\op}$-diagram of copresheaves
$$\cdots \xrightarrow{\underline{\Hom}(\Tate,-)}
\RigSH_{\tau}^{\eff,\,(\hyp)}(-;\Lambda)_*
\xrightarrow{\underline{\Hom}(\Tate,-)} 
\RigSH_{\tau}^{\eff,\,(\hyp)}(-;\Lambda)_*.$$
Given that the natural transformation 
$f^*\underline{\Hom}(\Tate,-)\to \underline{\Hom}(\Tate,-)\circ f^*$
is an equivalence for $f$ \'etale, we deduce that the 
presheaf $\RigSH^{(\hyp)}_{\tau}(-;\Lambda)^*$ 
can also be computed as the limit in
$\PSh(\Et/S;\CAT_{\infty})$ of the $\N^{\op}$-diagram of presheaves
$$\cdots \xrightarrow{\underline{\Hom}(\Tate,-)}
\RigSH_{\tau}^{\eff,\,(\hyp)}(-;\Lambda)^*
\xrightarrow{\underline{\Hom}(\Tate,-)} 
\RigSH_{\tau}^{\eff,\,(\hyp)}(-;\Lambda)^*.$$
Since $\RigSH_{\tau}^{\eff,\,(\hyp)}(-;\Lambda)^*$ was proven to be a 
$\tau$-(hyper)sheaf, this finishes the proof.
\end{proof}

\subsection{Compact generation}

$\empty$

\smallskip

In this subsection, we formulate conditions 
(in terms of $\Lambda$, $S$ and $\tau$) insuring that the 
$\infty$-category $\RigSH^{(\eff,\,\hyp)}_{\tau}(S;\Lambda)$
of rigid analytic motives over $S$ is compactly 
generated. Similar results in the algebraic setting 
were developed in \cite[\S 4.5.5]{ayoub-th2} and 
\cite[pages 29--30]{ayoub-etale}.

\begin{rmk}
\label{rmk:discrete-sheaf-of-modules}
Let $\mathcal{X}$ be an $\infty$-topos.
An abelian group object of $\mathcal{X}_{\leq 0}$ 
endowed with the structure of a $\pi_0\Lambda$-module
is called a discrete sheaf of $\pi_0\Lambda$-modules on $\mathcal{X}$.
The $n$-th cohomology group of $\mathcal{X}$ with coefficients in 
a discrete sheaf of $\pi_0\Lambda$-modules $\mathcal{F}$ 
is defined in \cite[Definition 7.2.2.14]{lurie} 
and will be denoted by 
$\H^n(\mathcal{X};\mathcal{F})$.
\symn{$\H^n$}
\end{rmk}

Recall the following notions. (Compare with
\cite[Definition 7.2.2.18]{lurie}.)

\begin{dfn}
\label{dfn:coh-dim-topos}
Let $\mathcal{X}$ be an $\infty$-topos.

\begin{enumerate}

\item[(1)] The $\Lambda$-cohomological dimension of 
an object $X\in \mathcal{X}$ is the smallest 
$d\in \N\sqcup\{-\infty,\infty\}$ 
such that for every discrete sheaf of $\pi_0\Lambda$-modules 
$\mathcal{F}$ on $\mathcal{X}_{/X}$, the cohomology groups
$\H^n(\mathcal{X}_{/X};\mathcal{F})$ vanish for $n>d$.
The global $\Lambda$-cohomological dimension of $\mathcal{X}$ 
is the $\Lambda$-cohomological dimension of a final object of 
$\mathcal{X}$.
\ncn{cohomological dimension!in an $\infty$-topos}

\item[(2)] The local $\Lambda$-cohomological dimension of 
$\mathcal{X}$ is the smallest $d\in \N\sqcup\{-\infty,\infty\}$
such that every object $X\in \mathcal{X}$ admits a cover
$(Y_i\to X)_i$ such that $Y_i$ is of 
$\Lambda$-cohomological dimension 
$\leq d$ for all $i$. (Recall that 
$(Y_i\to X)_i$ is a cover if $\coprod_i Y_i \to X$ is an 
effective epimorphism in the sense of \cite[\S 6.2.3]{lurie}.)
\ncn{cohomological dimension!local}

\end{enumerate}
\end{dfn}

\begin{rmk}
\label{rmk:coh-dim-hypercomp}
Keep the notation as in Definition 
\ref{dfn:coh-dim-topos}.
A discrete sheaf of $\pi_0\Lambda$-modules $\mathcal{F}$ 
on $\mathcal{X}_{/X}$ is a hypersheaf, i.e., belongs to 
$(\mathcal{X}_{/X})^{\hyp}\simeq (\mathcal{X}^{\hyp})_{/X^{\hyp}}$.
Thus, there are isomorphisms
$$\H^i(\mathcal{X}_{/X};\mathcal{F})\simeq 
\H^i(\mathcal{X}_{/X^{\hyp}};\mathcal{F})\simeq 
\H^i((\mathcal{X}^{\hyp})_{/X^{\hyp}};\mathcal{F}).$$
In particular, the $\Lambda$-cohomological dimension of an 
object $X$ is equal to 
the $\Lambda$-cohomological dimension of its hypercompletion $X^{\hyp}$
considered as an object of $\mathcal{X}$ or $\mathcal{X}^{\hyp}$.
Similarly, the global (resp. local)
$\Lambda$-cohomological dimensions
of $\mathcal{X}$ and $\mathcal{X}^{\hyp}$
coincide.
\end{rmk}

\begin{rmk}
\label{rmk:coh-dim-sites}
We define the local (resp. global) 
$\Lambda$-cohomological dimension of a site
$(\mathcal{C},\tau)$ to be the local (resp. global) 
$\Lambda$-cohomological dimension of the topos
$\Shv_{\tau}(\mathcal{C})$ (or, equivalently, 
$\Shv_{\tau}^{\hyp}(\mathcal{C})$). Similarly, we define 
the $\Lambda$-cohomological dimension of an 
object $X$ of a site $(\mathcal{C},\tau)$ 
to be the $\Lambda$-cohomological 
dimension of the image of $X$ in $\Shv_{\tau}(\mathcal{C})$
(or, equivalently, 
$\Shv_{\tau}^{\hyp}(\mathcal{C})$). 
By Lemma \ref{slice-topoi}, 
this coincides with the global 
$\Lambda$-cohomological dimension of the 
site $(\mathcal{C}_{/X},\tau)$.
\ncn{cohomological dimension!of a site}
\end{rmk}

We gather some well-known consequences of the finiteness of the
local $\Lambda$-cohomological dimension in the following statement.
(See Remark \ref{rmk:t-structure-connect-sheaves}.)

\begin{lemma}
\label{lem:pi-0-Lambda-coh-dim}
Let $(\mathcal{C},\tau)$ be a site of finite local 
$\Lambda$-cohomological dimension. 
\begin{enumerate}

\item[(1)] Postnikov towers in $\Shv^{\hyp}_{\tau}(\mathcal{C};\Lambda)$
converge, i.e., the obvious map
$$\mathcal{F}\to \lim_{n\in \N} 
\tau_{\leq n}\mathcal{F}$$
is an equivalence for every $\tau$-hypersheaf of 
$\Lambda$-modules $\mathcal{F}$ on $\mathcal{C}$.

\item[(2)] If $\mathcal{F}$ is a connective $\tau$-hypersheaf of 
$\Lambda$-modules on $\mathcal{C}$ and
$X\in \mathcal{C}$ is of $\Lambda$-cohomological dimension $\leq d$, 
then the $\Lambda$-module $\mathcal{F}(X)$ is $(-d)$-connective. 

\item[(3)] Assume that $\mathcal{C}$ is an ordinary 
category admitting fiber products and that every 
object of $\mathcal{C}$ is quasi-compact in the sense of 
\cite[Expos\'e VI, D\'efinitions 1.1]{SGAIV2}.
If $X\in \mathcal{C}$ 
is of finite $\Lambda$-cohomological dimension,
then the functor $\Shv^{\hyp}_{\tau}
(\mathcal{C};\Lambda)\to\Mod_{\Lambda}$,
$\mathcal{F}\mapsto \mathcal{F}(X)$ 
commutes with arbitrary colimits. 
In particular, $\Lambda_{\tau}(X)$ is a compact object of 
$\Shv^{\hyp}_{\tau}(\mathcal{C};\Lambda)$.

\end{enumerate}
\end{lemma}

\begin{proof}
We may replace $(\mathcal{C},\tau)$ with any site that gives 
rise to the same hypercomplete topos. Thus, we may assume that
every object of $\mathcal{C}$ has $\Lambda$-cohomological dimension 
$\leq d$. Property (2), for every object $X\in \mathcal{C}$,
follows from \cite[Proposition 4.5.58]{ayoub-th2}
when $(\mathcal{C},\tau)$ is an ordinary site and 
$\Lambda$ the unit spectrum. However, the proof of 
loc.~cit.~can be adapted without difficulty to our setting.
That proof gives also property (1). (Note that (1) can 
be deduced from (2), but usually these two properties are
proven together.)
Since $\mathcal{F}\mapsto \mathcal{F}(X)$ is an exact functor
between stable $\infty$-categories, it preserves pushouts. 
By \cite[Proposition 4.4.2.7]{lurie}, to prove property (3)
it is enough to show that this functor commutes with filtered colimits.
This follows from property (2) as in the proof 
of \cite[Corollaire 4.5.61]{ayoub-th2}.
(The extra conditions on $\mathcal{C}$ 
are used via \cite[Expos\'e VI, Corollaire 5.3]{SGAIV2} and can be 
substantially weakened.)

For a modern and more general treatment of this type of question,
we refer the reader to \cite[\S2]{cla-mat:etale-k-th}.
In particular, property (1) follows from 
\cite[Proposition 2.10]{cla-mat:etale-k-th}
(see also \cite[Example 2.11]{cla-mat:etale-k-th}).
Property (3) can be deduced from 
\cite[Proposition 2.23]{cla-mat:etale-k-th}.
Finally, we mention  
\cite[Proposition 7.2.1.10]{lurie},
which is obviously related to property (1).
\end{proof}

\begin{cor}
\label{cor:automatic-hypercomp}
Let $(\mathcal{C},\tau)$ be a site, and assume the following 
conditions:
\begin{enumerate}

\item[(1)] $\Lambda$ is eventually coconnective (i.e., 
its homotopy groups $\pi_i\Lambda$ 
vanish for $i$ big enough);

\item[(2)] $(\mathcal{C},\tau)$ has finite 
local $\Lambda$-cohomological dimension and
$\mathcal{C}$ is an ordinary category with fiber products;

\item[(3)] there exists a full subcategory 
$\mathcal{C}_0\subset \mathcal{C}$ stable under fiber products,
spanned by quasi-compact objects of finite $\Lambda$-cohomological
dimension, and such that every 
object of $\mathcal{C}$ admits a $\tau$-cover by 
objects of $\mathcal{C}_0$.

\end{enumerate}
Then every $\tau$-sheaf of $\Lambda$-modules on $\mathcal{C}$
is a $\tau$-hypersheaf, i.e., 
we have $\Shv_{\tau}^{\hyp}(\mathcal{C};\Lambda)=
\Shv_{\tau}(\mathcal{C};\Lambda)$.
\end{cor}

\begin{proof}
By Lemma \ref{lem:equi-of-sites-infty-topoi}, we may replace 
$\mathcal{C}$ with $\mathcal{C}_0$ and assume that 
every object of $\mathcal{C}$ is quasi-compact, quasi-separated 
and of finite $\Lambda$-cohomological dimension. 
For $X\in \mathcal{C}$, the $\tau$-sheaf 
$\Lambda_{\tau}(X)$ is hypercomplete since 
$\Lambda$ is eventually coconnective. Thus, it is enough to
show that $\tau$-hypersheaves are stable under colimits in 
$\Shv_{\tau}(\mathcal{C};\Lambda)$.
The result then follows from \cite[Proposition 2.23]{cla-mat:etale-k-th}
but we can also deduce it formally from Lemma 
\ref{lem:pi-0-Lambda-coh-dim} as follows.
Indeed, let $p:K\to \Shv_{\tau}^{\hyp}(\mathcal{C};\Lambda)$
be a diagram of $\tau$-hypersheaves of $\Lambda$-modules. 
The colimit of $p$ in $\Shv_{\tau}(\mathcal{C};\Lambda)$ 
is the $\tau$-sheafification of the colimit of $p$ in 
$\PSh(\mathcal{C};\Lambda)$. So it is enough to show that the 
colimit of $p$ in $\PSh(\mathcal{C};\Lambda)$ is already a 
$\tau$-hypersheaf. This follows immediately from 
Lemma \ref{lem:pi-0-Lambda-coh-dim}(3).
\end{proof}

We now give some estimates for the local and global 
$\Lambda$-cohomological dimensions of the various small sites
associated to a rigid analytic space.

\begin{lemma}
\label{lem:Nis-Krull-dim}
Let $X$ be a rigid analytic space of Krull dimension $\leq d$. 
The local $\Lambda$-cohomological dimension of $(\Etgr/X,\Nis)$
is $\leq d$. If $X$ is quasi-compact and quasi-separated, the
same is true for the global 
$\Lambda$-cohomological dimension.
\end{lemma}

\begin{proof}
Since every object of $\Etgr/X$ can be covered by quasi-compact
and quasi-separated rigid analytic spaces of Krull dimension
$\leq d$, it is enough to prove the assertion concerning the 
global $\Lambda$-cohomological dimension.
In particular, 
we may assume that $X$ is quasi-compact and quasi-separated. 
The site $(\Etgr/X,\Nis)$ is then equivalent to the limit of the
Nisnevich sites $(\Et/\mathcal{X}_{\sigma},\Nis)$,
for $\mathcal{X}\in \Mdl'(X)$ (see Remark
\ref{rmk:essential-form-sch}). 
It follows from \cite[Expos\'e VII, Th\'eor\`eme 5.7]{SGAIV2}
that the global $\Lambda$-cohomological dimension of the site 
$(\Etgr/X,\Nis)$
is smaller than the supremum of the global 
$\Lambda$-cohomological dimensions of the sites
$(\Et/\mathcal{X}_{\sigma},\Nis)$, for $\mathcal{X}\in \Mdl'(X)$.
But if $\mathcal{X}$ is a formal model of $X$ belonging to 
$\Mdl'(X)$,
the closed map $|X|\to |\mathcal{X}_{\sigma}|$ is surjective.
Thus, the dimension of $\mathcal{X}_{\sigma}$ 
is smaller than the dimension of $X$, and we conclude using 
\cite[Theorem 3.17]{cla-mat:etale-k-th}.
\end{proof}

\begin{dfn}
\label{dfn:virtual-coh-dim-field}
Let $G$ be a profinite group. The $\Lambda$-cohomological dimension 
of $G$ is the smallest $d\in \N\sqcup\{\infty\}$ such that, for 
every $\pi_0\Lambda$-module $M$ endowed with a continuous action of $G$, 
the cohomology groups $\H^i(G;M)$ vanish for $i>d$. The virtual 
$\Lambda$-cohomological dimension of $G$ is the infimum of the 
$\Lambda$-cohomological dimensions of the finite-index subgroups
of $G$. If $G$ admits a finite-index torsion-free
subgroup $H$, then 
the virtual $\Lambda$-cohomological dimension of $G$ is equal to the 
$\Lambda$-cohomological dimension of $H$.
(See \cite[Chapitre I, \S 3.3, Proposition 14$'$]{CohGal}.)

Let $k$ be a field with absolute Galois group $G_k$. 
The (virtual) $\Lambda$-cohomological dimension of $k$ is defined
to be the  (virtual) $\Lambda$-cohomological dimension of $G_k$.
\ncn{cohomological dimension!of a group}
\ncn{cohomological dimension!of a field}
\ncn{cohomological dimension!virtual}
\end{dfn}

\begin{rmk}
\label{rmk:about-cohgal}
Let $k$ be a field.
The following are classical facts about 
Galois cohomology.
\begin{enumerate}

\item[(1)] If the $\Lambda$-cohomological dimension of $k$ 
is different from its virtual 
$\Lambda$-cohomological dimension, then $k$ admits a real embedding
and $2$ is not invertible in $\pi_0\Lambda$.

\item[(2)] If $k$ has (virtual) $\Lambda$-cohomological dimension
$\leq d$ and $K/k$ is an extension of transcendence degree $\leq e$, 
then $K$ has (virtual) $\Lambda$-cohomological dimension
$\leq d+e$.

\item[(3)] Number fields have 
virtual $\Lambda$-cohomological dimension $\leq 3$, and
finite fields have $\Lambda$-cohomological dimension
$\leq 2$.

\end{enumerate}
Property (1) follows from 
\cite[Chapitre II, \S 4.1, Proposition 10$'$]{CohGal}.
Property (2) follows from 
\cite[Chapitre II, \S 4.2, Proposition 11]{CohGal}.
Property (3) follows from 
\cite[Chapitre II, \S 4.4, Proposition 13]{CohGal}.
\end{rmk}

\begin{dfn}
\label{dfn:Lambda-coh-dim-4}
Let $X$ be a scheme or a rigid analytic space. We denote by 
$\pvcd_{\Lambda}(X)\in \N\sqcup\{-\infty,\infty\}$ 
the supremum
of the virtual $\Lambda$-cohomological dimensions of 
the fields $\kappa(x)$ for $x\in |X|$. This number is called the 
punctual virtual $\Lambda$-cohomological dimension of $X$.
\symn{$\pvcd$}
\end{dfn}

\begin{lemma}
\label{lem:Lamba-tau-coh-pointwise}
Let $X$ be a rigid analytic space of Krull dimension $\leq d$ 
and of punctual virtual $\Lambda$-cohomological dimension $\leq e$.
Then, the local $\Lambda$-cohomological dimension of the site 
$(\Et/X,\et)$ is $\leq d+e$. The same is true for the global 
$\Lambda$-cohomological dimension if $X$ is quasi-compact and 
quasi-separated, and if the 
$\Lambda$-cohomological dimension of the residue field of 
every point of $X$ coincides with the virtual one.
\end{lemma}

\begin{proof}
Replacing $X$ by a suitable \'etale cover
(e.g., by $X[\frac{1}{2},\sqrt{-1}]\to X$ and 
$X[\frac{1}{3},\sqrt[3]{1}] \to X$), 
we may assume that the $\Lambda$-cohomological 
dimension of the residue field of 
each point of $X$ coincides with the virtual one.
We may also assume that $X$ is quasi-compact and quasi-separated.
Under these conditions, we will show that the global 
$\Lambda$-cohomological dimension of $(\Et/X,\et)$
is $\leq d+e$, which suffices to conclude.

Denote by $\pi:(\Et/X,\et) \to 
(\Etgr/X,\Nis)$ the obvious morphism of sites.
Given an \'etale sheaf $\mathcal{F}$ of $\pi_0\Lambda$-modules
on $\Et/X$, we denote by $\Rder\pi_*\mathcal{F}$ its (derived)
direct image. Using Lemma 
\ref{lem:Nis-Krull-dim}, we are reduced to showing that 
$\Rder\pi_*\mathcal{F}$ is $(-e)$-connective.
We check this on stalks at Nisnevich geometric rigid points of $X$
as in Construction \ref{cons:point-an-nis-et}.
Let $s\in S$ be a point and 
$t\to S$ a Nisnevich geometric rigid point over $s$. 
Thus, $t=\Spf(\kappa^+(t))$ with $\kappa^+(t)$ the 
adic completion of the Henselisation of $\kappa^+(s)$ at
a morphism $\Spec(\widetilde{\kappa}(t)) \to 
\Spec(\kappa^+(s))$
associated to a separable finite extension 
$\widetilde{\kappa}(t)/\widetilde{\kappa}(s)$.
It follows from Corollary 
\ref{cor:limit-etale-topi-rig-an}
that $(\Rder\pi_*\mathcal{F})_t$ is equivalent to 
$\Rder\Gamma_{\et}(t;(t\to S)^*\mathcal{F})$.
Thus, it is sufficient to show that the global 
$\Lambda$-cohomological dimension of $(\Et/t,\et)$
is smaller than $e$. Since $\kappa^+(t)$ is Henselian, 
every \'etale cover of $t$ can be refined by one of the 
form $\Spf(V)^{\rig} \to t$ where 
$V$ is the normalisation of $\kappa^+(t)$ 
in a finite separable extension of $\kappa(t)$.
Thus, the global cohomology of $(\Et/t,\et)$
coincides with the Galois cohomology of  
$\kappa(t)$. Since the field $\kappa(t)$ is the completion of 
an algebraic extension of $\kappa(s)$, we deduce that 
its $\Lambda$-cohomological dimension is $\leq e$ as needed. 
\end{proof}

The following is a corollary of the proof of Lemma 
\ref{lem:Lamba-tau-coh-pointwise}.

\begin{cor}
\label{cor:nis-etale-comparison}
Let $X$ be a rigid analytic space, and
let $\mathcal{F}$ be a discrete sheaf of $\Q$-vector spaces 
on $(\Et/X,\et)$. Then the natural map
$\H^*_{\Nis}(X;\mathcal{F}) \to \H^*_{\et}(X;\mathcal{F})$
is an isomorphism. 
\end{cor}

\begin{proof}
Arguing as in the proof of Lemma \ref{lem:Lamba-tau-coh-pointwise},
the result follows from the vanishing of the higher 
Galois cohomology groups with rational coefficients.
\end{proof}

\begin{cor}
\label{cor:fin-coh-dim-et-rig-space}
Let $X$ be a rigid analytic space of Krull dimension $\leq d$.
If $\Lambda$ is a $\Q$-algebra, then the local 
$\Lambda$-cohomological dimension of the site 
$(\Et/X,\et)$ is $\leq d$. 
If $X$ is quasi-compact and quasi-separated, the
same is true for the global 
$\Lambda$-cohomological dimension.
\end{cor}

\begin{dfn}
\label{dfn:lambda-tau-admiss}
Let $S$ be a scheme or a rigid analytic space. 
\begin{enumerate}

\item[(1)] We say that $S$ is $(\Lambda,\et)$-admissible if
there exists an open covering $(S_i)_i$ of $S$ such that 
each $S_i$ has finite Krull dimension and finite 
punctual virtual $\Lambda$-cohomological dimension.
For convenience, we also say that 
$S$ is $(\Lambda,\Nis)$-admissible when 
$S$ is locally of finite Krull dimension.
\ncn{schemes!$(\Lambda,\tau)$-admissible}
\ncn{rigid analytic spaces!$(\Lambda,\tau)$-admissible}

\item[(2)] If $2$ is not invertible in $\pi_0\Lambda$, 
we say that $S$ is $(\Lambda,\et)$-good if $\mathcal{O}(S)$ 
contains a primitive $n$-th root of unity for some $n\geq 3$.
For convenience, we agree that 
$S$ is always $(\Lambda,\tau)$-good if $2$ is invertible 
in $\pi_0\Lambda$ or if $\tau$ is the Nisnevich topology.
\ncn{schemes!$(\Lambda,\tau)$-good}
\ncn{rigid analytic spaces!$(\Lambda,\tau)$-good}

\end{enumerate}
\end{dfn}

\begin{rmk}
\label{rmk:good-virtual-coh-dim}
If $S$ is $(\Lambda,\et)$-good, then the $\Lambda$-cohomological dimension of the residue field of each of its points coincides with the virtual one.
This follows from Remark \ref{rmk:about-cohgal}.
\end{rmk}

\begin{lemma}
\label{lem:Lambda-coh-dim-val-ring}
Let $Y \to X$ be a morphism of
rigid analytic spaces which is locally of finite type, 
and let $y\in Y$ be a point 
with image $x\in X$. If the (virtual) 
$\Lambda$-cohomological dimension
of $\kappa(x)$ is finite, then so is the (virtual) 
$\Lambda$-cohomological 
dimension of $\kappa(y)$.
\end{lemma}

\begin{proof}
We use the fact that $\kappa(y)/\kappa(x)$ 
is topologically of finite type, i.e., 
that $\kappa(y)$ is the completion of a finite type 
extension of $\kappa(x)$.  
It follows that the absolute Galois group of 
$\kappa(y)$ can be identified with a closed subgroup of the 
absolute Galois group of a finite type extension of $\kappa(y)$. 
We then conclude using 
Remark \ref{rmk:about-cohgal}(2). Alternatively, 
one can deduce the result from
\cite[Lemma 2.8.4]{huber}.
\end{proof}

\begin{cor}
\label{cor:univers-Lambda-admis}
Let $\tau\in \{\Nis,\et\}$. Let $f:T \to S$ be a morphism 
of rigid analytic spaces which is locally of finite type. 
If $S$ is $(\Lambda,\tau)$-admissible, then so is $T$.
\end{cor}

\begin{proof}
This follows immediately from 
Lemma \ref{lem:Lambda-coh-dim-val-ring}.
\end{proof}

\begin{lemma}
\label{lem:auto-hypercomp-etale}
Let $\tau\in\{\Nis,\et\}$ and let $S$ be a $(\Lambda,\tau)$-admissible
rigid analytic space.
\begin{enumerate}

\item[(1)] (Case $\tau=\Nis$) Every Nisnevich sheaf of $\Lambda$-modules
on $\Etgr/S$ is a Nisnevich hypersheaf, i.e., we have 
$$\Shv_{\Nis}^{\hyp}(\Etgr/S;\Lambda)=\Shv_{\Nis}(\Etgr/S;\Lambda).$$
The same statement is true with ``$\Etgr/S$'' replaced with 
``$\Et/S$'' or ``$\RigSm/S$''.

\item[(2)] (Case $\tau=\et$) Assume that $\Lambda$ is eventually 
coconnective. Then every \'etale sheaf of $\Lambda$-modules 
on $\Et/S$ is an \'etale hypersheaf, i.e., we have 
$$\Shv_{\et}^{\hyp}(\Et/S;\Lambda)=
\Shv_{\et}(\Et/S;\Lambda).$$
The same statement is true with ``$\Et/S$'' replaced with 
``$\RigSm/S$''.

\end{enumerate}
\end{lemma}

\begin{proof}
If $\mathcal{F}$ is a $\tau$-sheaf of $\Lambda$-modules on 
$\RigSm/S$ whose restriction to 
$\Et/X$ (or $\Etgr/X$ if applicable) 
is a $\tau$-hypersheaf for every quasi-compact and 
quasi-separated $X\in \RigSm/S$, then $\mathcal{F}$ is
a $\tau$-hypersheaf. (Indeed, if this holds, the morphism
$\mathcal{F} \to \mathcal{F}^{\hyp}$ induces equivalences 
$\mathcal{F}(X) \simeq \mathcal{F}^{\hyp}(X)$ for every
$X\in \RigSm^{\qcqs}/S$, so it is itself an equivalence.) 
Therefore, using Corollary 
\ref{cor:univers-Lambda-admis},
it is enough to treat the cases of the small sites of $S$, 
with $S$ quasi-compact and quasi-separated.
The case of $(\Et/S,\et)$ follows then from 
Corollary \ref{cor:automatic-hypercomp} and Lemma
\ref{lem:Lamba-tau-coh-pointwise}.
The case of $(\Etgr/S,\Nis)$ needs a special treatment. 
For this, we remark that if $(X_{\alpha})_{\alpha}$ is
a cofiltered inverse system of quasi-compact and quasi-separated 
schemes of dimension $\leq d$ (with $d$ independent of $\alpha$),
then the proof of \cite[Theorem 3.17]{cla-mat:etale-k-th}
can be adapted to show that the site 
$\lim_{\alpha} (\Et/X_{\alpha},\Nis)$
is locally of homotopy dimension $\leq d$, which implies that the 
associated topos is hypercomplete by 
\cite[Corollary 7.2.1.12]{lurie}. Applying this to the 
inverse system $(\mathcal{S}_{\sigma})_{\mathcal{S}\in \Mdl'(S)}$
gives the result.
\end{proof}

\begin{prop}
\label{prop:automatic-hypercomp-motives}
Let $\tau\in\{\Nis,\et\}$ and let $S$ be a $(\Lambda,\tau)$-admissible
rigid analytic space. When $\tau$ is the \'etale topology, assume 
that $\Lambda$ is eventually coconnective. Then, we have
$$\RigSH_{\tau}^{(\eff),\,\hyp}(S;\Lambda)=
\RigSH_{\tau}^{(\eff)}(S;\Lambda).$$ 
\end{prop}

\begin{proof}
This follows immediately from Lemma
\ref{lem:auto-hypercomp-etale}.
\end{proof}

\begin{prop}
\label{prop:compact-sheaves-qcqs}
Let $\tau\in\{\Nis,\et\}$ and let $S$ be a rigid analytic space. 
\begin{enumerate}

\item[(1)] The $\infty$-category $\Shv_{\tau}(\RigSm/S;\Lambda)$
is compactly generated if $\tau$ is the Nisnevich topology or 
if $\Lambda$ is eventually coconnective.
A set of compact generators is given,
up to desuspension, by the $\Lambda_{\tau}(X)$ for $X\in \RigSm/S$ 
quasi-compact, quasi-separated and $(\Lambda,\tau)$-good.

\item[(2)] The $\infty$-category $\Shv^{\hyp}_{\tau}(\RigSm/S;\Lambda)$
is compactly generated if $S$ is $(\Lambda,\tau)$-admissible.
A set of compact generators is given, up to desuspension, 
by the $\Lambda_{\tau}(X)$ for $X\in \RigSm/S$ quasi-compact, quasi-separated and $(\Lambda,\tau)$-good.

\end{enumerate}
The above statements are also true with 
``$\RigSm/S$'' replaced with ``$\Et/S$'' and
``$\Etgr/S$'' when applicable (i.e., when $\tau$ is the Nisnevich 
topology).
\end{prop}

\begin{proof}
In each situation, 
we only need to show that $\Lambda_{\tau}(X)$ is a compact object
assuming that $X$ is quasi-compact and quasi-separated. 
The problem being local on $X$, we may actually 
assume that $X=\Spf(A)^{\rig}$ for an adic ring $A$
of principal ideal type.
Saying that $\Lambda_{\tau}(X)$ is compact is 
equivalent to saying that the functor 
$\mathcal{F}\mapsto \mathcal{F}(X)$ commutes with 
filtered colimits.  
This can be checked by first restricting to 
the small site of $X$. Therefore, we may replace $S$ 
by $X$ and assume that $S=\Spf(A)^{\rig}$ for an 
adic ring $A$. Moreover, it is enough to show the
versions of the above statements for $\Et/S$, when $\tau=\et$,
and for $\Etgr/S$, when $\tau=\Nis$. (Here we implicitly 
rely on Corollary \ref{cor:univers-Lambda-admis}.) 
We split the proof into two steps. (The reduction to $S=\Spf(A)^{\rig}$
is only needed in the second step.)

\paragraph*{Step 1} 
\noindent
Here we prove the second statement.
We concentrate on the \'etale topology; the case of the 
Nisnevich topology is similar. Thus, we need to show that 
$\Lambda_{\et}(X)$ is a compact object of 
$\Shv^{\hyp}_{\et}(\Et/S;\Lambda)$ when $X\in \Et/S$ is quasi-compact,
quasi-separated and $(\Lambda,\et)$-good.
This follows from combining 
Lemmas \ref{lem:pi-0-Lambda-coh-dim} 
and \ref{lem:Lamba-tau-coh-pointwise}, and using 
Remark \ref{rmk:good-virtual-coh-dim}.

\paragraph*{Step 2} 
\noindent
Here we prove the first statement. 
Let $\pi\in A$ be a generator of an ideal of definition.
We may write $A$ as the colimit of a cofiltered inductive system 
$(A_{\alpha})_{\alpha}$ where each $A_{\alpha}$ is an 
adic $\Z[[\pi]]$-algebra which is topologically of finite type.
Set $S_{\alpha}=\Spf(A_{\alpha})^{\rig}$.
Since the inclusion functor $\Prl_{\omega}\to \Prl$ 
commutes with filtered colimits by \cite[Proposition 5.5.7.6]{lurie},
it is enough by Lemma \ref{lem:sheaves-on-limit-etale}
below to show the first statement for each 
$S_{\alpha}$. Said differently, we may assume that 
$S$ is of finite type over $\Spf(\Z[[\pi]])^{\rig}$, and 
hence $(\Lambda,\tau)$-admissible. Since 
$\Lambda$ is eventually coconnective when $\tau=\et$, 
Lemma \ref{lem:auto-hypercomp-etale} implies that 
$\Shv_{\et}(\Et/S;\Lambda)$ is equivalent to 
$\Shv_{\et}^{\hyp}(\Et/S;\Lambda)$ and similarly for the 
small Nisnevich site. We may now use the first step to conclude.
\end{proof}

\begin{lemma}
\label{lem:sheaves-on-limit-etale}
Let $(\mathcal{S}_{\alpha})_{\alpha}$ be a cofiltered inverse system 
of quasi-compact and quasi-separated formal schemes 
with affine transition maps, and let $\mathcal{S}=\lim_{\alpha}
\mathcal{S}_{\alpha}$ be the limit of this system. 
We set $S_{\alpha}=\mathcal{S}^{\rig}_{\alpha}$ and 
$S=\mathcal{S}^{\rig}$. Then there is an equivalence 
\begin{equation}
\label{eq-lem:sheaves-on-limit-etale}
\underset{\alpha}{\colim}\,\Shv_{\et}(\Et/S_{\alpha};\Lambda)
\simeq 
\Shv_{\et}(\Et/S;\Lambda)
\end{equation}
in $\Prl$, where the colimit is also taken in $\Prl$.
A similar result is also true for the small Nisnevich sites.
\end{lemma}

\begin{proof}
We only discuss the \'etale case.
We have an equivalence of $\infty$-categories
\begin{equation}
\label{eq-prop:compact-sheaves-qcqs-1}
\underset{\alpha}{\colim}\,\PSh(\Et/S_{\alpha};\Lambda)
\simeq 
\PSh(\underset{\alpha}{\colim}\,
\Et/S_{\alpha};\Lambda)
\end{equation}
where the first colimit is taken in $\Prl$. 
(This is clear for $\mathcal{P}(-)$ 
instead of $\PSh(-;\Lambda)$ by the universal property 
of $\infty$-categories of presheaves, and 
we deduce the formula for $\PSh(-;\Lambda)$ using the equivalence 
$\PSh(-;\Lambda)\simeq \mathcal{P}(-)\otimes \Mod_{\Lambda}$.)
Using Remark \ref{rmk:effective-hypercover},
the fact that every cover in $\lim_{\alpha}(\Et/S_{\alpha},\et)$
is the image of a cover in $(\Et/S_{\alpha},\et)$ for some $\alpha$,
and the universal property of localisation given by
\cite[Proposition 5.5.4.20]{lurie}, 
we deduce from
\eqref{eq-prop:compact-sheaves-qcqs-1} an  
equivalence of $\infty$-categories
\begin{equation}
\label{eq-prop:compact-sheaves-qcqs-2}
\underset{\alpha}{\colim}\,\Shv_{\et}(\Et/S_{\alpha};\Lambda)
\simeq 
\Shv_{\et}(\underset{\alpha}{\colim}\,\Et/S_{\alpha};\Lambda)
\end{equation}
where the first colimit is taken in $\Prl$. 
On the other hand, by Corollary \ref{cor:limit-etale-topi-rig-an}, 
we have an equivalence of sites
$(\Et/S,\et)\simeq \lim_{\alpha} (\Et/S_{\alpha},\et)$.
Applying Lemma
\ref{lem:equi-of-sites-infty-topoi} 
we get an equivalence of $\infty$-categories
\begin{equation}
\label{eq-prop:compact-sheaves-qcqs-3}
\Shv_{\et}(\underset{\alpha}{\colim}\,\Et/S_{\alpha};\Lambda)
\simeq \Shv_{\et}(\Et/S;\Lambda).
\end{equation}
We conclude by combining 
\eqref{eq-prop:compact-sheaves-qcqs-2}
and 
\eqref{eq-prop:compact-sheaves-qcqs-3}.
\end{proof}

\begin{prop}
\label{prop:compact-shv-rigsm}
Let $\tau\in\{\Nis,\et\}$ and let $S$ be a rigid analytic space. 
\begin{enumerate}

\item[(1)] The $\infty$-category $\RigSH^{(\eff)}_{\tau}(S;\Lambda)$
is compactly generated if $\tau$ is the Nisnevich topology or 
if $\Lambda$ is eventually coconnective.
A set of compact generators is given,
up to desuspension and negative Tate twists when applicable, 
by the $\M^{(\eff)}(X)$ for $X\in \RigSm/S$ 
quasi-compact, quasi-separated and $(\Lambda,\tau)$-good. 

\item[(2)] The $\infty$-category $\RigSH^{(\eff),\,\hyp}_{\tau}(S;\Lambda)$
is compactly generated if $S$ is $(\Lambda,\tau)$-admissible.
A set of compact generators is given, up to desuspension 
and negative Tate twists when applicable, 
by the $\M^{(\eff)}(X)$ for $X\in \RigSm/S$ quasi-compact, 
quasi-separated and $(\Lambda,\tau)$-good.

\end{enumerate}
Moreover, under the stated assumptions, 
the monoidal $\infty$-category 
$\RigSH^{(\eff,\,\hyp)}_{\tau}(S;\Lambda)^{\otimes}$
belongs to $\CAlg(\Prl_{\omega})$ and, if $f:T\to S$ 
is a quasi-compact and quasi-separated 
morphism of rigid analytic spaces with $T$
assumed $(\Lambda,\tau)$-admissible in the hypercomplete case, 
the functor 
$f^*:\RigSH^{(\eff,\,\hyp)}_{\tau}(S;\Lambda)\to 
\RigSH^{(\eff,\,\hyp)}_{\tau}(T;\Lambda)$ is compact-preserving, 
i.e., belongs to $\Prl_{\omega}$.
\end{prop}

\begin{proof} 
Using Lemma \ref{lem:generation-rigsh}, we are left to show that
the objects $\M^{(\eff)}(X)$ are compact, for $X$ as in the statement.
In the effective case, this would follow from 
\cite[Corollary 5.5.7.3]{lurie} and Proposition
\ref{prop:compact-sheaves-qcqs} if we knew that 
$\RigSH^{\eff,\,(\hyp)}_{\tau}(S;\Lambda)$ is stable under filtered 
colimits in $\Shv^{(\hyp)}_{\tau}(\RigSm/S;\Lambda)$. But this is 
indeed the case by Proposition
\ref{prop:compact-sheaves-qcqs} and 
Remark \ref{rmk:on-B-1-local}.
The $\Tate$-stable case follows from the effective case using 
Remark \ref{rmk:symmetric-obj-}
and \cite[Proposition 5.5.7.6]{lurie}.
\end{proof}

\begin{rmk}
\label{rmk:compact-gen-alg-for}
A similar statement with a similar proof is  
also true for the $\infty$-category 
$\SH^{(\eff,\,\hyp)}_{\tau}(S;\Lambda)$ 
of algebraic motives over a scheme $S$, generalising 
\cite[Proposition 3.19]{ayoub-etale}. 
\end{rmk}

\subsection{Continuity, I. A preliminary result}

$\empty$

\smallskip

\label{subsect:continuity}

The goal of this subsection and the next one 
is to prove the continuity property for the functor 
$\RigSH_{\tau}^{(\eff)}(-;\Lambda)$
which, roughly speaking, asserts that this functor
transforms limits of certain cofiltered inverse systems of 
rigid analytic spaces into filtered colimits of presentable
$\infty$-categories. The precise statement is given in Theorem
\ref{thm:anstC} below.
(Note that we do not claim that $S$ is the limit of 
$(S_{\alpha})_{\alpha}$ in the categorical sense.)
Later, in Subsection 
\ref{sec:stalks}, we will generalise Theorem
\ref{thm:anstC}
to include more general inverse systems and a weaker notion
of limits; see Theorem \ref{thm:anstC-v2} below.

We let $\tau\in \{\Nis,\et\}$ be a topology
on rigid analytic spaces.

\begin{thm}
\label{thm:anstC}
\ncn{continuity}
Let $(\mathcal{S}_{\alpha})_{\alpha}$ be a cofiltered inverse system 
of quasi-compact and quasi-separated formal schemes 
with affine transition maps, and let $\mathcal{S}=\lim_{\alpha}
\mathcal{S}_{\alpha}$ be the limit of this system. 
We set $S_{\alpha}=\mathcal{S}^{\rig}_{\alpha}$ and 
$S=\mathcal{S}^{\rig}$. We assume one of the following two 
alternatives. 
\begin{enumerate}

\item[(1)] We work in the non-hypercomplete case. 

\item[(2)] We work in the hypercomplete case, and $S$
and the $S_{\alpha}$'s are $(\Lambda,\tau)$-admissible. 
When $\tau$ is the \'etale topology, we assume furthermore that
$\Lambda$ is eventually coconnective or that  
the numbers $\pvcd_{\Lambda}(S_{\alpha})$ 
are bounded independently of $\alpha$. (See
Definition \ref{dfn:Lambda-coh-dim-4}.)

\end{enumerate}
Then the obvious functor
\begin{equation}
\label{eq-thm:anstC}
\underset{\alpha}{\colim}\,
\RigSH^{(\eff,\,\hyp)}_{\tau}(S_{\alpha};\Lambda)
\to \RigSH^{(\eff,\,\hyp)}_{\tau}(S;\Lambda),
\end{equation}
where the colimit is taken in $\Prl$, is an equivalence.
\end{thm}

\begin{rmk}
\label{rmk:anstC-comp-gen}
Keep the notations and hypotheses as in Theorem
\ref{thm:anstC}.
Using \cite[Corollary 3.2.3.2]{lurie:higher-algebra}, 
we can upgrade \eqref{eq-thm:anstC}
into an equivalence
\begin{equation}
\label{eq-thm:anstC-mon}
\underset{\alpha}{\colim}\,\RigSH^{(\eff)}_{\tau}
(S_{\alpha};\Lambda)^{\otimes}
\simeq \RigSH^{(\eff)}_{\tau}(S;\Lambda)^{\otimes}
\end{equation}
in $\CAlg(\Prl)$, where the colimit is also taken in 
$\CAlg(\Prl)$.
\end{rmk}

\begin{rmk}
\label{rmk-thm:anstC-2cases}
The two alternatives considered in the statement of 
Theorem \ref{thm:anstC} have a nontrivial intersection
given as follows.
\begin{enumerate}

\item[(2$'$)] We work in the hypercomplete case and we assume that 
the $S_{\alpha}$'s and $S$ are $(\Lambda,\tau)$-admissible. 
When $\tau$ is the \'etale topology, we assume furthermore that
$\Lambda$ is eventually coconnective.

\end{enumerate}
Indeed, by Proposition
\ref{prop:automatic-hypercomp-motives}, we have in this case
$\RigSH^{(\eff),\,\hyp}_{\tau}(S_{\alpha};\Lambda)=
\RigSH^{(\eff)}_{\tau}(S_{\alpha};\Lambda)$, 
and similarly for $S$ in place 
of the $S_{\alpha}$'s. Said differently, the alternative (1)
covers the alternative (2) except when $\Lambda$ is not eventually 
coconnective, in which case we need a strong assumption
on the punctual virtual $\Lambda$-cohomological dimensions of the 
$S_{\alpha}$'s.
\end{rmk}

\begin{rmk}
\label{rmk-lem:sheaves-on-limit-etale}
Theorem \ref{thm:anstC} in the non-hypercomplete case is a 
motivic version of Lemma
\ref{lem:sheaves-on-limit-etale}. The conclusion 
of this lemma holds also in the hypercomplete case 
under the alternative (2)
as shown in corollary
\ref{cor:colim-small-et-site-2} 
below.
\end{rmk}

The proof of Theorem \ref{thm:anstC}
spans the entire subsection and the next one.
In fact, we will obtain this theorem 
as a combination of two other
results, namely Propositions 
\ref{prop:easy-equiv-S-alpha} and
\ref{prop:approx-rig-sh}, 
which are both interesting in their own right.
The proof of Proposition 
\ref{prop:approx-rig-sh}
will be given in Subsection 
\ref{subsect:continuity-2}.

\begin{nota}
\label{nota:motive-over-S-alpha}
Let $(S_{\alpha})_{\alpha}$ be a cofiltered inverse system of 
quasi-compact and quasi-separated rigid analytic spaces.
We define the $\infty$-category 
$\RigSH_{\tau}^{(\eff,\,\hyp)}((S_{\alpha})_{\alpha};\Lambda)$,
of rigid analytic motives over the rigid analytic pro-space 
$(S_{\alpha})_{\alpha}$, in the usual way from the {
limit site $\lim_{\alpha}(\RigSm/S_\alpha,\tau)$, that is, from the ordinary} category 
$$\RigSm/(S_{\alpha})_{\alpha}=
\underset{\alpha}{\colim}\, \RigSm/S_{\alpha}$$
endowed with the limit topology $\tau$. More precisely, one
repeats Definitions \ref{def:DAeff} and \ref{dfn:rigsh-stable}
with ``$\RigSm/S$'' replaced with ``$\RigSm/(S_{\alpha})_{\alpha}$''.
We denote also by 
$$\M^{(\eff)}:\RigSm/(S_{\alpha})_{\alpha}
\to \RigSH_{\tau}^{(\eff,\,\hyp)}((S_{\alpha})_{\alpha};\Lambda)$$
\symn{$\M^{(\eff)}$}
the obvious functor.
\symn{$\RigSH^{(\eff,\hyp)}$}
\symn{$\RigSm$}
\end{nota}

\begin{rmk}
\label{rmk:object-in-colim-RigSm}
\ncn{rigid analytic spaces!pro-objects}
Let $\Pro(\RigSpc)$ be the category of rigid analytic 
pro-spaces and consider the overcategory 
$\Pro(\RigSpc)/(S_{\alpha})_{\alpha}$
of $(S_{\alpha})_{\alpha}$-objects. 
There is a fully faithful embedding 
$$\RigSm/(S_{\alpha})_{\alpha}
\to \Pro(\RigSpc)/(S_{\alpha})_{\alpha}$$
and we will identify 
$\RigSm/(S_{\alpha})_{\alpha}$ with its essential image 
by this functor. Thus, we may think of 
an object of $\RigSm/(S_{\alpha})_{\alpha}$
as a pro-object $(X_{\alpha})_{\alpha\leq \alpha_0}$,
where $X_{\alpha_0}$ is a smooth rigid analytic 
$S_{\alpha_0}$-space and, for $\alpha\leq \alpha_0$, 
$X_{\alpha}\simeq X_{\alpha_0}\times_{S_{\alpha_0}}S_{\alpha}$.
If $(S_{\alpha})_{\alpha}$ is as in Theorem
\ref{thm:anstC},
given such a pro-object $(X_{\alpha})_{\alpha\leq \alpha_0}$, 
we denote by $X$ the rigid analytic $S$-space defined as follows. 
Assume first that there is a formal model 
$\mathcal{X}_{\alpha_0}$ of $X_{\alpha_0}$ over $\mathcal{S}_{\alpha_0}$.
Let $(\mathcal{X}_{\alpha})_{\alpha\leq \alpha_0}$ be the 
formal pro-scheme given by $\mathcal{X}_{\alpha}=\mathcal{X}_{\alpha_0}
\times_{\mathcal{S}_{\alpha_0}}\mathcal{S}_{\alpha}$. We set 
$X=\mathcal{X}^{\rig}$ where $\mathcal{X}=\lim_{\alpha\leq \alpha_0}\mathcal{X}_{\alpha}$.
This is independent of the choice of 
$\mathcal{X}_{\alpha_0}$ and the formation of $X$
is compatible with gluing rigid analytic 
$S_{\alpha_0}$-spaces along open immersions. 
Thus, the construction of $X$ can be extended to the general case 
where we do not assume the existence of a 
formal model for $X_{\alpha_0}$.
\end{rmk}

\begin{lemma}
\label{lem:compact-gen-rigsh-pro-obj}
Let $(S_{\alpha})_{\alpha}$ and $S$ be as in Theorem
\ref{thm:anstC} and assume that 
$S$ is $(\Lambda,\tau)$-admissible.
Then, the $\infty$-category 
$\Shv^{\hyp}_{\tau}(\RigSm/(S_{\alpha})_{\alpha};\Lambda)$
is compactly generated, up to desuspension, 
by the $\Lambda_{\tau}((X_{\alpha})_{\alpha\leq \alpha_0})$
with $X_{\alpha_0}$ quasi-compact, quasi-separated and 
$(\Lambda,\tau)$-good.
\end{lemma}

\begin{proof}
This can be shown by adapting the proof of Proposition
\ref{prop:compact-sheaves-qcqs}(2). The key point is to show that
$\lim_{\alpha\leq \alpha_0}(\Et/X_{\alpha},\tau)$
has finite local and global $\Lambda$-cohomological dimensions.
By Corollary \ref{cor:limit-etale-topi-rig-an}, this limit site
is equivalent 
to $(\Et/X,\tau)$.
Thus, we may use Lemma
\ref{lem:Lamba-tau-coh-pointwise} 
to conclude.
\end{proof}

\begin{prop}
\label{prop:easy-equiv-S-alpha}
Let $(S_{\alpha})_{\alpha}$ and $S$ be as in Theorem
\ref{thm:anstC} and assume one of the alternatives (1) or (2) of that theorem. Then the obvious functor
\begin{equation}
\label{eq-lem:easy-equiv-S-alpha}
\underset{\alpha}{\colim}\,\RigSH_{\tau}^{(\eff,\,\hyp)}
(S_{\alpha};\Lambda)
\to \RigSH_{\tau}^{(\eff,\,\hyp)}((S_{\alpha})_{\alpha};\Lambda),
\end{equation}
where the colimit is taken in $\Prl$, is an equivalence.
\end{prop}

\begin{proof}
We first work under the alternative (1), i.e., in the 
non-hypercomplete case. Here, the result is quite
straightforward. Arguing as in the proof of Lemma 
\ref{lem:sheaves-on-limit-etale}, we get an equivalence 
of $\infty$-categories
\begin{equation}
\label{eq-lem:easy-equiv-S-alpha-2}
\underset{\alpha}{\colim}\,
\Shv_{\tau}(\RigSm/S_{\alpha};\Lambda)
\simeq  \Shv_{\tau}(\RigSm/(S_{\alpha})_{\alpha};\Lambda),
\end{equation}
where the colimit is taken in $\Prl$.
Using the universal property of localisation given by
\cite[Proposition 5.5.4.20]{lurie}, we deduce from 
\eqref{eq-lem:easy-equiv-S-alpha-2} that 
\eqref{eq-lem:easy-equiv-S-alpha}
is an equivalence in the effective case.
We then deduce the $\Tate$-stable
case using Remark \ref{rmk:symmetric-obj-} 
and commutation of colimits with colimits.

Next, we work under the alternative (2). 
Arguing as before, we see that it is enough to 
prove the hypercomplete analogue of the equivalence
\eqref{eq-lem:easy-equiv-S-alpha-2}, i.e., it is enough to show that
\begin{equation}
\label{eq-lem:easy-equiv-S-alpha-3}
\underset{\alpha}{\colim}\,
\Shv^{\hyp}_{\tau}(\RigSm/S_{\alpha};\Lambda)
\to  \Shv^{\hyp}_{\tau}(\RigSm/(S_{\alpha})_{\alpha};\Lambda),
\end{equation}
is an equivalence.
It follows from Lemma \ref{lem:compact-gen-rigsh-pro-obj} 
that the functor \eqref{eq-lem:easy-equiv-S-alpha-3}
belongs to $\Prl_{\omega}$ and that it
takes a set of compact generators to a set of compact generators.
Thus, it remains to show that this functor is fully faithful 
on compact objects. Explicitly, we need 
to show the following assertion. Given two compact 
objects $\mathcal{M}$ and $\mathcal{N}$ in 
$\Shv_{\tau}^{\hyp}(\RigSm/S_{\alpha_0};\Lambda)$,
for some index $\alpha_0$, the natural map
\begin{equation}
\label{eq-lem:easy-equiv-S-alpha-4}
\underset{\alpha\leq \alpha_0}{\colim}\,
\Map(f_{\alpha\leq \alpha_0}^*\mathcal{M},
f_{\alpha \leq \alpha_0}^*\mathcal{N})
\to \Map(f_{\alpha_0}^*\mathcal{M},f_{\alpha_0}^*\mathcal{N})
\end{equation}
is an equivalence. Here $f_{\alpha\leq \alpha_0}:S_{\alpha}
\to S_{\alpha_0}$
and $f_{\alpha_0}:(S_{\alpha})_{\alpha}\to S_{\alpha_0}$ 
are the obvious morphisms.

{

Let $I$ be the indexing category of the inverse system
$(S_{\alpha})_{\alpha}$. 
We denote by
$\widetilde{S}:I \to \RigSpc$
the diagram of rigid analytic spaces defining the pro-object 
$(S_{\alpha})_{\alpha}$, i.e., 
sending $\alpha$ to $S_{\alpha}$. 
We define the site $(\RigSm/\widetilde{S},\tau)$ in the usual 
way, i.e., by adapting the beginning of \cite[\S 4.5.1]{ayoub-th2}.
We have a premorphism of sites (in the sense of 
\cite[D\'efinition 4.4.46]{ayoub-th2})
\begin{equation}
\label{eq-lem:easy-equiv-S-alpha-317}
\rho:(\RigSm/(S_{\alpha})_{\alpha},\tau) 
\to (\RigSm/\widetilde{S},\tau)
\end{equation}
induced by the functor $\RigSm/\widetilde{S}\to 
\RigSm/(S_{\alpha})_{\alpha}$ given by
$(\beta,X)\mapsto (X\times_{S_{\beta}}S_{\alpha})_{\alpha\leq \beta}$.
The inverse image functor $\rho^*$ is given, informally, by
$\rho^*(\mathcal{K})=\colim_{\beta}\,((S_{\alpha})_{\alpha} 
\to S_{\beta})^*\mathcal{K}_{\beta}$, 
where $\mathcal{K}_{\beta}$ is the restriction of 
$\mathcal{K}$ to $\RigSm/S_{\beta}$.
The inclusion $\RigSm/S_{\alpha_0}\subset \RigSm/\widetilde{S}$
induces a functor 
\begin{equation}
\label{eq-lem:easy-equiv-S-alpha-489}
\Shv^{(\hyp)}_{\tau}(\RigSm/S_{\alpha_0};\Lambda)
\to \Shv^{(\hyp)}_{\tau}(\RigSm/\widetilde{S};\Lambda).
\end{equation}
We may assume that $\mathcal{M}=\Lambda_{\tau}(X_{\alpha_0})$ with $X_{\alpha_0} \in \RigSm/S_{\alpha_0}$
quasi-compact, quasi-separated and $(\Lambda,\tau)$-good.
We let $\mathcal{R}$ be the image of $\mathcal{N}$ by the functor
\eqref{eq-lem:easy-equiv-S-alpha-489}.
Arguing as in the proof of \cite[Proposition 3.20]{ayoub-etale},
the assertion that \eqref{eq-lem:easy-equiv-S-alpha-4} is an equivalence
would follow if we show that the functor 
$$\rho^*:\PSh(\RigSm^{\qcqs}/\widetilde{S};\Lambda)
\to \PSh(\RigSm^{\qcqs}/(S_{\alpha})_{\alpha};\Lambda)$$
takes $\mathcal{R}$ to a presheaf $\rho^*(\mathcal{R})$
whose restriction to $\Et^{\qcqs}/(X_{\alpha\leq \alpha_0})_{\alpha}$ 
is a $\tau$-hypersheaf. This follows from Lemma 
\ref{lem:coh-limi-site-sga} below. (Compare with
\cite[Lemme 3.21]{ayoub-etale}.)}
\end{proof}

\begin{lemma}
\label{lem:coh-limi-site-sga}
Let $\widetilde{\mathcal{X}}:I \to \FSch$ be a diagram of quasi-compact
and quasi-separated formal schemes, with $I$ a cofiltered category,
and with affine transition morphisms.
Let $(\mathcal{X}_{\alpha})_{\alpha}$ be the associated pro-object
and $\mathcal{X}$ its limit. 
Set $\widetilde{X}=\widetilde{\mathcal{X}}{}^{\rig}$, 
$X_{\alpha}=\mathcal{X}^{\rig}_{\alpha}$ and $X=\mathcal{X}^{\rig}$.
Assume that the alternative (2) in Theorem 
\ref{thm:anstC} is satisfied with ``$(X_{\alpha})_{\alpha}$'' and 
``$X$'' instead of ``$(S_{\alpha})_{\alpha}$'' and ``$S$''.
{
Assume also that the $X_{\alpha}$'s are $(\Lambda,\tau)$-good.} 
Then the functor 
$$\rho^*:\PSh(\Et^{\qcqs}/\widetilde{X};\Lambda)
\to \PSh(\Et^{\qcqs}/(X_{\alpha})_{\alpha};\Lambda)$$ 
takes $\tau$-hypersheaves to $\tau$-hypersheaves. 
\end{lemma}

\begin{proof}
We split the proof into three steps. Below $\mathcal{K}$ is a 
$\tau$-hypersheaf of $\Lambda$-modules on $\Et^{\qcqs}/\widetilde{X}$.

\paragraph*{Step 1}
\noindent
We first deal with the case where $\Lambda$ is eventually 
coconnective. The proof in this case is similar to that of 
\cite[Lemme 3.21]{ayoub-etale}. First, one considers the case 
where $\mathcal{K}$ is discrete, i.e., is the Eilenberg--Mac 
Lane spectrum associated to an ordinary sheaf of 
$\pi_0\Lambda$-modules. This case follows from
\cite[Expos\'e VII, Th\'eor\`eme 5.7]{SGAIV2}.
By induction, one can then treat the case where $\mathcal{K}$ is 
bounded (i.e., where the discrete sheaves $\pi_i(\mathcal{K})$ vanish 
for $|i|$ big enough).
Finally, we deduce the general case from the bounded case
as follows.
A general $\mathcal{K}$ can be written as a colimit of 
objects of the form $\Lambda_{\tau}(\alpha_0,U)$, for 
$U\in \Et^{\qcqs}/X_{\alpha_0}$.
Since $\Lambda$ is eventually coconnective,
$\Lambda_{\tau}(\alpha_0,U)$ is bounded. The result for $\mathcal{K}$ follows then 
from the bounded case and Lemma
\ref{lem:pi-0-Lambda-coh-dim}(3) 
which implies that colimits in 
$\PSh(\Et^{\qcqs}/(X_{\alpha})_{\alpha};\Lambda)$
preserve $\tau$-hypersheaves.
(Here, we use that the site 
$(\Et^{\qcqs}/(X_{\alpha})_{\alpha},\tau)$
has finite local $\Lambda$-cohomological dimension as 
explained in the proof of 
Lemma \ref{lem:compact-gen-rigsh-pro-obj}.)

\paragraph*{Step 2}
\noindent
We next consider the case of the Nisnevich topology. 
The site $(\Et^{\qcqs}/(X_{\alpha})_{\alpha},\Nis)$ is equivalent
to $(\Et^{\qcqs}/X,\Nis)$.
Thus, by Lemma \ref{lem:auto-hypercomp-etale}(1), 
every Nisnevich sheaf on $\Et^{\qcqs}/(X_{\alpha})_{\alpha}$ 
is a Nisnevich hypersheaf.
Thus, to check that $\rho^*\mathcal{K}$ is a Nisnevich hypersheaf, it 
is enough to prove that $\rho^*\mathcal{K}$ has the Mayer--Vietoris 
property for the image in $\Et^{\qcqs}/(X_{\alpha})_{\alpha}$
of a Nisnevich square in $\Et^{\qcqs}/X_{\alpha}$,
for some $\alpha$. This is easily checked using 
exactness of filtered colimits on $\Mod_{\Lambda}$ and 
the formula $\rho^*\mathcal{K}=\colim_{\beta}\,((X_{\alpha})_{\alpha}
\to X_{\beta})^*\mathcal{K}_{\beta}$. The details are 
left to the reader.

\paragraph*{Step 3}
\noindent
We now treat the case of the \'etale topology assuming that
the numbers $\pvcd_{\Lambda}(X_{\alpha})$ 
are bounded independently of $\alpha$. In fact, since the 
$X_{\alpha}$'s are $(\Lambda,\et)$-good, there is a common bound $e$
for the $\Lambda$-cohomological dimensions of the residue fields of 
all the $X_{\alpha}$'s.

Denote by $\pi$ the morphism of sites of the form
$(\Et/(-),\et)\to (\Et/(-),\Nis)$ and by $\pi_*$ the induced
functor on $\infty$-categories of hypersheaves of $\Lambda$-modules.
Also, denote by
$$\rho^*_{\Nis}:\Shv^{\hyp}_{\Nis}(\Et^{\qcqs}/\widetilde{X};\Lambda)
\to \Shv^{\hyp}_{\Nis}(\Et^{\qcqs}/(X_{\alpha})_{\alpha};\Lambda)$$
the inverse image functor on Nisnevich hypersheaves. 
By the second step, $\rho^*_{\Nis}$ coincides
with $\rho^*$ on Nisnevich hypersheaves of $\Lambda$-modules.

By Lemma 
\ref{lem:pi-0-Lambda-coh-dim}(3),
the property that $\rho^*\mathcal{K}$ is an \'etale hypersheaf is stable 
by colimits in $\mathcal{K}$. Since $\mathcal{K}\simeq \colim_{n} \tau_{\geq -n}\mathcal{K}$,
we may assume that $\mathcal{K}$ is bounded from above, and even connective. By Lemma
\ref{lem:pi-0-Lambda-coh-dim}(1), we have an equivalence
$\mathcal{K}\simeq \lim_n\tau_{\leq n}\mathcal{K}$
yielding an equivalence 
$\pi_*\mathcal{K}\simeq \lim_n \pi_*\tau_{\leq n}\mathcal{K}$.
{
The proof of Lemma 
\ref{lem:Lamba-tau-coh-pointwise}
shows that $\pi_*\tau_{\leq n+1}\mathcal{K}
\to \pi_*\tau_{\leq n}\mathcal{K}$
induces an isomorphism on homotopy Nisnevich sheaves in degrees $\leq n-e$,
and the same is true for 
$\rho^*_{\Nis}\pi_*\tau_{\leq n+1}\mathcal{K}
\to \rho^*_{\Nis}\pi_*\tau_{\leq n}\mathcal{K}$.
Since $X$ and $X_{\alpha}$'s have finite Krull dimensions,
we deduce that the morphisms
$$\lim_n \pi_*\tau_{\leq n} \mathcal{K}
\to \pi_*\tau_{\leq m} \mathcal{K}
\quad \text{and} \qquad 
\lim_n \rho_{\Nis}\pi_*\tau_{\leq n} \mathcal{K}
\to \rho_{\Nis}\pi_*\tau_{\leq m} \mathcal{K}$$
induce isomorphisms on homotopy Nisnevich sheaves in degrees $\leq m-e$, 
for any integer $m$. It follows that the natural map
\begin{equation}
\label{eq-lem:coh-limi-site-sga-3}
\rho^*_{\Nis}\pi_*\mathcal{K} =
\rho^*_{\Nis}\lim_n\pi_*\tau_{\leq n}\mathcal{K}
\to \lim_n\rho_{\Nis}^*\pi_*\tau_{\leq n}\mathcal{K}
\end{equation}
induces isomorphisms on homotopy Nisnevich sheaves. Since 
both sides are Nisnevich hypersheaves, we deduce that 
\eqref{eq-lem:coh-limi-site-sga-3}
is an equivalence. Thus, we are left 
to show that $\rho^*_{\Nis}\pi_*\tau_{\leq n}\mathcal{K}$
is an \'etale hypersheaf for every $n$. 
This follows from the first step since 
$\tau_{\leq n}\mathcal{K}$ is naturally an \'etale hypersheaf of 
$\tau_{\leq n}\Lambda$-modules.}
\end{proof}

Lemma \ref{lem:coh-limi-site-sga}
has the following consequence which we state for 
completeness. 

\begin{cor}
\label{cor:colim-small-et-site-2}
Let $(S_{\alpha})_{\alpha}$ and $S$ be as in Theorem
\ref{thm:anstC} and assume one of the alternatives (1) or (2) of that theorem. Then the obvious functor
\begin{equation}
\label{eq-cor:colim-small-et-site-2}
\underset{\alpha}{\colim}\,\Shv_{\tau}^{(\hyp)}
(\Et/S_{\alpha};\Lambda)
\to \Shv_{\tau}^{(\hyp)}(\Et/S;\Lambda),
\end{equation}
where the colimit is taken in $\Prl$, is an equivalence.
\end{cor}

\begin{proof}
The non-hypercomplete case is already stated in 
Lemma \ref{lem:sheaves-on-limit-etale}. The hypercomplete case 
follows from Lemma \ref{lem:coh-limi-site-sga} by
arguing as in the proof of \cite[Proposition 3.20]{ayoub-etale}.
\end{proof}

The proof of Proposition
\ref{prop:easy-equiv-S-alpha},
adapted to the algebraic setting gives the following 
generalisation of \cite[Proposition 3.20]{ayoub-etale} and
\cite[Proposition C.12(4)]{hoyois-quadratic}.

\begin{prop}
\label{prop:cont-algebraic-29}
Let $(S_{\alpha})_{\alpha}$ be a cofiltered inverse system of 
quasi-compact and quasi-separated 
schemes with affine transition maps, 
and let $S=\lim_{\alpha}S_{\alpha}$
be the limit of this system. 
We assume one of the following two 
alternatives. 
\begin{enumerate}

\item[(1)] We work in the non-hypercomplete case. 

\item[(2)] We work in the hypercomplete case, and $S$ and 
the $S_{\alpha}$'s are $(\Lambda,\tau)$-admissible. 
When $\tau$ is the \'etale topology, we assume furthermore that
$\Lambda$ is eventually coconnective or that  
the numbers $\pvcd_{\Lambda}(S_{\alpha})$ 
are bounded independently of $\alpha$.

\end{enumerate}
Then the obvious functor
$$\underset{\alpha}{\colim}\,\SH_{\tau}^{(\eff,\,\hyp)}
(S_{\alpha};\Lambda)
\to \SH^{(\eff,\,\hyp)}_{\tau}(S;\Lambda),$$
where the colimit is taken in $\Prl$, is an equivalence.
\end{prop}

\begin{proof}
Indeed, in the algebraic setting, $\Sm^{\qcqs}/S$ is equivalent 
to $\colim_{\alpha}\,
\Sm^{\qcqs}/S_{\alpha}$.
\end{proof}

Theorem \ref{thm:anstC}
follows by combining Proposition 
\ref{prop:easy-equiv-S-alpha} 
and the next result.

\begin{prop}
\label{prop:approx-rig-sh}
Let $(\mathcal{S}_{\alpha})_{\alpha}$ be a cofiltered inverse system 
of quasi-compact and quasi-separated formal schemes 
with affine transition maps, and let $\mathcal{S}=\lim_{\alpha}
\mathcal{S}_{\alpha}$ be the limit of this system. 
We set $S_{\alpha}=\mathcal{S}^{\rig}_{\alpha}$ and 
$S=\mathcal{S}^{\rig}$.
Then the obvious functor
\begin{equation}
\label{eqn:funct-S-alpha-alpha-to-S}
\RigSH_{\tau}^{(\eff,\,\hyp)}((S_{\alpha})_{\alpha};\Lambda)
\to \RigSH_{\tau}^{(\eff,\,\hyp)}(S;\Lambda)
\end{equation}
is an equivalence.
\end{prop}

The proof of Proposition
\ref{prop:approx-rig-sh} is given in 
the next subsection.

\subsection{Continuity, II. Approximation up to homotopy}

$\empty$

\smallskip

\label{subsect:continuity-2}

The goal of this section is to prove Proposition
\ref{prop:approx-rig-sh}. The proof is similar to that of 
\cite[Proposition 4.5]{vezz-fw}, but some new ingredients
are necessary to deal with the generality considered in this paper.
We start with some reductions.

\begin{lemma}
\label{lem:reduc-effect-nis-non-hyp}
It is enough to prove Proposition
\ref{prop:approx-rig-sh} in the effective, non-hypercomplete case 
and for $\tau$ the Nisnevich topology. 
Said differently, it is enough to show that  
the obvious functor
\begin{equation}
\label{eq-lem:reduc-effect-nis-non-hyp}
\RigSH_{\Nis}^{\eff}((S_{\alpha})_{\alpha};\Lambda)
\to \RigSH_{\Nis}^{\eff}(S;\Lambda)
\end{equation}
is an equivalence.
\end{lemma}

\begin{proof}
The $\Tate$-stable
case follows from the effective case 
using Remark \ref{rmk:symmetric-obj-} 
and commutation of colimits with colimits.
Assume that 
\eqref{eq-lem:reduc-effect-nis-non-hyp}
is an equivalence, and let's show that 
\begin{equation}
\label{eq-lem:reduc-effect-nis-non-hyp-2}
\RigSH_{\tau}^{\eff,\,(\hyp)}((S_{\alpha})_{\alpha};\Lambda)
\to \RigSH_{\tau}^{\eff,\,(\hyp)}(S;\Lambda)
\end{equation}
is also an equivalence for $\tau\in \{\Nis,\et\}$.
There are three cases to consider: 
\begin{enumerate}

\item[(1)] the Nisnevich topology in the hypercomplete case;

\item[(2)] the \'etale topology in the non-hypercomplete case;

\item[(3)] the \'etale topology in the hypercomplete case.

\end{enumerate}
In each case, we will prove that the source and the target of 
\eqref{eq-lem:reduc-effect-nis-non-hyp-2}
are obtained from the source and the target of 
\eqref{eq-lem:reduc-effect-nis-non-hyp}
by localisation with respect to a set of morphisms
and its image by the equivalence
\eqref{eq-lem:reduc-effect-nis-non-hyp}, which suffices to conclude.
These sets consist respectively, up to desuspension,
of maps of the form
$\colim_{[n]\in \mathbf{\Delta}}\, 
\M^{\eff}((U_{n,\,\alpha})_{\alpha \leq \alpha_n}) 
\to \M^{\eff}((U_{-1,\,\alpha})_{\alpha \leq \alpha_{-1}})$ where 
$(U_{\bullet,\,\alpha})_{\alpha \leq \alpha_{\bullet}}$ is:
\begin{enumerate}

\item[(1)]  a hypercover in the limit site 
$\lim_{\alpha\leq \alpha_{-1}}(\Etgr/U_{-1,\,\alpha},\Nis)$;

\item[(2)] a {\v C}ech nerve associated to a cover in the limit site 
$\lim_{\alpha\leq \alpha_{-1}}(\Et/U_{-1,\,\alpha},\et)$;

\item[(3)] a hypercover in the limit site 
$\lim_{\alpha\leq \alpha_{-1}}(\Et/U_{-1,\,\alpha},\et)$.

\end{enumerate}
Localising the source of 
\eqref{eq-lem:reduc-effect-nis-non-hyp}
by one of these sets yield the source of 
\eqref{eq-lem:reduc-effect-nis-non-hyp-2}
by construction. We now show that localising 
the target of 
\eqref{eq-lem:reduc-effect-nis-non-hyp}
by the image of one of these sets yield the target of 
\eqref{eq-lem:reduc-effect-nis-non-hyp-2}.
This relies on the following two facts.
\begin{enumerate}

\item[(a)] Given an object $(Y_{\alpha})_{\alpha\leq \beta}$ 
in $\RigSm/(S_{\alpha})_{\alpha}$ and defining $Y$ as in 
Remark \ref{rmk:object-in-colim-RigSm}, we have an equivalence of sites
$$(\Et/Y,\tau) \simeq \lim_{\alpha\leq \beta}
(\Et/Y_{\alpha},\tau)$$
and similarly for ``$\Etgr$'' instead of ``$\Et$'' when applicable.

\item[(b)] Every $X\in \RigSm/S$ is locally for the analytic 
topology in the essential image $\RigSm'/S$ of the functor 
$\RigSm/(S_{\alpha})_{\alpha} \to \RigSm/S$.
In particular, we have an equivalence of sites
$$(\RigSm/S,\tau)\simeq (\RigSm'/S,\tau)$$
which is subject to Lemma
\ref{lem:equi-of-sites-infty-topoi}.
Thus, the $\infty$-category
$\RigSH^{\eff,\,(\hyp)}_{\tau}(S;\Lambda)$
can be defined using the site 
$(\RigSm'/S,\tau)$.

\end{enumerate}
Property (a) follows from Corollary 
\ref{cor:limit-etale-topi-rig-an} and
Remark \ref{rmk:analog-Nis-limit-etale-topi-rig-an}. 
To prove (b), we may assume that the inverse system
$(\mathcal{S}_{\alpha})_{\alpha}$ is affine, induced 
by an inductive system of adic rings $(A_{\alpha})_{\alpha}$
with colimit $A$, and that 
$X=\Spf(B)^{\rig}$ with $B$ a rig-\'etale adic 
$A\langle t\rangle$-algebra with $t=(t_1,\ldots, t_n)$
a system of coordinates. Then the result follows from 
Corollary \ref{cor:proj-limi-cal-E-A}.
\end{proof}

\begin{lemma}
\label{lem:reduction-anstC}
It is enough to prove that
\eqref{eq-lem:reduc-effect-nis-non-hyp}
is an equivalence assuming that the 
formal schemes $\mathcal{S}_{\alpha}$ are affine
of principal ideal type.
\end{lemma}

\begin{proof}
Without loss of generality, we may assume that there is a
final object $o$ in the indexing category of the inverse system
$(\mathcal{S}_{\alpha})_{\alpha}$. Replacing $\mathcal{S}_o$
by the blowup of an ideal of definition, and the 
$\mathcal{S}_{\alpha}$'s by their strict transforms, we may assume
that the $\mathcal{S}_{\alpha}$'s 
are locally of principal ideal type for every $\alpha$.
The presheaf $\RigSH_{\Nis}^{\eff}(-;\Lambda)$ 
has descent for the analytic topology by 
Theorem \ref{thm:hyperdesc}. Combining this with 
Proposition \ref{prop:easy-equiv-S-alpha} and 
\cite[Proposition 4.7.4.19]{lurie:higher-algebra}, 
we see that the problem is local on $\mathcal{S}_o$, which 
finishes the proof.
(Note that the condition for applying
\cite[Proposition 4.7.4.19]{lurie:higher-algebra}
is indeed satisfied by the base
change theorem for open immersions, a special case of 
the base change theorem for smooth morphisms; see Proposition
\ref{prop:6f1}.)
\end{proof}

We now introduce a notation that we keep using 
until the end of the proof of Proposition 
\ref{prop:approx-rig-sh}.

\begin{nota}
\label{not:closed-prospace}
Let $(\mathcal{S}_{\alpha})_{\alpha}$ be a cofiltered inverse 
system of affine formal schemes, and let $\mathcal{S}=\lim_{\alpha}\mathcal{S}_{\alpha}$. Denote by $\mathcal{S}'_{\alpha}$ 
the smallest closed 
formal subscheme of $\mathcal{S}_{\alpha}$ 
containing the image of $\mathcal{S} \to \mathcal{S}_{\alpha}$. 
(Said differently, $\mathcal{O}(\mathcal{S}'_{\alpha})$ 
is the quotient of $\mathcal{O}(\mathcal{S}_{\alpha})$ by the kernel of 
$\mathcal{O}(\mathcal{S}_{\alpha})\to \mathcal{O}(\mathcal{S})$
which is a closed ideal.) Then, we have a cofiltered inverse system 
of affine formal schemes $(\mathcal{S}'_{\alpha})_{\alpha}$
and a morphism $(\mathcal{S}'_{\alpha})_{\alpha}\to 
(\mathcal{S}_{\alpha})_{\alpha}$ of inverse systems
given by closed immersions and inducing an isomorphism 
$\lim_{\alpha}\mathcal{S}'_{\alpha}
\simeq \lim_{\alpha}\mathcal{S}_{\alpha}$
on the limit. 
\end{nota}

Although, in general, the pro-objects 
$(\mathcal{S}'_{\alpha})_{\alpha}$ and 
$(\mathcal{S}_{\alpha})_{\alpha}$ are not isomorphic,
we have the following.

\begin{lemma}
\label{lem:localisation-prospace}
Let $(\mathcal{S}_{\alpha})_{\alpha}$ be a cofiltered inverse 
system of affine formal schemes. 
Let $S_{\alpha}$ and $S'_{\alpha}$ be the rigid analytic 
spaces associated to $\mathcal{S}_{\alpha}$ and $\mathcal{S}'_{\alpha}$.
Then, the obvious functor 
\begin{equation}
\label{eq-lem:localisation-prospace}
\RigSH^{\eff}_{\Nis}((S_{\alpha})_{\alpha};\Lambda)
\to \RigSH^{\eff}_{\Nis}((S'_{\alpha})_{\alpha};\Lambda)
\end{equation}
is an equivalence.
\end{lemma}

\begin{proof}
It will be more convenient 
to use Proposition
\ref{prop:easy-equiv-S-alpha}
and prove that 
\begin{equation}
\label{eqn:colim-alpha-S-vs-Z}
\underset{\alpha}{\colim}\,\RigSH^{\eff}_{\Nis}(S_{\alpha};\Lambda)
\to \underset{\alpha}{\colim}\,
\RigSH^{\eff}_{\Nis}(S'_{\alpha};\Lambda)
\end{equation}
is an equivalence in $\Prl$.
We set $U_{\alpha}=S_{\alpha}\smallsetminus S'_{\alpha}$ and
denote by $j_{\alpha}:U_{\alpha}\to S_{\alpha}$ the obvious inclusion.
For each $\alpha$, 
$\RigSH_{\Nis}^{\eff}(S_{\alpha};\Lambda)\to 
\RigSH_{\Nis}^{\eff}(S'_{\alpha};\Lambda)$
is a localisation functor with respect to the 
class of maps of the form $0\to j_{\alpha,\,\sharp}M$ 
where $M\in \RigSH_{\Nis}^{\eff}(U_{\alpha};\Lambda)$.
This follows from the localisation theorem for 
rigid analytic motives; see Proposition \ref{prop:loc1}.
Moreover, by Lemma \ref{lem:generation-rigsh},
we may assume that $M$ is, up to desuspension, 
of the form $\M^{\eff}(X)$ with 
$X\in \RigSm/U_{\alpha}$ quasi-compact and quasi-separated.

It follows from the universal property of localisation 
(given by \cite[Proposition 5.5.4.20]{lurie})
that \eqref{eqn:colim-alpha-S-vs-Z}
is also a localisation functor with respect to 
the images of the maps $0\to j_{\alpha,\,\sharp}M$, with $M$ as above. 
Thus, it is enough to show that, for 
$X\in \RigSm/U_{\alpha}$ quasi-compact and quasi-separated,
there exists $\beta\leq \alpha$ such that 
$X\times_{S_{\alpha}}S_{\beta}=\emptyset$.
This follows from the fact that 
$X$ lies over a quasi-compact open subset 
$V\subset U_{\alpha}$ and that, for $\beta\leq \alpha$ small 
enough, we have $S_{\beta}\times_{S_{\alpha}}V=\emptyset$
by, for example, \cite[Chapter 0, Proposition 2.2.10]{fujiwara-kato}.
\end{proof}

\begin{nota}
\label{not:pro-rigspc-af-prime}
Let \sym{$\FSch_{\af,\,\pr}$} be the category of affine 
formal schemes of principal ideal type, 
and $\Pro(\FSch_{\af,\,\pr})$ the category of 
pro-objects in $\FSch_{\af,\,\pr}$. 
We have an idempotent endofunctor of
$\Pro(\FSch_{\af,\,\pr})$
given by $(\mathcal{S}_{\alpha})_{\alpha}\mapsto 
(\mathcal{S}_{\alpha}')_{\alpha}$. 
We define a new category \sym{$\Pro'(\FSch_{\af,\,\pr})$}, having the same 
objects as $\Pro(\FSch_{\af,\,\pr})$ and where morphisms are given by 
$$\begin{array}{rcl}
\Hom_{\Pro'(\FSch_{\af,\,\pr})}((\mathcal{T}_{\beta})_{\beta},(\mathcal{S}_{\alpha})_{\alpha}) & = & 
\Hom_{\Pro(\FSch_{\af,\,\pr})}((\mathcal{T}'_{\beta})_{\beta},(\mathcal{S}_{\alpha}')_{\alpha})\\
& \simeq & \Hom_{\Pro(\FSch_{\af,\,\pr})}((\mathcal{T}'_{\beta})_{\beta},(\mathcal{S}_{\alpha})_{\alpha}).
\end{array}$$
The obvious functor $\Pro(\FSch_{\af,\,\pr})^{\op}
\to \Pro'(\FSch_{\af,\,\pr})^{\op}$, given by the identity on objects, 
is a localisation functor and its right adjoint is given 
on objects by $(\mathcal{S}_{\alpha})_{\alpha}\mapsto 
(\mathcal{S}_{\alpha}')_{\alpha}$.
\end{nota}

\begin{cor}
\label{cor:ext-rigsh-pro-prime}
The functor 
$$\RigSH^{\eff}_{\Nis}((-)^{\rig};\Lambda):
\Pro(\FSch_{\af,\,\pr})^{\op}\to \Prl$$
extends uniquely to $\Pro'(\FSch_{\af,\,\pr})^{\op}$.
\end{cor}

\begin{proof}
Indeed, $\Pro(\FSch_{\af,\,\pr})^{\op}\to 
\Pro'(\FSch_{\af,\,\pr})^{\op}$
is a localisation functor and 
$\RigSH^{\eff}_{\Nis}((-)^{\rig};\Lambda)$
transforms the morphisms 
$(\mathcal{S}_{\alpha}')_{\alpha}\to 
(\mathcal{S}_{\alpha})_{\alpha}$ into equivalences by 
Lemma \ref{lem:localisation-prospace}.
Thus, the result follows from \cite[Proposition 5.2.7.12]{lurie}.
\end{proof}

\begin{rmk}
\label{rmk:const-rigsh-FSmRig-}
In the remainder of this subsection, we use the construction of 
$\RigSH_{\Nis}^{\eff}(S;\Lambda)$ as a localisation of the 
$\infty$-category of presheaves of $\Lambda$-modules
on $\FRigSm/\mathcal{S}$ as explained in 
Remark \ref{rmk:another-site-for-rigsh}. 
In fact, we will rather use the full subcategory of the latter, 
denoted by $\FRigSm_{\af,\,\pr}/\mathcal{S}$,
spanned by formal $\mathcal{S}$-schemes 
which are affine and of principal ideal type.
(If $\mathcal{S}$ is of principal ideal type and 
$\pi$ a generator of an ideal of definition, then 
the second condition is equivalent to having a
$\pi$-torsion-free structure sheaf.) We are free to do so 
since the obvious inclusion induces an equivalence of sites
$(\FRigSm/\mathcal{S},\rigNis)\simeq
(\FRigSm_{\af,\,\pr}/\mathcal{S},\rigNis)$.
We will also need the analogous fact for 
$\RigSH_{\Nis}^{\eff}((S_{\alpha})_{\alpha};\Lambda)$: 
It can be constructed as a localisation of the $\infty$-category of 
presheaves of $\Lambda$-modules on 
$$\FRigSm_{\af,\,\pr}/(\mathcal{S}_{\alpha})_{\alpha}=
\underset{\alpha}{\colim}\, \FRigSm_{\af,\,\pr}/\mathcal{S}_{\alpha}.$$
The above category will be endowed with the limit 
rig-Nisnevich topology
so that the resulting site is equivalent
to the one used in Notation \ref{nota:motive-over-S-alpha}
(with $\tau=\Nis$).
Moreover, \eqref{eq-lem:reduc-effect-nis-non-hyp}
is induced from the obvious functor 
$$\FRigSm_{\af,\,\pr}/(\mathcal{S}_{\alpha})_{\alpha}\to 
\FRigSm_{\af,\,\pr}/\mathcal{S}$$
by the naturality of the construction of 
categories of motives. 
\end{rmk}

\begin{rmk}
\label{rmk:object-in-colim-FRigSm}
(See Remark \ref{rmk:object-in-colim-RigSm}.)
Given a cofiltered inverse system 
of affine formal schemes of principal ideal type
$(\mathcal{S}_{\alpha})_{\alpha}$,
we denote by $\Pro(\FSch_{\af,\,\pr})/(\mathcal{S}_{\alpha})_{\alpha}$
the overcategory of $(\mathcal{S}_{\alpha})_{\alpha}$-objects. 
There is a fully faithful embedding 
\begin{equation}
\label{eq-not:rigsm-prime-1}
\FRigSm_{\af,\,\pr}/(\mathcal{S}_{\alpha})_{\alpha}
\to \Pro(\FSch_{\af,\,\pr})/(\mathcal{S}_{\alpha})_{\alpha}
\end{equation}
and we will identify 
$\FRigSm_{\af,\,\pr}/(\mathcal{S}_{\alpha})_{\alpha}$ with its essential
image by this functor. Thus, we may think of an object of
$\FRigSm_{\af,\,\pr}/(\mathcal{S}_{\alpha})_{\alpha}$ as a pro-object
$(\mathcal{X}_{\alpha})_{\alpha\leq \alpha_0}$,
where $\mathcal{X}_{\alpha_0}$ is a rig-smooth formal 
$\mathcal{S}_{\alpha_0}$-scheme and, for $\alpha\leq \alpha_0$, 
$\mathcal{X}_{\alpha}\simeq \mathcal{X}_{\alpha_0}
\times_{\mathcal{S}_{\alpha_0}}\mathcal{S}_{\alpha}/(0)^{\sat}$.
We set $\mathcal{S}=\lim_{\alpha}\mathcal{S}_{\alpha}$, 
and for an object 
$(\mathcal{X}_{\alpha})_{\alpha\leq \alpha_0}$ as before, we 
set $\mathcal{X}=\lim_{\alpha\leq \alpha_0}\mathcal{X}_{\alpha}$.
\end{rmk}

We now introduce a new category of formal pro-schemes over 
$(\mathcal{S}_{\alpha})_{\alpha}$ where, roughly speaking, 
the endofunctor introduced in Notation
\ref{not:pro-rigspc-af-prime}
becomes an equivalence. We will also consider the 
$\infty$-category of motives associated to this new category of 
formal pro-schemes, and use it to divide the sought after 
equivalence into two which are easier to establish.

\begin{nota}
\label{not:rigsm-prime}
Keep the assumptions as in Remark
\ref{rmk:object-in-colim-FRigSm}.
We denote by $\FRigSm'_{\af,\,\pr}/(\mathcal{S}_{\alpha})_{\alpha}$
the full subcategory of 
$\Pro'(\FSch_{\af,\,\pr})/(\mathcal{S}_{\alpha})_{\alpha}$
spanned by the objects which belong to the image of 
\eqref{eq-not:rigsm-prime-1}. More concretely, we have a functor 
\begin{equation}
\label{eq-not:rigsm-prime-2}
\FRigSm_{\af,\,\pr}/(\mathcal{S}_{\alpha})_{\alpha}
\to \FRigSm'_{\af,\,\pr}/(\mathcal{S}_{\alpha})_{\alpha}
\end{equation}
which is the identity on objects and such that, in the target, 
the set of morphisms from $(\mathcal{Y}_{\alpha})_{\alpha\leq \beta_0}$ 
to $(\mathcal{X}_{\alpha})_{\alpha\leq \alpha_0}$ 
is the set of morphisms from 
$(\mathcal{Y}'_{\alpha})_{\alpha\leq \beta_0}$ 
to $(\mathcal{X}'_{\alpha})_{\alpha\leq \alpha_0}$
over $(\mathcal{S}_{\alpha})_{\alpha}$.
\symn{$\FRigSm'_{\af,\,\pr}$}
\end{nota}

\begin{rmk}
\label{rmk:nisnevich-top-on-sm-prime}
Let $\FRigEt_{\af,\,\pr}/\mathcal{S}$ be the full subcategory 
of $\FRigSm_{\af,\,\pr}/\mathcal{S}$ spanned by rig-\'etale 
formal $\mathcal{S}$-schemes. Similarly, let 
$$\FRigEt_{\af,\,\pr}/(\mathcal{S}_{\alpha})_{\alpha}=
\underset{\alpha}{\colim}\,\FRigEt_{\af,\,\pr}/\mathcal{S}_{\alpha},$$
considered as a full subcategory of 
$\FRigSm_{\af,\,\pr}/(\mathcal{S}_{\alpha})_{\alpha}$,
and let $\FRigEt'_{\af,\,\pr}/(\mathcal{S}_{\alpha})_{\alpha}$
be its essential image by the functor 
\eqref{eq-not:rigsm-prime-2}.
The obvious functors
$$\FRigEt_{\af,\,\pr}/(\mathcal{S}_{\alpha})_{\alpha}
\to \FRigEt'_{\af,\,\pr}/(\mathcal{S}_{\alpha})_{\alpha}
\to \FRigEt_{\af,\,\pr}/\mathcal{S}$$
are equivalences of categories.
Indeed, it is so for their composition by 
Corollary \ref{cor:proj-limi-cal-E-A}, and the second functor 
is faithful. This allows us to define the rig-Nisnevich 
topology on 
$\FRigEt'_{\af,\,\pr}/(\mathcal{S}_{\alpha})_{\alpha}$, 
and more generally on 
$\FRigSm'_{\af,\,\pr}/(\mathcal{S}_{\alpha})_{\alpha}$
by replacing 
$(\mathcal{S}_{\alpha})_{\alpha}$ with a general object of the
latter category.
\symn{$\FRigEt_{\af,\,\pr}$}
\symn{$\FRigEt'_{\af,\,\pr}$}
\end{rmk}

\begin{prop}
\label{prop:extension-to-}
Let $(\mathcal{S}_{\alpha})_{\alpha}$ be a cofiltered inverse 
system of affine formal schemes of principal ideal type.
The functor 
$(\mathcal{X}_{\alpha})_{\alpha\leq \alpha_0}\mapsto 
\M^{\eff}((\mathcal{X}_{\alpha}^{\rig})_{\alpha\leq \alpha_0})$
extends naturally to a functor 
$$\M'^{\eff}(-):\FRigSm'_{\af,\,\pr}/(\mathcal{S}_{\alpha})_{\alpha}
\to \RigSH^{\eff}_{\Nis}((S_{\alpha})_{\alpha};\Lambda).$$
(As usual, we set $S_{\alpha}=\mathcal{S}_{\alpha}^{\rig}$.)
\symn{$\M'^{\eff}$}
\end{prop}

\begin{proof}
By Corollary \ref{cor:ext-rigsh-pro-prime},
there is a functor 
$$\RigSH^{\eff}_{\Nis}((-)^{\rig};\Lambda):
(\FRigSm'_{\af,\,\pr}/(\mathcal{S}_{\alpha})_{\alpha})^{\op}\to \Prl.$$
For every $(\mathcal{X}_{\alpha})_{\alpha\leq \alpha_0}$
in $\FRigSm'_{\af,\,\pr}/(\mathcal{S}_{\alpha})_{\alpha}$, with structure
morphism $f:(\mathcal{X}_{\alpha})_{\alpha\leq \alpha_0}\to 
(\mathcal{S}_{\alpha})_{\alpha}$, the associated inverse image functor 
$f^*$ admits a left adjoint $f_{\sharp}$. 
Moreover, the motive $\M^{\eff}(
(\mathcal{X}^{\rig}_{\alpha})_{\alpha\leq \alpha_0})$
is equivalent to $f_{\sharp}f^*\Lambda$.
Hence, the result follows by applying Lemma
\ref{lem:left-adjoint-to-final-obj} below. 
\end{proof}

\begin{lemma}
\label{lem:left-adjoint-to-final-obj}
Let $\mathcal{C}$ be an $\infty$-category and 
$\mathcal{F}:\mathcal{C}^{\op}\to \CAT_{\infty}$
a functor. Given a morphism $f:Y \to X$ in $\mathcal{C}$, 
we denote by $f^*:\mathcal{F}(X) \to \mathcal{F}(Y)$
the induced functor.
Assume that $\mathcal{C}$ admits a final object 
$\star$ and that for every object $X\in \mathcal{C}$, the 
functor $\pi_X^*$, associated to $\pi_X:X \to \star$, admits 
a left adjoint $\pi_{X,\,\sharp}$. Then, there is a functor 
$\mathcal{C} \to \Fun(\mathcal{F}(\star),\mathcal{F}(\star))$
sending $X\in \mathcal{C}$ to the endofunctor 
$\pi_{X,\,\sharp}\pi_X^*$ and a morphism $f:Y \to X$ to the composition of
$$\pi_{Y,\,\sharp}\pi_Y^* \simeq \pi_{Y,\,\sharp}f^*\pi_X^*
\xrightarrow{\eta} \pi_{Y,\,\sharp}f^*\pi_X^*\pi_{X,\,\sharp}\pi_X^*\simeq 
\pi_{Y,\,\sharp}\pi_Y^*\pi_{X,\,\sharp}\pi_X^*
\xrightarrow{\delta} \pi_{X,\,\sharp}\pi_X^*,$$
where $\eta$ is the unit of the adjunction $(\pi_{X,\,\sharp},\pi_X^*)$
and $\delta$ is the counit of the adjunction 
$(\pi_{Y,\,\sharp},\pi_Y^*)$.
\end{lemma}

\begin{proof}
Let $p:\mathcal{M} \to \mathcal{C}$ be the Cartesian fibration
associated to the functor $\mathcal{F}$ by Lurie's unstraightening construction \cite[\S 3.2]{lurie}. Since $\star$ is a final object
of $\mathcal{C}$, we have a natural transformation 
$\mathcal{F}(\star)_{\cst} \to \mathcal{F}$, 
where $\mathcal{F}(\star)_{\cst}:\mathcal{C}^{\op}\to \CAT_{\infty}$
is the constant functor with value $\mathcal{F}(\star)$.
This natural transformation induces a morphism of 
Cartesian fibrations
$$\xymatrix{\mathcal{F}(\star)\times \mathcal{C}
\ar[rr]^-G \ar[dr]_-q & & \mathcal{M} \ar[dl]^-p\\
& \mathcal{C}. & }$$
The fiber of $G$ over an object $X\in \mathcal{C}$ is 
the functor $\pi_X^*:\mathcal{F}(\star)\to \mathcal{F}(X)$,
which admits a left adjoint by assumption. By 
\cite[Proposition 7.3.2.6]{lurie:higher-algebra},
the functor $G$ admits a left adjoint 
$F$ relative to 
$\mathcal{C}$, in the sense of 
\cite[Definition 7.3.2.2]{lurie:higher-algebra}.
Thus, we have a commutative triangle 
$$\xymatrix{\mathcal{F}(\star)\times \mathcal{C}
\ar[dr]_-q & & \mathcal{M} \ar[dl]^-p \ar[ll]_-F\\
& \mathcal{C} & }$$
and a natural transformation $\id\to G\circ F$
over $\mathcal{C}$ which is a unit map. Moreover, the fiber of 
$F$ over an object $X\in \mathcal{C}$ is the functor 
$\pi_{X,\,\sharp}:\mathcal{F}(X)\to \mathcal{F}(\star)$.

Composing the endofunctor $F\circ G$ of 
$\mathcal{F}(\star)\times \mathcal{C}$ 
with the projection to $\mathcal{F}(\star)$ yields a functor
$\mathcal{F}(\star)\times \mathcal{C}\to \mathcal{F}(\star)$
and, by adjunction, a functor 
$\mathcal{C}\to \Fun(\mathcal{F}(\star),\mathcal{F}(\star))$.
We leave it to the reader to check that 
the latter satisfies the informal description given
in the statement.
\end{proof}

\begin{rmk}
\label{rmk-prop:extension-to-}
Let $(\mathcal{S}_{\alpha})_{\alpha}$ be a cofiltered inverse 
system of affine formal schemes of principal ideal type and 
$\mathcal{S}=\lim_{\alpha}\mathcal{S}_{\alpha}$.
We set $S_{\alpha}=\mathcal{S}^{\rig}_{\alpha}$ and 
$S=\mathcal{S}^{\rig}$.

\begin{enumerate}

\item[(1)] 
There is a commutative diagram
$$\xymatrix{\FRigSm_{\af,\,\pr}/(\mathcal{S}_{\alpha})_{\alpha}
\ar[r] \ar[dr]_-{\M^{\eff}((-)^{\rig})} & 
\FRigSm'_{\af,\,\pr}/(\mathcal{S}_{\alpha})_{\alpha} 
\ar[r] \ar[d]^-{\M'^{\eff}(-)} & \FRigSm_{\af,\,\pr}/\mathcal{S} 
\ar[d]^-{\M^{\eff}((-)^{\rig})} \\
& \RigSH^{\eff}_{\Nis}((S_{\alpha})_{\alpha};\Lambda)
\ar[r] & \RigSH^{\eff}_{\Nis}(S;\Lambda).\!}$$
This is not completely obvious. One needs to check 
that Lemma \ref{lem:left-adjoint-to-final-obj}
applied to the contravariant functor 
$\RigSH^{\eff}_{\Nis}((-)^{\rig};\Lambda)$ defined on 
$\FRigSm^{\af}/(\mathcal{S}_{\alpha})_{\alpha}$ and
$\FRigSm^{\af}/\mathcal{S}$
gives back the functor 
$\M^{\eff}((-)^{\rig})$.
To do so, one reduces to a similar question,
but for the contravariant functor $\RigSm/(-)^{\rig}$, 
which can be easily handled.

\item[(2)] It follows from the commutative triangle inside the
diagram in (1) that $\M'^{\eff}$ admits descent for the rig-Nisnevich 
topology, i.e., it takes a truncated hypercover for the 
rig-Nisnevich topology to a colimit diagram.
(See Remark \ref{rmk:nisnevich-top-on-sm-prime}.)

\item[(3)] By the universal properties of presheaf categories 
and localisation, the commutative diagram in (1) 
gives rise to a commutative diagram in $\Prl$:
$$\xymatrix{\RigSH^{\eff}_{\Nis}((S_{\alpha})_{\alpha};\Lambda)
\ar[r] \ar@{=}[dr] & \RigSH'^{\eff}_{\Nis}
((S_{\alpha})_{\alpha};\Lambda) 
\ar[d] \ar[dr] &  \\
& \RigSH^{\eff}_{\Nis}((S_{\alpha})_{\alpha};\Lambda)
\ar[r] & \RigSH^{\eff}_{\Nis}(S;\Lambda)}$$
where 
$\RigSH'^{\eff}_{\Nis}((S_{\alpha})_{\alpha};\Lambda)$
is defined from the site 
$(\FRigSm'_{\af,\,\pr}/(\mathcal{S}_{\alpha})_{\alpha},\rigNis)$
in the usual way, i.e., by adapting Definition 
\ref{def:DAeff}.
Thus, to finish the proof of Proposition
\ref{prop:approx-rig-sh}, it suffices to show 
Proposition \ref{prop:approx-rig-sh-prime} below.

\end{enumerate}
\end{rmk}

\begin{prop}
\label{prop:approx-rig-sh-prime}
Let $(\mathcal{S}_{\alpha})_{\alpha}$ be a cofiltered inverse 
system of affine formal schemes of principal ideal type and 
$\mathcal{S}=\lim_{\alpha}\mathcal{S}_{\alpha}$.
We set $S_{\alpha}=\mathcal{S}^{\rig}_{\alpha}$ and 
$S=\mathcal{S}^{\rig}$.
Then the obvious functor
\begin{equation}
\label{eqn:funct-S-alpha-alpha-to-S-prime}
\RigSH'^{\eff}_{\Nis}((S_{\alpha})_{\alpha};\Lambda)
\to \RigSH_{\Nis}^{\eff}(S;\Lambda)
\end{equation}
is an equivalence.
\end{prop}

\begin{nota}
\label{nota:A-A-alpha-A-infty-}
From now on, we fix a cofiltered inverse system 
$(\mathcal{S}_{\alpha})_{\alpha}$ of 
affine formal schemes of principal ideal type, 
and we let $\mathcal{S}=\lim_{\alpha}
\mathcal{S}_{\alpha}$. We define $(\mathcal{S}'_{\alpha})_{\alpha}$
as in Notation \ref{not:closed-prospace}, and we set 
$S_{\alpha}=\mathcal{S}_{\alpha}^{\rig}$, 
$S'_{\alpha}=\mathcal{S}'^{\rig}_{\alpha}$
and $S=\mathcal{S}^{\rig}$.
We set $A_{\alpha}=\mathcal{O}(\mathcal{S}_{\alpha})$, $A'_{\alpha}=\mathcal{O}(\mathcal{S}_{\alpha}')$ and
$A=\mathcal{O}(\mathcal{S})$. We identify $A'_{\alpha}$
with a subring of $A$ and set
$A'_{\infty}=\bigcup_{\alpha}\,A'_{\alpha}$
which is a dense subring of $A$.
We also assume that there is an element $\pi$,
which ``belongs'' to all the $A_{\alpha}$'s and generates 
an ideal of definition in each $A_{\alpha}$.
(This is not a restrictive assumption since it is clearly satisfied
when the indexing category of $(\mathcal{S}_{\alpha})_{\alpha}$
admits a final object.) Given 
$(\mathcal{X}_{\alpha})_{\alpha\leq \alpha_0}$
in $\FRigSm_{\af,\,\pr}/(\mathcal{S}_{\alpha})_{\alpha}$, 
we use similar notations: $B_{\alpha}=
\mathcal{O}(\mathcal{X}_{\alpha})$, 
$B'_{\alpha}=\mathcal{O}(\mathcal{X}'_{\alpha})$, 
$B=\mathcal{O}(\mathcal{X})$ and 
$B'_{\infty}=\bigcup_{\alpha\leq \alpha'}\,B'_{\alpha}$
which is a dense subring of $B$.
\end{nota}

\begin{rmk}
\label{rmk:compact-gen-A-infty-pol}
The $\infty$-category 
$\RigSH'^{\eff}_{\Nis}((S_{\alpha})_{\alpha};\Lambda)$
is compactly generated, up to desuspension, by 
$\M'^{\eff}((\mathcal{X}_{\alpha})_{\alpha\leq \alpha_0})$
where $(\mathcal{X}_{\alpha})_{\alpha\leq \alpha_0}$ 
belongs to $\FRigSm'_{\af,\,\pr}/(\mathcal{S}_{\alpha})_{\alpha}$. 
(This can be proven by adapting the proof of 
Proposition \ref{prop:compact-shv-rigsm}. The key point is that 
the small rig-Nisnevich site of 
$(\mathcal{X}_{\alpha})_{\alpha\leq \alpha_0}$
is equivalent to the small Nisnevich site of $X$; see 
Remark \ref{rmk:nisnevich-top-on-sm-prime}.)
Using Proposition \ref{prop:compact-shv-rigsm}, 
we deduce that the functor \eqref{eqn:funct-S-alpha-alpha-to-S-prime}
belongs to $\Prl_{\omega}$. 
This functor also sends a set of compact generators to 
a set of compact generators. 
Indeed, by Proposition
\ref{prop:on-etale-presentation}, a set of compact generators 
for $\RigSH_{\Nis}^{\eff}(S;\Lambda)$ is given, up to desuspension, 
by motives of smooth rigid $S$-affinoids $X=\Spf(B)^{\rig}$ 
with $B$ of the form
$$B=A\langle s_1,\ldots, s_m,t_1,\ldots, t_n\rangle/
(P_1,\ldots, P_n)^{\sat}$$
with $P_i\in A'_{\infty}[s_1,\ldots, s_m,t_1,\ldots, t_n]$ such that 
$\det(\partial P_i/\partial t_j)$
generates an open ideal in $B$.
Clearly, $\Spf(B)$ is in the image of 
$\FRigSm_{\af,\,\pr}/(\mathcal{S}_{\alpha})_{\alpha}
\to \FRigSm_{\af,\,\pr}/\mathcal{S}$.
In particular, to prove that the functor 
\eqref{prop:approx-rig-sh-prime}
is an equivalence, it remains to show that it is fully faithful. 
\end{rmk}

Before continuing with the proof, 
we recall the following two statements from \cite{vezz-fw}.

\begin{prop}
\label{prop:vezz-A1}
Let $R$ be an adic ring of principal ideal type and 
$\pi\in R$ a generator of an ideal of definition.
Let $s=(s_1,\ldots, s_m)$ 
and $t=(t_1,\ldots, t_n)$ be two systems of coordinates and let 
$P=(P_1,\ldots, P_n)$ be an $n$-tuple of polynomials 
in $R[s,t]$ with no constant term, i.e., 
such that $P|_{s=0,\,t=0}=(0,\ldots, 0)$. Assume also that 
$\det(\partial P_i/\partial t_j)|_{s=0,\,t=0}$ generates an open ideal in 
$R$. Then, there exists a unique $n$-tuple
$F=(F_1,\ldots, F_n)$ of formal power series in 
$(R[\pi^{-1}])[[s]]$
such that $P(s,F(s))=0$. Moreover, for $N$ large enough, 
the $F_i$'s belong to the subring
$R[[\pi^{-N}s]]$. 
\end{prop}

\begin{proof}
This is a slight generalisation of 
\cite[Proposition A.1]{vezz-fw}
and one can easily check that the proof 
of loc.~cit. still works in the present context. 
More precisely, instead of a Banach $K$-algebra, with 
$K$ a complete non-Archimedean field, as in loc.~cit., 
we consider the Banach ring $R[\pi^{-1}]$ 
endowed with the norm described in the proof of 
Proposition \ref{prop:approx-solution}.
(Note that $\det(\partial P_i/\partial t_j)|_{s=0,\,t=0}$ 
generates an open ideal
in $R$ if and only if it is invertible in $R[\pi^{-1}]$.)
\end{proof}

The previous statement has the following 
generalisation. (See \cite[Proposition A.2]{vezz-fw}.)

\begin{cor}
\label{cor:vezz-A2}
Let $R$ be an adic ring of principal ideal type and 
$\pi\in R$ a generator of an ideal of definition. 
Let $s=(s_1,\ldots, s_m)$ 
and $t=(t_1,\ldots, t_n)$ be two systems of coordinates,
let $a=(a_1,\ldots, a_m)$ 
and $b=(b_1,\ldots, b_n)$ be two tuples of elements in $R$, 
and let $P=(P_1,\ldots, P_n)$ be an $n$-tuple of polynomials 
in $R[s,t]$ such that $P|_{s=a,\,t=b}=(0,\ldots, 0)$. 
Assume also that $\det(\partial P_i/\partial t_j)|_{s=a,\,t=b}$ 
generates an open ideal in $R$. Then, there exists a unique $n$-tuple
$F=(F_1,\ldots, F_n)$ of formal power series in $(R[\pi^{-1}])[[s-a]]$
such that $P(s,F(s))=0$. Moreover, for $N$ large enough, 
the $F_i$'s belong to the subring $R[[\pi^{-N}(s-a)]]$. 
\end{cor}

We introduce some further notations.

\begin{nota}
\label{not:S-X-Y-A-B-C}
We fix two $\pi$-torsion-free
rig-smooth adic $A_{\alpha_0}$-algebras $B_{\alpha_0}$
and $C_{\alpha_0}$. For $\alpha\leq \alpha_0$, we set 
$B_{\alpha}=A_{\alpha}
\,\widehat{\otimes}_{A_{\alpha_0}} 
\,B_{\alpha_0}/(0)^{\sat}$, 
$C_{\alpha}=A_{\alpha}\,\widehat{\otimes}_{A_{\alpha_0}}
\,C_{\alpha_0}/(0)^{\sat}$,
$\mathcal{X}_{\alpha}=\Spf(B_{\alpha})$ and 
$\mathcal{Y}_{\alpha}=\Spf(C_{\alpha})$. 
Similarly, we set
$B=A\,\widehat{\otimes}_{A_{\alpha_0}}\,B_{\alpha_0}/(0)^{\sat}$,
$C=A\,\widehat{\otimes}_{A_{\alpha_0}}\,C_{\alpha_0}/(0)^{\sat}$, 
$\mathcal{X}=\Spf(B)$ and
$\mathcal{Y}=\Spf(C)$. 
We also denote by $B'_{\alpha}$, $B'_{\infty}$ and 
$\mathcal{X}'_{\alpha}$ as in Notation 
\ref{nota:A-A-alpha-A-infty-}, and we define similarly 
$C'_{\alpha}$, $C'_{\infty}$ and $\mathcal{Y}'_{\alpha}$.
Moreover, we assume that
$$B_{\alpha_0}=A_{\alpha_0}\langle s,t\rangle/(P)^{\sat}$$
with $s=(s_1,\ldots, s_m)$ and $t=(t_1,\ldots, t_n)$ 
two systems of coordinates, and $P=(P_1,\ldots, P_n)$ an 
$n$-tuple of polynomials in $A_{\alpha_0}[s,t]$ such that 
$\det(\partial P_i/\partial t_j)$ generates an open ideal 
of $A_{\alpha_0}$. 
\end{nota}

\begin{lemma}
\label{lem:approx-desc}
Given a morphism of formal schemes 
$f:\mathcal{Y}\to \mathcal{X}$, 
there exists an $\A^1$-homotopy 
$$H:\A^1_{\mathcal{Y}}=\Spf(C\langle \tau\rangle)
\to \mathcal{X}$$
from $f=H\circ i_0$ to a map $\tilde{f}=H\circ i_1$ such that 
$\tilde{f}:\mathcal{Y}\to \mathcal{X}$ descends to 
a unique map $\mathcal{Y}'_{\alpha}\to \mathcal{X}_{\alpha}$
for $\alpha\leq \alpha_0$ small enough.
\end{lemma}

\begin{proof}
Indeed, suppose that $f$ corresponds to a morphism of adic 
$A$-algebras $B\to C$ given by
$s_i\mapsto c_i$, for $1\leq i \leq m$, and 
$t_j\mapsto d_j$, for $1\leq j \leq n$, where the $c=(c_1,\ldots, c_m)$
and $d=(d_1,\ldots, d_n)$ are tuples of elements of $C$ 
satisfying $P(c,d)=0$.
Let $F=(F_1,\ldots, F_n)$ be the $n$-tuple of power series
in $C[\pi^{-1}][[s-c]]$ associated by Corollary
\ref{cor:vezz-A2} to the $n$-tuple of 
polynomials $P=(P_1,\ldots, P_n)$ 
(considered with coefficients in $C$ via the map
$A_{\alpha_0}\to C$) and their common zero $(c,d)$.
By the same corollary, for 
$\tilde{c}=(\tilde{c}_1,\ldots, \tilde{c}_m)$ 
an $m$-tuple of elements in $A$ close enough to $c$, the expressions
$F_i(c+(\tilde{c}-c)\cdot \tau)$ are well-defined elements of
$C\langle \tau\rangle$, and the assignment
$$s\mapsto c+(\tilde{c}-c)\cdot \tau,\quad 
t\mapsto F(c+(\tilde{c}-c)\cdot \tau)$$
gives rise to a map of $A$-algebras 
$B\to C\langle \tau\rangle$, and hence 
to a morphism $H:\A_{\mathcal{Y}}^1\to \mathcal{X}$ 
of formal schemes. By construction, 
$H\circ i_0=f$, and it remains to show that 
$\tilde{f}=H\circ i_1$ descends to a morphism 
$\mathcal{Y}'_{\alpha}\to 
\mathcal{X}_{\alpha}$ for a well-chosen $m$-tuple 
$\tilde{c}$. (The uniqueness is clear since $C'_{\alpha}\to C$
is injective.) 
This is the case when the $\tilde{c}_i$'s 
belong to the dense subring $C'_{\infty}=\bigcup_{\alpha\leq \alpha_0}\,
C'_{\alpha}$ of $C$. 
Indeed, refining $\alpha_0$, we may assume that the 
$\tilde{c}_i$'s belong to $C'_{\alpha_0}$. 
Consider the map $\mathcal{Y}'_{\alpha_0}\to  
\mathcal{S}_{\alpha_0}\times \A^m=
\Spf(A_{\alpha}\langle s\rangle)$ induced by $\tilde{c}$.
We have a rig-\'etale morphism $\mathcal{X}_{\alpha_0}\to 
\mathcal{S}_{\alpha_0}\times \A^m$ and the morphism
$\tilde{f}:\mathcal{Y}\to \mathcal{X}$ gives rise to 
a section $\sigma$ of the rig-\'etale projection 
$\mathcal{X}_{\alpha_0}
\times_{\mathcal{S}_{\alpha_0}\times \A^m,\,\tilde{c}}\mathcal{Y} 
\to \mathcal{Y}$. 
Then $\tilde{f}$ descends to 
a morphism $\mathcal{Y}'_{\alpha}\to \mathcal{X}_{\alpha}$
if and only if the section $\sigma$ descends to a section
of the rig-\'etale projection 
$(\mathcal{X}_{\alpha_0}\times_{\mathcal{S}_{\alpha_0}\times \A^m,\,\tilde{c}}\mathcal{Y}'_{\alpha}) \to 
\mathcal{Y}'_{\alpha}$. 
That this is true follows from Corollary
\ref{cor:proj-limi-cal-E-A}.
\end{proof}

\begin{cor}
\label{cor:for-constructing-homotopies}
Keep the notation as above.
Fix a system of coordinates $u=(u_1,\ldots,u_r)$ for 
$\A^r$.
Given a finite collection $f_1,\ldots,f_N$ in 
$\Hom_{\mathcal{S}}(\mathcal{Y}\times \A^r,\mathcal{X})$ 
we can find a collection $H_1,\ldots,H_N$ in 
$\Hom_{\mathcal{S}}(\mathcal{Y}\times \A^r\times \A^1, \mathcal{X})$ 
and some index $\alpha\leq \alpha_0$ such that:
\begin{enumerate}

\item[(1)] For all $1\leq k\leq N$, we have $f_k=H_k\circ i_0$ and 
the map $\tilde{f}_k=H_k\circ i_1$ descends to a unique map 
$\mathcal{Y}'_{\alpha}\times \A^r\to \mathcal{X}_{\alpha}$ over 
$\mathcal{S}_{\alpha}$.

\item[(2)] If $f_k\circ d_{i,\epsilon}=f_{k'}\circ d_{i,\epsilon}$ for some
$1\leq k,k'\leq N$ and some $(i,\epsilon)\in\{1,\ldots,r\}\times\{0,1\}$ 
then $H_k\circ d_{i,\epsilon}=H_{k'}\circ d_{i,\epsilon}$.

\item[(3)] If for some $1\leq k\leq N$ and some $\gamma\leq \alpha_0$ 
the map $f_k\circ d_{1,1}\in 
\Hom_{\mathcal{S}}(\mathcal{Y}\times \A^{r-1},\mathcal{X})$
comes from $\Hom_{\mathcal{S}_{\gamma}}(\mathcal{Y}'_{\gamma}
\times\A^{r-1},\mathcal{X}_{\gamma})$, then the homotopy 
$H_k\circ d_{1,1}\in
\Hom_{\mathcal{S}}(\mathcal{Y}\times \A^{r-1}\times \A^1,\mathcal{X})$ 
is constant, i.e., factors through the projection on
$\mathcal{Y}\times \A^{r-1}$.

\end{enumerate}
\end{cor}

\begin{proof}
Suppose that $f_k$ corresponds to a morphism of adic $A$-algebras
$B\to C\langle u \rangle$ given by $(s,t)\mapsto (c_k,d_k)$ 
where $c_k=(c_{k1},\ldots, c_{km})$ and $d_k=(d_{k1},\ldots, d_{kn})$
are tuples of elements of $C\langle u \rangle$ satisfying 
$P(c_k,d_k)=0$. By Lemma \ref{lem:approx-desc}, 
there are $n$-tuples of formal power
series $F_k=(F_{k1},\ldots,F_{kn})$ associated to the $f_k$'s
such that 
$$(s,t)\mapsto (c_k+(\tilde{c}_k-c)\cdot \tau,
F_k(c_k+(\tilde{c}_k-c_k)\cdot \tau))$$
defines a morphism $H_k:\mathcal{Y}\times \B^r\times \B^1
\to \mathcal{X}$ satisfying condition (1), 
for some $\alpha\leq \alpha_0$, when the $\tilde{c}_{ki}$'s are close 
enough to the $c_{ki}$'s and belong to the dense subring 
$C'_{\infty}\langle u\rangle=\bigcup_{\alpha\leq \alpha_0}C'_{\alpha}\langle u\rangle$ of $C\langle u\rangle$.

It remains to explain how to choose the $\tilde{c}_k$'s so
that the conditions (2) and (3) above are also satisfied.
To do so, we apply \cite[Proposition A.5]{vezz-fw}
to the $c_{ki}$'s. (This result of \cite{vezz-fw} 
is stated for Banach algebras over a non-Archimedean field
and a sequence of complete subalgebras, 
but holds more generally for Banach rings and a filtered 
family of complete subrings; and we apply it 
here to $C[\pi^{-1}]$ and the family 
$C'_{\alpha}[\pi^{-1}]$, for $\alpha \leq \alpha_0$.)
Thus we may find elements $\tilde{c}_{ki}\in 
C'_{\infty}\langle u\rangle$, which are arbitrary close to the $c_{ki}$'s, 
and satisfying the following properties:
\begin{enumerate}

\item[(2$'$)] If $c_k|_{u_i=\epsilon}=c_{k'}|_{u_i=\epsilon}$
for some
$1\leq k,k'\leq N$ and some $(i,\epsilon)\in\{1,\ldots,r\}\times\{0,1\}$ 
then $\tilde{c}_k|_{u_i=\epsilon}=\tilde{c}_{k'}|_{u_i=\epsilon}$.

\item[(3$'$)] If for some $1\leq k\leq N$ and some $\gamma\leq \alpha_0$,
$c_k|_{u_1=1}$ belongs to $C'_{\gamma}\langle u_2,\ldots,u_r\rangle$, 
then $\tilde{c}_k|_{u_1=1}=c_k|_{u_1=1}$.

\end{enumerate}
With these $\tilde{c}_{ki}$'s, it is easy to see that conditions
(2) and (3) are satisfied. Indeed, suppose that
$f_k\circ d_{i,\epsilon}=f_{k'}\circ d_{i,\epsilon}$ for some 
$i\in\{1,\ldots,r\}$ and $\epsilon\in\{0,1\}$. This means that 
$c_k|_{u_i=\epsilon}=c_{k'}|_{u_i=\epsilon}$ 
and $d_k|_{u_i=\epsilon}=d_{k'}|_{u_i=\epsilon}$;
we denote by $\bar{c}$ and $\bar{d}$ their respective common values. 
This implies that both $F_k|_{u_i=\epsilon}$ and 
$F_{k'}|_{u_i=\epsilon}$ are two $n$-tuples of formal power 
series $\bar{F}$ with coefficients in $C\langle u_2,\ldots, u_r\rangle$ converging around $\bar{c}$ and such that $P(s,\bar{F}(s))=0$ and
$\bar{F}(\bar{c})=\bar{d}$. By the uniqueness of such power series 
stated in Corollary 
\ref{cor:vezz-A2}, 
we conclude that they coincide. 
Moreover, by property (2$'$), we have 
$\tilde{c}_{k}|_{u_i=\epsilon}=
\tilde{c}_{k'}|_{u_i=\epsilon}$; we denote by $\bar{\tilde{c}}$
the common value. It follows that  
$$F_k(c_k+(\tilde{c}_{k}-c_k)\cdot \tau)|_{u_i=\epsilon}=
\bar{F}(\bar{c}+(\bar{\tilde{c}}-\bar{c})\cdot \tau)=
F_{k'}(c_{k'}+(\tilde{c}_{k'}-c_{k'})\cdot \tau)|_{\theta_r=\epsilon}$$
and thus $H_k\circ d_{i,\epsilon}=H_{k'}\circ d_{i,\epsilon}$ 
proving property (2). Property (3) follows immediately from 
property (3$'$) and the definition of $H_k$.
\end{proof}

\begin{proof}[Proof of Proposition \ref{prop:approx-rig-sh-prime}]
We split the argument into two steps.

\paragraph*{Step 1}
\noindent
Consider the $\A^1$-localisation functor 
$\Lder_{\A^1}$ on the $\infty$-categories 
of presheaves of $\Lambda$-modules
$$\PSh(\FRigSm'_{\af,\,\pr}/(\mathcal{S}_{\alpha})_{\alpha};\Lambda) 
\qquad \text{and}
\qquad \PSh(\FRigSm_{\af,\,\pr}/\mathcal{S};\Lambda).$$
For a presheaf $\mathcal{F}$ of $\Lambda$-modules, 
$\Lder_{\A^1}(\mathcal{F})$ is given by the 
colimit of the simplicial presheaf 
$\underline{\Hom}(\Delta^{\bullet},\mathcal{F})$ where
$\Delta^r$ refers to the $r$-th algebraic simplex and 
$$\underline{\Hom}(\Delta^r,\mathcal{F})(-)=
\mathcal{F}((-)\langle u_0,\ldots, u_r\rangle/(u_0+\cdots+u_r-1)).$$
Indeed, the map 
$\mathcal{F}\to \colim\,\underline{\Hom}(\Delta^{\bullet},\mathcal{F})$
is an $\A^1$-equivalence by
\cite[\S 2.3, Corollary 3.8]{mv-99}.
On the other hand, using \cite[\S 2.3, Proposition 3.4]{mv-99}
and the fact that the endofunctor 
$\underline{\Hom}(\Delta^1,-)$ preserves colimits,
we have equivalences
$$\colim\,
\underline{\Hom}(\Delta^{\bullet},\mathcal{F})
\simeq \colim\,
\underline{\Hom}(\Delta^{\bullet}\times \Delta^1,\mathcal{F})
\simeq 
\underline{\Hom}(\Delta^1,\colim\,
\underline{\Hom}(\Delta^{\bullet},\mathcal{F}))$$
showing that 
$\colim\,\underline{\Hom}(\Delta^{\bullet},\mathcal{F})$
is $\A^1$-local.

With $(\mathcal{X}_{\alpha})_{\alpha\leq \alpha_0}$ and
$(\mathcal{Y}_{\alpha})_{\alpha\leq \alpha_0}$ as in 
Notation \ref{not:S-X-Y-A-B-C}, we claim that the natural map
\begin{equation}
\label{eq:sing-b-1-approx}
\left(\Lder_{\A^1}\Lambda
((\mathcal{X}_{\alpha})_{\alpha\leq \alpha_0})\right)
((\mathcal{Y}_{\alpha})_{\alpha\leq \alpha_0})\to 
\left(\Lder_{\A^1}\Lambda(\mathcal{X})\right)
(\mathcal{Y})
\end{equation}
is an equivalence. By the commutation of colimits with 
tensor products, it is enough to prove this when 
$\Lambda$ is the sphere spectrum. (Here we use the explicit 
model for the $\A^1$-localisation recalled above.) 
Similarly, since tensoring with 
the Eilenberg--Mac Lane spectrum of $\Z$ is 
conservative on connective spectra, 
we reduce to prove this when $\Lambda$ is 
the (Eilenberg--Mac Lane 
spectrum associated to the) ring $\Z$. 
In this case, we may use another model for 
the $\A^1$-localisation functor $\Lder_{\A^1}$, namely 
the one taking $\mathcal{F}$ to the normalised complex 
associated to the cubical presheaf of complexes of abelian groups
$\underline{\Hom}(\A^{\bullet},\mathcal{F})$ where, 
as above, $\underline{\Hom}(\A^r,\mathcal{F})(-)=
\mathcal{F}((-)\langle u_1,\ldots, u_r\rangle)$.
(This is proven by adapting the method used for the simplicial presheaf 
$\underline{\Hom}(\Delta^{\bullet},\mathcal{F})$.
See also \cite[Th\'eor\`eme 2.23]{ayoub-h1}
for a closely related result.)
Thus, we are reduced to showing that the morphism of cubical
abelian groups 
\begin{equation}
\label{eq:sing-b-1-approx-2w}
\left(\underline{\Hom}(\A^{\bullet},
\Z((\mathcal{X}_{\alpha})_{\alpha\leq \alpha_0}))\right)
((\mathcal{Y}_{\alpha})_{\alpha\leq \alpha_0}) \to
\left(\underline{\Hom}(\A^{\bullet},\Z(\mathcal{X}))\right)
(\mathcal{Y})
\end{equation}
induces an isomorphism on the associated normalised complexes.
This follows from Corollary 
\ref{cor:for-constructing-homotopies}
by arguing as in \cite[Proposition 4.2]{vezz-fw}.
Note that, since $\Z((\mathcal{X}_{\alpha})_{\alpha\leq \alpha_0})$
is considered as a presheaf on 
$\FRigSm'_{\af,\,\pr}/(\mathcal{S}_{\alpha})_{\alpha}$, 
the elements of the left-hand side of 
\eqref{eq:sing-b-1-approx-2w} are linear combinations 
of $(\mathcal{S}_{\alpha})_{\alpha}$-morphisms of formal pro-schemes
from $(\mathcal{Y}'_{\alpha}\times \A^r)_{\alpha\leq \alpha_0}$ to 
$(\mathcal{X}_{\alpha})_{\alpha\leq \alpha_0}$.

\paragraph*{Step 2}
\noindent
Let $\phi:(\FRigSm_{\af,\,\pr}/\mathcal{S},\rigNis)
\to (\FRigSm'_{\af,\,\pr}/(\mathcal{S}_{\alpha})_{\alpha},\rigNis)$
be the premorphism of sites that gives rise to the 
adjunction 
$$\phi^*_{\mot}:\RigSH'^{\eff}_{\Nis}((S_{\alpha})_{\alpha};\Lambda) 
\rightleftarrows \RigSH_{\Nis}^{\eff}(S;\Lambda):\phi_{\mot,\,*}.$$
Our goal is to show that $\phi^*_{\mot}$ is an equivalence, and 
by Remark \ref{rmk:compact-gen-A-infty-pol} it remains to 
see that $\phi^*_{\mot}$ is fully faithful. 
We will prove that the unit morphism 
$\id\to \phi_{\mot,\,*}\phi_{\mot}^*$ is an equivalence. 
In order to do so, we note that the functor 
$$\phi_*:\PSh(\FRigSm_{\af,\,\pr}/\mathcal{S};\Lambda)
\to \PSh(\FRigSm'_{\af,\,\pr}/
(\mathcal{S}_{\alpha})_{\alpha};\Lambda)$$
preserves $(\A^1,\rigNis)$-local equivalences.
Preservation of $\rigNis$-local equivalences
follows immediately from Remark \ref{rmk:nisnevich-top-on-sm-prime}.
Preservation of $\A^1$-local equivalences is an easy consequence of 
the fact that $\A^1$ is an interval. (This is used 
to construct an explicit $\A^1$-homotopy between the 
identity of $\phi_*\Lambda((-)\times \A^1)$
and the endomorphism induced by the zero section.)
As a consequence, we are left to show that the morphism
$\mathcal{F} \to \phi_*\phi^*\mathcal{F}$ 
is an $(\A^1,\rigNis)$-local equivalence for
all presheaves of $\Lambda$-modules $\mathcal{F}$ on 
$\FRigSm'_{\af,\,\pr}/(\mathcal{S}_{\alpha})_{\alpha}$.
Since $\phi^*$ and $\phi_*$ commute with colimits, and since
$(\A^1,\rigNis)$-local equivalences are preserved by colimits, 
we may assume that $\mathcal{F}=
\Lambda((\mathcal{X}_{\alpha})_{\alpha\leq \alpha_0})$
with $(\mathcal{X}_{\alpha})_{\alpha\leq \alpha_0}$
as in Notation \ref{not:S-X-Y-A-B-C}. 
In this case, the morphism 
$\mathcal{F} \to \phi_*\phi^*\mathcal{F}$ 
can be rewritten as follows:
\begin{equation}
\label{eq:sing-b-1-approx-ue0}
\Lambda((\mathcal{X}_{\alpha})_{\alpha\leq \alpha_0})
\to \phi_*\Lambda(\mathcal{X}).
\end{equation}
We claim that this morphism 
is an $\A^1$-local equivalence.
Indeed, if we apply $\Lder_{\A^1}$ to 
\eqref{eq:sing-b-1-approx-ue0}
and if we evaluate at an object 
$(\mathcal{Y}_{\alpha})_{\alpha\leq \alpha_0}$
of $\FRigSm'_{\af,\,\pr}/
(\mathcal{S}_{\alpha})_{\alpha}$, 
we get precisely the map \eqref{eq:sing-b-1-approx}
which we know to be an equivalence.
\end{proof}

\subsection{Quasi-compact base change}

$\empty$

\smallskip

\label{subsect:quasi-compact-base-change}

We prove here the so-called 
quasi-compact base change theorem for rigid analytic 
motives. This will be obtained as an application 
of the continuity property for 
$\RigSH^{(\eff,\,\hyp)}_{\tau}(-;\Lambda)$
proved in Theorem \ref{thm:anstC}.
Our quasi-compact base change theorem can be compared with 
\cite[Proposition 4.4.1]{huber} and \cite[Theorems 5.3.1]{dJ-vdP}.

\begin{thm}[Quasi-compact base change]
\label{thm:general-base-change-thm}
\ncn{quasi-compact base change}
Consider a Cartesian square of rigid analytic spaces
$$\xymatrix{Y' \ar[r]^-{g'} \ar[d]^-{f'} & Y \ar[d]^-f\\
X' \ar[r]^-g & X}$$
with $f$ quasi-compact and quasi-separated. 
Let $\tau\in \{\Nis,\et\}$, and assume one of the 
following two alternatives. 
\begin{enumerate}

\item[(1)] We work in the non-hypercomplete case.
When $\tau$ is the \'etale topology, we assume furthermore that 
$\Lambda$ is eventually coconnective.

\item[(2)] 
We work in the hypercomplete case, and $X$, $X'$, $Y$ and $Y'$ are
$(\Lambda,\tau)$-admissible. When $\tau$ is the \'etale topology, 
we assume furthermore one of the following conditions:
\begin{itemize}

\item $\Lambda$ is eventually 
coconnective;

\item locally on $X$ and $X'$, one can find 
formal models $\mathcal{X}$ and $\mathcal{X}'$ such that 
$\mathcal{X}'$ is a limit of a cofiltered inverse system 
of finite type formal $\mathcal{X}$-schemes 
$(\mathcal{X}_{\alpha})_{\alpha}$ 
with affine transition morphisms and such that the numbers  
$\pvcd_{\Lambda}(\mathcal{X}^{\rig}_{\alpha})$ are  
bounded independently of $\alpha$. (For example, this 
holds if $g$ is locally of finite type.)

\end{itemize}

\end{enumerate}
Then, the commutative square
$$\xymatrix{\RigSH^{(\eff,\,\hyp)}_{\tau}(X;\Lambda) \ar[r]^-{f^*} 
\ar[d]^-{g^*} & \RigSH^{(\eff,\,\hyp)}_{\tau}(Y;\Lambda) \ar[d]^-{g'^*}\\
\RigSH^{(\eff,\,\hyp)}_{\tau}(X';\Lambda) \ar[r]^-{f'^*} & 
\RigSH^{(\eff,\,\hyp)}_{\tau}(Y';\Lambda) }$$
is right adjointable, i.e., the natural transformation 
$g^* \circ f_* \to f'_*\circ g'^*$ is an equivalence.
\end{thm}

\begin{proof}
Using Proposition \ref{prop:6f1}(3),
the problem is local on $X$ and $X'$.
In particular, we may assume 
that $X$ and $X'$ are quasi-compact and quasi-separated.
This implies the same for $Y$ and $Y'$.
We split the proof into two parts. In the first part, we assume that 
$g$ is of finite type and, in the second part, we explain how to 
remove this assumption.

\paragraph*{Part 1}
\noindent
Here we assume that $g$ is of finite type.
Since the problem is local on $X$ and $X'$, we may assume that 
$g$ factors as a closed immersion followed by a smooth morphism.  
Using the base change theorem for smooth morphisms
of Proposition \ref{prop:6f1}, we reduce to the 
case where $g$ is a closed immersion.
Thus, we may assume that $X=\Spf(A)^{\rig}$ and 
$X'=\Spf(A')^{\rig}$ where $A$ is an adic ring of principal 
ideal type and $A'$ a quotient of $A$ by a closed saturated ideal 
$I\subset A$. If $\pi\in A$ generates an ideal of definition, 
then $A'$ is the filtered colimit in the category of adic rings
of the $A$-algebras $A_{J,\,N}=A\langle J/\pi^N\rangle$
where $N\in \N$ and $J\subset I$ is a finitely generated ideal.

Set $\mathcal{X}=\Spf(A)$ and $\mathcal{X}'=\Spf(A')$.
Choose a formal model $\mathcal{Y}$ of $Y$ which is 
a formal $\mathcal{X}$-scheme and set $\mathcal{Y}'=\mathcal{Y}\times_{\mathcal{X}}\mathcal{X}'$. Let $K$
be the indexing category of the filtered inductive 
system $(A_{J,\,N})_{J,\,N}$, and write ``$\alpha$'' instead of 
``$J,N$'' for the objects of $K$. We denote by 
$o\in K$ the initial object (corresponding to $N=0$ and $J=(0)$).
Set $\mathcal{X}_{\alpha}=\Spf(A_{\alpha})$, 
$\mathcal{Y}_{\alpha}=
\mathcal{Y}\times_{\mathcal{X}}\mathcal{X}_{\alpha}$,
$X_{\alpha}=\mathcal{X}_{\alpha}^{\rig}$ and 
$Y_{\alpha}=\mathcal{Y}^{\alpha}_{\rig}$.
For $\alpha\to \beta$ in $K$, 
we have Cartesian squares of rigid analytic spaces
$$\xymatrix{Y_{\beta} \ar[r]^-{g'_{\beta\alpha}} \ar[d]^-{f_{\beta}} 
&  Y_{\alpha} \ar[d]^-{f_{\alpha}} \\
X_{\beta} \ar[r]^-{g_{\beta\alpha}} & X_{\alpha}}$$
where the horizontal arrows are open immersions.
(Note that $f_o=f$.) We deduce commutative 
squares of $\infty$-categories
\begin{equation}
\label{eq-thm:general-base-change-thm-1}
\begin{split}
\xymatrix{\RigSH^{(\eff,\,\hyp)}_{\tau}(X_{\alpha};\Lambda) 
\ar[r]^-{f_{\alpha}^*} 
\ar[d]^-{g_{\beta\alpha}^*} & 
\RigSH^{(\eff,\,\hyp)}_{\tau}(Y_{\alpha};\Lambda) 
\ar[d]^-{g'^*_{\beta\alpha}}\\
\RigSH^{(\eff,\,\hyp)}_{\tau}(X_{\beta};\Lambda) \ar[r]^-{f_{\beta}^*} & 
\RigSH^{(\eff,\,\hyp)}_{\tau}(Y_{\beta};\Lambda).\!}
\end{split}
\end{equation}
In fact, we have a functor 
$K \to \Fun(\Delta^1,\Prl_{\omega})$
sending $\alpha\in K$ to $f_{\alpha}^*$ and 
$\alpha\to \beta$ to the commutative square 
\eqref{eq-thm:general-base-change-thm-1}.
Moreover, since the squares 
\eqref{eq-thm:general-base-change-thm-1}
are right adjointable by Proposition \ref{prop:6f1}(3), 
this functor factors through
the sub-$\infty$-category 
$$\Fun^{\rm RAd}(\Delta^1,\Prl_{\omega})=
\Fun(\Delta^1,\Prl_{\omega})\cap 
\Fun^{\rm RAd}(\Delta^1,\CAT_{\infty}),$$
where $\Fun^{\rm RAd}(\Delta^1,\CAT_{\infty})$ is
the $\infty$-category introduced in 
\cite[Definition 4.7.4.16]{lurie:higher-algebra}.

Consider a colimit diagram 
$K^{\rhd}\to \Fun(\Delta^1,\Prl_{\omega})$
extending the one described above. Since all the 
$\infty$-categories we are considering are stable, 
Lemma
\ref{lem:colimits-of-right-adjointable}
below implies that 
this diagram factors also through the sub-$\infty$-category
$\Fun^{\rm RAd}(\Delta^1,\Prl_{\omega})$. Evaluating the functor 
$K^{\rhd}\to \Fun(\Delta^1,\Prl_{\omega})$
at the edge $o\to \infty$, where $\infty\in K^{\rhd}$ is the 
cone point, we obtain a commutative square in $\Prl_{\omega}$
$$\xymatrix@C=3pc{\RigSH^{(\eff,\,\hyp)}_{\tau}(X_o;\Lambda) 
\ar[r]^-{f_o^*} 
\ar[d]^-{\colim_{\alpha}\, g_{\alpha o}^*} & 
\RigSH^{(\eff,\,\hyp)}_{\tau}(Y_o;\Lambda) 
\ar[d]^-{\colim_{\alpha}\,g'^*_{\alpha o}}\\
\underset{\alpha}{\colim}\,
\RigSH^{(\eff,\,\hyp)}_{\tau}(X_{\alpha};\Lambda) 
\ar[r]^-{\colim_{\alpha}\,f_{\alpha}^*} & 
\underset{\alpha}{\colim}\,
\RigSH^{(\eff,\,\hyp)}_{\tau}(Y_{\alpha};\Lambda)}$$
which is right adjointable. 
By Theorem \ref{thm:anstC}, this square 
is equivalent to the one in the statement.

\paragraph*{Part 2}
\noindent
We now assume that 
$g$ is not necessarily of finite type. 
We may assume that $g$ is induced by a morphism 
$\Spf(A') \to \Spf(A)$ of affine formal schemes. Set 
$\mathcal{X}=\Spf(A)$ and $\mathcal{X}'=\Spf(A')$. 
Let $\mathcal{Y}$ be a quasi-compact and quasi-separated 
formal $\mathcal{X}$-scheme such that $Y=\mathcal{Y}^{\rig}$, and 
let $\mathcal{Y}'=\mathcal{Y}\times_{\mathcal{X}}\mathcal{X}'$ so that
$\mathcal{Y}'^{\rig}=Y'$.
Write $A'$ as a filtered colimit $A'=\colim_{\alpha}\, A_{\alpha}$
of finitely generated adic $A$-algebras $A_{\alpha}$. Set also
$\mathcal{X}_{\alpha}=\Spf(A_{\alpha})$, $\mathcal{Y}_{\alpha}=
\mathcal{Y}\times_{\mathcal{X}}\mathcal{X}_{\alpha}$, $X_{\alpha}=\mathcal{X}_{\alpha}^{\rig}$ and 
$Y_{\alpha}=\mathcal{Y}_{\alpha}^{\rig}$. 
If $\tau$ is the \'etale topology and $\Lambda$ is not 
eventually coconnective, we may assume that the numbers 
$\pvcd_{\Lambda}(X_{\alpha})$ are  
bounded independently of $\alpha$.

As in the first part of the proof, we have 
a diagram $K\to \Fun(\Delta^1,\Prl_{\omega})$
sending $\alpha\to \beta$ to squares of the form 
\eqref{eq-thm:general-base-change-thm-1}.
Since the morphisms 
$g_{\beta\alpha}:X_{\beta}\to X_{\alpha}$ are of finite type, 
these squares are right adjointable as shown in the 
first part of the proof. 
The result follows again by considering a colimit
diagram $K^{\rhd}\to \Fun(\Delta^1,\Prl_{\omega})$, and 
using Lemma \ref{lem:colimits-of-right-adjointable}
and Theorem \ref{thm:anstC}.
\end{proof}

The following lemma, which was used in the proof of 
Theorem \ref{thm:general-base-change-thm}, is well-known.
We include a proof for completeness. 
(Recall that we are using the notation \sym{$\Fun^{\rm RAd}$} 
following \cite[Definition 4.7.4.16]{lurie:higher-algebra}.)

\begin{lemma}
\label{lem:colimits-of-right-adjointable}
Let $K$ be a simplicial set. Let 
$\overline{\mathcal{C}}:
K^{\rhd}\to \Fun(\Delta^1,\Prl)$
be a colimit diagram and let $\mathcal{C}$ be its restriction 
to $K$. Assume the following conditions:
\begin{enumerate}

\item[(1)] $\mathcal{C}$ factors through 
$\Fun^{\rm RAd}(\Delta^1,\Prl)=\Fun(\Delta^1,\Prl)\cap 
\Fun^{\rm RAd}(\Delta^1,\CAT_{\infty})$;

\item[(2)] for every $s\in K$, the right adjoint to the functor 
$f_s:\mathcal{C}_0(s) \to \mathcal{C}_1(s)$, associated to 
$s$ by $\mathcal{C}$, is colimit-preserving.

\end{enumerate}
(Note that the second condition is satisfied if 
$f_s$ is compact-preserving, and the $\infty$-categories
$\mathcal{C}_0(s)$ and $\mathcal{C}_1(s)$ are stable
and compactly generated.)
Then, $\overline{\mathcal{C}}$ also factors through 
$\Fun^{\rm RAd}(\Delta^1,\Prl)$. Moreover, the resulting 
map $K^{\rhd} \to \Fun^{\rm RAd}(\Delta^1,\Prl)$
is a colimit diagram.
\end{lemma}

\begin{proof}
Using the equivalence $\Prl\simeq (\Prr)^{\op}$,
we deduce a limit diagram 
$$\overline{\mathcal{C}}{}':(K^{\op}){^\lhd}\to 
\Fun(\Delta^{1,\,\op},\Prr).$$
We denote by $\mathcal{C}'$ the restriction of $\overline{\mathcal{C}}{}'$
to $K^{\op}$.
Applying $\overline{\mathcal{C}}$ and 
$\overline{\mathcal{C}}{}'$ to an edge $e:s\to t$ in $K^{\rhd}$, 
we get the following commutative squares of 
$\infty$-categories
$$\xymatrix{\,\ar@{}[d]|(.6){\displaystyle 
\overline{\mathcal{C}}(e):} \\ \,}
\xymatrix{\mathcal{C}_0(s) \ar[d]\ar[r]^-{f_s} & 
\mathcal{C}_1(s)\ar[d]\\
\mathcal{C}_0(t) \ar[r]^-{f_t} & \mathcal{C}_1(t)}
\qquad 
\xymatrix{\,\ar@{}[d]|(.6){\displaystyle \text{and}} \\ \,}
\qquad
\xymatrix{\,\ar@{}[d]|(.6){\displaystyle 
\overline{\mathcal{C}}{}'(e):} \\ \,}
\xymatrix{\mathcal{C}_0(s) & 
\mathcal{C}_1(s) \ar[l]_-{g_s} \\
\mathcal{C}_0(t) \ar[u]
& \mathcal{C}_1(t),\! 
\ar[l]_-{g_t} \ar[u]}$$
where the functors in the second square are the 
right adjoints to the functors in the first square.
By condition (2) the functors $g_s$ admit right adjoints.
Moreover, the first square $\overline{\mathcal{C}}(e)$ 
is right adjointable 
if and only if the square $\overline{\mathcal{C}}{}'(e)$
is right adjointable. 
We can then reformulate the problem as follows: if
$\mathcal{C}'$ factors through 
$$\Fun^{\rm RAd}(\Delta^{1,\,\op},\Prr)=\Fun(\Delta^{1,\,\op},\Prr)\cap
\Fun^{\rm RAd}(\Delta^{1,\,\op},\CAT_{\infty}),$$ 
then the same holds true for $\overline{\mathcal{C}}{}'$
and the resulting map is a limit diagram.
Since limits in $\Prr$ are computed in $\CAT_{\infty}$
(by \cite[Theorem 5.5.3.18]{lurie}), this follows from 
\cite[Corollary 4.7.4.18(2)]{lurie:higher-algebra}.
\end{proof}

\begin{rmk}
\label{rmk-thm:general-base-change-thm}
Keep the notations and assumptions of Theorem
\ref{thm:general-base-change-thm}. The commutative square
$$\xymatrix{\Shv^{(\hyp)}_{\tau}(X;\Lambda) \ar[r]^-{f^*} 
\ar[d]^-{g^*} & \Shv^{(\hyp)}_{\tau}(Y;\Lambda) \ar[d]^-{g'^*}\\
\Shv^{(\hyp)}_{\tau}(X';\Lambda) \ar[r]^-{f'^*} & 
\Shv^{(\hyp)}_{\tau}(Y';\Lambda) }$$
is also right adjointable.
This is proven by the same method: instead of using
Theorem \ref{thm:anstC}, we use the much easier 
Corollary \ref{cor:colim-small-et-site-2}.
There is also an unstable version of this result, asserting that 
$$\xymatrix{\Shv^{(\hyp)}_{\tau}(X) \ar[r]^-{f^*} 
\ar[d]^-{g^*} & \Shv^{(\hyp)}_{\tau}(Y) \ar[d]^-{g'^*}\\
\Shv^{(\hyp)}_{\tau}(X') \ar[r]^-{f'^*} & 
\Shv^{(\hyp)}_{\tau}(Y') }$$
is right adjointable under some assumptions. This 
holds for instance when $\tau$ is the Nisnevich topology, and
$X$, $X'$, $Y$ and $Y'$ locally of finite Krull dimension.
When $\tau$ is the \'etale topology, we have a weaker result:
under the same assumption on the Krull dimensions, 
the base change morphism 
$g^*\circ f_*\to f'_* \circ g'^*$ is an isomorphism 
when evaluated at truncated \'etale sheaves and, in particular, 
at \'etale sheaves of sets. A proof of this can be 
obtained by adapting the proof of Theorem 
\ref{thm:general-base-change-thm}. 
Indeed, Corollary \ref{cor:colim-small-et-site-2}
is still true for the $\infty$-categories of $n$-truncated 
$\mathcal{S}$-valued sheaves $\Shv_{\tau}(-)_{\leq n}$. (In this 
case, there is no distinction between sheaves and hypersheaves.)
Similarly, if $h:T \to S$ is a quasi-compact morphism between
rigid analytic spaces 
locally of finite Krull dimension, the associated functor 
$h^*:\Shv_{\tau}(\Et/S)_{\leq n}\to \Shv_{\tau}(\Et/T)_{\leq n}$
belongs to $\Prl_{\omega}$.
\end{rmk}

\subsection{Stalks}

$\empty$

\smallskip

\label{sec:stalks}

In this subsection, we determine under some mild hypotheses 
the stalks of $\RigSH^{(\eff,\,\hyp)}_{\tau}(-;\Lambda)$, 
which is a $\tau$-(hyper)sheaf by Theorem \ref{thm:hyperdesc}. 
We then use this to generalise Theorem \ref{thm:anstC}.
We start with a general fact on presheaves with 
values in a compactly generated $\infty$-category.

\begin{prop}
\label{prop:stalk-equivalence}
Let $(\mathcal{C},\tau)$ be a site having enough points and 
let $\mathcal{V}$ be a compactly generated $\infty$-category.
For a morphism
$f:\mathcal{F}\to\mathcal{G}$ 
in $\PSh(\mathcal{C};\mathcal{V})$, 
the following conditions are equivalent:
\begin{enumerate}

\item[(1)] $\Lder_{\tau}(f):\Lder_{\tau}(\mathcal{F})
\to \Lder_{\tau}(\mathcal{G})$ is an equivalence in 
$\Shv^{\hyp}_{\tau}(\mathcal{C};\mathcal{V})$;

\item[(2)] $f_x:\mathcal{F}_x\to\mathcal{G}_x$ is an equivalence
in $\mathcal{V}$ for all $x$ in a conservative family 
of points of $(\mathcal{C},\tau)$.

\end{enumerate}
\end{prop}

\begin{proof}
By \cite[Proposition~2.5]{drew:motivic-hodge}, 
condition (1) holds if and only if, for all compact objects 
$A\in \mathcal{V}$, the maps of presheaves of spaces
$$\Map_{\mathcal{V}}(A,f):
\Map_{\mathcal{V}}(A,\mathcal{F}) \to 
\Map_{\mathcal{V}}(A,\mathcal{G})$$
induce equivalences after $\tau$-hypersheafification. 
This is the case if and only if for every $x$ as in (2), 
the induced maps on stalks 
$$\Map_{\mathcal{V}}(A,f)_x:
\Map_{\mathcal{V}}(A,\mathcal{F})_x \to 
\Map_{\mathcal{V}}(A,\mathcal{G})_x.$$
are equivalences. Since the $A$'s are compact and stalks are computed 
by filtered colimits, the above maps are equivalent to 
$$\Map_{\mathcal{V}}(A,f_x):
\Map_{\mathcal{V}}(A,\mathcal{F}_x) \to 
\Map_{\mathcal{V}}(A,\mathcal{G}_x).$$
Since $\mathcal{V}$ is compactly generated and 
$A$ varies among all compact objects, our condition 
is equivalent to asking that the maps
$f_x:\mathcal{F}_x\to \mathcal{G}_x$ are equivalences as needed.
\end{proof}

Later we use Proposition
\ref{prop:stalk-equivalence} with $\mathcal{V}=\Prl_{\omega}$. 
This is indeed possible by Proposition 
\ref{prop:Prlomega-comp-gen}
below, whose proof relies on two technical lemmas.
The first one is a variant of 
the characterisation of presentability given in 
\cite[Theorem 5.5.1.1(6)]{lurie}
which is certainly well-known. We provide an argument because 
we couldn't find a reference.

\begin{lemma}
\label{lem:characteri-presentable-infty-cat}
Let $\mathcal{C}$ be a locally small $\infty$-category 
admitting small colimits. Assume that there exists a 
regular cardinal $\kappa$ and a set $S\subset \mathcal{C}$ 
of $\kappa$-compact objects such that $\mathcal{C}$ 
coincides with its smallest full sub-$\infty$-category 
containing $S$ and stable under colimits.
Then $\mathcal{C}$ is $\kappa$-compactly generated
(in the sense of \cite[Definition 5.5.7.1]{lurie}).
\end{lemma}

\begin{proof}
The difference with \cite[Theorem 5.5.1.1(6)]{lurie}
is that we do not assume that every object of $\mathcal{C}$
is a colimit of a diagram with values in the full sub-$\infty$-category
spanned by $S$.

Let $\mathcal{E}\subset \mathcal{C}$
be the smallest sub-$\infty$-category of $\mathcal{C}$ 
containing $S$ and stable under $\kappa$-small colimits.
The $\infty$-category $\mathcal{E}$ can be constructed from $S$
by transfinite induction as follows. Let $\mathcal{E}_0$ 
be the full sub-$\infty$-category of $\mathcal{C}$ spanned
by $S$ and, for an ordinal $\nu>0$, let 
$\mathcal{E}_{\nu}$ be the full sub-$\infty$-category of 
$\mathcal{C}$ spanned by colimits of $\kappa$-small
diagrams in $\bigcup_{\mu<\nu}\mathcal{E}_{\mu}$.
Then $\mathcal{E}=\bigcup_{\nu<\kappa}\mathcal{E}_{\nu}$. 
This shows that $\mathcal{E}$ is essentially small 
and that every object of $\mathcal{E}$ is $\kappa$-compact
(by \cite[Corollary 5.3.4.15]{lurie}). 
By \cite[Proposition 5.3.5.11]{lurie}, the inclusion 
$\mathcal{E}\to \mathcal{C}$ extends uniquely to a functor 
$\phi:\Ind_{\kappa}(\mathcal{E})\to \mathcal{C}$
preserving $\kappa$-filtered colimits, and this functor 
is fully faithful.
In fact, by \cite[Proposition~5.3.6.2 and Example~5.3.6.8]{lurie},
$\Ind_{\kappa}(\mathcal{E})$ admits small colimits and 
the functor $\phi$ is colimit-preserving. Using that the essential image 
of $\phi$ contains $S$, we deduce that $\phi$ is an 
equivalence of $\infty$-categories. Since
$\Ind_{\kappa}(\mathcal{E})$ is presentable by 
\cite[Theorem 5.5.1.1]{lurie}, this 
finishes the proof.
(Note that $\Ind_{\kappa}(\mathcal{E})$ is 
$\kappa$-accessible by definition, see
\cite[Definition 5.4.2.1]{lurie}.)
\end{proof}

\begin{lemma}
\label{lem:criterion-compact-gen}
Let $\mathcal{C}$ and $\mathcal{D}$ be $\infty$-categories
such that $\mathcal{C}$ is compactly generated and 
$\mathcal{D}$ admits small colimits. Assume that there is a 
functor $G:\mathcal{D} \to \mathcal{C}$ 
with the following properties:
\begin{enumerate}

\item[(1)] it admits a left adjoint;

\item[(2)] it is conservative;

\item[(3)] it commutes with filtered colimits.

\end{enumerate}
Then $\mathcal{D}$ is compactly generated. 
Moreover, if $F$ is a left adjoint to $G$, then $F$ takes a 
set of compact generators of $\mathcal{C}$ to a set of compact 
generators of $\mathcal{D}$.
\end{lemma}

\begin{proof}
Since $G$ commutes with filtered colimits, the functor 
$F$ takes a compact object of $\mathcal{C}$ to a compact object 
of $\mathcal{D}$. Let $\mathcal{C}_0$ be the 
full sub-$\infty$-category of $\mathcal{C}$ spanned by 
compact objects, and let $\mathcal{D}'\subset \mathcal{D}$
be the smallest sub-$\infty$-category containing $F(\mathcal{C}_0)$
and stable under colimits. By Lemma 
\ref{lem:characteri-presentable-infty-cat}, 
$\mathcal{D}'$ is compactly generated
since $\mathcal{C}_0$ is essentially small. 
Thus, it suffices to show that the inclusion functor 
$U:\mathcal{D}'\to \mathcal{D}$ is an equivalence. 
By \cite[Corollary 5.5.2.9 \& Remark 5.5.2.10]{lurie}, 
the functor $U$ admits a right adjoint $V$
and it is enough to show that 
$V$ is conservative. This follows from the hypothesis that 
$G$ is conservative. Indeed, we have $G\simeq G'\circ V$
where $G'$ is right adjoint to the functor 
$F':\mathcal{C}\to \mathcal{D}'$ induced by $F$
(which exists by \cite[Corollary 5.5.2.9]{lurie}).
\end{proof}

\begin{prop}
\label{prop:Prlomega-comp-gen}
The $\infty$-category $\Prl_{\omega}$ is compactly generated.
\end{prop}

\begin{proof}
This is probably well-known, but we couldn't find a reference.
We include a proof here for completeness. 
Denote by $\Cat_{\infty}^{\rm rex,\, idem}$ the 
sub-$\infty$-category of $\Cat_{\infty}$ whose objects are the
idempotent complete small $\infty$-categories
admitting finite colimits and whose morphisms are
the right exact functors. 
By \cite[Lemma 5.3.2.9(1)]{lurie:higher-algebra}, the functor 
$\mathcal{C}\mapsto \Ind_{\omega}(\mathcal{C})$ induces an 
equivalence of $\infty$-categories between 
$\Cat_{\infty}^{\rm rex,\,idem}$ and $\Prl_{\omega}$.
Thus, it is enough to show that $\Cat_{\infty}^{\rm rex,\,idem}$
is compactly generated.
Since $\Prl_{\omega}$ admits small colimits by 
\cite[Proposition 5.5.7.6]{lurie}, 
the same is true for $\Cat_{\infty}^{\rm rex,\,idem}$
which is moreover obviously locally small.

We will show that 
$\Cat_{\infty}^{\rm rex,\,idem}$ 
is compactly generated by applying 
Lemma \ref{lem:criterion-compact-gen}
to the inclusion functor 
$\Cat_{\infty}^{\rm rex,\,idem} \to \Cat_{\infty}$.
First, note that $\Cat_{\infty}$ is compactly generated. Indeed, 
$\Cat_{\infty}$ is the $\infty$-category associated to the 
combinatorial simplicial model category $\Set_{\Delta}^+$ of 
marked simplicial sets where the cofibrations are generated by 
monomorphisms with compact domain and codomain, and where 
fibrant objects are stable by filtered colimits. 
(See \cite[Propositions 3.1.3.7 \& 3.1.4.1, \& Theorem 3.1.5.1]{lurie}.)
We now check that the inclusion functor  
$\Cat_{\infty}^{\rm rex,\,idem} \to \Cat_{\infty}$ 
satisfies properties (1)--(3) of 
Lemma \ref{lem:criterion-compact-gen}. 
Property (1) follows from  \cite[Corollary 5.3.6.10]{lurie}.
Property (2) is obvious: an inverse of a right exact equivalence 
of $\infty$-categories is right exact. For property (3), we need to show 
the following: given a filtered diagram in 
$\Cat_{\infty}^{\rm rex,\,idem}$,
its colimit computed in $\Cat_{\infty}$
admits finite colimits and is idempotent complete. 
The first property follows from \cite[Proposition 5.5.7.11]{lurie}. 
The second property follows from 
\cite[Corollary 4.4.5.21]{lurie}.\footnote{Corollary 4.4.5.21 
can be found in the electronic version of \cite{lurie}
on the author's webpage, but not in the published version.}
\end{proof}

We record the following lemma for later use.

\begin{lemma}
\label{lem:sheaf-conserv-point}
Let $(\mathcal{C},\tau)$ be a site and let 
$\mathcal{F}:\mathcal{C}^{\op}\to \CAT_{\infty}$ 
be a presheaf on $\mathcal{C}$. Set 
$\mathcal{E}=\lim_{\mathcal{C}^{\op}}\mathcal{F}$. 
(If $\mathcal{C}$ admits a final object $\star$, 
then $\mathcal{E}\simeq \mathcal{F}(\star)$.)
Given an object $X\in \mathcal{C}$, we denote by 
$A\mapsto A_X$ the obvious functor $\mathcal{E} \to \mathcal{F}(X)$.
\begin{enumerate}

\item[(1)] Assume that $\mathcal{F}$ is a $\tau$-(hyper)sheaf.
Then, for $A,B\in \mathcal{E}$, the presheaf on 
$\mathcal{C}$, given informally by 
$X \mapsto \Map_{\mathcal{F}(X)}(A_X,B_X)$, 
is a $\tau$-(hyper)sheaf.

\item[(2)] Assume that $\mathcal{F}$ is a $\tau$-hypersheaf
and that the limit diagram $(\mathcal{C}^{\rhd})^{\op}\to 
\CAT_{\infty}$ extending $\mathcal{F}$ 
factors through $\Prl_{\omega}$.
Assume also that $(\mathcal{C},\tau)$ 
admits a conservative family of points $(x_i)_i$.
Then, the family of functors $(\mathcal{E} \to 
\mathcal{F}_{x_i})_i$, where the stalks $\mathcal{F}_{x_i}$ 
are computed in $\Prl_{\omega}$, is conservative. 

\end{enumerate}
\end{lemma}

\begin{proof}
We denote by 
$M:(\CAT_{\infty})_{\partial \Delta^1/}\to \mathcal{S}$
the copresheaf corepresented by 
$\partial\Delta^1\to \Delta^1$. The functor $M$ commutes with 
limits and admits the following informal description.
It sends an $\infty$-category 
$\mathcal{Q}$ together with a functor 
$q:\partial \Delta^1\to \mathcal{Q}$ to the mapping space
$\Map_{\mathcal{Q}}(q(0),q(1))$. 
This is indeed a consequence of 
\cite[Proposition 1.2]{Dugger-Spivak-Mapping}.

To give a precise construction of the 
presheaf described informally in (1), we consider 
$\mathcal{E}$ as an object of 
$(\CAT_{\infty})_{\partial\Delta^1/}$ using the functor
$e:\partial\Delta^1\to \mathcal{E}$ mapping $0$ to $A$ and 
$1$ to $B$. By the definition of 
$\mathcal{E}$, the presheaf $\mathcal{F}$ lifts to a 
$(\CAT_{\infty})_{\mathcal{E}/}$-valued presheaf $\mathcal{F}'$. 
The functor $e$ gives rise to a functor 
$$(\CAT_{\infty})_{\mathcal{E}/}
\to (\CAT_{\infty})_{\partial\Delta^1/}$$
and we denote by $\mathcal{F}''$ 
the $(\CAT_{\infty})_{\partial\Delta^1/}$-valued presheaf
obtained from $\mathcal{F}'$ by composing with this functor.
By construction, $\mathcal{F}''$ is a lift of $\mathcal{F}$
admitting the following informal description. 
It sends an object $X\in \mathcal{C}$ to the 
$\infty$-category $\mathcal{F}(X)$ together with the functor 
$\partial\Delta^1
\to \mathcal{F}(X)$ mapping $0$ to $A_X$ and $1$ to $B_X$. 
The presheaf $X\mapsto \Map_{\mathcal{F}(X)}(A_X,B_X)$ in (1)
is then defined to be $M\circ \mathcal{F}''$.
That said, the conclusion of assertion (1) is now clear. 
Indeed, the projection 
$(\CAT_{\infty})_{\partial\Delta^1/}\to \CAT_{\infty}$ preserves and 
detects limits by \cite[Proposition 1.2.13.8]{lurie} and,
as mentioned above, the functor $M$ is limit-preserving. 
Thus, the conclusion follows from Remark
\ref{recollection-hypersheaves}(1).

Given a point $x$ of $(\mathcal{C},\tau)$, we denote by 
$A\mapsto A_x$ the functor $\mathcal{E}\to \mathcal{F}_x$.
To prove the second assertion, we fix a morphism 
$f:A \to B$ in $\mathcal{E}$ inducing equivalences 
$A_{x_i}\simeq B_{x_i}$ for all $i$. 
We need to prove that $f$ is an equivalence. Since $\mathcal{E}$
is compactly generated, it is enough to show that 
$f$ induces an equivalence 
$\Map_{\mathcal{E}}(C,A) \to \Map_{\mathcal{E}}(C,B)$
for every compact object $C\in \mathcal{E}$. 
The compositions with the $f_X$'s, for $X\in \mathcal{C}$, 
induce a morphism of presheaves
\begin{equation}
\label{eq-lem:sheaf-conserv-point}
(X\mapsto \Map_{\mathcal{F}(X)}(C_X,A_X))
\to (X\mapsto \Map_{\mathcal{F}(X)}(C_X,B_X)),
\end{equation}
whose construction we leave to the reader.
By assertion (1),  this is actually a 
morphism of $\tau$-hypersheaves.
Thus, to conclude, it is enough to show that the morphism
\eqref{eq-lem:sheaf-conserv-point} induces
equivalences on stalks at $x_i$ for every $i$. 
Since $C$ is compact, the stalk at $x_i$ of this morphism
is given by the map $\Map_{\mathcal{F}_{x_i}}(C_{x_i},A_{x_i})
\to \Map_{\mathcal{F}_{x_i}}(C_{x_i},B_{x_i})$
which is indeed an equivalence since 
$A_{x_i}\simeq B_{x_i}$.
\end{proof}

By Theorem \ref{thm:hyperdesc},
the $\Prl$-valued presheaf $\RigSH^{(\eff),\,\hyp}_{\tau}(-;\Lambda)$
has $\tau$-hyperdescent. Therefore, it is particularly 
useful to determine its stalks.
The next theorem shows that, under some mild hypotheses, 
these stalks can also be understood as 
$\infty$-categories of rigid analytic motives
over rigid points (in the sense of 
Definition \ref{dfn:point-like-rigid-analytic-space}).

\begin{thm}
\label{thm:etst}
Let $S$ be a rigid analytic space and let 
$\overline{s} \to S$ be an algebraic rigid point of $S$. 
(See Remark \ref{rmk:tau-geom-pt-of-S}.) Let
$\tau\in \{\Nis,\et\}$, and assume one of the following two alternatives.
\begin{enumerate}

\item[(1)] We work in the non-hypercomplete case.

\item[(2)] We work in the hypercomplete case and  
$S$ is $(\Lambda,\tau)$-admissible.

\end{enumerate}
Then there is an equivalence of $\infty$-categories
$$\RigSH_{\tau}^{(\eff,\,\hyp)}(-;\Lambda)_{\overline{s}}\simeq 
\RigSH^{(\eff,\,\hyp)}_{\tau}(\overline{s};\Lambda),$$
where the left-hand side is the stalk of 
$\RigSH_{\tau}^{(\eff,\,\hyp)}(-;\Lambda)$ 
at $\overline{s}$, i.e., the colimit, taken in $\Prl$,
of the diagram $(\overline{s} \to U\to S)
\mapsto \RigSH_{\tau}^{(\eff,\,\hyp)}(U;\Lambda)$ with $U\in \Et/S$.
\end{thm}

\begin{proof}
We need to show that the obvious functor 
$$\underset{\overline{s} \to U \to S}{\colim}\,
\RigSH^{(\eff,\,\hyp)}_{\tau}(U;\Lambda)
\to \RigSH^{(\eff,\,\hyp)}_{\tau}(\overline{s};\Lambda)$$
is an equivalence. The question being local on $S$ around
the image of $\overline{s}$, we may assume that 
$S$ is quasi-compact and quasi-separated. In particular, 
$S$ admits a formal model $\mathcal{S}$. The functor 
$$(\Spf(\kappa^+(\overline{s}))\to \mathcal{U} \to \mathcal{S})
\mapsto (\overline{s} \to \mathcal{U}^{\rig} \to S),$$
with $\mathcal{U}$ affine and rig-\'etale over $\mathcal{S}$, 
is cofinal. Moreover, by Lemma
\ref{lem:limit-of-formal-schemes}, 
we have a canonical isomorphism of formal schemes
$$\Spf(\kappa^+(\overline{s}))\simeq 
\lim_{\Spf(\kappa^+(\overline{s}))\to \mathcal{U} \to \mathcal{S}}\,
\mathcal{U}.$$
The result follows now from Theorem 
\ref{thm:anstC}. Indeed, if $S$ is $(\Lambda,\tau)$-admissible 
then so are $\overline{s}$ and every \'etale rigid analytic 
$S$-space $U$. (For $\overline{s}$, use that the absolute Galois
group of $\kappa(\overline{s})$ is a closed subgroup
of the absolute Galois group of $\kappa(s)$; for $U$, use 
Corollary \ref{cor:univers-Lambda-admis}.)
Moreover, by the proof of Lemma \ref{lem:Lambda-coh-dim-val-ring}, 
we have the inequality 
$\pvcd_{\Lambda}(U)\leq \pvcd_{\Lambda}(S)$,
and, since $S$ is quasi-compact, 
the $(\Lambda,\tau)$-admissibility of $S$ implies that 
$\pvcd_{\Lambda}(S)$ is finite. 
\end{proof}

\begin{rmk}
\label{rmk:analytic-point-fiber-}
Theorem \ref{thm:etst}
applies in the case of a rigid point $s\to S$ 
associated to a point $s\in |S|$. In this case, the stalk 
$\RigSH_{\tau}^{(\eff,\,\tau)}(-;\Lambda)_s$ has a simpler description: 
it is the colimit, taken in $\Prl$,
of the diagram $U \mapsto \RigSH_{\tau}^{(\eff,\,\hyp)}(U;\Lambda)$,
where $U$ runs over the open neighbourhoods $U\subset S$ of $s$.
Indeed, every \'etale neighbourhood $s\to T \to S$ of $s$ in $S$ 
can be refined by an open neighbourhood. (This follows from 
Corollary \ref{cor:proj-limi-cal-E-A} and 
Lemma \ref{lem:limit-of-formal-schemes}(1).) Similarly, if
$\overline{s} \to S$ is a $\Nis$-geometric rigid point as in 
Construction \ref{cons:point-an-nis-et}(1), 
we may restrict in the description of the stalk in Theorem 
\ref{thm:etst} to those \'etale neighbourhoods $U$ 
admitting good reduction.
\end{rmk}

\begin{cor}
\label{cor:conserv-fiber-funct}
Let $S$ be a rigid analytic space. Assume one of the following 
two alternatives.
\begin{enumerate}

\item[(1)] We work in the non-hypercomplete case, and 
$S$ is locally of finite Krull dimension. 
When $\tau$ is the \'etale topology, we assume furthermore 
that $\Lambda$ is eventually coconnective

\item[(2)] We work in the hypercomplete case, and 
$S$ is $(\Lambda,\tau)$-admissible.

\end{enumerate}
Then, the functors
$$\RigSH^{(\eff,\,\hyp)}_{\tau}(S;\Lambda)
\to \RigSH^{(\eff,\,\hyp)}_{\tau}(s;\Lambda),$$
for $s\in S$, are jointly conservative. 
\end{cor}

\begin{proof}
Let $\Op/S$ denote the category of open subspaces of $S$
endowed with the analytic topology. By Theorem 
\ref{thm:hyperdesc},
$\RigSH^{(\eff,\,\hyp)}_{\tau}(-;\Lambda)$ is a hypersheaf on 
$\Op/S$. (In the non-hypercomplete case, we use  
\cite[Theorem 3.12]{cla-mat:etale-k-th}
and \cite[Corollary 7.2.1.12]{lurie}
which insure that a sheaf on $\Op/S$ is 
automatically a hypersheaf.)
Moreover, by Proposition
\ref{prop:compact-sheaves-qcqs}, 
this presheaf takes values in $\Prl_{\omega}$.
The result follows now from Lemma
\ref{lem:sheaf-conserv-point} 
and Theorem
\ref{thm:etst}. 
\end{proof}

\begin{rmk}
\label{rmk:conserv-points-sh}
The algebraic analogue of Corollary 
\ref{cor:conserv-fiber-funct} is also true:
given a scheme $S$ and assuming one of the alternatives 
of this corollary, the functors 
$\SH^{(\eff,\,\hyp)}_{\tau}(S;\Lambda)
\to \SH^{(\eff,\,\hyp)}_{\tau}(s;\Lambda)$, for $s\in |S|$,
are jointly conservative. This can be deduced from
Proposition \ref{prop:loc1} by arguing as in the proof of
\cite[Corollary 14]{framed-localis}.
\end{rmk}

Our next goal is to upgrade Theorem \ref{thm:anstC} 
to a motivic analogue of \cite[Proposition 2.4.4]{huber};
see Theorem \ref{thm:anstC-v2} below.
We first introduce, following
\cite[Definition 2.4.2 \& Remark 2.4.5]{huber},
a notion of weak limit in the category of rigid analytic spaces. 

\begin{dfn}
\label{dfn:limit-rigid-an-spc}
\ncn{rigid analytic spaces!weak limit}
Let $(S_{\alpha})_{\alpha}$ be a cofiltered inverse system of 
rigid analytic spaces, with quasi-compact and quasi-separated 
transition maps. 
Let $S$ be a rigid analytic space endowed
with a map of pro-objects 
$(f_{\alpha})_{\alpha}:S \to (S_{\alpha})_{\alpha}$, i.e., with an 
element $(f_{\alpha})_{\alpha} \in \lim_{\alpha}\Hom(S,S_{\alpha})$.
We say that $S$ is a weak limit of $(S_{\alpha})_{\alpha}$ 
and write $S\sim \lim_{\alpha}S_{\alpha}$ if the following 
two conditions are satisfied:
\symn{$\sim$}
\begin{enumerate}

\item[(1)] the map $|S| \to \lim_{\alpha}|S_{\alpha}|$ is a homeomorphism;

\item[(2)] for every $s\in |S|$ with images 
$s_{\alpha}\in |S_{\alpha}|$, the morphism
$$\underset{\alpha}{\colim}\,\kappa^+(s_{\alpha})\to 
\kappa^+(s),$$
where the colimit is taken in the category of adic rings, 
is an isomorphism.

\end{enumerate}
\end{dfn}

\begin{exm}
\label{exm:pro-formal-schemes-rig}
Let $(\mathcal{S}_{\alpha})_{\alpha}$ be a 
cofiltered inverse system of formal schemes with affine transition 
maps and let $\mathcal{S}=\lim_{\alpha}\mathcal{S}_{\alpha}$ be its limit.
Set $S=\mathcal{S}^{\rig}$ and $S_{\alpha}=\mathcal{S}_{\alpha}^{\rig}$.
Then $S$ is a weak limit of $(S_{\alpha})_{\alpha}$.
Indeed, condition (1) follows from 
commutation of limits with limits, see Notation
\ref{not:visualisation-}.
The point is that any admissible blowup of $\mathcal{S}$ 
can be obtained as the strict transform of $\mathcal{S}$ 
with respect to an admissible blowup of an $\mathcal{S}_{\alpha}$ 
for some $\alpha$. Condition (2) follows from 
Lemma \ref{lem:limit-of-formal-schemes}(1).
\end{exm}

\begin{exm}
\label{exm:limit-open-neigh}
Let $X$ be a rigid analytic space and 
$Z\subset X$ a closed subspace. Let $(U_{\alpha})_{\alpha}$ 
be an inverse system of open neighbourhoods of $Z$ in $X$ 
such that, locally at every point of $Z$, this inverse system 
is cofinal in the system of all neighbourhoods of $Z$ in $X$. 
(When $X$ is quasi-compact, this is equivalent to saying that 
$(U_{\alpha})_{\alpha}$ is cofinal in the system of all 
neighbourhoods of $Z$ in $X$.) Then, $Z$ is a weak limit of 
$(U_{\alpha})_{\alpha}$. Indeed, condition (2) is obvious and,
for condition (1), we need to show that 
$|Z|=\bigcap_{\alpha}|U_{\alpha}|$. This follows easily from the 
fact that $|X|$ is a valuative topological space
(in the sense of \cite[Chapter 0, Definition 2.3.1]{fujiwara-kato}) 
and that $|Z|\subset |X|$ is stable by generisation. 
\end{exm}

The following lemma can be compared with 
\cite[Remark 2.4.3(i)]{huber}. See also the proof of 
\cite[Proposition 7.16]{scholze}.

\begin{lemma}
\label{lem:on-limit-rigid-an-spc}
Keep the notation as in Definition 
\ref{dfn:limit-rigid-an-spc} and consider the following 
variants of conditions (1) and (2):
\begin{enumerate}

\item[(1$'$)] the $f_{\alpha}$'s are quasi-compact 
and quasi-separated, and 
the map $|S| \to \lim_{\alpha}|S_{\alpha}|$
is a bijection;

\item[(2$'$)] for every $s\in |S|$ with images 
$s_{\alpha}\in |S_{\alpha}|$, 
the induced morphism of fields
$$\underset{\alpha}{\colim}\,\kappa(s_{\alpha}) \to \kappa(s)$$
has dense image.

\end{enumerate}
Then, conditions (2) and (2$'$) are equivalent. Moreover, if condition 
(2) is satisfied, then conditions (1) and (1$'$) are equivalent.  
\end{lemma}

\begin{proof}
We identify $\kappa^+(s_{\alpha})$ with a subring of 
$\kappa^+(s)$ and $\kappa(s_{\alpha})$ with a subfield 
of $\kappa(s)$. We may assume that there is an element
$\pi\in \kappa^+(s)$ which belongs to 
all the $\kappa^+(s_{\alpha})$'s and generates an ideal of definition
in each one of them. If (2) is satisfied, then 
$\kappa^+(s)$ is the $\pi$-adic completion of 
$\bigcup_{\alpha}\kappa^+(s_{\alpha})$, which implies that 
$\bigcup_{\alpha}\kappa(s_{\alpha})$ is dense in $\kappa(s)$.
Conversely, if (2$'$) is satisfied, then 
$\kappa^+(s)$ is the Hausdorff completion of 
$\kappa^+(s)\cap \bigcup_{\alpha}\kappa(s_{\alpha})$.
Then condition (2) follows from the following equalities 
$\pi^n\kappa^+(s)\cap \bigcup_{\alpha}
\kappa(s_{\alpha})=\bigcup_{\alpha} \pi^n\kappa^+(s_{\alpha})$
which are easily checked using the valuation on $\kappa(s)$.

Clearly, (1) implies (1$'$). 
We next assume that (2) is satisfied, and show that 
(1$'$) implies (1). Using that the $f_{\alpha}$'s and the transition 
morphisms of the inverse system $(S_{\alpha})_{\alpha}$ 
are quasi-compact 
and quasi-separated, we may reduce to the case where 
$S$ and all the $S_{\alpha}$'s are quasi-compact and quasi-separated. 
By \cite[\href{https://stacks.math.columbia.edu/tag/09XU}
{Lemma 09XU}]{stacks-project},
it is then enough to show that the bijection 
$|S|\simeq \lim_{\alpha}|S_{\alpha}|$
detects generisations.
Given $s\in |S|$ with images $s_{\alpha}\in |S_{\alpha}|$, 
the generisations of $s$ are the points of 
$\Spf(\kappa^+(s))$ while the generisations of 
$(s_{\alpha})_{\alpha}$ are the points of 
$\lim_{\alpha}\Spf(\kappa^+(s_{\alpha}))$.
Thus, condition (2) implies that $f$ induces a bijection between
the generisations of $s$ and those of $(s_{\alpha})_{\alpha}$.
\end{proof}

The following can be compared with 
\cite[Remark 2.4.3(ii)]{huber} and
\cite[Proposition 7.16]{scholze}.

\begin{lemma}
\label{lem:limit-an-spc-base-change}
Let $(S_{\alpha})_{\alpha}$ be a cofiltered inverse system of 
rigid analytic spaces, with quasi-compact and quasi-separated 
transition maps, and admitting a weak limit $S$.
Let $X$ be a rigid analytic $S_{\alpha_0}$-space for some 
index $\alpha_0$. Then $X\times_{S_{\alpha_0}}S$ is a weak limit of 
$(X\times_{S_{\alpha_0}}S_{\alpha})_{\alpha\leq \alpha_0}$.
\end{lemma}

\begin{proof}
We reduce easily to the case where $S$, the $S_{\alpha}$'s and 
$X$ are quasi-compact and quasi-separated. We will check that 
condition (1$'$) of Lemma 
\ref{lem:on-limit-rigid-an-spc}
and condition (2) of Definition
\ref{dfn:limit-rigid-an-spc}
are satisfied by the maps 
$X\times_{S_{\alpha_0}}S \to X\times_{S_{\alpha_0}}S_{\alpha}$, for 
$\alpha\leq \alpha_0$.
A point of $|X\times_{S_{\alpha_0}}S|$ corresponds to a 
point $s\in |S|$ and a point of $|X\times_{S_{\alpha_0}}s|$
mapping to the closed point of $|s|$. Using a similar description 
for the points of the $|X\times_{S_{\alpha_0}}S_{\alpha}|$'s, 
condition (1$'$) and (2) follow from the following assertion: given
$s\in |S|$ with images $s_{\alpha}\in |S_{\alpha}|$, 
$X\times_{S_{\alpha_0}}s$ is a weak limit of 
$(X\times_{S_{\alpha_0}}s)_{\alpha\leq \alpha_0}$.
To prove this assertion, choose a formal model 
$\mathcal{X} \to \mathcal{S}_{\alpha_0}$ 
of $X\to S_{\alpha_0}$ and use Example 
\ref{exm:pro-formal-schemes-rig} and
the isomorphism of formal schemes 
$\mathcal{X}\times_{\mathcal{S}_{\alpha_0}}\Spf(\kappa^+(s))
\simeq \lim_{\alpha\leq \alpha_0} 
\mathcal{X}\times_{\mathcal{S}_{\alpha_0}}\Spf(\kappa^+(s_{\alpha}))$.
\end{proof}

\begin{thm}
\label{thm:anstC-v2}
\ncn{continuity}
Let $(S_{\alpha})_{\alpha}$ be a cofiltered inverse system of 
rigid analytic spaces, with quasi-compact and quasi-separated 
transition maps, and admitting a weak limit $S$.
Let $\tau\in \{\Nis,\et\}$, and assume one of the following two 
alternatives.
\begin{enumerate}

\item[(1)] We work in the non-hypercomplete case, and $S$ and the 
$S_{\alpha}$'s are locally of finite Krull dimension.
When $\tau$ is the \'etale topology, we assume furthermore 
that $\Lambda$ is eventually coconnective.

\item[(2)] We work in the hypercomplete case, and 
$S$ and the $S_{\alpha}$'s are $(\Lambda,\tau)$-admissible
(see Definition
\ref{dfn:lambda-tau-admiss}).
When $\tau$ is the \'etale topology, we assume furthermore 
that $\Lambda$ is eventually coconnective or that,
for every $s\in |S|$ with images $s_{\alpha}\in |S_{\alpha}|$, 
the $\Lambda$-cohomological dimensions of the residue fields
$\kappa(s_{\alpha})$ are bounded independently of $\alpha$.

\end{enumerate}
Then, the obvious functor 
\begin{equation}
\label{eq-thm:anstC-v2-1}
\underset{\alpha}{\colim}\,\RigSH^{(\eff,\,\hyp)}_{\tau}
(S_{\alpha};\Lambda) \to \RigSH^{(\eff,\,\hyp)}_{\tau}(S;\Lambda),
\end{equation}
where the colimit is taken in $\Prl$, is an equivalence.
\end{thm}

\begin{proof}
Let $U_{\alpha_0,\,\bullet} \to S_{\alpha_0}$ be 
a hypercover of $S_{\alpha_0}$ in the analytic topology 
with $U_{\alpha_0,\,n}$ a disjoint union of a family 
$(U_{\alpha_0,\,n,\,i})_{i\in I_n}$ of open subspaces 
of $S_{\alpha_0}$. Set $U_{\alpha,\,n,\,i}=
U_{\alpha_0,\,n,\,i}\times_{S_{\alpha_0}}S_{\alpha}$
and $U_{n,\,i}=U_{\alpha_0,\,n,\,i}\times_{S_{\alpha_0}}S$.
We have hypercovers $U_{\alpha,\,\bullet} \to S_{\alpha}$ and 
$U_{\bullet} \to S$ with $U_{\alpha,\,n}=\coprod_{i\in I_n}
U_{\alpha,\,n,\,i}$ and similarly for $U_n$.
By \cite[Proposition 4.7.4.19]{lurie:higher-algebra},
there is an equivalence of $\infty$-categories
\begin{equation}
\label{eq-thm:anstC-v2-lim-colim}
\underset{\alpha}{\colim}\, \lim_{[n]\in \mathbf{\Delta}}
\prod_{i\in I_n}
\RigSH^{(\eff,\,\hyp)}_{\tau}(U_{\alpha,\,n,\,i};\Lambda)
\simeq \lim_{[n]\in \mathbf{\Delta}}
\prod_{i\in I_n}
\underset{\alpha}{\colim}\,\RigSH^{(\eff,\,\hyp)}_{\tau}
(U_{\alpha,\,n,\,i};\Lambda).
\end{equation}
The right adjointability of the squares that is needed for 
\cite[Proposition 4.7.4.19]{lurie:higher-algebra}
holds by the base change theorem for open immersions, 
which is a special case of Proposition
\ref{prop:6f1}(3).
The presheaf $\RigSH^{(\eff,\,\hyp)}_{\tau}(-;\Lambda)$ 
admits descent for the hypercovers $U_{\bullet}\to S$ and
$U_{\alpha,\,\bullet}\to S_{\alpha}$ by
Theorem \ref{thm:hyperdesc}. (In the non-hypercomplete case, 
we use the assumption that $S$ and the $S_{\alpha}$'s 
have locally finite Krull dimension so that descent implies 
hyperdescent by \cite[Theorem 3.12]{cla-mat:etale-k-th}
and \cite[Corollary 7.2.1.12]{lurie}.) Therefore, the equivalence 
\eqref{eq-thm:anstC-v2-lim-colim}
shows that it is enough to prove the theorem for 
the inverse systems $(U_{\alpha,\,n,\,i})_{\alpha\leq \alpha_0}$.
In particular, we may assume that the 
$S_{\alpha}$'s are quasi-compact and quasi-separated.

Denote by $\Op^{\qcqs}/S$ the category of 
quasi-compact and quasi-separated open subspaces of $S$, 
and similarly for other rigid analytic spaces.
Given that $\Op^{\qcqs}/S=
\colim_{\alpha}\,\Op^{\qcqs}/S_{\alpha}$, 
there exists a $\Prl$-valued presheaf $\mathcal{R}$ on 
$\Op^{\qcqs}/S$ given by 
$$\mathcal{R}(U)=\underset{\alpha\geq \alpha_0}{\colim}\,
\RigSH^{(\eff,\,\hyp)}_{\tau}(U_{\alpha};\Lambda)$$
for any $U_{\alpha_0}\in \Op^{\qcqs}/S_{\alpha_0}$ such that
$U=U_{\alpha_0}\times_{S_{\alpha_0}}S$. (As usual, we set
$U_{\alpha}=U_{\alpha_0}\times_{S_{\alpha_0}}S_{\alpha}$.)
Moreover, we have a morphism of $\Prl$-valued presheaves
$$\phi:\mathcal{R} \to \RigSH^{(\eff,\,\hyp)}_{\tau}(-;\Lambda)$$
on $\Op^{\qcqs}/S$.
Since $S$ belongs to $\Op^{\qcqs}/S$, it suffices to 
show that $\phi$ is an equivalence of presheaves. 
We will achieve this by showing the following two properties:
\begin{enumerate}

\item[(1)] $\mathcal{R}$ and $\RigSH^{(\eff,\,\hyp)}_{\tau}(-;\Lambda)$ 
are hypersheaves on $\Op^{\qcqs}/S$ for the analytic topology;

\item[(2)] $\phi$ induces an equivalence on stalks for the analytic 
topology at every point $s\in |S|$.

\end{enumerate}
This suffices indeed by Propositions 
\ref{prop:stalk-equivalence} and 
\ref{prop:Prlomega-comp-gen}, since the presheaves 
$\mathcal{R}$ and $\RigSH^{(\eff,\,\hyp)}_{\tau}(-;\Lambda)$ 
on $\Op^{\qcqs}/S$ take values in 
$\Prl_{\omega}$ by Proposition 
\ref{prop:compact-shv-rigsm}.

First, we prove (1). That 
$\RigSH^{(\eff,\,\hyp)}_{\tau}(-;\Lambda)$ 
is a hypersheaf on $\Op^{\qcqs}/S$ was mentioned above.
To handle the case of $\mathcal{R}$, we use again 
\cite[Theorem 3.12]{cla-mat:etale-k-th}
and \cite[Corollary 7.2.1.12]{lurie}
which insure that a sheaf on $\Op^{\qcqs}/S$ is 
automatically a hypersheaf.
Thus, it is enough to show that $\mathcal{R}$ 
admits descent for truncated hypercovers 
$U_{\bullet}$ in $\Op^{\qcqs}/S$. 
We may assume that $U_{-1}=S$. Every such hypercover, is the 
inverse image of a truncated hypercover $U_{\alpha_0,\,\bullet}$
with $U_{\alpha_0,-1}=S_{\alpha_0}$. We may then use the equivalence
\eqref{eq-thm:anstC-v2-lim-colim} 
to conclude.

Next, we prove (2). Fix $s\in |S|$ with images 
$s_{\alpha}\in |S_{\alpha}|$. Since every 
quasi-compact and quasi-separated open neighbourhood of $s$ is 
the inverse image of a quasi-compact and quasi-separated open neighbourhood of $s_{\alpha}$, for $\alpha$ small enough, 
the functor $\phi_s$ can be rewritten as follows:
$$\underset{\alpha}{\colim}\,
\RigSH^{(\eff,\,\hyp)}_{\tau}(-;\Lambda)_{s_{\alpha}}
\to \RigSH^{(\eff,\,\hyp)}_{\tau}(-;\Lambda)_s.$$
Using Theorem \ref{thm:etst} (and Remark 
\ref{rmk:analytic-point-fiber-}), this functor is equivalent to
$$\underset{\alpha}{\colim}\,
\RigSH^{(\eff,\,\hyp)}_{\tau}(s_{\alpha};\Lambda)
\to \RigSH^{(\eff,\,\hyp)}_{\tau}(s;\Lambda).$$
By Theorem \ref{thm:anstC}, the latter is an equivalence.
\end{proof}

\subsection{(Semi-)separatedness}

$\empty$

\smallskip

\label{subsect:separatedness-}

In this subsection, we discuss two basic properties of the 
functor $\RigSH^{(\eff,\,\hyp)}_{\tau}(-;\Lambda)$, namely 
semi-separatedness and separatedness.

\begin{dfn}
\label{dfn:totally-insep-morph}
Let $e:X'\to X$ be a morphism of rigid analytic spaces.
\begin{enumerate}

\item[(1)] 
We say that $e$ is radicial if $|e|:|X'|\to |X|$ is injective and, 
for every $x'\in |X'|$ with image $x\in |X|$, the residue field
$\kappa(x')$ contains a dense purely inseparable extension of $\kappa(x)$.
\ncn{rigid analytic spaces!radicial morphism}

\item[(2)] We say that $e$ is a universal homeomorphism if 
it is quasi-compact, quasi-separated, surjective and radicial.
(See Remark \ref{rmk:on-faithfully-radicial} below.)
\ncn{rigid analytic spaces!universal homeomorphism}

\end{enumerate}
\end{dfn}

\begin{rmk}
\label{rmk:on-faithfully-radicial}
$\empty$

\begin{enumerate}

\item[(1)] Radicial morphisms and universal homeomorphisms 
are stable under base change. 

\item[(2)] If $e:X'\to X$ is a universal homeomorphism, 
then $|e|:|X'|\to |X|$ is a quasi-compact and quasi-separated 
bijection which detects generisation. By 
\cite[\href{https://stacks.math.columbia.edu/tag/09XU}
{Lemma 09XU}]{stacks-project}, this implies that
$|e|:|X'|\to |X|$ is a homeomorphism of topological spaces.
Moreover, by (1), this property is preserved by base change,
which explains our terminology.

\item[(3)] A morphism of schemes $e:X'\to X$ 
is called a universal homeomorphism
if every base change of $e$ induces a homeomorphism on the 
underlying topological spaces. By 
\cite[Chapitre IV, Corollaire 18.12.13]{EGAIV4}, this is equivalent 
to saying that $e$ is entire, surjective and radicial.

\end{enumerate}
\end{rmk}

\begin{lemma}
\label{lem:faithfully-radicial-etale-topos}
Let $e:X'\to X$ be a universal homeomorphism 
of rigid analytic spaces. The induced morphism 
$e:(\Et/X',\tau) \to (\Et/X,\tau)$
is an equivalence of sites,
{
i.e., induces an equivalence between the associated ordinary topoi,}
for $\tau\in \{\an,\Nis,\et\}$.
In particular, we have an equivalence of $\infty$-categories
$\Shv_{\tau}^{(\hyp)}(\Et/X';\Lambda)
\simeq \Shv_{\tau}^{(\hyp)}(\Et/X;\Lambda)$.
\end{lemma}

\begin{proof}
The second assertion follows from the first one using
Lemma \ref{lem:equi-of-sites-infty-topoi}. To prove the 
first assertion, we need to show that the unit 
$\id\to e_*e^*$ and counit $e^*e_*\to \id$ are equivalences
on $\tau$-sheaves of sets (i.e., on discrete $\tau$-sheaves). 
For $x\in |X|$, we have a morphism of 
sites $(\Et/x,\tau)\to (\Et/X,\tau)$, and we denote by 
$x^*$ the associated inverse image functor. Then, 
the functors $x^*$, for $x\in |X|$, are jointly conservative 
on $\tau$-sheaves of sets. The same discussion is equally valid 
for points of $X'$. Thus, we are left to show that 
the natural transformations $x^*\to x^*e_*e^*$ and
$x'^*e^*e_*\to x'^*$ are equivalences on $\tau$-sheaves of sets
for all $x\in |X|$ and $x'\in |X'|$. 
Assuming that $x$ is the image of $x'$, these natural transformations
are equivalent to 
$x^*\to e_{x,\,*}e_x^*x^*$ and 
$e_x^*e_{x,\,*}x'^*\to x'^*$, where $e_x:x' \to x$
is the obvious morphism. This follows from 
Remark \ref{rmk-thm:general-base-change-thm}
and the fact that the morphism $x'\to X'\times_X x$ identifies
$x'$ with $(X'\times_X x)_{\red}$.
Thus, we are reduced to prove the lemma for rigid 
points. Since $\kappa(x')$ contains a dense 
purely inseparable extension of $\kappa(x)$, 
the functor $\Et/x\to \Et/x'$ is an equivalence
of categories which respects the analytic, Nisnevich 
and \'etale topologies. 
\end{proof}

\begin{rmk}
\label{rmk:universal-home-topo}
Lemma \ref{lem:faithfully-radicial-etale-topos}
admits a variant for universal homeomorphisms of schemes which 
is well-known, see \cite[Expos\'e VIII, Th\'eor\`eme 1.1]{SGAIV2}.
\end{rmk}

\begin{cor}
\label{cor:cocartesian-square-}
Let $e:S'\to S$ be a universal homeomorphism of 
rigid analytic spaces. Then, for $\tau\in \{\Nis,\et\}$,
we have a coCartesian square in $\Prl$
$$\xymatrix{\RigSH^{(\eff)}_{\Nis}(S;\Lambda) \ar[r]^-{e^*} 
\ar[d] & 
\RigSH^{(\eff)}_{\Nis}(S';\Lambda) \ar[d]\\
\RigSH^{(\eff,\,\hyp)}_{\tau}(S;\Lambda) \ar[r]^-{e^*} 
& \RigSH^{(\eff,\,\hyp)}_{\tau}(S';\Lambda).}$$ 
Said differently, $\RigSH^{(\eff)}_{\Nis}(S';\Lambda)
\to \RigSH^{(\eff,\,\hyp)}_{\tau}(S';\Lambda)$ is a localisation 
functor with respect to the image by $e^*$ of 
morphisms of the form 
$\colim_{[n]\in \mathbf{\Delta}}\M^{(\eff)}(U_{\bullet})
\to \M^{(\eff)}(U_{-1})$, and their desuspensions and 
negative Tate twists when applicable, with $U_{\bullet}$ a 
$\tau$-hypercover in $\RigSm/S$ which we assume 
to be truncated in the non-hypercomplete case.
\end{cor}

\begin{proof}
Using Remark \ref{rmk:symmetric-obj-}, one reduces easily to 
the effective case. From the construction, one sees immediately that  
$\RigSH^{\eff}_{\Nis}(S';\Lambda)
\to \RigSH^{\eff,\,(\hyp)}_{\tau}(S';\Lambda)$
is the localisation functor with respect to morphisms of the form 
$\alpha'^*_n\mathcal{F}' \to 
\alpha'^*_n\mathcal{G}'$ where:
\begin{itemize}

\item $\alpha'_n:(\RigSm/S',\tau) \to (\Et/\B^n_{S'},\tau)$
is the premorphism of sites given by the obvious functor;

\item $\mathcal{F}'\to \mathcal{G}'$ is a morphism in 
$\Shv_{\Nis}(\Et/\B^n_{S'};\Lambda)$ inducing an equivalence in 
$\Shv_{\tau}^{(\hyp)}(\Et/\B^n_{S'};\Lambda)$.

\end{itemize}
For example, $\mathcal{F}'\to \mathcal{G}'$ could be 
$\colim_{[n]\in \mathbf{\Delta}}\,\Lambda_{\Nis}(U'_{\bullet})\to 
\Lambda_{\Nis}(U'_{-1})$ with $U'_{\bullet}$ a $\tau$-hypercover in 
$(\Et/\B^n_{S'},\tau)$ which is truncated in the 
non-hypercomplete case. The result follows now from 
the commutative square
$$\xymatrix{\Shv_{\Nis}(\Et/\B^n_S;\Lambda)
\ar[r]_-{\sim}^-{e^*} \ar[d]_-{\alpha_n^*} & 
\Shv_{\Nis}(\Et/\B^n_{S'};\Lambda) \ar[d]^-{\alpha'^*_n} \\
\Shv_{\Nis}(\RigSm/S;\Lambda) \ar[r]^-{e^*} & 
\Shv_{\Nis}(\RigSm/S';\Lambda)}$$
and Lemma \ref{lem:faithfully-radicial-etale-topos} 
which insures that the upper horizontal arrow
is an equivalence of $\infty$-categories respecting 
$\tau$-local equivalences 
(in both the hypercomplete and non-hypercomplete cases).
\end{proof}

\begin{thm}[Semi-separatedness]
\label{thm:seprig}
\ncn{semi-separatedness}
Let $\tau\in\{\Nis,\et\}$.
Let $e:X'\to X$ be a universal homeomorphism 
of rigid analytic spaces.
Assume that $X$ has locally finite Krull dimension. 
Assume also that every prime number is invertible in either  
$\mathcal{O}_X$ or $\pi_0\Lambda$. Then the functor 
$$e^*:\RigSH^{(\hyp)}_{\tau}(X;\Lambda)
\to \RigSH^{(\hyp)}_{\tau}(X';\Lambda)$$
is an equivalence of $\infty$-categories.
\end{thm}

\begin{proof}
By Corollary \ref{cor:cocartesian-square-}, 
we may assume that $\tau$ is the Nisnevich topology.
Since $X$ and $X'$ are locally of finite Krull dimension,
we are automatically working in the non-hypercomplete case
by Proposition \ref{prop:automatic-hypercomp-motives}. 
We need to show that the unit $\id \to e_*e^*$ and 
the counit $e^*e_*\to \id$ are equivalences.
By Corollary \ref{cor:conserv-fiber-funct}, 
it is enough to show that the natural transformations
$x^*\to x^*e_*e^*$ and $x'^*e^*e_* \to x'^*$
are equivalences for all points $x\in |X|$ and $x'\in |X'|$.
(Here, we denote by $x$ the morphism of rigid analytic spaces
$x\to X$ associated to the point $x\in |X|$, and similarly 
for $x'$.)
Assuming that $x$ is the image of $x'$, 
these natural transformations are equivalent to 
$x^*\to e_{x,\,*}e_x^*x^*$ and 
$e_x^*e_{x,\,*}x'^*\to x'^*$, where $e_x:x' \to x$
is the obvious morphism. This follows from Theorem
\ref{thm:general-base-change-thm}
and the fact that the morphism $x'\to X'\times_X x$ identifies
$x'$ with $(X'\times_X x)_{\red}$.
Thus, we are reduced to prove the result for the morphism 
$e_x:x'\to x$ of rigid points. 
Moreover, we can write $x'\sim \lim_{\alpha} x_{\alpha}$
with $(x_{\alpha})_{\alpha}$ the cofiltered inverse system of 
rigid analytic $x$-points such that $\kappa(x_{\alpha})$ 
is a finite purely inseparable extension of $\kappa(x)$ 
contained in $\kappa(x')$. Using Theorem \ref{thm:anstC-v2}, 
we reduce to show that $e^*$ is an equivalence 
for a morphism of rigid points 
$e:x'\to x$ such that $\kappa(x')/\kappa(x)$ is a finite
purely inseparable extension.

Arguing as in \cite[Sous-lemme 1.4]{ayoub-etale}, 
we see that $e^*e_*\simeq \id$. Thus, we only need to check 
that $\id \to e_*e^*$ is an equivalence. Since 
$e^*$ and $e_*$ commute with colimits (by Proposition
\ref{prop:compact-shv-rigsm}), it is enough
to show that $\id \to e_*e^*$ is an equivalence when 
applied to a set of compact generators. 
Such a set is given, up to 
desuspension and negative Tate twists, 
by objects of the form $f_{\sharp}\Lambda$ with 
$f:\Spf(A)^{\rig} \to x$ where $A$ a rig-smooth $\kappa^+(x)$-adic
algebra. Set $A'=A\,\widehat{\otimes}_{\kappa^+(x)}\kappa^+(x')$, 
and let $e':\Spf(A')^{\rig}\to \Spf(A)^{\rig}$ and $f':\Spf(A')^{\rig}\to x'$ be the obvious morphisms. 
Using Propositions \ref{prop:6f1} and 
\ref{prop:base-change-finite-and-projection}(2),
we have equivalences 
$e_*e^*f_{\sharp} \simeq e_*f'_{\sharp}e'^*\simeq f_{\sharp}e'_*e'^*$.
Thus, to finish the proof, we only need to show that 
$\Lambda \to e'_*e'^*\Lambda$ 
is an equivalence in $\RigSH_{\Nis}(\Spf(A)^{\rig};\Lambda)$.
Recall that there is a morphism of $\Prl$-valued presheaves 
$$\An^*:\SH_{\Nis}(-;\Lambda) \to \RigSH_{\Nis}((-)^{\an};\Lambda)$$
on $\Sch^{\lft}/U$, with $U=\Spec(A[\pi^{-1}])$ 
where $\pi\in \kappa^+(x)$ a generator of an ideal of definition.
Calling $e'':\Spec(A'[\pi^{-1}]) \to \Spec(A[\pi^{-1}])$
the obvious morphism, we have, by Proposition
\ref{prop:ran-star-com-f-lower-star-prelim},
equivalences $\An^* e''_*e''^*\simeq e'_*e'^*\An^*$.
Thus, it is enough to show that 
$\Lambda \to e''_*e''^*\Lambda$ 
is an equivalence in $\SH_{\Nis}(\Spec(A[\pi^{-1}]);\Lambda)$.
This follows from Theorem \ref{thm:sepalgebr} below.
\end{proof}

\begin{thm}
\label{thm:sepalgebr}
Let $\tau\in\{\Nis,\et\}$.
Let $e:X'\to X$ be a universal homeomorphism of schemes.
Assume that every prime number is invertible in either  
$\mathcal{O}_X$ or $\pi_0\Lambda$. Then the functor 
$$e^*:\SH^{(\hyp)}_{\tau}(X;\Lambda)
\to \SH^{(\hyp)}_{\tau}(X';\Lambda)$$
is an equivalence of $\infty$-categories.
\end{thm}

\begin{proof}
Using the algebraic analogue of Corollary 
\ref{cor:cocartesian-square-}, we may assume that 
$\tau$ is the Nisnevich topology and we may work in 
the non-hypercomplete case. Then, the statement is 
\cite[Theorem 2.1.1]{perfection-sh}. Alternatively, 
we may remark that the proof of 
\cite[Th\'eor\`eme 3.9]{ayoub-etale}
can be extended easily to the case of $\SH_{\Nis}(-;\Lambda)$. 
We explain this below.

The problem is local on $X$, so we may assume that $X$ is affine. By 
\cite[\href{https://stacks.math.columbia.edu/tag/0EUJ}
{Lemma 0EUJ}]{stacks-project}, $X'$ is the limit 
of a cofiltered inverse system of finitely presented 
$X$-schemes $(X'_{\alpha})_{\alpha}$, with $X'_{\alpha}\to X$ universal 
homeomorphisms. Using Proposition \ref{prop:cont-algebraic-29},
we thus reduce to the case where $e$ is assumed to be 
of finite presentation. In this case, writing  
$X$ as the limit of a cofiltered inverse system $(X_{\alpha})_{\alpha}$
consisting of $\Z$-schemes which are essentially of finite type,
the scheme $X'$ is the limit of 
$(X_{\alpha_0}'\times_{X_{\alpha_0}}X_{\alpha})_{\alpha\leq \alpha_0}$
for a finite universal homeomorphism $X_{\alpha_0}'\to X_{\alpha_0}$.
Using Proposition \ref{prop:cont-algebraic-29} again
and base change for finite morphisms, we reduce to the case where
$X$ is of finite type over $\Z$. In conclusion, we may assume that
$X$ has finite Krull dimension and that $X'\to X$ is finite.

Arguing as in the beginning of the proof of Theorem
\ref{thm:seprig}, and using Remark \ref{rmk:conserv-points-sh} 
instead of Corollary \ref{cor:conserv-fiber-funct} and 
base change for finite morphisms instead of 
Theorem \ref{thm:general-base-change-thm}, we 
reduce to the case where $X$ is the spectrum of a field $K$,
and $X'$ the spectrum of a finite purely inseparable 
extension $K'/K$. If $K$ has characteristic zero, then $K=K'$
and there is nothing left to prove. So, we may assume that 
$K$ has positive characteristic $p$. We then write
$\Spec(K)$ as the limit of a cofiltered inverse system of finite 
type $\F_p$-schemes $(X_{\alpha})_{\alpha}$ 
and $\Spec(K')$ as the limit of 
$(X_{\alpha_0}'\times_{X_{\alpha_0}}X_{\alpha})_{\alpha\leq \alpha_0}$
for a finite universal
homeomorphism $X_{\alpha_0}'\to X_{\alpha_0}$. Thus, as before, 
we are finally reduced to treat the case where
$X$ and $X'$ are of finite type over $\F_p$.
This case follows from \cite[Th\'eor\`eme 1.2]{ayoub-etale}. Indeed, 
the condition $(\mathbf{SS}_p)$ of loc.~cit. is satisfied 
for $\SH_{\Nis}(-;\Lambda)$, when $p$ is invertible in $\pi_0\Lambda$,
as shown in \cite[Annexe C]{ayoub-etale}.
In fact, in loc.~cit., this is stated explicitly in 
\cite[Th\'eor\`eme C.1]{ayoub-etale}
for $\DA^{\hyp}_{\et}(-;\Lambda)$, but the proofs apply 
also to $\SH_{\Nis}(-;\Lambda)$. Indeed, the main point is 
to show that elevation to the power $p^n$ on the multiplicative group
$\G_{\rm m}$ induces an autoequivalence of $\M(\G_{\rm m})$
in $\SH(\F_p;\Lambda)$; see \cite[Lemme C.4]{ayoub-etale}. 
This follows from the fact that 
elevation to the power $m$ on $\G_{\rm m}$
induces the endomorphism of
$\Lambda(1)$ given by 
multiplication by the element 
$m_{\epsilon}=\sum_{i=1}^m\langle(-1)^{i-1}\rangle$
in ${\rm K}_0^{\rm MW}(\F_p)$; see \cite[Lemma 3.14]{A1-alg-topol}.
That this element is invertible in the endomorphism ring of 
$\Lambda(1)$ when $m=p^n$ is proven in 
\cite[Lemma 2.2.8]{perfection-sh}.
\end{proof}

\begin{rmk}
\label{rmk:finitness-cond-semi-sep}
In the statement of Theorem
\ref{thm:seprig}, we made the assumption that the rigid analytic space
$X$ has locally finite Krull dimension, whereas the analogous 
assumption was not necessary for Theorem 
\ref{thm:sepalgebr}. This is because we do not know if the analogue of
\cite[\href{https://stacks.math.columbia.edu/tag/0EUJ}
{Lemma 0EUJ}]{stacks-project}
holds for rigid analytic spaces. This is indeed the only 
obstacle for removing the assumption on the Krull dimension in 
Theorem \ref{thm:seprig}. Said differently, semi-separatedness 
for rigid analytic motives holds for a universal homeomorphism
$e:X' \to X$ when, locally on $X$, this morphism can be obtained 
as a weak limit of a cofiltered inverse system of universal 
homeomorphisms $(e_{\alpha}:X'_{\alpha} \to X_{\alpha})_{\alpha}$
where the $X_{\alpha}$'s have finite Krull dimension.
\end{rmk}

\begin{prop}[Separatedness]
\label{prop:separatedness-rig}
\ncn{separatedness}
Let $f:Y\to X$ be a morphism of rigid analytic spaces.
Assume that $X$ is $(\Lambda,\et)$-admissible, and that for 
every point $x\in |X|$, there is a point $y\in |Y|$ mapping to 
$x$ and such that $\kappa(y)$ contains a dense algebraic extension of
$\kappa(x)$. Then the functor 
$$f^*:\RigSH^{(\eff),\,\hyp}_{\et}(X;\Lambda)\to 
\RigSH^{(\eff),\,\hyp}_{\et}(Y;\Lambda)$$
is conservative. 
\end{prop}

\begin{proof}
Using Corollary \ref{cor:conserv-fiber-funct}, 
we reduce to the case of rigid points.
More precisely, we need to prove that 
a morphism $f:y\to x$ of rigid points, with
$\kappa(y)$ containing a dense algebraic extension of
$\kappa(x)$, induces a conservative functor 
$$f^*:\RigSH^{(\eff),\,\hyp}_{\et}(x;\Lambda)\to 
\RigSH^{(\eff),\,\hyp}_{\et}(y;\Lambda).$$
To do so, we may obviously replace $y$ by any rigid
$x$-point $y'$ admitting an $x$-morphism $y'\to y$. 
Since the completion of a separable closure of 
$\kappa(x)$ is algebraically closed, we may 
take for $y'$ a rigid $x$-point 
$\overline{x}$ as in Construction
\ref{cons:point-an-nis-et}(2):
$\kappa(\overline{x})$ is the completion of a separable closure
$\overline{\kappa}(x)$ of $\kappa(x)$ and 
$\kappa^+(\overline{x})$ is the completion of a valuation 
ring $\overline{\kappa}{}^+(x) \subset \overline{\kappa}(x)$
extending $\kappa^+(x)$. In this case, we have 
$\overline{x}\sim\lim_{\alpha} x_{\alpha}$
where $(x_{\alpha})_{\alpha}$ is the inverse system of rigid
$x$-points such that $\kappa(x_{\alpha})$ is a finite subextension  
of $\overline{\kappa}(x)/\kappa(x)$. By Theorem 
\ref{thm:etst}, we have an equivalence:
$$\RigSH^{(\eff),\,\hyp}_{\et}(-;\Lambda)_{\overline{x}}
\simeq \RigSH^{(\eff),\,\hyp}_{\et}(\overline{x};\Lambda)$$
where the left-hand side is the stalk of 
$\RigSH^{(\eff),\,\hyp}_{\et}(-;\Lambda)$ 
at the point $\overline{x}$ of the site $(\Et/x,\et)$. 
Since this point is conservative, we deduce from
Lemma \ref{lem:sheaf-conserv-point}(2)
that the functor 
$$\RigSH_{\et}^{(\eff),\,\hyp}(x;\Lambda) 
\to 
\RigSH_{\et}^{(\eff),\,\hyp}(\overline{x};\Lambda)$$
is conservative, as needed.
\end{proof}

\begin{cor}
\label{cor:separatedness-rig-finite-type}
Let $X$ be a $(\Lambda,\et)$-admissible rigid analytic space, 
and let $f:Y \to X$ be a locally of finite type surjective
morphism. Then the functor 
$$f^*:\RigSH^{(\eff),\,\hyp}_{\et}(X;\Lambda)\to 
\RigSH^{(\eff),\,\hyp}_{\et}(Y;\Lambda)$$
is conservative. 
\end{cor}

\begin{proof}
For every point $x\in |X|$, we may find a point $y\in |Y|$ 
mapping to $x$ and such that $\kappa(y)/\kappa(x)$ is a finite
extension. (This follows from 
\cite[Chapter II, Proposition 8.2.6]{fujiwara-kato}
by a standard argument.)
Thus, the result is a particular case of Proposition 
\ref{prop:separatedness-rig}.
\end{proof}

\begin{cor}
\label{cor:effective-semi-separatedness}
Let $e:X'\to X$ be a universal homeomorphism of rigid analytic 
spaces, and assume that $X$ is $(\Lambda,\et)$-admissible. 
Then, the functor 
$$e^*:\RigSH^{(\eff),\,\hyp}_{\et}(X;\Lambda)
\to \RigSH^{(\eff),\,\hyp}_{\et}(X';\Lambda)$$
is an equivalence of $\infty$-categories.
\end{cor}

\begin{proof}
The morphism $(X')_{\red}\to (X'\times_X X')_{\red}$ is a 
closed immersion and a universal homeomorphism, hence it is an 
isomorphism. Arguing as in \cite[Sous-lemme 1.4]{ayoub-etale}, 
we deduce that $e^*e_*\simeq \id$. Since
$e^*$ is conservative by Proposition
\ref{prop:separatedness-rig}, the result follows.
\end{proof}

\begin{rmk}
\label{rmk:on-effective-semi-sep}
Of course, the $\Tate$-stable case of Corollary 
\ref{cor:effective-semi-separatedness} 
is already covered by Theorem
\ref{thm:seprig} under weaker assumptions. 
The content of this corollary 
is that semi-separatedness holds also for effective
\'etale motives.
It is worth noting that the algebraic analogue of this
result is unknown. 
\end{rmk}

\begin{rmk}
\label{rmk:rigda-rigdm-effective}
Corollary \ref{cor:effective-semi-separatedness}
can be used to improve on the main result 
of \cite{vezz-DADM}. Indeed, given a rigid variety 
$B$ over a non-Archimedean field $K$, 
Corollary \ref{cor:effective-semi-separatedness}
implies that $\RigDA^{\eff,\,(\hyp)}_{\et}(B;\Q)$
is equivalent to the $\infty$-category 
$\RigDA^{\eff,\,(\hyp)}_{{\rm Frob}\et}(B^{\rm Perf};\Q)$
introduced in \cite[Definition 3.5]{vezz-DADM}. Thus, 
assuming that $B$ is normal, 
\cite[Theorem 4.1]{vezz-DADM} can be stated more naturally as 
an equivalence of $\infty$-categories
$$\RigDA^{\eff,\,(\hyp)}_{\et}(B;\Q)\simeq 
\RigDM^{\eff,\,(\hyp)}_{\et}(B;\Q).$$
In fact, this equivalence can be obtained more directly 
by arguing as in the proof of loc.~cit., without mentioning 
the $\infty$-category 
$\RigDA^{\eff,\,(\hyp)}_{{\rm Frob}\et}(B^{\rm Perf};\Q)$.
We leave the details to the interested reader.
\end{rmk}

\subsection{Rigidity}

$\empty$

\smallskip

\label{subsect:rigidity-rigid-motive}

Here, we discuss the rigidity property for rigid analytic motives.
Rigidity is the property that the $\infty$-category of 
torsion \'etale motives over a base is equivalent to 
the $\infty$-category of torsion \'etale sheaves on the small 
\'etale site of the same base.
Rigidity for rigid analytic motives was obtained in 
\cite[Theorem 2.1]{bamb-vezz} for $\RigDA_{\et}^{\hyp}(S;\Lambda)$,
with $S$ of finite type over a non-Archimedean field and 
$\Lambda$ an ordinary torsion ring.
Rigidity in the algebraic setting was obtained in 
\cite[Th\'eor\`eme 4.1]{ayoub-etale} for $\DA_{\et}^{\hyp}(-;\Lambda)$,
with $\Lambda$ an ordinary torsion ring, 
and in \cite[Theorem 6.6]{bachmann-rigidity} for 
$\SH_{\et}^{\hyp}(-;\Lambda)$, with $\Lambda$ the sphere spectrum.
In the recent preprint 
\cite{bachmann-rigidity-remarks}, Bachmann proved rigidity
for effective motives and removed all finiteness assumptions
on the base scheme.
We shall revisit these results in this subsection, 
mainly following \cite{bachmann-rigidity,bachmann-rigidity-remarks}. 

\begin{nota}
\label{not:p-torsion-stable-infty-cat}
Let $\mathcal{C}$ be a stable presentable 
$\infty$-category and $\ell$ a prime number. 
An object $A$ of $\mathcal{C}$ is said to be $\ell$-nilpotent if the 
zero object of $\mathcal{C}$ is a colimit of the $\N$-diagram
$$A \xrightarrow{\ell\cdot \id} A \xrightarrow{\ell\cdot \id} 
A \xrightarrow{\ell\cdot \id} \cdots.$$
An object $A$ of $\mathcal{C}$ is said to be 
$\ell$-complete if the zero object of $\mathcal{C}$ is a limit of
the $\N^{\op}$-diagram
$$\cdots \xrightarrow{\ell\cdot \id} A \xrightarrow{\ell\cdot \id} 
A \xrightarrow{\ell\cdot \id} A.$$
We denote by $\mathcal{C}_{\ellnil}\subset \mathcal{C}$
and $\mathcal{C}_{\ellcpl}\subset \mathcal{C}$
the sub-$\infty$-categories spanned by 
$\ell$-nilpotent and $\ell$-complete objects
respectively. Given an object $A$ of $\mathcal{C}$, we denote 
by $A/\ell^n$ the cofiber of the map
$$\ell^n\cdot \id:A \to A.$$
Since multiplication by $\ell^{2n}$ is zero on 
$A/\ell^n$, it is both $\ell$-nilpotent and $\ell$-complete.
\ncn{nilpotent@$\ell$-nilpotent}
\ncn{complete@$\ell$-complete}
\symn{$(-)_{\ellnil}$}
\symn{$(-)_{\ellcpl}$}
\end{nota}

We gather a few facts concerning the notions of 
$\ell$-nilpotent and $\ell$-complete objects
in the following remark. We refer the reader to 
\cite[Part II, Chapter 7]{lurie:SAG}
where these notions are developed in greater 
generality. See also \cite[\S 2.1]{bachmann-rigidity}.

\begin{rmk}
\label{rmk:on-ell-nilp-compl}
Let $\mathcal{C}$ be a stable presentable $\infty$-category and 
$\ell$ a prime number. We denote by $\mathcal{C}[\ell^{-1}]$ 
the full sub-$\infty$-category of $\mathcal{C}$ spanned by those
objects for which multiplication by $\ell$ is an equivalence.
\symn{$(-)[\ell^{-1}]$}
\begin{enumerate}

\item[(1)] The $\infty$-category $\mathcal{C}_{\ellnil}$ is stable,
presentable and generated under colimits by the 
objects of the form $A/\ell^n$, for $A\in \mathcal{C}$. 
The inclusion functor $\mathcal{C}_{\ellnil}\to \mathcal{C}$
commutes with colimits and finite limits. 
If $\mathcal{C}$ is compactly generated, then so it is
$\mathcal{C}_{\ellnil}$.

\item[(2)] The $\infty$-category 
$\mathcal{C}_{\ellcpl}$ is the localisation of
$\mathcal{C}$ with respect to the maps $0\to A$, for 
$A\in \mathcal{C}[\ell^{-1}]$. We denote by 
$(-)^{\wedge}_{\ell}:\mathcal{C}\to \mathcal{C}_{\ellcpl}$ the {left} adjoint to the inclusion functor.
This is called the $\ell$-completion functor.
\symn{$(-)^{\wedge}_{\ell}$}

\item[(3)] The $\ell$-completion functor induces an equivalence of 
$\infty$-categories
$$(-)^{\wedge}_{\ell}:\mathcal{C}_{\ellnil}
\xrightarrow{\sim} \mathcal{C}_{\ellcpl}.$$
In particular, we see that $\mathcal{C}_{\ellcpl}$ is stable,
presentable and generated under colimits by the 
objects of the form $A/\ell^n$, for $A\in \mathcal{C}$.
If $\mathcal{C}$ is compactly generated, then so is
$\mathcal{C}_{\ellcpl}$.

\item[(4)] If $\mathcal{C}$ underlies a presentable symmetric monoidal 
$\infty$-category $\mathcal{C}^{\otimes}$, then there is an 
essentially unique morphism
$\mathcal{C}^{\otimes} \to \mathcal{C}^{\otimes}_{\ellcpl}$
in $\CAlg(\Prl)$ 
whose underlying functor is 
$(-)^{\wedge}_{\ell}:\mathcal{C}\to \mathcal{C}_{\ellcpl}$.

\item[(5)]
Suppose that $\mathcal{C}$ is given as a colimit 
in $\Prl$ of an inductive system 
$(\mathcal{C}_{\alpha})_{\alpha}$ of stable presentable 
$\infty$-categories. Then $\mathcal{C}[\ell^{-1}]$
is also the colimit of the inductive system 
$(\mathcal{C}_{\alpha}[\ell^{-1}])_{\alpha}$ in $\Prl$.
(This uses the fact that a colimit in $\Prl$ can be computed 
as a limit in $\Prr$.)
In particular, $\mathcal{C}[\ell^{-1}]$
is generated under colimits by the images of the functors 
$\mathcal{C}_{\alpha}[\ell^{-1}]\to \mathcal{C}[\ell^{-1}]$. 
It follows from (2) and the universal property of localisations
(see \cite[Proposition 5.5.4.20]{lurie})
that $\mathcal{C}_{\ellcpl}$ is the colimit in $\Prl$ of the 
inductive system $(\mathcal{C}_{\alpha,\,\ellcpl})_{\alpha}$.
Using (3), we deduce that 
$\mathcal{C}_{\ellnil}$ is also the  
colimit in $\Prl$ of the inductive system
$(\mathcal{C}_{\alpha,\,\ellnil})_{\alpha}$.

\end{enumerate}
\end{rmk}

\begin{thm}[Rigidity]
\label{thm:rigrig}
\ncn{rigidity}
Let $S$ be a rigid analytic space and $\ell$ a prime number
which is invertible in $\widetilde{\kappa}(s)$ for every
$s\in |S|$. Assume one of the following two alternatives.
\begin{enumerate}

\item[(1)] We work in the non-hypercomplete case and $\Lambda$ 
is eventually coconnective.

\item[(2)] We work in the hypercomplete case.

\end{enumerate}
Then the obvious functor 
\begin{equation}
\label{eq-thm:rigrig-1}
\Shv_{\et}^{(\hyp)}(\Et/S;\Lambda)_{\ellcpl} \to 
\RigSH^{(\eff,\,\hyp)}_{\et}(S;\Lambda)_{\ellcpl}
\end{equation}
is an equivalence of $\infty$-categories. 
(The same is true with ``$\ellnil$'' instead of 
``$\ellcpl$''.)
\end{thm}

We also have the algebraic analogue of Theorem \ref{thm:rigrig} 
which can be stated as follows.

\begin{thm}
\label{thm:rigrig-algebraic}
Let $S$ be a scheme and $\ell$ a prime number
which is invertible on $S$. Assume one of the following 
two alternatives.
\begin{enumerate}

\item[(1)] We work in the non-hypercomplete case and $\Lambda$ 
is eventually coconnective.

\item[(2)] We work in the hypercomplete case.

\end{enumerate}
Then the obvious functor 
\begin{equation}
\label{eq-thm:rigrig-algebraic}
\Shv_{\et}^{(\hyp)}(\Et/S;\Lambda)_{\ellcpl} \to 
\SH^{(\eff,\,\hyp)}_{\et}(S;\Lambda)_{\ellcpl}
\end{equation}
is an equivalence of $\infty$-categories.
(The same is true with ``$\ellnil$'' instead of 
``$\ellcpl$''.)
\end{thm}

\begin{proof}
We first consider the alternative (1). We may assume that 
$S$ is affine and given as the limit of a cofiltered inverse system 
$(S_{\alpha})_{\alpha}$ of affine schemes of finite type over $\Z$. 
By the algebraic analogue of Lemma 
\ref{lem:sheaves-on-limit-etale} and
Proposition \ref{prop:cont-algebraic-29}, 
it is enough to prove the conclusion for the $S_{\alpha}$'s.
Thus, me may assume that $S$ is of finite type over 
$\Z$ and hence $(\Lambda,\et)$-admissible. 
By the algebraic analogue of Lemma 
\ref{lem:auto-hypercomp-etale}(2) and Proposition
\ref{prop:automatic-hypercomp-motives-algebraic} below, 
we are then automatically 
working in the hypercomplete case. This means that we only need
to consider the alternative (2). 
In that case, the result is essentially 
\cite[Theorem 6.6]{bachmann-rigidity}
improved in \cite[Theorem 3.1]{bachmann-rigidity-remarks}.
\end{proof}

\begin{rmk}
\label{rmk:reduction-rigidity-alt-2}
Arguing as above, we only need to prove Theorem \ref{thm:rigrig} 
under the second alternative.
Indeed, by Lemma \ref{lem:sheaves-on-limit-etale} and
Theorem \ref{thm:anstC} we may assume that $S$ is 
$(\Lambda,\et)$-admissible. In this case, there is no
distinction between the hypercomplete and the non-hypercomplete 
cases by Lemma \ref{lem:auto-hypercomp-etale}(2) and
Proposition \ref{prop:automatic-hypercomp-motives}.
\end{rmk}

{
Our proof of Theorem \ref{thm:rigrig} follows the arguments 
in \cite{bachmann-rigidity,bachmann-rigidity-remarks}
and relies on some of the key steps in loc.~cit.
We start with a reduction to the $(\Lambda,\et)$-admissible case. 

\begin{lemma}
\label{lemma:reduction-rigidity-to-admissible}
To prove Theorem \ref{thm:rigrig},
we may work in the hypercomplete case and assume that 
$S$ is $(\Lambda,\et)$-admissible.
\end{lemma}

\begin{proof}
We said already that it is enough to work under the second alternative.
Assume that Theorem \ref{thm:rigrig}
is known in the hypercomplete case when the base
is $(\Lambda,\et)$-admissible. 
To prove the theorem in general, 
we argue as in the proof of 
\cite[Theorem 3.1]{bachmann-rigidity-remarks}.
We may assume that 
$S=\mathcal{S}^{\rig}$ where $\mathcal{S}$ is an affine formal 
scheme given as the limit of an affine formal pro-scheme 
$(\mathcal{S}_{\alpha})_{\alpha}$ such that the 
$\mathcal{S}_{\alpha}$'s are of finite type over 
$\Z[\ell^{-1}][[\pi]]$. We set 
$S_{\alpha}=\mathcal{S}_{\alpha}^{\rig}$; these are 
$(\Lambda,\et)$-admissible rigid analytic spaces. By 
Lemma \ref{lem:sheaves-on-limit-etale}, 
Theorem \ref{thm:anstC} and Remark 
\ref{rmk:on-ell-nilp-compl}(5),  
we have a commutative square 
$$\xymatrix{\colim_{\alpha}\Shv_{\Nis}(S_{\alpha};\Lambda)_{\ellcpl}
\ar[r]^-{\sim} \ar[d] & \Shv_{\Nis}(S;\Lambda)_{\ellcpl} \ar[d] \\
\colim_{\alpha} \RigSH^{(\eff)}_{\Nis}(S_{\alpha};\Lambda)_{\ellcpl}
\ar[r]^-{\sim} & \RigSH^{(\eff)}_{\Nis}(S;\Lambda)_{\ellcpl}}$$
where the horizontal arrows are
equivalences of $\infty$-categories.
It follows that in the analogous commutative square
$$\xymatrix{\colim_{\alpha}\Shv^{\hyp}_{\et}(S_{\alpha};\Lambda)_{\ellcpl}
\ar[r] \ar[d]^-{\sim} & \Shv^{\hyp}_{\et}(S;\Lambda)_{\ellcpl} \ar[d] \\
\colim_{\alpha} \RigSH^{(\eff),\,\hyp}_{\et}(S_{\alpha};\Lambda)_{\ellcpl}
\ar[r] & \RigSH^{(\eff),\,\hyp}_{\et}
(S;\Lambda)_{\ellcpl},\,}$$
the horizontal arrows are localisation functors, whereas, 
by assumption, the left vertical arrow is an equivalence. 
This shows that 
\begin{equation}
\label{eq-lemma:reduction-rigidity-to-admissible-7}
\Shv^{\hyp}_{\et}(S;\Lambda)_{\ellcpl}
\to \RigSH^{(\eff),\,\hyp}_{\et}
(S;\Lambda)_{\ellcpl}
\end{equation}
is a localisation functor. 
To finish the proof, it remains to see that 
\eqref{eq-lemma:reduction-rigidity-to-admissible-7}
is conservative.
Given a geometric rigid point 
$\overline{s} \to S$, we have a commutative square
$$\xymatrix{\Shv^{\hyp}_{\et}(S;\Lambda)_{\ellcpl}
\ar[r] \ar[d] &  \RigSH^{(\eff),\,\hyp}_{\et}
(S;\Lambda)_{\ellcpl} \ar[d]\\
\Shv^{\hyp}_{\et}(\overline{s};\Lambda)_{\ellcpl}
\ar[r]^-{\sim} &  \RigSH^{(\eff),\,\hyp}_{\et}
(\overline{s};\Lambda)_{\ellcpl},\,}$$
and the bottom arrow is an equivalence, again by assumption, 
since $\overline{s}$ is $(\Lambda,\et)$-admissible. 
This proves that the functor
$(-)_{\overline{s}}:\Shv^{\hyp}_{\et}(S;\Lambda)_{\ellcpl}
\to (\Mod_{\Lambda})_{\ellcpl}$ factors through
\eqref{eq-lemma:reduction-rigidity-to-admissible-7}.
We conclude using Propositions
\ref{prop:enough-points-rig-an}
and 
\ref{prop:stalk-equivalence}.
\end{proof}

}

We now introduce some notations.

\begin{nota}
\label{nota:lambda-ell-etale-sheaf}
Let $S$ be a rigid analytic space. The 
$\ell$-completion of the constant \'etale sheaf $\Lambda\in 
\Shv^{(\hyp)}_{\et}(\Et/S;\Lambda)$ will be denoted simply 
by $\Lambda_{\ell}$. This is the unit object of 
$\Shv^{(\hyp)}_{\et}(\Et/S;\Lambda)_{\ellcpl}$ endowed 
with its natural monoidal structure.
We denote by 
$$\iota_S^*:\Shv^{(\hyp)}_{\et}(\Et/S;\Lambda)
\to \RigSH_{\et}^{\eff,\,(\hyp)}(S;\Lambda)$$
the obvious functor, and by $\iota_{S,\,*}$ its right adjoint.
Similarly, we denote by 
$$\iota_{S,\,\ell}^*:\Shv^{(\hyp)}_{\et}(\Et/S;\Lambda)_{\ellcpl}
\to \RigSH_{\et}^{\eff,\,(\hyp)}(S;\Lambda)_{\ellcpl}$$
the functor induced by $\iota_S^*$ on $\ell$-completed 
objects, and by
$\iota_{S,\,\ell,\,*}$ its right adjoint. We denote by
$$\Sigma_{\Tate,\,\ell}^{\infty}:
\RigSH^{\eff,\,(\hyp)}_{\et}(S;\Lambda)_{\ellcpl}
\to \RigSH^{(\hyp)}_{\et}(S;\Lambda)_{\ellcpl}$$
the functor induced by $\Sigma_{\Tate}^{\infty}$ 
on $\ell$-completed objects, 
and by \sym{$\Omega^{\infty}_{\Tate,\,\ell}$} its right adjoint. 
(See Definition \ref{dfn:rigsh-stable}.) 
The functor 
\eqref{eq-thm:rigrig-1}
is given by $\iota^*_{S,\,\ell}$ in the effective case and by 
$\Sigma_{\Tate,\,\ell}^{\infty}\circ 
\iota^*_{S,\,\ell}$ in the $\Tate$-stable case.
These notations apply also when $S$ is a scheme.
\symn{$\iota^*, \iota_*$}
\symn{$\Sigma_{\Tate,\,\ell}^{\infty}$}
\end{nota}

Recall that 
$\U^1_S$ is the relative unit sphere over the rigid analytic space $S$.
(See Notation \ref{not:relative-ball}(3).)

\begin{lemma}
\label{lem:conseq-algeb-rig}
Let $S$ be a rigid analytic space and $\ell$ a prime number
which is invertible in $\kappa(s)$ for every $s\in |S|$. 
There is a $\otimes$-invertible object $\Lambda_{\ell}(1)$ 
in $\Shv^{\hyp}_{\et}(S;\Lambda)_{\ellcpl}$
together with a morphism
\begin{equation}
\label{eq-lem:conseq-algeb-rig-4}
\sigma:\Lambda_{\ell} \to \Lambda_{\ell}(1)[1]
\end{equation}
in $\Shv^{\hyp}_{\et}(\Et/\U^1_S;\Lambda)_{\ellcpl}$
endowed with a trivialisation (i.e., a homotopy 
to the null morphism) over the unit section
$1_S\subset \U^1_S$.
Moreover, the induced morphism 
$\sigma:\Tate^{\wedge}_{\ell} \to 
\iota^*_{S,\,\ell}(\Lambda_{\ell}(1)[1])$
is an equivalence in $\RigSH^{\eff,\,\hyp}_{\et}(S;\Lambda)_{\ellcpl}$.
\end{lemma}

\begin{proof}
We may construct $\Lambda_{\ell}(1)$ and $\sigma:\Lambda_{\ell}\to
\Lambda_{\ell}(1)[1]$ locally on $S$ provided that the 
construction is compatible with base change. 
Assume that $S=\Spf(A)^{\rig}$ with 
$A$ an adic ring. Let $I\subset A$ be an ideal of definition and set 
$U=\Spec(A)\smallsetminus \Spec(A/I)$.
We denote by $\Lambda_{\ell}(1)
\in \Shv^{\hyp}_{\et}(\Et/U;\Lambda)_{\ellcpl}$
the $\otimes$-invertible object obtained from  
the one introduced in \cite[Definition 3.9]{bachmann-rigidity}
by extension of scalars to $\Lambda$. Also, let 
$\sigma:\Lambda_{\ell} \to \Lambda_{\ell}(1)[1]$
be the morphism in $\Shv_{\et}^{\hyp}(\Et/\A^1_U\smallsetminus 
0_U;\Lambda)_{\ellcpl}$ obtained from the one
introduced in \cite[Definition 3.13]{bachmann-rigidity}
by extension of scalars to $\Lambda$.
As explained in the beginning of \cite[\S 6]{bachmann-rigidity},
a trivialisation of $\sigma$ above $1_S$ gives rise to a morphism 
$\Tate^{\wedge}_{\ell} \to \iota_{U,\,\ell}^*\Lambda_{\ell}(1)[1]$
in $\SH^{\eff,\,\hyp}_{\et}(U;\Lambda)_{\ellcpl}$.
As explained in the beginning of the proof of 
\cite[Theorem 3.1]{bachmann-rigidity-remarks}, 
this morphism is an equivalence 
(see also \cite[Theorem 6.5]{bachmann-rigidity}
in the $\Tate$-stable case).
The lemma follows now from the existence of a commutative
square of stable presentable $\infty$-categories
$$\xymatrix{\Shv_{\et}^{\hyp}(\Et/U;\Lambda) \ar[r]^-{\iota_U^*} 
\ar[d] & 
\SH^{\eff,\,\hyp}_{\et}(U;\Lambda) \ar[d]
\\
\Shv_{\et}^{\hyp}(\Et/S;\Lambda) \ar[r]^-{\iota_S^*} & 
\RigSH^{\eff,\,\hyp}_{\et}(S;\Lambda) 
}$$
where the vertical arrows are induced by the analytification functor. 
\end{proof}

{
\begin{cor}
\label{cor:ell-completed-RigSH-eff-non-eff}
Let $S$ be a rigid analytic space and $\ell$ a prime number
which is invertible in $\kappa(s)$ for every $s\in |S|$.
Then the obvious functor 
$$\Sigma_{\Tate,\,\ell}^{\infty}:
\RigSH^{\eff,\,\hyp}_{\et}(S;\Lambda)_{\ellcpl}
\to \RigSH^{\eff,\,\hyp}_{\et}(S;\Lambda)_{\ellcpl}$$
is an equivalence of $\infty$-categories.
(The same is true with ``$\ellnil$'' instead of 
``$\ellcpl$''.)
\end{cor}

\begin{proof}
Indeed, by Remarks 
\ref{rmk:symmetric-obj-} 
and
\ref{rmk:on-ell-nilp-compl}(5),
$\RigSH^{\hyp}_{\et}(S;\Lambda)_{\ellcpl}$ is
the colimit in $\Prl$ of the $\N$-diagram whose 
transition maps are given by 
tensoring with $\Tate^{\wedge}_{\ell}$ in 
$\RigSH^{\eff,\,\hyp}_{\et}(S;\Lambda)_{\ellcpl}$.
The result follows since $\Tate^{\wedge}_{\ell}$
is $\otimes$-invertible by Lemma 
\ref{lem:conseq-algeb-rig}.
\end{proof}

}

\begin{lemma}
\label{lem:fully-faith-rigidity}
Let $S$ be a $(\Lambda,\et)$-admissible rigid analytic space 
and $\ell$ a prime number which is invertible in 
$\widetilde{\kappa}(s)$ for every $s\in |S|$. 
Then the obvious functor
\begin{equation} 
\label{eq-prop:fully-faith-rigidity-1}
\Shv_{\et}^{\hyp}(\Et/S;\Lambda)_{\ellcpl} \to 
\RigSH^{(\eff),\,\hyp}_{\et}(S;\Lambda)_{\ellcpl}
\end{equation}
is fully faithful.
(The same is true with ``$\ellnil$'' instead of 
``$\ellcpl$''.)
\end{lemma}

\begin{proof}
By Corollary
\ref{cor:ell-completed-RigSH-eff-non-eff},
we only need to treat the effective case.
The functor 
$$\iota_S^*:\PSh(\Et/S;\Lambda)
\to \PSh(\RigSm/S;\Lambda)$$
is fully faithful and its right adjoint
commutes with \'etale hypersheafification. 
It follows that the induced functor on \'etale hypersheaves
$$\iota_S^*:\Shv_{\et}^{\hyp}(\Et/S;\Lambda)
\to \Shv_{\et}^{\hyp}(\RigSm/S;\Lambda)$$
is also fully faithful, and the same is true for the induced
functor on $\ell$-complete objects 
$$\iota_{S,\,\ell}^*:\Shv_{\et}^{\hyp}(\Et/S;\Lambda)_{\ellcpl}
\to \Shv_{\et}^{\hyp}(\RigSm/S;\Lambda)_{\ellcpl}.$$
We claim that the functor $\iota_{S,\,\ell}^*$ 
takes values in the sub-$\infty$-category 
$\RigSH^{\eff,\,\hyp}_{\et}(S;\Lambda)_{\ellcpl}$
spanned by $\B^1$-local objects; this would finish the proof. 
Indeed, let $\mathcal{F}$ be an $\ell$-complete 
\'etale hypersheaf of $\Lambda$-modules on $\Et/S$.
Saying that $\iota_{S,\,\ell}^*\mathcal{F}$ is $\B^1$-local
is equivalent to saying that for every $X\in \RigSm/S$, the map 
$\Gamma(X;\mathcal{F}|_X) \to \Gamma(\B^1_X;\mathcal{F}|_{\B^1_X})$
is an equivalence. (Here, we denote by $\mathcal{F}|_X$ 
the $\ell$-complete inverse image of $\mathcal{F}$ along the morphism 
$X \to S$, and similarly for $\mathcal{F}|_{\B^1_X}$.)
Since $X$ is $(\Lambda,\et)$-admissible, the claim 
follows from Lemma \ref{lem:B-1-invariance-etale-coh}(1) below
(see also \cite[Example 0.1.1(2)]{huber}).
\end{proof}

\begin{lemma}
\label{lem:B-1-invariance-etale-coh}
Let $X$ be a $(\Lambda,\et)$-admissible rigid analytic space 
and $\ell$ a prime number which is invertible in 
$\widetilde{\kappa}(x)$ for every $x\in |X|$. 
Let $p:\B^1_X \to X$ be the obvious 
projection and let $\mathcal{F}$ be an $\ell$-complete \'etale hypersheaf 
on $\Et/X$. Then the map $\mathcal{F} \to 
p_*p^*\mathcal{F}$ is an equivalence. 
\end{lemma}

\begin{proof}
It is enough to prove the results
on the stalks for all geometric algebraic rigid points 
$\overline{x} \to X$. Using Remark
\ref{rmk-thm:general-base-change-thm}, 
we reduce to show the following. Given a geometric 
rigid point $s=\Spf(V)^{\rig}$ and an 
$\ell$-complete \'etale hypersheaf of $\Lambda$-modules
$\mathcal{F}$ on $\Et/s$,
the map $\mathcal{F}(s) 
\to \Gamma(\B^1_s;\mathcal{F}|_{\B^1_s})$
is an equivalence.
Using Lemmas \ref{lem:pi-0-Lambda-coh-dim} and 
\ref{lem:Lamba-tau-coh-pointwise}, we reduce to the case where 
$\mathcal{F}$ is bounded. By an easy induction, 
we reduce to the case where $\mathcal{F}$ is discrete, 
and we may then assume that $\mathcal{F}$ is an ordinary 
\'etale sheaf of $\Z/\ell^n$-modules. 
The site $(\Et/s,\et)$ is equivalent to 
$(\FRigEt/\Spf(V),\riget)$ and, since $s$ is geometric, 
it is also equivalent to $(\Et/\Spec(V'),\et)$,
where $V'=V/\sqrt{(\pi)}$ with $\pi$ a generator 
of an ideal of definition of $V$. Thus, we may consider 
$\mathcal{F}$ as an ordinary \'etale sheaf on 
$\FRigEt/\Spf(V)$ and on $\Et/\Spec(V')$.
We then have equivalences:
\begin{equation}
\label{eq-lem:B-1-invariance-etale-coh-1}
\Rder\Gamma_{\et}(\B^1_s;\mathcal{F}|_{\B^1_s})\simeq 
\Rder\Gamma_{\riget}(\A^1_V;\mathcal{F}|_{\A^1_V})
\simeq \Rder\Gamma_{\et}(\A^1_{V'};
i^*j_*(\Z/\ell^n) \otimes_{\Z/\ell^n}
\mathcal{F}|_{\A^1_{V'}}).
\end{equation}
Here $i$ denotes the closed immersion $\Spec(V')
\to \Spec(V)$ and its base changes, and $j$ denotes 
the open complement of $i$ and its base changes.
The second equivalence in 
\eqref{eq-lem:B-1-invariance-etale-coh-1}
follows from \cite[Corollary 3.5.16]{huber}. 
(More precisely, we reduce to the case where $\mathcal{F}$ 
is of the form 
$i'_*\Z/\ell^n$ with $i':\Spec(V'') \to \Spec(V')$ 
a closed immersion, and we remark that \cite[Corollary 3.5.16]{huber}
is still valid if we replace the closed point of $\Spec(V)$ 
by a closed subscheme contained in $\Spec(V')$.) 
Using the smooth base change theorem in \'etale cohomology
\cite[Expos\'e XVI, Th\'eor\`eme 1.1]{SGAIV3} 
and the fact that the fraction field of 
$V$ is algebraically closed, 
we deduce that $i^*j_*\Z/\ell^n \simeq \Z/\ell^n$
on $\A^1_{V'}$. Thus, the last term in 
\eqref{eq-lem:B-1-invariance-etale-coh-1}
is equivalent to $\Rder\Gamma_{\et}(\A^1_{V'};
\mathcal{F}|_{\A^1_{V'}})$
which, by homotopy invariance of \'etale cohomology
\cite[Expos\'e XV, Corollaire 2.2]{SGAIV3}, is 
equivalent to $\mathcal{F}(V')\simeq \mathcal{F}(s)$.
This proves that
$\mathcal{F}(s)$ is indeed equivalent to 
$\Rder\Gamma(\B^1_s;\mathcal{F}|_{\B^1_s})$ as needed.
\end{proof}

\begin{proof}[Proof of Theorem \ref{thm:rigrig}]
Using Lemmas
\ref{lemma:reduction-rigidity-to-admissible}
and 
\ref{lem:fully-faith-rigidity},
it remains to see that the 
functor \eqref{eq-prop:fully-faith-rigidity-1}
is essentially surjective (still under the assumption that 
$S$ is $(\Lambda,\et)$-admissible).
Moreover, it is enough to do so in the $\Tate$-stable case,
by Corollary 
\ref{cor:ell-completed-RigSH-eff-non-eff}.
We follow the argument used in the proof of 
\cite[Theorem 2.1]{bamb-vezz}.

The question being local on $S$, we may assume that 
$S=\Spf(A)^{\rig}$ with $A$ an adic ring of 
principal ideal type. Let $\pi\in A$ be a generator of an 
ideal of definition and set $U=\Spec(A[\pi^{-1}])$. 
It is enough to show that the image of the functor 
\eqref{eq-prop:fully-faith-rigidity-1}, 
in the $\Tate$-stable case, 
contains a set of generators of 
$\RigSH^{\hyp}_{\et}(S;\Lambda)_{\ellcpl}$.
Such a set of generators is given, up to shift and Tate twists, 
by $\M(V)/\ell^n$ where $n\in \N$ and 
$V=\Spf(B)^{\rig}$ with $B$ a rig-\'etale adic $A$-algebra
satisfying the conclusion of Proposition 
\ref{prop:formal-compl-smooth}. Thus, there exists a smooth 
affine $U$-scheme $X$ and an open immersion 
$v:V \to X^{\an}$. Since we are allowed to replace $V$ 
by the components of an analytic hypercover, we 
may assume that $\Omega_{X/U}$ is free.
Fix a projective compactification $j:X\to P$ over $U$
and denote by $f:X \to U$ and $p:P\to U$ the structural morphisms.
Thus, we have a commutative diagram
$$\xymatrix{V \ar[r]^-v \ar[dr]_-g  
\ar@/^2pc/[rr]^-{\overline{v}} & X^{\an} \ar[r]^-{j^{\an}} 
\ar[d]^-{f^{\an}} & P^{\an} 
\ar[dl]^-{p^{\an}} \\
& S. &}$$
The motive $\M(V)$ is equivalent to $g_{\sharp}\Lambda\simeq 
f^{\an}_{\sharp} v_{\sharp}\Lambda$. 
Using Corollary \ref{cor:6ffalg},
we see that $\M(V)$ is equivalent, up to shift and Tate twist, to
$f^{\an}_! v_{\sharp}\Lambda\simeq p^{\an}_!
j^{\an}_{\sharp} v_{\sharp}\Lambda \simeq p_*
\overline{v}_{\sharp}\Lambda$.

Using Lemmas \ref{lem:conseq-algeb-rig}
and \ref{lem:fully-faith-rigidity}, the image of the 
functor \eqref{eq-prop:fully-faith-rigidity-1}, 
in the $\Tate$-stable case,
is closed under shift and Tate twists. Therefore, 
it remains to see that the latter image
contains $p^{\an}_* \overline{v}_{\sharp}\Lambda/\ell^n$.
Clearly, $\overline{v}_{\sharp}\Lambda/\ell^n$ belongs to the 
image of 
$$\Sigma_{\Tate,\,\ell}^{\infty} \circ 
\iota_{P^{\an},\,\ell}^*:
\Shv_{\et}^{\hyp}(\Et/P^{\an};\Lambda)_{\ellcpl}
\to \RigSH^{\hyp}_{\et}(P^{\an};\Lambda)_{\ellcpl}.$$
Thus, it is enough to show that the natural transformation 
$\Sigma_{\Tate,\,\ell}^{\infty}\circ 
\iota_{S,\,\ell}^* \circ p^{\an}_*
\to p^{\an}_*\circ \Sigma_{\Tate,\,\ell}^{\infty} \circ 
\iota_{P^{\an},\,\ell}^*$
is an equivalence. (The first $p^{\an}_*$ is the direct image functor 
on \'etale hypersheaves, and the second $p^{\an}_*$ is the direct 
image functor on rigid analytic motives.) 
Using Corollary 
\ref{cor:ell-completed-RigSH-eff-non-eff},
it is enough to show that the natural 
transformation $\iota_{S,\,\ell}^* \circ p^{\an}_* \to 
p^{\an}_* \circ  \iota_{P^{\an},\,\ell}^*$
is an equivalence. Given an $\ell$-complete \'etale hypersheaf 
$\mathcal{F}$ on $\Et/P^{\an}$, the evaluation of
$\iota_{S,\,\ell}^*p^{\an}_*\mathcal{F} \to 
p^{\an}_*\iota_{P^{\an}}^*\mathcal{F}$ on a smooth 
rigid analytic $S$-space $Y$
is given by
$$\Gamma(Y;g^*p^{\an}_*\mathcal{F})\to \Gamma(Y\times_S 
P^{\an}, g'^*\mathcal{F})=
\Gamma(Y;p'_*g'^*\mathcal{F})$$
where $p'$ and $g'$ are as in the Cartesian square
$$\xymatrix{Y\times_S P^{\an} \ar[r]^-{g'} \ar[d]^-{p'}
& P^{\an} \ar[d]^-{p^{\an}} \\
Y \ar[r]^g & S.\!}$$
The result follows now from the quasi-compact base change theorem, 
see Remark \ref{rmk-thm:general-base-change-thm}.
\end{proof}

\section{Rigid analytic motives as modules in formal motives}

\label{sect:formal-mot-module-main-thm}

This section contains one of the key results of the paper
which, roughly speaking, gives a description of the functor 
$\RigSH^{(\hyp)}_{\tau}(-;\Lambda)$ 
in terms of the functor $\SH^{(\hyp)}_{\tau}(-;\Lambda)$.
This can be considered as a vast generalisation of 
\cite[Scholie 1.3.26]{ayoub-rig}. In fact, we prefer 
to work with the functor 
$\FSH^{(\hyp)}_{\tau}(-;\Lambda)$, sending a formal scheme 
to the $\infty$-category of formal motives, instead of the functor 
$\SH^{(\hyp)}_{\tau}(-;\Lambda)$, but this is a merely aesthetic
difference, by Theorem \ref{thm:formal-mot-alg-mot}. 
For a precise form of the description alluded to, 
we refer the reader to Theorems \ref{thm:main-thm-}
and \ref{thm:compute-chi}.

We start by recalling the definition and the basic properties
of the $\infty$-category $\FSH^{(\eff,\,\hyp)}_{\tau}
(\mathcal{S};\Lambda)$ of formal motives 
over a formal scheme $\mathcal{S}$.

\subsection{Formal and algebraic motives}

$\empty$

\smallskip

\label{subsect:formal-and-alg-mot}

Recall that we denote by $\FSch$ the category of formal schemes
and that, given a formal scheme $\mathcal{S}$, we denote 
by $\FSm/\mathcal{S}$ the category of smooth formal 
$\mathcal{S}$-schemes. 
(Notations \ref{not:categ-formal-scheme} and 
\ref{not:big-sites}.)
The $\infty$-category of formal motives over a formal scheme 
is constructed as in Definitions \ref{def:DAeff} and 
\ref{dfn:rigsh-stable}. 

We fix a formal scheme $\mathcal{S}$ and $\tau\in \{\Nis,\et\}$. 

\begin{dfn}
\label{def:DAeff-form}
Let $\FSH^{\eff,\,(\hyp)}_{\tau}(\mathcal{S};\Lambda)$ 
be the full sub-$\infty$-category 
of $\Shv_{\tau}^{(\hyp)}(\FSm/\mathcal{S};\Lambda)$
spanned by those objects which are local with respect to 
the collection of maps of the form 
$\Lambda_{\tau}(\A^1_{\mathcal{X}}) \to
\Lambda_{\tau}(\mathcal{X})$, for $\mathcal{X}\in
\FSm/\mathcal{S}$, and their desuspensions. Let
\begin{equation}
\label{eq-def:DAeff-form-1}
\Lder_{\A^1}:\Shv^{(\hyp)}_{\tau}(\FSm/\mathcal{S};\Lambda)
\to \FSH^{\eff,\,(\hyp)}_{\tau}(\mathcal{S};\Lambda)
\end{equation}
be the left adjoint to the obvious inclusion.
This is called the $\A^1$-localisation functor.
Given a smooth formal $\mathcal{S}$-scheme $\mathcal{X}$, we set
$\M^{\eff}(X)=\Lder_{\A^1}(\Lambda_{\tau}(\mathcal{X}))$. 
This is the effective motive of $\mathcal{X}$.
\symn{$\FSH^{\eff,\,(\hyp)}$}
\symn{$\Lder_{\A^1}$}
\symn{$\M^{\eff}$}
\end{dfn}

\begin{rmk}
\label{rmk:mon-str-fsh-eff}
By \cite[Proposition 2.2.1.9]{lurie:higher-algebra}, 
$\FSH^{\eff,\,(\hyp)}_{\tau}(\mathcal{S};\Lambda)$ underlies a unique
monoidal $\infty$-category
$\FSH^{\eff,\,(\hyp)}_{\tau}(\mathcal{S};\Lambda)^{\otimes}$ such that  
$\Lder_{\A^1}$ lifts to a monoidal functor.
Moreover, this monoidal $\infty$-category is 
presentable, i.e., belongs to $\CAlg(\Prl)$.
\end{rmk}

\begin{dfn}
\label{dfn:fsh-stable}
Let $\Tate_{\mathcal{S}}$ (or simply \sym{$\Tate$} if $\mathcal{S}$ 
is clear from the context) be the image by $\Lder_{\A^1}$ of the 
cofiber of the split inclusion
$\Lambda_{\tau}(\mathcal{S})\to
\Lambda_{\tau}(\A^1_{\mathcal{S}}\smallsetminus 0_{\mathcal{S}})$
induced by the unit section. With the notation of 
\cite[Definition 2.6]{robalo:k-theory-bridge}, we set 
\begin{equation}
\label{eq-def:formDAeff-2}
\FSH_{\tau}^{(\hyp)}(\mathcal{S};\Lambda)^{\otimes}
=\FSH^{\eff,\,(\hyp)}_{\tau}(\mathcal{S};\Lambda)^{\otimes}
[\Tate_{\mathcal{S}}^{-1}].
\end{equation}
More precisely, there is a morphism
$\Sigma^{\infty}_{\Tate}:
\FSH_{\tau}^{\eff,\,(\hyp)}(\mathcal{S};\Lambda)^{\otimes}
\to \FSH^{(\hyp)}_{\tau}(\mathcal{S};\Lambda)^{\otimes}$
in $\CAlg(\Prl)$, sending 
$\Tate_{\mathcal{S}}$ to a $\otimes$-invertible object, and which is 
initial for this property. We denote by 
$\Omega^{\infty}_{\Tate}:\FSH^{(\hyp)}_{\tau}(\mathcal{S};\Lambda)
\to \FSH^{\eff,\,(\hyp)}_{\tau}(\mathcal{S};\Lambda)$
the right adjoint to $\Sigma_{\Tate}^{\infty}$.
Given a smooth formal $\mathcal{S}$-scheme $\mathcal{X}$, 
we set $\M(\mathcal{X})=\Sigma^{\infty}_{\Tate}\M^{\eff}(\mathcal{X})$. This is the motive of $X$.
\symn{$\FSH^{(\hyp)}$}
\symn{$\Sigma^{\infty}_{\Tate}$}
\symn{$\Omega^{\infty}_{\Tate}$}
\symn{$\M$}
\end{dfn}

\begin{dfn}
\label{dfn:form-mot}
Objects of $\FSH^{(\hyp)}_{\tau}(\mathcal{S};\Lambda)$ 
are called \nc{formal motives}
over $\mathcal{S}$. We will denote by $\Lambda$ 
(or $\Lambda_{\mathcal{S}}$ if we need
to be more precise) 
the monoidal unit of $\FSH^{(\hyp)}_{\tau}(\mathcal{S};\Lambda)$. 
For any $n\in\N$, we denote by \sym{$\Lambda(n)$} the image of 
$\Tate_{\mathcal{S}}^{\otimes n}[-n]$ by $\Sigma^{\infty}_{\Tate}$, 
and by $\Lambda(-n)$ the $\otimes$-inverse of $\Lambda(n)$.
For $n\in \Z$, we denote by 
$M\mapsto M(n)$ the \nc{Tate twist} given by tensoring with $\Lambda(n)$.
\end{dfn}

\begin{rmk}
\label{rmk:on-formal-motives}
$\empty$

\begin{enumerate}

\item[(1)] Remark \ref{rmk:symmetric-obj-} applies also in the case 
of formal motives: the $\infty$-category 
$\FSH^{(\hyp)}_{\tau}(\mathcal{S};\Lambda)$ underlying 
\eqref{eq-def:formDAeff-2}
is equivalent to the colimit in $\Prl$ 
of the $\N$-diagram whose transition maps are given by 
tensoring with $\Tate_{\mathcal{S}}$ in $\FSH^{\eff,\,(\hyp)}_{\tau}
(\mathcal{S};\Lambda)$.

\item[(2)] When $\Lambda$ is the Eilenberg--Mac Lane spectrum associated
to an ordinary ring, also denoted by $\Lambda$, the category 
$\FSH_{\tau}^{(\eff,\,\hyp)}(\mathcal{S};\Lambda)$ 
is more commonly denoted 
by $\FDA^{(\eff,\,\hyp)}_{\tau}(\mathcal{S};\Lambda)$. Also, when 
$\tau$ is the Nisnevich topology, we sometimes 
drop the subscript ``$\Nis$''.
\symn{$\FDA^{(\eff,\,\hyp)}$}

\item[(3)] Just as in Remark \ref{rmk:rigda-model-cat}, 
there is a more traditional description of the $\infty$-category
$\FSH^{(\eff,\,\hyp)}_{\tau}(\mathcal{S};\Lambda)$ 
using the language of model categories. This is the approach taken in 
\cite[\S 1.4.2]{ayoub-rig}.

\item[(4)] If $S$ is an ordinary scheme considered as a formal scheme in 
the obvious way, i.e., such that the zero ideal 
is an ideal of definition, then the $\infty$-category 
$\FSH^{(\eff,\,\hyp)}_{\tau}(S;\Lambda)$ is the usual $\infty$-category
$\SH^{(\eff,\,\hyp)}_{\tau}(S;\Lambda)$ of algebraic motives over $S$.
More generally, by Theorem 
\ref{thm:formal-mot-alg-mot}
below, the $\infty$-categories introduced in 
Definitions \ref{def:DAeff-form} and 
\ref{dfn:fsh-stable} are always equivalent to 
$\infty$-categories of algebraic motives.

\end{enumerate}
\end{rmk}

\begin{lemma}
\label{lem:generation-fsh}
The monoidal $\infty$-category 
$\FSH^{(\eff,\,\hyp)}_{\tau}(\mathcal{S};\Lambda)^{\otimes}$
is presentable and its underlying $\infty$-category is 
generated under colimits, and up to desuspension and negative 
Tate twists when applicable, by the motives 
$\M^{(\eff)}(\mathcal{X})$ with
$\mathcal{X}\in \FSm/\mathcal{S}$ quasi-compact and quasi-separated.
\end{lemma}

\begin{proof}
See the proof of Lemma
\ref{lem:generation-rigsh}.
\end{proof}

\begin{prop}
\label{prop:basic-functorial-fsh}
The assignment $\mathcal{S}\mapsto \FSH_{\tau}^{(\eff,\,\hyp)}
(\mathcal{S};\Lambda)^{\otimes}$
extends naturally into a functor
\begin{equation}
\label{eq-cor:compactness-fda}
\FSH_{\tau}^{(\eff,\,\hyp)}(-;\Lambda)^{\otimes}:
\FSch^{\op}\to \CAlg(\Prl).
\end{equation}
\end{prop}

\begin{proof}
We refer to \cite[\S 9.1]{robalo} for the construction of 
an analogous functor in the algebraic setting.
\end{proof}

\begin{nota}
\label{nota:f-stars-for-form}
Let $f:\mathcal{Y}\to \mathcal{X}$ 
be a morphism of formal schemes. 
The image of $f$ by 
\eqref{eq-cor:compactness-fda}
is the inverse image functor 
$$f^*:\FSH^{(\eff,\,\hyp)}_{\tau}(\mathcal{X};\Lambda)
\to 
\FSH^{(\eff,\,\hyp)}_{\tau}(\mathcal{Y};\Lambda)$$
which has the structure of a monoidal functor.
Its right adjoint $f_*$ is the direct image functor.
It has the structure of a right-lax monoidal functor.
(See Lemma \ref{lem:adj-monoidal-cat} below.)
\end{nota}

\begin{nota}
\label{nota:special-fiber-sigma-star}
Recall that we denote by $\mathcal{X}_{\sigma}$ 
the special fiber of a formal scheme $\mathcal{X}$. 
(See Notation \ref{not:special-fiber}.)
The functor $\mathcal{X}\mapsto \mathcal{X}_{\sigma}$ 
induces a functor $(-)_{\sigma}:\FSm/\mathcal{S}\to
\Sm/\mathcal{S}_{\sigma}$ which is continuous for the 
topology $\tau$. By the functoriality of the construction 
of $\infty$-categories of motives, we deduce an adjunction 
\begin{equation}
\label{eqn:adj-FDA-DA}
\sigma^*:\FSH_{\tau}^{(\eff,\,\hyp)}(\mathcal{S};\Lambda)
\rightleftarrows \SH_{\tau}^{(\eff,\,\hyp)}
(\mathcal{S}_{\sigma};\Lambda):\sigma_*.
\end{equation}
In fact, modulo the identification of Remark
\ref{rmk:on-formal-motives}(4), $\sigma^*$ is simply the inverse image 
functor associated to the morphism of formal schemes
$\mathcal{X}_{\sigma} \to \mathcal{X}$. 
\symn{$\sigma^*, \sigma_*$}
\end{nota}

\begin{thm}
\label{thm:formal-mot-alg-mot}
The functors $\sigma^*$ and $\sigma_*$ in 
\eqref{eqn:adj-FDA-DA}
are equivalences of $\infty$-categories.
\end{thm}

\begin{proof}
This is \cite[Corollaires 1.4.24 \& 1.4.29]{ayoub-rig} under
the assumption that $\mathcal{S}$ is of finite type
over $\Spf(k^{\circ})$, with $k^{\circ}$ a complete valuation 
ring of height $\leq 1$. However, this assumption is not 
used in the proofs of these results. 
\end{proof}

\begin{rmk}
\label{rmk:6-oper-formal-sch}
Let $f:\mathcal{Y} \to \mathcal{X}$ be a morphism of 
formal schemes. Modulo the equivalences of Theorem 
\ref{thm:formal-mot-alg-mot}, the operations $f^*$ and $f_*$
coincide with the operations $f^*_{\sigma}$ and $f_{\sigma,*}$
associated to the morphism of schemes 
$f_{\sigma}:\mathcal{Y}_{\sigma}\to \mathcal{X}_{\sigma}$.
When $f_{\sigma}$ is locally of finite type, we denote by
$f_!$ and $f^!$ the operations on formal motives corresponding 
to the operations $f_{\sigma,!}$ and $f_{\sigma}^!$ on 
algebraic motives modulo the equivalences of Theorem 
\ref{thm:formal-mot-alg-mot} (in the $\Tate$-stable case). Similarly, 
if $f_{\sigma}$ is smooth, we denote by $f_{\sharp}$ the 
operation corresponding to $f_{\sigma,\,\sharp}$.
\end{rmk}

\begin{nota}
\label{nota:generic-fiber-eta-star}
Recall that we denote by $\mathcal{X}^{\rig}$ the generic 
fiber of a formal scheme $\mathcal{X}$. (See 
Notation \ref{not:rigid-analytic-spaces}.)
The functor $\mathcal{X}\mapsto \mathcal{X}^{\rig}$ 
induces a functor $(-)^{\rig}:\FSm/\mathcal{S}
\to \RigSm/\mathcal{S}^{\rig}$ which is continuous for the 
topology $\tau$. By the functoriality of the construction 
of $\infty$-categories of motives, we deduce an adjunction 
\begin{equation}
\label{eq:adj-FDA-RigDA}
\xi_{\mathcal{S}}:\FSH_{\tau}^{(\eff,\,\hyp)}(\mathcal{S};\Lambda)
\rightleftarrows \RigSH_{\tau}^{(\eff,\,\hyp)}
(\mathcal{S}^{\rig};\Lambda):\chi_{\mathcal{S}}.
\end{equation}
Composing with the equivalences of Theorem
\ref{thm:formal-mot-alg-mot}, we get also an equivalent adjunction 
\begin{equation}
\label{eq:adj-DA-RigDA}
\xi_{\mathcal{S}}:\SH_{\tau}^{(\eff,\,\hyp)}(\mathcal{S}_{\sigma};\Lambda)
\rightleftarrows \RigSH_{\tau}^{(\eff,\,\hyp)}
(\mathcal{S}^{\rig};\Lambda):\chi_{\mathcal{S}}.
\end{equation}
These adjunctions will play an important role in this section.
\symn{$\xi$}\symn{$\chi$}
\end{nota}

\begin{prop}
\label{prop:xi-presheaf}
The functors $\xi_{\mathcal{S}}$, for $\mathcal{S}\in \FSch$,
are part of a morphism of $\CAlg(\Prl)$-valued presheaves
\begin{equation}
\label{eq-prop:xi-presheaf-1}
\xi:\FSH_{\tau}^{(\eff,\,\hyp)}(-;\Lambda)^{\otimes}
\to \RigSH_{\tau}^{(\eff,\,\hyp)}((-)^{\rig};\Lambda)^{\otimes}
\end{equation}
on $\FSch$. In particular, the functors $\xi_{\mathcal{S}}$ 
are monoidal and commute with the inverse image functors. Moreover, if 
$f:\mathcal{T}\to \mathcal{S}$ 
is a smooth morphism in $\FSch$, the natural transformation 
$$f^{\rig}_{\sharp}\circ \xi_{\mathcal{T}} \to \xi_{\mathcal{S}}
\circ f_{\sharp}$$
is an equivalence.
\end{prop}

\begin{proof}
One argues as in \cite[\S 9.1]{robalo} for  
the first assertion. The second assertion is clear.
\end{proof}

In the rest of this subsection we use the above constructions
to produce a convenient conservative family of functors for the 
$\infty$-category $\RigSH^{(\eff,\,\hyp)}_{\tau}(S;\Lambda)$, 
for $S$ a rigid analytic space.
This family is rather big: it is indexed by formal models of smooth 
rigid analytic $S$-spaces. For a better result, we refer the reader
to Corollary \ref{cor:chi-non-trivial-generators} below.
We start by recording the following general fact.

\begin{prop}
\label{conservativity-generation}
Let $(F_i:\mathcal{C}_i \to \mathcal{D})_i$ 
be a small family of functors in 
$\Prl$ having the same target $\mathcal{D}$. 
Let $G_i$ be the right adjoint of $F_i$. 
Then the following conditions are equivalent:
\begin{enumerate}

\item[(1)] the family $(G_i:\mathcal{D}\to \mathcal{C}_i)_{i\in I}$
is conservative;

\item[(2)] $\mathcal{D}$ is generated under colimits by 
objects of the form $F_i(A)$, with $A\in \mathcal{C}_i$.

\end{enumerate}
\end{prop}

\begin{proof}
Assume first that (2) is satisfied. Let $f:X\to Y$ be a map in 
$\mathcal{D}$ such that $G_i(f)$ is an equivalence for every $i$.
We want to show that $f$ is an equivalence.  
To do so, consider the full sub-$\infty$-category 
$\mathcal{D}_0\subset \mathcal{D}$
spanned by objects $E$ such that $\Map_{\mathcal{D}}(E,X)
\to \Map_{\mathcal{D}}(E,Y)$ is an equivalence.
Clearly, $\mathcal{D}_0$ is stable under arbitrary colimits and
contains the images of the $F_i$'s. By (2), it follows that 
$\mathcal{D}_0=\mathcal{D}$, and thus $f$ is an equivalence 
by the Yoneda lemma.

We now assume that (1) is satisfied.
Denote by $\mathcal{D}'\subset \mathcal{D}$ the smallest full
sub-$\infty$-category containing the images of the $F_i$'s and 
stable under arbitrary colimits. We need to show that 
$\mathcal{D}'=\mathcal{D}$. We claim that 
the $\infty$-category $\mathcal{D}'$
is presentable. Indeed, as the $F_i$'s are colimit-preserving
and the $\mathcal{C}_i$'s are presentable, $\mathcal{D}'$ is 
the smallest sub-$\infty$-category of $\mathcal{D}$
stable under colimits and containing a certain small set of 
objects (namely the union of images of sets of generators 
for the $\mathcal{C}_i$'s). These objects are $\kappa$-compact
for $\kappa$ large enough. Thus, our claim follows 
from Lemma \ref{lem:characteri-presentable-infty-cat}.
Using \cite[Corollary 5.5.2.9]{lurie}, 
we may thus consider the right adjoint $\rho$
to the inclusion functor $\mathcal{D}'\to \mathcal{D}$.
Fix an object $X\in \mathcal{D}$. We will show that
$\rho(X)\to X$ is an equivalence, which will finish the proof.
Since the $G_i$'s form
a conservative family, it is enough to show that 
the maps $G_i(\rho(X)) \to G_i(X)$ are equivalences.
By the Yoneda lemma, it is enough to show that the maps 
$$\Map_{\mathcal{C}_i}(A,G_i(\rho(X)))\to 
\Map_{\mathcal{C}_i}(A,G_i(X))$$ 
are equivalences for all $A\in \mathcal{C}_i$.
By adjunction, these maps are equivalent to 
$$\Map_{\mathcal{D}}(F_i(A),\rho(X))\to 
\Map_{\mathcal{D}}(F_i(A),X),$$
which are equivalences since the $F_i(A)$'s belong to $\mathcal{D}'$.
\end{proof}

\begin{prop}
\label{prop:chi-trivial-generators} 
Let $S$ be a rigid analytic space. 
For every $U\in \RigSm^{\qcqs}/S$, denote by 
$f_U:U \to S$ the structural morphism and choose a formal model 
$\mathcal{U}$ of $U$. 
Then, the functors
$$\chi_{\mathcal{U}}\circ f_U^*:
\RigSH^{(\eff,\,\hyp)}_{\tau}(S;\Lambda)
\to \FSH^{(\eff,\,\hyp)}_{\tau}(\mathcal{U};\Lambda),$$
for $U\in \RigSm^{\qcqs}/S$, 
form a conservative family. In fact, the same is true
if we restrict to those $U$'s admitting affine formal models
of principal ideal type.
\end{prop}

\begin{proof}
The functor $\chi_{\mathcal{U}}\circ f_U^*$
has a left adjoint $f_{U,\,\sharp}\circ \xi_{\mathcal{U}}$
sending the monoidal unit of 
$\FSH^{(\eff,\,\hyp)}_{\tau}(\mathcal{U};\Lambda)$
to $\M^{(\eff)}(U)$. We conclude by 
Lemma \ref{lem:generation-rigsh} and Proposition
\ref{conservativity-generation}.
\end{proof}

\subsection{Descent, continuity and stalks, I. 
The case of formal motives}

$\empty$

\smallskip

In this subsection, we gather a few basic properties of the 
functor $\mathcal{S}\mapsto 
\FSH_{\tau}^{(\eff,\,\hyp)}(\mathcal{S};\Lambda)$,
$f\mapsto f^*$, from Proposition 
\ref{prop:basic-functorial-fsh}. 
We fix a topology $\tau\in \{\Nis,\et\}$.

\begin{prop}
\label{prop:FDA-hypersheaf}
The contravariant functor 
$$\mathcal{S}\mapsto \FSH^{(\eff,\,\hyp)}_{\tau}(\mathcal{S};\Lambda), 
\quad f\mapsto f^*$$ 
defines a $\tau$-(hyper)sheaf on $\FSch$ with values in $\Prl$.
\end{prop}

\begin{proof}
The proof is similar to that of Theorem \ref{thm:hyperdesc}. 
It suffices to prove that for every formal scheme $\mathcal{S}$, 
the functor 
$$\FSH_{\tau}^{(\eff,\,\hyp)}(-;\Lambda):
(\Et/\mathcal{S})^{\op} \to \Prl,$$
is a $\tau$-(hyper)sheaf. One reduces, 
by an essentially formal argument, to show that the functor 
$$\Shv_{\tau}^{(\hyp)}(\FSm/-;\Lambda):
(\Et/\mathcal{S})^{\op} \to \Prl$$
is a $\tau$-(hyper)sheaf, and this follows from 
Corollary \ref{cor:sheafy-sheaves}. The formal argument alluded to
can be found in the proof of Theorem
\ref{thm:hyperdesc}, and we will not repeat it here.
\end{proof}

A formal scheme $\mathcal{S}$ is said to be 
$(\Lambda,\tau)$-admissible (resp. $(\Lambda,\tau)$-good)
if the scheme $\mathcal{S}_{\sigma}$ is 
$(\Lambda,\tau)$-admissible (resp. $(\Lambda,\tau)$-good)
in the sense of Definition 
\ref{dfn:lambda-tau-admiss}.
\ncn{formal schemes!$(\Lambda,\tau)$-admissible}
\ncn{formal schemes!$(\Lambda,\tau)$-good
}
\begin{prop}
\label{prop:automatic-hypercomp-motives-algebraic}
Let $\tau\in\{\Nis,\et\}$ and let $\mathcal{S}$ 
be a $(\Lambda,\tau)$-admissible formal scheme. 
When $\tau$ is the \'etale topology, assume 
that $\Lambda$ is eventually coconnective. Then, we have
$$\FSH_{\tau}^{(\eff),\,\hyp}(\mathcal{S};\Lambda)=
\FSH_{\tau}^{(\eff)}(\mathcal{S};\Lambda).$$ 
\end{prop}

\begin{proof}
This is proven in the same way as Proposition
\ref{prop:automatic-hypercomp-motives}.
\end{proof}

\begin{prop}
\label{prop:compact-shv-forsm}
Let $\mathcal{S}$ be a formal scheme. 
\begin{enumerate}

\item[(1)] 
The $\infty$-category $\FSH^{(\eff)}_{\tau}(\mathcal{S};\Lambda)$
is compactly generated if $\tau$ is the Nisnevich topology or 
if $\Lambda$ is eventually coconnective.
A set of compact generators is given,
up to desuspension and negative Tate twists when applicable, 
by the $\M^{(\eff)}(\mathcal{X})$ for $\mathcal{X}\in \FSm/\mathcal{S}$ 
quasi-compact, quasi-separated and $(\Lambda,\tau)$-good. 

\item[(2)] The $\infty$-category $\FSH^{(\eff),\,\hyp}_{\tau}
(\mathcal{S};\Lambda)$
is compactly generated if $\mathcal{S}$ is $(\Lambda,\tau)$-admissible.
A set of compact generators is given, up to desuspension 
and negative Tate twists when applicable, 
by the $\M^{(\eff)}(\mathcal{X})$ for $\mathcal{X}\in \FSm/\mathcal{S}$ 
quasi-compact, quasi-separated and $(\Lambda,\tau)$-good.

\end{enumerate}
Moreover, under the stated assumptions, the monoidal $\infty$-category 
$\FSH^{(\eff,\,\hyp)}_{\tau}(\mathcal{S};\Lambda)^{\otimes}$
belongs to $\CAlg(\Prl_{\omega})$ and, if $f:\mathcal{T}\to 
\mathcal{S}$ is a quasi-compact and quasi-separated 
morphism of formal schemes with $\mathcal{T}$
assumed $(\Lambda,\tau)$-admissible in the hypercomplete case, 
the functor 
$f^*:\FSH^{(\eff,\,\hyp)}_{\tau}(\mathcal{S};\Lambda)\to 
\FSH^{(\eff,\,\hyp)}_{\tau}(\mathcal{T};\Lambda)$ is compact-preserving, 
i.e., belongs to $\Prl_{\omega}$.
\end{prop}

\begin{proof}
This is proven in the same way as Proposition
\ref{prop:compact-shv-rigsm}.
\end{proof}

Given a formal scheme $\mathcal{S}$, we write 
``$\pvcd_{\Lambda}(\mathcal{S})$'' instead of 
``$\pvcd_{\Lambda}(\mathcal{S}_{\sigma})$''; 
see Definition \ref{dfn:Lambda-coh-dim-4}.
Our next statement is an analogue of Theorem 
\ref{thm:anstC} for formal motives. 
\symn{$\pvcd$}

\begin{prop}
\label{prop:anstC-FDA}
Let $(\mathcal{S}_{\alpha})_{\alpha}$ be a 
cofiltered inverse system of quasi-compact and quasi-separated 
formal schemes with affine transition maps, 
and let $\mathcal{S}=\lim_{\alpha}\mathcal{S}_{\alpha}$
be the limit of this system. We assume one of the following two 
alternatives. 
\begin{enumerate}

\item[(1)] We work in the non-hypercomplete case. 

\item[(2)] We work in the hypercomplete case, and $\mathcal{S}$ and 
the $\mathcal{S}_{\alpha}$'s are $(\Lambda,\tau)$-admissible. 
When $\tau$ is the \'etale topology, we assume furthermore that
$\Lambda$ is eventually coconnective or that  
the numbers $\pvcd_{\Lambda}(\mathcal{S}_{\alpha})$ 
are bounded independently of $\alpha$.

\end{enumerate}
Then the obvious functor
$$\underset{\alpha}{\colim}\,\FSH_{\tau}^{(\eff,\,\hyp)}
(\mathcal{S}_{\alpha};\Lambda)
\to \FSH^{(\eff,\,\hyp)}_{\tau}(\mathcal{S};\Lambda),$$
where the colimit is taken in $\Prl$, is an equivalence.
\end{prop}

\begin{proof}
This follows immediately from Proposition
\ref{prop:cont-algebraic-29} and Theorem
\ref{thm:formal-mot-alg-mot}.
\end{proof}

We will use Proposition 
\ref{prop:anstC-FDA}
to compute the stalks of $\FSH^{(\eff,\,\hyp)}_{\tau}(-;\Lambda)$
for the topology $\rig\text{-}\tau$ on $\FSch$.
(See Corollary \ref{cor:rig-topol-formal-schemes}). 
We first describe a conservative family of points 
for this topology.

\begin{rmk}
\label{rmk:rigid-point-formal-scheme}
Let $\mathcal{S}$ be a formal scheme. A rigid point of 
$\mathcal{S}$ is a morphism $\mathfrak{s}:\Spf(V) \to \mathcal{S}$ 
where $V$ is an adic valuation ring of principal ideal type.
We sometimes also denote by $\mathfrak{s}$ the formal scheme
$\Spf(V)$.
The assignment 
$(\Spf(V)\to \mathcal{S})\mapsto (\Spf(V)^{\rig}\to \mathcal{S}^{\rig})$
is an equivalence of groupoids between rigid points of
$\mathcal{S}$ and those of $\mathcal{S}^{\rig}$. (See Remark  
\ref{rmk:tau-geom-pt-of-S}.) We will say that a rigid 
point $\mathfrak{s}:\Spf(V) \to \mathcal{S}$ is algebraic 
(resp. $\tau$-geometric)
if the associated rigid point of $\mathcal{S}^{\rig}$
is algebraic (resp. $\tau$-geometric). See 
Remarks \ref{rmk:algebraic-morphism-etale-etc} and 
\ref{rmk:tau-geom-pt-of-S}, and Definition 
\ref{dfn:geometric-points-tau}.
\ncn{formal schemes!rigid points}
\ncn{rigid points!of formal schemes}
\end{rmk}

\begin{prop}
\label{prop:enough-points-rig-tau-20u}
Let $\mathcal{S}$ be a formal scheme.
We denote by $\FRigEt/\mathcal{S}$ the category of 
rig-\'etale formal $\mathcal{S}$-schemes. 
Then, the site $(\FRigEt/\mathcal{S},\rig\text{-}\tau)$
admits a conservative family of points
indexed by $\tau$-geometric algebraic rigid points 
$\mathfrak{s}=\Spf(V) \to \mathcal{S}$. To such a rigid point
$\mathfrak{s}$, the associated topos-theoretic point is given by 
$$\mathcal{F}\mapsto \mathcal{F}_{\mathfrak{s}}=
\underset{\Spf(V) \to \mathcal{U} \to \mathcal{S}}{\colim}\,
\mathcal{F}(\mathcal{U})$$
where the colimit is over rig-\'etale neighbourhoods $\mathcal{U}$ 
of $\mathfrak{s}$. Moreover, one may restrict to those 
rigid points of $\mathcal{S}^{\rig}$
as in Construction \ref{cons:point-an-nis-et}.
\symn{$\FRigEt$}
\end{prop}

\begin{proof}
This follows from Corollary \ref{cor:rig-topol-formal-schemes} 
and Proposition \ref{prop:enough-points-rig-an}.
\end{proof}

\begin{prop}
\label{prop:FDAstalks}
Let $\mathcal{S}$ be a formal scheme and let 
$\mathfrak{s}\to \mathcal{S}$ be an algebraic rigid point 
of $\mathcal{S}$. 
Assume one of the following two alternatives.
\begin{enumerate}

\item[(1)] We work in the non-hypercomplete case.

\item[(2)] We work in the hypercomplete case, and  
$\mathcal{S}$ and $\mathcal{S}^{\rig}$ are $(\Lambda,\tau)$-admissible.
When $\tau$ is the \'etale topology, we assume furthermore 
that $\Lambda$ is eventually coconnective or that the numbers 
$\pvcd_{\Lambda}(\mathcal{S}')$, for admissible blowups 
$\mathcal{S}' \to \mathcal{S}$, 
are bounded independently of $\mathcal{S}'$. 

\end{enumerate}
Then there is an equivalence of $\infty$-categories
$$\FSH^{(\eff,\,\hyp)}_{\tau}(-;\Lambda)_{\mathfrak{s}}
\simeq \FSH^{(\eff,\,\hyp)}_{\tau}(\mathfrak{s};\Lambda)$$
where the left-hand side is the stalk of 
$\FSH_{\tau}^{(\eff,\,\hyp)}(-;\Lambda)$ 
at $\mathfrak{s}$, i.e., the colimit, taken in $\Prl$,
of the diagram $(\mathfrak{s} \to \mathcal{U}\to \mathcal{S})
\mapsto \FSH_{\tau}^{(\eff,\,\hyp)}(\mathcal{U};\Lambda)$ with 
$\mathcal{U}\in \FRigEt/\mathcal{S}$.
\end{prop}

\begin{proof}
This follows from Proposition
\ref{prop:anstC-FDA}. Indeed, the condition that 
$\mathcal{S}^{\rig}$ is $(\Lambda,\tau)$-admissible implies that 
$\mathfrak{s}$ is $(\Lambda,\tau)$-admissible. Moreover, if 
the numbers $\pvcd_{\Lambda}(\mathcal{S}')$
are bounded independently of $\mathcal{S}'$
for admissible blowups $\mathcal{S}'\to \mathcal{S}$, 
then the same is true
for the numbers $\pvcd_{\Lambda}(\mathcal{U})$ for the saturated 
rig-\'etale neighbourhoods $\mathfrak{s} \to \mathcal{U}\to \mathcal{S}$.
\end{proof}

\subsection{Statement of the main result}

$\empty$

\smallskip

\label{subsect:statement-main-result}

Let $\mathcal{S}$ be a formal scheme. 
By Proposition \ref{prop:xi-presheaf}, 
we have a monoidal functor
$$\xi^{\otimes}_{\mathcal{S}}:
\FSH^{(\eff,\,\hyp)}_{\tau}(\mathcal{S};\Lambda)^{\otimes}
\to \RigSH^{(\eff,\,\hyp)}_{\tau}
(\mathcal{S}^{\rig};\Lambda)^{\otimes}.$$ 
From Corollary \ref{cor:adjunction-on-com-alge} below, 
we deduce that $\chi_{\mathcal{S}}\Lambda$ 
underlies a commutative algebra in the monoidal 
$\infty$-category 
$\FSH^{(\eff,\hyp)}_{\tau}(\mathcal{S};\Lambda)^{\otimes}$, 
which we also denote by $\chi_{\mathcal{S}}\Lambda$.
Moreover, the functor $\chi_{\mathcal{S}}$ admits a factorization 
$$\RigSH^{(\eff,\,\hyp)}_{\tau}(\mathcal{S}^{\rig};\Lambda)
\xrightarrow{\widetilde{\chi}_{\mathcal{S}}}
\FSH^{(\eff,\,\hyp)}_{\tau}(\mathcal{S};\chi \Lambda) 
\xrightarrow{\rm ff} \FSH^{(\eff,\,\hyp)}_{\tau}(\mathcal{S};\Lambda),$$
where $\FSH^{(\eff,\,\hyp)}_{\tau}(\mathcal{S};\chi \Lambda)$ 
is the $\infty$-category of 
$\chi_{\mathcal{S}}\Lambda$-modules in 
$\FSH^{(\eff,\,\hyp)}_{\tau}(\mathcal{S};\Lambda)^{\otimes}$
and ${\rm ff}$ is the forgetful functor.
The functor $\widetilde{\chi}_{\mathcal{S}}$ admits a left adjoint 
$$\widetilde{\xi}_{\mathcal{S}}:
\FSH^{(\eff,\,\hyp)}_{\tau}(\mathcal{S};\chi\Lambda) \to 
\RigSH^{(\eff,\,\hyp)}_{\tau}(\mathcal{S}^{\rig};\Lambda)$$
that sends a $\chi_{\mathcal{S}}\Lambda$-module $M$ to 
$\xi_{\mathcal{S}}(M)\otimes_{\xi_{\mathcal{S}}
\chi_{\mathcal{S}}\Lambda}\Lambda$.
It will be important for us to know that the functors 
$\widetilde{\xi}_{\mathcal{S}}$, for $\mathcal{S}\in \FSch$, 
are part of a morphism 
$$\widetilde{\xi}{}^{\otimes}:\FSH^{(\eff,\,\hyp)}_{\tau}
(-;\chi\Lambda)^{\otimes}
\to \RigSH^{(\eff,\,\hyp)}_{\tau}((-)^{\rig};\Lambda)^{\otimes}$$
in the $\infty$-category 
$\PSh(\FSch;\CAlg(\Prl))$
of presheaves on $\FSch$ valued in 
$\CAlg(\Prl)$. The construction of $\widetilde{\xi}{}^{\otimes}$
will be carried in Subsection \ref{subsect:const-tilde-xi}
below. Before stating the main result of this section,
we introduce the following assumptions. 
\symn{$\FSH^{(\eff,\,\hyp)}(-;\chi \Lambda)$}
\symn{$\widetilde{\chi}$}
\symn{$\widetilde{\xi}$}

\begin{assu}
\label{assu:for-main-thm}
We assume (at least) one of the following four alternatives:
\begin{enumerate}

\item[(i)] $\tau$ is the Nisnevich topology;

\item[(ii)] $\pi_0\Lambda$ is a $\Q$-algebra;

\item[(iii)] we work in the non-hypercomplete case,
$\Lambda$ is eventually coconnective and every prime number 
which is not invertible in $\pi_0\Lambda$ is invertible 
on every formal scheme we consider;

\item[(iv)] we work in the hypercomplete case, 
every formal scheme we consider is
$(\Lambda,\tau)$-admissible and its generic fiber is also 
$(\Lambda,\tau)$-admissible, and every prime number 
which is not invertible in $\pi_0\Lambda$ is invertible 
on every formal scheme we consider.

\end{enumerate}
Moreover, under one of the alternatives (iii) or (iv), 
when we write ``$\FSch$'', we actually mean the full 
subcategory of formal schemes satisfying the properties in 
(iii) or (iv) respectively.
\end{assu}

\begin{assu}
\label{assu:for-main-thm-2}
We assume that $\tau$ is the \'etale topology and that
one of the two alternatives (iii) or (iv)
above is satisfied.
\end{assu}

\begin{thm}
\label{thm:main-thm-}
$\empty$

\begin{enumerate}

\item[(1)] We work under Assumption 
\ref{assu:for-main-thm}.
Given a formal scheme $\mathcal{S}$, the functor 
$$\widetilde{\xi}_{\mathcal{S}}:\FSH^{(\hyp)}_{\tau}
(\mathcal{S};\chi\Lambda) \to 
\RigSH^{(\hyp)}_{\tau}(\mathcal{S}^{\rig};\Lambda)$$
is fully faithful.

\item[(2)] We work under Assumption 
\ref{assu:for-main-thm-2}. The morphism of $\CAlg(\Prl)$-valued presheaves
$$\widetilde{\xi}{}^{\otimes}:
\FSH^{(\hyp)}_{\et}(-;\chi\Lambda)^{\otimes}
\to \RigSH^{(\hyp)}_{\et}((-)^{\rig};\Lambda)^{\otimes}$$
exhibits $\RigSH^{(\hyp)}_{\et}((-)^{\rig};\Lambda)^{\otimes}$ 
as the rig-\'etale sheaf associated to 
$\FSH^{(\hyp)}_{\et}(-;\chi\Lambda)^{\otimes}$.

\end{enumerate}
\end{thm}

\begin{rmk}
\label{rmk:on-effective-version-main-thm}
Our proof of Theorem \ref{thm:main-thm-}
relies crucially on $\Tate$-stability.
Therefore, we do not expect this theorem to hold 
for the effective $\infty$-categories of motives.
\end{rmk}

{
\begin{rmk}
\label{rmk:main-on-rig}
One can reformulate Theorem 
\ref{thm:main-thm-}(2) as an equivalence between functors defined on 
rigid analytic spaces. Indeed, by Corollary 
\ref{cor:rig-topol-formal-schemes}, 
we have an equivalence of sites 
$$(\RigSpc^{\qcqs},\et) \xrightarrow{\sim} (\FSch^{\qcqs},\riget).$$
Moreover, the left Kan extension of the $\CAlg(\Prl)$-valued
presheaf
$\FSH^{(\hyp)}_{\et}(-;\chi\Lambda)^{\otimes}$ 
along the functor $(-)^{\rig}:\FSch^{\qcqs} \to \RigSpc^{\qcqs}$
is easily seen to be given by 
\begin{equation}
\label{eq-rmk:main-on-rig-1}
S\mapsto \underset{\mathcal{S}\in \Mdl(S)}{\colim}\,
\FSH^{(\hyp)}_{\et}(\mathcal{S};\chi\Lambda)^{\otimes}.
\end{equation}
(See Notation \ref{nota:Mdl}.) 
Thus, Theorem \ref{thm:main-thm-}(2) implies that the 
morphism of 
$\CAlg(\Prl)$-valued
presheaves given by 
$$\underset{\mathcal{S}\in \Mdl(S)}{\colim}\,
\FSH^{(\hyp)}_{\et}(\mathcal{S};\chi\Lambda)^{\otimes}
\to \RigSH^{(\hyp)}_{\et}(S;\Lambda)^{\otimes}$$
exhibits $\RigSH^{(\hyp)}_{\et}(-;\Lambda)^{\otimes}$
as the \'etale sheafification of the 
$\CAlg(\Prl)$-valued presheaf 
\eqref{eq-rmk:main-on-rig-1}.
\end{rmk}}

\subsection{Construction of 
\texorpdfstring{$\widetilde{\xi}{}^{\otimes}$}{xi}}

$\empty$

\smallskip

\label{subsect:const-tilde-xi}

We denote by \sym{$\Fin$} the category of finite pointed sets.
Up to isomorphism, the objects of $\Fin$ are the pointed sets
$\langle n\rangle=\{1,\ldots, n\}\cup\{*\}$, for $n\in \N$.
For $1\leq i \leq n$, we denote by 
$\rho^i:\langle n\rangle \to \langle 1\rangle$ 
the unique map such that $(\rho^i)^{-1}(1)=\{i\}$.
Recall that a symmetric monoidal $\infty$-category
is a coCartesian fibration $\mathcal{C}^{\otimes}\to \Fin$
such that the induced functor 
$(\rho^i_!)_i:\mathcal{C}_{\langle n\rangle}
\to \prod_{1\leq i \leq n} \mathcal{C}_{\langle 1\rangle}$
is an equivalence for all $n\geq 0$.
We usually write ``$\mathcal{C}_{\langle n\rangle}$'' 
instead of ``$\mathcal{C}_{\langle n\rangle}^{\otimes}$''
to denote the fiber of $\mathcal{C}^{\otimes}\to \Fin$
at $\langle n\rangle$.
The $\infty$-category $\mathcal{C}_{\langle 1\rangle}$ 
is called the underlying $\infty$-category of $\mathcal{C}^{\otimes}$
and is denoted by $\mathcal{C}$.
Recall also that a monoidal functor is a morphism of coCartesian fibrations between symmetric monoidal $\infty$-categories,
i.e., a functor over $\Fin$ which preserves coCartesian edges.
\symn{$(-)_{\langle n\rangle}$}

We remind the reader that ``monoidal'' always means 
``symmetric monoidal'' in this paper.
We denote by $\CAlg(\CAT_{\infty})$ the 
$\infty$-category of (possibly large) 
monoidal $\infty$-categories and monoidal functors between them. 
The following lemma is well-known.

\begin{lemma}
\label{lem:adj-monoidal-cat}
Let $F^{\otimes}:\mathcal{C}^{\otimes} \to \mathcal{D}^{\otimes}$ 
be a monoidal functor between monoidal $\infty$-categories. 
Then the following conditions are 
equivalent.
\begin{enumerate}

\item[(1)] The underlying functor $F$ admits a right adjoint
$G:\mathcal{D}\to\mathcal{C}$;

\item[(2)] The functor $F^{\otimes}$ admits a right adjoint $G^{\otimes}$ 
making the following triangle commutative
$$\xymatrix{\mathcal{C}^{\otimes}\ar[dr]_-p & & 
\mathcal{D}^{\otimes}\ar[dl]^-q \ar[ll]_-{G^{\otimes}}\\
& \Fin& }$$
with $p$ and $q$ the defining coCartesian fibrations.

\end{enumerate}
Moreover, if these conditions are satisfied, we have the following 
two extra properties.
\begin{enumerate}

\item[(a)] The natural transformations 
$$p\to p\circ G^{\otimes}\circ F^{\otimes}=p 
\qquad \text{and} \qquad
q=q\circ F^{\otimes}\circ G^{\otimes}\to q,$$
induced by the 
unit and the counit of the adjunction 
$(F^{\otimes},G^{\otimes})$, are the identity natural transformations
of $p$ and $q$.

\item[(b)] The functor $G^{\otimes}$ is a
right-lax monoidal functor (i.e., preserves coCartesian edges over 
the arrows $\rho^i:\langle n\rangle \to \langle 1 \rangle$ for 
$1\leq i \leq n$) and its underlying functor 
$G_{\langle 1\rangle}$ is equivalent to $G$.

\end{enumerate}
\end{lemma}

\begin{proof}
This is contained in 
\cite[Propositions 7.3.2.5 \& 7.3.2.6, 
\& Corollary 7.3.2.7]{lurie:higher-algebra}.
We also remark that property (a) is automatic. In fact, 
more generally, every invertible natural transformation of $p$ 
is the identity, and similarly for $q$.
\end{proof}

\begin{cor}
\label{cor:adjunction-on-com-alge}
Let $F^{\otimes}:\mathcal{C}^{\otimes} \to \mathcal{D}^{\otimes}$ 
be a monoidal functor between monoidal $\infty$-categories, and 
assume that $F$ admits a right adjoint $G$. 
Then the induced functor 
$$\CAlg(F):\CAlg(\mathcal{C})
\to \CAlg(\mathcal{D})$$ 
admits also a right adjoint, which is given by $\CAlg(G)$.
\end{cor}

\begin{proof}
Let $p:\mathcal{C}^{\otimes} \to \Fin$ and 
$q:\mathcal{D}^{\otimes} \to \Fin$ be the defining
coCartesian fibrations. Recall that $\CAlg(\mathcal{C})$ 
is the full sub-$\infty$-category of ${\rm Sect}(p)=\Fun(\Fin,\mathcal{C}^{\otimes})\times_{\Fun(\Fin,\,\Fin)}{\id_{\Fin}}$
spanned by those sections of $p$ sending the arrows 
$\rho^i:\langle n\rangle \to \langle 1\rangle$, for $1\leq i\leq n$, 
to coCartesian edges, and similarly for 
$\CAlg(\mathcal{D})$. It follows that 
$F^{\otimes}$ and $G^{\otimes}$ induce functors 
$\CAlg(F)$ and $\CAlg(G)$, and that the 
unit and counit of the adjunction 
$(F^{\otimes},G^{\otimes})$ define natural transformations
$$\id\to  \CAlg(G)\circ \CAlg(F)
\qquad \text{and} \qquad 
\CAlg(F)\circ \CAlg(G) \to \id$$
satisfying the usual identities up to homotopy.
\end{proof}

We now start our construction of $\widetilde{\xi}{}^{\otimes}$.
By Proposition \ref{prop:xi-presheaf}, we have a morphism 
$$\xi^{\otimes}:\FSH^{(\eff,\,\hyp)}_{\tau}(-;\Lambda)^{\otimes}
\to \RigSH^{(\eff,\,\hyp)}_{\tau}((-)^{\rig};\Lambda)^{\otimes}$$
in the $\infty$-category $\Fun(\FSch^{\op},\CAlg(\CAT_{\infty}))$.
The formation of $\infty$-categories of commutative algebras
gives a functor $\CAlg(-):\CAlg(\CAT_{\infty}) \to \CAT_{\infty}$.
Applying this functor to $\xi^{\otimes}$ yields a morphism
$$\CAlg(\xi):
\CAlg(\FSH^{(\eff,\,\hyp)}_{\tau}(-;\Lambda)) \to 
\CAlg(\RigSH^{(\eff,\,\hyp)}_{\tau}((-)^{\rig};\Lambda))$$
in the $\infty$-category $\Fun(\FSch^{\op},\CAT_{\infty})$. 
Applying Lurie's unstraightening construction \cite[\S 3.2]{lurie}
to this morphism, we get a commutative triangle
$$\xymatrix{\Xi_0 \ar[rr]^-F \ar[dr]_-{p_0} & & \Xi_1 \ar[ld]^-{p_1}\\
& \FSch^{\op} &}$$
where $p_0$ and $p_1$ are coCartesian fibrations classified by 
$$\CAlg(\FSH_{\tau}^{(\eff,\,\hyp)}(-;\Lambda))
\qquad \text{and} \qquad 
\CAlg(\RigSH_{\tau}^{(\eff,\,\hyp)}((-)^{\rig};\Lambda)),$$ 
and $F$ is the functor induced by $\CAlg(\xi)$.
By Corollary
\ref{cor:adjunction-on-com-alge}, 
the fibers of $F$ admit right adjoints. More precisely, for 
$\mathcal{S}\in \FSch$, the functor $F_{\mathcal{S}}=
\CAlg(\xi_{\mathcal{S}})$
admits a right adjoint, which is given by 
$\CAlg(\chi_{\mathcal{S}})$. 
(Note that $\chi^{\otimes}_{\mathcal{S}}$ 
is a right-lax monoidal functor.)
Applying \cite[Proposition 7.3.2.6]{lurie:higher-algebra}, 
we deduce that $F$ admits a right adjoint $G$ 
making the following triangle 
$$\xymatrix{\Xi_0  \ar[dr]_-{p_0} & & \ar[ll]_-G \Xi_1 \ar[ld]^-{p_1}\\
& \FSch^{\op} &}$$
commutative and such that, for every $\mathcal{S}\in \FSch$, the
functor $G_{\mathcal{S}}$ is equivalent to $\CAlg(\chi_{\mathcal{S}})$.

We now consider the $\infty$-categories 
${\rm Sect}(p_0)$ and ${\rm Sect}(p_1)$ of sections of 
$p_0$ and $p_1$. The functor $G$ induces a functor 
$G':{\rm Sect}(p_1) \to {\rm Sect}(p_0)$. We have 
an obvious object $\one\in {\rm Sect}(p_1)$, such that 
$\one_{\mathcal{S}}\in \CAlg(\RigSH_{\tau}^{(\eff,\,\hyp)}
(\mathcal{S}^{\rig};\Lambda))$ is the initial algebra
for every $\mathcal{S}\in \FSch$. We set:
$$\mathcal{A}=G'(\one).$$
By construction, $\mathcal{A}$ 
is a section of the coCartesian fibration $p_0$ such that
$\mathcal{A}_{\mathcal{S}}$ is equivalent to 
$\chi_{\mathcal{S}}\Lambda$ considered 
as an object of $\CAlg(\FSH_{\tau}^{(\eff,\,\hyp)}
(\mathcal{S};\Lambda))$.
For a morphism $f:\mathcal{T}\to \mathcal{S}$ of formal schemes, 
the induced morphism $\mathcal{A}_{\mathcal{S}}\to 
\mathcal{A}_{\mathcal{T}}$ in $\Xi_0$ corresponds to a morphism 
$f^*\mathcal{A}_{\mathcal{S}} \to \mathcal{A}_{\mathcal{T}}$.
This is the morphism induced by the natural transformation 
$f^*\circ \chi_{\mathcal{S}} \to \chi_{\mathcal{T}}
\circ f^{\rig,\,*}$ which one obtains by 
adjunction from the equivalence 
$f^{\rig,\,*}\circ \xi_S\simeq \xi_T \circ f^*$.
The following fact, which we record for later use, 
follows easily from this description.

\begin{lemma}
\label{lem:adjointable-cal-A-}
Let $f:\mathcal{T}\to \mathcal{S}$ be a morphism of formal schemes. 
For $f$ to be sent to a
$p_0$-coCartesian edge by $\mathcal{A}$, it suffices that
the commutative square 
$$\xymatrix{\FSH_{\tau}^{(\eff,\,\hyp)}(\mathcal{S};\Lambda) 
\ar[r]^-{\xi_{\mathcal{S}}} 
\ar[d]^-{f^*} & \RigSH_{\tau}^{(\eff,\,\hyp)}
(\mathcal{S}^{\rig};\Lambda) \ar[d]^-{f^{\rig,\,*}}\\
\FSH_{\tau}^{(\eff,\,\hyp)}(\mathcal{T};\Lambda) 
\ar[r]^-{\xi_{\mathcal{T}}} & \RigSH_{\tau}^{(\eff,\,\hyp)}
(\mathcal{T}^{\rig};\Lambda)}$$ 
is right adjointable. This happens when $f$ is smooth.
\end{lemma}

\begin{proof}
Only the last assertion requires a proof. 
If $f$ is smooth, then there is a commutative square
$$\xymatrix{\FSH_{\tau}^{(\eff,\,\hyp)}(\mathcal{T};\Lambda) 
\ar[r]^-{\xi_{\mathcal{T}}} 
\ar[d]^-{f_{\sharp}} & \RigSH_{\tau}^{(\eff,\,\hyp)}
(\mathcal{T}^{\rig};\Lambda) \ar[d]^-{f^{\rig}_{\sharp}}\\
\FSH_{\tau}^{(\eff,\,\hyp)}(\mathcal{S};\Lambda) 
\ar[r]^-{\xi_{\mathcal{S}}} & \RigSH_{\tau}^{(\eff,\,\hyp)}
(\mathcal{S}^{\rig};\Lambda)}$$
by Proposition \ref{prop:xi-presheaf}.
The natural transformation
$f^*\circ \chi_{\mathcal{S}} \to \chi_{\mathcal{T}}\circ 
f^{\rig,\,*}$ deduced from the square of the statement via 
the adjunctions $(\xi_{\mathcal{S}},\chi_{\mathcal{S}})$ 
and $(\xi_{\mathcal{T}},\chi_{\mathcal{T}})$ coincides with the
natural equivalence deduced from the above square 
via the adjunctions
$(\xi_{\mathcal{S}}\circ f_{\sharp},f^*\circ \chi_{\mathcal{S}})$ 
and $(f^{\rig}_{\sharp}\circ \xi_{\mathcal{T}}, 
\chi_{\mathcal{T}}\circ f^{\rig,\,*})$.
\end{proof}

Before going further, we need a small digression
about algebras and modules in general monoidal 
$\infty$-categories. Let $\mathcal{C}^{\otimes}$ 
be a monoidal $\infty$-category and 
$p:\mathcal{C}^{\otimes} \to \Fin$ the defining  
coCartesian fibration. By 
\cite[\S 3.3.3]{lurie:higher-algebra}, 
we may associate to $\mathcal{C}^{\otimes}$ a functor 
\begin{equation}
\label{eqn:mod-C-otimes-to-Fin-CAlg}
f:\Mod(\mathcal{C})^{\otimes}
\to \Fin\times \CAlg(\mathcal{C})
\end{equation}
such that, for each commutative algebra 
$A$ of $\mathcal{C}^{\otimes}$, the induced functor 
\begin{equation}
\label{eqn:mod-C-otimes-to-Fin-CAlg-2}
f_A:\Mod_A(\mathcal{C})^{\otimes}=
\Mod(\mathcal{C})^{\otimes}\times_{\CAlg(\mathcal{C})}\{A\}
\to \Fin
\end{equation}
makes $\Mod_A(\mathcal{C})^{\otimes}$ into an $\infty$-operad.
This is the $\infty$-operad of $A$-modules, which is a monoidal 
$\infty$-category whenever $\mathcal{C}$ admits enough colimits,
and these colimits are compatible with the monoidal structure.
We recall below the construction of the simplicial set 
$\Mod(\mathcal{C})^{\otimes}$ which is a particular case
of \cite[Construction 3.3.3.1]{lurie:higher-algebra}.

\begin{cons}
\label{cons:Mod-C-otimes-}
Recall that a map $\gamma:\langle m\rangle \to \langle n\rangle$
is said to be inert (resp. semi-inert) 
if the induced map $\gamma^{-1}(\{1,\ldots,n\})\to 
\{1,\ldots,n\}$ is a bijection (resp. an injection). 
The map $\gamma$ is said to be null if its image is the base-point
of $\langle n\rangle$.
Let ${\rm K}\subset \Fun(\Delta^1,\Fin)$
be the full subcategory spanned by the semi-inert maps.
We have two obvious functors 
$e_0,e_1:{\rm K} \to \Fin$ induced by the inclusions
$\{0\},\{1\}\subset \Delta^1$. Given $\langle m\rangle\in \Fin$, 
a morphism $\delta$ in the fiber $e_0^{-1}(\langle m\rangle)$
of $e_0$ at $\langle m\rangle$
is said to be inert if the map $e_1(\delta)$, which belongs
to $\Fin$, is inert.

We define a simplicial set $\Mod(\mathcal{C})^{\otimes}$ 
as follows. Giving a map $\Delta^n \to \Mod(\mathcal{C})^{\otimes}$
is equivalent to giving a map $\Delta^n \to \Fin$, and a
functor $\Delta^n\times_{\Fin,\,e_0}{\rm K} \to \mathcal{C}^{\otimes}$
making the triangle
$$\xymatrix{\Delta^n\times_{\Fin,\,e_0}{\rm K} \ar[r] 
\ar[dr]_-{e_1\circ {\rm pr}_{\rm K}} & 
\mathcal{C}^{\otimes} \ar[d]^-p\\
& \Fin}$$
commutative and such that the following condition is satisfied.
For every vertex $\{i\}\subset \Delta^n$, the induced functor
$\{i\}\times_{\Fin,\,e_0}{\rm K}\to \mathcal{C}^{\otimes}$ 
takes an inert map to a $p$-coCartesian morphism.

There is a full inclusion $\Fin\times \Fin \to {\rm K}$,
sending a pair of objects to the null morphism between them,
which is a section to $(e_0,e_1)$.
This induces the functor \eqref{eqn:mod-C-otimes-to-Fin-CAlg}.
That the functor \eqref{eqn:mod-C-otimes-to-Fin-CAlg-2} 
defines an $\infty$-operad is a particular case 
of \cite[Theorem 3.3.3.9]{lurie:higher-algebra}.
According to \cite[Theorem 4.5.3.1]{lurie:higher-algebra},
the functor \eqref{eqn:mod-C-otimes-to-Fin-CAlg}
is a coCartersian fibration when 
$\mathcal{C}$ admits geometric realisations 
which are moreover compatible with the monoidal structure.
In this case, the functor \eqref{eqn:mod-C-otimes-to-Fin-CAlg-2}
is also a coCartesian fibration and thus the 
$\infty$-operad $\Mod_A(\mathcal{C})^{\otimes}$ is a 
monoidal $\infty$-category. (This is also stated 
explicitly in \cite[Theorems 4.5.2.1]{lurie:higher-algebra}.)
\end{cons}

\begin{rmk}
\label{rmk:functoriality-mod-otimes}
It follows from Construction
\ref{cons:Mod-C-otimes-} that $\Mod(-)^{\otimes}$ 
defines a functor from $\CAlg(\CAT_{\infty})$ to 
$\CAT_{\infty}$ endowed with a natural transformation
$f:\Mod(-)^{\otimes} \to \Fin\times \CAlg(-)$.
In fact, Construction \ref{cons:Mod-C-otimes-}
shows more: $\Mod(-)^{\otimes}$ and $f$ naturally extend
to a larger $\infty$-category of monoidal $\infty$-categories
where the morphisms are given by right-lax monoidal functors.
\end{rmk}

Now, we go back to the situation we are interested in.
We start again with our morphism $\xi^{\otimes}$ 
in $\Fun(\FSch^{\op},\CAlg(\CAT_{\infty}))$. Applying 
the functors $\Mod(-)^{\otimes}$ and $\CAlg(-)$, 
we obtain a commutative square in $\Fun(\FSch^{\op},\CAT_{\infty})$:
$$\xymatrix{\Mod(\FSH^{(\eff,\,\hyp)}_{\tau}(-;\Lambda))^{\otimes} 
\ar[rr]^-{\Mod(\xi)^{\otimes}} 
\ar[d]_-{f_0} && \Mod(\RigSH_{\tau}^{(\eff,\,\hyp)}
((-)^{\rig};\Lambda))^{\otimes} \ar[d]^-{f_1}\\
\Fin\times \CAlg(\FSH_{\tau}^{(\eff,\,\hyp)}(-;\Lambda)) 
\ar[rr]^-{\CAlg(\xi)}  && 
\Fin\times \CAlg(\RigSH_{\tau}^{(\eff,\,\hyp)}((-)^{\rig};\Lambda)).\!}$$
Applying Lurie's unstraightening construction \cite[\S 3.2]{lurie}, we get 
a commutative diagram
$$\xymatrix{\mathfrak{M}^{\otimes}_0 \ar[rr]^-{H^{\otimes}} 
\ar[d]_-{q_0} & & \mathfrak{M}^{\otimes}_1 \ar[d]^-{q_1}\\
\Fin\times \Xi_0 \ar[rr]^-F \ar[dr]_-{p_0} & & 
\Fin\times \Xi_1 \ar[dl]^-{p_1} \\
& \Fin\times \FSch^{\op}.\! &}$$
The functors $p_0$, $p_1$, $q_0$, $q_1$, $p_0\circ q_0$
and $p_1\circ q_1$ are coCartesian fibrations. 
Indeed, for $p_0$ and $p_1$, this is by construction.
For the remaining functors, this follows from the Lemma
\ref{lem:cocart-fibration-M-Xi} below
and \cite[Proposition 2.4.2.3(3)]{lurie}.

\begin{lemma}
\label{lem:cocart-fibration-M-Xi}
Let $\mathcal{C}$ be an $\infty$-category 
and $\mathcal{E}^{\otimes}:\mathcal{C} \to \CAlg(\CAT_{\infty})$ a
functor. Consider the commutative triangle 
$$\xymatrix{\mathcal{M}^{\otimes} \ar[rr]^-r \ar[dr] & & 
\Fin\times \mathcal{D} \ar[dl]\\
& \mathcal{C} & }$$
obtained by applying Lurie's unstraightening construction 
\cite[\S 3.2]{lurie}
to the morphism
$$\Mod(\mathcal{E}(-))^{\otimes} \to 
\Fin\times \CAlg(\mathcal{E}(-))$$
in $\Fun(\mathcal{C},\CAT_{\infty})$.
We assume the following conditions:
\begin{itemize}

\item for every $X \in \mathcal{C}$, the $\infty$-category 
$\mathcal{E}(X)$ admits geometric realisations and 
these are compatible with the monoidal structure;

\item for every morphism $f:X \to Y$, the induced functor 
$\mathcal{E}(f)$ commutes with geometric realisations.

\end{itemize}
Then $r$ is a coCartesian fibration.
\end{lemma}

\begin{proof} 
By \cite[Theorem 4.5.3.1]{lurie:higher-algebra}, the morphism
$r_X:\mathcal{M}^{\otimes}_X\to \Fin\times \mathcal{D}_X$ is a coCartesian 
fibration for every $X\in \mathcal{C}$. 
Using \cite[Proposition 2.4.2.11]{lurie}, 
we deduce that $r$ is a locally 
coCartesian fibration. By \cite[Proposition 2.4.2.8]{lurie},
it remains to check that locally $r$-coCartesian morphisms
are stable under composition. Consider a commutative triangle in 
$\Fin\times \mathcal{D}$ that we depict informally as
$$\xymatrix{(\langle n_0\rangle,X_0,R_0) 
\ar[rr]^-{(\gamma_{02},f_{02},\phi_{02})} 
\ar[dr]_-{(\gamma_{01},f_{01},\phi_{01})\,\,\,\,} & & 
(\langle n_2\rangle,X_2,R_2)\\
& (\langle n_1\rangle,X_1,R_1). 
\ar[ur]_-{\,\,\,\,(\gamma_{12},f_{12},\phi_{12})}}$$
Here $X_i$, for $0\leq i \leq 2$,
are objects of $\mathcal{C}$ and 
$f_{ij}:X_i\to X_j$,
for $0\leq i <j\leq 2$,
are morphisms of $\mathcal{C}$,
each $R_i$ is a commutative algebra in 
$\mathcal{E}(X_i)$
and each $\phi_{ij}:\mathcal{E}(f_{ij})(R_i)\to R_j$
is a morphism of commutative algebras in 
$\mathcal{E}(X_j)$,
and the $\gamma_{ij}$'s are maps in $\Fin$.
From this triangle, we deduce a triangle of $\infty$-categories
$$\xymatrix@C=.5pc@R=2.3pc{
\Mod_{R_0}(\mathcal{E}(X_0))_{\langle n_0\rangle}
\ar[rr]^-{(\gamma_{02},f_{02},\phi_{02})_!} 
\ar[dr]_-{(\gamma_{01},f_{01},\phi_{01})_!\,\,\,\,} & & 
\Mod_{R_2}(\mathcal{E}(X_2))_{\langle n_2\rangle}\\
& \Mod_{R_1}(\mathcal{E}(X_1))_{\langle n_1\rangle}
\ar[ur]_-{\,\,\,\,(\gamma_{12},f_{12},\phi_{12})_!}}$$
and we need to show that this triangle commutes up to equivalence.
Using that the $\mathcal{E}(f_{ij})$'s commute with the tensor
product of modules, one reduces easily to the case where $n_0=n_1=n_2=1$
and $\gamma_{ij}$ are the identity maps. 
We are then left to check that 
$$\mathcal{E}(f_{12})(\mathcal{E}(f_{01})(-)
\otimes_{\mathcal{E}(f_{01})(R_0)}R_1)
\otimes_{\mathcal{E}(f_{12})(R_1)}R_2 
\simeq 
\mathcal{E}(f_{02})(-)\otimes_{\mathcal{E}(f_{02})
(R_0)}R_2,$$
which follows again from the fact that 
the $\mathcal{E}(f_{ij})$'s commute with the tensor
product of modules.
\end{proof}

Recall that we have constructed a section 
$\mathcal{A}:\FSch^{\op} \to \Xi_0$ together with a morphism
$F\mathcal{A} \to \one$. Using Lemma 
\ref{lem:cocart-fibration-M-Xi} and 
\cite[Proposition 2.4.2.3(2)]{lurie}, 
we get coCartesian fibrations 
$$\begin{array}{rcl}
\Phi_0 & = & \mathfrak{M}_0^{\otimes}\times_{\Xi_0,\,\one \to \mathcal{A}}
(\Delta^1\times \FSch^{\op}) \to 
\Delta^1\times \Fin\times \FSch^{\op},\\
& \vspace{-.3cm} &\\
\Phi_1 & = & 
\mathfrak{M}_1^{\otimes}\times_{\Xi_1,\,\one \to F
\mathcal{A} \to \one} (\Delta^2\times \FSch^{\op}) \to 
\Delta^2\times \Fin\times \FSch^{\op},
\end{array}$$
and a morphism 
$\Phi_0 \to \Phi_1\times_{\Delta^2}\Delta^{\{0,1\}}$
induced by $H^{\otimes}$.
Let us pause and describe informally what we have constructed.
For $\mathcal{S}\in \FSch$, the coCartesian fibration
$(\Phi_0)_{\mathcal{S}}\to \Delta^1\times \Fin$
is classified by the monoidal functor 
$-\otimes_{\Lambda} \chi\Lambda:
\FSH_{\tau}^{(\eff,\,\hyp)}(\mathcal{S};\Lambda)^{\otimes} \to 
\FSH_{\tau}^{(\eff,\,\hyp)}(\mathcal{S};\chi\Lambda)^{\otimes}$. 
Similarly, the coCartesian fibration
$(\Phi_1)_{\mathcal{S}}\to \Delta^2\times \Fin$
is classified by the commutative triangle 
$$\xymatrix@C=6pc{\RigSH_{\tau}^{(\eff,\,\hyp)}
(\mathcal{S}^{\rig};\Lambda)^{\otimes} 
\ar@{=}[dr]
\ar[r]^-{-\otimes_{\Lambda} \xi\chi\Lambda} 
& \RigSH_{\tau}^{(\eff,\,\hyp)}(\mathcal{S}^{\rig};
\xi\chi\Lambda)^{\otimes} \ar[d]^-{-\otimes_{\xi\chi\Lambda}\Lambda}\\
& \RigSH_{\tau}^{(\eff,\,\hyp)}
(\mathcal{S}^{\rig};\Lambda)^{\otimes}.\!}$$
Finally, applying Lurie's straightening construction 
\cite[\S 3.2]{lurie}, we get 
the following commutative diagram in the $\infty$-category
$\Fun(\FSch^{\op},\CAlg(\CAT_{\infty}))$:
$$\xymatrix@C=4pc{\FSH_{\tau}^{(\eff,\,\hyp)}
(-;\Lambda)^{\otimes} \ar[r]^-{-\otimes_{\Lambda} \chi\Lambda} 
\ar[d]^-{\xi^{\otimes}} & 
\FSH_{\tau}^{(\eff,\,\hyp)}(-;\chi\Lambda)^{\otimes} 
\ar[d]^-{\xi^{\otimes}} & \\
\RigSH_{\tau}^{(\eff,\,\hyp)}(-;\Lambda)^{\otimes} 
\ar[r]^-{-\otimes_{\Lambda} \xi\chi\Lambda} 
\ar@/_2pc/@{=}[rr]_-{\empty} & \RigSH_{\tau}^{(\eff,\,\hyp)}
(-;\xi\chi\Lambda)^{\otimes} 
\ar[r]^-{-\otimes_{\xi\chi\Lambda}\Lambda} & 
\RigSH_{\tau}^{(\eff,\,\hyp)}(-;\Lambda)^{\otimes}.}$$
The morphism $\widetilde{\xi}{}^{\otimes}$ 
is then defined as the composition of
$$\widetilde{\xi}{}^{\otimes}:
\FSH_{\tau}^{(\eff,\,\hyp)}(-;\chi\Lambda)^{\otimes} 
\xrightarrow{\xi^{\otimes}} \RigSH_{\tau}^{(\eff,\,\hyp)}
(-;\xi\chi\Lambda)^{\otimes}
\xrightarrow{-\otimes_{\xi\chi\Lambda}\Lambda}
\RigSH_{\tau}^{(\eff,\,\hyp)}(-;\Lambda)^{\otimes}.$$
\symn{$\widetilde{\xi}^{\otimes}$}

\subsection{Descent, continuity and stalks, II. 
The case of \texorpdfstring{$\chi\Lambda$}{chi}-modules}

$\empty$

\smallskip

\label{subsect:desc--FDA-chi}

We gather here a few basic properties of the 
functor $\FSH^{(\eff,\,\hyp)}_{\tau}(-;\chi\Lambda)^{\otimes}$
and the natural transformation $\widetilde{\xi}{}^{\otimes}$ 
constructed in Subsection
\ref{subsect:const-tilde-xi}.

\begin{prop}
\label{prop:modules-hyperdescent}
The contravariant functor 
$$\mathcal{S}\mapsto \FSH^{(\eff,\,\hyp)}_{\tau}(\mathcal{S};
\chi\Lambda), \quad f\mapsto f^*$$ 
defines a $\tau$-(hyper)sheaf on $\FSch$ 
with values in $\Prl$.
\end{prop}

\begin{proof}
Fix an internal hypercover $\mathcal{U}_{\bullet}$ in the site 
$(\FSch,\tau)$, with $\mathcal{U}_n\to \mathcal{U}_{-1}$ \'etale 
for every $n\in \N$, and which we assume to be truncated in the 
non-hypercomplete case. We need to show that 
$$\FSH_{\tau}^{(\eff,\hyp)}(\mathcal{U}_{\bullet};\chi\Lambda):
\mathbf{\Delta}_+=\mathbf{\Delta}^{\lhd}\to \CAT_{\infty}$$
is a limit diagram. To do so, we use the fact that 
$\FSH_{\tau}^{(\eff,\,\hyp)}(\mathcal{U}_{\bullet};\Lambda)$
is a limit diagram (by Proposition \ref{prop:FDA-hypersheaf}) 
and exhibit a natural transformation
\begin{equation}
\label{eq-prop:modules-hyperdescent}
\FSH_{\tau}^{(\eff,\,\hyp)}(\mathcal{U}_{\bullet};\chi\Lambda)
\to \FSH_{\tau}^{(\eff,\,\hyp)}(\mathcal{U}_{\bullet};\Lambda)
\end{equation}
satisfying the hypotheses of 
\cite[Corollary 5.2.2.37]{lurie:higher-algebra}.
To do so, we start with the obvious natural transformation 
$$-\otimes_{\Lambda}\chi\Lambda:
\FSH_{\tau}^{(\eff,\,\hyp)}(-;\Lambda)\to 
\FSH_{\tau}^{(\eff,\,\hyp)}(-;\chi\Lambda),$$
that we restrict to $\Et/\mathcal{U}_{-1}$, 
and consider the morphism of coCartesian fibrations
$$\xymatrix{\mathcal{F} \ar[rr]^-{F} \ar[dr]_-p & & \mathcal{G} 
\ar[dl]^-q\\
& \Et/\mathcal{U}_{-1} &}$$
associated to this natural transformation by Lurie's 
unstraightening construction \cite[\S 3.2]{lurie}.
Fiberwise, $F$ admits right adjoints. By 
\cite[Proposition 7.3.2.6]{lurie:higher-algebra}, we deduce that 
$F$ admits a right adjoint $G:\mathcal{G}\to \mathcal{F}$
making the triangle 
$$\xymatrix{\mathcal{F} \ar[dr]_-p & & \mathcal{G} \ar[ll]_-{G}
\ar[dl]^-q\\
& \Et/\mathcal{U}_{-1} &}$$
commutative and 
which is fiberwise given by the forgetful functor.
We claim that $G$ is in fact a morphism of coCartesian fibrations, 
i.e., takes a $q$-coCartesian edge to a $p$-coCartesian edge,
and thus determines a natural transformation 
\begin{equation}
\label{eq-prop:modules-hyperdescent-2}
\FSH_{\tau}^{(\eff,\,\hyp)}(-;\chi\Lambda)\to 
\FSH_{\tau}^{(\eff,\,\hyp)}(-;\Lambda)
\end{equation}
on $\Et/\mathcal{U}_{-1}$ 
given objectwise by the forgetful functor.
To prove this, we need to check that the square
$$\xymatrix@C=3pc{\FSH_{\tau}^{(\eff,\,\hyp)}(\mathcal{V};\Lambda) 
\ar[r]^-{-\otimes_{\Lambda}\,\chi\Lambda} \ar[d]^-{e^*}
& \FSH_{\tau}^{(\eff,\,\hyp)}(\mathcal{V};\chi\Lambda)
\ar[d]^-{e^*} \\
\FSH_{\tau}^{(\eff,\,\hyp)}(\mathcal{V}';\Lambda) 
\ar[r]^-{-\otimes_{\Lambda}\,\chi\Lambda} & 
\FSH_{\tau}^{(\eff,\,\hyp)}(\mathcal{V}';\chi\Lambda)}$$
is right adjointable for every map 
$e:\mathcal{V} \to \mathcal{V}'$ in $\Et/\mathcal{U}_{-1}$.
This follows from Lemma
\ref{lem:adjointable-cal-A-} 
which implies that $e^*\chi_{\mathcal{V}}\Lambda
\to \chi_{\mathcal{V}'}\Lambda$ is an equivalence.
That said, we define 
\eqref{eq-prop:modules-hyperdescent}
to be the restriction of 
\eqref{eq-prop:modules-hyperdescent-2}.
That the hypotheses of 
\cite[Lemma 5.2.2.37]{lurie:higher-algebra}
are satisfied is clear:
\begin{itemize}

\item hypothesis (1) of loc.~cit.{} follows from
Proposition \ref{prop:FDA-hypersheaf};

\item hypothesis (2) of loc.~cit.{} follows from 
\cite[Corollary~4.2.3.2]{lurie:higher-algebra};

\item hypothesis (3) of loc.~cit.{} is clear since the 
$\infty$-categories $\FSH_{\tau}^{(\eff,\,\hyp)}(U_n;\Lambda)$ are presentable;

\item hypothesis (4) of loc.~cit.{}, and more generally 
the right adjointability of the squares
$$\xymatrix{\FSH_{\tau}^{(\eff,\,\hyp)}(\mathcal{V};
\chi\Lambda)\ar[r]^-{e^*} 
\ar[d] & \FSH_{\tau}^{(\eff,\,\hyp)}(\mathcal{V}';\chi\Lambda) 
\ar[d]\\
\FSH_{\tau}^{(\eff,\,\hyp)}(\mathcal{V};\Lambda)\ar[r]^-{e^*} 
& \FSH_{\tau}^{(\eff,\,\hyp)}(\mathcal{V}';\Lambda),\!}$$
for $e:\mathcal{V}'\to \mathcal{V}$ in 
$\Et/\mathcal{U}_{-1}$, is clear by 
construction.

\end{itemize}
This completes the proof.
\end{proof}

\begin{lemma}
\label{lem:tilde-xi-otimes-prl}
The natural transformation 
$$\widetilde{\xi}{}^{\otimes}:
\FSH^{(\eff,\,\hyp)}_{\tau}(-;\chi\Lambda)^{\otimes} \to 
\RigSH^{(\eff,\,\hyp)}_{\tau}((-)^{\rig};\Lambda)^{\otimes}$$
is a morphism in $\Fun(\FSch^{\op},\CAlg(\Prl))$.
Moreover, in the following two cases, 
if we restrict this natural transformation to
the subcategory $\mathcal{V}\subset \FSch$, we get a morphism in 
$\Fun(\mathcal{V}^{\op},\CAlg(\Prl_{\omega}))$.
\begin{enumerate}

\item[(1)] We work in the non-hypercomplete case and,
if $\tau$ is the \'etale topology, we assume that 
$\Lambda$ is eventually coconnective. In this case, we
may take $\mathcal{V}$ to be the wide subcategory
of $\FSch$ consisting of quasi-compact morphisms. 

\item[(2)] We work in the hypercomplete case. In this case, 
$\mathcal{V}$ is the subcategory whose objects are those
formal schemes $\mathcal{S}$ such that $\mathcal{S}^{\rig}$
is $(\Lambda,\tau)$-admissible and whose morphisms 
are the quasi-compact and quasi-separated ones. 

\end{enumerate}
\end{lemma}

\begin{proof}
By \cite[Theorem 3.4.4.2]{lurie:higher-algebra}, 
$\FSH_{\tau}^{(\eff,\,\hyp)}(\mathcal{S};\chi\Lambda)^{\otimes}$ is a
presentable monoidal $\infty$-category for every 
$\mathcal{S}\in \FSch$. Moreover, the image of 
$-\otimes_{\Lambda}\chi\Lambda:\FSH_{\tau}^{(\eff,\,\hyp)}
(\mathcal{S};\Lambda) \to \FSH_{\tau}^{(\eff,\,\hyp)}
(\mathcal{S};\chi\Lambda)$ generates 
$\FSH_{\tau}^{(\eff,\,\hyp)}(\mathcal{S};\chi\Lambda)$ 
by colimits. This follows from Proposition
\ref{conservativity-generation} since 
the right adjoint to $-\otimes_{\Lambda}\chi\Lambda$ is conservative
by \cite[Corollary 4.2.3.2]{lurie:higher-algebra}.
By \cite[Corollary 3.4.4.6]{lurie:higher-algebra}, this right adjoint 
also preserves all colimits, which implies that 
$-\otimes_{\Lambda}\chi\Lambda$ preserves compact objects. 
In particular, we see that 
$\FSH^{(\eff,\,\hyp)}_{\tau}
(\mathcal{S};\chi\Lambda)$ is compactly generated 
when $\FSH^{(\eff,\,\hyp)}_{\tau}
(\mathcal{S};\Lambda)$ is. Thus, the second part of the statement 
follows easily from Propositions
\ref{prop:compact-shv-rigsm} and 
\ref{prop:compact-shv-forsm}.
\end{proof}

Our next goal is to prove the continuity property for 
$\FSH^{(\eff,\,\hyp)}_{\tau}(-;\chi\Lambda)$.
\begin{thm}
\label{thm:anstC-FDA-chi}
\ncn{continuity}
Let $(\mathcal{S}_{\alpha})_{\alpha}$ be a cofiltered inverse system 
of quasi-compact and quasi-separated formal schemes 
with affine transition maps, and let $\mathcal{S}=\lim_{\alpha}
\mathcal{S}_{\alpha}$ be the limit of this system. 
We assume one of the following two alternatives. 
\begin{enumerate}

\item[(1)] We work in the non-hypercomplete case. When $\tau$ 
is the \'etale topology, we assume furthermore that $\Lambda$ is 
eventually coconnective.

\item[(2)]
We work in the hypercomplete case, and $\mathcal{S}$ and
$\mathcal{S}^{\rig}$ as well as 
the $\mathcal{S}_{\alpha}$'s and the
$\mathcal{S}^{\rig}_{\alpha}$'s are $(\Lambda,\tau)$-admissible. 
When $\tau$ is the \'etale topology, we assume furthermore that
$\Lambda$ is eventually coconnective or that  
the numbers $\pvcd_{\Lambda}(\mathcal{S}_{\alpha})$ 
and $\pvcd_{\Lambda}(\mathcal{S}^{\rig}_{\alpha})$ 
are bounded independently of $\alpha$.

\end{enumerate}
Then the obvious functor
\begin{equation}
\label{eq-thm:anstC-FDA-chi}
\underset{\alpha}{\colim}\,
\FSH^{(\eff,\,\hyp)}_{\tau}(\mathcal{S}_{\alpha};\chi\Lambda)
\to \FSH^{(\eff,\,\hyp)}_{\tau}(\mathcal{S};\chi\Lambda),
\end{equation}
where the colimit is taken in $\Prl$, is an equivalence.
\end{thm}

\begin{rmk}
\label{rmk:first-alter-cont-chi}
Compared to the analogous statements for rigid analytic 
and formal motives (see Theorem
\ref{thm:anstC} and Proposition
\ref{prop:anstC-FDA}), we have to assume, in the non-hypercomplete 
case, that $\Lambda$ is eventually coconnective when 
$\tau$ is the \'etale topology. This is due to Lemma 
\ref{lem:colimi-chi-}
below, that we were only able to prove under this extra assumption
which insures the compact generation of the $\infty$-categories
of $\chi\Lambda$-modules in formal motives. 
\end{rmk}

We will obtain Theorem
\ref{thm:anstC-FDA-chi} as a consequence of Theorem
\ref{thm:anstC} and Proposition \ref{prop:anstC-FDA}.
To do so, we need some $\infty$-categorical facts. 
We start with the following result, which is well-known 
but for which we couldn't find a reference.

\begin{lemma}
\label{lem:ff-is-colim-preserv-Mod-}
Let $\mathcal{C}^{\otimes}$ be a monoidal $\infty$-category
admitting colimits which are compatible
with the monoidal structure. Then, the forgetful functor 
${\rm ff}:\Mod(\mathcal{C})\to \mathcal{C}$ 
commutes with filtered colimits.
\end{lemma}

\begin{proof}
By \cite[Theorem 4.5.3.1]{lurie:higher-algebra}, 
we have a coCartesian fibration 
$\Mod(\mathcal{C})\to \CAlg(\mathcal{C})$. 
By \cite[Corollary 3.4.4.6(2)]{lurie:higher-algebra}, 
for every $A\in \CAlg(\mathcal{C})$, the $\infty$-category
$\Mod_A(\mathcal{C})$ admits colimits and the forgetful 
functor ${\rm ff}_A:\Mod_A(\mathcal{C})\to \mathcal{C}$
is colimit-preserving. Also, the base change
functor $\Mod_A(\mathcal{C})\to \Mod_B(\mathcal{C})$, 
associated to a morphism $A\to B$ in $\CAlg(\mathcal{C})$,
is colimit-preserving since it admits a right adjoint.
Moreover, by \cite[Corollaries 3.2.3.2 \& 3.2.3.3]{lurie:higher-algebra},
the $\infty$-category $\CAlg(\mathcal{C})$ 
admits colimits and the forgetful functor 
$\CAlg(\mathcal{C}) \to \mathcal{C}$ preserves 
the filtered ones.
Using {\cite[Proposition 4.3.1.5(2) \& 
Corollary 4.3.1.11]{lurie}},
we deduce that $\Mod(\mathcal{C})$ admits 
colimits and that they are computed as follows. 
Let $p:K\to \Mod(\mathcal{C})$ be a diagram and let 
$q:K\to \CAlg(\mathcal{C})$ be the diagram obtained 
by composing with the forgetful functor. 
Let $A_{\infty}\in \CAlg(\mathcal{C})$ be a colimit of 
$q$ and let $p':K \to \Mod_{A_{\infty}}(\mathcal{C})$
be a diagram endowed with a morphism $p\to p'$ 
in $\Mod(\mathcal{C})^K$ given by coCartesian edges. 
(See the beginning of the proof of \cite[Corollary 4.3.1.11]{lurie}.)
Then, the colimit of $p$ is equivalent to the colimit of 
$p'$ computed in $\Mod_{A_{\infty}}(\mathcal{C})$.

Now assume that $K$ is a filtered partially ordered set, 
and let $L$ be the subset of $K\times K$ 
consisting of those pairs $(i,j)$ with $i\leq j$.
We endow $L$ with the induced order. 
Consider the commutative square
$$\xymatrix{K \ar[r]^-p \ar[d] & \Mod(\mathcal{C}) \ar[d]\\
L \ar[r]^-{\tilde{q}} & \CAlg(\mathcal{C}),\!}$$
where the vertical left arrow is the diagonal map given by 
$i\mapsto (i,i)$ and $\tilde{q}$ is the diagram obtained 
by composing $q$ with the map $L\to K$ given by $(i,j)\mapsto j$.
Let $\tilde{p}:L\to \Mod(\mathcal{C})$
be the relative left Kan extension (in the sense of 
\cite[Definition 4.3.2.2]{lurie}). 
Setting $A_i=q(i)$ and $M_i=p(i)$, we have informally
$\tilde{p}(i,j)=A_j\otimes_{A_i}M_i$.
The diagrams $p$ and $\tilde{p}$ have the same colimits,
so it is enough to show that ${\rm ff}(\colim\,\tilde{p})
\simeq \colim\, {\rm ff}\circ \tilde{p}$.
Now, a colimit over $L$ can be computed as a double 
colimit
$$\underset{(i,j)\in L}{\colim}\simeq \underset{i\in K}{\colim}\;
\underset{j\in K_{i/}}{\colim}.$$
Moreover, since the diagram $i\mapsto \colim_{j\in K_{i/}}\,
\tilde{p}(i,-)$ lands in $\Mod_{A_{\infty}}(\mathcal{C})$,
its colimit commutes with ${\rm ff}_{A_{\infty}}$ as mentioned 
above. Thus, it is enough to prove the statement for the 
diagrams $\tilde{p}(i,-):K_{i/}\to \Mod(\mathcal{C})$.
Said differently, we may assume that $p$ takes an
edge of $K$ to a coCartesian edge of the coCartesian 
fibration $\Mod(\mathcal{C})\to \CAlg(\mathcal{C})$.

We may assume that $K$ has an initial object 
$o\in K$. We have a natural transformation 
between the following two functors 
$\Mod_{A_o}(\mathcal{C}) \to \mathcal{C}$.
\begin{enumerate}

\item[(1)] The first one sends $M\in \Mod_{A_o}(\mathcal{C})$ 
to the colimit in $\mathcal{C}$ 
of the diagram $i\mapsto {\rm ff}_{A_i}(A_i\otimes_{A_o}M)$.

\item[(2)] The second one sends $M\in \Mod_{A_o}(\mathcal{C})$
to ${\rm ff}_{A_{\infty}}(A_{\infty}\otimes_{A_o}M)$. 

\end{enumerate}
We want to show that this natural transformation is an 
equivalence. (Together with the description of colimits in $\Mod(\mathcal{C})$ given at the beginning, this would complete the proof.)
To do so, we remark that the two functors above
are colimit-preserving. Using 
\cite[Proposition 4.7.3.14]{lurie:higher-algebra}, 
we reduce to show that this natural transformation is 
an equivalence on $A_o$-modules of the form $A_o\otimes M$, 
with $M\in \mathcal{C}$. In this case, we have to show that the
morphism 
$$\underset{i\in K}{\colim}\; {\rm ff}(A_i\otimes M)
\to {\rm ff}(A_{\infty}\otimes M)$$
is an equivalence. This is clear since 
$\CAlg(\mathcal{C})\to \mathcal{C}$ commutes with filtered colimits.
\end{proof}

Before stating the next $\infty$-categorical result, 
we introduce some notation.
Let $\mathcal{C}$ be an $\infty$-category
and $\mathcal{E}^{\otimes}:\mathcal{C}\to \CAlg(\Prl)$
a functor. Consider the commutative triangle
$$\xymatrix{\mathcal{M}^{\otimes} \ar[rr]^-r \ar[dr]_-q
& & \Fin\times \mathcal{D}
\ar[dl]^-{\id\times p} \\
& \Fin\times \mathcal{C}}$$
obtained by applying Lurie's unstraightening construction 
\cite[\S 3.2]{lurie}
to the functor sending $X\in \mathcal{C}$ 
to the commutative triangle 
$$\xymatrix{\Mod(\mathcal{E}(X))^{\otimes} \ar[rr] \ar[dr] & & \Fin\times \CAlg(\mathcal{E}(X)) \ar[dl]\\
& \Fin. & }$$
By Lemma \ref{lem:cocart-fibration-M-Xi} and 
\cite[Proposition 2.4.2.3(3)]{lurie}, 
the maps $p$, $q$ and $r$ are all coCartesian 
fibrations.
Assume that we are given a section 
$A$ of the coCartesian fibration
$p:\mathcal{D} \to \mathcal{C}$, and consider 
$\mathcal{M}^{\otimes}_A=
\mathcal{M}^{\otimes}\times_{\mathcal{D},\,A}\mathcal{C}$.
The obvious functor $\mathcal{M}^{\otimes}_A\to 
\Fin\times \mathcal{C}$ is a coCartesian fibration.
By Lurie's straightening construction 
\cite[\S 3.2]{lurie}, it determines a functor 
$$\Mod_A(\mathcal{E})^{\otimes}:\mathcal{C} \to \CAlg(\Prl).$$
For proving Theorem
\ref{thm:anstC-FDA-chi}, 
we will use the following general result.

\begin{lemma}
\label{lem:colim-mod-categ}
Assume that $\mathcal{C}$ 
is filtered and set $\mathcal{E}^{\otimes}_{\infty}=
\colim_{\mathcal{C}}\,\mathcal{E}^{\otimes}$. (Here and below, the 
colimit is taken in $\CAlg(\Prl)$.) Let
$\widetilde{A}:\mathcal{C}\to \CAlg(\mathcal{E}_{\infty})$
be the composition of the section $A$ with the obvious functor
$\mathcal{D} \to \CAlg(\mathcal{E}_{\infty})$, and set 
$A_{\infty}=\colim\, \widetilde{A}$.
Then there is an equivalence 
\begin{equation}
\label{eq-lem:colim-mod-categ}
\underset{\mathcal{C}}{\colim}\,\Mod_A(\mathcal{E})^{\otimes} \simeq \Mod_{A_{\infty}}
(\mathcal{E}_{\infty})^{\otimes}.
\end{equation}
\end{lemma}

\begin{proof}
By \cite[Corollary 3.2.3.2]{lurie:higher-algebra}, 
the forgetful functor $\CAlg(\Prl)\to \Prl$ 
detects filtered colimits. Therefore, it is enough to 
prove that 
$$\colim\,\Mod_A(\mathcal{E}) \to \Mod_{A_{\infty}}
(\mathcal{E}_{\infty})$$
is an equivalence, where the colimit is taken in $\Prl$.
By \cite[Corollary 4.5.1.6]{lurie:higher-algebra}, the 
$\infty$-category $\Mod_{A(c)}(\mathcal{E}(c))$ is equivalent 
to the $\infty$-category ${\rm LMod}_{A(c)}(\mathcal{E}(c))$
of left-$A(c)$-modules, for every $c\in \mathcal{C}$, 
and similarly for $\Mod_{A_{\infty}}(\mathcal{E}_{\infty})$.
In fact, \cite[Corollary 4.5.1.6]{lurie:higher-algebra}
shows also that the functor 
$\Mod_A(\mathcal{E}):\mathcal{C}\to \Prl$
is equivalent to the functor 
${\rm LMod}_A(\mathcal{E}):\mathcal{C}\to \Prl$
which is constructed similarly as above. More explicitly, 
one applies Lurie's unstraightening construction 
\cite[\S 3.2]{lurie} to the functor
sending $c\in \mathcal{C}$ to the functor
${\rm LMod}(\mathcal{E}(c)) \to {\rm Alg}(\mathcal{E}(c))$ 
(see \cite[Definition 4.2.1.13 \& 
Example 4.2.1.18]{lurie:higher-algebra}) to get a morphism 
of coCartesian fibrations
$$\xymatrix{\mathcal{M}'^{\otimes} \ar[rr]^-{r'} \ar[dr]_-{q'}
& & \mathcal{D}'
\ar[dl]^-{p'} \\
& \mathcal{C}.\!}$$
Then, the functor ${\rm LMod}_A(\mathcal{E})$ is obtained by 
applying Lurie's straightening construction 
\cite[\S 3.2]{lurie} to the coCartesian fibration 
$\mathcal{M}_A'=\mathcal{M}'\times_{\mathcal{D}',\,A}\mathcal{C}
\to \mathcal{C}$. That said, we are left to show that 
\begin{equation}
\label{eq-lem:colim-mod-categ-2}
\underset{\mathcal{C}}{\colim}\,{\rm LMod}_A(\mathcal{E}) 
\to {\rm LMod}_{A_{\infty}}
(\mathcal{E}_{\infty})
\end{equation}
is an equivalence, where the colimit is taken in $\Prl$.
Using the functor $\widehat{\Theta}:{\rm Pr}^{\rm Alg} 
\to {\rm Pr}^{\rm Mod}$ of
\cite[Construction 4.8.3.24 
\& Notation 4.8.5.10]{lurie:higher-algebra}
and the forgetful functor ${\rm ff}:{\rm Pr}^{\rm Mod} \to \Prl$,
we may rewrite \eqref{eq-lem:colim-mod-categ-2} as
\begin{equation}
\label{eq-lem:colim-mod-categ-4}
\underset{\mathcal{C}}{\colim}\, 
{\rm ff}\circ \widehat{\Theta}(\mathcal{E},A) \to
{\rm ff}\circ \widehat{\Theta}(\mathcal{E}_{\infty},A_{\infty}).
\end{equation}
We give below an informal description of the objects we have just introduced and refer the reader to loc.~cit. for the precise definitions:
\begin{itemize}

\item ${\rm Pr}^{\rm Alg}$ is the $\infty$-category
whose objects are pairs $(\mathcal{X}^{\otimes},R)$ consisting
of a presentable monoidal $\infty$-category 
$\mathcal{X}^{\otimes}$ and an associative algebra $R\in 
{\rm Alg}(\mathcal{X})$;

\item ${\rm Pr}^{\rm Mod}\simeq {\rm LMod}(\Prl)$ 
is the $\infty$-category whose objects 
are pairs $(\mathcal{X}^{\otimes},\mathcal{Y})$ consisting 
of a presentable monoidal $\infty$-category 
$\mathcal{X}^{\otimes}$ and an $\mathcal{X}^{\otimes}$-module
$\mathcal{Y}$ in $\Prlmon$;

\item $\widehat{\Theta}$ sends $(\mathcal{X}^{\otimes},R)$ to 
$(\mathcal{X}^{\otimes},\Mod_R(\mathcal{X}))$ and 
${\rm ff}$ sends $(\mathcal{X}^{\otimes},\mathcal{Y})$
to $\mathcal{Y}$;

\item $(\mathcal{E},A)$ denotes the functor $\mathcal{C}
\to {\rm Pr}^{\rm Alg}$
given informally by $c\mapsto (\mathcal{E}(c),A(c))$.

\end{itemize}
By Lemma
\ref{lem:ff-is-colim-preserv-Mod-}, 
the functor ${\rm ff}$ commutes with filtered colimits. Using
\cite[Theorem 4.8.5.11]{lurie:higher-algebra}
and \cite[Proposition 4.4.2.9]{lurie}, we deduce that 
$\widehat{\Theta}$ commutes also with filtered colimits. 
Since $\colim_{\mathcal{C}}\,(\mathcal{E},A)\simeq 
(\mathcal{E}_{\infty},A_{\infty})$, this proves that 
\eqref{eq-lem:colim-mod-categ-4} is an equivalence.
\end{proof}

Using Proposition \ref{prop:anstC-FDA}, 
Lemma \ref{lem:colim-mod-categ}
and the construction of the functor 
$\FSH_{\tau}^{(\eff,\,\hyp)}(-;\chi\Lambda)$, we see that Theorem 
\ref{thm:anstC-FDA-chi}
is a consequence of the following lemma.

\begin{lemma}
\label{lem:colimi-chi-}
With the notation and assumptions of Theorem
\ref{thm:anstC-FDA-chi}, we have an equivalence
$$\underset{\alpha}{\colim}\,
f_{\alpha}^*\,\chi_{\mathcal{S}_{\alpha}}\Lambda
\to \chi_{\mathcal{S}}\Lambda$$
in $\FSH_{\tau}^{(\eff,\,\hyp)}(\mathcal{S};\Lambda)$,
where $f_{\alpha}:\mathcal{S}\to \mathcal{S}_{\alpha}$ 
is the obvious map.
\end{lemma}

\begin{proof}
Under the assumptions of Theorem 
\ref{thm:anstC-FDA-chi}, 
Theorem \ref{thm:anstC} and Proposition 
\ref{prop:anstC-FDA}
provide us with 
equivalences in $\Prl_{\omega}$
\begin{equation}
\label{eq-lem:colimi-chi-rig-2390i}
\underset{\alpha}{\colim}\,
\RigSH^{(\eff,\,\hyp)}_{\tau}(\mathcal{S}_{\alpha}^{\rig};\Lambda)
\simeq 
\RigSH^{(\eff,\,\hyp)}_{\tau}(\mathcal{S}^{\rig};\Lambda)
\end{equation}
\begin{equation}
\label{eq-lem:colimi-chi-for-2390i}
\text{and} \quad \underset{\alpha}{\colim}\,
\FSH^{(\eff,\,\hyp)}_{\tau}(\mathcal{S}_{\alpha};\Lambda)
\simeq 
\FSH^{(\eff,\,\hyp)}_{\tau}(\mathcal{S};\Lambda),
\end{equation}
where the colimits are also taken in $\Prl_{\omega}$.
(See Propositions 
\ref{prop:compact-shv-rigsm}
and
\ref{prop:compact-shv-forsm}.) 
In particular, the $\infty$-category
$\FSH_{\tau}^{(\eff,\,\hyp)}(\mathcal{S};\Lambda)$ 
is compactly generated and it suffices to show that a compact object 
$M$ in this $\infty$-category 
induces an equivalence
\begin{equation}
\label{eq-lem:colimi-chi-}
\Map_{\FSH_{\tau}^{(\eff,\,\hyp)}(\mathcal{S};\Lambda)}
(M,\underset{\alpha}{\colim}\,
f_{\alpha}^*\,\chi_{\mathcal{S}_{\alpha}}\Lambda)
\to \Map_{\FSH_{\tau}^{(\eff,\,\hyp)}(\mathcal{S};\Lambda)}(M,
\chi_{\mathcal{S}}\Lambda).
\end{equation}
For $\beta\leq \alpha$, we denote by $f_{\beta\alpha}:
\mathcal{S}_{\beta}\to \mathcal{S}_{\alpha}$ the transition 
map in the inverse system
$(\mathcal{S}_{\alpha})_{\alpha}$.
Since $M$ is compact,
there exists an index $\rho$ and a compact object 
$M_{\rho}\in \FSH_{\tau}^{(\eff,\,\hyp)}(\mathcal{S}_{\rho};\Lambda)$ 
such that $M\simeq f^*_{\rho}M_{\rho}$. We have canonical equivalences:
$$\begin{array}{rcl}
\Map_{\FSH_{\tau}^{(\eff,\,\hyp)}(\mathcal{S};\,\Lambda)}
(M,\underset{\alpha}{\colim}\,
f_{\alpha}^*\,\chi_{\mathcal{S}_{\alpha}}\Lambda) & 
\overset{(1)}{\simeq} & 
\underset{\alpha}{\colim}\,\Map_{\FSH_{\tau}^{(\eff,\,\hyp)}
(\mathcal{S};\,\Lambda)}
(M,f_{\alpha}^*\,\chi_{\mathcal{S}_{\alpha}}\Lambda)\\
& \overset{(2)}{\simeq} & 
\underset{\alpha\leq \rho}{\colim}\,
\underset{\beta\leq \alpha}{\colim}\,
\Map_{\FSH_{\tau}^{(\eff,\,\hyp)}(\mathcal{S}_{\beta};\,\Lambda)}
(f_{\beta\rho}^*M_{\rho},f_{\beta\alpha}^*\,
\chi_{\mathcal{S}_{\alpha}}\Lambda)\\
& \overset{(3)}{\simeq}  &
\underset{\beta\leq \rho}{\colim}\,
\Map_{\FSH_{\tau}^{(\eff,\,\hyp)}(\mathcal{S}_{\beta};\,\Lambda)}
(f_{\beta\rho}^*M_{\rho},\chi_{\mathcal{S}_{\beta}}\Lambda)\\
& \overset{(4)}{\simeq}  &
\underset{\beta\leq \rho}{\colim}\,\Map_{\RigSH_{\tau}^{(\eff,\,\hyp)}
(\mathcal{S}^{\rig}_{\beta};\,\Lambda)}
(f_{\beta\rho}^{\rig,\,*}\xi_{\mathcal{S}_{\rho}}M_{\rho},\Lambda)\\
& \overset{(5)}{\simeq}  &\Map_{\RigSH_{\tau}^{(\eff,\,\hyp)}
(\mathcal{S}^{\rig};\,\Lambda)}
(f_{\rho}^{\rig,\,*}\xi_{\mathcal{S}_{\rho}}M_{\rho},\Lambda)\\
& \overset{(6)}{\simeq}  &\Map_{\FSH_{\tau}^{(\eff,\,\hyp)}
(\mathcal{S};\,\Lambda)}(M,\chi_{\mathcal{S}}\Lambda)
\end{array}$$
where
\begin{enumerate}

\item[(1)] follows from the assumption that $M$ is compact,

\item[(2)] follows from the fact that the colimit in
\eqref{eq-lem:colimi-chi-for-2390i}
is taken in $\Prl_{\omega}$,

\item[(3)] follows from the cofinality of the diagonal map 
$\beta \mapsto (\beta\leq \beta)$,

\item[(4)] follows from the adjunction 
$(\xi_{\mathcal{S}_{\beta}},\chi_{\mathcal{S}_{\beta}})$
and the commutation $\xi_{\mathcal{S}_{\beta}} 
f^*_{\beta\rho}\simeq f^{\rig,\,*}_{\beta\rho}
\xi_{\mathcal{S}_{\rho}}$,

\item[(5)] follows from the fact that the colimit in
\eqref{eq-lem:colimi-chi-rig-2390i}
is taken in $\Prl_{\omega}$,

\item[(6)] follows from the commutation 
$f_{\rho}^{\rig,\,*}\xi_{\mathcal{S}_{\rho}}\simeq 
\xi_{\mathcal{S}} f_{\rho}^*$ and the adjunction
$(\xi_{\mathcal{S}},\chi_{\mathcal{S}})$.

\end{enumerate}
It is easy to see that the composition of the above equivalences 
coincide with the map 
\eqref{eq-lem:colimi-chi-}.
\end{proof}

\begin{rmk}
\label{rmk:extension-lemma-colimi-chi-}
Lemma \ref{lem:colimi-chi-} admits a useful extension as follows. 
Keep the notation and assumptions of Theorem
\ref{thm:anstC-FDA-chi}. Let $I$ be the indexing category
of the inverse system $(\mathcal{S}_{\alpha})_{\alpha}$
and let $\alpha\mapsto N_{\alpha}$ be a section of the 
coCartesian fibration associated to the functor 
$I^{\op}\to \CAT_{\infty}$, 
$\alpha \mapsto \RigSH^{(\eff,\,\hyp)}_{\tau}
(\mathcal{S}_{\alpha}^{\rig};\Lambda)$.
Let $N\in \RigSH^{(\eff,\,\hyp)}_{\tau}(\mathcal{S}^{\rig};\Lambda)$ 
be the colimit of the $f^{\rig,\,*}_{\alpha}N_{\alpha}$'s.
Then there is an equivalence 
$$\underset{\alpha}{\colim}\,
f_{\alpha}^*\,\chi_{\mathcal{S}_{\alpha}}N_{\alpha}
\xrightarrow{\sim} \chi_{\mathcal{S}}N$$
in $\FSH_{\tau}^{(\eff,\,\hyp)}(\mathcal{S};\Lambda)$.
This is shown using exactly the same reasoning as in the proof of 
Lemma \ref{lem:colimi-chi-}.
\end{rmk}

We finish this subsection with a computation of the 
stalks of $\FSH^{(\eff,\,\hyp)}_{\tau}(-;\chi\Lambda)$
for the topology $\rig\text{-}\tau$ on $\FSch$.

\begin{thm}
\label{thm:FDAstalks-chi}
Let $\mathcal{S}$ be a formal scheme and let 
$\mathfrak{s}\to \mathcal{S}$ be an algebraic rigid point 
of $\mathcal{S}$. 
Assume one of the following two alternatives.
\begin{enumerate}

\item[(1)] We work in the non-hypercomplete case and, if $\tau$ is the 
\'etale topology, we assume that $\Lambda$ is eventually 
coconnective.

\item[(2)] We work in the hypercomplete case, and  
$\mathcal{S}$ and $\mathcal{S}^{\rig}$ are $(\Lambda,\tau)$-admissible.
When $\tau$ is the \'etale topology, we assume furthermore 
that $\Lambda$ is eventually coconnective or that the numbers 
$\pvcd_{\Lambda}(\mathcal{S}')$, for admissible blowups 
$\mathcal{S}' \to \mathcal{S}$, 
are bounded independently of $\mathcal{S}'$.

\end{enumerate}
Then there is an equivalence of $\infty$-categories
$$\FSH^{(\eff,\,\hyp)}_{\tau}(-;\chi\Lambda)_{\mathfrak{s}}
\simeq \FSH^{(\eff,\,\hyp)}_{\tau}(\mathfrak{s};\chi\Lambda)$$
where the left-hand side is the stalk of 
$\FSH_{\tau}^{(\eff,\,\hyp)}(-;\chi\Lambda)$ 
at $\mathfrak{s}$, i.e., the colimit, taken in $\Prl$,
of the diagram $(\mathfrak{s} \to \mathcal{U}\to \mathcal{S})
\mapsto \FSH_{\tau}^{(\eff,\,\hyp)}(\mathcal{U};\chi\Lambda)$ with 
$\mathcal{U}\in \FRigEt/\mathcal{S}$.
\end{thm}

\begin{proof}
This follows from Theorem \ref{thm:anstC-FDA-chi}.
Indeed, the condition that 
$\mathcal{S}^{\rig}$ is $(\Lambda,\tau)$-admissible implies that 
$\mathfrak{s}$ is $(\Lambda,\tau)$-admissible. Moreover, if 
the numbers $\pvcd_{\Lambda}(\mathcal{S}')$
are bounded independently of $\mathcal{S}'$
for admissible blowups $\mathcal{S}'\to \mathcal{S}$, 
then the same is true
for the numbers $\pvcd_{\Lambda}(\mathcal{U})$ for the saturated 
rig-\'etale neighbourhoods $\mathfrak{s} \to \mathcal{U}\to \mathcal{S}$.
\end{proof}

\subsection{Proof of the main result, I. Fully faithfulness}

$\empty$

\smallskip

\label{subsect:proof-of-main-thm}

Our goal in this subsection is to prove the first part of 
Theorem \ref{thm:main-thm-} concerning the fully faithfulness 
of the functor $\widetilde{\xi}_{\mathcal{S}}$. 
(The second part of this theorem will be proved in the next
subsection.) A key ingredient 
is a projection formula for the functor $\chi_{\mathcal{S}}$
as in the following statement. This projection formula
is also a key ingredient in the proof of the extended 
proper base change theorem for rigid analytic motives, see
Theorem \ref{thm:prop-base} below.

\begin{thm}
\label{thm:proj-form-chi-xi}
\ncn{projection formula}
We work under Assumption 
\ref{assu:for-main-thm}. 
Let $\mathcal{S}$ be a formal scheme and set
$S=\mathcal{S}^{\rig}$.
Then, for $M\in \RigSH^{(\hyp)}_{\tau}(S;\Lambda)$
and $N\in \FSH^{(\hyp)}_{\tau}(\mathcal{S};\Lambda)$, the obvious map
\begin{equation}
\label{eq-thm:proj-form-chi-xi-1}
\chi_{\mathcal{S}}(M)\otimes N\to 
\chi_{\mathcal{S}}(M\otimes \xi_{\mathcal{S}}(N))
\end{equation}
is an equivalence.
\end{thm}

We first prove the following reduction.

\begin{lemma}
\label{lem:alternative-iii-vs-iv-3i}
To prove Theorem \ref{thm:proj-form-chi-xi},
it is enough to consider the alternatives (i), (ii) and (iv) 
of Assumptions \ref{assu:for-main-thm}. Moreover, when working under the
alternative (iv), we may assume the following extra conditions:
\begin{enumerate}

\item[(1)] $\tau$ is the \'etale topology;

\item[(2)] $\Lambda$ is the Eilenberg--Mac Lane spectrum 
associated to the ring $\Z/\ell$, with 
$\ell$ a prime number invertible on $\mathcal{S}$;

\item[(3)] $M$ and $N$ are compact objects.

\end{enumerate}
\end{lemma}

\begin{proof}
We split the proof into two parts.

\paragraph*{Part 1}
\noindent
Here we show that the conclusion of
Theorem \ref{thm:proj-form-chi-xi} holds under (iii) if it  
holds under (iv).

We work under the alternative (iii).
The problem is local on $\mathcal{S}$. Thus, 
we may assume that $\mathcal{S}$ is affine, given as a limit
of a cofiltered inverse system $(\mathcal{S}_{\alpha})_{\alpha}$ of
affine formal schemes such that the $\mathcal{S}_{\alpha}$'s 
and their generic fibers 
$S_{\alpha}=\mathcal{S}^{\rig}_{\alpha}$ are 
$(\Lambda,\tau)$-admissible. By Theorem \ref{thm:anstC} and 
Proposition \ref{prop:anstC-FDA}, we have equivalences
$$\underset{\alpha}{\colim}\,
\RigSH_{\tau}(S_{\alpha};\Lambda)
\simeq 
\RigSH_{\tau}(S;\Lambda)\quad \text{and}\quad
\underset{\alpha}{\colim}\,
\FSH_{\tau}(\mathcal{S}_{\alpha};\Lambda)
\simeq 
\FSH_{\tau}(\mathcal{S};\Lambda)$$
in $\Prl$, and the colimits are taken in $\Prl$.
Using that the tensor product of $\Prl$ commutes
with filtered colimits, we deduce an equivalence
$$\underset{\alpha}{\colim}\,
\left(\RigSH_{\tau}(S_{\alpha};\Lambda) \otimes 
\FSH_{\tau}(\mathcal{S}_{\alpha};\Lambda)\right)
\simeq 
\RigSH_{\tau}(S;\Lambda) \otimes \FSH_{\tau}(\mathcal{S};\Lambda).$$
Since the functors 
$\xi_{\mathcal{S}_{\alpha}}$ and $\chi_{\mathcal{S}_{\alpha}}$
belong to $\Prl$ and are in adjunction, and since 
$\xi_{\mathcal{S}}$ is the colimit of the 
$\xi_{\mathcal{S}_{\alpha}}$'s, we deduce that 
$\chi_{\mathcal{S}}$ is the colimit of the 
$\chi_{\mathcal{S}_{\alpha}}$'s. 
(Here we use Propositions \ref{prop:compact-shv-rigsm} and
\ref{prop:compact-shv-forsm}
to view $\xi_{\mathcal{S}}$ and the $\xi_{\mathcal{S}_{\alpha}}$'s
as functors in $\Prl_{\omega}$ with colimit-preserving right adjoints.)
Considering $\chi_{\mathcal{S}}(-)\otimes (-)$ and 
$\chi_{\mathcal{S}}(-\otimes \xi_{\mathcal{S}}(-))$ as functors from
$\RigSH_{\tau}(S;\Lambda) \otimes \FSH_{\tau}(\mathcal{S};\Lambda)$
to $\FSH_{\tau}(\mathcal{S};\Lambda)$, and similarly 
with ``$\mathcal{S}_{\alpha}$'' instead of ``$\mathcal{S}$'', 
it follows that the natural transformation 
$\chi_{\mathcal{S}}(-)\otimes (-)\to 
\chi_{\mathcal{S}}(-\otimes \xi_{\mathcal{S}}(-))$
is the colimit of the natural transformations
$\chi_{\mathcal{S}_{\alpha}}(-)\otimes (-)\to 
\chi_{\mathcal{S}_{\alpha}}(-\otimes \xi_{\mathcal{S}_{\alpha}}(-))$.
This reduces us to treat the case where $\mathcal{S}$ and 
$S$ are $(\Lambda,\tau)$-admissible. But in this case, we have 
$$\RigSH_{\tau}(S;\Lambda)\simeq 
\RigSH^{\hyp}_{\tau}(S;\Lambda)
\qquad \text{and} \qquad 
\FSH_{\tau}(\mathcal{S};\Lambda)\simeq 
\FSH^{\hyp}_{\tau}(\mathcal{S};\Lambda)$$
by Propositions \ref{prop:automatic-hypercomp-motives} and 
\ref{prop:automatic-hypercomp-motives-algebraic}. 
Therefore, this case is covered by the 
alternative (iv).

\paragraph*{Part 2}
\noindent 
Here we assume that the conclusion of
Theorem \ref{thm:proj-form-chi-xi} holds under (i) and (ii), and we 
show that we may assume conditions (1), (2) and (3) when proving 
Theorem \ref{thm:proj-form-chi-xi} under (iv).

Assume the alternative (iv). If $\tau$ is the Nisnevich 
topology, then there is nothing to prove since 
Theorem \ref{thm:proj-form-chi-xi} holds under (i). 
Thus, we may assume that $\tau$ is the \'etale topology. 
By Propositions \ref{prop:compact-shv-rigsm} and
\ref{prop:compact-shv-forsm}, the $\infty$-categories 
$\RigSH^{\hyp}_{\et}(S;\Lambda)$ and 
$\FSH^{\hyp}_{\et}(\mathcal{S};\Lambda)$
are compactly
generated, and the functor $\chi_{\mathcal{S}}$
commutes with colimits (since its left adjoint 
is compact-preserving). This will be used freely in the 
discussion below.

Let $M_{\Q}=M\otimes \Q$ and $N_{\Q}=N\otimes \Q$
be the rationalisations of $M$ and $N$, and let 
$M_{\rm tor}$ and $N_{\rm tor}$ be the fibers of 
$M\to M_{\Q}$ and $N\to N_{\Q}$.
Since Theorem \ref{thm:proj-form-chi-xi} 
holds under the alternative (ii), we deduce 
that the morphism \eqref{eq-thm:proj-form-chi-xi-1}
becomes an equivalence if we replace $M$ by $M_{\Q}$ or 
$N$ by $N_{\Q}$. Thus, it remains to show that the morphism
\eqref{eq-thm:proj-form-chi-xi-1}
becomes an equivalence if we replace $M$ and $N$ by 
$M_{\rm tor}$ and $N_{\rm tor}$.
Now, $M_{\rm tor}$ is a coproduct of $\ell$-nilpotent objects, 
where $\ell$ varies among the prime numbers 
which are not invertible in $\pi_0\Lambda$, and 
similarly for $N_{\rm tor}$. Moreover, every 
$\ell$-nilpotent object is a colimit of compact $\ell$-nilpotent 
objects. Thus, it is enough to show that the
morphism \eqref{eq-thm:proj-form-chi-xi-1}
is an equivalence when $M$ and $N$ are 
$\ell$-nilpotent compact objects.

By Theorems 
\ref{thm:rigrig},
\ref{thm:rigrig-algebraic} and 
\ref{thm:formal-mot-alg-mot}, 
we have equivalences of $\infty$-categories
$$\Shv_{\et}^{\hyp}(\Et/S;\Lambda)_{\ellnil}
\simeq \RigSH^{\hyp}_{\et}(S;\Lambda)_{\ellnil}
\qquad\text{and} \qquad 
\Shv_{\et}^{\hyp}(\Et/\mathcal{S};\Lambda)_{\ellnil}
\simeq \FSH^{\hyp}_{\et}(\mathcal{S};\Lambda)_{\ellnil}.$$
We denote by $M_0$ and $N_0$
the objects of $\Shv_{\et}^{\hyp}(\Et/S;\Lambda)_{\ellnil}$ 
and $\Shv_{\et}^{\hyp}(\Et/\mathcal{S};\Lambda)_{\ellnil}$
corresponding to $M$ and $N$ by these equivalences.
It is enough to show that 
\begin{equation}
\label{eq-thm:proj-form-chi-xi-1-lemma-reduc}
\chi_{\mathcal{S}}(M_0)\otimes_{\Lambda} N_0\to 
\chi_{\mathcal{S}}(M_0\otimes_{\Lambda} \xi_{\mathcal{S}}(N_0))
\end{equation}
is an equivalence. 
(Here $\xi_{\mathcal{S}}$ is the inverse image functor associated to
the morphism of sites $(\Et/S,\et)\to (\Et/\mathcal{S},\et)$
given by $(-)^{\rig}$, and $\chi_{\mathcal{S}}$ is its right adjoint.)
Since $M_0$ and $N_0$ are compact, 
they are eventually connective. It follows from Lemmas 
\ref{lem:pi-0-Lambda-coh-dim} and 
\ref{lem:Lamba-tau-coh-pointwise}
(and the analogue of the latter for schemes)
that we have equivalences
$$\begin{array}{rcl}
\chi_{\mathcal{S}}(M_0)\otimes_{\Lambda} N_0 & \simeq & 
{\displaystyle 
\lim_r \chi_{\mathcal{S}}(M_0\otimes_{\Lambda} 
\tau_{\leq r}\Lambda)\otimes_{\Lambda} N_0}\\
\chi_{\mathcal{S}}(M_0\otimes_{\Lambda} \xi_{\mathcal{S}}(N_0))
& \simeq & {\displaystyle \lim_r
\chi_{\mathcal{S}}((M_0\otimes_{\Lambda} \tau_{\leq r}\Lambda)
\otimes_{\Lambda} \xi_{\mathcal{S}}(N_0))}.
\end{array}$$
Thus, it is enough to show that 
\eqref{eq-thm:proj-form-chi-xi-1-lemma-reduc}
becomes an equivalence if we replace $M_0$ by 
$M_0\otimes_{\Lambda} \tau_{\leq r}\Lambda$. 
The latter, being a compact object of 
$\Shv_{\et}^{\hyp}(\Et/S;\tau_{\leq r}\Lambda)$,
is eventually connective and coconnective. Thus, if we momentarily 
renounce on having $M_0$ compact, which we do, 
we may assume that $M_0$ is eventually connective and coconnective. 
By an easy induction, we may even assume that $M_0$ is 
in the heart of $\Shv_{\et}^{\hyp}(\Et/S;\Lambda)$ and
that $\ell$ acts by $0$ on $M_0$, i.e., $M_0$ is an 
ordinary \'etale sheaf of $\pi_0\Lambda/\ell$-modules.

Furthermore, we may take $N_0=\Lambda_{\et}(\mathcal{U})/\ell$, 
with $\mathcal{U}$ an 
\'etale formal $\mathcal{S}$-scheme, since the objects of this form
and their desuspensions generate 
$\Shv^{\hyp}_{\et}(\Et/\mathcal{S};\Lambda)$
under colimits.
In this case, we have
$$\begin{array}{rcl}
\chi_{\mathcal{S}}(M_0)\otimes_{\Lambda} N_0 & \simeq &  
\chi_{\mathcal{S}}(M_0)\otimes_{\Z} \Z_{\et}(\mathcal{U})/\ell\\
& \simeq & \chi_{\mathcal{S}}(M_0)\otimes_{\Z/\ell} 
(\Z_{\et}(\mathcal{U})/\ell \oplus \Z_{\et}(\mathcal{U})/\ell[1]),\\
& &\\
\chi_{\mathcal{S}}(M_0\otimes_{\Lambda} \xi_{\mathcal{S}}(N_0)) & \simeq &  
\chi_{\mathcal{S}}(M_0\otimes_{\Z} 
\xi_{\mathcal{S}}(\Z_{\et}(\mathcal{U})/\ell))\\
& \simeq & \chi_{\mathcal{S}}(M_0\otimes_{\Z/\ell} 
\xi_{\mathcal{S}}(\Z_{\et}(\mathcal{U})/\ell \oplus \Z_{\et}(\mathcal{U})/\ell[1])).
\end{array}$$
This shows that we may assume that $\Lambda=\Z/\ell$ as claimed. 
It remains to replace $M_0$ by a compact \'etale sheaf of 
$\Z/\ell$-modules to finish the proof.
\end{proof}

To prove Theorem \ref{thm:proj-form-chi-xi}, 
we need some preliminaries. We start by introducing 
a new $\infty$-category of motives. Let $\mathcal{S}$ be 
a formal scheme and fix a topology $\tau\in \{\Nis,\et\}$.

\begin{dfn}
\label{dfn:arrow-site}
We define the $\infty$-category
$\overline{\FSH}{}^{(\eff,\,\hyp)}_{\tau}(\mathcal{S};\Lambda)$
by repeating Definitions \ref{def:DAeff-form} and 
\ref{dfn:fsh-stable}
while replacing $\FSm/\mathcal{S}$ with the category
$\FRigSm/\mathcal{S}$ of rig-smooth formal $\mathcal{S}$-schemes {(see Definition
\ref{dfn:rig-smooth} and
Remark \ref{rmk:rig-topol-})}.
\symn{$\overline{\FSH}{}^{(\eff,\,\hyp)}$}
\end{dfn}

\begin{rmk}
\label{rmk:two-functors-overline-FSH}
There are functors relating 
$\overline{\FSH}{}_{\tau}^{(\eff,\,\hyp)}(\mathcal{S};\Lambda)$
to other $\infty$-categories of motives considered before.
Below, we set as usual $S=\mathcal{S}^{\rig}$. 
\begin{enumerate}

\item[(1)] The inclusion functor $\iota_{\mathcal{S}}:\FSm/\mathcal{S}
\to \FRigSm/\mathcal{S}$ induces an adjunction
$$\iota_{\mathcal{S}}^*:\FSH^{(\eff,\,\hyp)}_{\tau}(\mathcal{S};\Lambda)
\rightleftarrows 
\overline{\FSH}{}^{(\eff,\,\hyp)}_{\tau}(\mathcal{S};\Lambda):
\iota_{\mathcal{S},\,*}.$$
The functor $\iota_{\mathcal{S},\,*}$
is induced by the restriction functor along $\iota_{\mathcal{S}}$, 
and the functor $\iota_{\mathcal{S}}^*$ is fully faithful
and underlies a monoidal functor.
\symn{$\iota^*, \iota_*$}

\item[(2)] The functor 
$(-)^{\rig}:\FRigSm/\mathcal{S} \to \RigSm/S$
induces an adjunction 
$$\overline{\xi}_{\mathcal{S}}:
\overline{\FSH}{}_{\tau}^{(\eff,\,\hyp)}
(\mathcal{S};\Lambda)
\rightleftarrows
\RigSH_{\tau}^{(\eff,\,\hyp)}(S;\Lambda):
\overline{\chi}_{\mathcal{S}}.$$
By Remark \ref{rmk:another-site-for-rigsh},
$\overline{\xi}_{\mathcal{S}}$ is a localisation functor,
and $\overline{\chi}_{\mathcal{S}}$ is fully faithful 
and identifies the $\infty$-category 
$\RigSH_{\tau}^{(\eff,\,\hyp)}(S;\Lambda)$ 
with the full sub-$\infty$-category of 
$\overline{\FSH}{}_{\tau}^{(\eff,\,\hyp)}
(\mathcal{S};\Lambda)$ spanned by those objects admitting 
$\rig\text{-}\tau$-(hyper)descent.
\symn{$\overline{\xi}$}

\end{enumerate}
Clearly, we have natural equivalences
$\xi_{\mathcal{S}}\simeq \overline{\xi}_{\mathcal{S}}\circ 
\iota^*_{\mathcal{S}}$ and $\chi_{\mathcal{S}}\simeq 
\iota_{\mathcal{S},\,*}\circ \overline{\chi}_{\mathcal{S}}$.
\end{rmk}

We record the following lemma for later use.

\begin{lemma}
\label{lem:iota--monoidal}
The functor $\iota_{\mathcal{S},\,*}$ underlies a monoidal functor
\begin{equation}
\label{eq-lem:iota--monoidal-1}
\iota_{\mathcal{S},\,*}^{\otimes}:
\overline{\FSH}{}^{(\eff,\,\hyp)}_{\tau}(\mathcal{S};\Lambda)^{\otimes}
\to \FSH^{(\eff,\,\hyp)}_{\tau}(\mathcal{S};\Lambda)^{\otimes}
\end{equation}
which belongs to $\CAlg(\Prl)$.
\end{lemma}

\begin{proof}
The functor 
\begin{equation}
\label{eq-lem:iota--monoidal-2}
\iota_{\mathcal{S},\,*}:\PSh(\FRigSm/\mathcal{S};\Lambda)
\to \PSh(\FSm/\mathcal{S};\Lambda)
\end{equation}
underlies a monoidal functor $\iota_{\mathcal{S},\,*}^{\otimes}$
and admits a right adjoint. (Recall that the tensor product 
on presheaves is given objectwise, see Remarks
\ref{rmk:monoid-str-shv-general}
and 
\ref{rmk:Lambda-tau-X}.)
Moreover, it commutes with the $\tau$-(hyper)sheafification 
functor. Indeed, restricting to the small sites 
$(\Et/\mathcal{X},\tau)$, for $\mathcal{X}$ in 
$\FSm/\mathcal{S}$ (resp. $\FRigSm/\mathcal{S}$),
detects $\tau$-(hyper)sheaves and $\tau$-local equivalences
in the hypercomplete and non-hypercomplete cases.
It follows that the functor 
\eqref{eq-lem:iota--monoidal-2} 
induces a left adjoint functor 
\begin{equation}
\label{eq-lem:iota--monoidal-3}
\iota_{\mathcal{S},\,*}:
\Shv_{\tau}^{(\hyp)}(\FRigSm/\mathcal{S};\Lambda)
\to \Shv_{\tau}^{(\hyp)}(\FSm/\mathcal{S};\Lambda)
\end{equation}
underlying a monoidal functor. 
Moreover, for $\mathcal{X}$ a smooth formal $\mathcal{S}$-scheme, 
we have an equivalence 
$\iota_{\mathcal{S},\,*}(\Lambda_{\tau}(\mathcal{X}))\simeq
\Lambda_{\tau}(\mathcal{X})$. Using 
\cite[Proposition 5.5.4.20]{lurie}, 
it follows that \eqref{eq-lem:iota--monoidal-3} 
preserves $\A^1$-local equivalences inducing a left adjoint functor 
\begin{equation}
\label{eq-lem:iota--monoidal-4}
\iota_{\mathcal{S},\,*}:
\overline{\FSH}{}_{\tau}^{\eff,\,(\hyp)}(\mathcal{S};\Lambda)
\to \FSH_{\tau}^{\eff,\,(\hyp)}(\mathcal{S};\Lambda)
\end{equation}
underlying a monoidal functor.
This functor sends $\Tate_{\mathcal{S}}$ to $\Tate_{\mathcal{S}}$, 
and induces a left adjoint functor 
\begin{equation}
\label{eq-lem:iota--monoidal-6}
\iota_{\mathcal{S},\,*}:
\overline{\FSH}{}_{\tau}^{(\hyp)}(\mathcal{S};\Lambda)
\to \FSH_{\tau}^{(\hyp)}(\mathcal{S};\Lambda)
\end{equation}
underlying a monoidal functor.
From the above discussion, we see that 
the functors \eqref{eq-lem:iota--monoidal-4} 
and \eqref{eq-lem:iota--monoidal-6} are right 
adjoint to the functors $\iota_{\mathcal{S}}^*$
in Remark \ref{rmk:two-functors-overline-FSH}(1), 
finishing the proof. 
\end{proof}

\begin{rmk}
\label{rmk:inverse-image-overline-FSH}
There is also an obvious functorial dependence of 
$\overline{\FSH}{}_{\tau}^{(\eff,\,\hyp)}
(\mathcal{S};\Lambda)$ on the 
formal scheme $\mathcal{S}$. 
A morphism of formal schemes 
$f:\mathcal{T}\to \mathcal{S}$ 
induces an inverse image functor 
$$f^*:\overline{\FSH}{}_{\tau}^{(\eff,\,\hyp)}
(\mathcal{S};\Lambda)
\to 
\overline{\FSH}{}_{\tau}^{(\eff,\,\hyp)}
(\mathcal{T};\Lambda)$$
which is a left adjoint and underlies a monoidal functor. 
Moreover, we have natural equivalences
$$f^*\circ \iota_{\mathcal{S}}^*\simeq \iota_{\mathcal{T}}^*\circ f^*
\qquad \text{and}\qquad
f^{\rig,\,*}\circ \overline{\xi}_{\mathcal{S}}\simeq
\overline{\xi}_{\mathcal{T}}\circ f^*.$$
When $f$ is rig-smooth, $f^*$ admits a left adjoint 
$f_{\sharp}$ and there is a natural equivalence
$\overline{\xi}_{\mathcal{S}}\circ f_{\sharp}\simeq
f^{\rig}_{\sharp}\circ \overline{\xi}_{\mathcal{T}}$.
If $f$ is smooth, we also have a natural equivalence
$\iota_{\mathcal{S}}^*\circ f_{\sharp}\simeq f_{\sharp}\circ 
\iota_{\mathcal{T}}^*$.
\end{rmk}

We now state the main technical result needed for 
proving Theorem \ref{thm:proj-form-chi-xi}.

\begin{prop}
\label{prop:tech-for-full}
Let $\mathcal{S}$ be a formal scheme
and set $S=\mathcal{S}^{\rig}$.
Let $M$ and $N$ be objects of $\RigSH^{(\hyp)}_{\tau}(S;\Lambda)$
and $\FSH^{(\hyp)}_{\tau}(\mathcal{S};\Lambda)$ respectively.
We work under one the alternatives (i), (ii) or (iv) 
of Assumption \ref{assu:for-main-thm} and, when working 
under (iv), we assume the conditions (1), (2) and (3) of
Lemma \ref{lem:alternative-iii-vs-iv-3i}.
Then, the obvious morphism
\begin{equation}
\label{eq-prop:tech-for-full}
\overline{\chi}_{\mathcal{S}}(M) \otimes \iota_{\mathcal{S}}^*(N)
\to \overline{\chi}_{\mathcal{S}}(M\otimes \xi_{\mathcal{S}}(N))
\end{equation}
is an equivalence in 
$\overline{\FSH}{}^{(\hyp)}_{\tau}(\mathcal{S};\Lambda)$.
\end{prop}

We first explain how Theorem \ref{thm:proj-form-chi-xi}
follows from Proposition \ref{prop:tech-for-full}.

\begin{proof}[Proof of Theorem \ref{thm:proj-form-chi-xi}]
By Lemma \ref{lem:alternative-iii-vs-iv-3i}, 
we may work under one the alternatives (i), (ii) or (iv) of 
Assumption \ref{assu:for-main-thm}, and assume 
the conditions (1), (2) and (3) of
Lemma \ref{lem:alternative-iii-vs-iv-3i} when working under (iv).
Then, we have a chain of equivalences 
$$\begin{array}{rcl}
\chi_{\mathcal{S}}(M)\otimes N & \overset{(1)}{\simeq} & 
\iota_{\mathcal{S},\,*}(\overline{\chi}_{\mathcal{S}}(M))
\otimes \iota_{\mathcal{S},\,*}(\iota_{\mathcal{S}}^*(N))\\
& \overset{(2)}{\simeq} & \iota_{\mathcal{S},\,*}(\overline{\chi}_{\mathcal{S}}(M)
\otimes \iota_{\mathcal{S}}^*(N))\\
& \overset{(3)}{\simeq} & \iota_{\mathcal{S},\,*}(\overline{\chi}_{\mathcal{S}}(M \otimes \xi_{\mathcal{S}}(N)))\\
& \overset{(4)}{\simeq} & 
\chi_{\mathcal{S}}(M \otimes \xi_{\mathcal{S}}(N))
\end{array}$$
where
\begin{enumerate}

\item[(1)] follows from the equivalence $\chi_{\mathcal{S}}\simeq 
\iota_{\mathcal{S},\,*}\circ \overline{\chi}_{\mathcal{S}}$ 
and the fully faithfulness of 
$\iota_{\mathcal{S}}^*$,

\item[(2)] follows from Lemma \ref{lem:iota--monoidal},

\item[(3)] follows from Proposition 
\ref{prop:tech-for-full},

\item[(4)] follows from the equivalence $\chi_{\mathcal{S}}\simeq 
\iota_{\mathcal{S},\,*}\circ \overline{\chi}_{\mathcal{S}}$.

\end{enumerate}
It is easy to see that the composition of the above equivalences 
coincides with the natural morphism
$\chi_{\mathcal{S}}(M)\otimes N \to \chi_{\mathcal{S}}(M\otimes \xi_{\mathcal{S}}(N))$.
\end{proof}

\begin{proof}[Proof of Proposition \ref{prop:tech-for-full}]
The morphism \eqref{eq-prop:tech-for-full}
is given by the following composition
$$\begin{array}{rcl}
\overline{\chi}_{\mathcal{S}}(M) \otimes \iota_{\mathcal{S}}^*(N)
& \xrightarrow{(1)} & 
\overline{\chi}_{\mathcal{S}}\overline{\xi}_{\mathcal{S}}
(\overline{\chi}_{\mathcal{S}}(M) \otimes \iota_{\mathcal{S}}^*(N))\\
& \overset{(2)}{\simeq} & \overline{\chi}_{\mathcal{S}}
(\overline{\xi}_{\mathcal{S}}\overline{\chi}_{\mathcal{S}}(M) 
\otimes \overline{\xi}_{\mathcal{S}}\iota_{\mathcal{S}}^*(N))\\
& \overset{(3)}{\simeq} & \overline{\chi}_{\mathcal{S}}(M\otimes \xi_{\mathcal{S}}(N))
\end{array}$$
where the equivalence $(2)$ follows from the fact that 
$\overline{\xi}_{\mathcal{S}}$ is monoidal, and the 
equivalence (3) follows from the fact that 
$\overline{\chi}_{\mathcal{S}}$ is fully faithful
and the equivalence $\xi_{\mathcal{S}}\simeq 
\overline{\xi}_{\mathcal{S}}\circ \iota^*_{\mathcal{S}}$.
Thus, to prove the proposition, it remains to show that 
the morphism (1) is an equivalence. This would follows 
if the object 
$E=\overline{\chi}_{\mathcal{S}}(M) \otimes \iota_{\mathcal{S}}^*(N)$ 
belongs to the image of the functor $\overline{\chi}_{\mathcal{S}}$. 
Recall that the latter identifies $\RigSH^{(\hyp)}_{\tau}(S;\Lambda)$
with the full sub-$\infty$-category of 
$\overline{\FSH}{}^{(\hyp)}_{\tau}(\mathcal{S};\Lambda)$
spanned by those objects admitting $\rig\text{-}\tau$-(hyper)descent.
Thus, we need to show that $E$ is local with respect to 
morphisms of the form
\begin{equation}
\label{eq-prop:tech-for-full-11}
\underset{[n]\in\mathbf{\Delta}}{\colim}\,
\M(\mathcal{U}_{\bullet}) \to \M(\mathcal{U}_{-1}),
\end{equation}
and their desuspensions and negative Tate twists, 
where $\mathcal{U}_{\bullet}$ is a $\rig\text{-}\tau$-hypercover 
which we assume to be truncated in the non-hypercomplete case.
(Here $\mathcal{U}_{-1}$ is a rig-smooth formal $\mathcal{S}$-scheme
and $\mathcal{U}_n$, for $n\in \N$, are rig-\'etale over 
$\mathcal{U}_{-1}$.)
Since $M$ and $N$ are general objects of 
$\RigSH^{(\hyp)}_{\tau}(S;\Lambda)$ and 
$\FSH^{(\hyp)}_{\tau}(\mathcal{S};\Lambda)$,
it is enough to show that $E$ is local with respect to 
\eqref{eq-prop:tech-for-full-11}
without worrying about desuspensions and 
negative Tate twists. By a standard argument, the case of a 
$\rig\text{-}\tau$-hypercover $\mathcal{U}$ 
follows if we can treat the cases of a 
$\rig\text{-}\tau$-hypercover $\mathcal{U}'$ 
refining $\mathcal{U}$ and its base change to each of the 
$\mathcal{U}_n$'s. Using the description of 
$\rig\text{-}\tau$-covers given in Remark 
\ref{rmk:rig-topol-} and Proposition
\ref{prop:descr-rig-etale},
we may thus assume that 
$\mathcal{U}_{\bullet}$ satisfies the following,
according to the cases $\tau=\Nis$ 
and $\tau=\et$.
\begin{enumerate}

\item[$(\Nis)$] The morphism of formal simplicial schemes 
$\mathcal{U}_{\bullet}\to \mathcal{U}_{-1}$ (here $\bullet\geq 0$)
factors through an admissible blowup 
$\widetilde{\mathcal{U}}_{-1} \to \mathcal{U}_{-1}$
and the resulting morphism 
$\mathcal{U}_{\bullet}\to 
\widetilde{\mathcal{U}}_{-1}$ is a Nisnevich hypercover of 
$\widetilde{\mathcal{U}}_{-1}$ which is 
truncated in the non-hypercomplete case.

\item[$(\et)$] The morphism of formal simplicial schemes 
$\mathcal{U}_{\bullet}\to \mathcal{U}_{-1}$ (here $\bullet\geq 0$)
factors through an admissible blowup 
$\widetilde{\mathcal{U}}_{-1} \to \mathcal{U}_{-1}$
and the resulting morphism 
$\mathcal{U}_{\bullet}\to \widetilde{\mathcal{U}}_{-1}$
factors as
$$\mathcal{U}_{\bullet} \xrightarrow{(2)} 
\widetilde{\mathcal{U}}_{\bullet}
\xrightarrow{(1)} \widetilde{\mathcal{U}}_{-1}$$
where $(1)$ is a Nisnevich hypercover of 
$\widetilde{\mathcal{U}}_{-1}$ which is 
truncated in the non-hypercomplete case
and (2) is a relative hypercover
for the topology generated by finite rig-\'etale coverings
(in the sense of Definition 
\ref{dfn:finite-rig-etale}(3)) which is also
truncated in the non-hypercomplete case.

\end{enumerate}
We denote by ``\sym{$\rig\fet$}'' the
topology on formal schemes generated by  
finite rig-\'etale coverings.
Since $E$ admits Nisnevich (hyper)descent by construction, 
we see that the result would follow if we can prove the 
following two properties {(where we denote by
$\M:\FRigSm/\mathcal{S}
\to \overline{\FSH}{}^{(\hyp)}_{\tau}(\mathcal{S};\Lambda)$
the ``associated motive'' functor as in Definitions
\ref{dfn:rigsh-stable} and 
\ref{dfn:fsh-stable})}:
\symn{$\M$}
\begin{enumerate}

\item[(A)] $E$ is local with respect to morphisms
$\M(\mathcal{V}) \to \M(\mathcal{U})$, where
$\mathcal{V}\to \mathcal{U}$ is an admissible blowup;

\item[(B)] if $\tau$ is the \'etale topology, then $E$ is local with 
respect to morphisms of the form
\begin{equation}
\label{eq-prop:tech-for-full-111}
\underset{[n]\in\mathbf{\Delta}}{\colim}\,
\M(\mathcal{V}_{\bullet}) \to \M(\mathcal{V}_{-1}),
\end{equation}
where $\mathcal{V}_{\bullet}$ is a hypercover 
for the topology $\rig\fet$,
which we assume to be truncated in the non-hypercomplete case.

\end{enumerate}
We split the rest of the proof into several parts.
In the first part, we prove property (A). In the second part, we 
establish a preliminary fact for proving property (B). 
In the remaining parts,
we prove property (B) assuming one of the alternatives (ii)
or (iv) in Assumption \ref{assu:for-main-thm}.

\paragraph*{Part 1}
\noindent
Here we prove property (A). We start by introducing some notations.
We denote by $f:\mathcal{U}\to \mathcal{S}$ the structural morphism
and by $e:\mathcal{V}\to \mathcal{U}$ the admissible blowup,
and we set $g=f\circ e$.
Since $\M(\mathcal{U})=f_{\sharp}\Lambda$ and 
$\M(\mathcal{V})=g_{\sharp}\Lambda$ 
(see Remark \ref{rmk:inverse-image-overline-FSH}), 
it is enough to show that the obvious morphism
$$\Map_{\overline{\FSH}{}^{(\hyp)}_{\tau}(\mathcal{U};\,\Lambda)}
(\Lambda, f^*E) \to 
\Map_{\overline{\FSH}{}^{(\hyp)}_{\tau}(\mathcal{V};\,\Lambda)}
(\Lambda, g^*E)$$
is an equivalence. This map can be identified with 
$$\Map_{\FSH^{(\hyp)}_{\tau}(\mathcal{U};\,\Lambda)}
(\Lambda, \iota_{\mathcal{U},\,*}f^*E) \to 
\Map_{\FSH^{(\hyp)}_{\tau}(\mathcal{V};\,\Lambda)}
(\Lambda, \iota_{\mathcal{V},\,*}g^*E)$$
which is induced by a morphism 
$\iota_{\mathcal{U},\,*}f^*E\to e_*\iota_{\mathcal{V},\,*}g^*E$
in $\FSH^{(\hyp)}_{\tau}(\mathcal{U};\Lambda)$, and it is enough
to show that the latter is an equivalence.
We have a chain of equivalences
$$\begin{array}{rcl}
\iota_{\mathcal{U},\,*}f^*E & = & 
\iota_{\mathcal{U},\,*} f^*(\overline{\chi}_{\mathcal{S}}(M)
\otimes \iota_{\mathcal{S}}^*(N))\\
& \overset{(1)}{\simeq} & (\iota_{\mathcal{U},\,*}f^*\overline{\chi}_{\mathcal{S}}(M))
\otimes (\iota_{\mathcal{U},\,*}f^*\iota_{\mathcal{S}}^*(N))\\
& \overset{(2)}{\simeq} & 
\chi_{\mathcal{U}}(f^{\rig,\,*}(M))\otimes f^*(N)
\end{array}$$
where (1) follows from 
Lemma \ref{lem:iota--monoidal}
and (2) follows from the natural equivalences 
$$f^*\circ \overline{\chi}_{\mathcal{S}}\simeq 
\overline{\chi}_{\mathcal{U}}\circ f^{\rig,\,*},\quad
\iota_{\mathcal{U},\,*}\circ \overline{\chi}_{\mathcal{U}}
\simeq \chi_{\mathcal{U}}, \quad f^*\circ \iota_{\mathcal{S}}^*
\simeq \iota_{\mathcal{U}}^* \circ f^* \quad \text{and} \quad
\iota_{\mathcal{U},\,*}\circ \iota_{\mathcal{U}}^*\simeq \id.$$
The same applies with ``$\mathcal{V}$'' and ``$g$'' instead of 
``$\mathcal{U}$'' and ``$f$''.
Thus, we are left to show that the morphism
$$\chi_{\mathcal{U}}(f^{\rig,\,*}(M))\otimes f^*(N)
\to e_*(\chi_{\mathcal{V}}(g^{\rig,\,*}(M))\otimes g^*(N))$$
is an equivalence. Since $e_{\sigma}$ is a projective morphism, 
we may use Theorem \ref{thm:formal-mot-alg-mot}
and the projective projection formula for algebraic motives
(see \cite[Th\'eor\`eme 2.3.40]{ayoub-th1}
and Proposition
\ref{prop:base-change-finite-and-projection}(1)
in the rigid analytic setting) to rewrite the above morphism as
$$\chi_{\mathcal{U}}(f^{\rig,\,*}(M))\otimes f^*(N)
\to e_*(\chi_{\mathcal{V}}(g^{\rig,\,*}(M)))\otimes f^*(N).$$
The result follows now from the commutation 
$e_*\circ \chi_{\mathcal{V}}\simeq \chi_{\mathcal{U}}\circ
e^{\rig}_*$ and the fact that $e^{\rig}:\mathcal{V}^{\rig}\to
\mathcal{U}^{\rig}$ is an isomorphism (which implies that 
$e^{\rig}_*\circ g^{\rig,\,*}\simeq f^{\rig,\,*}$). 

\paragraph*{Part 2}
\noindent 
Until the end of the proof, 
$\tau$ will be the \'etale topology. In this part, 
we formulate a property which implies  
property (B) for a fixed hypercover $\mathcal{V}_{\bullet}$;
see property (B$'$) below.

For $n\geq -1$, we denote by $g_n:\mathcal{V}_n\to \mathcal{S}$
and $e_n:\mathcal{V}_n\to \mathcal{V}_{-1}$
the obvious morphisms.
As in the first part, we need to prove that 
$$\Map_{\FSH^{(\hyp)}_{\tau}(\mathcal{V}_{-1};\,\Lambda)}
(\Lambda, \iota_{\mathcal{V}_{-1},\,*}g_{-1}^*E) \to 
\lim_{[n]\in \mathbf{\Delta}}
\Map_{\FSH^{(\hyp)}_{\tau}(\mathcal{V}_n;\,\Lambda)}
(\Lambda, \iota_{\mathcal{V}_n,\,*}g_n^*E)$$
is an equivalence. As explained in the first part, we have an
equivalence
$$\iota_{\mathcal{V}_n,\,*}g_n^*E\simeq \chi_{\mathcal{V}_n}
(g_n^{\rig,\,*}(M))\otimes g_n^*(N),$$
and it is enough to prove that 
$$\chi_{\mathcal{V}_{-1}}
(g_{-1}^{\rig,\,*}(M))\otimes g_{-1}^*(N)
\to \lim_{[n]\in \mathbf{\Delta}}
e_{n,\,*}\left(\chi_{\mathcal{V}_n}(g^{\rig,\,*}_n(M))\otimes 
g_n^*(N)\right)$$
is an equivalence in $\FSH^{(\hyp)}_{\et}(\mathcal{V}_{-1};\Lambda)$.
Since $e_{n,\,\sigma}$ is a finite morphism, we may use Theorem \ref{thm:formal-mot-alg-mot}
and the projective projection formula for algebraic motives
to rewrite the above morphism as 
$$\begin{array}{rcl}
\chi_{\mathcal{V}_{-1}}
(g_{-1}^{\rig,\,*}(M))\otimes g_{-1}^*(N)
& \to & {\displaystyle \lim_{[n]\in \mathbf{\Delta}}
\left(e_{n,\,*}(\chi_{\mathcal{V}_n}(g^{\rig,\,*}_n(M)))\otimes 
g_{-1}^*(N)\right)}\\
& \simeq & {\displaystyle \lim_{[n]\in \mathbf{\Delta}}
\left(\chi_{\mathcal{V}_{-1}}(e^{\rig}_{n,\,*}(g^{\rig,\,*}_n(M)))\otimes 
g_{-1}^*(N)\right)}\\
& \simeq & {\displaystyle \lim_{[n]\in \mathbf{\Delta}}
\left(\chi_{\mathcal{V}_{-1}}(e^{\rig}_{n,\,*}e^{\rig,\,*}_n(g^{\rig,\,*}_{-1}(M)))\otimes g_{-1}^*(N)\right).}
\end{array}$$
Since $g^{\rig,\,*}_{-1}(M)$ belongs to 
$\RigSH^{(\hyp)}_{\et}(\mathcal{V}_{-1}^{\rig};\Lambda)$, 
it admits (hyper)descent with respect to 
$\mathcal{V}_{\bullet}^{\rig}$. 
Using that $\chi_{\mathcal{V}_{-1}}$ is a right adjoint
functor, we deduce that the morphism 
$$\chi_{\mathcal{V}_{-1}}
(g_{-1}^{\rig,\,*}(M))
\to
\lim_{[n]\in \mathbf{\Delta}}
\chi_{\mathcal{V}_{-1}}(e^{\rig}_{n,\,*}e^{\rig,\,*}_n(g^{\rig,\,*}_{-1}(M))).$$
is an equivalence. Thus, we see that property (B) follows from 
the following property:
\begin{enumerate}

\item[(B$'$)] Set
$A^{\bullet}=\chi_{\mathcal{V}_{-1}}(e^{\rig}_{\bullet,\,*}e^{\rig,\,*}_{\bullet}(g^{\rig,\,*}_{-1}(M)))$ and $B=g^*_{-1}(N)$.
Then, the obvious morphism
\begin{equation}
\label{eq-prop:tech-for-full-111-3}
(\lim_{[n]\in \mathbf{\Delta}} A^n)\otimes B \to 
\lim_{[n]\in \mathbf{\Delta}} (A^n\otimes B)
\end{equation}
is an equivalence in $\FSH^{(\hyp)}_{\et}(\mathcal{V}_{-1};\Lambda)$.

\end{enumerate}

\paragraph*{Part 3}
\noindent 
Here we prove property (B) 
assuming that $\pi_0\Lambda$ is a $\Q$-algebra. 

For a formal scheme $\mathcal{X}$, the site 
$(\FRigEt/\mathcal{X}, \rig\fet)$ has zero global and local 
$\Lambda$-cohomological dimensions. 
Indeed, let $\mathcal{F}$ be an ordinary $\rig\fet$-sheaf 
of $\Q$-vector spaces on $\FRigEt/\mathcal{X}$.
For every finite rig-\'etale covering 
$\mathcal{X}''\to \mathcal{X}'$ in $\FRigEt/\mathcal{X}$, 
there is a normalised transfer 
map $\mathcal{F}(\mathcal{X}'') \to \mathcal{F}(\mathcal{X}')$
which is a section to the restriction map. (This map can be 
constructed $\rig\fet$-locally on 
$\mathcal{X}'$, and thus we may assume that 
$\mathcal{X}''$ is isomorphic to a finite coproduct of
copies of $\mathcal{X}'$.) 
Using these 
normalised transfer maps, one can show that 
the {\v C}ech cohomology of $\mathcal{X}$ with values in $\mathcal{F}$
vanishes in degrees $\geq 1$. More precisely, given a finite 
rig-\'etale cover $\mathcal{X}'\to \mathcal{X}$, one can build,
using the normalised transfer maps, a contracting homotopy 
from $\mathcal{F}(\mathcal{X}'_{\bullet})$, where 
$\mathcal{X}'_{\bullet}$ is the {\v C}ech nerve of 
$\mathcal{X}'\to \mathcal{X}$, to the constant simplicial complex
$\mathcal{F}(\mathcal{X})$. We leave the easy details to the reader. 
By Corollary \ref{cor:automatic-hypercomp}, 
it follows that every $\rig\fet$-sheaf of 
$\Lambda$-modules on $\FRigEt/\mathcal{X}$ is automatically a 
$\rig\fet$-(hyper)sheaf. (Indeed, although $\Lambda$ is
not assumed to be eventually coconnective, the condition that 
$\pi_0\Lambda$ is a $\Q$-algebra implies that there exists a morphism
of commutative ring spectra $\Q\to \Lambda$, and thus we may replace
$\Lambda$ by $\Q$ in order to apply 
Corollary \ref{cor:automatic-hypercomp}.)

By the above discussion, it is enough
to check property (B) when $\mathcal{V}_{\bullet}$ 
is the {\v C}ech nerve associated to a finite rig-\'etale covering 
$e_0:\mathcal{V}_0\to \mathcal{V}_{-1}$. 
Moreover, we may assume that the formal $\mathcal{V}_{-1}$-scheme 
$\mathcal{V}_0$ admits an action of a finite group $G$
which is simply transitive on the geometric fibers of 
$e_0^{\rig}:\mathcal{V}_0^{\rig} \to \mathcal{V}_{-1}^{\rig}$. 
The {\v C}ech nerve $\mathcal{V}_{\bullet}$ can be refined
by the following $\rig\fet$-hypercover
\begin{equation}
\label{eq-prop:tech-for-full-111-2}
\xymatrix@C=1.3pc{\cdots \mathcal{V}_0\times G\times G 
\ar[r] \ar@<-.4pc>[r] \ar@<.4pc>[r] &
\mathcal{V}_0\times G \ar@<-.2pc>[r] \ar@<.2pc>[r]& \mathcal{V}_0
\ar[r] & \mathcal{V}_{-1}.}
\end{equation}
Since the latter has the same form when base-changed to 
each $\mathcal{V}_n$, we are left to prove property (B)
with \eqref{eq-prop:tech-for-full-111-2} 
instead of $\mathcal{V}_{\bullet}$. As explained in the second part, 
it suffices to prove property (B$'$) for 
\eqref{eq-prop:tech-for-full-111-2}. 
In this case, the cosimplicial object $A^{\bullet}$
defines an action of $G$ on 
$A^0\in \FSH^{(\hyp)}_{\et}(\mathcal{V}_{-1};\Lambda)$, 
and we may rewrite
\eqref{eq-prop:tech-for-full-111-3} as
$$(A^0)^G \otimes B\to (A^0\otimes B)^G.$$
That this is an equivalence follows from the fact that
taking the ``$G$-invariant subobject'' in a $\Q$-linear 
$\infty$-category
is equivalent to taking the image of the projector 
$|G|^{-1}\sum_{g\in G} g$.

\paragraph*{Part 4}
\noindent
Here we prove property (B)
under the alternative (iv) and assuming conditions (1), (2) and (3)
of Lemma \ref{lem:alternative-iii-vs-iv-3i}.

By Theorems 
\ref{thm:rigrig},
\ref{thm:rigrig-algebraic} and 
\ref{thm:formal-mot-alg-mot}, 
we have equivalences of $\infty$-categories
$$\Shv_{\et}^{\hyp}(\Et/\mathcal{X}^{\rig};\Z/\ell) \simeq
\RigSH^{\hyp}_{\et}(\mathcal{X}^{\rig};\Z/\ell)
\qquad \text{and}\qquad 
\Shv_{\et}^{\hyp}(\Et/\mathcal{X};\Z/\ell) \simeq
\FSH^{\hyp}_{\et}(\mathcal{X};\Z/\ell)$$
for every formal $\mathcal{S}$-scheme $\mathcal{X}$. 
Let $M_0\in \Shv_{\et}^{\hyp}(\Et/\mathcal{X}^{\rig};\Z/\ell)$
and $N_0\in \Shv_{\et}^{\hyp}(\Et/\mathcal{X};\Z/\ell)$
be the \'etale hypersheaves corresponding to $M$ and $N$ 
by these equivalences.
Set $A^{\bullet}_0=
\chi_{\mathcal{V}_{-1}}(e^{\rig}_{\bullet,\,*}e^{\rig,\,*}_{\bullet}(g^{\rig,\,*}_{-1}(M_0)))$ and $B_0=g^*_{-1}(N_0)$.
We need to prove that 
\begin{equation}
\label{eq-prop:tech-for-full-111-22}
(\lim_{[n]\in \mathbf{\Delta}} A^n_0)\otimes B_0 \to 
\lim_{[n]\in \mathbf{\Delta}} (A_0^n\otimes B_0)
\end{equation}
is an equivalence in $\Shv^{\hyp}_{\et}(\Et/\mathcal{V}_{-1};\Z/\ell)$.
We will do this by proving that
\eqref{eq-prop:tech-for-full-111-22} 
induces an equivalence at 
every geometric point $v\to \mathcal{V}_{-1,\,\sigma}$.
Since $M_0$ and $N_0$ are compact, $A^{\bullet}_0$ and 
$A^{\bullet}_0\otimes B_0$ 
are eventually connective and coconnective 
as cosimplicial objects, i.e., uniformly in the cosimplicial degree.
Since the homotopy limit of a cosimplicial object in 
complexes of $\Z/\ell$-modules can be computed using 
the total complex of the associated double complex,  
this implies that
$$(\lim_{[n]\in \mathbf{\Delta}}A^n_0)_v\simeq 
\lim_{[n]\in \mathbf{\Delta}}(A^n_0)_v
\qquad \text{and} \qquad
(\lim_{[n]\in \mathbf{\Delta}}(A^n_0\otimes B))_v\simeq 
\lim_{[n]\in \mathbf{\Delta}}((A^n_0)_v \otimes B_v).$$
Thus, the fiber of 
\eqref{eq-prop:tech-for-full-111-22} 
at $v$ can be identified with the map
$$(\lim_{[n]\in \mathbf{\Delta}} (A^n_0)_v)\otimes (B_0)_v \to 
\lim_{[n]\in \mathbf{\Delta}} ((A_0^n)_v\otimes (B_0)_v).$$
That the latter is an equivalence follows from the fact that 
$(B_0)_v$ is a perfect complex of $\Z/\ell$-modules
(which is a consequence of the assumption that $N$ is compact).
\end{proof}

The method used for proving Theorem \ref{thm:proj-form-chi-xi}
can be also used to prove the following result. 

\begin{prop}
\label{prop:chi-commute-direct-sums}
We work under Assumption 
\ref{assu:for-main-thm}. 
Let $\mathcal{S}$ be a formal scheme and set
$S=\mathcal{S}^{\rig}$. The functor 
$\chi_{\mathcal{S}}:\RigSH^{(\hyp)}_{\tau}(S;\Lambda)
\to \FSH^{(\hyp)}_{\tau}(\mathcal{S};\Lambda)$
preserves colimits. 
\end{prop}

\begin{proof}
This is clear under the alternatives (iii) and (iv)
which imply that $\xi_{\mathcal{S}}$ belongs to $\Prl_{\omega}$
by Propositions 
\ref{prop:compact-shv-rigsm} and
\ref{prop:compact-shv-forsm}.
Thus, it is enough to consider the alternatives (i) and (ii).

By Lemma \ref{lem:iota--monoidal}, 
the functor $\iota_{\mathcal{S},\,*}$ preserves colimits. 
Since $\chi_{\mathcal{S}}=\iota_{\mathcal{S},\,*}
\circ \overline{\chi}_{\mathcal{S}}$, it is enough to show
that the functor $\overline{\chi}_{\mathcal{S}}$ preserves colimits.
The latter is fully faithful with essential image 
the full-subcategory of 
$\overline{\FSH}{}^{(\hyp)}_{\tau}(\mathcal{S};\Lambda)$
spanned by those objects admitting
$\rig\text{-}\tau$-(hyper)descent. Thus, it is 
enough to show that the property of 
admitting $\rig\text{-}\tau$-(hyper)descent
is preserved under colimits.

Let $E:I \to \overline{\FSH}{}^{(\hyp)}_{\tau}(\mathcal{S};\Lambda)$
be a diagram with colimit $E(\infty)$ and such that 
$E(\alpha)$ admits $\rig\text{-}\tau$-(hyper)descent
for every $\alpha\in I$. We need to show that 
$E(\infty)$ admits $\rig\text{-}\tau$-(hyper)descent. 
As in the proof of Proposition 
\ref{prop:tech-for-full}, we reduce to show the following 
two properties:
\begin{enumerate}

\item[(A)] $E(\infty)$ is local with respect to morphisms
$\M(\mathcal{V}) \to \M(\mathcal{U})$, where
$\mathcal{V}\to \mathcal{U}$ is an admissible blowup;

\item[(B)] if $\tau$ is the \'etale topology, then $E(\infty)$ 
is local with respect to morphisms of the form
$$\underset{[n]\in\mathbf{\Delta}}{\colim}\,
\M(\mathcal{V}_{\bullet}) \to \M(\mathcal{V}_{-1}),$$
where $\mathcal{V}_{\bullet}$ is a hypercover 
for the topology $\rig\fet$,
which we assume to be truncated in the non-hypercomplete case.

\end{enumerate}
We split the rest of the proof into two parts. 

\paragraph*{Part 1}
\noindent 
Here we prove property (A). 
We start by introducing some notations.
We denote by $f:\mathcal{U}\to \mathcal{S}$ the structural morphism
and by $e:\mathcal{V}\to \mathcal{U}$ the admissible blowup,
and we set $g=f\circ e$.
We need to show that the obvious morphism 
$$\Map_{\overline{\FSH}{}^{(\hyp)}_{\tau}(\mathcal{U};\,\Lambda)}
(\Lambda, f^*E(\infty)) \to 
\Map_{\overline{\FSH}{}^{(\hyp)}_{\tau}(\mathcal{V};\,\Lambda)}
(\Lambda, g^*E(\infty))$$
is an equivalence. As in the first part of the proof of Proposition 
\ref{prop:tech-for-full}, it is enough to show that 
$$\iota_{\mathcal{U},\,*} f^*E(\infty)
\to e_*\iota_{\mathcal{V},\,*} g^*E(\infty)$$
is an equivalence. Since the objects 
$E(\alpha)$ admit $\rig\text{-}\tau$-(hyper)descent, for $\alpha\in I$, 
we deduce that the morphisms 
$$\iota_{\mathcal{U},\,*} f^*E(\alpha)
\to e_*\iota_{\mathcal{V},\,*} g^*E(\alpha)$$
are equivalences. Since the functors $f^*$, $g^*$, 
$\iota_{\mathcal{U},\,*}$ and 
$\iota_{\mathcal{V},\,*}$ preserve colimits (see 
Lemma \ref{lem:iota--monoidal}), it suffices to show that 
the functor $e_*:\FSH^{(\hyp)}_{\tau}(\mathcal{V};\Lambda)
\to \FSH^{(\hyp)}_{\tau}(\mathcal{U};\Lambda)$ 
preserves colimits. By Theorem 
\ref{thm:formal-mot-alg-mot}, it is equivalent to show that 
the functor 
$e_{\sigma,\,*}:\SH^{(\hyp)}_{\tau}(\mathcal{V}_{\sigma};\Lambda)
\to \SH^{(\hyp)}_{\tau}(\mathcal{U}_{\sigma};\Lambda)$ 
preserves colimits. This follows from the fact that 
$e_{\sigma}$ is projective which implies that 
$e_{\sigma,\,*} \simeq e_{\sigma,\,!}$ admits a right 
adjoint $e_{\sigma}^!$; see \cite[Th\'eor\`eme 1.7.17]{ayoub-th1}.

\paragraph*{Part 2}
\noindent 
Here we prove property (B). 
In particular, we work under the alternative (ii) 
and assume that $\tau$ is the \'etale topology. 

For $n\geq -1$, we denote by $g_n:\mathcal{V}_n\to \mathcal{S}$
and $e_n:\mathcal{V}_n\to \mathcal{V}_{-1}$
the obvious morphisms.
As in the second part of the proof of Proposition 
\ref{prop:tech-for-full}, we need to show that 
$$\iota_{\mathcal{V}_{-1},\,*}g_{-1}^*E(\infty)
\to \lim_{[n]\in \mathbf{\Delta}}
e_{n,\,*}\iota_{\mathcal{V}_n,\,*}g_n^*E(\infty)$$
is an equivalence. Since the objects 
$E(\alpha)$ admit $\riget$-(hyper)descent, for $\alpha\in I$, 
we deduce that the morphisms 
$$\iota_{\mathcal{V}_{-1},\,*}g_{-1}^*E(\alpha)
\to \lim_{[n]\in \mathbf{\Delta}}
e_{n,\,*}\iota_{\mathcal{V}_n,\,*}g_n^*E(\alpha)$$
are equivalences. For $n\geq -1$, the functors 
$g^*_n$, $\iota_{\mathcal{V}_n,\,*}$ and 
$e_{n,\,*}$ commute with colimits. 
(For the second one, we use Lemma \ref{lem:iota--monoidal}
and, for the third one, we use that $e_{n,\,\sigma}$ is 
finite which implies that $e_{n,\,\sigma,\,*}\simeq
e_{n,\,\sigma,\,!}$ admits a 
right adjoint $e_{n,\,\sigma}^!$; see
\cite[Th\'eor\`eme 1.7.17]{ayoub-th1}.)
Therefore, it is enough to show that the 
obvious morphism
\begin{equation}
\label{eq-prop:chi-commute-direct-sums-5}
\underset{\alpha\in I}{\colim}\,\lim_{[n]\in \mathbf{\Delta}}
e_{n,\,*}\iota_{\mathcal{V}_n,\,*}g_n^*E(\alpha)
\to \lim_{[n]\in \mathbf{\Delta}}\underset{\alpha\in I}{\colim}\,
e_{n,\,*}\iota_{\mathcal{V}_n,\,*}g_n^*
E(\alpha)
\end{equation}
is an equivalence. Now, as explained in the third part of 
the proof of Proposition \ref{prop:tech-for-full}, 
we may assume from the beginning that $\mathcal{V}_{\bullet}$
is of the form \eqref{eq-prop:tech-for-full-111-2}.
In this case, the morphism 
\eqref{eq-prop:chi-commute-direct-sums-5}
can be rewritten as follows:
$$\underset{\alpha\in I}{\colim}\,
(e_{0,\,*}\iota_{\mathcal{V}_0,\,*}g_0^*E(\alpha))^G
\to (\underset{\alpha\in I}{\colim}\,
e_{0,\,*}\iota_{\mathcal{V}_0,\,*}g_0^*
E(\alpha))^G.$$
That this is an equivalence follows from the fact that
taking the ``$G$-invariant subobject'' in a $\Q$-linear 
$\infty$-category
is equivalent to taking the image of the projector 
$|G|^{-1}\sum_{g\in G} g$.
\end{proof}

With Theorem \ref{thm:proj-form-chi-xi} and 
Proposition \ref{prop:chi-commute-direct-sums}
at hand, we can prove the first assertion in Theorem
\ref{thm:main-thm-}.

\begin{proof}[Proof of Theorem
\ref{thm:main-thm-}(1)]
We need to show that the unit map
$\id \to \widetilde{\chi}_{\mathcal{S}}\circ 
\widetilde{\xi}_{\mathcal{S}}$
is an equivalence. Clearly, $\widetilde{\xi}_{\mathcal{S}}$ 
preserves colimits and the same is true for 
$\widetilde{\chi}_{\mathcal{S}}$
by Proposition \ref{prop:chi-commute-direct-sums}
combined with \cite[Corollary 3.4.4.6(2)]{lurie:higher-algebra}.
It is thus enough to show that 
$M\to \widetilde{\chi}_{\mathcal{S}}
\widetilde{\xi}_{\mathcal{S}}M$ is an equivalence 
for $M$ varying in a set of objects generating 
$\FSH^{(\hyp)}_{\tau}(\mathcal{S};\chi\Lambda)$ under colimits.
Thus, we may assume that $M$ is a free $\chi_{\mathcal{S}}\Lambda$-module, i.e., that $M\simeq \chi_{\mathcal{S}}(\Lambda)
\otimes N$ for some $N\in \FSH^{(\hyp)}_{\tau}(\mathcal{S};\Lambda)$.
In this case, the unit map coincides with the obvious map
$\chi_{\mathcal{S}}(\Lambda)
\otimes N\to \chi_{\mathcal{S}}\xi_{\mathcal{S}}(N)$
which is an equivalence by Theorem 
\ref{thm:proj-form-chi-xi}.
\end{proof}

\subsection{Proof of the main result, II. Sheafification}

$\empty$

\smallskip

\label{subsect:proof-of-main-thm-second}

Our goal in this subsection is to prove the second part of 
Theorem \ref{thm:main-thm-}. Using \cite[Corollaries 3.2.2.5 \&
3.2.3.2]{lurie}, this is equivalent to proving the following statement.

\begin{thm}
\label{thm:sheafifi-FSH-chi-}
We work under Assumption 
\ref{assu:for-main-thm-2}.
The morphism of $\Prl$-valued presheaves
$$\widetilde{\xi}:
\FSH^{(\hyp)}_{\et}(-;\chi\Lambda)
\to \RigSH^{(\hyp)}_{\et}((-)^{\rig};\Lambda)$$
exhibits $\RigSH^{(\hyp)}_{\et}((-)^{\rig};\Lambda)$ 
as the rig-\'etale sheaf associated to 
$\FSH^{(\hyp)}_{\et}(-;\chi\Lambda)$.
\end{thm}

\begin{rmk}
\label{rmk:sheaf-gives-hypersheaf-main-thm}
In the hypercomplete case, 
Theorem \ref{thm:sheafifi-FSH-chi-},
combined with Theorem
\ref{thm:hyperdesc}, shows
that the \'etale sheafification of 
$\FSH^{\hyp}_{\et}(-;\chi\Lambda)$
is already an \'etale hypersheaf. 
\end{rmk}

\begin{rmk}
\label{rmk:on-sieve-formal-schemes}
Let $\mathcal{S}$ be a formal scheme.
\begin{enumerate}

\item[(1)] Recall that 
a sieve $H\subset \mathcal{S}$
is a sub-presheaf of $\mathcal{S}$ considered as 
a presheaf on $\FSch$. A formal $H$-scheme 
is a formal $\mathcal{S}$-scheme such that the structural
morphism $\mathcal{T}\to \mathcal{S}$ factors through $H$.
We say that $H$ is generated by a family
$(\mathcal{S}_i \to \mathcal{S})_i$ if $H$ is equal to the 
union of the images of the morphisms $\mathcal{S}_i \to \mathcal{S}$
considered as morphisms of presheaves on $\FSch$.
Equivalently, $H$ is the smallest sieve of $\mathcal{S}$ 
such that the $\mathcal{S}_i$'s are formal $H$-schemes.
\ncn{sieves}
\ncn{sieves!generated}

\item[(2)] We say that a sieve $H\subset \mathcal{S}$ 
is a rig-\'etale sieve if 
the inclusion $H\subset \mathcal{S}$ becomes an isomorphism 
after rig-\'etale sheafification.  
Equivalently, $H$ contains the sieve generated by 
a rig-\'etale cover of $\mathcal{S}$. (Of course, this also makes
sense for any other topology.)
\ncn{sieves!rig-\'etale}

\end{enumerate}
\end{rmk}

We will need the following definition.

\begin{dfn}
\label{dfn:potential-good-reduction-cover}
Let $\mathcal{S}$ be a formal scheme.
\begin{enumerate}

\item[(1)] A formal $\mathcal{S}$-scheme $\mathcal{U}$ is said 
to be nearly smooth (resp. \'etale) if, locally on $\mathcal{U}$, 
it is of finite type and 
there exists a finite morphism $\mathcal{U}'\to \mathcal{U}$ 
from a smooth (resp. \'etale) formal $\mathcal{S}$-scheme $\mathcal{U}'$
inducing an isomorphism $\mathcal{U}'^{\rig}\simeq \mathcal{U}^{\rig}$ 
on generic fibers. 
\ncn{formal schemes!nearly smooth}
\ncn{formal schemes!nearly \'etale}

\item[(2)] Let $H\subset \mathcal{S}$ be a sieve.
A formal $\mathcal{S}$-scheme $\mathcal{U}$ is said 
to be $H$-potentially nearly smooth (resp. \'etale) if
$\mathcal{U}\times_{\mathcal{S}}\mathcal{T}$
is nearly smooth (resp. \'etale) over $\mathcal{T}$ for every formal 
$H$-scheme $\mathcal{T}$. If $H$ is generated by a family 
$(\mathcal{S}_i \to \mathcal{S})_i$, it is enough to ask that 
$\mathcal{U}\times_{\mathcal{S}}\mathcal{S}_i$
is nearly smooth (resp. \'etale) over $\mathcal{S}_i$ for every $i$.
\ncn{formal schemes!potentially nearly smooth}
\ncn{formal schemes!potentially nearly \'etale}

\item[(3)] A formal $\mathcal{S}$-scheme $\mathcal{U}$ is said 
to be potentially nearly smooth (resp. \'etale) if it 
is $H$-potentially nearly smooth (resp. \'etale) for some
rig-\'etale sieve $H\subset \mathcal{S}$.

\end{enumerate}
As usual, we say that a morphism of formal 
schemes $\mathcal{T}\to \mathcal{S}$ is ($H$-potentially, 
potentially) nearly smooth if the formal $\mathcal{S}$-scheme
$\mathcal{T}$ is so.
\end{dfn} 

\begin{rmk}
\label{rmk:composition-nearly-smooth-etale}
It follows immediately from the definition that
the class of nearly smooth (resp. \'etale) morphisms
is stable under base change and composition. 
Similarly, the class of potentially nearly smooth 
(resp. \'etale) morphisms
is stable under base change. It follows from  
Proposition \ref{prop:compose-potential-smooth-etale} below that 
the class of potentially nearly \'etale morphisms 
is also stable under composition if we restrict to
quasi-compact and quasi-separated formal schemes. 
However, this is not the case 
for the class of potentially nearly smooth morphisms.
\end{rmk}

We gather a few properties concerning the notion
of (potentially) nearly \'etale morphisms in the 
following proposition.

\begin{prop}
\label{lem:nearly-etale-is-rig-etale}
$\empty$

\begin{enumerate}

\item[(1)] A nearly \'etale morphism of formal schemes 
is rig-\'etale.

\item[(2)] Let $f:\mathcal{T}\to \mathcal{S}$ be a potentially 
nearly \'etale morphism of formal schemes. 
Then, there exists a rig-\'etale cover $g:\mathcal{T}'\to \mathcal{T}$
such that $f\circ g$ is rig-\'etale.\footnote{It is plausible that 
$f$ itself is rig-\'etale, but we didn't strive to prove this 
since we do not need it.}

\item[(3)] A quasi-compact and quasi-separated 
rig-\'etale morphism of formal schemes
is potentially nearly \'etale.

\end{enumerate}
\end{prop}

\begin{proof}
Assertion (1) is clear. Indeed, the notion of rig-\'etaleness 
is local for the rig topology (see Definition 
\ref{dfn:new-etale-rig}(2)) and a finite morphism 
$\mathcal{U}'\to \mathcal{U}$
as in Definition 
\ref{dfn:potential-good-reduction-cover}(1)  
is a rig cover.

We now prove (2). 
By assumption, there is
a rig-\'etale cover $e:\mathcal{S}'\to \mathcal{S}$ such that 
$f':\mathcal{T}'=\mathcal{T}\times_{\mathcal{S}}\mathcal{S}'
\to \mathcal{S}'$ is nearly \'etale. 
By (1), we know that 
$f'$ is rig-\'etale. If follows that 
$e\circ f':\mathcal{T}'\to \mathcal{S}$ is 
also rig-\'etale. Now, remark that $g:\mathcal{T}'\to \mathcal{T}$,
which is a base change of $e$, is a rig-\'etale cover. 
This proves the second assertion.

It remains to prove (3).
Let $f:\mathcal{T} \to \mathcal{S}$ be a quasi-compact and 
quasi-separated rig-\'etale morphism. 
Our goal is to show that $f$ is potentially nearly \'etale. 
The problem is local on $\mathcal{S}$ for the 
rig-\'etale topology and, since
$f$ is quasi-compact and quasi-separated, 
it is local for the Zariski topology on 
$\mathcal{T}$. Thus, we may assume that $\mathcal{S}=\Spf(A)$, with 
$A$ an adic ring of principal ideal type, and 
$\mathcal{T}=\Spf(B)$, with $B$ a rig-\'etale adic 
$A$-algebra such that the zero ideal of $B$ is saturated. 
We fix a generator $\pi\in A$ of an ideal of definition.

We will show that every algebraic geometric rigid point 
$\mathfrak{s}:\Spf(V)\to \mathcal{S}$ admits a rig-\'etale 
neighbourhood $\mathcal{U}_{\mathfrak{s}}$ such that 
$\mathcal{T}\times_{\mathcal{S}}\mathcal{U}_{\mathfrak{s}}$ 
is nearly \'etale over $\mathcal{U}_{\mathfrak{s}}$. 
This suffices to conclude.

Fix $\mathfrak{s}$ as above.
Consider the rig-\'etale $V$-algebra 
$W=V\,\widehat{\otimes}_A\, B/(0)^{\sat}$. 
Arguing as in the proof of Proposition
\ref{prop:descr-rig-etale}, we see that $\Spf(W)$ is the 
completion of a quasi-finite affine flat $V$-scheme, 
necessarily of finite presentation by 
\cite[Chapter 0, Corollary 9.2.8]{fujiwara-kato}.
From Zariski's main theorem \cite[Chapitre IV, 
Th\'eor\`eme 8.12.6]{EGAIV3},
we deduce that $\Spf(W)$ is an open formal subscheme of 
$\Spf(W')$ where $W'$ is a finite flat $V$-algebra.
Moreover since $V[\pi^{-1}]$ is an algebraically 
closed field it follows that $W'[\pi^{-1}]$ is 
a finite direct product of copies of $V[\pi^{-1}]$. 
Replacing $\mathcal{S}$ with a rig-\'etale neighbourhood
of $\mathfrak{s}$ and $\mathcal{T}$ with an open 
covering, we may assume that $W$ is the completion 
of a localisation of $W'$, i.e., there exists
$u\in W'$ which is invertible in 
$W'[\pi^{-1}]$ and such that 
$W$ is the completion of $W'[u^{-1}]$.

Using that $W'[\pi^{-1}]$ is a direct product 
of copies of $V[\pi^{-1}]$, we may find a morphism of 
$V$-algebras
$$V[t]/((t-a_1)\cdots (t-a_r)) \to W',$$
inducing an isomorphism after inverting $\pi$,
where the $a_i$'s belong to $V$ and such that 
two distinct $a_i$'s differ additively 
by an invertible element of 
$V[\pi^{-1}]$.
We may extend this morphism into a presentation
\begin{equation}
\label{eq-lem:nearly-etale-is-rig-etale-7}
V\langle t,s_1,\ldots, s_m\rangle/((t-a_1)\cdots (t-a_r),
\pi^Ns_1-P_1,\ldots, \pi^Ns_m-P_m)^{\sat}\simeq W'
\end{equation}
where $N\in \N$ is large enough and the $P_i$'s are
polynomials in $V[t]$.
The left-hand side of the isomorphism 
\eqref{eq-lem:nearly-etale-is-rig-etale-7}
gives a presentation of the rig-\'etale $V$-algebra $W'$
as in Definition \ref{dfn:new-etale-rig}.
Using Proposition
\ref{prop:on-etale-presentation} and Lemma
\ref{lem:limit-of-formal-schemes}, 
we may assume that the $a_i$'s and the coefficients
of the $P_j$'s belong to the image of the map
\begin{equation}
\label{eq-lem:nearly-etale-is-rig-etale-9}
\underset{\Spf(V) \to \mathcal{U} \to \mathcal{S}}{\colim}
\mathcal{O}(\mathcal{U}) \to V,
\end{equation}
where the colimit is over affine rig-\'etale
neighbourhoods of $\mathfrak{s}$ in $\mathcal{S}$.
Similarly, we may assume that $u\in W'$ 
is the image of a polynomial $Q\in A[t,s_1,\ldots, s_m]$
with coefficients in the image of 
\eqref{eq-lem:nearly-etale-is-rig-etale-9}.
Thus, we may find a rig-\'etale neighbourhood 
$\mathcal{U}_{\mathfrak{s}}=\Spf(A_{\mathfrak{s}})$ 
of $\mathfrak{s}$ and lifts $\widetilde{a}_i$'s,
$\widetilde{P}_j$'s and $\widetilde{Q}$ 
to $A_{\mathfrak{s}}$
of the $a_i$'s, $P_j$'s and $Q$. 
We then set
$$C'_{\mathfrak{s}}=A_{\mathfrak{s}}\langle t,s_1,\ldots, s_m\rangle/((t-\widetilde{a}_1)\cdots (t-\widetilde{a}_r),
\pi^Ns_1-\widetilde{P}_1,\ldots, \pi^Ns_m-\widetilde{P}_m)^{\sat}$$
$$\text{and} \qquad C_{\mathfrak{s}}=C'_{\mathfrak{s}}\langle v\rangle/(v\cdot \widetilde{Q}-1).$$
Refining $\mathcal{U}_{\mathfrak{s}}$, we may assume that 
two $\widetilde{a}_i$'s differ by an invertible element of
$A_{\mathfrak{s}}[\pi^{-1}]$. This insures that 
$C'_{\mathfrak{s}}$ is a rig-\'etale $A_{\mathfrak{s}}$-algebra. 
By construction, we have an isomorphism
$$V\,\widehat{\otimes}_{A_{\mathfrak{s}}}\,C_{\mathfrak{s}}/(0)^{\sat}
\simeq W \simeq V\,\widehat{\otimes}_A\,B/(0)^{\sat}.$$
Using Corollary 
\ref{cor:proj-limi-cal-E-A}, we may refine $\mathcal{U}_{\mathfrak{s}}$
and assume that 
$$C_{\mathfrak{s}}\simeq 
A_{\mathfrak{s}}\,\widehat{\otimes}_A\,B/(0)^{\sat}.$$
Therefore, to conclude, it is enough to see that 
$\Spf(C'_{\mathfrak{s}})$ is nearly \'etale over 
$\Spf(A_{\mathfrak{s}})$ 
for $\mathcal{U}_{\mathfrak{s}}$ sufficiently
small. After refining 
$\mathcal{U}_{\mathfrak{s}}$ if necessary, we may 
assume that the classes of the $\widetilde{P}_i$'s in the ring
$A_{\mathfrak{s}}[t]/((t-\widetilde{a}_1)\cdots (t-\widetilde{a}_r))$,
divided by $\pi^N$, are algebraic over this ring. (Indeed, the $P_i$'s 
satisfy the analogous property.) In this case, the claim
is clear since the normalisation of $C'_{\mathfrak{s}}$ in 
$C'_{\mathfrak{s}}[\pi^{-1}]$ is then a finite 
direct product of copies of the normalisation of
$A_{\mathfrak{s}}$ in $A_{\mathfrak{s}}[\pi^{-1}]$.
\end{proof}

\begin{prop}
\label{prop:compose-potential-smooth-etale}
Let $\mathcal{T} \to \mathcal{S}$ be a quasi-compact and
quasi-separated potentially nearly \'etale morphism
of formal schemes. Let $\mathcal{V}$ be a
potentially nearly smooth formal $\mathcal{T}$-scheme. 
Then $\mathcal{V}$ is also potentially nearly smooth
as a formal $\mathcal{S}$-scheme.
\end{prop}

\begin{proof}
The problem is local on $\mathcal{S}$ for the rig-\'etale topology.
Thus, we may assume that $\mathcal{S}$ and $\mathcal{T}$ 
are quasi-compact and quasi-separated, and that the morphism 
$\mathcal{T}\to \mathcal{S}$ is nearly \'etale. 
The problem is also local on $\mathcal{T}$. Thus, 
we may assume that there is a finite morphism 
$\mathcal{T}_1 \to \mathcal{T}$ from an \'etale 
formal $\mathcal{S}$-scheme
$\mathcal{T}_1$ inducing an isomorphism on generic fibers.
It is clearly enough to show that the formal $\mathcal{S}$-scheme
$\mathcal{T}_1\times_{\mathcal{T}}\mathcal{V}$ is potentially
nearly smooth over $\mathcal{S}$. Thus, we may replace 
$\mathcal{T}$ with $\mathcal{T}_1$ and $\mathcal{V}$ with
$\mathcal{T}_1\times_{\mathcal{T}}\mathcal{V}$, and assume that 
$\mathcal{T}\to \mathcal{S}$ is \'etale. 
Let $\mathcal{T}'\to \mathcal{T}$ be a rig-\'etale cover such that 
$\mathcal{V}\times_{\mathcal{T}}\mathcal{T}'$ is nearly smooth
over $\mathcal{T}'$. By Lemma 
\ref{lem:refining-etale-cover-over-etale-spaces} below, 
there is a rig-\'etale cover $\mathcal{S}'\to \mathcal{S}$ and 
and a morphism of formal $\mathcal{T}$-schemes 
$\mathcal{T}\times_{\mathcal{S}}\mathcal{S}'\to \mathcal{T}'$.
We claim that the formal $\mathcal{S}'$-scheme
$\mathcal{V}\times_{\mathcal{S}}\mathcal{S}'$
is nearly smooth. Indeed, we have an isomorphism
$\mathcal{V}\times_{\mathcal{S}}\mathcal{S}'\simeq
\mathcal{V}\times_{\mathcal{T}}(\mathcal{T}\times_{\mathcal{S}}\mathcal{S}')$
and the formal $\mathcal{T}\times_{\mathcal{S}}\mathcal{S}'$-scheme
$\mathcal{V}\times_{\mathcal{T}}(\mathcal{T}\times_{\mathcal{S}}\mathcal{S}')$ is nearly smooth since it is a base change of 
the formal $\mathcal{T}'$-scheme $\mathcal{V}\times_{\mathcal{T}}\mathcal{T}'$. The structural morphism of 
the formal $\mathcal{S}'$-scheme
$\mathcal{V}\times_{\mathcal{S}}\mathcal{S}'$
is thus the composition of two nearly smooth morphisms
$$\mathcal{V}\times_{\mathcal{T}}(\mathcal{T}\times_{\mathcal{S}}\mathcal{S}')
\to \mathcal{T}\times_{\mathcal{S}}\mathcal{S}' \to \mathcal{S}'.$$
This finishes the proof since nearly smooth morphisms are 
preserved under composition.
\end{proof}

\begin{lemma}
\label{lem:refining-etale-cover-over-etale-spaces}
Let $\mathcal{T} \to \mathcal{S}$ be a quasi-compact and quasi-separated
\'etale morphism of formal schemes, and let $\mathcal{T}'\to \mathcal{T}$ 
be a rig-\'etale cover. Then there exists
a rig-\'etale cover $\mathcal{S}'\to \mathcal{S}$ 
and a morphism of $\mathcal{T}$-schemes
$\mathcal{T}\times_{\mathcal{S}}\mathcal{S}'\to \mathcal{T}'$.
\end{lemma}

\begin{proof}
This is proven in the same manner as Corollary
\ref{cor:refining-etale-cover-open}. Given an algebraic 
geometric rigid point 
$\mathfrak{s} \to \mathcal{S}$, we consider 
$\mathfrak{t}=\mathfrak{s}\times_{\mathcal{S}}\mathcal{T}$.
This is a quasi-compact and quasi-separated 
\'etale formal $\mathfrak{s}$-scheme.
Thus $\mathfrak{t}$ is a disjoint union of quasi-compact 
open formal subschemes of $\mathfrak{s}$. 
In particular, the morphism 
$\mathfrak{t} \to \mathcal{T}$ factors through $\mathcal{T}'$.
We then use Corollary 
\ref{cor:limit-etale-topi-rig-an} and Lemma
\ref{lem:limit-of-formal-schemes}
to conclude.
\end{proof}

\begin{dfn}
\label{dfn:good-sieve-}
Let $\mathcal{S}$ be a formal scheme. 
\begin{enumerate}

\item[(1)] Let $K\subset \mathcal{S}$ be a sieve.
A rig-\'etale sieve $H\subset \mathcal{S}$
is said to be $K$-potentially nearly \'etale 
if it can be generated by a family 
$(\mathcal{S}_i \to \mathcal{S})_i$
consisting of rig-\'etale morphisms which are
$K$-potentially nearly \'etale.
\ncn{sieves!potentially nearly \'etale}

\item[(2)] A rig-\'etale sieve $H\subset \mathcal{S}$
is said to be potentially nearly \'etale if it is $K$-potentially 
nearly \'etale for some rig-\'etale sieve $K\subset \mathcal{S}$.

\end{enumerate}
\end{dfn}

\begin{cor}
\label{cor:existence-of-potent-good-red-covers}
Let $\mathcal{S}$ be a quasi-compact and quasi-separated formal scheme. 
Let $H\subset \mathcal{S}$ be a rig-\'etale sieve. Then, we may refine
$H$ by a rig-\'etale sieve which is potentially 
nearly \'etale.
\end{cor}

\begin{proof}
After refinement, we may assume that $H$ is generated 
by a rig-\'etale cover
$(\mathcal{S}_i \to \mathcal{S})_{i\in I}$ where $I$ is finite and 
every $\mathcal{S}_i$ is a quasi-compact and quasi-separated
rig-\'etale formal $\mathcal{S}$-scheme. 
By Proposition \ref{lem:nearly-etale-is-rig-etale}(3), 
each $\mathcal{S}_i$ is $K_i$-potentially 
nearly \'etale over $\mathcal{S}$ for some rig-\'etale 
sieve $K_i\subset \mathcal{S}$. It follows that $H$ 
is $K$-potentially nearly \'etale, with $K=\cap_i K_i$ which is
a rig-\'etale sieve since $I$ is finite.
\end{proof}

\begin{nota}
\label{nota:xi-tilde-H-}
$\empty$

\begin{enumerate}

\item[(1)] Given a presheaf of sets $H$ on $\FSch$, we denote by 
$\FSH^{(\hyp)}_{\et}(H;\chi\Lambda)$ the object of $\Prl$ obtained
by evaluating on $H$ the right Kan extension of 
$\FSH^{(\hyp)}_{\et}(-;\chi\Lambda)$ along the
Yoneda embedding 
$\FSch^{\op}\to \mathcal{P}(\FSch)^{\op}$. 
We define similarly $\RigSH^{(\hyp)}_{\et}(H^{\rig};\Lambda)$. 
\symn{$\FSH^{(\hyp)}$}
\symn{$\RigSH^{(\hyp)}$}

\item[(2)] Let $\mathcal{S}$ be a formal scheme and $H\subset \mathcal{S}$
a rig-\'etale sieve. We denote by 
\begin{equation}
\label{eq-nota:xi-tilde-H-1}
\widetilde{\xi}_H:\FSH^{(\hyp)}_{\et}(H;\chi\Lambda)
\to \RigSH^{(\hyp)}_{\et}(\mathcal{S}^{\rig};\Lambda)
\end{equation}
the functor obtained 
by evaluating on $H$ the right Kan extension of
\sym{$\widetilde{\xi}$} and then  
composing with the equivalence 
$\RigSH^{(\hyp)}_{\et}(H^{\rig};\Lambda)
\simeq 
\RigSH^{(\hyp)}_{\et}(\mathcal{S}^{\rig};\Lambda)$
provided by Theorem \ref{thm:hyperdesc}.

\end{enumerate}
\end{nota}

\begin{nota}
\label{nota:motive-poten-good-reduction-}
Let $\mathcal{S}$ be a formal scheme and $H\subset \mathcal{S}$
a rig-\'etale sieve. We denote by 
$$\RigSH^{(\hyp)}_{\et}(\mathcal{S}^{\rig};
\Lambda)_{\langle H\rangle}\subset 
\RigSH^{(\hyp)}_{\et}(\mathcal{S}^{\rig};
\Lambda)$$
the full sub-$\infty$-category generated under colimits,
desuspensions and negative Tate twists by motives of the form
$\M(\mathcal{U}^{\rig})$ where $\mathcal{U}$ is a 
formal $\mathcal{S}$-scheme which is $H$-potentially 
nearly smooth.
\symn{$\RigSH(-)_{\langle -\rangle}$}
\end{nota}

\begin{prop}
\label{prop:FSH-chi-on-cirble}
We work under Assumption \ref{assu:for-main-thm-2}.
Let $\mathcal{S}$ be a quasi-compact and quasi-separated formal scheme
and let $H\subset \mathcal{S}$ be a rig-\'etale sieve.
\begin{enumerate}

\item[(1)] The functor \eqref{eq-nota:xi-tilde-H-1} 
is fully faithful and its essential image contains 
$\RigSH^{(\hyp)}_{\et}(\mathcal{S}^{\rig};
\Lambda)_{\langle H\rangle}$.

\item[(2)] Assume that the sieve $H$ is $K$-potentially 
nearly \'etale for a rig-\'etale sieve $K\subset \mathcal{S}$. 
Then the essential image of 
\eqref{eq-nota:xi-tilde-H-1} is contained in 
$\RigSH^{(\hyp)}_{\et}(\mathcal{S}^{\rig};
\Lambda)_{\langle K\rangle}$.

\end{enumerate}
\end{prop}

\begin{proof}
Up to equivalences, the functor
$\widetilde{\xi}_H$ is given by 
$$\lim_{\mathcal{T}\to H}
\FSH^{(\hyp)}_{\et}(\mathcal{T};\chi\Lambda)
\to \lim_{\mathcal{T}\to H}
\RigSH^{(\hyp)}_{\et}(\mathcal{T}^{\rig};\Lambda),$$
where the limit is over the category of formal $H$-schemes.
Since limits in $\CAT_{\infty}$ preserve fully faithful embeddings,
Theorem \ref{thm:main-thm-}(1), proved in 
Subsection \ref{subsect:proof-of-main-thm}, 
implies that the functor $\widetilde{\xi}_H$ is fully faithful. 
Moreover, an object 
$M\in \RigSH^{(\hyp)}_{\et}(\mathcal{S}^{\rig};\Lambda)$
belongs to the essential image of $\widetilde{\xi}_H$
if and only if, for every $e:\mathcal{T} \to \mathcal{S}$
factoring through $H$, $e^{\rig,\,*}M$ belongs to the essential
image of $\widetilde{\xi}_{\mathcal{T}}$. 
This shows the first assertion. Indeed, if $\mathcal{U}$ 
is a formal $\mathcal{S}$-scheme which is $H$-potentially 
nearly smooth and $e:\mathcal{T}\to \mathcal{S}$ as before, then 
$e^{\rig,\,*}\M(\mathcal{U}^{\rig})\simeq 
\M((\mathcal{U}\times_{\mathcal{S}}\mathcal{T})^{\rig})$
is a colimit of objects of the form
$\M(\mathcal{V}^{\rig})\simeq \xi_{\mathcal{T}}\M(\mathcal{V})$
where $\mathcal{V}$ is a smooth formal $\mathcal{T}$-scheme
admitting a finite morphism to an open formal subscheme of 
$\mathcal{U}\times_{\mathcal{S}}\mathcal{T}$ 
which induces an isomorphism on generic fibers.
(Recall that such $\mathcal{V}$'s exist locally on the 
nearly smooth formal $\mathcal{T}$-scheme
$\mathcal{U}\times_{\mathcal{S}}\mathcal{T}$.)

To prove the second assertion, we assume that 
$H$ is generated by a rig-\'etale cover 
$(\mathcal{S}_i \to \mathcal{S})_i$ such that the 
formal $\mathcal{S}$-schemes
$\mathcal{S}_i$ are rig-\'etale and $K$-potentially nearly \'etale.
We want to show that the essential image of 
$\widetilde{\xi}_H$ is contained in 
$\RigSH^{(\hyp)}_{\et}(\mathcal{S}^{\rig};\Lambda)_{\langle K\rangle}$.
Let $M$ be in the essential image of $\widetilde{\xi}_H$.
Let $\mathcal{T}=\coprod_i\mathcal{S}_i$ and form the {\v C}ech
nerve $\mathcal{T}_{\bullet}$. Denote by 
$e_n:\mathcal{T}_n\to \mathcal{S}$
the obvious morphism. Then
$$\underset{[n]\in \mathbf{\Delta}}{\colim}\;
e^{\rig}_{n,\,\sharp}e_n^{\rig,\,*}M\to M$$
is an equivalence. (Indeed, by the projection formula, 
the simplicial object 
$e^{\rig}_{\bullet,\,\sharp}e_{\bullet}^{\rig,\,*}M$ 
is equivalent to 
$\M(\mathcal{T}^{\rig}_{\bullet})\otimes M$ and
$\mathcal{T}^{\rig}_{\bullet}\to \mathcal{S}^{\rig}$
is a truncated \'etale hypercover of $\mathcal{S}^{\rig}$.)
Since $e_n^{\rig,\,*}M$ belongs to the essential image of 
$\widetilde{\xi}_{\mathcal{T}_n}$, it is enough to show 
that the essential image of 
$e_{n,\,\sharp}^{\rig}\circ \widetilde{\xi}_{\mathcal{T}_n}$
is contained in 
$\RigSH^{(\hyp)}_{\et}(\mathcal{S}^{\rig};\Lambda)_{\langle K\rangle}$.
This would follow if we can prove that for every 
smooth formal $\mathcal{T}_n$-scheme $\mathcal{V}$
the formal $\mathcal{S}$-scheme $\mathcal{V}$ is 
$K$-potentially nearly smooth. This is a direct consequence of the 
definitions (and also a special case of Proposition
\ref{prop:compose-potential-smooth-etale}).
\end{proof}

Recall that $\Lder_{\riget}$ denotes the 
rig-\'etale sheafification functor.
In particular, 
$\Lder_{\riget}\FSH^{(\hyp)}_{\et}(-;\chi\Lambda)$
is the rig-\'etale sheaf associated to the 
$\Prl$-valued presheaf 
$\FSH^{(\hyp)}_{\et}(-;\chi\Lambda)$.

\begin{prop}
\label{prop:rig-fet-sheaf-FSH-chi}
We work under Assumption \ref{assu:for-main-thm-2}.
Let $\mathcal{S}$ be a quasi-compact and quasi-separated 
formal scheme. Then the functor 
\begin{equation}
\label{eq-cor:rig-fet-sheaf-FSH-chi-1}
\Lder_{\riget}\FSH^{(\hyp)}_{\et}(\mathcal{S};\chi\Lambda)
\to 
\RigSH^{(\hyp)}_{\et}(\mathcal{S}^{\rig};
\Lambda)
\end{equation}
is fully faithful with essential image 
the full sub-$\infty$-category generated under colimits,
desuspensions and negative Tate twists by motives of the form
$\M(\mathcal{U}^{\rig})$ where $\mathcal{U}$ is a 
formal $\mathcal{S}$-scheme which is potentially 
nearly smooth. In fact, we can restrict to those 
$\mathcal{U}$'s which are smooth over a quasi-compact and
quasi-separated rig-\'etale formal $\mathcal{S}$-scheme.
\end{prop}

\begin{proof}
We split the proof into three steps.

\paragraph*{Step 1}
\noindent
Let $\Lder^1_{\riget}$ be the endofunctor on presheaves
over $\FSch$ described informally as follows. Given a formal 
scheme $\mathcal{S}$ and a presheaf 
$\mathcal{F}$ with values in an $\infty$-category admitting 
limits and colimits, we have 
$$\Lder^1_{\riget}(\mathcal{F})(\mathcal{S})=
\underset{H\subset \mathcal{S}}{\colim}\,
\overline{\mathcal{F}}(H)$$
where $\overline{\mathcal{F}}$ is the right Kan extension along 
the Yoneda embedding and 
the colimit is over the rig-\'etale sieves $H\subset \mathcal{S}$.
For a precise construction of such an endofunctor, 
we refer the reader to \cite[Construction 6.2.2.9 \&
Remark 6.2.2.12]{lurie}.\footnote{This can be found 
in the electronic version of \cite{lurie}
on the author's webpage, but not in the published version.} 
(In loc.~cit., this is done for 
presheaves with values in $\mathcal{S}$, but the construction
makes sense for more general presheaves.)
\symn{$\Lder^1_\tau$}

Let $\mathcal{S}$ be a quasi-compact and quasi-separated 
formal scheme. 
Let ${\rm Sv}(\mathcal{S})$ be the set of 
rig-\'etale sieves of $\mathcal{S}$ ordered by containment 
and let ${\rm Sv}'(\mathcal{S})$ be the subset of 
${\rm Sv}(\mathcal{S})\times {\rm Sv}(\mathcal{S})$, 
endowed with the induced order, 
consisting of those pairs $(H,K)$ such that 
$H$ is $K$-potentially nearly \'etale. We have two projections 
${\rm Sv}'(\mathcal{S}) \to {\rm Sv}(\mathcal{S})$
which are cofinal by Corollary 
\ref{cor:existence-of-potent-good-red-covers}
and \cite[Theorem 4.1.3.1]{lurie}.
By Proposition \ref{prop:FSH-chi-on-cirble},
every pair $(H,K)\in{\rm Sv}'(\mathcal{S})$ gives rise to a
sequence of fully faithful embeddings
$$\RigSH^{(\hyp)}_{\et}(\mathcal{S}^{\rig};\Lambda)_{\langle H \rangle}
\to \FSH^{(\hyp)}_{\et}(H;\chi\Lambda)
\to \RigSH^{(\hyp)}_{\et}(\mathcal{S}^{\rig};\Lambda)_{\langle K \rangle}
\to \FSH^{(\hyp)}_{\et}(K;\chi\Lambda),$$
in which we identified $\FSH^{(\hyp)}_{\et}(H;\chi\Lambda)$ with its essential image under $\widetilde{\xi}_H$, 
and similarly for $K$ instead of $H$.
Passing to the colimit over ${\rm Sv}'(\mathcal{S})$
and using the cofinality of the two projections 
${\rm Sv}'(\mathcal{S}) \to {\rm Sv}(\mathcal{S})$, we obtain
an equivalence in $\Prl$:
$$\Lder^1_{\riget}\FSH^{(\hyp)}_{\et}(\mathcal{S};\chi\Lambda)
\simeq \underset{H\subset \mathcal{S}}{\colim}\,
\RigSH^{(\hyp)}_{\et}(\mathcal{S}^{\rig};\Lambda)_{\langle H \rangle}.$$
Since the sub-$\infty$-categories 
$\RigSH^{(\hyp)}_{\et}(\mathcal{S}^{\rig};\Lambda)_{\langle H \rangle}$
are generated under colimits by a set of compact generators of 
$\RigSH^{(\hyp)}_{\et}(\mathcal{S}^{\rig};\Lambda)$, 
it follows immediately that the induced functor 
$$\widetilde{\xi}{}^{\,1}_{\mathcal{S}}:\Lder_{\riget}^1\FSH^{(\hyp)}_{\et}(\mathcal{S};\chi\Lambda)
\to 
\RigSH^{(\hyp)}_{\et}(\mathcal{S}^{\rig};
\Lambda)$$
is fully faithful 
with essential image 
the full sub-$\infty$-category generated under colimits,
desuspensions and negative Tate twists by motives of the form
$\M(\mathcal{U}^{\rig})$ where $\mathcal{U}$ is a 
formal $\mathcal{S}$-scheme which is potentially 
nearly smooth.
\symn{$\widetilde{\xi}^{\,1}$}

\paragraph*{Step 2}
\noindent
Here, we prove that
$\Lder^1_{\riget}\FSH^{(\hyp)}_{\et}(-;\chi\Lambda)$,
restricted to $\FSch^{\qcqs}$, is already a rig-\'etale sheaf.
This will prove the statement except for the last sentence.
 
We argue as in the proof of Proposition
\ref{prop:FSH-chi-on-cirble}.
Let $H\subset \mathcal{S}$ be a rig-\'etale sieve 
generated by a finite family $(\mathcal{S}_i \to \mathcal{S})_i$
such that the $\mathcal{S}_i$'s are quasi-compact and rig-\'etale
over $\mathcal{S}$. We consider the functor 
$$\widetilde{\xi}{}^{\,1}_H:
\Lder^1_{\riget}\FSH^{(\hyp)}_{\et}(H;\chi\Lambda)
\to \RigSH^{(\hyp)}_{\et}(\mathcal{S}^{\rig};\Lambda)$$
defined as in Notation \ref{nota:xi-tilde-H-}(2).
This is a fully faithful functor with essential image 
the sub-$\infty$-category
spanned by those $M\in \RigSH^{(\hyp)}_{\et}
(\mathcal{S}^{\rig};\Lambda)$
such that $e^{\rig,\,*}M$ belongs to the essential image of 
$\widetilde{\xi}{}^{\,1}_{\mathcal{T}}$
for every $e:\mathcal{T}\to \mathcal{S}$ factoring through $H$.
Our goal is to show that $\widetilde{\xi}{}^{\,1}_{\mathcal{S}}$
and $\widetilde{\xi}{}^{\,1}_H$ have the same essential image.

Let $\mathcal{T}=\coprod_i \mathcal{S}_i$ and form the {\v C}ech
nerve $\mathcal{T}_{\bullet}$ associated to $\mathcal{T}\to
\mathcal{S}$. Let $e_n:\mathcal{T}_n \to \mathcal{S}$ be the 
obvious morphism. Let $M$ be in the essential image of 
$\widetilde{\xi}{}^{\,1}_H$.
We have an equivalence 
$$\underset{[n]\in \mathbf{\Delta}}{\colim}\;
e^{\rig}_{n,\,\sharp}e_n^{\rig,\,*}M\to M.$$
Therefore, it is enough to show that 
$e^{\rig}_{n,\,\sharp}e_n^{\rig,\,*}M$ belongs to
the essential image of $\widetilde{\xi}{}^{\,1}_{\mathcal{S}}$. 
Using the description of the essential image of 
$\widetilde{\xi}{}^{\,1}_H$ given above, it suffices to show that
$e^{\rig}_{n,\,\sharp}$ takes the essential image of 
$\widetilde{\xi}{}^{\,1}_{\mathcal{T}_n}$
to the  essential image of 
$\widetilde{\xi}{}^{\,1}_{\mathcal{S}}$.
This follows from the description of the essential images
of $\widetilde{\xi}{}^{\,1}_{\mathcal{S}}$ and 
$\widetilde{\xi}{}^{\,1}_{\mathcal{T}_n}$ given above,
and the fact that a potentially nearly smooth 
formal $\mathcal{T}_n$-scheme is also potentially nearly smooth as a formal
$\mathcal{S}$-scheme which follows from Propositions
\ref{lem:nearly-etale-is-rig-etale}(3) and 
\ref{prop:compose-potential-smooth-etale}.

\paragraph*{Step 3}
\noindent
It remains to show the last assertion in the statement,
concerning the generators under colimits of the essential image of 
\eqref{eq-cor:rig-fet-sheaf-FSH-chi-1}. Let 
$\mathcal{C}$ be the sub-$\infty$-category of 
$\RigSH^{(\hyp)}_{\et}(\mathcal{S}^{\rig};\Lambda)$
generated under colimits, desuspension and negative Tate 
twists by $\M(\mathcal{V}^{\rig})$, with 
$\mathcal{V}$ smooth over a rig-\'etale formal 
$\mathcal{S}$-scheme. We want to show that 
$\mathcal{C}$ coincides with the essential image of 
\eqref{eq-cor:rig-fet-sheaf-FSH-chi-1}.
By the previous steps, it is enough to show that 
$\M(\mathcal{U}^{\rig})\in \mathcal{C}$ 
for every potentially nearly 
smooth formal $\mathcal{S}$-scheme $\mathcal{U}$.
Let $\mathcal{T} \to \mathcal{S}$ be a rig-\'etale cover 
such that $\mathcal{U}\times_{\mathcal{S}}\mathcal{T}$
is nearly smooth over $\mathcal{T}$. Let $\mathcal{T}_{\bullet}$
be the {\v C}ech nerve associated to $\mathcal{T} \to \mathcal{S}$. 
Since 
$$\M(\mathcal{U}^{\rig})\simeq 
\underset{[n]\in\mathbf{\Delta}}{\colim}\,
\M((\mathcal{U}\times_{\mathcal{S}}\mathcal{T}_n)^{\rig}),$$
it is enough to show that $\M((\mathcal{U}\times_{\mathcal{S}}
\mathcal{T}_n)^{\rig})\in \mathcal{C}$ for every $n\in\N$.
The problem is local on 
$\mathcal{U}\times_{\mathcal{S}}\mathcal{T}_n$. 
Since the latter is nearly smooth, we are reduced to show that 
$\M(\mathcal{V}^{\rig}) \in \mathcal{C}$ if 
$\mathcal{V}$ is a formal $\mathcal{T}_n$-scheme admitting a 
finite morphism $\mathcal{V}'\to \mathcal{V}$
inducing an isomorphism $\mathcal{V}'^{\rig}\simeq \mathcal{V}^{\rig}$
and such that $\mathcal{V}'$ is smooth over $\mathcal{T}_n$.
This is clear since $\M(\mathcal{V}'^{\rig})\in \mathcal{C}$
by construction.
\end{proof}

\begin{cor}
\label{cor:limits-S-L-rig-et-FSH}
Let $(\mathcal{S}_{\alpha})_{\alpha}$ be a cofiltered inverse 
system of quasi-compact and quasi-separated 
formal schemes with affine transition morphisms, and let 
$\mathcal{S}=\lim_{\alpha}\mathcal{S}_{\alpha}$.
Assume one of the following conditions.
\begin{enumerate}

\item[(1)] We work under the alternative (iii) of Assumption
\ref{assu:for-main-thm}.

\item[(2)] We work under the alternative (iv) of Assumption
\ref{assu:for-main-thm}. We assume furthermore that
$\Lambda$ is eventually coconnective or that  
the numbers $\pvcd_{\Lambda}(\mathcal{S}^{\rig}_{\alpha})$ 
are bounded independently of $\alpha$.

\end{enumerate}
Then, we have an equivalence in $\Prl$:
$$\underset{\alpha}{\colim}\;
\Lder_{\riget}\FSH_{\et}
(\mathcal{S}_{\alpha};\chi\Lambda)\simeq 
\Lder_{\riget}\FSH_{\et}(\mathcal{S};\chi\Lambda).$$
\end{cor}

\begin{proof}
This follows from Theorem
\ref{thm:anstC}, Proposition
\ref{prop:rig-fet-sheaf-FSH-chi}
and the following assertion.
Given a rig-\'etale formal $\mathcal{S}$-scheme $\mathcal{T}$
and a smooth formal $\mathcal{T}$-scheme $\mathcal{V}$,
we can find, locally for the rig topology 
on $\mathcal{T}$ and $\mathcal{V}$, an index $\alpha_0$,
a rig-\'etale formal $\mathcal{S}_{\alpha_0}$-scheme
$\mathcal{T}_{\alpha_0}$, 
a smooth formal $\mathcal{T}_{\alpha_0}$-scheme 
$\mathcal{V}_{\alpha_0}$, and isomorphisms of formal
$\mathcal{S}$-schemes 
$$\mathcal{T}/(0)^{\sat}
\simeq \lim_{\alpha\leq \alpha_0}\mathcal{T}_{\alpha}/(0)^{\sat}
\qquad \text{and} \qquad
\mathcal{V}/(0)^{\sat}\simeq 
\lim_{\alpha\leq \alpha_0}\mathcal{V}_{\alpha}/(0)^{\sat}.$$
(As usual, for $\alpha\leq \alpha_0$, we set 
$\mathcal{T}_{\alpha}=\mathcal{T}_{\alpha_0}\times_{\mathcal{S}_{\alpha_0}}\mathcal{S}_{\alpha}$ 
and similarly for $\mathcal{V}_{\alpha}$.)
To prove this assertion, we may assume that the 
$\mathcal{S}_{\alpha}=\Spf(A_{\alpha})$'s are affine, 
that $\mathcal{T}=\Spf(B)$ 
with $B$ adic rig-\'etale over $A=\colim_{\alpha}\,A_{\alpha}$
and admitting a presentation as in 
Definition \ref{dfn:new-etale-rig}, 
and $\mathcal{V}=\Spf(C)$ 
with $C$ an adic $B$-algebra \'etale over 
$B\langle t_1,\ldots, t_m\rangle$. Then, the result follows 
easily from Corollary \ref{cor:proj-limi-cal-E-A}.
\end{proof}

\begin{rmk}
\label{rmk:reduction-for-main-thm-iv}
Recall that our goal in this subsection is to prove 
Theorem \ref{thm:sheafifi-FSH-chi-}. This is equivalent
to the statement that  
the morphism of rig-\'etale $\Prl$-valued sheaves
\begin{equation}
\label{eq-rmk:reduction-for-main-thm-iv}
\Lder_{\riget}\widetilde{\xi}:
\Lder_{\riget}\FSH^{(\hyp)}_{\et}(-;\chi\Lambda)
\to 
\RigSH^{(\hyp)}_{\et}((-)^{\rig};\Lambda)
\end{equation}
is an equivalence under Assumption \ref{assu:for-main-thm-2}.
Clearly, it is enough to do so after restricting
\eqref{eq-rmk:reduction-for-main-thm-iv} 
to affine formal schemes. Every affine formal scheme
is the limit of a cofiltered inverse system of $(\Lambda,\et)$-admissible
affine formal schemes with 
$(\Lambda,\et)$-admissible generic fiber. Thus, when working under the 
alternative (iii) of Assumption \ref{assu:for-main-thm}, 
Theorem \ref{thm:anstC} and Corollary \ref{cor:limits-S-L-rig-et-FSH}
allow us to restrict  
\eqref{eq-rmk:reduction-for-main-thm-iv}
further to the subcategory of $(\Lambda,\et)$-admissible
affine formal schemes with
$(\Lambda,\et)$-admissible generic fiber. 
By Propositions \ref{prop:automatic-hypercomp-motives} and 
\ref{prop:automatic-hypercomp-motives-algebraic},
we are then automatically working under
the alternative (iv) of Assumption \ref{assu:for-main-thm}. Said differently, to prove 
Theorem \ref{thm:sheafifi-FSH-chi-} we may work from this point 
onwards under the alternative (iv)
of Assumption \ref{assu:for-main-thm}. In particular, since 
we only consider formal schemes with finite dimensional generic fibers, 
\eqref{eq-rmk:reduction-for-main-thm-iv}
is a morphism of rig-Nisnevich hypersheaves.
(See the proof of Lemma
\ref{lem:auto-hypercomp-etale}.)
As a consequence, it is enough to show that 
\eqref{eq-rmk:reduction-for-main-thm-iv} 
induces equivalences on the stalks for the rig-Nisnevich 
topology. Using Theorem \ref{thm:etst}
and the analogous statement for $\Lder_{\riget}\FSH^{(\hyp)}_{\et}(-;\chi\Lambda)$ which follows in the same way from 
Corollary \ref{cor:limits-S-L-rig-et-FSH}, we are left 
to show the following statement.
\end{rmk}

\begin{prop}
\label{prop:grgen}
Let $s$ be a rigid point and set 
$\mathfrak{s}=\Spf(\kappa^+(s))$. Assume the following 
conditions:
\begin{enumerate}

\item[(1)] every prime number 
is invertible either in $\kappa^+(s)$ or in $\pi_0\Lambda$;

\item[(2)] when working in the non-hypercomplete case, 
$\Lambda$ is eventually coconnective.

\end{enumerate}
Then, $\RigSH^{(\hyp)}_{\et}(s;\Lambda)$ is 
generated under colimits, desuspension and negative 
Tate twists by motives of the form
$\M(\mathcal{U}^{\rig})$ with 
$\mathcal{U}$ smooth over a rig-\'etale
formal $\mathfrak{s}$-scheme
(or, equivalently, by the motives 
$M(U)$ with $U$ smooth with good reduction over an 
\'etale rigid analytic $s$-space).
\end{prop}

\begin{proof}
This is a generalisation of 
\cite[Theorem 2.5.34]{ayoub-rig}, and we will adapt 
the proof of loc.~cit. to our situation. 
Let $\mathcal{C}(s)$ be the 
sub-$\infty$-category of $\RigSH^{(\hyp)}_{\et}(s;\Lambda)$ 
generated under colimits, desuspension and negative Tate twists 
by motives of the form $\M(\mathcal{U}^{\rig})$,
with $\mathcal{U}$ smooth over a rig-\'etale formal 
$\mathfrak{s}$-scheme.
Note that $\mathcal{C}(s)$ is equally generated by 
motives of the form $\M(U)$,
with $U$ smooth with good reduction over an \'etale rigid analytic 
$s$-space.
Our goal is to show that $\mathcal{C}(s)$ is equal to
$\RigSH_{\et}^{(\hyp)}(s;\Lambda)$.
We divide the proof into several steps.

\paragraph*{Step 1}
\noindent 
Here we show that it is enough to prove the proposition 
under the following assumptions:
\begin{itemize}

\item $\pi_0\Lambda$ is a $\Q$-algebra;

\item $\kappa(s)$ is algebraically closed and 
$\kappa^+(s)$ has finite height.

\end{itemize}
In particular, $s$ is $(\Lambda,\et)$-admissible and 
we will be working in the hypercomplete case.

Indeed, we can find a cofiltered inverse system of 
rigid points $(s_{\alpha})_{\alpha}$ with 
$s\sim \lim_{\alpha}s_{\alpha}$ such that 
the valuation rings $\kappa^+(s_{\alpha})$ have finite
ranks and the fields $\kappa(s_{\alpha})$ 
have finite virtual $\Lambda$-cohomological dimensions. 
We set $\mathfrak{s}_{\alpha}=\Spf(\kappa^+(s_{\alpha}))$
so that $\mathfrak{s}=\lim_{\alpha}\mathfrak{s}_{\alpha}$.
Our goal is to prove that 
$\mathcal{C}(s)=\RigSH^{(\hyp)}_{\et}(s;\Lambda)$
and, by Lemma \ref{lem:generation-rigsh}, it is enough to show 
that $\M(\mathcal{V}^{\rig})\in \mathcal{C}(s)$
for $\mathcal{V}$ a rig-smooth formal $\mathfrak{s}$-scheme. 
Moreover, we may assume that $\mathcal{V}=\Spf(A)$ 
where $A$ is an adic $\kappa^+(s)$-algebra 
which is rig-\'etale over $\kappa^+(s)\langle t_1,\cdots, t_m\rangle$.
Thus, using Corollary
\ref{cor:proj-limi-cal-E-A}, there is an index $\alpha$ 
and a rig-smooth formal $\mathfrak{s}_{\alpha}$-scheme
$\mathcal{V}_{\alpha}$ such that $\mathcal{V}^{\rig}=
\mathcal{V}_{\alpha}^{\rig}\times_{s_{\alpha}}s$. 
Since $\mathcal{C}(s)$ contains the image of 
$\mathcal{C}(s_{\alpha})$ by the inverse image functor along
$s\to s_{\alpha}$, we see that it is enough to show that 
$\M(\mathcal{V}_{\alpha}^{\rig})\in \mathcal{C}(s_{\alpha})$. 
Thus, we may replace $s$ by $s_{\alpha}$ and assume that 
$s$ is $(\Lambda,\et)$-admissible.
In particular, by Proposition
\ref{prop:automatic-hypercomp-motives},
the non-hypercomplete case is then 
covered by the hypercomplete case. Also, 
the $\infty$-category $\RigSH^{\hyp}_{\et}(s;\Lambda)$
is compactly generated by Proposition
\ref{prop:compact-shv-rigsm}.

Next, we explain how to reduce to the case where
$\pi_0\Lambda$ is a $\Q$-algebra. Let 
$M\in \RigSH^{\hyp}_{\et}(s;\Lambda)$ 
and consider the cofiber sequence 
$M\to M_{\Q}\to M_{\rm tor}$
where $M_{\Q}=M\otimes \Q$ is the rationalisation of $M$.
The motive $M_{\rm tor}$ is a direct coproduct of $\ell$-nilpotent 
motives $M_{\ell}$ for $\ell$ non invertible in $\pi_0\Lambda$. 
By Theorem \ref{thm:rigrig}, 
we have an equivalence of $\infty$-categories
$$\Shv_{\et}^{\hyp}(\Et/s;\Lambda)_{\ellnil}
\simeq \RigSH^{\hyp}_{\et}(s;\Lambda)_{\ellnil}.$$
This implies that $M_{\ell}$ belongs to the sub-$\infty$-category
of $\RigSH^{\hyp}_{\et}(s;\Lambda)$ generated under colimits 
by motives of the form $\M(U)$, where $U$ is an \'etale
rigid analytic $s$-space. This show that 
$M_{\rm tor}$ belongs to $\mathcal{C}(s)$, and 
we are left to show that $M_{\Q}$ belongs to $\mathcal{C}(s)$.
To do so, we may replace $\Lambda$ with $\Lambda_{\Q}$ 
and assume that $\pi_0\Lambda$ is a $\Q$-algebra.

It remains to explain how to reduce to the case where 
$\kappa(s)$ is algebraically closed.
Let $\kappa^+(\overline{s})$ be the adic completion of a
valuation ring extending $\kappa^+(s)$ inside a 
separable closure of $\kappa(s)$, and let $\kappa(\overline{s})$
be the fraction field of $\kappa^+(\overline{s})$. 
This defines a geometric algebraic point $\overline{s}$
over $s$ as in Construction \ref{cons:point-an-nis-et}(2). 
We have 
$\overline{s}\sim \lim_{\alpha} \overline{s}_{\alpha}$
where $(\overline{s}_{\alpha})_{\alpha}$ is the cofiltered inverse system 
of rigid points such that $\kappa(\overline{s}_{\alpha})/\kappa(s)$ 
is a finite separable extension contained in $\kappa(\overline{s})$.
Using Theorem \ref{thm:anstC} and arguing as above, 
we have an equivalence in $\Prl_{\omega}$:
\begin{equation}
\label{eqn:continuity-cal-C-rig-point}
\mathcal{C}(\overline{s})\simeq \underset{\alpha}{\colim}\,
\mathcal{C}(\overline{s}_{\alpha}).
\end{equation}
Denote by $e:\overline{s} \to s$, $e_{\alpha}:\overline{s}
\to \overline{s}_{\alpha}$
and $r_{\alpha}:\overline{s}_{\alpha}\to s$ 
the obvious morphisms.
Consider a compact motive 
$M\in \RigSH^{\hyp}_{\et}(s;\Lambda)$ and assume that we 
know that $e^*M\in \mathcal{C}(\overline{s})$.
Since $e^*M$ is compact, the equivalence
\eqref{eqn:continuity-cal-C-rig-point}
implies that there exists $\alpha_0$ and 
a compact object $N\in \mathcal{C}(\overline{s}_{\alpha_0})$ 
such that $e^*M=e_{\alpha_0}^*N$.
In particular, the two compact objects
$r_{\alpha_0}^*M$ and $N$ of 
$\RigSH^{\hyp}_{\et}(\overline{s}_{\alpha_0};\Lambda)$
become equivalent when pulled back to $\overline{s}$. 
By Theorem \ref{thm:anstC}, they actually become equivalent 
when pulled back to $\overline{s}_{\alpha}$, for $\alpha\leq 
\alpha_0$ sufficiently small. This shows that 
$r_{\alpha}^*M$ belongs to $\mathcal{C}(\overline{s}_{\alpha})$.
We now conclude as in the second step of the proof of 
Proposition \ref{prop:rig-fet-sheaf-FSH-chi}: using the {\v C}ech 
nerve associated to $\overline{s}_{\alpha} \to s$,
we reduce to show that, for $n\geq 1$, 
$$M\otimes \M(\overbrace{\overline{s}_{\alpha}\times_s\cdots\times_s
\overline{s}_{\alpha}}^{n\,\text{times}})\simeq 
(r_{\alpha,\,\sharp}r_{\alpha}^*M)\otimes 
\M(\overbrace{\overline{s}_{\alpha}\times_s\cdots\times_s
\overline{s}_{\alpha}}^{n-1\,\text{times}}),$$
belongs to $\mathcal{C}(s)$ which is clear.

\paragraph*{Step 2}
\noindent
In the remainder of the proof, we work under the two 
assumptions introduced in the first step.
We set $K=\kappa(s)$, $V=\kappa^+(s)$ and we fix
$\pi\in V$ a generator of an ideal of definition.
We set $\eta=\Spec(K)$ and use a subscript 
``$\eta$'' to denote the fiber at $\eta$ of a $V$-scheme.
By Lemma \ref{lem:generation-rigsh},
the $\infty$-category $\RigSH_{\et}^{\hyp}(s;\Lambda)$ 
is generated under colimits 
by the motives $\M(Y)$, for $Y\in \RigSm^{\rm qc}/s$,
and their desuspensions and negative Tate twists.
We will show that $\M(Y)\in \mathcal{C}(s)$ by induction 
on the relative dimension $d$ of $|Y|$ over $|s|$. 
The case of relative dimension zero is clear
because $Y$ is then \'etale over $s$. 
In general, the problem is local on $Y$.
Thus, by Proposition
\ref{prop:formal-compl-smooth},
we may assume that 
$Y$ is the $\pi$-adic completion $\widehat{P}$
of a $V$-scheme $P$ of finite presentation and
generically smooth.
Replacing $P$ with the Zariski closure of $P_{\eta}$, 
we may also assume that $P$ is flat over $V$. 

\paragraph*{Step 3}
\noindent 
(This is analogous to the second step in the proof of 
\cite[Th\'eor\`eme 2.5.34]{ayoub-rig}.) 
In this step, we will prove the following preliminary assertion.
Let $E\subset P$ be a closed subscheme, generically of 
codimension $\geq 1$, and let $Z=\widehat{E}{}^{\rig}$
considered as a closed rigid analytic subspace of $Y$.
Then the relative motive
$\M(Y/Y\smallsetminus Z)$, defined as the cofiber of 
$\M(Y\smallsetminus Z) \to \M(Y)$,
belongs to $\mathcal{C}(s)$. 
The proof of this uses the induction on the 
relative dimension of $|Y|$ over $|s|$, and
we will argue by 
a second induction on the dimension of $E_{\eta}$. 
The base case for the second induction is when 
$E_{\eta}$ is empty: the relative motive is then zero 
and the claim is obvious.
Let $E'\subset E$ be the closure of the singularity locus of 
$E_{\eta}$ and $Z'=\widehat{E}{}'^{\rig}$. 
Since $\kappa(s)$ is algebraically closed and hence perfect, 
$E'_{\eta}$ has codimension $\geq 1$ in $E_{\eta}$. 
By the second induction, we may 
assume that $\M(Y/Y\smallsetminus Z')$ belongs to 
$\mathcal{C}(s)$. We are thus left to show that 
$\M(Y\smallsetminus Z'/Y\smallsetminus Z)$
belongs to $\mathcal{C}(s)$.
The rigid analytic space $Y\smallsetminus Z'$ is not necessarily 
quasi-compact, but we may write it as a filtered union of 
quasi-compact opens $Y_{\alpha}=(\widehat{P}_{\alpha})^{\rig}$ 
where $P_{\alpha}$ are open subschemes of admissible blowups of 
$P$, not meeting the closure of $E'_{\eta}$.
Thus, we are left to show that $\M(Y_{\alpha}/Y_{\alpha}\smallsetminus Z_{\alpha})$ belongs to $\mathcal{C}(s)$ 
with $Z_{\alpha}$ the generic fiber of the 
formal completion of $E_{\alpha}=E\times_P P_{\alpha}$. 
Replacing $Y$ with $Y_{\alpha}$ and $E$ with $E_{\alpha}$, 
we are thus reduced to showing that $\M(Y/Y\smallsetminus Z)$ belongs 
$\mathcal{C}(s)$ under the assumption that $E_{\eta}$ is smooth.

As usual, we may also assume that $P$ is affine, and that 
$E$ is flat over $V$. Now, assume we are given a finite type morphism 
$e:\widetilde{P} \to P$ and a closed subscheme 
$\widetilde{E}\subset \widetilde{P}$
with the following properties:
\begin{itemize}

\item $e_{\eta}$ is \'etale, $\widetilde{E}\subset e^{-1}(E)$ and
$\widetilde{E}_{\eta}=e_{\eta}^{-1}(E_{\eta})$;

\item the induced morphism $\widetilde{E} \to E$ 
is proper and an isomorphism 
$\widetilde{E}_{\eta}\simeq E_{\eta}$ on generic fibers.

\end{itemize}
Then, letting $\widetilde{Y}$ and $\widetilde{Z}$ be the 
generic fibers of the $\pi$-adic completions of $\widetilde{P}$ and 
$\widetilde{E}$, we have, by \'etale excision, an isomorphism
$\M(\widetilde{Y}/\widetilde{Y}\smallsetminus \widetilde{Z})\simeq 
\M(Y/Y\smallsetminus Z)$. Using this principle twice, we may 
assume that $P$ is isomorphic to $E\times \A^c$, for some $c\geq 1$,
and that $E\subset P$ is the zero section.
In this case, the relative motive 
$\M(Y/Y\smallsetminus Z)$ is isomorphic to $\M(Z)(c)[2c]$, and we 
may conclude using the induction on the relative dimension of $Y$.

\paragraph*{Step 4}
\noindent
(This is analogous to the third step in the proof of 
\cite[Th\'eor\`eme 2.5.34]{ayoub-rig}.)
In this step, we show that we may assume $P$ 
to be ``poly-stable''.
By means of \cite[Lemma 9.2]{berk-contr}, 
applied to some compactification of $P$, we may find
a proper surjective morphism $e:Q\to P$
with the following properties:
\begin{itemize}

\item there is a finite group $G$ acting on the $P$-scheme $Q$,
a dense open subscheme $L\subset P_{\eta}$ with inverse image
$M=e^{-1}(L)$ dense in $Q_{\eta}$, and such that 
$M \to M/G$ is a finite \'etale Galois cover with group $G$
and $M/G \to L$ is a universal homeomorphism;

\item the projection $Q\to \Spec(V)$ factors as a composition of 
$$Q=Q_d\xrightarrow{f_d}
Q_{d-1}{\to}\ldots\to Q_1\xrightarrow{f_1} Q_0=\Spec(V)$$
and, for every $1\leq i \leq d$, the 
morphism $f_i$ decomposes, \'etale locally on the source and 
the target, as 
\begin{equation}
\label{eqn:uva}
\Spec(B) \xrightarrow{\text{\'etale}}
\Spec(A[u,v]/(uv-a)) \to \Spec(A)
\end{equation}
with $A$ a flat $V$-algebra of finite type, $u$ and $v$
two indeterminates, and $a\in A$ invertible in 
$A[\pi^{-1}]$. 

\end{itemize}
In particular, we see
that the $f_i$'s have relative dimension $1$ and that  
the $(f_i)_{\eta}$'s are smooth.

Let $E\subset P$ be the closure of $P_{\eta}\smallsetminus L$ in $P$
and $F\subset Q$ the closure of $Q_{\eta}\smallsetminus M$ in $Q$.
By the second step, it is enough to prove that 
$\M(\widehat{P}{}^{\rig}\smallsetminus \widehat{E}{}^{\rig})$
belongs to $\mathcal{C}(s)$. By Lemma 
\ref{lem:direct-summand-finite-cover} below, 
$\M(\widehat{P}{}^{\rig}\smallsetminus \widehat{E}{}^{\rig})$ 
is a direct summand of 
$\M(\widehat{Q}{}^{\rig}\smallsetminus \widehat{F}{}^{\rig})$
and it is enough to see that the latter is in 
$\mathcal{C}(s)$. Using the second step again, we see that it is 
enough to show that $\M(\widehat{Q}{}^{\rig})$
belongs to $\mathcal{C}(s)$.
Thus, replacing $P$ with $Q$ and 
$Y$ with $\widehat{Q}{}^{\rig}$, we may assume that 
the projection $P\to \Spec(V)$ can be factored as a composition 
\begin{equation}
\label{eqn:fi}
P=P_d\xrightarrow{f_d}
P_{d-1}{\to}\ldots\to P_1\xrightarrow{f_1} P_0=\Spec(V)
\end{equation}
with $f_i$ given, \'etale locally on the source and 
the target, by 
\eqref{eqn:uva}.

\paragraph*{Step 5} 
\noindent
We now conclude the proof.
We argue by induction on the number of integers $i\in \{1,\ldots, d\}$
such that $f_i$ is not smooth.
If all the $f_i$'s are smooth, then the formal scheme
$\widehat{P}$ is smooth over $\Spf(V)$ and 
$\M(Y)\in \mathcal{C}(s)$ by construction. 
Now suppose that at least one of the $f_i$'s 
is not smooth. Arguing as in 
\cite[page 332]{ayoub-rig},\footnote{We remind the reader that
the page references to \cite{ayoub-rig} correspond to the 
published version.} we may assume that 
$f_d:P_d \to P_{d-1}$ is not smooth. The problem is local 
for the \'etale topology on $Y$. (More precisely, 
if $Y_{\bullet} \to Y$ is a truncated \'etale hypercover 
then it is enough 
to prove that $\M(Y_n)\in \mathcal{C}(s)$ for $n\geq 0$.)
Therefore, we may assume that a factorization as in 
\eqref{eqn:uva}
exists globally for $f_d$, i.e., that $f_d$ is a composition of 
$$P=P_d \xrightarrow{\text{\'etale}} P_{d-1}[u,v]/(uv-a)
\to P_{d-1}$$
for some $a\in \mathcal{O}(P_{d-1})$
which is invertible in $\mathcal{O}((P_{d-1})_{\eta})$.
Arguing by \'etale excision as in \cite[page 333]{ayoub-rig}, 
we conclude that it suffices to treat the case where 
$P=P_{d-1}[u,v]/(uv-a)$.

We set $R=P_{d-1}$. By the induction on the relative 
dimension of $|Y|\to |s|$, we know that 
$\M(\widehat{R}{}^{\rig})$ belongs to $\mathcal{C}(s)$.
Consider the blowup $e:W \to R[u]$ of the ideal $(a,u)$. 
Since $a$ is invertible on $R_{\eta}$, $e_{\eta}$ is an isomorphism
and $\widehat{W}{}^{\rig}\simeq \widehat{R}{}^{\rig}\times \B^1$.
Moreover, $W$ admits a Zariski cover given by 
$P=R[u,v]/(uv-a)$ and $P'=R[u,w]/(aw-u)\simeq R[w]$ 
intersecting at 
$P''=R[u,v,v^{-1}]/(uv-a)\simeq R[v,v^{-1}]$.
Thus, we have a cofiber sequence
$$\M(\widehat{R}{}^{\rig}\times \U^1)
\to \M(\widehat{P}{}^{\rig})\oplus 
\M(\widehat{R}{}^{\rig}\times \B^1)
\to \M(\widehat{R}{}^{\rig}\times \B^1)$$
showing that $\M(Y)$ is isomorphic to 
$\M(\widehat{R}{}^{\rig})\oplus 
\M(\widehat{R}{}^{\rig})(1)[1]$. 
This finishes the proof. 
\end{proof}

\begin{lemma}
\label{lem:direct-summand-finite-cover}
Let $S$ be a rigid analytic space, $f:Y \to X$ a morphism of 
smooth rigid analytic $S$-spaces and
$G$ a finite group acting on the rigid analytic $X$-space $Y$.
Assume that $Y \to Y/G$ is a finite \'etale cover and that 
$Y/G \to X$ is a universal homeomorphism. 
Assume also that the order of $G$ is invertible in 
$\pi_0\Lambda$ and that every prime number is invertible 
either in $\mathcal{O}(X)$ or in $\pi_0\Lambda$. Then, in
the $\infty$-category $\RigSH^{(\hyp)}_{\et}(S;\Lambda)$, 
the morphism $\M(Y) \to \M(X)$ induced by $f$
exhibits $\M(X)$ as the image of the 
projector $|G|^{-1}\sum_{g\in G}g$ acting on $\M(Y)$.
\end{lemma}

\begin{proof}
Let $\pi_X:X \to S$ and $\pi_Y:Y \to S$ be the structural morphisms.
Since
$\M(X)=\pi_{X,\,\sharp}\pi_X^*\Lambda$, 
there is an equivalence of copresheaves 
$$\Map_{\RigSH^{(\hyp)}_{\et}(S;\,\Lambda)}(\M(X),-)
\simeq \Map_{\RigSH^{(\hyp)}_{\et}(S;\,\Lambda)}
(\Lambda,\pi_{X,\,*}\pi_X^*(-)),$$
and similarly for $Y$ instead of $X$.
Thus, by Yoneda's lemma, it is enough to show that, for 
every $M\in \RigSH^{(\hyp)}_{\et}(S;\Lambda)$, 
the obvious morphism
$\pi_{X,\,*}\pi_X^*M\to \pi_{Y,\,*}\pi_Y^*M$
exhibits $\pi_{X,\,*}\pi_X^*M$ as the image of the projector 
$|G|^{-1}\sum_{g\in G}g$ acting on $\pi_{Y,\,*}\pi_Y^*M$.
Set $X'=Y/G$ and let $\pi_{X'}:X'\to S$ be the structural morphism. 
By \'etale descent, the image of the projector 
$|G|^{-1}\sum_{g\in G}g$ acting on $\pi_{Y,\,*}\pi_Y^*M$
is equivalent to $\pi_{X',\,*}\pi_{X'}^*M$. 
Thus, we need to show that the natural transformation 
$\pi_{X,\,*}\pi_X^*\to \pi_{X',\,*}\pi_{X'}^*$
is an equivalence. This follows from the fact that the unit 
morphism $\id \to e_*e^*$ is an equivalence, 
which is a consequence of Theorem \ref{thm:sepalgebr}.
\end{proof}

Now that we have completed the proof of 
Theorem \ref{thm:sheafifi-FSH-chi-}, we record the following 
generalisation of Proposition \ref{prop:grgen}.

\begin{cor}
\label{cor-prop:grgen}
Let $S$ be a rigid analytic space. Assume the following 
conditions:
\begin{enumerate}

\item[(1)] every prime
number is invertible either in every $\kappa^+(s)$ for 
$s\in |S|$ or in $\pi_0\Lambda$;

\item[(2)] when working in the non-hypercomplete case, $\Lambda$
is eventually coconnective. 

\end{enumerate}
Then $\RigSH^{(\hyp)}_{\et}(S;\Lambda)$
is generated under colimits, desuspension and negative 
Tate twists by the motives $\M(U)$ with 
$U$ smooth with good reduction over an \'etale rigid analytic 
$S$-space.
\end{cor}

\begin{proof}
The problem is local on $S$. Thus, we may assume that 
$S=\Spf(A)^{\rig}$ with $A$ an adic ring.
We may write $A$ as the colimit in the category of adic 
rings of a filtered direct system $(A_{\alpha})_{\alpha}$
such that the $\mathcal{S}_{\alpha}=\Spf(A_{\alpha})$ and 
the $S_{\alpha}=\Spf(A_{\alpha})^{\rig}$ are 
$(\Lambda,\et)$-admissible. 
Arguing as in the first step of the proof of 
Proposition \ref{prop:grgen}, we see that it is enough 
to prove the corollary for each $S_{\alpha}$. 
Said differently, we may assume that 
$\mathcal{S}=\Spf(A)$ and $S$ are $(\Lambda,\et)$-admissible.
By Theorem \ref{thm:sheafifi-FSH-chi-},
we have an equivalence
$$\Lder_{\riget}\FSH^{(\hyp)}_{\et}(\mathcal{S};\chi\Lambda)
\simeq \RigSH^{(\hyp)}_{\et}(S;\Lambda).$$
We may now conclude using 
Proposition \ref{prop:rig-fet-sheaf-FSH-chi}.
\end{proof}

\begin{cor}
\label{cor:chi-non-trivial-generators} 
Let $S$ be a rigid analytic space and assume the conditions
(1) and (2) of Corollary \ref{cor-prop:grgen}.
For every $U\in \Et^{\qcqs}/S$, denote by 
$f_U:U \to S$ the structural morphism and choose a formal model 
$\mathcal{U}$ of $U$. Then, the functors
$$\chi_{\mathcal{U}}\circ f_U^*:
\RigSH^{(\hyp)}_{\et}(S;\Lambda)
\to \FSH^{(\hyp)}_{\et}(\mathcal{U};\Lambda),$$
for $U\in \Et^{\qcqs}/S$, 
form a conservative family. In fact, the same is true
if we restrict to those $U$'s admitting affine formal models
of principal ideal type.
\end{cor}

\begin{proof}
This follows immediately from Proposition \ref{conservativity-generation}
and Corollary \ref{cor-prop:grgen}.
\end{proof}

We end the subsection with the following statement.

\begin{thm}
\label{thm:rqrem}
We assume that $\tau$ is the \'etale topology
and work under one of the alternatives (ii), (iii) and (iv) 
of Assumption \ref{assu:for-main-thm}.
Let $s$ be a geometric rigid point and set 
$\mathfrak{s}=\Spf(\kappa^+(s))$. Then  
$$\widetilde{\xi}_{\mathfrak{s}}:
\FSH^{(\hyp)}_{\et}(\mathfrak{s};\chi\Lambda)
\to \RigSH^{(\hyp)}_{\et}(s;\Lambda)$$
is an equivalence of $\infty$-categories.
\end{thm}

\begin{proof}
When working under (iii) or (iv), this is a direct consequence 
of Theorem \ref{thm:main-thm-}(2) and the fact that 
every rig-\'etale cover of $\mathfrak{s}$ splits.
In the generality considered in the statement, 
we argue as follows.
The functor 
$\widetilde{\xi}_{\mathfrak{s}}$ is fully faithful
by Theorem \ref{thm:main-thm-}(1).
Since this functor preserves colimits, it remains to see
that its image  generates $\RigSH^{(\hyp)}_{\et}(s;\Lambda)$
under colimits. This follows from Proposition 
\ref{prop:grgen} and the fact that an \'etale rigid analytic 
$s$-space is a coproduct of open subspaces.
\end{proof}

\subsection{Complement}

$\empty$

\smallskip

\label{subsect:complement}

Theorem \ref{thm:main-thm-}
is especially useful if we have a handle on the 
commutative algebras $\chi_{\mathcal{S}}\Lambda$, for
$\mathcal{S}\in \FSch$. Our goal 
in this subsection is to obtain a purely algebro-geometric description
of these commutative algebras, i.e., one that does not involve rigid analytic geometry. In order to do so, we need to assume that 
$\tau$ is the \'etale topology; the case of the Nisnevich 
topology seems to require techniques of resolution
of singularities which are stronger than what is available.

Given a formal scheme $\mathcal{S}$, we will implicitly 
identify the $\infty$-categories 
$\SH^{(\eff,\,\hyp)}_{\tau}(\mathcal{S}_{\sigma};\Lambda)$
and $\FSH^{(\eff,\,\hyp)}_{\tau}(\mathcal{S};\Lambda)$
by means of Theorem \ref{thm:formal-mot-alg-mot}. In particular, 
$\chi_{\mathcal{S}}\Lambda$ will be considered as a commutative
algebra in $\SH^{(\eff,\,\hyp)}_{\tau}(\mathcal{S}_{\sigma};\Lambda)$.
Our goal is to prove Theorem \ref{thm:compute-chi} below. 
The proof will occupy most of the subsection, and it is 
inspired by the proof of \cite[Th\'eor\`eme
1.3.38]{ayoub-rig}.

\begin{thm}
\label{thm:compute-chi}
Let $B$ be a scheme, $B_{\sigma}\subset B$ 
a closed subscheme locally of finite presentation
up to nilimmersion, and $B_{\eta}\subset B$ its open complement.
Consider the functor 
$$\chi_B:\SH^{(\hyp)}_{\et}(B_{\eta};\Lambda)
\to \SH^{(\hyp)}_{\et}(B_{\sigma};\Lambda)$$
given by $\chi_B=i^*\circ j_*$, where 
$i:B_{\sigma}\to B$ and $j:B_{\eta}\to B$
are the obvious immersions.
Assume that every prime number is invertible either 
in $\pi_0\Lambda$ or in $\mathcal{O}(B)$.
Assume one of the following alternatives.
\symn{$\chi$}
\begin{enumerate}

\item[(1)] We work in the non-hypercomplete case and
$\Lambda$ is eventually coconnective;

\item[(2)] We work in the hypercomplete case and $B$ is 
$(\Lambda,\et)$-admissible.

\end{enumerate}
Let $\widehat{B}$ be the formal completion of $B$ at 
$B_{\sigma}$. (Note that $B_{\sigma} = \widehat{B}_{\sigma}$
up to nilimmersion.) Then, there is an equivalence
$\chi_B\Lambda\simeq \chi_{\widehat{B}}\Lambda$ of commutative algebras
in $\SH^{(\hyp)}_{\et}(B_{\sigma};\Lambda)$.
\end{thm}

\begin{rmk}
\label{rmk:chi-B-computable}
One has a good handle on the motive $\chi_B\Lambda$ 
in many situations. For example, if 
$B$ is regular and $B_{\sigma}$ is a principal regular divisor in $B$, 
then $\chi_B\Lambda\simeq \Lambda\oplus \Lambda(-1)[-1]$.
This follows from absolute purity; see Corollary 
\ref{cor:computing-chi-alg} below. More generally, 
absolute purity can be used to give a precise description 
of $\chi_B\Lambda$ when 
$B_{\sigma}$ is a normal crossing divisor of a regular scheme $B$. 
In general, assuming that $B$ is quasi-excellent, 
one can access $\chi_B\Lambda$ using
techniques of resolution of singularities 
to reduce to the case where $B$ is regular and 
$B_{\sigma}$ is a normal crossing divisor. 
In fact, these techniques will also be used in the proof of 
Theorem \ref{thm:compute-chi}.
\end{rmk}

{
\begin{rmk}\label{rmk:k((pi))}
Let $k$ be a field of characteristic zero having finite 
virtual $\Lambda$-cohomological dimension. 
In the non-hypercomplete case, assume that 
$\Lambda$ is eventually coconnective. Let 
$K$ be the discretely valued field $k((\pi))$
and $R\subset K$ its valuation ring. 
For $n\in \N^{\times}$, we denote by $K_n=K[\pi^{1/n}]$
the finite extension of $K$ obtained by adjoining an $n$-th root of 
unity, and $R_n\subset K_n$ its valuation ring. 
Also, we let $K_{\infty}$ be the completion of 
$\bigcup_{n\in \N^{\times}} K_n$ and $R_{\infty}\subset K_{\infty}$ 
its valuation ring. Using Theorem \ref{thm:compute-chi}
(and Remark \ref{rmk:chi-B-computable}), 
we obtain canonical equivalences of commutative algebras
$$\chi_{R_n}\Lambda \simeq q_{n,\,*}\Lambda$$
where $q_n:T_n \to \Spec(k)$
is the structural projection of the $1$-dimensional torus
$T_n\simeq \G_{\rm m}$ given by $\Spec(k[\pi^{1/n},\pi^{-1/n}])$. 
It follows formally that we have an 
equivalence of $\infty$-categories
$$\SH^{(\hyp)}_{\et}(k,\chi_{R_n}\Lambda)\simeq 
\uSH^{(\hyp)}_{\et}(T_n;\Lambda),$$
where $\uSH^{(\hyp)}_{\et}(T_n;\Lambda)$ is the full sub-$\infty$-category
of $\SH^{(\hyp)}_{\et}(T_n;\Lambda)$ generated under colimits 
by the image of the functor $q_n^*$. Letting $n$ go to $\infty$, 
we obtain an equivalence of $\infty$-categories
\begin{equation}
\label{eq-rmk:k((pi))-1}
\SH^{(\hyp)}_{\et}(k,\chi_{R_{\infty}}\Lambda)\simeq 
\uSH^{(\hyp)}_{\et}(T_{\infty};\Lambda),
\end{equation}
where $T_{\infty}$ is the pro-torus given by the spectrum of 
$\bigcup_{n\in \N^{\times}}k[\pi^{1/n},\pi^{-1/n}]$
and $\uSH^{(\hyp)}_{\et}(T_{\infty};\Lambda)$ is 
defined similarly as for the $T_n$'s. Now, assume furthermore that 
$k$ is algebraically closed. Then the valued field $K_{\infty}$ is 
also algebraically closed. 
Combining the equivalence 
\eqref{eq-rmk:k((pi))-1}
with Theorem \ref{thm:rqrem}, 
we obtain an equivalence 
of $\infty$-categories
\begin{equation}
\label{eq-rmk:k((pi))-3}
\uSH^{(\hyp)}_{\et}(T_{\infty};\Lambda)
\simeq \RigSH^{(\hyp)}_{\et}(K_{\infty};\Lambda).
\end{equation}
Moreover, one can check that this equivalence is given by the 
composition of
$$\uSH^{(\hyp)}_{\et}(T_{\infty};\Lambda)
\subset \SH^{(\hyp)}_{\et}(T_{\infty};\Lambda)
\to \SH^{(\hyp)}_{\et}(K_{\infty};\Lambda)
\xrightarrow{\An^*}
\RigSH^{(\hyp)}_{\et}(K_{\infty};\Lambda).$$
In fact, by Galois descent, one can show that the equivalence 
\eqref{eq-rmk:k((pi))-3}
is also true without assuming that $k$ is algebraically closed.
We obtain in this way a weak version of 
\cite[Scholie 1.3.26(1)]{ayoub-rig}
(for the \'etale topology and after replacing $K$ with $K_{\infty}$).
See also \cite[Th\'eor\`eme 2.5.75]{ayoub-rig} for a similar
statement for motives with transfers.
\symn{$\uSH^{(\hyp)}$}
\end{rmk}}

Our first task is to construct a morphism of commutative 
algebras $\chi_B \Lambda \to \chi_{\widehat{B}}\Lambda$
which we will eventually prove to be an equivalence. 
In order to do so, we need a digression on the notion of 
rigid analytic schemes, generalising 
\cite[D\'efinition 1.4.1]{ayoub-rig}.

\begin{dfn}
\label{dfn:rigid-analytic-scheme}
A rigid analytic scheme $S$ is a triple $(S_{\eta},\widehat{S},\iota_S)$
consisting of a rigid analytic space $S_{\eta}$, 
called the generic fiber of $S$, 
a formal scheme $\widehat{S}$, called the completion of $S$, 
and an open immersion $\iota_S:\widehat{S}{}^{\rig}
\to S_{\eta}$. (We think of $S$ as obtained from 
$S_{\eta}$ and $\widehat{S}$ by gluing along $\widehat{S}{}^{\rig}$.)
Given a rigid analytic scheme $S$, we set
$S_{\sigma}=\widehat{S}_{\sigma}$ and call it the special 
fiber of $S$. A morphism of rigid analytic schemes 
$f:T\to S$ is a pair of morphisms $(f_{\eta},\widehat{f})$,
where $f_{\eta}:T_{\eta} \to S_{\eta}$ is a morphism of 
rigid analytic spaces and $\widehat{f}:\widehat{T}
\to \widehat{S}$ is a morphism of formal schemes, and such that 
$\iota_S\circ \widehat{f}_{\eta}=f_{\eta}\circ \iota_T$.
The morphism $f$ is said to be \'etale (resp. smooth)
if both $f_{\eta}$ and $\widehat{f}$ are \'etale (resp. smooth).
\ncn{rigid analytic schemes}
\end{dfn}

\begin{nota}
\label{nota:categ-of-rig-analytic-schemes}
We denote by \sym{$\RigSch$} the category of rigid analytic schemes.
Given a rigid analytic scheme $S$, we denote by 
$\RigSch/S$ the overcategory of rigid analytic $S$-schemes
and $\Et/S$ (resp. $\RigSm/S$) its full subcategory 
consisting of \'etale (resp. smooth) objects.
\symn{$\Et$}
\symn{$\RigSm$}
\end{nota}

\begin{rmk}
\label{rmk:rigid-analytic-sch-and-functors}
$\empty$

\begin{enumerate}

\item[(1)] We have a fully faithful embedding 
$\RigSpc \to \RigSch$ sending a rigid analytic 
space $S$ to the triple $(S,\emptyset,\emptyset\to S)$.
We will identify $\RigSpc$ with its essential image in $\RigSch$.

\item[(2)] We have a fully faithful embedding 
$\FSch \to \RigSch$ sending a formal scheme $\mathcal{S}$ 
to the triple $(\mathcal{S}^{\rig},\mathcal{S},
\id_{\mathcal{S}^{\rig}})$.
We will identify $\FSpc$ with its essential image in $\RigSch$.

\end{enumerate}
\end{rmk}

\begin{rmk}
\label{rmk:closed-immersion-and-complementation-rigsch}
A morphism $j$ of rigid analytic schemes is 
said to be a closed (resp. an open) immersion if both
$j_{\eta}$ and $\widehat{j}$ are closed (resp. open) immersions.
Given a closed immersion $Z \to S$ of rigid analytic schemes, 
the complement $S\smallsetminus Z$ is defined 
to be the rigid analytic scheme given by the triple
$(S_{\eta}\smallsetminus Z_{\eta},\widehat{S}\smallsetminus 
\widehat{Z},\iota_{S\smallsetminus Z})$
where $\iota_{S\smallsetminus Z}$ is obtained by restriction and 
corestriction from $\iota_S$. We have an obvious 
open immersion $S\smallsetminus Z \to S$.

We warn the reader about the following notation clash: 
given a closed immersion of formal schemes 
$\mathcal{Z}\to \mathcal{S}$, then 
``$\mathcal{S}\smallsetminus\mathcal{Z}$'' can mean two different things.
It can mean the open formal subscheme of $\mathcal{S}$
supported on the open subset $|\mathcal{S}|\smallsetminus |\mathcal{Z}|$
of $|\mathcal{S}|$. 
It can also mean the rigid analytic scheme obtained as the complement 
of $\mathcal{Z}$ in $\mathcal{S}$ considered as
rigid analytic schemes.
Each time there is a risk of confusion, we will specify if 
the complementation is taken in the category of 
formal schemes or the category of rigid analytic schemes.
\end{rmk}

Next, we generalise Construction 
\ref{cons:analytification-}.

\begin{cons}
\label{cons:analytification-scheme-3}
Let $B$ be a scheme, $B_{\sigma}\subset B$ 
a closed subscheme locally of finite presentation
up to nilimmersion, and 
$B_{\eta}\subset B$ its open complement.
There exists an analytification functor 
\begin{equation}
\label{eq-cons:analytification-scheme-37}
(-)^{\an}:\Sch^{\lft}/B \to \RigSch/\widehat{B}
\end{equation}
which is uniquely determined by the following two
properties.
\symn{$(-)^{\an}$}
\begin{enumerate}

\item[(1)] It is compatible with gluing along open immersions.

\item[(2)] For a separated finite type $B$-scheme $X$ with an open 
immersion $X\to \overline{X}$ into a proper $B$-scheme, and 
complement $Y=\overline{X}\smallsetminus X$, we have 
\begin{equation}
\label{eq-cons:scheme-analytification-}
X^{\an}=\widehat{\overline{X}}\smallsetminus 
(\widehat{Y})
\end{equation}
where, for a $B$-scheme $W$, $\widehat{W}$ is the formal 
completion of $W$ at $W_{\sigma}=W\times_B B_{\sigma}$.

\end{enumerate}
We stress that in \eqref{eq-cons:scheme-analytification-} 
the complement is taken in the category of
rigid analytic schemes.
\end{cons}

\begin{rmk}
\label{rmk:on-analytif-rigsch}
Keep the notation of Construction
\ref{cons:analytification-scheme-3}. 
The functor 
\eqref{eq-cons:analytification-scheme-37}
commutes with finite limits, and
preserves \'etale and smooth morphisms, 
closed immersions 
and complementary open immersions, as well as proper morphisms.
For $X\in \Sch^{\lft}/B$, we have a canonical isomorphism
$(X^{\an})_{\eta}\simeq (X_{\eta})^{\an}$
so there is no ambiguity in writing 
``$X^{\an}_{\eta}$''. 
The formal completions of $X$ and $X^{\an}$ 
are canonically isomorphic, i.e.,
$\widehat{X^{\an}}\simeq \widehat{X}$, 
and we have isomorphisms
$(X^{\an})_{\sigma}\simeq X_{\sigma}\simeq (X_{\sigma})^{\an}$
up to nilimmersions.
\end{rmk}

\begin{dfn}
\label{dfn:etale-topol-rigspc}
Let $(f_i:S_i \to S)_i$ be a family of \'etale morphisms of 
rigid analytic schemes. We say that this family is an \'etale 
(resp. Nisnevich) cover
if both families $(f_{i,\,\eta}:S_{i,\,\eta}\to S_{\eta})_i$ 
and $(\widehat{f}_i:\widehat{S}_i \to \widehat{S})_i$ are
\'etale (resp. Nisnevich) covers. The topology generated by 
\'etale (resp. Nisnevich) covers is called the \'etale 
(resp. Nisnevich) topology and is denoted by ``\sym{$\et$}''
(resp. ``\sym{$\Nis$}'').
\end{dfn}

\begin{nota}
\label{nota:tate-ball-over-}
Let $X$ be a rigid analytic scheme. 
We denote by $\B^n_X$ the relative $n$-dimensional ball
given by the triple $(\B^n_{X_{\eta}},\A^n_{\widehat{X}},
\id_{\B^n}\times \iota_X)$.
Similarly, we denote by $\U^1_X\subset \B^1_X$
the relative unit circle given by the triple
$(\U^1_{X_{\eta}},\A^1_{\widehat{X}}\smallsetminus 0_{\widehat{X}},
\id_{\U^1}\times \iota_X)$.
\symn{$\B^n$}
\symn{$\U^1$}
\end{nota}

\begin{dfn}
\label{dfn:motives-over-rigid-analytic-schemes}
Given a rigid analytic scheme $S$, we define the 
monoidal $\infty$-category 
of rigid analytic motives 
$\RigSH^{(\eff,\,\hyp)}_{\tau}(S;\Lambda)^{\otimes}$ 
from the smooth \'etale site 
$(\RigSm/S,\tau)$ using the interval $\B^1_S$ and the 
motive of $\U^1_S$ pointed by the unit section, 
just as in Definitions 
\ref{def:DAeff} and \ref{dfn:rigsh-stable}.
\symn{$\RigSH^{(\eff,\,\hyp)}$}
\end{dfn}

\begin{rmk}
\label{rmk:extension-to-motives-over-rig-scheme}
Many of the results that we have established for $\infty$-categories
of motives over rigid analytic spaces hold true for 
$\infty$-categories
of motives over rigid analytic schemes, and often the proof we gave
can be read in the context of rigid analytic schemes.
This is the case for instance for 
Proposition \ref{prop:6f1}. Moreover,  
Proposition \ref{prop:loc1} holds true 
for rigid analytic schemes, except that the proof of 
the localisation property requires some extra arguments.
These extra arguments can be found in the proof of 
\cite[Proposition 1.4.21]{ayoub-rig}.
Proposition 
\ref{prop:homotopic-functor} also extends: with the notation of 
Construction 
\ref{cons:analytification-scheme-3}, 
the contravariant functor 
$$X \mapsto \RigSH^{(\hyp)}_{\tau}(X^{\an};\Lambda),
\quad f\mapsto f^{\an,\,*}$$ 
from $\Sch^{\lft}/B$ to $\Prl$
is a stable homotopical functor in the sense that it satisfies
the $\infty$-categorical versions of the
properties (1)--(6) listed in \cite[\S 1.4.1]{ayoub-th1}.
\end{rmk}

Keep the notation as in Construction 
\ref{cons:analytification-scheme-3}.
Given a $B$-scheme $X$ which is locally of finite type,
the analytification functor 
\eqref{eq-cons:analytification-scheme-37}
induces a premorphism of sites
\begin{equation}
\label{eqn:analytific-funct-ran-sch}
\An_X:(\RigSm/X^{\an},\tau) \to (\Sm/X,\tau).
\end{equation}
By the functoriality of the construction of the 
$\infty$-categories of motives, \eqref{eqn:analytific-funct-ran-sch}
induces a functor
\begin{equation}
\label{eqn:analytific-funct-ran-star-sch}
\An_X^*:\SH^{(\eff,\,\hyp)}_{\tau}(X;\Lambda)
\to \RigSH^{(\eff,\,\hyp)}_{\tau}(X^{\an};\Lambda).
\end{equation}
(This generalises the functor 
\eqref{eqn:analytific-funct-ran-star}.)
Given a morphism 
$f:Y \to X$ in $\Sch^{\lft}/B$, there is an  
equivalence
$f^{\an,\,*}\circ \An^*_X\simeq \An_Y^*\circ f^*$.
In fact, the generalisation of Proposition 
\ref{prop:analytif-presheaf}
holds true: we have a morphism of 
$\CAlg(\Prl)$-valued presheaves
\begin{equation}
\label{eq-prop:analytif-presheaf-1-rigsch}
\SH_{\tau}^{(\eff,\,\hyp)}(-;\Lambda)^{\otimes}
\to \RigSH_{\tau}^{(\eff,\,\hyp)}((-)^{\an};\Lambda)^{\otimes}
\end{equation}
on $\Sch^{\lft}/B$.
Also, note that if $Z$ is a $B_{\sigma}$-scheme which is locally 
of finite type, then 
$\An^*_Z$ is an equivalence of $\infty$-categories.
\symn{$\An^*$}

\begin{nota}
\label{nota:chi-B-specialisation-system}
Let $B$ be a scheme, $B_{\sigma}\subset B$ 
a closed subscheme locally of finite presentation, and 
$B_{\eta}\subset B$ its open complement.
\begin{enumerate}

\item[(1)] Given a $B$-scheme $X$, 
we set $X_{\sigma}=X\times_B B_{\sigma}$ and 
$X_{\eta}=X\times_B B_{\eta}$, and we define the 
functor 
\begin{equation}
\label{eq-nota:chi-B-specialisation-system}
\chi_X:\SH^{(\eff,\,\hyp)}_{\tau}(X_{\eta};\Lambda)
\to \SH^{(\eff,\,\hyp)}_{\tau}(X_{\sigma};\Lambda)
\end{equation}
as in the statement of Theorem 
\ref{thm:compute-chi}.
More precisely, we denote by $i:X_{\sigma}\to X$ and 
$j:X_{\eta}\to X$ the obvious inclusions, and set 
$\chi_X=i^*\circ j_*$.
\symn{$\chi$}

\item[(2)] Given a rigid analytic $\widehat{B}$-scheme $X$, 
we define the functor 
\begin{equation}
\label{eq-nota:chi-B-specialisation-system-rigsch}
\chi_X:\RigSH^{(\eff,\,\hyp)}_{\tau}(X_{\eta};\Lambda)
\to \SH^{(\eff,\,\hyp)}_{\tau}(X_{\sigma};\Lambda)
\end{equation}
similarly.
More precisely, we denote by $i:X_{\sigma}\to X$ and 
$j:X_{\eta}\to X$ the obvious inclusions, and set
$\chi_X=i^*\circ j_*$.
\end{enumerate}
\end{nota}

\begin{rmk}
\label{rmk:chi-X-specialisation-system}
In the $\Tate$-stable case, 
the collection of functors $\{\chi_X\}_X$, for $X\in \Sch/B$, 
is part of a specialisation system in the sense of 
\cite[D\'efinition 3.1.1]{ayoub-th2}.
In fact, this specialisation system is considered in 
\cite[Exemple 3.1.4]{ayoub-th2} where it is called 
the \nc{canonical specialisation system}. 
Similarly, the collection of functors 
$\{\chi_{X^{\an}}\circ \An_{X_{\eta}}^*\}_X$, 
for $X\in \Sch/B$, 
is part of a specialisation system; see
\cite[Proposition 1.4.41]{ayoub-rig}.
There are natural transformations 
\begin{equation}
\label{eqn:morphism-of-specialisation-systems}
\rho_X:\chi_X \to \chi_{X^{\an}}\circ \An_{X_{\eta}}^*,
\end{equation}
given by the composition of
$$\chi_X=i^*\circ j_* \simeq \An_{X_{\sigma}}^*\circ i^*\circ j_* 
\simeq i^{\an,\,*}\circ \An_X^* \circ j_*  \to i^{\an,\,*} \circ
j^{\an}_* \circ \An_{X_{\eta}}^* 
\simeq \chi_{X^{\an}}\circ \An_{X_{\eta}}^*,$$
which
are part of a morphism of specialisation systems;
see \cite[Lemme 1.4.42]{ayoub-rig}.
\symn{$\rho$}
\end{rmk}

\begin{rmk}
\label{rmk:independence-of-B-}
The natural transformation $\rho_X$
is independent of $B$ in the following way.
Let $B'\in \Sch^{\lft}/B$ and $X\in \Sch^{\lft}/B'$. 
Then we have two natural transformations 
``$\chi_X\to \chi_{X^{\an}}\circ \An^*_{X_{\eta}}$'',
one associated with $X$ considered as a $B$-scheme
and one associated with $X$ considered as a $B'$-scheme.
We claim that these two natural transformations 
are equivalent. To explain how, we write momentarily 
$\chi_{(X/B)^{\an}}$, $\An^*_{X_{\eta}/B}$, etc., 
to stress the dependency on the scheme $B$.
There is a canonical isomorphism 
$$(X/B')^{\an}\simeq 
(X/B)^{\an}\times_{(B'/B)^{\an}}\widehat{B}{}',$$
and hence an open immersion of rigid analytic 
$\widehat{B}$-schemes $\iota:(X/B')^{\an} \to (X/B)^{\an}$
inducing an isomorphism on special fibers.
Moreover, we have natural equivalences
$$\chi_{(X/B)^{\an}}\simeq \chi_{(X/B')^{\an}}\circ 
\iota_{\eta}^* \qquad \text{and} \qquad 
\An^*_{X_{\eta}/B'}\simeq \iota_{\eta}^*\circ \An^*_{X_{\eta}/B}.$$
Modulo these equivalences, the two natural transformations 
``$\chi_X\to \chi_{X^{\an}}\circ \An^*_{X_{\eta}}$''
give the same natural transformation 
$\chi_X \to \chi_{(X/B')^{\an}}\circ \iota_{\eta}^*\circ
\An^*_{X_{\eta}/B}$.
\end{rmk}

\begin{lemma}
\label{lem:description-chi-rigsch}
Let $X$ be a rigid analytic $\widehat{B}$-scheme.
The functor \eqref{eq-nota:chi-B-specialisation-system-rigsch}
is equivalent to the composition of
$$\RigSH^{(\eff,\,\hyp)}_{\tau}(X_{\eta};\Lambda)
\xrightarrow{\iota_X^*} 
\RigSH^{(\eff,\,\hyp)}_{\tau}(\widehat{X}{}^{\rig};\Lambda)
\xrightarrow{\chi_{\widehat{X}}} \SH^{(\eff,\,\hyp)}_{\tau}(X_{\sigma};\Lambda),$$
where $\chi_{\widehat{X}}$ is the functor introduced in Notation
\ref{nota:generic-fiber-eta-star}. 
\end{lemma}

\begin{proof}
For the sake of clarity, we will momentarily 
write ``$\chi'_{\widehat{X}}$'' instead of 
``$\chi_{\widehat{X}}$'' for the functor 
introduced in Notation
\ref{nota:generic-fiber-eta-star} and use 
``$\chi_{\widehat{X}}$''
to denote the functor introduced in 
Notation 
\ref{nota:chi-B-specialisation-system}(2)
with $\widehat{X}$
considered as a rigid analytic $\widehat{B}$-scheme
via the fully faithful embedding $\FSch\to \RigSch$.

We have an equivalence
$\chi_X \simeq \chi_{\widehat{X}}\circ \iota_X^*$
which follows from the fact that $(\iota_X)_{\sigma}$ is the 
identification $\widehat{X}_{\sigma}\simeq X_{\sigma}$.
Thus, to prove the lemma, it is enough to show that the
two functors 
$$\chi_{\widehat{X}},\;\;\chi'_{\widehat{X}}:
\RigSH^{(\eff,\,\hyp)}_{\tau}(\widehat{X}_{\eta};\Lambda)
\to \SH^{(\eff,\,\hyp)}_{\tau}(X_{\sigma};\Lambda)$$
are equivalent. (Note that $\widehat{X}_{\eta}=\widehat{X}{}^{\rig}$; here 
we use ``$\widehat{X}_{\eta}$'' because we want to think 
about $\widehat{X}$ as a rigid analytic scheme via the fully faithful
embedding of Remark
\ref{rmk:rigid-analytic-sch-and-functors}(2).) 
In order to do that, we remark that the 
base change functor $\RigSm/\widehat{X} \to \Sm/X_{\sigma}$
factors as follows
$$\RigSm/\widehat{X} \xrightarrow{\widehat{(-)}} 
\FSm/\widehat{X} \xrightarrow{(-)_{\sigma}} \Sm/X_{\sigma}.$$
We deduce immediately from the construction of the 
$\infty$-categories of motives that the inverse image functor
$i^*:\RigSH^{(\eff,\,\hyp)}_{\tau}(\widehat{X};\Lambda)
\to \SH^{(\eff,\,\hyp)}_{\tau}(X_{\sigma};\Lambda)$
is the composition of 
$$\RigSH^{(\eff,\,\hyp)}_{\tau}(\widehat{X};\Lambda)
\xrightarrow{\widehat{(-)}{}^*}
\FSH^{(\eff,\,\hyp)}_{\tau}(\widehat{X};\Lambda)
\xrightarrow{\sigma^*} 
\SH^{(\eff,\,\hyp)}_{\tau}(X_{\sigma};\Lambda)$$
where $\sigma^*$ is the equivalence of 
Theorem \ref{thm:formal-mot-alg-mot} and 
$\widehat{(-)}{}^*$ is the functor 
that takes the motive of a rigid analytic 
$\widehat{X}$-scheme to the motive of its formal completion.
The formal completion functor $\widehat{(-)}$ 
is right adjoint to the obvious inclusion 
$\inc:\FSm/\widehat{X} \to \RigSch/\widehat{X}$.
It follows that $\widehat{(-)}{}^*$ is right adjoint
to the functor
$$\inc^*:\FSH^{(\eff,\,\hyp)}_{\tau}(\widehat{X};\Lambda)
\to \RigSH^{(\eff,\,\hyp)}_{\tau}(\widehat{X};\Lambda).$$ 
This means that we have an equivalence 
$\widehat{(-)}{}^*\simeq \inc_*$.
In conclusion, we see that 
$\chi_{\widehat{X}}$ is equivalent to the composition of 
$$\RigSH^{(\eff,\,\hyp)}_{\tau}(\widehat{X}_{\eta};\Lambda)
\xrightarrow{j_*}
\RigSH^{(\eff,\,\hyp)}_{\tau}(\widehat{X};\Lambda)
\xrightarrow{\inc_*}
\FSH^{(\eff,\,\hyp)}_{\tau}(\widehat{X};\Lambda)
\xrightarrow{\sigma^*}
\SH^{(\eff,\,\hyp)}_{\tau}(X_{\sigma};\Lambda).$$
Since $j^*\circ \inc^*$ is clearly equivalent to the functor 
$\xi_{\widehat{X}}$ from Notation 
\ref{nota:generic-fiber-eta-star}, the result follows.
\end{proof}

\begin{cor}
\label{cor:chi-widehat-B}
The functor $\chi_{\widehat{B}}$ obtained by 
taking $X=\widehat{B}$ in Notation
\ref{nota:chi-B-specialisation-system}(2) 
coincides with the functor 
$\chi_{\widehat{B}}$ obtained by taking 
$\mathcal{S}=\widehat{B}$ in Notation 
\ref{nota:generic-fiber-eta-star}.
\end{cor}

From Corollary \ref{cor:chi-widehat-B}, we see that 
Theorem \ref{thm:compute-chi}
follows from the following statement.

\begin{thm}
\label{thm:morphism-of-specialisation-equival}
Let $B$ be a scheme, $B_{\sigma}\subset B$ 
a closed subscheme locally of finite presentation up to nilimmersion, 
and $B_{\eta}\subset B$ its open complement.
Assume that every prime number is invertible either 
in $\pi_0\Lambda$ or in $\mathcal{O}(B)$.
Assume one of the following alternatives.
\begin{enumerate}

\item[(1)] We work in the non-hypercomplete case and
$\Lambda$ is eventually coconnective;

\item[(2)] We work in the hypercomplete case and $B$ is 
$(\Lambda,\et)$-admissible.

\end{enumerate}
Then, for every $X\in \Sch^{\lft}/B$, the natural transformation 
$\rho_X:\chi_X \to \chi_{X^{\an}}\circ \An_{X_{\eta}}^*$, 
between functors from $\SH^{(\hyp)}_{\et}(X_{\eta};\Lambda)$
to $\SH^{(\hyp)}_{\et}(X_{\sigma};\Lambda)$,
is an equivalence.
\end{thm}

We start by proving a reduction.

\begin{lemma}
\label{lem:reduct-for-iso-specialisation-systems}
To prove Theorem 
\ref{thm:morphism-of-specialisation-equival}, 
we may assume that $\Lambda$ is eventually coconnective 
and that $B$ is essentially of finite type over 
$\Spec(\Z)$. In particular, there is no need to 
distinguish the non-hypercomplete and the hypercomplete cases.
\end{lemma}

\begin{proof}
We first explain how to reduce to the case where
$\Lambda$ is eventually coconnective.
For this, we only need to consider the alternative (2). 
It follows from Propositions
\ref{prop:compact-shv-rigsm} and
\ref{prop:compact-shv-forsm} that 
$\rho_X$ is a natural transformation between 
colimit-preserving functors between compactly generated 
categories. Thus, it is enough to prove that 
$\chi_XM \to \chi_{X^{\an}}\An_{X_{\eta}}^*M$
is an equivalence for $M\in \SH^{\hyp}_{\et}(X_{\eta};\Lambda)$
compact. 
Arguing as in the second part of the proof of Lemma
\ref{lem:alternative-iii-vs-iv-3i}, we reduce to  
the following two cases:
\begin{itemize}

\item $\pi_0\Lambda$ is a $\Q$-algebra;

\item $M$ is $\ell$-nilpotent for a prime $\ell$ invertible on $B$.

\end{itemize}
In the first case, we may replace $\Lambda$ by $\Q$
and assume that $\Lambda$ is eventually coconnective
as claimed.
In the second case, let 
$M_0\in \Shv^{\hyp}_{\et}(\Et/X_{\eta};\Lambda)_{\ell}$
be the object corresponding to $M$ by the equivalence 
$$\Shv^{\hyp}_{\et}(\Et/X_{\eta};\Lambda)_{\ellnil}
\simeq 
\SH^{\hyp}_{\et}(X_{\eta};\Lambda)_{\ellnil}$$
provided by Theorem 
\ref{thm:rigrig-algebraic}. 
Using also Theorem
\ref{thm:rigrig}, we reduce to show that 
$\chi_XM_0 \to \chi_{X^{\an}}\An_{X_{\eta}}^*M_0$
is an equivalence. (Here the functors $\chi_X$, 
$\chi_{X^{\an}}$ and $\An_{X_{\eta}}^*$ are 
defined on \'etale hypersheaves of $\Lambda$-modules 
by the same formulas as their motivic versions.)
Using Lemma
\ref{lem:pi-0-Lambda-coh-dim}, 
one obtains equivalences 
$$\chi_XM_0\simeq \lim_r\chi_X(M_0\otimes_{\Lambda}
\tau_{\leq r}\Lambda)
\qquad \text{and} \qquad 
\chi_{X^{\an}}\An_{X_{\eta}}^*M_0
\simeq \lim_r 
\chi_{X^{\an}}\An_{X_{\eta}}^*
(M_0\otimes_{\Lambda}\tau_{\leq r}\Lambda).$$
(Indeed, as $M_0$ is compact, the inverse system 
$(M_0\otimes_{\Lambda}\tau_{\leq r}\Lambda)_r$ 
consists of eventually coconnective \'etale sheaves and
is eventually constant on homotopy sheaves.)
This shows that we may replace $M$ and $\Lambda$ by 
$M\otimes_{\Lambda}\tau_{\leq r}\Lambda$ and 
$\tau_{\leq r}\Lambda$, and assume that $\Lambda$ 
is eventually coconnective as claimed.

We now assume that $\Lambda$ is eventually coconnective
and explain how to reduce to the case where $B$ is essentially 
of finite type over $\Spec(\Z)$.
By Propositions \ref{prop:automatic-hypercomp-motives}
and \ref{prop:automatic-hypercomp-motives-algebraic}, 
the alternative (2) is covered by the alternative (1). 
By Remark \ref{rmk:independence-of-B-}, 
we only need to consider the case $X=B$.
The problem is local on $B$, so we may assume that 
$B$ is affine given as a limit of a cofiltered inverse system
$(B_{\alpha})_{\alpha}$ of affine schemes 
which are essentially of finite type over 
$\Z$. We may also assume that there are closed subschemes 
$B_{\alpha,\,\sigma}\subset B_{\alpha}$
such that, for every $\beta\leq \alpha$, 
$B_{\beta,\,\sigma}$ is the inverse image of $B_{\alpha,\,\sigma}$, 
and $B_{\sigma}$ is the limit of the inverse system 
$(B_{\alpha,\,\sigma})_{\alpha}$. 
Set $B_{\alpha,\,\eta}=B_{\alpha}\smallsetminus B_{\alpha,\,\sigma}$
so that $B_{\eta}$ is the limit of the inverse system 
$(B_{\alpha,\,\eta})_{\alpha}$. Let
$i_{\alpha}:B_{\alpha,\,\sigma}\to B_{\alpha}$ and 
$j_{\alpha}:B_{\alpha,\,\eta}\to B_{\alpha}$ 
be the obvious immersions, and let
$f_{\alpha}:B\to B_{\alpha}$ and 
$f_{\beta\alpha}:B_{\beta} \to B_{\alpha}$
be the obvious morphisms. 

We need to show that $\chi_BM\to 
\chi_{\widehat{B}}\,\An^*_{X_{\eta}}M$ is an equivalence 
for all $M\in \SH_{\et}(B_{\eta};\Lambda)$.
Since the three functors 
$\chi_B$, $\chi_{\widehat{B}}$ and $\An^*_{X_{\eta}}$
commute with colimits, we may assume that 
$M$ is compact.
By Proposition \ref{prop:cont-algebraic-29}, 
we have an equivalence
$$\SH_{\et}(B;\Lambda)\simeq \underset{\alpha}{\colim}\,
\SH_{\et}(B_{\alpha};\Lambda)$$
in $\Prl$, and similarly for $B_{\sigma}$ and $B_{\eta}$.
Since $M\in \SH_{\et}(B_{\eta};\Lambda)$ 
is assumed compact, we may find an index $\alpha_0$,
a compact object $M_{\alpha_0}\in 
\SH_{\et}(B_{\alpha_0,\,\eta};\Lambda)$
and an equivalence $f_{\alpha_0,\,\eta}^*M_{\alpha_0} \simeq M$.
We set $M_{\alpha}=f_{\alpha\alpha_0,\,\eta}^*M$.
With this, we have an equivalence 
$$j_*M\simeq \underset{\alpha\leq \alpha_0}
{\colim}\; f^*_{\alpha} j_{\alpha,\,*}M_{\alpha}.$$
(It is not totally obvious how to construct such an equivalence.
One needs to argue as in the proof of Lemma \ref{lem:colimi-chi-};
see also Remark \ref{rmk:extension-lemma-colimi-chi-}.) 
Applying $i^*$, we deduce an equivalence
\begin{equation}
\label{eq-lem:reduct-for-iso-specialisation-systems-7}
\chi_BM\simeq \underset{\alpha\leq \alpha_0}{\colim}\;
f_{\alpha}^*\chi_{B_{\alpha}}M_{\alpha}.
\end{equation}
Similarly, by Remark \ref{rmk:extension-lemma-colimi-chi-}
and using Corollary 
\ref{cor:chi-widehat-B},
we have an equivalence 
\begin{equation}
\label{eq-lem:reduct-for-iso-specialisation-systems-11}
\chi_{\widehat{B}}\,\An^*_{B_{\eta}}M
\simeq \underset{\alpha\leq \alpha_0}{\colim}\;
f_{\alpha}^{\an,\,*}\chi_{\widehat{B}_{\alpha}}
\An^*_{B_{\alpha,\,\eta}} M_{\alpha}.
\end{equation}
Therefore, it is enough to show that 
$\chi_{B_{\alpha}}M_{\alpha}
\to \chi_{\widehat{B}_{\alpha}} 
\An^*_{B_{\alpha,\,\eta}}M_{\alpha}$
is an equivalence.
In particular, we may assume that $B$ is quasi-excellent 
and $(\Lambda,\et)$-admissible. 
In this case, since $\Lambda$ is eventually coconnective, 
we are automatically working in the hypercomplete case
by Propositions \ref{prop:automatic-hypercomp-motives} and
\ref{prop:automatic-hypercomp-motives-algebraic}.
This finishes the proof.
\end{proof}

Our next task is to prove the following 
weak version of Theorem 
\ref{thm:morphism-of-specialisation-equival}
(which we are able to justify even when $\tau$ is the Nisnevich
topology).

\begin{prop}
\label{prop-thm:morphism-of-specialisation-equival-weak}
Let $B$ be a quasi-excellent $(\Lambda,\tau)$-admissible 
scheme, $B_{\sigma}\subset B$ 
a closed subscheme, and 
$B_{\eta}\subset B$ its open complement.
If $\tau$ is the \'etale topology, 
assume that every prime number is invertible either 
in $\pi_0\Lambda$ or in $\mathcal{O}(B)$.
Then, there is a natural transformation 
$\chi_{\widehat{B}}\circ \An^*_{B_{\eta}}\to \chi_B$, 
between functors from $\SH^{\hyp}_{\tau}(B_{\eta};\Lambda)$
to $\SH^{\hyp}_{\tau}(B_{\sigma};\Lambda)$,
which is a section to the natural transformation 
$\rho_B$, i.e., such that the composition of 
$$\chi_{\widehat{B}}\circ \An^*_{B_{\eta}}\to \chi_B
\xrightarrow{\rho_B} 
\chi_{\widehat{B}}\circ \An^*_{B_{\eta}}$$
is the identity.

\end{prop}

To prove Proposition 
\ref{prop-thm:morphism-of-specialisation-equival-weak}
we need a digression. (Compare with 
\cite[page 112]{ayoub-rig}.\footnote{We remind the reader that
the page references to \cite{ayoub-rig} correspond to the 
published version.})

\begin{cons}
\label{cons:arrow-site+}
Let $\mathcal{S}$ be a formal scheme. 
We denote by $\FRigSm_{\af}/\mathcal{S}$ 
the full subcategory of $\FSch/\mathcal{S}$ 
spanned by rig-smooth formal $\mathcal{S}$-schemes 
which are affine. Consider the functor  
\begin{equation}
\label{eq-dfn:arrow-site-1+}
\mathfrak{D}_{\mathcal{S}}:
\FRigSm_{\af}/\mathcal{S} \to \Sch
\end{equation}
sending an affine formal scheme $\Spf(A)$ over 
$\mathcal{S}$ to the scheme $\Spec(A)$.
Consider also the two related functors 
$\mathfrak{D}_{\mathcal{S},\,\sigma}$ and 
$\mathfrak{D}_{\mathcal{S},\,\eta}$
between the same categories, sending 
an affine formal scheme $\Spf(A)$ over 
$\mathcal{S}$ to the schemes $\Spf(A)_{\sigma}$ and 
$\Spec(A)\smallsetminus \Spf(A)_{\sigma}$
respectively. 
We consider $\mathfrak{D}_{\mathcal{S}}$, $\mathfrak{D}_{\mathcal{S},\,\sigma}$ and $\mathfrak{D}_{\mathcal{S},\,\eta}$
as diagrams of schemes and define the smooth $\tau$-sites
$$(\Sm/\mathfrak{D}_{\mathcal{S}},\tau),\quad 
(\Sm/\mathfrak{D}_{\mathcal{S},\,\sigma},\tau) \quad \text{and} \quad 
(\Sm/\mathfrak{D}_{\mathcal{S},\,\eta},\tau)$$
as in \cite[\S 4.5.1]{ayoub-th2}.
To fix the notation, let us recall that an object of 
$\Sm/\mathfrak{D}_{\mathcal{S}}$ 
is a pair $(\mathcal{U},V)$ consisting of an object
$\mathcal{U}\in \FRigSm_{\af}/\mathcal{S}$ and 
a smooth $\mathcal{O}(\mathcal{U})$-scheme $V$. The topology 
$\tau$ on $\Sm/\mathfrak{D}_{\mathcal{S}}$
is generated by families of the form
$((\id_{\mathcal{U}},e_i):
(\mathcal{U},V_i)\to (\mathcal{U},V))_i$ 
where the family $(e_i)_i$ is a cover for the topology $\tau$.
\symn{$\FRigSm_{\af}$}
\symn{$\mathfrak{D}$}

The $\infty$-category $\SH^{(\eff,\,\hyp)}_{\tau}
(\mathfrak{D}_{\mathcal{S}};\Lambda)$ is constructed from the site
$(\Sm/\mathfrak{D}_{\mathcal{S}},\tau)$, using the interval
$\A^1$ and the motive of $\A^1\smallsetminus 0$ pointed by 
the unit section, as in Definitions \ref{def:DAeff} and 
\ref{dfn:rigsh-stable} (or Definition 
\ref{def:DAeff-form} and \ref{dfn:fsh-stable}), and similarly 
for $\mathfrak{D}_{\mathcal{S},\,\sigma}$ and 
$\mathfrak{D}_{\mathcal{S},\,\eta}$. (For a construction using the 
language of model categories, see \cite[\S 4.5.2]{ayoub-th2}.)
We note here that $\A^1$ (resp. $\A^1\smallsetminus 0$) 
is considered as a presheaf of sets on 
$\Sm/\mathfrak{D}_{\mathcal{S}}$, sending 
$(\mathcal{U},V)$ to $\mathcal{O}(V)$ 
(resp. $\mathcal{O}^{\times}(V)$).
This presheaf is not representable unless $\mathcal{S}$ 
is affine, but the Cartesian product with this presheaf
preserves representable presheaves. (For instance, 
we have $\A^1\times (\mathcal{U},V)=(\mathcal{U},\A^1_V)$.)
We have morphisms of diagrams of schemes
$\mathfrak{i}:\mathfrak{D}_{\mathcal{S},\,\sigma} \to 
\mathfrak{D}_{\mathcal{S}}$ and 
$\mathfrak{j}:\mathfrak{D}_{\mathcal{S},\,\eta} \to 
\mathfrak{D}_{\mathcal{S}}$, 
and we define the functor 
\begin{equation}
\label{eq-dfn:arrow-site-13+}
\chi_{\mathfrak{D}_{\mathcal{S}}}:
\SH^{(\eff,\,\hyp)}_{\tau}(\mathfrak{D}_{\mathcal{S},\,\eta};\Lambda)
\to \SH^{(\eff,\,\hyp)}_{\tau}
(\mathfrak{D}_{\mathcal{S},\,\sigma};\Lambda)
\end{equation}
to be the composite $\mathfrak{i}^*\circ \mathfrak{j}_*$.

Similarly, consider the functor 
\begin{equation}
\label{eq-dfn:arrow-site-9+}
\mathfrak{D}^{\an}_{\mathcal{S}}:
\FRigSm_{\af}/\mathcal{S} \to \RigSch
\end{equation}
sending an affine formal scheme $\Spf(A)$ over 
$\mathcal{S}$ to $\Spf(A)$ considered as a rigid analytic scheme.
Consider also the related functor
$\mathfrak{D}^{\an}_{\mathcal{S},\,\eta}$
between the same categories, sending 
an affine formal scheme $\Spf(A)$ over 
$\mathcal{S}$ to the rigid analytic space 
$\Spf(A)^{\rig}$. 
We consider $\mathfrak{D}^{\an}_{\mathcal{S}}$ and 
$\mathfrak{D}^{\an}_{\mathcal{S},\,\eta}$
as diagrams of rigid analytic schemes 
and define the smooth $\tau$-sites
$(\RigSm/\mathfrak{D}^{\an}_{\mathcal{S}},\tau)$ and 
$(\RigSm/\mathfrak{D}^{\an}_{\mathcal{S},\,\eta},\tau)$
as in \cite[\S 4.5.1]{ayoub-th2}.
The $\infty$-category $\RigSH^{(\eff,\,\hyp)}_{\tau}
(\mathfrak{D}^{\an}_{\mathcal{S}};\Lambda)$ 
is constructed from the site
$(\RigSm/\mathfrak{D}^{\an}_{\mathcal{S}},\tau)$, 
using the interval $\B^1$ and the motive of 
$\U^1$ pointed by the unit section, as in 
Definitions \ref{def:DAeff} and 
\ref{dfn:rigsh-stable}, and similarly for 
$\mathfrak{D}^{\an}_{\mathcal{S},\,\eta}$.
We have morphisms of diagrams of rigid analytic schemes
$\mathfrak{i}^{\an}:\mathfrak{D}_{\mathcal{S},\,\sigma} \to 
\mathfrak{D}_{\mathcal{S}}^{\an}$ and 
$\mathfrak{j}^{\an}:\mathfrak{D}^{\an}_{\mathcal{S},\,\eta} \to 
\mathfrak{D}_{\mathcal{S}}^{\an}$, 
and we define the functor 
\begin{equation}
\label{eq-dfn:arrow-site-137+}
\chi_{\mathfrak{D}^{\an}_{\mathcal{S}}}:
\RigSH^{(\eff,\,\hyp)}_{\tau}
(\mathfrak{D}^{\an}_{\mathcal{S},\,\eta};\Lambda)
\to \SH^{(\eff,\,\hyp)}_{\tau}
(\mathfrak{D}_{\mathcal{S},\,\sigma};\Lambda)
\end{equation}
to be the composite 
$\mathfrak{i}^{\an,\,*}\circ \mathfrak{j}^{\an}_*$.
The analytification functor induces functors 
$$\An^*_{\mathfrak{D}_{\mathcal{S}}}:
\SH^{(\eff,\,\hyp)}_{\tau}(\mathfrak{D}_{\mathcal{S}};\Lambda)
\to \RigSH^{(\eff,\,\hyp)}_{\tau}
(\mathfrak{D}^{\an}_{\mathcal{S}};\Lambda) \qquad \text{and}$$
$$\An^*_{\mathfrak{D}_{\mathcal{S},\,\eta}}:
\SH^{(\eff,\,\hyp)}_{\tau}(\mathfrak{D}_{\mathcal{S},\,\eta};\Lambda)
\to \RigSH^{(\eff,\,\hyp)}_{\tau}
(\mathfrak{D}^{\an}_{\mathcal{S},\,\eta};\Lambda).$$
We may then define a natural transformation 
\begin{equation}
\label{eq-dfn:arrow-site-13708+}
\rho_{\mathfrak{D}_{\mathcal{S}}}:
\chi_{\mathfrak{D}_{\mathcal{S}}}\to 
\chi_{\mathfrak{D}^{\an}_{\mathcal{S}}}\circ 
\An^*_{\mathfrak{D}_{\eta}}
\end{equation}
as in Remark \ref{rmk:chi-X-specialisation-system}.
\end{cons}

\begin{rmk}
\label{rmk:diag-frak-D-formal-schemes}
The functor 
\eqref{eq-dfn:arrow-site-9+}
factors through the subcategory 
$\FSch\subset \RigSch$ and defines a 
diagram of formal schemes that we denote by 
$\mathfrak{D}^{\form}_{\mathcal{S}}$. 
As in Construction 
\ref{cons:arrow-site+}, we can define an 
$\infty$-category $\FSH^{(\eff,\,\hyp)}_{\tau}
(\mathfrak{D}^{\form}_{\mathcal{S}};\Lambda)$
of formal motives over 
$\mathfrak{D}^{\form}_{\mathcal{S}}$
using the smooth site 
$(\FSm/\mathfrak{D}^{\form}_{\mathcal{S}},\tau)$.
Moreover, we have an equivalence of 
$\infty$-categories
$$\sigma^*:\FSH^{(\eff,\,\hyp)}_{\tau}
(\mathfrak{D}^{\form}_{\mathcal{S}};\Lambda)
\xrightarrow{\sim} \SH^{(\eff,\,\hyp)}_{\tau}
(\mathfrak{D}_{\mathcal{S},\,\sigma};\Lambda)$$
as in Theorem 
\ref{thm:formal-mot-alg-mot}.
\symn{$\mathfrak{D}^{\form}$}
\end{rmk}

\begin{lemma}
\label{lem:chi-mathfrak-D-an}
The functor 
$\chi_{\mathfrak{D}_{\mathcal{S}}^{\an}}$ 
coincides with the composition of 
$$\RigSH^{(\eff,\,\hyp)}_{\tau}
(\mathfrak{D}^{\an}_{\mathcal{S},\,\eta};\Lambda) 
\xrightarrow{\chi_{\mathfrak{D}^{\form}_{\mathcal{S}}}}
\FSH^{(\eff,\,\hyp)}_{\tau}
(\mathfrak{D}^{\form}_{\mathcal{S}};\Lambda)
\xrightarrow{\sigma^*}
\SH^{(\eff,\,\hyp)}_{\tau}
(\mathfrak{D}_{\mathcal{S},\,\sigma};\Lambda)$$
where $\chi_{\mathfrak{D}^{\form}_{\mathcal{S}}}$
is the restriction along the functor  
$(-)^{\rig}:
\FSm/\mathfrak{D}^{\form}_{\mathcal{S}}
\to 
\RigSm/\mathfrak{D}^{\an}_{\mathcal{S},\,\eta}$
sending a pair 
$(\mathcal{U},\mathcal{V})$ to 
$(\mathcal{U},\mathcal{V}^{\rig})$.
\end{lemma}

\begin{proof}
This is diagrammatic version of Lemma 
\ref{lem:description-chi-rigsch}
which is proven in the same way.
\end{proof}

\begin{rmk}
\label{rmk:diagonal-functors-}
There are five diagonal functors 
emanating from $\FRigSm_{\af}/\mathcal{S}$ 
and taking values in the categories 
$\Sm/\mathfrak{D}_{\mathcal{S}}$, 
$\Sm/\mathfrak{D}_{\mathcal{S},\,\sigma}$,
$\Sm/\mathfrak{D}_{\mathcal{S},\,\eta}$,
$\RigSm/\mathfrak{D}^{\an}_{\mathcal{S}}$ and 
$\RigSm/\mathfrak{D}^{\an}_{\mathcal{S},\,\eta}$.
These functors will be denoted respectively by 
\sym{$\diag$}, $\diag_{\sigma}$, $\diag_{\eta}$, $\diag^{\an}$ 
and $\diag^{\an}_{\eta}$. 
They send an affine formal scheme 
$\mathcal{U}=\Spf(A)$ over $\mathcal{S}$ to the pairs
$(\mathcal{U},\Spec(A))$, $(\mathcal{U},\mathcal{U}_{\sigma})$,
$(\mathcal{U},\Spec(A)\smallsetminus \mathcal{U}_{\sigma})$,
$(\mathcal{U},\mathcal{U})$ and $(\mathcal{U},\mathcal{U}^{\rig})$
respectively. We now concentrate on the case of 
${\rm diag}$, but what we are going to say 
can be adapted to the remaining four diagonal functors. 
The functor ${\rm diag}$ induces an adjunction 
$$\diag^*:
\PSh(\FRigSm_{\af}/\mathcal{S};\Lambda)
\rightleftarrows 
\PSh(\Sm/\mathfrak{D}_{\mathcal{S}};\Lambda):
\diag_*$$
where $\diag_*$ is the restriction functor.
As in Remark \ref{rmk:rigda-model-cat},
we denote by $T_{\mathcal{S}}$ (instead of $\Tate_{\mathcal{S}}$)
the cofiber of the split inclusion of $\Lambda(\mathcal{S})
\to \Lambda(\A^1_{\mathcal{S}}\smallsetminus 0_{\mathcal{S}})$
(without $\tau$-(hyper)sheafification), and similarly 
for $T_{\mathfrak{D}_{\mathcal{S}}}$. (Here $\mathcal{S}$ and 
$\A^1_{\mathcal{S}}\smallsetminus 0_{\mathcal{S}}$ are considered
as presheaves of sets on $\FRigSm_{\af}/\mathcal{S}$
which are not necessarily representable.)
Noting that $\diag_*(T_{\mathfrak{D}_{\mathcal{S}}})
\simeq T_{\mathcal{S}}$, we may extend the above adjunction 
to $T$-spectra:
$$\diag^*:
\Spect_T(\PSh(\FRigSm_{\af}/\mathcal{S};\Lambda))
\rightleftarrows 
\Spect_T(\PSh(\Sm/\mathfrak{D}_{\mathcal{S}};\Lambda)):
\diag_*.$$
Here, by abuse of notation, we write 
$\Spect_T(\PSh(-;\Lambda))$ for the 
$\infty$-category associated to the simplicial category
$\Spect_T(\PSh_{\Delta}(-;\Lambda))$
endowed with its levelwise global model structure; 
compare with Remark \ref{rmk:rigda-model-cat}.
We have the following equivalences
$$\diag_{\sigma,\,*}\simeq \diag_*\circ \mathfrak{i}_*
\simeq \diag^{\an}_*\circ \mathfrak{i}^{\an}_*, \quad
\diag_{\eta,\,*}\simeq \diag_*\circ \mathfrak{j}_*
\quad \text{and} \quad
\diag^{\an}_{\eta,\,*}\simeq \diag^{\an}_*\circ \mathfrak{j}^{\an}_*.$$
Moreover, there are natural equivalences
$\An_{\mathfrak{D}_{\mathcal{S}}}^*\circ \diag^*
\simeq \diag^{\an,\,*}$ and 
$\An_{\mathfrak{D}_{\mathcal{S},\,\eta}}^*\circ \diag_{\eta}^*
\simeq \diag_{\eta}^{\an,\,*}$
inducing natural transformations 
\begin{equation}
\label{eq-rmk:diagonal-functors-}
\diag_* \to \diag^{\an}_*\circ \An_{\mathfrak{D}_{\mathcal{S}}}^*
\qquad \text{and} \qquad 
\diag_{\eta,\,*} \to \diag^{\an}_{\eta,\,*}\circ 
\An_{\mathfrak{D}_{\mathcal{S},\,\eta}}^*.
\end{equation}
\end{rmk}

\begin{lemma}
\label{lem:diag-eta-on-representable}
Below, we consider $\diag_{\eta,\,*}$ and $\diag_{\eta,\,*}^{\an}$
as ordinary functors on ordinary categories of presheaves of sets.
Given a rigid analytic space $W$ over $\mathcal{S}^{\rig}$, 
we denote also by $W$ the presheaf of sets on 
$\FRigSm_{\af}/\mathcal{S}$ 
given by $W(\mathcal{X})=\Hom_{\mathcal{S}^{\rig}}
(\mathcal{X}^{\rig},W)$.
\begin{enumerate}

\item[(1)] Let $(\mathcal{U},V)$ be an object of 
$\Sm/\mathfrak{D}_{\mathcal{S},\,\eta}$ which we identify with the 
presheaf of sets it represents. 
Denote by $V^{\an}$ the analytification of $V$ 
with respect to the adic ring $\mathcal{O}(\mathcal{U})$.
Then, there is a morphism of presheaves of sets
\begin{equation}
\label{eq-lem:diag-eta-on-representable-1}
\diag_{\eta,\,*}(\mathcal{U},V)\to V^{\an}
\end{equation}
which induces an isomorphism after sheafification for the rig topology.

\item[(2)] Let $(\mathcal{U},V)$ be an object of 
$\RigSm/\mathfrak{D}^{\an}_{\mathcal{S},\,\eta}$ 
which we identify with the presheaf of sets it represents. 
Then, there is a morphism of presheaves of sets
\begin{equation}
\label{eq-lem:diag-eta-on-representable-2}
\diag^{\an}_{\eta,\,*}(\mathcal{U},V)\to V
\end{equation}
which induces an isomorphism after sheafification for the rig topology.

\end{enumerate}
\end{lemma}

\begin{proof}
We only prove the first part, which is slightly more interesting. 
Set $A=\mathcal{O}(\mathcal{U})$ and 
let $\mathcal{T}=\Spf(B)$ be a rig-smooth affine formal 
$\mathcal{S}$-scheme. A section of 
$\diag_{\eta,\,*}(\mathcal{U},V)$ on $\mathcal{T}$ 
is a pair $(f,g)$
consisting of a morphism of formal $\mathcal{S}$-schemes 
$f:\mathcal{T}\to \mathcal{U}$ and a morphism of 
schemes $g:\Spec(B)\smallsetminus \mathcal{T}_{\sigma} \to V$
over $\Spec(A)\smallsetminus \mathcal{U}_{\sigma}$.
This gives rise to a section of the
$(\Spec(B)\smallsetminus \mathcal{T}_{\sigma})$-scheme
$V\times_{\Spec(A)}\Spec(B)$ and, by analytification over
$\mathcal{T}$, to a morphism $\mathcal{T}^{\rig}\to V^{\an}
\times_{\mathcal{S}^{\rig}}\mathcal{T}^{\rig}$. This defines the 
morphism of presheaves \eqref{eq-lem:diag-eta-on-representable-1}. 
It remains to see that this morphism 
induces an equivalence on stalks for the rig topology.
To do so, we evaluate
\eqref{eq-lem:diag-eta-on-representable-1}
on a rig point $\mathfrak{t}=\Spf(R)$ over
$\mathcal{S}$, with $R$ an adic valuation ring with fraction
field $K$. We may replace $\mathcal{S}$ with 
$\mathfrak{t}$ and assume that $V$ is a smooth $K$-scheme. 
The question being local, we may assume that $V$ 
is compactifiable over $R$ and fix an open immersion 
$V\to \overline{V}$ into a proper $R$-scheme $\overline{V}$.
In this case, the evaluation of 
\eqref{eq-lem:diag-eta-on-representable-1}
on $\mathfrak{t}$ is the obvious map between
\begin{enumerate}

\item[(1)] the set of $K$-points $x:\Spec(K)\to V$;

\item[(2)] the set of $R$-points $\mathfrak{x}:\Spf(R) \to \widehat{\overline{V}}$
such that there exists an admissible blowup 
$\overline{V}{}'\to \overline{V}$ with the property that 
the lift $\mathfrak{x}':\Spf(R)\to 
\widehat{\overline{V}}{}'$ of $\mathfrak{x}$ factors through 
the complement of the special fiber of the Zariski closure of 
$\overline{V}{}_{\eta}'\smallsetminus V$ in 
$\overline{V}{}'$. (See Construction
\ref{cons:analytification-}.)

\end{enumerate}
To give a morphism of formal $R$-schemes 
$\mathfrak{x}:\Spf(R) \to \widehat{\overline{V}}$
is equivalent to giving a morphism of $R$-schemes
$\widetilde{x}:\Spec(R)\to \overline{V}$, and the condition in 
(2) corresponds to the condition that 
$\widetilde{x}$ sends $\Spec(K)$ to $V$.
Hence, the set described in (2) can be identified with
\begin{enumerate}

\item[(2$'$)] the set of $R$-points $\widetilde{x}:\Spec(R) \to \overline{V}$
sending $\Spec(K)$ to $V$.

\end{enumerate}
That the obvious map between (1) and (2$'$) is a bijection
is clear. (Note that the existence of this map follows from the 
valuative criterion of properness but, once the existence 
of this map is granted, it is clearly a bijection.)
\end{proof}

Recall that the weak equivalences 
of the stable $(\B^1,\tau)$-local
model structure are called the 
stable $(\B^1,\tau)$-local equivalences; see
Remark \ref{rmk:rigda-model-cat}. 
Similarly, we have the notions of stable 
$(\A^1,\tau)$-local equivalences and 
stable $(\A^1,\rig\text{-}\tau)$-local equivalences.
For later use, we record the following result.

\begin{lemma}
\label{lem:diag-exact}
$\empty$

\begin{enumerate}

\item[(1)]
The functor 
$${\rm diag}_{\eta,\,*}:
\Spect_T(\PSh(\Sm/\mathfrak{D}_{\mathcal{S},\,\eta};\Lambda))
\to \Spect_T(\PSh(\FRigSm_{\af}/\mathcal{S};\Lambda))$$
takes a stable $(\A^1,\tau)$-local equivalence to 
a stable $(\A^1,\rig\text{-}\tau)$-local equivalence. 

\item[(2)]
The functor 
$${\rm diag}^{\an}_{\eta,\,*}:
\Spect_T(\PSh(\RigSm/\mathfrak{D}^{\an}_{\mathcal{S},\,\eta};\Lambda))
\to \Spect_T(\PSh(\FRigSm_{\af}/\mathcal{S};\Lambda))$$
takes a stable $(\B^1,\tau)$-local equivalence to 
a stable $(\A^1,\rig\text{-}\tau)$-local equivalence.

\end{enumerate}
\end{lemma}

\begin{proof}
We only treat the first part; the second part is proven in the 
same way.
The functor ${\rm diag}_{\eta,\,*}$ commutes with colimits.
Thus, by \cite[Proposition 5.5.4.20]{lurie}, it is enough to show that 
${\rm diag}_{\eta,\,*}$ transforms the following types of morphisms
\begin{enumerate}

\item[(1)] $\colim_{[n]\in \mathbf{\Delta}} \,
\Lambda(\mathcal{U},V_n) \to 
\Lambda(\mathcal{U},V_{-1})$, 
where $V_{\bullet}$ is a $\tau$-hypercover,

\item[(2)] $\Lambda(\mathcal{U},V) \to 
\Lambda(\mathcal{U},\A^1_V)$,

\item[(3)] a morphism of $T$-spectra $F \to F'$ such that 
$F_n \to F'_n$ is an equivalence for $n$ large enough,

\end{enumerate}
into $(\A^1,\rig\text{-}\tau)$-local equivalences,
for (1) and (2), and into 
stable $(\A^1,\rig\text{-}\tau)$-local equivalences,
for (3). The case of (3) is obvious, so we only need to
discuss morphisms of type (1) and (2).

In (1) and (2) above, $\mathcal{U}$ is an affine formal scheme
which is rig-smooth over $\mathcal{S}$. We set 
$U=\Spec(\mathcal{O}(\mathcal{U}))$, $U_{\sigma}=\mathcal{U}_{\sigma}$
and $U_{\eta}=U\smallsetminus U_{\sigma}$. Then
$V$ and the $V_n$'s, for $n\geq -1$, are smooth $U_{\eta}$-schemes. 
By Lemma \ref{lem:diag-eta-on-representable}(1),
${\rm diag}_{\eta,\,*}$ takes morphisms of 
type (1) and (2) to morphisms which are
rig-locally equivalent to 
\begin{enumerate}

\item[(1$'$)] $\colim_{[n]\in \mathbf{\Delta}} \,
\Lambda(V^{\an}_n) \to 
\Lambda(V^{\an}_{-1})$,

\item[(2$'$)] $\Lambda(V^{\an}) \to 
\Lambda((\A^1_V)^{\an})$,

\end{enumerate}
where we use the notation introduced in 
aforementioned lemma. By Remark
\ref{rmk:another-site-for-rigsh}, it is enough to show 
that (1$'$) and (2$'$) are 
$(\B^1,\tau)$-equivalences in 
$\PSh(\RigSm/\mathcal{S}^{\rig};\Lambda)$
which is obvious. 
\end{proof}

We now state the main technical result needed for 
proving Proposition
\ref{prop-thm:morphism-of-specialisation-equival-weak}.
(Compare with \cite[Th\'eor\`eme 1.3.37]{ayoub-rig}.)

\begin{prop}
\label{prop:tech-for-full+}
Let $B$ be a quasi-excellent $(\Lambda,\tau)$-admissible 
scheme, $B_{\sigma}\subset B$ 
a closed subscheme locally of finite presentation, and 
$B_{\eta}\subset B$ its open complement.
If $\tau$ is the \'etale topology, 
assume that every prime number is invertible either 
in $\pi_0\Lambda$ or in $\mathcal{O}(B)$.
\begin{enumerate}

\item[(1)] Consider the commutative diagram of 
diagrams of schemes
$$\xymatrix{\mathfrak{D}_{\widehat{B},\,\eta} 
\ar[r]^-{\mathfrak{j}} \ar[d]^-{\mathfrak{u}_{\eta}} & 
\mathfrak{D}_{\widehat{B}} 
\ar[d]^-{\mathfrak{u}} &
\mathfrak{D}_{\widehat{B},\,\sigma}\ar[l]_-{\mathfrak{i}} 
\ar[d]^-{\mathfrak{u}_{\sigma}} \\
B_{\eta} \ar[r]^-j & B & B_{\sigma}.\! \ar[l]_-i}$$
Then, the composite functor 
\begin{equation}
\label{eq-prop:tech-for-full+-113}
\diag_{\sigma,\,*}\circ \mathfrak{i}^* \circ \mathfrak{j}_*
\circ \mathfrak{u}_{\eta}^*:
\SH^{\hyp}_{\tau}(B_{\eta};\Lambda)
\to 
\Spect_T(\PSh(\FRigSm_{\af}/\widehat{B};\Lambda))
\end{equation}
takes values in 
$\RigSH^{\hyp}_{\tau}(\widehat{B}{}^{\rig};\Lambda)$
considered as the full sub-$\infty$-category 
of the target of
\eqref{eq-prop:tech-for-full+-113} 
spanned by those objects which are 
stably $(\A^1,\rig\text{-}\tau)$-local.

\item[(2)] Consider the commutative diagram of 
diagrams of rigid analytic schemes
$$\xymatrix{\mathfrak{D}^{\an}_{\widehat{B},\,\eta} 
\ar[r]^-{\mathfrak{j}^{\an}} \ar[d]^-{\mathfrak{u}^{\an}_{\eta}} & 
\mathfrak{D}^{\an}_{\widehat{B}} 
\ar[d]^-{\mathfrak{u}^{\an}} &
\mathfrak{D}_{\widehat{B},\,\sigma}\ar[l]_-{\mathfrak{i}^{\an}} 
\ar[d]^-{\mathfrak{u}_{\sigma}} \\
\widehat{B}{}^{\rig} \ar[r]^-{j^{\an}} & \widehat{B} & B_{\sigma}.\! \ar[l]_-{i^{\an}}}$$
Then, the composite functor 
\begin{equation}
\label{eq-prop:tech-for-full+-115}
\diag^{\an}_{\sigma,\,*}\circ \mathfrak{i}^{\an,\,*}\circ 
\mathfrak{j}^{\an}_*
\circ \mathfrak{u}_{\eta}^{\an,\,*}:
\RigSH^{\hyp}_{\tau}(\widehat{B}{}^{\rig};\Lambda)
\to 
\Spect_T(\PSh(\FRigSm_{\af}/\widehat{B};\Lambda))
\end{equation}
takes values in 
$\RigSH^{\hyp}_{\tau}(\widehat{B}{}^{\rig};\Lambda)$
considered as the full sub-$\infty$-category 
of the target of
\eqref{eq-prop:tech-for-full+-115}
spanned by those objects which are 
stably $(\A^1,\rig\text{-}\tau)$-local.
Moreover, the induced endofunctor of 
$\RigSH^{\hyp}_{\tau}(\widehat{B}{}^{\rig};\Lambda)$
is equivalent to the identity functor.

\end{enumerate}
\end{prop}

\begin{proof}
We start with part (2) which is easier. 
Let $\diag^{\form}:\FRigSm_{\af}/\widehat{B}
\to \FSm/\mathfrak{D}^{\form}_{\widehat{B}}$
be the diagonal functor sending an affine formal
scheme $\mathcal{U}$ to the pair $(\mathcal{U},\mathcal{U})$,
and let $\diag^{\form}_*$ be constructed as in 
Remark \ref{rmk:diagonal-functors-}. 
We have an equivalence ${\rm diag}^{\form}_*\circ \sigma_*
\simeq {\rm diag}_{\sigma,\,*}$, where 
$\sigma_*$ is restriction along the functor 
$(-)_{\sigma}:
\FSm/\mathfrak{D}^{\form}_{\widehat{B}}
\to \Sm/\mathfrak{D}_{\widehat{B},\,\sigma}$.
By Lemma \ref{lem:chi-mathfrak-D-an} and 
Theorem \ref{thm:formal-mot-alg-mot}, the 
composite functor 
\eqref{eq-prop:tech-for-full+-115}
is equivalent to the composite functor
\begin{equation}
\label{eq-prop:tech-for-full+-1157}
\diag^{\form}_*\circ \chi_{\mathfrak{D}^{\form}_{\widehat{B}}}
\circ \mathfrak{u}_{\eta}^{\an,\,*}:
\RigSH^{\hyp}_{\tau}(\widehat{B}{}^{\rig};\Lambda)
\to 
\Spect_T(\PSh(\FRigSm_{\af}/\widehat{B};\Lambda)).
\end{equation}
Now, $\chi_{\mathfrak{D}^{\form}_{\widehat{B}}}$ is 
restriction along the functor
$(-)^{\rig}:\FSm/\mathfrak{D}^{\form}_{\widehat{B}}
\to \RigSm/\mathfrak{D}^{\an}_{\widehat{B}{}^{\rig}}$
and $\mathfrak{u}_{\eta}^*$ is restriction along 
along the functor
$\RigSm/\mathfrak{D}^{\an}_{\widehat{B}{}^{\rig}}
\to \RigSm/\widehat{B}{}^{\rig}$
sending a pair $(\mathcal{U},\mathcal{V})$ 
to $\mathcal{V}^{\rig}$. It follows that the composite functor
\eqref{eq-prop:tech-for-full+-1157}
is restriction along the functor 
$(-)^{\rig}:\FRigSm_{\af}/\widehat{B}
\to \RigSm/\widehat{B}{}^{\rig}$.
The claim now follows from Remark
\ref{rmk:another-site-for-rigsh}.

We now concentrate on part (1). 
We fix an object $M\in \SH^{\hyp}_{\tau}(B_{\eta};\Lambda)$.
Our goal is to show that $\diag_{\sigma,\,*}
\mathfrak{i}^*\mathfrak{j}_*\mathfrak{u}_{\eta}^*M$
belongs to the full sub-$\infty$-category 
\begin{equation}
\label{eq-prop:tech-for-full+-1156}
\RigSH^{\hyp}_{\tau}(\widehat{B}{}^{\rig};\Lambda)
\subset \Spect_T(\PSh(\FRigSm_{\af}/\widehat{B};\Lambda)).
\end{equation}
The proof of this is similar to the proof of Proposition 
\ref{prop:tech-for-full} and, instead of repeating 
large portions of that proof we will refer to it when possible.
It follows from Propositions
\ref{prop:compact-shv-rigsm} and
\ref{prop:compact-shv-forsm} that 
the sub-$\infty$-category 
\eqref{eq-prop:tech-for-full+-1156}
is closed under colimits and that the functors
$\diag_{\sigma,\,*}$, $\mathfrak{i}^*$, $\mathfrak{j}_*$
and $\mathfrak{u}_{\eta}^*$ are colimit-preserving. 
Thus, we may assume that $M$ is compact. We split the proof
into several steps.

\paragraph*{Step 1}
\noindent
Arguing as in the second part of the proof of Lemma
\ref{lem:alternative-iii-vs-iv-3i}, 
we may assume one of the following alternatives:
\begin{enumerate}

\item[(1)] $\tau$ is the Nisnevich topology;

\item[(2)] $\pi_0\Lambda$ is a $\Q$-algebra;

\item[(3)] $\tau$ is the \'etale topology and
$M$ is $\ell$-nilpotent for a prime $\ell$ invertible on $B$. 

\end{enumerate}
Moreover, we claim that under the alternative (3), 
we may assume that $\Lambda$ is eventually coconnective. 
To prove this, let $M_0\in \Shv^{\hyp}_{\et}(\Et/B_{\eta};\Lambda)_{\ellnil}$ be the object corresponding to $M$ 
by the equivalence 
$$\Shv^{\hyp}_{\et}(\Et/B_{\eta};\Lambda)_{\ellnil}
\simeq 
\SH^{\hyp}_{\et}(B_{\eta};\Lambda)_{\ellnil}$$
provided by Theorem 
\ref{thm:rigrig-algebraic}.
Then, as a $T$-spectrum, $M$ is given at level $m$ 
by $\iota_{B_{\eta}}^*M_0(m)[m]$,
where $\iota_{B_{\eta}}^*$ is as in Notation
\ref{nota:lambda-ell-etale-sheaf}.
(See \cite[Corollary 4.9]{ayoub-etale} in the case where 
$\Lambda$ is an Eilenberg--Mac Lane spectrum; the general case can 
be treated similarly.) Similarly, as a $T$-spectrum,
$\mathfrak{i}^*\mathfrak{j}_*\mathfrak{u}_{\eta}^*M$
is given at level $m$ by 
$\iota_{\mathfrak{D}_{\widehat{B},\,\sigma}}^*
\mathfrak{i}^*\mathfrak{j}_*\mathfrak{u}_{\eta}^*M_0(m)[m]$.
Using this and 
Lemma \ref{lem:pi-0-Lambda-coh-dim}, 
one deduces an equivalence 
$$\diag_{\sigma,\,*}
\mathfrak{i}^*\mathfrak{j}_*\mathfrak{u}_{\eta}^*M
\simeq \lim_r 
\diag_{\sigma,\,*}
\mathfrak{i}^*\mathfrak{j}_*\mathfrak{u}_{\eta}^*(M\otimes_{\Lambda}
\tau_{\leq r}\Lambda).$$
Since the sub-$\infty$-category 
\eqref{eq-prop:tech-for-full+-1156} is stable under limits,
we deduce that it is enough to prove the result for 
$M\otimes_{\Lambda} \tau_{\leq r}\Lambda$ .
This proves our claim.

In conclusion, when $\tau$ is the \'etale 
topology, we may assume that $\Lambda$ is 
eventually coconnective. (Indeed, if $\pi_0\Lambda$ is a 
$\Q$-algebra, there is a morphism $\Q\to \Lambda$ and
we may replace $\Lambda$ by $\Q$.)

\paragraph*{Step 2}
\noindent
From now on, we set $E=\diag_{\sigma}^*
\mathfrak{i}^*\mathfrak{j}_*\mathfrak{u}_{\eta}^*M$
and, for $m\in \N$, we denote by 
$E_m$ the $m$-th level of the $T$-spectrum $E$. 
In this step, we show that $E$ admits levelwise hyperdescent 
for the rig-Nisnevich topology. 
Arguing as in the beginning of the proof of Proposition
\ref{prop:tech-for-full}, we need to show that 
$E_m$ has descent for every rig-Nisnevich hypercover 
$\mathcal{U}_{\bullet}$ in $\FRigSm_{\af}/\widehat{B}$
admitting a morphism of augmented simplicial formal schemes
$\widetilde{\mathcal{U}}_{\bullet}
\to \mathcal{U}_{\bullet}$ such that:
\begin{itemize}

\item $\widetilde{\mathcal{U}}_{\bullet}$ is a Nisnevich 
hypercover;

\item $\widetilde{\mathcal{U}}_{-1} \to \mathcal{U}_{-1}$ is an 
admissible blowup;

\item $\widetilde{\mathcal{U}}_n \to \mathcal{U}_n$
is an isomorphism for $n\geq 0$.

\end{itemize}
In particular, we see that $\widetilde{\mathcal{U}}_n$ is 
affine except possibly when $n=-1$.
For $n\geq -1$, we set $U_n=\Spec(\mathcal{O}(\mathcal{U}_n))$
and, for $n\geq 0$, we set $\widetilde{U}_n=U_n$.
Since $\widetilde{\mathcal{U}}_{-1}\to \mathcal{U}_{-1}$ is an 
admissible blowup, it is the formal completion of a unique
blowup $e:\widetilde{U}_{-1}\to U_{-1}$ with center supported on 
$\mathcal{U}_{-1,\,\sigma}\subset U_{-1}$.
For $n\geq -1$, we set $U_{n,\,\sigma}=\mathcal{U}_{n,\,\sigma}$, 
$\widetilde{U}_{n,\,\sigma}=\widetilde{\mathcal{U}}_{n,\,\sigma}$,
$U_{n,\,\eta}=U_n\smallsetminus U_{n,\,\sigma}$ and 
$\widetilde{U}_{n,\,\eta}=\widetilde{U}_n\smallsetminus U_{n,\,\sigma}$.
We denote by $u_n:U_n\to B$ and $\widetilde{u}_n:\widetilde{U}_n\to B$
the obvious morphisms.

Since $M$ can be shifted and twisted, 
it suffices to prove that the map
$$\Map(\Lambda(\mathcal{U}_{-1}),E_0)\to 
\lim_{[n]\in\mathbf{\Delta}}
\Map(\Lambda(\mathcal{U}_n),E_0)$$
is an equivalence, where the mapping spaces are taken in 
$\PSh(\FRigSm_{\af}/\widehat{B};\Lambda)$.
Looking at the definition of $E_0$, we see that this map
is equivalent to 
\begin{equation}
\label{eq-prop:tech-for-full+-11568}
\begin{array}{rcl}
\Map_{\SH^{\hyp}_{\tau}(U_{-1,\,\sigma};\,\Lambda)}
(\Lambda,\chi_{U_{-1}} u_{-1,\,\eta}^*M)
& \to & {\displaystyle \lim_{[n]\in \mathbf{\Delta}}
\Map_{\SH^{\hyp}_{\tau}(U_{n,\,\sigma};\,\Lambda)}
(\Lambda,\chi_{U_n} u_{n,\,\eta}^*M)}\\
& = & {\displaystyle \lim_{[n]\in \mathbf{\Delta}}
\Map_{\SH^{\hyp}_{\tau}(\widetilde{U}_{n,\,\sigma};\,\Lambda)}
(\Lambda,\chi_{\widetilde{U}_n} \widetilde{u}_{n,\,\eta}^*M).}
\end{array}
\end{equation}
For $n\geq 0$, we let $v_n:\widetilde{U}_n \to \widetilde{U}_{-1}$
be the obvious morphism. Since $B$ is quasi-excellent, 
the $v_n$'s are regular morphisms. By Lemma 
\ref{lem:regular-base-change} below, the morphism
$$\chi_{\widetilde{U}_n} \widetilde{u}_{n,\,\eta}^*M
\to v_{n,\,\sigma}^*
\chi_{\widetilde{U}_{-1}} \widetilde{u}_{-1,\,\eta}^*M$$
is an equivalence. Therefore, the left-hand side in
\eqref{eq-prop:tech-for-full+-11568} 
is equivalent to 
$$\lim_{[n]\in \mathbf{\Delta}}
\Map_{\SH^{\hyp}_{\tau}(\widetilde{U}_{n,\,\sigma};\,\Lambda)}
(\Lambda,v_{n,\,\sigma}^*
\chi_{\widetilde{U}_{-1}} \widetilde{u}_{-1,\,\eta}^*M).$$
Since $\widetilde{U}_{\bullet,\,\sigma}$ is a Nisnevich hypercover,
the latter is equivalent to 
$\Map_{\SH^{\hyp}_{\tau}(\widetilde{U}_{-1,\,\sigma};\,\Lambda)}
(\Lambda,\chi_{\widetilde{U}_{-1}} u_{-1,\,\eta}^*M)$.
Thus, we are left to show that the morphism
$$\chi_{U_{-1}}u^*_{-1,\,\eta}M\to 
e_{\sigma,\,*}\chi_{\widetilde{U}_{-1}}
\widetilde{u}^*_{-1,\,\eta}M$$
is an equivalence. This follows from the projective base change 
theorem and the fact that $e_{\eta}$ is an isomorphism.

\paragraph*{Step 3}
\noindent
In this step and the next one, 
we assume that $\tau$ is the \'etale topology and we 
prove that $E$ admits levelwise
hyperdescent for the rig-\'etale topology. 
By the second step, we already know that $E$ admits levelwise 
hyperdescent for the rig-Nisnevich topology. Thus,
arguing as in the beginning of the proof of Proposition
\ref{prop:tech-for-full}, it remains to 
show that $E$ has levelwise descent 
for the topology $\rig\fet$.

In this step, we deal with the case
where $\pi_0\Lambda$ is a $\Q$-algebra. 
As explained in the third part of the proof of Proposition
\ref{prop:tech-for-full}, we only need to show that 
$E$ has levelwise descent for a $\rig\fet$-hypercover
of the form 
\begin{equation}
\label{eq-prop:tech-for-full+-115681}
\xymatrix@C=1.3pc{\cdots \mathcal{V}_0\times G\times G 
\ar[r] \ar@<-.4pc>[r] \ar@<.4pc>[r] &
\mathcal{V}_0\times G \ar@<-.2pc>[r] \ar@<.2pc>[r]& \mathcal{V}_0
\ar[r] & \mathcal{V}_{-1}.}
\end{equation}
where $\mathcal{V}_{-1}$ is an affine rig-smooth 
formal $\widehat{B}$-scheme and
$\mathcal{V}_0\to \mathcal{V}_{-1}$ is a finite rig-\'etale
covering admitting an action of a finite group $G$
which is simply transitive on the geometric fibers of 
$\mathcal{V}_0^{\rig} \to \mathcal{V}_{-1}^{\rig}$. 
For $n\in \{-1,0\}$, we set 
$V_n=\Spec(\mathcal{O}(\mathcal{V}_n))$, 
$V_{n,\,\sigma}=\mathcal{V}_{n,\,\sigma}$ and $V_{n,\,\eta}=V_n\smallsetminus V_{n,\,\sigma}$. 
We also denote by $v_{-1}:V_{-1} \to B$, $v_0:V_0\to B$ 
and $e:V_0\to V_{-1}$ the obvious morphisms.
For later use, we note that $e_{\eta}:V_{0,\,\eta}\to V_{-1,\,\eta}$ 
is a finite \'etale cover
admitting an action of $G$
which is simply transitive on geometric fibers.

Since $M$ can be shifted and twisted, 
it suffices to prove that the map
$$\Map(\Lambda(\mathcal{V}_{-1}),E_0)\to 
\Map(\Lambda(\mathcal{V}_0),E_0)^G$$
is an equivalence, where the mapping spaces are taken in 
$\PSh(\FRigSm_{\af}/\widehat{B};\Lambda)$.
Looking at the definition of $E_0$, we see that this map
is equivalent to
$$\Map_{\SH^{\hyp}_{\et}(V_{-1,\,\sigma};\,\Lambda)}(\Lambda,
\chi_{V_{-1}}v_{-1,\,\eta}^*M)
\to \Map_{\SH^{\hyp}_{\et}(V_{0,\,\sigma};\,\Lambda)}(\Lambda,
\chi_{V_0}v_{0,\,\eta}^*M)^G.$$
Thus, it is enough to show that 
$$\chi_{V_{-1}}v_{-1,\,\eta}^*M 
\to  (e_{\sigma,\,*}
\chi_{V_0}v_{0,\,\eta}^*M)^G 
\simeq  (\chi_{V_{-1}}e_{\eta,\,*}v_{-1,\,\eta}^*M)^G$$
is an equivalence. (The equivalence above follows from the proper
base change theorem and the fact that $e$ is finite.)
Taking the ``$G$-invariant subobject'' in a $\Q$-linear 
$\infty$-category
is equivalent to taking the image of the projector 
$|G|^{-1}\sum_{g\in G} g$, and hence it commutes with 
the functor $\chi_{V_{-1}}$. 
Thus, it is enough to show that $v^*_{-1,\,\eta}M_0
\to (e_{\eta,\,*}e_{\eta}^*v^*_{-1,\,\eta}M_0)^G$
is an equivalence, which follows from 
\'etale descent in $\SH^{\hyp}_{\et}(V_{-1,\,\eta};\Lambda)$.

\paragraph*{Step 4}
\noindent
Here we complete the proof that $E$ 
admits levelwise hyperdescent for the rig-\'etale topology. 
By the first and the third steps, we may assume that 
$M$ is $\ell$-nilpotent and that $\Lambda$ is 
eventually coconnective. Let $M_0\in \Shv^{\hyp}_{\et}(\Et/B_{\eta};\Lambda)_{\ellnil}$ be the object corresponding to $M$ 
by the equivalence 
$$\Shv^{\hyp}_{\et}(\Et/B_{\eta};\Lambda)_{\ellnil}
\simeq 
\SH^{\hyp}_{\et}(B_{\eta};\Lambda)_{\ellnil}$$
provided by Theorem 
\ref{thm:rigrig-algebraic}.
As in the third step, it suffices to show descent for the 
$\rig\fet$-hypercover
\eqref{eq-prop:tech-for-full+-115681}
and it is enough to prove that 
$$\chi_{V_{-1}}v_{-1,\,\eta}^*M
\to (\chi_{V_{-1}}e_{\eta,\,*}v_{-1,\,\eta}^*M)^G$$
is an equivalence. Using Theorem 
\ref{thm:rigrig-algebraic}, we may as well prove that 
$$\chi_{V_{-1}}v_{-1,\,\eta}^*M_0
\to (\chi_{V_{-1}}e_{\eta,\,*}v_{-1,\,\eta}^*M_0)^G$$
is an equivalence. Since $\Lambda$ is eventually coconnective
and $M_0$ is compact, we deduce that the \'etale sheaf 
$M_0$ is also eventually coconnective.  
Taking the ``$G$-invariant subobject''
commutes with direct images and, if we restrict to 
eventually coconnective \'etale sheaves, it also commute with 
inverse images. (The latter assertion can be proven 
using an explicit model for the $G$-invariant
functor; see the fourth part of the proof of 
Proposition \ref{prop:tech-for-full}
for a similar argument.) Thus, as in the previous step, 
it is enough to show that $v^*_{-1,\,\eta}M_0
\to (e_{\eta,\,*}e_{\eta}^*v^*_{-1,\,\eta}M_0)^G$
is an equivalence, which follows from 
\'etale descent in $\Shv^{\hyp}_{\et}(\Et/V_{-1,\,\eta};\Lambda)$.

\paragraph*{Step 5}
\noindent
In this last step, we check that $E$ is levelwise $\A^1$-invariant 
and an $\Omega$-spectrum. Since $M$ can be shifted, it is enough 
to show that the maps 
\begin{equation}
\label{eq-prop:tech-for-full+-1156810}
\begin{array}{l}
\Map(\Lambda(\mathcal{U}),E_m)
\to \Map(\Lambda(\A^1_{\mathcal{U}}),E_m),\\
\vspace{-.3cm} \\
\Map(\Lambda(\mathcal{U}),E_m)
\to \fiber\{\Map(\Lambda(\A^1_{\mathcal{U}}\smallsetminus 
0_{\mathcal{U}}),E_{m+1}) \xrightarrow{1^*} 
\Map(\Lambda(\mathcal{U}),E_{m+1})\}
\end{array}
\end{equation}
are equivalences for every $\mathcal{U}\in \FRigSm_{\af}/\widehat{B}$.

Set $U=\Spec(\mathcal{O}(\mathcal{U}))$, 
$U_{\sigma}=\mathcal{U}_{\sigma}$ and 
$U_{\eta}=U\smallsetminus U_{\sigma}$. 
Let $\mathcal{V}$ be an affine smooth formal $\mathcal{U}$-scheme,
and set $V=\Spec(\mathcal{O}(\mathcal{V}))$, 
$V_{\sigma}=\mathcal{V}_{\sigma}$ and 
$V_{\eta}=V\smallsetminus V_{\sigma}$.
Denote by $u:U\to B$ and $g:V\to U$ the obvious morphisms. 
Then we have equivalences
$$\begin{array}{rcl}
\Map(\Lambda(\mathcal{U}),E_m) & \simeq & 
\Map_{\SH^{\hyp}_{\tau}(U_{\sigma};\,\Lambda)}(\Lambda(-m),
\chi_U u_{\eta}^*M),\\
&  &\\
\Map(\Lambda(\mathcal{V}),E_m) & \simeq &  
\Map_{\SH^{\hyp}_{\tau}(V_{\sigma};\,\Lambda)}(\Lambda(-m),
\chi_V g_{\eta}^*u_{\eta}^*M),\\
& \overset{(1)}{\simeq} & 
\Map_{\SH^{\hyp}_{\tau}(V_{\sigma};\,\Lambda)}(\Lambda(-m),
g_{\sigma}^*\chi_Uu_{\eta}^*M),\\
& \overset{(2)}{\simeq} & 
\Map_{\SH^{\hyp}_{\tau}(U_{\sigma};\,\Lambda)}(\Lambda(-m),
g_{\sigma,\,*}g_{\sigma}^*\chi_Uu_{\eta}^*M).
\end{array}$$
The equivalence (1) follows from Lemma
\ref{lem:regular-base-change} 
below and the fact that $g$ is regular.
The equivalence (2) follows by adjunction. 
Letting $p:\A^1_{U_{\sigma}}\to U_{\sigma}$
and $q:\A^1_{U_{\sigma}}\smallsetminus 0_{U_{\sigma}}\to 
U_{\sigma}$ be the obvious projections, we 
deduce that the maps \eqref{eq-prop:tech-for-full+-1156810}
are equivalent to the following ones:
$$\begin{array}{l}
\Map_{\SH^{\hyp}_{\tau}(U_{\sigma};\,\Lambda)}
(\Lambda(-m),\chi_Uu_{\eta}^*M)
\to \Map_{\SH^{\hyp}_{\tau}(U_{\sigma};\,\Lambda)}
(\Lambda(-m),p_*p^*\chi_Uu_{\eta}^*M),\\
\vspace{-.3cm} \\
\Map_{\SH^{\hyp}_{\tau}(U_{\sigma};\,\Lambda)}
(\Lambda(-m),\chi_Uu_{\eta}^*M)
\to \Map_{\SH^{\hyp}_{\tau}(U_{\sigma};\,\Lambda)}
(\Lambda(-m-1),\fiber\{1^*:q_*q^*\chi_Uu_{\eta}^*M
\to \chi_Uu_{\eta}^*M\})
\end{array}$$
which are clearly equivalences as needed.
\end{proof}

The following lemma was used in the proof of Proposition
\ref{prop:tech-for-full+}. We prove it in a greater generality 
than needed because of its potential usefulness. 

\begin{lemma}[Regular base change]
\label{lem:regular-base-change}
\ncn{regular base change}
Consider a Cartesian square of schemes
$$\xymatrix{Y' \ar[r]^-{g'} \ar[d]^-{f'} & Y\ar[d]^-f\\
X' \ar[r]^-g & X}$$
with $X$ locally noetherian, $g$ regular,
and $f$ quasi-compact and quasi-separated.
Assume one of the following alternatives:
\begin{enumerate}

\item[(1)] we work in the non-hypercomplete case and, when $\tau$
is the \'etale topology, we assume furthermore 
that $\Lambda$ is eventually coconnective;

\item[(2)] we work in the hypercomplete case, and the schemes 
$X$, $X'$, $Y$ and $Y'$ are $(\Lambda,\tau)$-admissible.

\end{enumerate}
Then, the natural transformation 
$g^*\circ f_*\to f'_*\circ g'^*$, between functors from 
$\SH^{(\eff,\,\hyp)}_{\tau}(Y;\Lambda)$ to
$\SH^{(\eff,\,\hyp)}_{\tau}(X';\Lambda)$,
is an equivalence.
\end{lemma}

\begin{proof}
This is a generalisation of \cite[Corollary 1.A.4]{ayoub-rig}
and, as in loc.~cit., 
its proof consists in reducing to the smooth base change 
theorem using Popescu's theorem on regular algebras and 
Proposition \ref{prop:cont-algebraic-29}. However, here we need
an extra argument to reduce to the case where
$\Lambda$ is eventually coconnective so that 
Proposition \ref{prop:cont-algebraic-29}
applies.
The problem being local on $X$, $X'$ and $Y$, we may assume that 
$X$, $X'$, $Y$ and $Y'$ are affine. (This uses the hypothesis that
$f$ is quasi-compact and quasi-separated.)
By Proposition
\ref{prop:compact-shv-forsm}, the $\infty$-category 
$\SH^{(\eff,\,\hyp)}_{\tau}(X;\Lambda)$
is compactly generated, and similarly for 
$X'$, $Y$ and $Y'$. By the same proposition, 
the functors $f_*$ and $f'_*$ are colimit-preserving, 
and thus belong to $\Prl$. (The same is obviously true for 
$g^*$ and $g'^*$.)

We first prove the lemma under the alternative (1).
By \cite[Theorem 1.8]{popescu-86}, 
the $X$-scheme $X'$ is a limit of a cofiltered inverse
system $(X'_{\alpha})_{\alpha}$ of smooth affine $X$-schemes. 
For each $\alpha$, consider a Cartesian square
$$\xymatrix{Y'_{\alpha} \ar[r]^-{g'_{\alpha}} 
\ar[d]^-{f'_{\alpha}} & Y\ar[d]^-f \\
X'_{\alpha} \ar[r]^-{g_{\alpha}} & X.}$$
By the smooth base change theorem, we have commutative
squares in $\Prl$
$$\xymatrix{\SH^{(\eff)}_{\tau}(Y'_{\alpha};\Lambda) 
\ar[d]^-{f'_{\alpha,\,*}} & 
\SH^{(\eff)}_{\tau}(Y;\Lambda) \ar[d]^-{f_*} 
\ar[l]_-{g'^*_{\alpha}}\\
\SH^{(\eff)}_{\tau}(X'_{\alpha};\Lambda) & \SH^{(\eff)}_{\tau}(X;\Lambda)
\ar[l]_-{g^*_{\alpha}} }$$
Taking the colimit in $\Prl$ of these squares yields 
a commutative square expressing that 
$g^*\circ f_*$ is equivalent to $f'_*\circ g'^*$ as needed.
(This is actually not obvious; one needs to argue as in the 
proof of Theorem \ref{thm:general-base-change-thm}. 
We leave the details to the reader.)

Next, we prove the lemma under the alternative (2). 
Using Proposition
\ref{prop:automatic-hypercomp-motives-algebraic}, we may conclude 
using the lemma under the alternative (1) if $\tau$ is the 
Nisnevich topology or if
$\Lambda$ is eventually coconnective and, more generally,
if $\Lambda$ is an algebra over an eventually coconnective 
commutative ring spectrum. In particular, we may assume that 
$\tau$ is the \'etale topology, and the 
result holds if $\pi_0\Lambda$ is a $\Q$-algebra.
Arguing as in the second part of the proof of Lemma 
\ref{lem:alternative-iii-vs-iv-3i}, it remains to prove that 
$g^*f_*M\to f'_*g'^*M$ is an equivalence when
$M\in \SH^{(\eff),\,\hyp}_{\et}(Y;\Lambda)_{\ellnil}$, 
for some prime $\ell$ invertible on $X$. 
Moreover, we may assume that $M$ is compact. 
By Theorem \ref{thm:rigrig-algebraic}, it is 
enough to show that 
$g^*f_*M_0\to f'_*g'^*M_0$ is an equivalence
for $M_0\in \Shv^{\hyp}_{\et}(\Et/Y;\Lambda)_{\ellnil}$.
Using 
Lemma \ref{lem:pi-0-Lambda-coh-dim}, one deduces equivalences
$$g^*f_*M_0\simeq \lim_r 
g^*f_*(M_0\otimes_{\Lambda}\tau_{\leq r}\Lambda) 
\qquad \text{and} \qquad 
f'_*g'^*M_0\simeq \lim_r
f'_*g'^*(M_0\otimes_{\Lambda}\tau_{\leq r}\Lambda).$$
Thus, we may replace $M$ and $\Lambda$ with 
$M\otimes_{\Lambda}\tau_{\leq r}\Lambda$
and $\tau_{\leq r}\Lambda$. We are then automatically working 
under the alternative (1), and the result follows.
\end{proof}

\begin{proof}[Proof of Proposition
\ref{prop-thm:morphism-of-specialisation-equival-weak}]
We have a commutative square of natural transformations 
\begin{equation}
\label{eq-prop-thm:morphism-of-specialisation-equival-weak-1}
\begin{split}
\xymatrix{\diag_{\eta,\,*}\circ \mathfrak{u}_{\eta}^*
\ar[r]^-{\alpha} \ar[d]^-{\beta'} &
\diag_{\sigma,\,*}\circ \mathfrak{i}^*\circ \mathfrak{j}_*\circ 
\mathfrak{u}_{\eta}^*
\ar[d]^-{\beta} \\ 
\diag^{\an}_{\eta,\,*}\circ \mathfrak{u}^{\an,\,*}_{\eta}
\circ \An^*_{B_{\eta}}
\ar[r]^-{\alpha'} &
\diag_{\sigma,\,*}\circ \mathfrak{i}^{\an,\,*}\circ \mathfrak{j}^{\an}_*
\circ \mathfrak{u}^{\an,\,*}_{\eta}\circ \An_{B_{\eta}}^*.}
\end{split}
\end{equation}
The natural transformation $\alpha$ is obtained from 
$\mathfrak{j}_*\to \mathfrak{i}_*\circ \mathfrak{i}^*\circ \mathfrak{j}_*$
by applying $\diag_*$, and similarly for the natural 
transformation $\alpha'$.
The natural transformation $\beta$ is deduced from 
\eqref{eq-dfn:arrow-site-13708+}
(with $\mathcal{S}=\widehat{B}$).
Finally, the natural transformation 
$\beta'$ is deduced from the second natural transformation in 
\eqref{eq-rmk:diagonal-functors-} (with $\mathcal{S}=\widehat{B}$)
and the equivalence $\mathfrak{u}_{\eta}^{\an,\,*}\circ 
\An^*_{B_{\eta}}\simeq \An^*_{\mathfrak{D}_{\widehat{B},\,\eta}}
\circ \mathfrak{u}^*_{\eta}$.

We claim that the natural transformation 
$\beta\circ \alpha$ is given by stable $(\A^1,\rig\text{-}\tau)$-local
equivalences. We will prove this by showing that 
$\alpha'$ is an equivalence and that 
$\beta'$ is given by 
stable $(\A^1,\rig\text{-}\tau)$-local
equivalences. We then use this to finish the proof of 
the proposition. We split the remainder of the proof into 
three steps accordingly.

\paragraph*{Step 1}
\noindent 
Here we prove that 
$\alpha'$ is an equivalence. In fact, even the natural transformation 
$$\diag^{\an}_{\eta,\,*}
\to \diag_{\sigma,\,*} \circ \mathfrak{i}^{\an,\,*}\circ 
\mathfrak{j}_*^{\an}=\diag_{\sigma,\,*}\circ 
\chi_{\mathfrak{D}^{\an}_{\widehat{B}}}$$
is an equivalence. Indeed, by Lemma \ref{lem:chi-mathfrak-D-an} and 
Theorem \ref{thm:formal-mot-alg-mot}, we have an equivalence
$$\diag_{\sigma,\,*}\circ 
\chi_{\mathfrak{D}^{\an}_{\widehat{B}}}\simeq 
\diag^{\form}_*\circ \chi_{\mathfrak{D}^{\form}_{\widehat{B}}}.$$
(See the beginning of the proof of Proposition
\ref{prop:tech-for-full+}.)
Thus, we need to show that the natural transformation 
$$\diag^{\an}_{\eta,\,*}
\to \diag^{\form}_*\circ \chi_{\mathfrak{D}^{\form}_{\widehat{B}}}$$
is an equivalence. This follows from 
the equality $\diag^{\an}_{\eta}=(-)^{\rig}\circ \diag^{\form}$
and the fact that 
$\chi_{\mathfrak{D}^{\form}_{\widehat{B}}}$
is restriction along the functor 
$(-)^{\rig}:\FSm/\mathfrak{D}_{\widehat{B}}^{\form}
\to 
\FSm/\mathfrak{D}_{\widehat{B},\,\eta}^{\rm an}$.

\paragraph*{Step 2}
\noindent 
Here we prove that 
$\beta'$ is given by stable $(\A^1,\rig\text{-}\tau)$-local
equivalences. Since all the functors composing the source and
the target of $\beta'$ are colimit-preserving and since 
stable $(\A^1,\rig\text{-}\tau)$-local
equivalences are preserved by colimits, it is enough to show that
$$\beta'_M:\diag_{\eta,\,*}\mathfrak{u}^*_{\eta}M
\to \diag^{\an}_{\eta,\,*}
\mathfrak{u}^{\an,\,*}_{\eta}\An^*_{B_{\eta}}M$$
is a stable $(\A^1,\rig\text{-}\tau)$-local
equivalence when $M$ is of the form 
$\Lder_{\A^1,\,\tau,\,\st}\Sus_T^m\Lambda(X)$ for $n\in\N$ and $X\in\Sm/B_{\eta}$.
(Here, $\Lder_{\A^1,\,\tau,\,\st}$ is the stable 
$(\A^1,\tau)$-localisation functor and 
$\Sus_T^m$ is the left adjoint sending a $T$-spectrum
to its $m$-th level.) We have an equivalence
$$\mathfrak{u}^*_{\eta}M\simeq 
\Lder_{\A^1,\,\tau,\,\st}\mathfrak{u}^*_{\eta}\Sus_T^m\Lambda(X)$$
where, on the right-hand side, 
$\mathfrak{u}^*_{\eta}:\Spect_T(\PSh(\Sm/B_{\eta};\Lambda))
\to \Spect_T(\PSh(\Sm/\mathfrak{D}_{\widehat{B},\,\eta};\Lambda))$ 
is the inverse image functor on $T$-spectra of
presheaves of $\Lambda$-modules.
Using Lemma \ref{lem:diag-exact}(1), we deduce 
a stable $(\A^1,\rig\text{-}\tau)$-local equivalence
$$\diag_{\eta,\,*}\mathfrak{u}^*_{\eta}
\Sus^m_T\Lambda(X) \to 
\diag_{\eta,\,*}\mathfrak{u}^*_{\eta}M.$$
Similarly, we have $\An^*_{B_{\eta}}M\simeq 
\Lder_{\B^1,\,\tau,\,\st}\Sus^m_T\Lambda(X^{\an})$.
Arguing as before and using 
Lemma \ref{lem:diag-exact}(2), we deduce a stable 
$(\A^1,\rig\text{-}\tau)$-local equivalence
$$\diag^{\an}_{\eta,\,*}\mathfrak{u}^{\an,\,*}_{\eta}
\Sus^m_T\Lambda(X^{\an})
\to \diag^{\an}_{\eta,\,*}
\mathfrak{u}^{\an,\,*}_{\eta}\An^*_{B_{\eta}}M.$$
The result follows now by remarking that the obvious morphism 
$$\diag_{\eta,\,*}\mathfrak{u}^*_{\eta}
\Sus^m_T\Lambda(X)
\to 
\diag^{\an}_{\eta,\,*}\mathfrak{u}^{\an,\,*}_{\eta}
\Sus^m_T\Lambda(X^{\an})$$
is an isomorphism.

\paragraph*{Step 3}
\noindent
We are now ready to finish the proof of the proposition. 
By Proposition \ref{prop:tech-for-full+}(1), 
the functor $\diag_{\sigma,\,*}\circ 
\mathfrak{i}^*\circ \mathfrak{j}_*\circ 
\mathfrak{u}_{\eta}^*$
takes values in the $\infty$-subcategory spanned 
by stably $(\A^1,\rig\text{-}\tau)$-local objects. 
Therefore, $\alpha$ factors through the functor 
$\Lder_{\A^1,\,\rig\text{-}\tau,\,\st}\circ
\diag_{\eta,\,*}\circ \mathfrak{u}_{\eta}^*$
and the composition of 
$$\Lder_{\A^1,\,\rig\text{-}\tau,\,\st}\circ
\diag_{\eta,\,*}\circ \mathfrak{u}_{\eta}^*
\xrightarrow{\widetilde{\alpha}}
\diag_{\sigma,\,*}\circ \mathfrak{i}^*\circ \mathfrak{j}_*\circ 
\mathfrak{u}_{\eta}^*
\xrightarrow{\beta}
\diag_{\sigma,\,*}\circ \mathfrak{i}^{\an,\,*}\circ \mathfrak{j}^{\an}_*
\circ \mathfrak{u}^{\an,\,*}_{\eta}\circ \An_{B_{\eta}}^*$$
is given by stable $(\A^1,\rig\text{-}\tau)$-local 
equivalences (by the first and second steps). 
Since the source and the target of this 
composition take values in the $\infty$-subcategory spanned 
by stably $(\A^1,\rig\text{-}\tau)$-local objects
(by Proposition \ref{prop:tech-for-full+}(2) for the target), 
this composition is in fact a natural equivalence. 
Thus, we have shown that $\beta$ admits a section. 
Applying the restriction functor 
$$\mathfrak{r}_*:
\Spect_T(\PSh(\FRigSm_{\af}/\widehat{B};\Lambda))
\to \Spect_T(\PSh(\FSm_{\af}/\widehat{B};\Lambda))$$
to $\beta$, we deduce a natural transformation 
$$\mathfrak{r}_*(\beta):
\mathfrak{r}_*\circ 
\diag_{\sigma,\,*}\circ 
\mathfrak{i}^*\circ \mathfrak{j}_*\circ 
\mathfrak{u}_{\eta}^*
\to \mathfrak{r}_*\circ \diag_{\sigma,\,*}\circ 
\mathfrak{i}^{\an,\,*}\circ \mathfrak{j}^{\an}_*
\circ \mathfrak{u}^{\an,\,*}_{\eta}\circ \An_{B_{\eta}}^*$$
admitting a section. We claim that this natural transformation 
is equivalent to $\rho_B:\chi_B\to \chi_{\widehat{B}}\circ 
\An^*_{B_{\eta}}$. We only explain how to identify 
$\mathfrak{r}_*\circ 
\diag_{\sigma,\,*}\circ 
\mathfrak{i}^*\circ \mathfrak{j}_*\circ 
\mathfrak{u}_{\eta}^*$ with $\chi_B$;
the identification of $\mathfrak{r}_*\circ \diag_{\sigma,\,*}\circ 
\mathfrak{i}^{\an,\,*}\circ \mathfrak{j}^{\an}_*
\circ \mathfrak{u}^{\an,\,*}_{\eta}$ with $\chi_{\widehat{B}}$
is similar and easier.

Denote by $\mathfrak{D}^{\rm sm}_{\widehat{B}}$
the diagram of schemes obtained by restricting the functor 
$\mathfrak{D}_{\widehat{B}}$ to the subcategory 
$\FSm_{\af}/\widehat{B}\subset \FRigSm_{\af}/\widehat{B}$.
Define $\mathfrak{D}^{\rm sm}_{\widehat{B},\,\sigma}$ and 
$\mathfrak{D}^{\rm sm}_{\widehat{B},\,\eta}$ similarly
and denote by 
$$\mathfrak{i}^{\rm sm}:
\mathfrak{D}^{\rm sm}_{\widehat{B},\,\sigma}
\to \mathfrak{D}^{\rm sm}_{\widehat{B}}
\qquad \text{and} \qquad 
\mathfrak{j}^{\rm sm}:\mathfrak{D}^{\rm sm}_{\widehat{B},\,\eta}
\to \mathfrak{D}^{\rm sm}_{\widehat{B}}$$
the obvious inclusions. We also consider the diagonal functor 
${\rm diag}^{\rm sm}_{\sigma}:
\FSm/\widehat{B}\to \Sm/\mathfrak{D}^{\rm sm}_{\widehat{B},\,\sigma}$
sending a formal scheme $\mathcal{U}$ to the pair 
$(\mathcal{U},\mathcal{U}_{\sigma})$. With these notations, we
have an equivalence
$$\mathfrak{r}_*\circ \diag_{\sigma,\,*}\circ 
\mathfrak{i}^*\circ \mathfrak{j}_*\circ 
\mathfrak{u}_{\eta}^*\simeq \diag^{\rm sm}_{\sigma,\,*}\circ 
\mathfrak{i}^{{\rm sm},\,*}\circ \mathfrak{j}^{\rm sm}_*\circ 
\mathfrak{u}^{{\rm sm},\,*}_{\eta}.$$
Now, remark that the diagram of schemes 
$\mathfrak{D}^{\rm sm}_{\widehat{B}}$ takes values in 
regular $B$-schemes. By Lemma \ref{lem:regular-base-change}, 
we deduce an equivalence
$$\mathfrak{u}_{\sigma}^{{\rm sm},\,*} \circ \chi_B=
\mathfrak{u}_{\sigma}^{{\rm sm},\,*} \circ i^* \circ j_*
\simeq 
\mathfrak{i}^{{\rm sm},\,*}\circ \mathfrak{j}^{\rm sm}_*\circ 
\mathfrak{u}^{{\rm sm},\,*}_{\eta}.$$
We conclude by remarking that $\diag_{\sigma,\,*}^{\rm sm}\circ 
\mathfrak{u}_{\sigma}^{{\rm sm},\,*}$ is equivalent 
to the identity functor.
\end{proof}

We are now almost ready to finish the proof of Theorem 
\ref{thm:morphism-of-specialisation-equival}, but
we still need two results which are of independent interest.
The following is a version of 
\cite[Proposition 2.2.27(2)]{ayoub-th1}
with integral coefficients.

\begin{prop}
\label{prop:compact-gener-alteration-direct-image}
Let $B$ be a $(\Lambda,\et)$-admissible 
scheme, $B_{\sigma}\subset B$ 
a closed subscheme, and $B_{\eta}\subset B$ its open complement.
Assume one of the following alternatives:
\begin{itemize}

\item $B$ is quasi-compact and quasi-excellent of characteristic zero;

\item $B$ is of finite type over a quasi-compact and 
quasi-excellent scheme of dimension $\leq 1$.

\end{itemize}
Assume that every prime number is invertible either 
in $\pi_0\Lambda$ or in $\mathcal{O}(B)$.
Then, the $\infty$-category $\SH^{\hyp}_{\et}(B_{\eta};\Lambda)$
is compactly generated, up to desuspension and Tate twists, 
by motives of the form $f_{\eta,\,*}\Lambda$,
where $f:X \to B$ is a proper morphism with $X$ regular
and such that $X_{\sigma}$ is a normal crossing divisor.
\end{prop}

\begin{proof}
By \cite[Theorem 1.1]{temkin-excellent} and 
\cite[Theorem 5.13]{dejong-curves}, given a finite type 
$B$-scheme $X$ with $X_{\eta}$ integral and dense in $X$, 
we may find a proper morphism $e:X'\to X$ such that:
\begin{enumerate}

\item[(1)] $X'$ is regular and 
$X'_{\sigma}$ is a strict normal crossing 
divisor of $X'$;

\item[(2)] $X'_{\eta}$ is integral and dense in $X'$, 
and $X'\to X$ is dominant and generically finite;

\item[(3)] there exists a finite group $G$ acting on the $X$-scheme 
$X'$ and a dense open $U\subset X_{\eta}$ 
with inverse image $U'\subset X'_{\eta}$, such that the morphism 
$U'\to U$ factors as a finite \'etale Galois cover 
$U'\to U'/G$ with group $G$ and a universal homeomorphism 
$U'/G \to U$.

\end{enumerate}
Now, let $\mathcal{T}$ (resp. $\mathcal{T}'$)
be the smallest full sub-$\infty$-category of 
$\SH^{\hyp}_{\et}(B_{\eta};\Lambda)$ closed under 
colimits, desuspension and Tate twists, and containing 
the motives of the form $f_{\eta,\,*}\Lambda$,
where $f:X \to B$ is a proper morphism 
(resp. a proper morphism with $X$ regular
and $X_{\sigma}$ a normal crossing divisor).
By \cite[Lemme 2.2.23]{ayoub-th1}, we have
$\mathcal{T}=\SH^{\hyp}_{\et}(B_{\eta};\Lambda)$,
and it is enough to show that $\mathcal{T}\subset \mathcal{T}'$.
Said differently, we need to show that 
$f_{\eta,\,*}\Lambda\in \mathcal{T}'$ for any proper morphism 
$f:X \to B$. We argue by induction on the dimension of 
$X_{\eta}$.

Given a dense open immersion $j:U\to X_{\eta}$, we have an equivalence
\begin{equation}
\label{eq-prop:compact-gener-alteration-direct-image}
(f_{\eta,\,*}\Lambda \in \mathcal{T}') \Leftrightarrow
(f_{\eta,\,*}j_!\Lambda\in \mathcal{T}')
\end{equation}
by the induction hypothesis and the localisation property.
Thus, given a proper morphism $e_1:X_1\to X$ 
such that $e_1^{-1}(U)$ is dense in $X_{1,\,\eta}$ and 
$e_1^{-1}(U) \simeq U$, we may replace 
$X$ with $X_1$. Applying this to the normalisation of 
$X$, we reduce to the case where $X$ is integral
and $X_{\eta}$ dense in $X$.

Now, let $e:X'\to X$, $G$, $U$ and $U'$ be as in (1)--(3)
above. Set $f'=f\circ e$, and denote by $j:U\to X_{\eta}$ 
and $j':U'\to X'_{\eta}$ the obvious inclusions.
Then $f'_{\eta,\,*}\Lambda\in \mathcal{T}'$ by definition
and $f'_{\eta,\,*}j'_!\Lambda \in \mathcal{T}'$
by the equivalence
\eqref{eq-prop:compact-gener-alteration-direct-image},
for $X'$ instead of $X$, which is also valid under the 
induction hypothesis since $X'_{\eta}$ 
has the same dimension as $X_{\eta}$.
Moreover, by the equivalence
\eqref{eq-prop:compact-gener-alteration-direct-image},
we only need to show that 
$f_{\eta,\,*}j_!\Lambda \in \mathcal{T}'$.
Since $\mathcal{T}'$ is closed under colimits, it is enough 
to show that 
$$f_{\eta,\,*}j_!\Lambda\simeq \underset{G}{\colim}\; 
f'_{\eta,\,*}j'_!\Lambda$$
where $f'_{\eta,\,*}j'_!\Lambda$ is endowed with the $G$-action 
induced from the action of $G$ on $X'$.
Let $u:U'\to U$ and $v:U'/G \to U$ be the obvious morphisms.
Since $e$ is proper, we have $f'_{\eta,\,*}j'_!\Lambda\simeq
f_{\eta,\,*} j_! u_*\Lambda$. Since $f_{\eta,\,*}$ and
$j_!$ commute with colimits, we have 
$$\underset{G}{\colim} \; f'_{\eta,\,*}j'_!\Lambda
\simeq f_{\eta,\,*}j_!  (\underset{G}{\colim} \; u_*\Lambda).$$ 
Thus, we are left to show that 
$\Lambda \to \colim_G\,u_*\Lambda$ is an 
equivalence. By \'etale descent, we have 
$v_*\Lambda\simeq \colim_G\,u_*\Lambda$
and by Theorem
\ref{thm:sepalgebr} 
we have $\Lambda\simeq v_*\Lambda$. 
This finishes the proof. 
\end{proof}

The following is a generalisation of 
\cite[Th\'eor\`eme 7.4]{ayoub-etale}.

\begin{prop}
\label{prop:absolute-purity-SH-}
Let $S$ be a regular $(\Lambda,\et)$-admissible scheme 
and assume that every prime number is invertible either 
in $\pi_0\Lambda$ or in $\mathcal{O}(S)$.
Let 
$$\xymatrix{T' \ar[r]^-{s'} \ar[d]^-{t'} & T \ar[d]^-t\\
S' \ar[r]^-s & S}$$
be a transversal square of closed immersions in the
sense of \cite[D\'efinition 7.2]{ayoub-etale}.
Then, the morphism 
$s'^*t^!\Lambda \to t'^!s^*\Lambda$ 
is an equivalence in $\SH^{\hyp}_{\et}(T';\Lambda)$.
\end{prop}

\begin{proof}
More generally, given a $\Lambda$-module $M\in \Mod_{\Lambda}$, 
we will prove that $s'^*t^! M \to t'^!s^*M$ 
is an equivalence. Since the functors $s^*$, $t^!$, $s'^*$
and $t'^*$ are colimit-preserving, we may assume that 
$M$ is compact. When $\Lambda$ is the Eilenberg--Mac Lane spectrum 
associated to an ordinary ring, this is 
\cite[Th\'eor\`eme 7.4]{ayoub-etale}. 
It follows that the proposition is known if $\pi_0\Lambda$ 
is a $\Q$-algebra or, said differently, if we replace 
$M$ by $M_{\Q}=M\otimes \Q$. Thus, we are left to treat the case
where $M$ is $\ell$-nilpotent for a prime $\ell$
invertible on $S$. We may apply Theorem 
\ref{thm:rigrig-algebraic} and work with the 
$\infty$-categories of \'etale sheaves 
$\Shv^{\hyp}_{\et}(\Et/(-);\Lambda)_{\ellnil}$
instead of $\SH^{\hyp}_{\et}(-;\Lambda)$.
We have equivalences 
$$t^!M \simeq \lim_r t^!(M\otimes_{\Lambda}\tau_{\leq r}\Lambda)
\qquad \text{and} \qquad
t'^! M \simeq \lim_r t'^!(M\otimes_{\Lambda}\tau_{\leq r}\Lambda).$$
Since $S$ is $(\Lambda,\et)$-admissible and $M$ is compact, 
Lemma \ref{lem:pi-0-Lambda-coh-dim} 
implies that the inverse system 
$(t^!(M\otimes_{\Lambda}\tau_{\leq r}\Lambda))_r$ in 
$\Shv_{\et}^{\hyp}(\Et/T;\Lambda)$
is eventually constant on homotopy sheaves.
It follows that 
$$s'^*t^!M \simeq \lim_r 
s'^*t^!(M\otimes_{\Lambda}\tau_{\leq r}\Lambda).$$
Thus, it is enough to prove that the maps
$$s'^*t^!(M\otimes_{\Lambda}\tau_{\leq r}\Lambda)
\to t'^!s^*(M\otimes_{\Lambda}\tau_{\leq r}\Lambda)$$
are equivalences. Said differently, we may assume that 
$\Lambda$ is eventually coconnective. 
By an easy induction, we reduce to the case where 
$\Lambda$ is the Eilenberg--Mac Lane spectrum associated 
to $\Z/\ell$. (See the proof of Lemma 
\ref{lem:alternative-iii-vs-iv-3i}.) 
In this case, the result is proven in 
\cite[Proposition 7.8]{ayoub-etale}
as a consequence of Gabber's absolute purity 
\cite[Expos\'e XVI, Th\'eor\`eme 3.1.1]{Travaux-Gabber}.
\end{proof}

\begin{cor}
\label{cor:computing-chi-alg}
Let $B$ be a $(\Lambda,\et)$-admissible scheme, $B_{\sigma}\subset B$ 
a closed subscheme, and $B_{\eta}\subset B$ its open complement.
Below, we use Notation \ref{nota:chi-B-specialisation-system}.
\begin{enumerate}

\item[(1)] Assume that $B$ is regular and that $B_{\sigma}$ is a 
regular subscheme of codimension $c$ defined as the vanishing locus
of a global regular sequence $a_1,\ldots, a_c\in \mathcal{O}(B)$.
Then, we have equivalences 
$$i^!\Lambda\simeq \Lambda(-c)[-2c] \qquad \text{and} \qquad
\chi_B\Lambda\simeq \Lambda\oplus \Lambda(-c)[-2c+1]$$
in $\SH^{\hyp}_{\et}(B_{\sigma};\Lambda)$.

\item[(2)] Assume that $B$ is regular and that 
$B_{\sigma}$ is a strict normal crossing 
divisor. Let $D\subset B_{\sigma}$ be an irreducible component 
and $D^{\circ}$ the intersection of $D$ with the regular locus
of $(B_{\sigma})_{\red}$.
Let $u:D^{\circ} \to D$ and $v:D \to B_{\sigma}$ be
the obvious inclusions. The morphism 
$$v^*\chi_B\Lambda \to u_*u^*v^*\chi_B\Lambda$$
is an equivalence in $\SH^{\hyp}_{\et}(D;\Lambda)$.

\end{enumerate}
\end{cor}

\begin{proof}
For the first assertion, we consider the commutative diagram
with Cartesian squares
$$\xymatrix{B_{\eta} \ar[r]^-j \ar[d]^-{a_{\eta}} 
& B\ar[d]^-a & 
B_{\sigma} \ar[l]_-i \ar[d]^-{a_{\sigma}} \\
\A^c_B\smallsetminus 0_B \ar[r]^-{j_0} & \A^c_B & \ar[l]_{i_0} B,\!}$$
where $a$ is the section of $\A^c_B \to B$ induced by the 
$c$-tuple $(a_1,\ldots, a_c)$ and $i_0$ is the zero section. 
By Proposition \ref{prop:absolute-purity-SH-}, 
we have equivalences 
$i^!\Lambda\simeq a_{\sigma}^*i_0^!\Lambda$ 
and $\chi_B\Lambda\simeq a_{\sigma}^*\chi_{\A^c_B}\Lambda$,
which enable us to conclude.

We now pass to the second assertion. 
Since the problem is local over $B$, we may assume that 
$(B_{\sigma})_{\red}$ 
is defined by an equation of the form
$a_1\cdots a_c=0$, where $a_1,\ldots, a_c$ is a regular sequence.
Consider the commutative diagram
with Cartesian squares
$$\xymatrix{B_{\eta} \ar[r]^-j \ar[d]^-{a_{\eta}} 
& B\ar[d]^-a & 
B_{\sigma} \ar[l]_-i \ar[d]^-{a_{\sigma}} \\
U \ar[r]^-{j'} & \A^c_B & \ar[l]_{i'} E,\!}$$
where $E$ is defined by the equation $t_1\cdots t_c=0$,
with $(t_1,\ldots, t_c)$ a system of coordinates on $\A^c$,
and $U=\A^c_B\smallsetminus E$.
For $I\subset \{1,\ldots, c\}$ nonempty, we let 
$D_I \subset B_{\sigma}$ and $H_I\subset E$
be the closed subschemes defined by the equations
$\prod_{i\in I} a_i=0$ and $\prod_{i\in I} t_i=0$
respectively. We have transversal squares
$$\xymatrix{D_I \ar[r]^-{i_I} \ar[d]^-{a_I} & B\ar[d]^-a\\
H_I \ar[r]^-{i'_I} & \A^c_B.\!}$$
By Proposition \ref{prop:absolute-purity-SH-},
we deduce equivalences $a_I^*i'^!_I\Lambda\simeq i_I^!\Lambda$. 
Since $i^!\Lambda$ and 
$i'^!\Lambda$ can be built from the $i_I^!\Lambda$'s 
and the $i_I'^!\Lambda$'s using the same recipe, we deduce that the 
obvious map $a_{\sigma}^*i'^!\Lambda\to i^!\Lambda$ 
is an equivalence. It follows that 
$$a_{\sigma}^*\chi_{\A^c_B}\Lambda
\to \chi_B\Lambda$$
is also an equivalence. 
We may assume that $D=D_1$. We set $H=H_1$ and define
$H^{\circ}$ as in the statement.
We also let $v':H\to E$ and $u':H^{\circ} \to H$
be the obvious inclusions.
By \cite[Th\'eor\`eme 3.3.11]{ayoub-th2}, the obvious map
$$v'^*\chi_{\A^c_B}\Lambda
\to u'_*u'^*v'^*\chi_{\A^c_B}\Lambda$$
is an equivalence. We have a commutative diagram
$$\xymatrix{a_1^*v'^*\chi_{\A^c_B}\Lambda \ar[r] \ar[d]^-{\sim} 
\ar[r]^-{\sim} & v^*a_{\sigma}^*\chi_{\A^c_B}\Lambda \ar[r]^-{\sim} & v^*\chi_B\Lambda\ar[d]\\
a_1^*u'_*u'^*v'^*\chi_{\A^c_B}\Lambda\ar[r] & u_*u^*v^*a_{\sigma}^*
\chi_{\A^c_B}\Lambda \ar[r]^-{\sim} & u_*u^*v^*\chi_B\Lambda.\!}$$
So, we are left to show that the morphism
$a_1^*u'_*u'^*v'^*\chi_{\A^c_B}\Lambda\to u_*u^*v^*a_{\sigma}^*
\chi_{\A^c_B}\Lambda$ is an equivalence. 
For this, we remark that 
$u'^*v'^*\chi_{\A^c_B}\Lambda \simeq \Lambda\oplus \Lambda(-1)[-1]$
in $\SH^{\hyp}_{\et}(H^{\circ};\Lambda)$, and that this morphism
is equivalent to 
$$a^*_1u'_*(\Lambda \oplus \Lambda(-1)[-1]) \to 
u_*(\Lambda \oplus \Lambda(-1)[-1]).$$ 
Thus, it remains to show that 
$b^*z'^!\Lambda \to z^!\Lambda$
is an equivalence, with $z:D\smallsetminus D^{\circ} \to D$, 
$z':H\smallsetminus H^{\circ} \to H$ and
$b:D\smallsetminus D^{\circ} \to H\smallsetminus H^{\circ}$
the obvious morphisms. This is proven in the same way we 
proved above that $a^*_{\sigma}i'^!\Lambda \to i^!\Lambda$ was an 
equivalence.
\end{proof}

We are finally ready to conclude.

\begin{proof}[Proof of Theorem 
\ref{thm:morphism-of-specialisation-equival}]
By Lemma 
\ref{lem:reduct-for-iso-specialisation-systems}, 
we may assume that $B$ is essentially of finite type over 
$\Spec(\Z)$ and work in the hypercomplete case.
Since the source and target of $\rho_X$ consist of 
colimit-preserving functors, it is enough to prove
that $\chi_X M \to \chi_{X^{\an}}\An_{X_{\eta}}^*M$
is an equivalence when $M$ belongs to set of 
compact generators of $\SH^{\hyp}_{\et}(X_{\eta};\Lambda)$.
By Proposition \ref{prop:compact-gener-alteration-direct-image}, 
we may assume that $M=f_{\eta,\,*}\Lambda$
where $f:Y \to X$ is a proper morphism such that 
$Y$ is regular and $Y_{\sigma}$ is a normal crossing divisor. 
By the proper base change theorem, we have equivalences
$$\chi_X f_{\eta,\,*}\Lambda\simeq 
f_{\sigma,\,*}\chi_Y\Lambda \qquad \text{and} \qquad
\chi_{X^{\an}}\An_{X_{\eta}}^* f_{\eta,\,*}\Lambda
\simeq f_{\sigma,\,*} \chi_{Y^{\an}}\An_{Y_{\eta}}^*\Lambda
\simeq f_{\sigma,\,*} \chi_{Y^{\an}}\Lambda.$$
Thus, replacing $X$ with $Y$, we may assume that 
$X$ is regular and $X_{\sigma}$ a strict normal crossing divisor
and, in this case, we only need to show that 
$\chi_X\Lambda\to \chi_{X^{\an}}\Lambda$ 
is an equivalence.
By Proposition \ref{prop-thm:morphism-of-specialisation-equival-weak},
this morphism admits a section, and thus 
$\chi_{X^{\an}}\Lambda$
is the image of a projector $p$ of 
$\chi_X\Lambda$. We need to prove that $p$ is the identity,
and it is enough to show this after restriction to each 
irreducible component of $X_{\sigma}$. Using Corollary 
\ref{cor:computing-chi-alg}(2), it is enough to do so 
after restricting to the regular locus of $X_{\sigma}$. 
Said differently, we may assume that 
$X_{\sigma}$ is a regular divisor.

From now on, we assume that $X$ is regular 
and that $X_{\sigma}$ is 
a regular divisor defined by the zero locus of 
$a\in \mathcal{O}(X)$.
We denote by $p$ the projector of $\chi_X$ provided by 
Proposition 
\ref{prop-thm:morphism-of-specialisation-equival-weak}.
Our goal is to show that $p$ acts on $\chi_X\Lambda
\simeq \Lambda\oplus \Lambda(-1)[-1]$
by the identity, and it is enough to show that 
$p$ is an equivalence. First, note that we have 
a commutative square
$$\xymatrix{\Lambda \ar[r] \ar@{=}[d] & \chi_X\Lambda  
\ar[d]^-p \\
\Lambda \ar[r] & \chi_X\Lambda}$$
since $p$ is an algebra endomorphism of $\chi_X\Lambda$. 
(Indeed, the section constructed in Proposition
\ref{prop-thm:morphism-of-specialisation-equival-weak}
respects the natural right-lax monoidal structures.)
Thus, with respect to the decomposition 
$\chi_X\Lambda
\simeq \Lambda\oplus \Lambda(-1)[-1]$,
$p$ is given by a triangular matrix
$$p=\begin{pmatrix}
1 & r\\
0 & q
\end{pmatrix}.$$
We will show that $q$ is the identity of $\Lambda(-1)[-1]$. 
To do so, we consider the morphism
$\Lambda \to \Lambda(1)[1]$ in $\SH^{\hyp}_{\et}(X_{\eta};\Lambda)$
corresponding to $a\in \mathcal{O}^{\times}(X_{\eta})$, i.e., 
induced by the section 
$a:X_{\eta} \to \A^1_{X_{\eta}}\smallsetminus 
0_{X_{\eta}}$.
Applying $\chi_X$ and then 
$p:\chi_X \to \chi_X$ yields a commutative square
$$\xymatrix@R=3pc{\Lambda\oplus \Lambda(-1)[-1] 
\ar[r]^-{\begin{pmatrix}\scriptstyle{0\,0} \\
\vspace{-.7cm}\\
\scriptstyle{1\,0}\end{pmatrix}} 
\ar[d]^-{\begin{pmatrix}\scriptstyle{1\,r} \\
\vspace{-.7cm}\\
\scriptstyle{0\,q}\end{pmatrix}} & \Lambda(1)[1] \oplus \Lambda
\ar[d]^-{\begin{pmatrix}\scriptstyle{1\,r} \\
\vspace{-.7cm}\\
\scriptstyle{0\,q}\end{pmatrix}} \\
\Lambda\oplus \Lambda(-1)[-1] 
\ar[r]^-{\begin{pmatrix}\scriptstyle{0\,0} \\
\vspace{-.7cm}\\
\scriptstyle{1\,0}\end{pmatrix}} & \Lambda(1)[1] \oplus \Lambda.\!}$$
This forces $q$ to be the identity, as needed.
\end{proof}

\section{The six-functor formalism for rigid analytic motives}

\label{sec:6ff}

In this section, we develop the six-functor formalism
for rigid analytic motives, getting rid of the quasi-projectivity 
assumption imposed in \cite[\S 1.4]{ayoub-rig}. The 
key step in doing so is to prove an extended 
proper base change theorem for rigid analytic motives; 
see Theorem \ref{thm:prop-base} below.
An important particularity in the rigid analytic setting is the 
existence of canonical compactifications (aka., Huber compactifications).
We will not make use of these compactifications in defining the 
exceptional direct image functors, but see Theorem
\ref{thm:f-!-using-huber-compactification} below.

\subsection{Extended proper base change theorem}

$\empty$

\smallskip

\label{subsect:proper-base-change}

Our goal in this subsection is to prove a general
extended proper base change theorem for rigid analytic motives;
see Theorem \ref{thm:prop-base} below. 
This will be achieved by reducing to the usual proper base change
theorem for algebraic motives. A 
compatibility property for the functors $\chi_{\mathcal{S}}$,
for $\mathcal{S}\in \FSch$, and the operations 
$f_{\sharp}$, for $f$ smooth, plays a key role in this 
reduction; it is given in Theorem 
\ref{thm:chijsharp} below which we deduce quite easily
from Theorem \ref{thm:proj-form-chi-xi}
(which was a key step in proving 
Theorem \ref{thm:main-thm-}).
We start by a well-known generalisation of some facts contained in
\cite[Scholie 1.4.1]{ayoub-th1}.

\begin{prop}
\label{prop:prop-base-proformal}
\ncn{proper base change}
\ncn{smooth base change}
Consider a Cartesian square in $\FSch$
$$\xymatrix{\mathcal{Y}' \ar[r]^-{g'} \ar[d]^-{f'}
& \mathcal{Y} \ar[d]^-f\\
\mathcal{X}' \ar[r]^-g & \mathcal{X}}$$
with $f$ proper. 
\begin{enumerate}

\item[(1)] The 
commutative square 
$$\xymatrix{\FSH^{(\hyp)}_{\tau}(\mathcal{X};\Lambda) \ar[r]^-{f^*}
\ar[d]^-{g^*} & 
\FSH^{(\hyp)}_{\tau}(\mathcal{Y};\Lambda) \ar[d]^-{g'^*}\\
\FSH^{(\hyp)}_{\tau}(\mathcal{X}';\Lambda)
\ar[r]^-{f'^*} & \FSH^{(\hyp)}_{\tau}(\mathcal{Y}';\Lambda)}$$
is right adjointable, i.e., the natural transformation 
$g^*\circ f_* \to f'_*\circ g'^*$
is an equivalence.

\item[(2)] If $g$ is smooth, the commutative square
$$\xymatrix{\FSH^{(\hyp)}_{\tau}(\mathcal{X}';\Lambda) \ar[r]^-{f'^*}
\ar[d]^-{g_{\sharp}} & 
\FSH^{(\hyp)}_{\tau}(\mathcal{Y}';\Lambda) \ar[d]^-{g'_{\sharp}}\\
\FSH^{(\hyp)}_{\tau}(\mathcal{X};\Lambda)
\ar[r]^-{f^*} & \FSH^{(\hyp)}_{\tau}(\mathcal{Y};\Lambda)}$$
is right adjointable, i.e., the natural transformation 
$g_{\sharp}\circ f'_* \to f_*\circ g'_{\sharp}$
is an equivalence.

\end{enumerate}
\end{prop}

\begin{proof}
By Theorem \ref{thm:formal-mot-alg-mot}, 
we reduce to show the statement for a Cartesian 
square in $\Sch$
$$\xymatrix{Y' \ar[r]^{g'} \ar[d]_-{f'} & Y \ar[d]^-f\\
X' \ar[r]^g & X}$$
with $f$ proper. When $f$ is projective, this is covered by 
\cite[Scholie 1.4.1]{ayoub-th1}; see also 
\cite[Proposition 3.5]{ayoub-etale}. 
The passage from the projective to the proper case is a 
well-known procedure, that we revisit here because
we don't know a reference in the generality we are considering. 
(Under noetherianness 
assumptions, an argument can be found in the proof of 
\cite[Proposition 2.3.11(2)]{cd}.)

The question is local on $X$, so we may assume that 
$X$ is quasi-compact and quasi-separated. 
Using a covering of $Y$ by finitely many affine open subschemes, 
assertion (1) (resp. assertion (2)) 
follows if we can prove that the natural transformation
$$g^*\circ f_* \circ v_{\sharp}\to 
f'_*\circ g'^*\circ v_{\sharp}
\qquad \text{(resp. }
g_{\sharp} \circ f'_*\circ v'_{\sharp}
\to f_*\circ g'_{\sharp}\circ v'_{\sharp}\text{)}$$
is an equivalence for every open immersion 
$v:V \to Y$ with base change $v':V'\to Y'$.
Letting $g'':V'\to V$ be the base change of $g'$, 
this natural transformation can be rewritten as follows:
$$g^*\circ (f_* \circ v_{\sharp})\to 
(f'_*\circ v'_{\sharp})\circ g''^*
\qquad \text{(resp. }
g_{\sharp} \circ (f'_*\circ v'_{\sharp})
\to (f_*\circ v_{\sharp})\circ g''_{\sharp}\text{)}.$$
By the refined version of Chow's lemma given in
\cite[Corollary 2.6]{nagata-deligne},
we may find a blowup $e:Z\to Y$, with centre
disjoint from $V$, such that $h=f\circ e$ is a projective
morphism. Let $w:V \to Z$ be the open immersion such that 
$v=e\circ w$. Set $Z'=Z\times_Y Y'$ and let  
$e':Z'\to Y'$, $h':Z'\to X'$ and $w':V'\to Z'$
be the base change of $e$, $h$ and $w$ along $g$.
Using \cite[Scholie 1.4.1]{ayoub-th1}, 
we have natural equivalences 
$v_{\sharp}\simeq e_*\circ w_{\sharp}$ and 
$v'_{\sharp}\simeq e'_*\circ w'_{\sharp}$.
Thus, we may rewrite the above natural transformation as follows:
$$g^*\circ (h_* \circ w_{\sharp})\to 
(h'_*\circ w'_{\sharp})\circ g''^*
\qquad \text{(resp. }
g_{\sharp} \circ (h'_*\circ w'_{\sharp})
\to (h_*\circ w_{\sharp})\circ g''_{\sharp}\text{)}.$$
Thus, we may replace $f$ and $f'$ by $h$ and $h'$,
thereby reducing the general case to the case of a 
projective morphism.
\end{proof}

\begin{lemma}
\label{lem:direct-image-in-prl}
Let $f:\mathcal{Y} \to \mathcal{X}$ 
be a proper morphism of formal schemes. Then, 
the functor 
$$f_*:\FSH^{(\hyp)}_{\tau}(\mathcal{Y};\Lambda) \to 
\FSH^{(\hyp)}_{\tau}(\mathcal{X};\Lambda)$$
is colimit-preserving and thus admits a right adjoint.
\end{lemma}

\begin{proof}
By Theorem \ref{thm:formal-mot-alg-mot}, 
we reduce to show the statement for a proper morphism
of schemes $f:Y \to X$.
When $f$ is projective, this follows from 
\cite[Th\'eor\`eme 1.7.17]{ayoub-th1}. In general, we may assume that 
$X$ is quasi-compact and quasi-separated, and reduce to show 
that $f_*\circ v_{\sharp}$ is colimit-preserving for every 
open immersion $v:V\to Y$ with $V$ affine. Then, we use 
the refined version of Chow's lemma given in
\cite[Corollary 2.6]{nagata-deligne},
to find a blowup $Y'\to Y$ with centre disjoint from $V$ and 
such that $Y'\to X$ is projective. We conclude using the
equivalence $f_*\circ v_{\sharp}\simeq f'_*\circ v'_{\sharp}$
where $f':Y'\to X$ and $v':V \to Y'$ are the obvious morphisms.
\end{proof}

Our main task in this subsection is to prove 
a variant of Proposition \ref{prop:prop-base-proformal}
for rigid analytic motives. (A version of 
Proposition \ref{prop:prop-base-proformal}(a) 
holds true in the rigid analytic setting even without  
assuming that $f$ is proper but under some mild technical 
assumptions; see Theorem 
\ref{thm:general-base-change-thm}. We will explain below 
how to remove these technical assumptions
when $f$ is assumed to be proper.) 
A key ingredient is provided by the following theorem.

\begin{thm}
\label{thm:chijsharp}
We work under Assumption 
\ref{assu:for-main-thm}.
Let $f:\mathcal{T} \to \mathcal{S}$ 
be a smooth morphism of formal schemes. 
The commutative square
$$\xymatrix{\FSH^{(\hyp)}_{\tau}(\mathcal{T};\Lambda) 
\ar[r]^-{\xi_{\mathcal{T}}} \ar[d]^-{f_{\sharp}} & 
\RigSH^{(\hyp)}_{\tau}(\mathcal{T}^{\rig};\Lambda) 
\ar[d]^-{f^{\rig}_{\sharp}}\\ 
\FSH^{(\hyp)}_{\tau}(\mathcal{S};\Lambda) 
\ar[r]^-{\xi_{\mathcal{S}}}  & 
\RigSH^{(\hyp)}_{\tau}(\mathcal{S}^{\rig};\Lambda)}$$
is right adjointable, i.e., the induced natural transformation
$f_{\sharp}\circ \chi_{\mathcal{T}} \to \chi_{\mathcal{S}} 
\circ f^{\rig}_{\sharp}$
is an equivalence.
\end{thm}

\begin{proof}
We split the proof into two steps. In the first one, we consider 
the case where $f$ is an open immersion and, in the second one, 
we treat the general case.

\paragraph*{Step 1}
\noindent 
Here we treat the case of 
an open immersion $j:\mathcal{U} \to \mathcal{S}$.
For $M\in \RigSH^{(\hyp)}_{\tau}(\mathcal{S}^{\rig};\Lambda)$, we have a 
commutative diagram
$$\xymatrix{\chi_{\mathcal{S}}(M)\otimes 
j_{\sharp}\Lambda \ar[r]^-{(1)} \ar[d]^-{\sim} 
& \chi_{\mathcal{S}}(M\otimes \xi_{\mathcal{S}} j_{\sharp}\Lambda) 
\ar[r]^-{\sim} & \chi_{\mathcal{S}}(M\otimes j^{\rig}_{\sharp}\Lambda) \ar[d]^-{\sim}\\
j_{\sharp}j^*\chi_{\mathcal{S}}M \ar[r]^-{\sim} & 
j_{\sharp}\,\chi_{\mathcal{U}}\,j^{\rig,\,*}M
\ar[r]^-{(2)} &
\chi_{\mathcal{S}}j^{\rig}_{\sharp}j^{\rig,\,*}M,\!}$$
where all the arrows, except the labeled ones, 
are equivalences for obvious reasons.  
By Theorem \ref{thm:proj-form-chi-xi}, the morphism (1)
is also an equivalence, and hence the same is true for 
the morphism (2).
Thus, the natural transformation $j_{\sharp}\circ \chi_{\mathcal{U}}
\to \chi_{\mathcal{S}}\circ j^{\rig}_{\sharp}$
becomes an equivalence when applied to the functor 
$j^{\rig,\,*}$. Since the latter is essentially surjective, 
the result follows.

\paragraph*{Step 2} 
\noindent
Here we treat the general case.
Clearly, the problem is local on $\mathcal{S}$. 
We claim that it is also local on $\mathcal{T}$.
Indeed, let $(u_i:\mathcal{T}_i \to \mathcal{T})_i$
be an open covering of $\mathcal{T}$.
The $\infty$-category 
$\RigSH^{(\hyp)}_{\tau}(\mathcal{T}^{\rig};\Lambda)$
is generated under colimits by the images of the functors 
$u^{\rig}_{i,\,\sharp}$. Clearly, the functors 
$f_{\sharp}$ and $f_{\sharp}^{\rig}$ are colimit-preserving.
By Proposition
\ref{prop:chi-commute-direct-sums}, the same is true for 
$\chi_{\mathcal{T}}$ and $\chi_{\mathcal{S}}$.
Thus, it is enough to prove that the natural transformations
$f_{\sharp}\circ \chi_{\mathcal{T}}\circ u^{\rig}_{i,\,\sharp}
\to \chi_{\mathcal{S}}\circ f^{\rig}_{\sharp}\circ u^{\rig}_{i,\,\sharp}$
are equivalences. Using the first step, this natural transformation 
is equivalent to $(f\circ u_i)_{\sharp}\circ 
\chi_{\mathcal{T}_i} \to \chi_{\mathcal{S}}\circ 
(f\circ u_i)_{\sharp}^{\rig}$ which brings us to prove 
the theorem for the morphisms 
$f\circ u_i$. This proves our claim.

The problem being local on $\mathcal{T}$ and $\mathcal{S}$, 
we may assume that there is a closed immersion 
$i:\mathcal{T} \to \A^n_{\mathcal{S}}$. 
We may also assume that there is an \'etale neighbourhood 
of $\mathcal{T}$ in $\A^n_{\mathcal{S}}$ which is 
isomorphic to an \'etale neighbourhood of the zero section 
$\mathcal{T}\to \A^m_{\mathcal{T}}$ (where $m$ is the codimension
of the immersion $i$). Thus, letting  
$p:\A^n_{\mathcal{S}}\to \mathcal{S}$ be the obvious projection,
we have natural equivalences
$$p_{\sharp}\circ i_*\simeq f_{\sharp}(m)[2m]
\qquad \text{and} \qquad p_{\sharp}^{\rig}\circ i^{\rig}_*\simeq f^{\rig}_{\sharp}(m)[2m].$$
Moreover, the following diagram is commutative
$$\xymatrix{p_{\sharp} \circ i_*\circ \chi_{\mathcal{T}}
\ar[r]^-{\sim} \ar[d]^-{\sim} 
& p_{\sharp}\circ \chi_{\A^n_{\mathcal{S}}}\circ i^{\rig}_*
\ar[r] & \chi_{\mathcal{S}} \circ p^{\rig}_{\sharp}\circ 
i^{\rig}_* \ar[d]^-{\sim} \\
f_{\sharp}\circ \chi_{\mathcal{T}}(m)[2m] \ar[rr] & & 
\chi_{\mathcal{S}}\circ f^{\rig}_{\sharp}(m)[2m].\!}$$
This shows that it suffices to treat the case of the projection 
$p:\A^n_{\mathcal{S}}\to \mathcal{S}$. 

Let $j:\A^n_{\mathcal{S}}\to \P^n_{\mathcal{S}}$ be an open 
immersion into the relative projective space of dimension $n$ 
and let $q:\P^n_{\mathcal{S}}\to \mathcal{S}$ be the obvious 
projection. The morphism 
$p_{\sharp}\circ \chi_{\A^n_{\mathcal{S}}}
\to \chi_{\mathcal{S}}\circ p^{\rig}_{\sharp}$
is equivalent to the composition of 
$$q_{\sharp} \circ j_{\sharp}\circ \chi_{\A^n_{\mathcal{S}}}
\to q_{\sharp} \circ \chi_{\P^n_{\mathcal{S}}}
\circ j^{\rig}_{\sharp} \to 
\chi_{\mathcal{S}}
\circ q^{\rig}_{\sharp} \circ j^{\rig}_{\sharp}$$
and the first morphism is an equivalence by the first step. 
Thus, we are left to treat the case of $q:\P^n_{\mathcal{S}}
\to \mathcal{S}$. By 
\cite[Th\'eor\`eme 1.7.17]{ayoub-th1} and Corollary 
\ref{cor:6ffalg}, we have equivalences
$$q_{\sharp}\simeq q_*\circ \Th(\Omega_q) 
\qquad \text{and} \qquad 
q^{\rig}_{\sharp}\simeq q^{\rig}_*\circ \Th(\Omega_{q^{\rig}}),$$
and the following square
$$\xymatrix{q_{\sharp} \circ \chi_{\P^n_{\mathcal{S}}} \ar[r]^-{\sim} 
\ar[d] & q_* \circ \Th(\Omega_q) \circ \chi_{\P^n_{\mathcal{S}}} 
\ar[d]^-{\sim}\\ 
\chi_{\mathcal{S}}\circ q^{\rig}_{\sharp} \ar[r]^-{\sim} & 
\chi_{\mathcal{S}} \circ q^{\rig}_*\circ \Th(\Omega_{q^{\rig}})}$$
is commutative.
This finishes the proof.
\end{proof}

Here is the main result of this subsection.

\begin{thm}[Extended proper base change]
\label{thm:prop-base}
\ncn{extended proper base change}
Consider a Cartesian square in $\RigSpc$
$$\xymatrix{Y' \ar[r]^-{g'} \ar[d]^-{f'}
& Y \ar[d]^-f\\
X' \ar[r]^-g & X}$$
with $f$ proper.
\begin{enumerate}

\item[(1)] The commutative square 
$$\xymatrix{\RigSH^{(\hyp)}_{\tau}(X;\Lambda) \ar[r]^-{f^*}
\ar[d]^-{g^*} & 
\RigSH^{(\hyp)}_{\tau}(Y;\Lambda) \ar[d]^-{g'^*}\\
\RigSH^{(\hyp)}_{\tau}(X';\Lambda)
\ar[r]^-{f'^*} & \RigSH^{(\hyp)}_{\tau}(Y';\Lambda)}$$
is right adjointable, i.e., the natural transformation 
$g^*\circ f_* \to f'_*\circ g'^*$
is an equivalence.

\item[(2)] If $g$ is smooth, the commutative square
$$\xymatrix{\RigSH^{(\hyp)}_{\tau}(X';\Lambda) \ar[r]^-{f'^*}
\ar[d]^-{g_{\sharp}} & 
\RigSH^{(\hyp)}_{\tau}(Y';\Lambda) \ar[d]^-{g'_{\sharp}}\\
\RigSH^{(\hyp)}_{\tau}(X;\Lambda)
\ar[r]^-{f^*} & \FSH^{(\hyp)}_{\tau}(Y;\Lambda)}$$
is right adjointable, i.e., the natural transformation 
$g_{\sharp}\circ f'_* \to f_*\circ g'_{\sharp}$
is an equivalence.

\end{enumerate}
\end{thm}

\begin{proof}
The question is local on $X$ and $X'$. Thus, we may assume that 
$X$ and $X'$ are quasi-compact and quasi-separated. 
We split the proof into three steps. The first two steps concern
part (2): in the first step we show that it is enough 
to treat the case where $g$ has good reduction, and in the second 
step we prove part (2) while working in the 
non-hypercomplete case and assuming that 
$\tau$ is the Nisnevich topology. Finally, in the third step, we 
use what we learned in the second step to prove the theorem 
in complete generality.

\paragraph*{Step 1}
\noindent
Here, we assume that part (2) is known when 
$g$ has good reduction and we explain how to 
deduce it in general. 
The problem being local on $X'$, we may assume that 
our Cartesian square is the composition of two 
Cartesian squares
$$\xymatrix{Y' \ar[r]^-{e'} \ar[d]^-{f'} & Y_1 \ar[r]^-{h'} 
\ar[d]^-{f_1} & Y\ar[d]^-f\\
X' \ar[r]^-e & X_1 \ar[r]^-h & X}$$
where $e$ is \'etale and $h$ is smooth with good reduction. 
(For instance, we may assume that $h$ is the projection 
of a relative ball.)
By assumption, part (2) is known for the right square, so 
it remains to prove it for the left square. 
Said differently, we may assume that $g$ is \'etale. 
Using Lemma \ref{lem:local-structure-of-rig-etale-morphisms}
below, we reduce further to the case where 
$g$ is finite \'etale. 
In this case, there is a natural equivalence 
$g_{\sharp} \simeq g_*$ 
constructed as follows. Consider the Cartesian square
$$\xymatrix{X'\times_X X' \ar[r]^-{\pr_2}\ar[d]_-{\pr_1}
& X' \ar[d]^-g\\
X' \ar[r]^-g & X,}$$
and the diagonal embedding $\Delta:X' \to X'\times_XX'$
which is a clopen immersion. 
Since $g$ is locally projective, we may use Proposition 
\ref{prop:base-change-finite-and-projection}(2)
which implies that the natural transformation 
$$g_{\sharp}\circ \pr_{1,\,*} \to g_*\circ 
\pr_{2,\,\sharp}$$
is an equivalence. 
Applying this equivalence to the functor $\Delta_{\sharp}\simeq \Delta_*$,
we get the equivalence $g_{\sharp}\simeq g_*$. 
Similarly, we have an equivalence 
$g'_{\sharp}\simeq g'_*$.
Moreover, modulo
these equivalences, the natural transformation 
$g_{\sharp}\circ f'_*\to f_*\circ g'_{\sharp}$
coincides with the obvious equivalence 
$g_*\circ f'_*\simeq f_*\circ g'_*$. 
This proves the claimed reduction.

\paragraph*{Step 2}
\noindent
We now prove part (2) of the statement under
Assumption \ref{assu:for-main-thm} so that we can use 
Theorem \ref{thm:chijsharp}. (More precisely, 
we will assume that all the formal models 
used below satisfy this 
assumption.) In the third step we explain 
how to get rid of this assumption.

The problem being local on $X$ and $X'$, 
we may also assume that $f$ is the generic fiber of a 
proper morphism $\widetilde{f}:\mathcal{Y} \to 
\mathcal{X}$ in $\FSch$ and that 
$g$ is the generic fiber of a smooth morphism 
$\widetilde{g}:\mathcal{X}'\to \mathcal{X}$
of formal schemes (since $g$ can be assumed to have 
good reduction, by the first step). We form a Cartesian square
$$\xymatrix{\mathcal{Y}' \ar[r]^-{\widetilde{g}{}'} \ar[d]^-{\widetilde{f}{}'} & \mathcal{Y} 
\ar[d]^-{\widetilde{f}}\\
\mathcal{X}' \ar[r]^-{\widetilde{g}} & \mathcal{X}.\!}$$
For every quasi-compact and quasi-separated smooth 
rigid analytic $X$-space $L$, with structural morphism 
$p_L:L\to X$, choose a formal model $\mathcal{L}$
which is a finite type formal $\mathcal{X}$-scheme.
By Proposition \ref{prop:chi-trivial-generators}, 
when $L$ varies, the functors 
$$\chi_{\mathcal{L}}\circ p_L^*:\RigSH^{(\hyp)}_{\tau}(X;\Lambda)
\to \FSH^{(\hyp)}_{\tau}(\mathcal{L};\Lambda)$$ 
form a conservative family.
Therefore, it is enough to show that the natural transformation
$$\chi_{\mathcal{L}}\circ p_L^* \circ 
g_{\sharp}\circ f'_* \to 
\chi_{\mathcal{L}}\circ p_L^*\circ f_* \circ 
g'_{\sharp}$$
is an equivalence for each $p_L:L\to X$ and $\mathcal{L}$ as above. 
Letting $f_L$, $f'_L$, $g_L$ and $g'_L$ be the base change 
of the morphisms $f$, $f'$, $g$ and $g'$ along $p_L:L\to X$, 
and using Proposition 
\ref{prop:6f1}, 
we reduce to show that
the natural transformation  
$$\chi_{\mathcal{L}}\circ  
g_{L,\,\sharp}\circ f'_{L,\,*} \to 
\chi_{\mathcal{L}} \circ f_{L,\,*} \circ 
g'_{L,\,\sharp}$$
is an equivalence. Thus, replacing $X$ with $L$ and 
$\mathcal{X}$ with $\mathcal{L}$, we may
concentrate on the natural transformation 
$$\chi_{\mathcal{X}}\circ  
g_{\sharp}\circ f'_* \to 
\chi_{\mathcal{X}} \circ f_* \circ 
g'_{\sharp}.$$
Using Theorem
\ref{thm:chijsharp}, 
we can rewrite 
this natural transformation as follows:
$$\widetilde{g}_{\sharp}\circ \widetilde{f}{}'_*\circ \chi_{\mathcal{Y}'}
\to \widetilde{f}_*\circ \widetilde{g}{}'_{\sharp}\circ 
\chi_{\mathcal{Y}'}.$$
We now conclude using Proposition 
\ref{prop:prop-base-proformal}(2).

\paragraph*{Step 3}
\noindent 
In this step, we will prove the 
theorem in complete generality. By Theorem
\ref{thm:general-base-change-thm} and the second 
step, the theorem is known for the $\infty$-categories
$\RigSH_{\Nis}(-;\Lambda)$, i.e., when $\tau$ is the Nisnevich 
topology and we work in the non-hypercomplete case.
This will be our starting point.
(Of course, by the second step, 
the theorem is known more generally, e.g., 
when $\tau$ is the Nisnevich 
topology and we work in the hypercomplete case, 
but this will not be used below.)

For a rigid analytic space $S$, the functor 
$\Lder_S:\RigSH_{\Nis}(S;\Lambda)\to \RigSH^{(\hyp)}_{\tau}(S;\Lambda)$
is a localisation functor with respect to the set 
$\mathcal{H}_S$ consisting of maps of the form
$\colim_{[n]\in\mathbf{\Delta}}\;\M(T_n)\to \M(T_{-1})$,
and their desuspensions and negative Tate twists, where 
$T_{\bullet}$ is a $\tau$-hypercover which is assumed to be
truncated in the non-hypercomplete case. We
claim that the functor 
$$f_*:\RigSH_{\Nis}(Y;\Lambda) \to \RigSH_{\Nis}(X;\Lambda)$$ 
takes $\mathcal{H}_Y$-equivalences
to $\mathcal{H}_X$-equivalences, and that the same is true for
$f'_*$. Assuming this claim, one has equivalences 
$\Lder_X\circ f_*\simeq f_*\circ \Lder_Y$, and similarly for $f'_*$.
Since the functors $\Lder_Y$ and $\Lder_{Y'}$ are essentially 
surjective on objects, it suffices to prove that the natural 
transformations 
$$g^*\circ f_*\circ \Lder_Y\to f'_*\circ g'^*\circ \Lder_Y
\quad\text{and}\quad
g_{\sharp}\circ f'_*\circ \Lder_{Y'}\to 
f_*\circ g'_{\sharp}\circ \Lder_{Y'}$$
are equivalences. Thus, using our claim and the obvious 
analogous commutations for $g^*$, $g_{\sharp}$, $g'^*$ and 
$g'_{\sharp}$, the above natural transformations are 
equivalent to 
$$\Lder_{X'}\circ g^*\circ f_* \to \Lder_{X'}\circ
f'_*\circ g'^*
\quad\text{and}\quad
\Lder_X \circ g_{\sharp}\circ f'_* \to 
\Lder_X\circ f_*\circ g'_{\sharp},$$
and the result follows. 

It remains to prove our claim, and it is enough to consider the case of 
$f$ (which is a general proper morphism).
Using a covering of $Y$ by finitely many affine open
subspaces, we see that it suffices to show that 
$f_*\circ v_{\sharp}$ takes $\mathcal{H}_V$-equivalences
to $\mathcal{H}_X$-equivalences
for every open immersion $v:V\to Y$ such that 
$V$ admits a locally closed immersion into a relative projective 
space $P\simeq \P^n_X$ over $X$. (For what we mean by a locally
closed immersion, see Definition 
\ref{dfn:smooth-etale-proper-finite-rig}.
For the existence of a cover by open subspaces 
with the required property,
see the proof of Proposition \ref{prop:exist-compac}(2) below.)
Let $U\subset P$ be an open subspace containing 
$V$ as a closed subset. Set $Q=Y\times_X P$, $W=V\times_X U$,
$W_1=V\times_X P$ and
$W_2=Y\times_XU$.
Thus, $Q$ is a proper rigid analytic $X$-space, and 
$W$, $W_1$ and $W_2$ are open subspaces of $Q$ 
containing $Y$, via the diagonal embedding 
$Y\to Q$, as a closed subset. 
We have a commutative diagram of immersions with Cartesian squares
$$\xymatrix{V \ar[dr]^-t \ar@{=}[r] \ar@{=}[d] 
& V \ar[dr]^-{t_2} & &\\
V \ar[dr]_-{t_1}& W \ar[r]^-{e_2} \ar[d]^-{e_1} \ar[dr]|-w & W_2 \ar[d]^-{w_2}\\
& W_1 \ar[r]^-{w_1} & Q.}$$
Using Proposition
\ref{prop:loc1}(4), we obtain equivalences
$e_{1,\,\sharp} \circ t_* \simeq t_{1,\,*}$ and
$e_{2,\,\sharp} \circ t_* \simeq t_{2,\,*}$.
Applying this to $w_{1,\,\sharp}$ and $w_{2,\,\sharp}$, 
we obtain equivalences
\begin{equation}
\label{eq-thm:prop-base-1}
w_{1,\,\sharp}\circ t_{1,\,*}\simeq w_{\sharp}\circ t_*\simeq 
w_{2,\,\sharp}\circ t_{2,\,*}.
\end{equation}
Now, consider the commutative diagram with a Cartesian square
$$\xymatrix{V \ar@{=}[rd] \ar[r]^-{t_1} & W_1 
\ar[r]^-{w_1} \ar[d]^-{q'} & Q \ar[d]^-q\\
 & V \ar[r]^-v & Y.\!}$$
By the second step, we deduce equivalences of functors from
$\RigSH_{\Nis}(V;\Lambda)$ to $\RigSH_{\Nis}(Y;\Lambda)$:
$$v_{\sharp}\simeq v_{\sharp}\circ q'_*\circ t_{1,\,*}\simeq 
q_*\circ w_{1,\,\sharp} \circ t_{1,\,*}\simeq 
q_*\circ w_{\sharp} \circ t_*.$$
Thus, it will be enough to show that the functor 
$f_*\circ q_*\circ w_{\sharp}\circ t_*$
takes $\mathcal{H}_V$-equivalences to 
$\mathcal{H}_X$-equivalences. 
Next, consider the commutative diagram with Cartesian squares
$$\xymatrix{V \ar[r]^-{t_2} \ar@{=}[d] 
& W_2 \ar[r]^-{w_2} \ar[d]^-{h'} & Q \ar[d]^-h\\
V \ar[r]^-s & U \ar[r]^-u & P.\!}$$
By the second step, we 
we deduce equivalences of functors from
$\RigSH_{\Nis}(V;\Lambda)$ to $\RigSH_{\Nis}(Y;\Lambda)$:
$$u_{\sharp}\circ s_*\simeq u_{\sharp} \circ h'_*\circ t_{2,\,*}
\simeq h_*\circ w_{2,\,\sharp}\circ t_{2,\,*}\simeq 
h_*\circ w_{\sharp}\circ t_*.$$
Since $p\circ h=f\circ q$ with $p:P\to X$ the structural 
projection of the relative projective space $P$, 
we are left to show that 
$p_*\circ u_{\sharp}\circ s_*$ 
takes $\mathcal{H}_V$-equivalences to 
$\mathcal{H}_X$-equivalences.
This is actually true for each of the functors 
$p_*$, $u_{\sharp}$ and $s_*$.
For the first one, we use the equivalence 
$p_*\simeq p_{\sharp}\circ \Th^{-1}(\Omega_p)$ 
provided by Corollary 
\ref{cor:6ffalg}. For the second one, this is clear, and 
for the third one, this follows from 
Lemma \ref{lem:i-lower-star-sheaf}.
\end{proof}

The following lemma was used in the first step of the proof of 
Theorem \ref{thm:prop-base}.

\begin{lemma}
\label{lem:local-structure-of-rig-etale-morphisms}
Let $f:T\to S$ be an \'etale morphism of rigid analytic spaces.
Then, locally on $S$ and $T$, we may find a commutative triangle
$$\xymatrix{T \ar[r]^-j \ar[dr]_-f & T' \ar[d]^-{f'}\\
& S}$$
where $j$ is an open immersion and $f'$ is a finite \'etale morphism.
\end{lemma}

\begin{proof}
This is a well-known fact. In the generality we are considering here, 
it can be proven by adapting the argument used in proving
Proposition \ref{lem:nearly-etale-is-rig-etale}(3). 
More precisely, it is enough to show that a rig-\'etale morphism 
of formal schemes $f:\mathcal{T}\to \mathcal{S}$ is locally, for the 
rig topology on $\mathcal{S}$ and $\mathcal{T}$, the composition
of an open immersion and a finite rig-\'etale morphism. 
We argue locally around a rigid point
$\mathfrak{s}:\Spf(V) \to \mathcal{S}$
corresponding to $s\in |\mathcal{S}^{\rig}|$.
As in the proof of 
Proposition \ref{lem:nearly-etale-is-rig-etale}(3),
we may assume that the formal scheme  
$\mathfrak{s}\times_{\mathcal{S}}\mathcal{T}/(0)^{\sat}$ 
is the formal spectrum of the $\pi$-adic completion 
of an algebra of the form 
\begin{equation}
\label{eq-lem:local-structure-of-rig-etale-morphisms-1}
V\langle t,s_1,\ldots, s_m\rangle/(R,
\pi^Ns_1-P_1,\ldots, \pi^Ns_m-P_m)^{\sat}[Q^{-1}]
\end{equation}
where $R\in V[t]$ is a monic polynomial 
which is separable over $V[\pi^{-1}]$, and $Q\in V[t,s_1,\ldots, s_m]$.
(The polynomial $R$ is the analogue of the polynomial 
$(t-a_1)\cdots (t-a_r)$ in 
\eqref{eq-lem:nearly-etale-is-rig-etale-7}.
Here, since $V[\pi^{-1}]$ is not algebraically closed, our polynomial 
$R$ will not split in general.)
The remainder of the argument is identical to the one used in 
the proof of Proposition \ref{lem:nearly-etale-is-rig-etale}(3). 
\end{proof}

The following is a corollary of the proof of Theorem 
\ref{thm:prop-base}.

\begin{cor}
\label{cor:f-lower-star-commute-colimits-rigan}
Let $f:Y \to X$ be a proper morphism of rigid analytic spaces. 
Then, the functor 
$$f_*:\RigSH^{(\hyp)}_{\tau}(Y;\Lambda) \to 
\RigSH^{(\hyp)}_{\tau}(X;\Lambda)$$
is colimit-preserving and thus admits a right adjoint.
\end{cor}

\begin{proof}
This is true for the functor  
$$f_*:\RigSH_{\Nis}(Y;\Lambda)\to \RigSH_{\Nis}(X;\Lambda)$$
by Proposition \ref{prop:compact-shv-rigsm}(1).
The result in general follows from the fact that this functor
takes $\mathcal{H}_Y$-equivalences to $\mathcal{H}_X$-equivalences
as shown in the third step of the proof of Theorem
\ref{thm:prop-base}.
\end{proof}

We end this subsection
by establishing the projection formula for direct images
along proper morphisms.

\begin{prop}
\label{prop:proper-projective-formula-weakly-propro}
\ncn{proper projection formula}
$\empty$

\begin{enumerate}

\item[(1)] Let $f:\mathcal{Y} \to \mathcal{X}$ be a proper
morphism of formal schemes. For 
$M\in \FSH^{(\hyp)}_{\tau}(\mathcal{X};\Lambda)$
and $N\in \FSH^{(\hyp)}_{\tau}(\mathcal{Y};\Lambda)$, 
the morphism 
$$M\otimes f_*N \to f_*(f^*M \otimes N)$$
is an equivalence.

\item[(2)] Let $f:Y \to X$ be a proper
morphism of rigid analytic spaces. For 
$M\in \RigSH^{(\hyp)}_{\tau}(X;\Lambda)$
and $N\in \RigSH^{(\hyp)}_{\tau}(Y;\Lambda)$, 
the morphism 
$$M\otimes f_*N \to f_*(f^*M \otimes N)$$
is an equivalence.

\end{enumerate}
\end{prop}

\begin{proof}
We only prove the second part. The proof of the first part 
is similar: in the argument below, 
use Proposition \ref{prop:prop-base-proformal} and 
Lemma \ref{lem:direct-image-in-prl}
instead of Theorem \ref{thm:prop-base},
and Corollary
\ref{cor:f-lower-star-commute-colimits-rigan}.

The functor 
$f_*$ is colimit-preserving by Corollary
\ref{cor:f-lower-star-commute-colimits-rigan}. 
Hence, it is enough
to prove the result when $M$ varies in a set of 
generators under colimits for the $\infty$-category
$\RigSH^{(\hyp)}_{\tau}(X;\Lambda)$.
Thus, we may assume that $M=g_{\sharp}\Lambda$
where $g:X'\to X$ is a smooth morphism.
We form the Cartesian square
$$\xymatrix{Y' \ar[r]^-{g'} \ar[d]^-{f'} & Y \ar[d]^-f\\
X' \ar[r]^-g & X.\!}$$
By Proposition \ref{prop:6f1}(2),
we have natural equivalences
$$M\otimes (-)\simeq g_{\sharp}\circ g^*(-)
\qquad \text{and} \qquad 
(f^*M)\otimes (-) \simeq 
g'_{\sharp} \circ g'^*(-).$$
Modulo these equivalences, the morphism of the statement is the 
composition of
$$g_{\sharp}g^*f_*N\to g_{\sharp}f'_*g'^*N
\to f_*g'_{\sharp}g'^*N.$$
The result follows now from Theorem \ref{thm:prop-base}.
\end{proof}

Recall that an object in a monoidal $\infty$-category $\mathcal{C}^\otimes$ is strongly dualisable if it is so as an object of the homotopy category of $\mathcal{C}$ endowed with 
the induced monoidal structure. The following is
a well-known consequence of the projection 
formula for proper direct images.

\begin{cor}
\label{cor:strong-dual-}
\ncn{duality}
$\empty$

\begin{enumerate}

\item[(1)] Let $f:\mathcal{Y} \to \mathcal{X}$ 
be a smooth and proper morphism of formal schemes. 
Then $f_{\sharp}\Lambda$ 
is strongly dualisable in the monoidal $\infty$-category
$\FSH^{(\hyp)}_{\tau}(\mathcal{X};\Lambda)^{\otimes}$
and its dual is $f_*\Lambda$.

\item[(2)] Let $f:Y \to X$ be a smooth and proper morphism of 
rigid analytic spaces. Then $f_{\sharp}\Lambda$ 
is strongly dualisable in the monoidal $\infty$-category
$\RigSH^{(\hyp)}_{\tau}(X;\Lambda)^{\otimes}$
and its dual is $f_*\Lambda$. 

\end{enumerate}
\end{cor}

\begin{proof}
We only treat the case of rigid analytic motives.
We need to show that there is an equivalence between the
endofunctors $\underline{\Hom}(f_{\sharp}\Lambda,-)$ 
and $(f_*\Lambda) \otimes-$.
We have natural equivalences
$$\underline{\Hom}(f_{\sharp}\Lambda,-)\overset{(1)}{\simeq} 
f_*f^*(-)
\simeq f_*(\Lambda \otimes f^*M)
\overset{(2)}{\simeq} 
f_*(\Lambda)\otimes M$$
where $(1)$ is deduced by adjunction from the 
smooth projection formula 
$f_{\sharp}\Lambda \otimes -\simeq f_{\sharp}\circ f^*(-)$
(see Proposition \ref{prop:6f1}(2)) and (2) is deduced 
from Proposition 
\ref{prop:proper-projective-formula-weakly-propro}(2).
\end{proof}

\subsection{Weak compactifications}

$\empty$

\smallskip

\label{subsect:huber-compact}

In this subsection, we discuss the notion of a weak compactification 
of a rigid analytic $S$-space. For us, it will be enough 
to know that weak compactifications exist locally. 
We will also briefly discuss Huber's 
compactifications.

\begin{dfn}
\label{dfn:weak-compact}
\ncn{rigid analytic spaces!weak compactification of a morphism}
Let $f:Y\to X$ be a morphism of rigid analytic spaces. 
A weak compactification of $f$ is a commutative triangle
\begin{equation}
\label{eq-dfn:weak-compact}
\begin{split}
\xymatrix{Y \ar[r]^i \ar[dr]_-f & W
\ar[d]^-h \\
& X}
\end{split}
\end{equation}
of rigid analytic spaces, where $i$ is a locally closed 
immersion and $h$ a proper morphism.
(See Definition \ref{dfn:smooth-etale-proper-finite-rig}.)
By abuse of language, we say that $h$ is
a weak compactification of $f$ or that $W$ is a 
weak compactification of $Y$.
We define the category of weak compactifications of $f$ 
to be the full subcategory of $(\RigSpc/X)_{f/}$
spanned by the weak compactifications of $f$.
We say that $f$ is weakly compactifiable if it admits a  
weak compactification. (Clearly, for $f$ to be weakly compactifiable, 
it is necessary that $f$ is separated and locally of finite type.)
\end{dfn}

\begin{prop}
\label{prop:exist-compac}
Let $f:Y\to X$ be a morphism of rigid analytic spaces. 
\begin{enumerate}

\item[(1)] The category of weak compactifications of $f$ 
has fiber products and equalizers. 
In particular, when $f:Y \to X$ is weakly compactifiable,
this category is cofiltered.

\item[(2)] Assume that $f$ is locally of finite type. 
Then, locally on $Y$, $f$ is weakly compactifiable.

\end{enumerate}
\end{prop}

\begin{proof}
The first part follows from standard properties of 
proper morphisms and locally closed immersions.
For the second part, since the question is local on $Y$, 
we may assume that $f$ 
factors through an open subspace $U\subset X$ and that $Y\to U$ is the 
generic fiber of a finite type morphism 
$\mathcal{Y}\to \mathcal{U}$ between affine formal schemes. 
In this case, we may factor $f$ as the composition of
$$Y \xrightarrow{s} \B^N_U \xrightarrow{u} \P^N_X \xrightarrow{p} X$$
where $s$ is a closed immersion, $u$ the obvious open immersion 
and $p$ the obvious projection. 
\end{proof}

We will need a short digression concerning 
the notion of relative interior.

\begin{dfn}
\label{dfn:interior-comp}
Let $f:X \to W$ be a morphism between rigid analytic spaces.
Let $V\subset W$ be an open subspace. 
We say that $X$ maps into the interior of $V$ relatively to $W$
and write $f(X)\Subset_W V$ if the closure of $f(|X|)$ in 
$|W|$ is contained in $|V|$.
\symn{$\Subset$}
\end{dfn}

\begin{rmk}
\label{rmk-on-dfn:interior-comp}
Often we use Definition \ref{dfn:interior-comp}
when $f$ is a locally closed immersion. In this case, we write simply 
``$X\Subset_W V$'' instead of ``$f(X)\Subset_W V$''.
\end{rmk}

Below, we use freely the fact that the underlying topological 
space of a rigid analytic space is valuative in the sense of 
\cite[Chapter 0, Definition 2.3.1]{fujiwara-kato}.

\begin{lemma}
\label{lem:closure-locally-closed-subset}
Let $f:X \to W$ be a morphism between 
quasi-compact and quasi-separated rigid analytic spaces.
A point of $|W|$ belongs to $\overline{f(|X|)}$
if and only if its maximal generisation belongs to $f(|X|)$.
Moreover, we have the equalities:
\begin{equation}
\label{eq-lem:closure-locally-closed-subset}
\overline{f(|X|)}=\bigcap_{f(X)\,\Subset_W\, V}|V|=
\bigcap_{f(X)\,\Subset_W\, V}\overline{|V|}.
\end{equation}
\end{lemma}

\begin{proof}
The first assertion follows from \cite[Chapter 0, 
Theorem 2.2.26]{fujiwara-kato}
and the fact that $f(|X|)$ is stable under generisation.
It follows that $\overline{f(|X|)}$ is 
also stable under generisation, which implies that 
$\overline{f(|X|)}$ is the intersection of its open neighbourhoods. 
(Indeed, if $w\in |W|$ does not belong to 
$\overline{f(|X|)}$, we have
$\overline{\{w\}}\cap \overline{f(|X|)}=\emptyset$.)
This gives the first equality in 
\eqref{eq-lem:closure-locally-closed-subset}. 
The second equality follows from 
\cite[Chapter 0, Proposition 2.3.7]{fujiwara-kato}.
\end{proof}

\begin{lemma}
\label{lem:Subset-Subset}
Let $f:X \to W$ be a morphism between 
quasi-compact and quasi-separated rigid analytic spaces. 
Let $V\subset W$ be an open subspace such that 
$f(X)\Subset_W V$. There exists an open subspace 
$V'\subset W$ such that $f(X)\Subset_W V'$ and 
$V'\Subset_W V$.
\end{lemma}

\begin{proof}
By Lemma 
\ref{lem:closure-locally-closed-subset},
we have 
$$\overline{f(|X|)}=\bigcap_{f(X)\Subset_W V'} \overline{|V'|}
\subset |V|.$$
By \cite[Chapter 0, Corollary 2.2.12]{fujiwara-kato}, 
there exists a quasi-compact open subspace $V'\subset W$ 
with $f(X)\Subset_W V'$ such that 
$\overline{|V'|}\subset |V|$ as needed.
\end{proof}

We now discuss Huber's compactifications. We will freely use 
results and notations from Subsection 
\ref{subsect:rel-adic-sp}. We start with a definition.

\begin{dfn}
\label{dfn:universally-stably-uniform}
$\empty$

\begin{enumerate}

\item[(1)] A Tate ring $A$ is said to be universally uniform if 
every finitely generated Tate $A$-algebra is uniform.
(Recall that a finitely generated Tate 
$A$-algebra is a quotient of
$A\langle t\rangle=A_0\langle t\rangle[\pi^{-1}]$
where $t=(t_1,\ldots, t_n)$ is a system of coordinates, $A_0\subset A$
a ring of definition and $\pi\in A$ a topologically nilpotent 
unit contained in $A_0$.)
In particular, a universally uniform Tate ring is also 
stably uniform in the sense of \cite[pages 30--31]{buzz-ver}.
A Tate affinoid ring $R$ is said to be universally uniform
if $R^{\pm}$ is universally uniform.
\ncn{Tate ring!universally uniform}
\ncn{Tate affinoid ring!universally uniform}

\item[(2)] A universally uniform adic space 
is a uniform adic space (as in Definition
\ref{dfn:uniform-adic-space}) which is locally 
isomorphic to $\Spa(A)$, where $A$ is a universally uniform 
Tate affinoid ring.
\ncn{adic spaces!universally uniform}
\end{enumerate} 
\end{dfn}

\begin{nota}
\label{not:adic-finite-type=}
$\empty$

\begin{enumerate}

\item[(1)] Let $S$ be a universally uniform adic space.
We denote by $\Adic/S$ the category of uniform adic 
$S$-spaces. We denote by $\Adic^{\lft}/S$ 
(resp. $\Adic^{\sft}/S$) the full subcategory of $\Adic/S$
spanned by those adic $S$-spaces which are locally of finite type
(resp. which are separated
of finite type).
\symn{$\Adic$}
\symn{$\Adic^{\lft}$}
\symn{$\Adic^{\sft}$}

\item[(2)] Let $S$ be a rigid analytic space. We denote by 
$\RigSpc^{\lft}/S$ (resp. $\RigSpc^{\sft}/S$) 
the full subcategory of $\RigSpc/S$ spanned by those 
rigid analytic $S$-spaces which are locally of finite type 
(resp. which are separated of finite type).
\symn{$\RigSpc^{\lft}$}
\symn{$\RigSpc^{\sft}$}

\item[(3)] Let $S$ be a universally uniform adic space.
By Corollary \ref{cor:from-Huber-to-FK}, $S$ determines 
a rigid analytic space which we denote also by $S$, and we have
equivalences of categories
$\Adic^{\lft}/S\simeq \RigSpc^{\lft}/S$ and 
$\Adic^{\sft}/S\simeq \RigSpc^{\sft}/S$.

\end{enumerate}
\end{nota}

\begin{nota}
\label{nota:Tate-ring-B-c}
Let $A$ be a Tate affinoid ring and $B$ a Tate affinoid $A$-algebra. 
We define a new Tate affinoid $A$-algebra 
$B_{\rm c}=(B_{\rm c}^{\pm},B_{\rm c}^+)$
by setting $B_{\rm c}^{\pm}=B^{\pm}$ and letting $B_{\rm c}^+$ to be the 
integral closure of the subring $A^++B^{\circ\circ}\subset B$.
\symn{$(-)_{\rm c}$}
\end{nota}

The following theorem is due to Huber.

\begin{thm}
\label{thm:huber-compact}
Let $S$ be a quasi-compact and quasi-separated 
universally uniform adic space. 
There is a functor $\Adic^{\rm sft}/S \to 
\Fun({\Delta^1},\Adic/S)$
sending a separated adic $S$-space of finite type $X$ 
to an open immersion $j_X:X \to X^{\rm c}$ over $S$
satisfying the following properties.
\begin{enumerate}

\item[(1)]  Every point of $|X^{\rm c}|$ is a specialisation of 
a point of $|X|$. 
Moreover, for every $x\in |X|$ and 
every valuation ring $V\subset \kappa^+(x)$ containing 
$\kappa^+(s')$ for a specialisation $s'\in |S|$ 
of the image of $x$ in $|S|$, there exists a unique point 
$x'\in |X^{\rm c}|$ which is a specialisation of 
$x$ and such that $\kappa^+(x')=V$.

\item[(2)] The morphism 
$\mathcal{O}_{X^{\rm c}}\to j_{X,\,*}\mathcal{O}_X$
is an isomorphism.

\item[(3)] \emph{(Compatibility with base change)}
If $S' \to S$ is an open immersion, then the morphism
$$j_X\times_S S':X\times_S S'\to X^{\rm c}\times_S S'$$ 
coincides with 
$j_{X'}: X'\to X'^{\rm c}$ where 
$X'$ is the adic $S'$-space $X\times_S S'$. 

\item[(4)] If $S=\Spa(A)$ and $X=\Spa(B)$, then $X^{\rm c}=\Spa(B_{\rm c})$.

\end{enumerate}
\end{thm}

\begin{proof}
This is essentially contained in \cite[Theorem 5.1.5]{huber}.
In loc.~cit., it is assumed that adic spaces satisfy one of the 
conditions in \cite[(1.1.1)]{huber}, but this is only needed 
to insure universal sheafyness. Here, we use instead
universal uniformness and \cite[Theorem 7]{buzz-ver}.
\end{proof}

In the next proposition, we denote a uniform adic space 
and the associated rigid analytic space by the same symbol.
(This is an abuse of notation justified by Corollary 
\ref{cor:from-Huber-to-FK}.)

\begin{prop}
\label{prop:comparison-compactification}
Let $S$ be a quasi-compact and quasi-separated 
universally uniform adic space, and 
let $X$ be a separated adic $S$-space of finite type.
Let $i:X \to W$ be a weak compactification of $X$
over $S$. Then, $X^{\rm c}$ is naturally a weak limit of 
the rigid analytic pro-space $(V)_{X\Subset_W V}$
in the sense of Definition
\ref{dfn:limit-rigid-an-spc}.
\end{prop}

\begin{proof}
By the universal property of Huber's compactifications 
(see \cite[Theorem 5.1.5]{huber}), the locally closed immersion 
$i$ extends to a morphism $i':X^{\rm c}\to W$. Since 
$i'(|X^{\rm c}|)$ is contained in the closure 
of $|X|$ in $|W|$, there is a natural map 
\begin{equation}
\label{eq-prop:comparison-compactification-45}
X^{\rm c}\to (V)_{X\Subset_W V},
\end{equation}
and we need to prove that it exhibits $X^{\rm c}$ as a weak limit of 
$(V)_{X\Subset_W V}$.

We first check that the map 
\begin{equation}
\label{eq-prop:comparison-compactification-47}
|X^{\rm c}|\to \lim_{X\Subset_W V}|V|
\end{equation}
is a bijection. Injectivity is clear since each locally closed immersion
$X\to V$, with $X\Subset_W V$, induces an injection 
$|X^{\rm c}|\to |V^{\rm c}|$ and the map $|X^{\rm c}|\to |V|$
factors this injection. For surjectivity, we use Lemma
\ref{lem:closure-locally-closed-subset} 
which implies that every point $v\in \lim_{X\Subset_W V}|V|$
is a specialisation of a point $x\in |X|$. Thus, we have
$\kappa(v)=\kappa(x)$ and $\kappa^+(s)\subset
\kappa^+(v)\subset \kappa^+(x)$. By Theorem 
\ref{thm:huber-compact}(1), 
the valuation ring $\kappa^+(v)\subset \kappa(x)$ determines 
a point of $|X^{\rm c}|$ which is necessarily sent to $v$ by
\eqref{eq-prop:comparison-compactification-45} since 
$W\to S$ is separated.

It remains to see that for every point
$x$ of $|X^{\rm c}|$ with image $v$ in $\lim_{X\Subset_W V}|V|$
the map $\kappa(v) \to \kappa(x)$ has dense image. 
In fact, we have $\kappa(v)\simeq \kappa(x)$.
To prove this, 
we may assume that $x$ belongs to $|X|$, since the residue field of 
$x$ is equal to the residue field of its maximal generisation
and similarly for $v$.
The claimed result is then clear since $X \to (V)_{X\Subset_W V}$
is given by locally closed immersions.
\end{proof}

\subsection{The exceptional functors, I. Construction}

$\empty$

\smallskip

\label{subsect:exceptional-functors}

In this subsection, we define the exceptional functors 
$f_!$ and $f^!$ associated with a morphism $f$ of rigid analytic 
spaces which is locally of finite type, and establish some of their basic 
properties. We start with the easy case of a locally closed immersion.

\begin{lemma}
\label{defn:j-lower-!-for-locally-closed-immer}
Let $i:Z \to X$ be a locally closed immersion of rigid analytic
spaces. Let $U\subset X$ be an open neighbourhood of $Z$
in which $Z$ is closed. Denote by $s:Z \to U$ and $j:U \to X$
the obvious immersions. Then, the composite functor
$$j_{\sharp}\circ s_*:\RigSH^{(\hyp)}_{\tau}(Z;\Lambda)
\to \RigSH^{(\hyp)}_{\tau}(X;\Lambda)$$
is independent of the choice of $U$ and we denote it by $i_!$.
\end{lemma}

\begin{proof}
Let $U'\subset U$ be an open neighbourhood of $Z$. 
Let $s':Z\to U'$ and $u:U'\to U$ 
be the obvious immersions. We need to show that 
$u_{\sharp}\circ s'_*\simeq s_*$.
To do so, we use the Cartesian square
$$\xymatrix{Z \ar@{=}[r] \ar[d]^-{s'} 
& Z\ar[d]^-s\\
U' \ar[r]^-u & U}$$
and Proposition \ref{prop:loc1}(4).
\end{proof}

\begin{lemma}
\label{lem:composition-i-!-locally-closed}
Let $s:Y \to X$ and $t:Z \to Y$ be locally closed immersions
of rigid analytic spaces. There is an equivalence
$(s\circ t)_!\simeq s_!\circ t_!$
of functors from $\RigSH^{(\hyp)}_{\tau}(Z;\Lambda)$ 
to $\RigSH^{(\hyp)}_{\tau}(X;\Lambda)$.
\end{lemma}

\begin{proof}
Indeed, let $U\subset X$ be an open neighbourhood of 
$Y$ in which $Y$ is closed, and let 
$V\subset U$ be an open neighbourhood of $Z$ in which $Z$ is closed. 
Set $W=Y\cap V$. Consider the commutative diagram with 
a Cartesian square
$$\xymatrix{Z\ar[r]^-e \ar[dr] \ar@/^2pc/[rr]^-t & 
W \ar[r]^-{w} \ar[d]^-d & Y\ar[d]^c \ar[dr]^-s \\
& V \ar[r]^-v & U \ar[r]^-u & X.\!}$$
Using Proposition \ref{prop:loc1}(4), 
we have natural equivalences 
$$u_{\sharp}\circ c_*\circ w_{\sharp}\circ e_*
\simeq u_{\sharp}\circ v_{\sharp}\circ d_*\circ e_*
\simeq (u\circ v)_{\sharp}\circ (d\circ e)_*$$
as needed.
\end{proof}

\begin{prop}
\label{prop:base-change-for-i-!-locally closed}
Consider a Cartesian square of rigid analytic spaces
$$\xymatrix{T\ar[d]^-{f'} \ar[r]^-{i'} & Y \ar[d]^-f\\
Z \ar[r]^-i & X}$$
where $i$ is a locally closed immersion.
\begin{enumerate}

\item[(1)] There is a natural equivalence 
$f^* \circ i_!\simeq  i'_!\circ f'^*$
between functors from 
$\RigSH^{(\hyp)}_{\tau}(Z;\Lambda)$ 
to $\RigSH^{(\hyp)}_{\tau}(Y;\Lambda)$.

\item[(2)] Assume that
$f$ is a proper morphism. 
There is a natural equivalence 
$f_*\circ i'_!\simeq i_!\circ f'_*$
between functors from 
$\RigSH^{(\hyp)}_{\tau}(T;\Lambda)$ 
to $\RigSH^{(\hyp)}_{\tau}(X;\Lambda)$.

\end{enumerate}
\end{prop}

\begin{proof}
Part (1) follows from Proposition \ref{prop:loc1}(4) and
part (2) follows from Theorem \ref{thm:prop-base}.
\end{proof}

We next discuss the case of weakly compactifiable morphisms.

\begin{dfn}
\label{dfn:f-!-up-down}
\ncn{exceptional functors}
Let $f:Y \to X$ be a weakly compactifiable morphism
of rigid analytic spaces. Choose a weak compactification 
$$\xymatrix{Y \ar[r]^-i \ar[dr]_-f
& W \ar[d]^-h \ar[d]\\
& X}$$
of $f$ and define the functor 
$$f_!:\RigSH^{(\hyp)}_{\tau}(Y;\Lambda) \to 
\RigSH^{(\hyp)}_{\tau}(X;\Lambda)$$ 
by setting $f_!=h_*\circ i_!$. 
It follows from Corollary
\ref{cor:f-lower-star-commute-colimits-rigan}
that the functor $f_!$ is colimit-preserving;
we denote by $f^!$ its right
adjoint. The functors $f_!$ and $f^!$ are called the exceptional 
direct and inverse image functors.
\symn{$(-)_{\shriek}, (-)^{\shriek}$}
\end{dfn}

\begin{lemma}
\label{lem:f-!-well-def}
Keep the notations as in Definition 
\ref{dfn:f-!-up-down}. The functor $f_!$ is independent of the choice of 
the weak compactification of $f$.
\end{lemma}

\begin{proof}
Let $i':Y \to W'$ be a second weak compactification of $f$
and denote by $h':W'\to X$ the structural projection.
Without loss of generality, we may assume that $W'$ 
is finer than $W$. Let $U\subset W$ be an open neighbourhood of 
$Y$ in which $Y$ is closed, and let $U'\subset W'$ be the
inverse image of $U$. We then have a commutative diagram 
with a Cartesian square
$$\xymatrix@R=1pc{
& U' 
\ar[r]^-{j'} \ar[dd]^-{g'} & W' \ar[dd]^-g 
\ar[dr]^-{h'} \\
Y \ar[rd]_-s \ar[ru]^-{s'} & & & X.\\
& U 
\ar[r]^-{j} & W \ar[ur]_-h}$$
We need to compare 
$h_*\circ j_{\sharp} \circ s_*$
with 
$h'_*\circ j'_{\sharp} \circ s'_*$.
We have natural transformations
$$h_*\circ j_{\sharp} \circ s_*\simeq 
h_*\circ j_{\sharp} \circ g'_*\circ s'_*
\to h_*\circ g_*\circ j'_{\sharp}\circ s'_*\simeq 
h'_*\circ j'_{\sharp}\circ s'_*$$
where the middle one is an equivalence by 
Theorem \ref{thm:prop-base}.
\end{proof}

\begin{exm}
\label{exm:f!-2-cases}
Using Lemma \ref{lem:f-!-well-def}
and a well-chosen weak compactification, 
we obtain the following particular cases.

\begin{enumerate}

\item[(1)] If $j:U\to X$ is an open immersion, 
then $j_!\simeq j_{\sharp}$ and $j^!\simeq j^*$.

\item[(2)] If $f:Y \to X$ is proper, then $f_!\simeq f_*$.

\end{enumerate}
\end{exm}

\begin{rmk}
\label{rmk:what-we-are-goind-to-do-f-!-}
At this point we have constructed, for each weakly 
compactifiable morphism $f:Y\to X$ of rigid analytic spaces, 
a functor $f_!:\RigSH^{(\hyp)}_{\tau}(Y;\Lambda)\to
\RigSH^{(\hyp)}_{\tau}(X;\Lambda)$. Due to the choice 
of a weak compactification involved in the construction, 
it is not clear why $f\mapsto f_!$ would be functorial in any sense.
The main goal of the remainder of this subsection is to prove that 
in fact it is, as long as we restrict to morphisms between weakly compactifiable rigid analytic spaces over a fixed base. (Note that morphisms between such spaces are automatically weakly compactifiable, so that our construction applies.)
\end{rmk}

\begin{nota}
\label{nota:rigspc-wc-}
Let $S$ be a rigid analytic space. 
\begin{enumerate}

\item[(1)] We denote by $\RigSpc^{\wc}/S\subset \RigSpc/S$ 
the full subcategory of weakly compactifiable rigid analytic 
$S$-spaces. Recall that, by definition, $\RigSpc^{\wc}/S$ is 
contained in $\RigSpc^{\lft}/S$ (see Notation
\ref{not:adic-finite-type=}(2)) and that every object in 
$\RigSpc^{\lft}/S$ is locally isomorphic to an object of 
$\RigSpc^{\wc}/S$ by Proposition 
\ref{prop:exist-compac}.
\symn{$\RigSpc^{\wc}$}

\item[(2)] We denote by $\RigSpc^{\proper}/S\subset \RigSpc/S$
the full subcategory of proper rigid analytic $S$-spaces.
\symn{$\RigSpc^{\proper}$}

\item[(3)] We denote by $\WComp/S$
the category whose objects are pairs $(X,W)$
where $X$ is a rigid analytic 
$S$-space and $W$ a weak compactification of $X$.
There are functors 
$$\mathfrak{d}_S:\WComp/S \to \RigSpc^{\wc}/S
\quad \text{and} \quad
\mathfrak{w}_S:\WComp/S \to \RigSpc^{\proper}/S$$
sending a pair $(X,W)$ to $X$ and $W$
respectively.
\symn{$\WComp$}
\symn{$\mathfrak{d}$}
\symn{$\mathfrak{w}$}

\end{enumerate}
\end{nota}

\begin{prop}
\label{prop:composing-f-!}
Let $S$ be a rigid analytic space. There is a functor 
\begin{equation}
\label{eq-prop:composing-f-!}
\RigSH^{(\hyp)}_{\tau}(-;\Lambda)_!:\RigSpc^{\wc}/S
\to \Prl
\end{equation}
sending an object $X$ to the $\infty$-category
$\RigSH^{(\hyp)}_{\tau}(X;\Lambda)$ and a morphism 
$f$ to the functor $f_!$.
\end{prop}

We fix a rigid analytic space $S$. The functor 
\eqref{eq-prop:composing-f-!}
will be constructed below and the fact that it 
extends the functors in Definition \ref{dfn:f-!-up-down}
is proven in Lemma 
\ref{lem:left-kan-extension-properties-first}.
We start by constructing a similar functor 
defined on $\WComp/S$.

\begin{nota}
\label{nota:subcat-rigsh-!-wcomp}
Given an object $(X,W)$ in $\WComp/S$, we denote by 
$\RigSH^{(\hyp)}_{\tau}((X,W);\Lambda)_!$ the full 
sub-$\infty$-category 
of $\RigSH^{(\hyp)}_{\tau}(W;\Lambda)$ spanned by the 
essential image of the fully faithful embedding
\begin{equation}
\label{eq-nota:subcat-rigsh-!-wcomp-1}
i_!:\RigSH^{(\hyp)}_{\tau}(X;\Lambda)
\to \RigSH^{(\hyp)}_{\tau}(W;\Lambda),
\end{equation}
where $i:X\to W$ is the given locally closed immersion.
\symn{$\RigSH(-)_{\shriek}$}
\end{nota}

\begin{prop}
\label{prop:h-star-preserves-!-}
Let $(f,h):(X',W') \to (X,W)$ be a morphism in 
$\WComp/S$. 

\begin{enumerate}

\item[(1)] The functor
$$h_*:\RigSH^{(\hyp)}_{\tau}(W';\Lambda)\to 
\RigSH^{(\hyp)}_{\tau}(W;\Lambda)$$
takes 
$\RigSH^{(\hyp)}_{\tau}((X',W');\Lambda)_!$
into 
$\RigSH^{(\hyp)}_{\tau}((X,W);\Lambda)_!$, and induces a functor 
\begin{equation}
\label{eq-lem:h-star-preserves-!-1}
(f,h)_!:\RigSH^{(\hyp)}_{\tau}((X',W');\Lambda)_!
\to 
\RigSH^{(\hyp)}_{\tau}((X,W);\Lambda)_!.
\end{equation}

\item[(2)] There is a commutative square
$$\xymatrix{\RigSH^{(\hyp)}_{\tau}(X';\Lambda) \ar[r] 
\ar[d]^-{f_!} & \RigSH^{(\hyp)}_{\tau}((X',W');\Lambda)_!
\ar[d]^-{(f,h)_!}\\
\RigSH^{(\hyp)}_{\tau}(X;\Lambda)\ar[r] & 
\RigSH^{(\hyp)}_{\tau}((X,W);\Lambda)_!}$$
where the horizontal arrows are equivalences.

\item[(3)] If $f$ is an isomorphism, 
then $(f,h)_!$
is an equivalence of $\infty$-categories.

\end{enumerate}
\end{prop}

\begin{proof}
Consider the commutative diagram with a Cartesian square
$$\xymatrix{X' \ar[r]^u \ar[dr]_-f \ar@/^2pc/[rr]^-{i'} 
 & V \ar[r]^-v \ar[d]^-{h'} & W'\ar[d]^-h\\
& X \ar[r]^-i & W.\!}$$
By Lemma \ref{lem:composition-i-!-locally-closed},
we have $i'_!\simeq v_!\circ u_!$. 
Thus, the essential image of $i'_!$ is contained in the essential 
image of $v_!$. On the other hand, by
Proposition \ref{prop:base-change-for-i-!-locally closed}(2), we have 
$h_*\circ v_!\simeq i_!\circ h'_*$. Thus, $h_*$ takes the 
essential image of $v_!$ into the essential image of $i_!$,
which proves the first statement. 

Next, we verify the second statement. 
Note that $V$ is a weak compactification of $X'$ over $X$.
Thus, by Lemma \ref{lem:f-!-well-def}, we have 
$f_!\simeq h'_*\circ u_!$. Using Proposition 
\ref{prop:base-change-for-i-!-locally closed}(2) again, we 
obtain natural equivalences
$$i_!\circ f_!\simeq i_!\circ h'_*\circ u_!\simeq 
h_*\circ v_! \circ u_!\simeq h_*\circ i'_!.$$
This gives the commutative square in the second statement. 
Finally, the third statement follows from the second one
using Lemma \ref{lem:f-!-well-def}.
\end{proof}

\begin{nota}
\label{not:ordinary-functor-invers-direct-star}
We will denote by 
\begin{equation}
\label{eq-prop:composing-f-*-upper}
\RigSH^{(\hyp)}_{\tau}(-;\Lambda)^*:\RigSpc^{\op} \to \Prl
\end{equation}
the functor from Proposition
\ref{prop:basic-functorial-rigsh} 
(in the $\Tate$-stable case and after 
forgetting the monoidal structure) and by 
\begin{equation}
\label{eq-prop:composing-f-*-lower}
\RigSH^{(\hyp)}_{\tau}(-;\Lambda)_*:\RigSpc \to \Prr
\end{equation}
the functor deduced from
\eqref{eq-prop:composing-f-*-upper} 
using the equivalence $(\Prl)^{\op}\simeq \Prr$.
By Corollary 
\ref{cor:f-lower-star-commute-colimits-rigan}, 
the restriction of 
\eqref{eq-prop:composing-f-*-lower}
to $\RigSpc^{\proper}/S$ yields a $\Prl$-valued functor. 
In particular, we have a functor 
\begin{equation}
\label{eq-rmk:composing-mathfrak-d-S-yield-Prl}
\RigSH^{(\hyp)}_{\tau}(\mathfrak{w}_S(-);\Lambda)_*:
\WComp/S \to \Prl.
\end{equation}
\end{nota}

To go further, we need the following well-known general lemma.

\begin{lemma}
\label{lem:fully-faith-sub-functor}
Let $B$ be a simplicial set and let 
$\mathcal{C}:B \to \CAT_{\infty}$ be a diagram of 
$\infty$-categories.
Assume that, for each vertex $x\in B_0$, we are given a full
sub-$\infty$-category $\mathcal{C}'(x)\subset \mathcal{C}(x)$.
Assume also that, for every edge $e\in B_1$, the functor 
$\mathcal{C}(e_0) \to \mathcal{C}(e_1)$ takes 
$\mathcal{C}'(e_0)$ into $\mathcal{C}'(e_1)$.
Then, there exists a diagram 
$\mathcal{C}':B\to \CAT_{\infty}$ and a natural transformation 
$\mathcal{C}'\to \mathcal{C}$ such that for every 
edge $e\in B_1$, $\mathcal{C}'(e)$ is equivalent to the functor 
induced from $\mathcal{C}(e)$ on the sub-$\infty$-categories
$\mathcal{C}'(e_0)$ and $\mathcal{C}'(e_1)$.
\end{lemma}

\begin{proof}
By Lurie's unstraightening \cite[\S 3.2]{lurie}, 
one reduces to prove an analogous statement for coCartesian 
fibrations which is easy and left to the reader.
\end{proof}

By Proposition \ref{prop:h-star-preserves-!-}(1), 
we may apply Lemma \ref{lem:fully-faith-sub-functor}
to the functor 
\eqref{eq-rmk:composing-mathfrak-d-S-yield-Prl}
and the full sub-$\infty$-categories introduced in 
Notation \ref{nota:subcat-rigsh-!-wcomp}.
This yields a functor 
\begin{equation}
\label{eq-prop:composing-f-!-127}
\RigSH^{(\hyp)}_{\tau}((-,-);\Lambda)_!:
\WComp/S\to \Prl.
\end{equation}
(The fact that this functor lands in $\Prl$, and not just in 
$\CAT_\infty$, follows from Corollary 
\ref{cor:f-lower-star-commute-colimits-rigan}
together with Proposition 
\ref{prop:h-star-preserves-!-}(2).)
By left Kan extension 
along the functor $\mathfrak{d}_S$, we obtain from
\eqref{eq-prop:composing-f-!-127}
a functor 
\begin{equation}
\label{eq-prop:composing-f-!-31}
\RigSH^{(\hyp)}_{\tau}(-;\Lambda)_!:
\RigSpc^{\wc}/S\to \Prl.
\end{equation}
The following lemma shows that this left Kan 
extension behaves as we want it to.
\symn{$\RigSH(-)_{\shriek}$}

\begin{lemma}
\label{lem:left-kan-extension-properties-first}
The obvious natural transformation 
\begin{equation}
\label{eq-lem:left-kan-extension-properties-first-1}
\RigSH^{(\hyp)}_{\tau}((-,-);\Lambda)_!\to 
\RigSH^{(\hyp)}_{\tau}(-;\Lambda)_!\circ \mathfrak{d}_S
\end{equation}
is an equivalence. In particular, the functor 
\eqref{eq-prop:composing-f-!-31} sends a morphism
$f:Y \to X$ in $\RigSpc^{\wc}/S$
to the functor $f_!$ of Definition 
\ref{dfn:f-!-up-down}.
\end{lemma}

\begin{proof}
Given an object $X\in \RigSpc^{\wc}/S$, there is an 
equivalence in $\Prl$:
\begin{equation}
\label{eq-lem:left-kan-extension-properties-first-2}
\RigSH^{(\hyp)}_{\tau}(X;\Lambda)_!
\simeq \underset{(Y,W)\in (\WComp/S)_{/X}}{\colim}\;
\RigSH^{(\hyp)}_{\tau}((Y,W);\Lambda)_!
\end{equation}
where the category $(\WComp/S)_{/X}$ consists of pairs 
$(Y,W)$ with $Y$ a rigid analytic $X$-space and 
$W$ a compactification of $Y$ over $S$.
Fix a weak compactification $P$ of $X$ over $S$.
There is an obvious forgetful functor 
$$\alpha:(\WComp/S)_{/(X,P)} \to (\WComp/S)_{/X}$$
admitting a right adjoint $\beta$ given by 
$(Y,W)\mapsto (Y,W\times_SP)$. Moreover, it follows from 
Proposition 
\ref{prop:h-star-preserves-!-}(3), 
that the counit of the adjunction $\alpha\circ \beta\to \id$
induces an equivalence between the functor 
$$\RigSH^{(\hyp)}_{\tau}((-,-);\Lambda)_!:
(\WComp/S)_{/X} \to \Prl$$
and its composition with the endofunctor 
$\alpha\circ \beta$ of 
$(\WComp/S)_{/X}$. Since $\beta$ is right adjoint to $\alpha$, 
composition with $\beta$ is equivalent to left Kan extension 
along $\alpha$. This implies that the colimit in
\eqref{eq-lem:left-kan-extension-properties-first-2} 
is equivalent to 
$$\underset{(Y,W)\in (\WComp/S)_{/(X,P)}}{\colim}\;
\RigSH^{(\hyp)}_{\tau}((Y,W);\Lambda)_!\simeq 
\RigSH^{(\hyp)}_{\tau}((X,P);\Lambda)_!$$
since $(X,P)$ is the final object of 
$(\WComp/S)_{/(X,P)}$. This proves the lemma.
\end{proof}

\begin{cor}
\label{cor:restriction-Open-rigsh-!}
Let $X$ be a weakly compactifiable rigid analytic 
$S$-space, and let $\Op/X$ be the category of 
open subspaces of $X$. Then, the functors 
\begin{equation}
\label{eq-lem:restriction-Open-rigsh-!-1}
\RigSH^{(\hyp)}_{\tau}(-;\Lambda)_!:\Op/X
\to \Prl
\quad\text{and}\quad
\RigSH^{(\hyp)}_{\tau}(-;\Lambda)^*:(\Op/X)^{\op}
\to \Prr
\end{equation}
are exchanged by the equivalence 
$(\Prl)^{\op}\simeq \Prr$.
\symn{$\Op$}
\end{cor}

\begin{proof}
Let $P$ be a weak compactification of $X$. 
Then, for every open subspace $U\subset X$, 
$P$ is also a weak compactification of $U$. 
Thus, we have a functor $\Op/X\to 
\WComp/S$ given by $U\mapsto (U,P)$. 
Therefore, by Lemma \ref{lem:left-kan-extension-properties-first},
the first functor in 
\eqref{eq-lem:restriction-Open-rigsh-!-1}
is equivalent to the functor given by
$U\mapsto \RigSH^{(\hyp)}_{\tau}((U,P);\Lambda)_!$.
It is immediate from the construction of 
\eqref{eq-prop:composing-f-!-127}
that this functor is equivalent to the one sending 
an open immersion $u:U\to X$ to the essential image of the 
fully faithful embedding $u_{\sharp}$. This proves the 
corollary.
\end{proof}

\begin{rmk}
\label{rmk:highlighting-dependence-S-!}
Using the equivalence $(\Prl)^{\op}\simeq \Prr$, 
the functor \eqref{eq-prop:composing-f-!} gives rise to a functor 
\begin{equation}
\label{eq-prop:composing-f-!-upper}
\RigSH^{(\hyp)}_{\tau}(-;\Lambda)^!:(\RigSpc^{\wc}/S)^{\op}
\to \Prr
\end{equation}
sending a morphism $f$ to the functor $f^!$.
\symn{$(-)^{\shriek}$}
\end{rmk}

\begin{prop}
\label{prop:f-up-!-sheaf}
The functor \eqref{eq-prop:composing-f-!-upper}
is a $\Prr$-valued sheaf for the analytic topology.
\end{prop}

\begin{proof}
It is enough to show that, for every 
$X\in \RigSpc^{\wc}/S$, the restriction of
\eqref{eq-prop:composing-f-!-upper}
to $\Op/X$ is a sheaf for the analytic topology. 
This follows from Corollary \ref{cor:restriction-Open-rigsh-!}
and Theorem \ref{thm:hyperdesc}. (Indeed, the inclusion functors
$\Prl\to \CAT_{\infty}$ and $\Prr\to \CAT_{\infty}$ are
limit-preserving by \cite[Proposition 5.5.3.13 \& 
Theorem 5.5.3.18]{lurie}.) 
\end{proof}

\begin{cor}
\label{cor:ext-f-!-all-ft-mor}
There is a unique extension of 
\eqref{eq-prop:composing-f-!}
into a functor 
\begin{equation}
\label{eq-cor:ext-f-!-all-ft-mor}
\RigSH^{(\hyp)}_{\tau}(-;\Lambda)_!:\RigSpc^{\lft}/S \to \Prl
\end{equation}
such that the following condition is satisfied. The functor
\begin{equation}
\label{eq-cor:ext-f-!-all-ft-mor-2}
\RigSH^{(\hyp)}_{\tau}(-;\Lambda)^!:
(\RigSpc^{\lft}/S)^{\op} \to \Prr,
\end{equation}
obtained from \eqref{eq-cor:ext-f-!-all-ft-mor}
using the equivalence $(\Prl)^{\op}\simeq \Prr$,
is a $\Prr$-valued sheaf for the analytic topology. 
\end{cor}

\begin{proof}
This follows from Proposition
\ref{prop:f-up-!-sheaf} using Lemma
\ref{lem:equi-of-sites-infty-topoi}.
Indeed, a $\Prr$-valued $\tau$-sheaf on a site 
$(\mathcal{C},\tau)$
is equivalent to a limit-preserving functor on 
$\Shv_{\tau}(\mathcal{C})^{\op}$; see Definition
\ref{hypersheaves}.
\end{proof}

\begin{rmk}
\label{rmk:do-not-know-rigsh-!-general}
At this point, it is unclear that the $\infty$-category 
$\RigSH^{(\hyp)}_{\tau}(X;\Lambda)_!$ is equivalent to the 
$\infty$-category $\RigSH^{(\hyp)}_{\tau}(X;\Lambda)$ 
for a general object $X\in \RigSpc^{\lft}/S$. This will be 
proven in Subsection \ref{subsect:exceptional-functors-II};
see Corollary \ref{cor:rigsh-!-equiv-rigsh-star}
below. When $X$ is weakly compactifiable, this is already
stated in Proposition \ref{prop:composing-f-!}.
\end{rmk}

We end this subsection with the following result
relating our approach to the one in \cite[\S 5.2]{huber}.

\begin{thm}
\label{thm:f-!-using-huber-compactification}
Let $X$ and $Y$ be quasi-compact and quasi-separated 
uniform adic spaces, and let 
$f:Y \to X$ be a weakly compactifiable morphism of 
rigid analytic spaces. Let $f^{\rm c}:Y^{\rm c}\to X$
be the projection of Huber's compactification of 
$Y$ over $X$, and $j:Y\to Y^{\rm c}$ the obvious inclusion.
Assume one of the following two alternatives.
\begin{enumerate}

\item[(1)] We work in the non-hypercomplete case, and $X$ is 
locally of finite Krull dimension.
When $\tau$ is the \'etale topology, we assume furthermore 
that $\Lambda$ is eventually coconnective.

\item[(2)] We work in the hypercomplete case, and 
$X$ is $(\Lambda,\tau)$-admissible
(see Definition \ref{dfn:lambda-tau-admiss}).

\end{enumerate}
Then, the functor $f_!$ of Definition 
\ref{dfn:f-!-up-down} coincides with the composite functor
$f^{\rm c}_*\circ j_{\sharp}$.
\end{thm}

\begin{proof}
Fix a weak compactification 
$W$ of $Y$ over $X$, and let 
$h:W\to X$ and $i:Y\to W$ be the given morphisms.
The morphism 
$i$ extends to a morphism $i':Y^{\rm c} \to W$.
We have $f^{\rm c}_*\simeq h_* {\circ} i'_*$.
Thus, we only need to show that there is an equivalence
$i'_*\circ j_{\sharp}\simeq i_!$.
The Cartesian square
$$\xymatrix{Y \ar[r]^-j \ar@{=}[d] & Y^{\rm c} \ar[d]^-{i'}\\
Y \ar[r]^-i & W}$$
and Proposition
\ref{prop:base-change-for-i-!-locally closed}(1)
give an equivalence $i'^*\circ i_!\simeq j_!=j_{\sharp}$.
Thus, it is enough to show that the morphism 
$$i_! \to i'_*\circ i'^*\circ i_!$$
is an equivalence. By Proposition 
\ref{prop:comparison-compactification}, 
$Y^{\rm c}$ is the weak limit of the rigid analytic 
pro-space $(V)_{Y\Subset_W V}$.
It follows from Theorem 
\ref{thm:anstC-v2} that there is an equivalence of 
$\infty$-categories
$$\underset{Y\Subset_W V}{\colim}\;
\RigSH^{(\hyp)}_{\tau}(V;\Lambda)
\to \RigSH^{(\hyp)}_{\tau}(Y^{\rm c};\Lambda).$$
Arguing as in the proof of Lemma \ref{lem:colimi-chi-}
(see also Remark \ref{rmk:extension-lemma-colimi-chi-}),
we deduce an equivalence 
$$\underset{Y\Subset_W V}{\colim}\,
r_{V,\,*} \circ r_V^* \simeq i'_*\circ i'^*$$
where, for an open subspace $U\subset W$, $r_U:U\to W$ 
denotes the obvious inclusion.
Therefore, it is enough to prove that 
$$i_!\to r_{V,\,*}\circ r_V^*\circ i_!$$
is an equivalence for every $V\subset W$ such that $Y\Subset_W V$.
Letting $Q$ be the open subspace of $W$ with underlying 
topological space $|Q|=|W|\smallsetminus \overline{|Y|}$, we have 
$W=V\cup Q$.
So it suffices to prove that 
$$r_V^* \circ i_!\to  r_V^*\circ r_{V,\,*}\circ r_V^*\circ i_!
\quad \text{and} \qquad 
r_Q^* \circ i_!\to 
r_Q^*\circ r_{V,\,*}\circ r_V^*\circ i_!$$
are equivalences. For the first one, we use that 
$r_V^*\circ r_{V,\,*}\simeq \id$. For the second one, 
we use Proposition
\ref{prop:base-change-for-i-!-locally closed}(1)
and the fact that $Y\times_WQ=\emptyset$, 
which imply that the source and the target of the natural 
transformation are the zero functor.
\end{proof}

\begin{rmk}
\label{rmk:huber-method-local}
Theorem \ref{thm:f-!-using-huber-compactification}
can be extended to separable morphisms of finite type which are
not assumed to be weakly compactifiable. Indeed, one can 
construct a variant of the functor 
\eqref{eq-prop:composing-f-!}
using Huber's compactifications (instead of weak compactifications)
and show that this new functor coincides with 
\eqref{eq-prop:composing-f-!} on $\RigSpc^{\wc}/S$ and 
gives rise to a sheaf for the analytic topology 
via the equivalence $(\Prl)^{\op}\simeq \Prr$.
We will not pursue this further in this paper, and 
leave it to the interested reader.
\end{rmk}

\subsection{The exceptional functors, II. Exchange}

$\empty$

\smallskip

\label{subsect:exceptional-functors-II}

The goal of this subsection is to prove Theorem 
\ref{thm:exist-functo-exceptional-image} below
and derive a few consequences. This theorem can be seen
as a strengthening of Corollary 
\ref{cor:ext-f-!-all-ft-mor} 
and gives a way to encapsulate the coherence properties
of the exchange equivalences between the ordinary inverse (resp. 
direct) image functors and the exceptional direct (resp. inverse)
image functors. It should be mentioned that Theorem 
\ref{thm:exist-functo-exceptional-image}
is not the best possible statement one could hope for.
For a better statement, we refer to Theorem 
\ref{thm:corresp-6-funct} below
whose proof relies unfortunately on unproven 
claims in \cite{Gait-Rozen-I} concerning 
$(\infty,2)$-categories. However, Theorem 
\ref{thm:exist-functo-exceptional-image}
is probably good enough in practice.

\begin{nota}
\label{not:arrows-categories-rigspc}
Given a simplicial set $B$ and a diagram 
$\mathcal{C}:B \to \CAT_{\infty}$, we
denote by $\int_B\mathcal{C}\to B$ 
a coCartesian fibration classified by $\mathcal{C}$.
When $B$ is an ordinary category and $\mathcal{C}$ takes values
in the sub-$\infty$-category of $\CAT_{\infty}$ 
spanned by ordinary categories, we take for 
$\int_B\mathcal{C}$ the ordinary category given by the Grothendieck
construction. In particular, objects of $\int_B\mathcal{C}$
are represented by pairs $(b,c)$ where $b\in B$ and 
$c\in \mathcal{C}(b)$.
\symn{$\int$}
\end{nota}

\begin{thm}
\label{thm:exist-functo-exceptional-image}
There are functors 
\begin{equation}
\label{eq-thm:exist-functo-exceptional-image-1}
\begin{array}{rcl}
\RigSH^{(\hyp)}_{\tau}(-;\Lambda)_!^* & : & 
{\displaystyle \int_{\RigSpc^{\op}}\RigSpc^{\lft} \to \Prl}\\
& \vspace{-.3cm} & \\
\RigSH^{(\hyp)}_{\tau}(-;\Lambda)_*^! & : &
{\displaystyle
\left(\int_{\RigSpc^{\op}}\RigSpc^{\lft}\right)^{\op} \to \Prr}
\end{array}
\end{equation}
which are exchanged by the equivalence 
$(\Prl)^{\op}\simeq \Prr$ and which admit
the following informal description.
\symn{$\RigSH(-)_*^{\shriek}$}
\symn{$\RigSH(-)^*_{\shriek}$}
\begin{itemize}

\item These functors send an object $(S,X)$, with $S$ a
rigid analytic space and $X$ an object of $\RigSpc^{\lft}/S$,
to the $\infty$-category $\RigSH^{(\hyp)}_{\tau}(X;\Lambda)$.

\item These functors send an arrow $(g,f):(S,Y) \to (T,X)$, 
consisting of morphisms 
$g:T \to S$ and $f:T\times_S Y \to X$, to the functors
$f_!\circ g'^*$ and $g'_*\circ f^!$ respectively,
with $g':T\times_S Y \to Y$ the base change of $g$.

\end{itemize}
Moreover, the functors in
\eqref{eq-thm:exist-functo-exceptional-image-1}
satisfy the following properties.
\begin{enumerate}

\item[(1)] The ordinary functors 
\begin{equation}
\label{eq-thm:exist-functo-exceptional-image-1.6}
\begin{array}{rcl}
\RigSH^{(\hyp)}_{\tau}(-;\Lambda)^* & : &
\RigSpc^{\op}\to \Prl\\
& \vspace{-.3cm} &\\
\RigSH^{(\hyp)}_{\tau}(-;\Lambda)_* & : &
\RigSpc^{\op}\to \Prr
\end{array}
\end{equation}
(as in Notation
\ref{not:ordinary-functor-invers-direct-star}) 
are obtained from the functors in 
\eqref{eq-thm:exist-functo-exceptional-image-1}
by composition with the diagonal functor $\RigSpc^{\op}\to 
\int_{\RigSpc^{\op}}\RigSpc^{\lft}$, given by 
$S\mapsto (S,S)$.

\item[(2)] For a rigid analytic space $S$, the functors
\begin{equation}
\label{eq-thm:exist-functo-exceptional-image-1.81}
\begin{array}{rcl}
\RigSH^{(\hyp)}_{\tau}(-;\Lambda)_! & : & 
\RigSpc^{\lft}/S \to \Prl \\
& \vspace{-.3cm} & \\
\RigSH^{(\hyp)}_{\tau}(-;\Lambda)^! & : & 
\RigSpc^{\lft}/S \to \Prr
\end{array}
\end{equation}
(as in Corollary \ref{cor:ext-f-!-all-ft-mor})
are obtained from the functors in
\eqref{eq-thm:exist-functo-exceptional-image-1}
by restriction to $\RigSpc^{\lft}/S$.

\end{enumerate}
\end{thm}

To construct the functors in
\eqref{eq-thm:exist-functo-exceptional-image-1},
we start with the functor
\begin{equation}
\label{eq-thm:exist-functo-exceptional-image-3}
\RigSH^{(\hyp)}_{\tau}(-;\Lambda)^{*,*}:
\int_{\RigSpc^{\op}} (\RigSpc^{\proper})^{\op} \to \Prl
\end{equation}
admitting the following informal description.
\begin{itemize}

\item It sends a pair $(S,X)$, with $S$ a
rigid analytic space and $X$ an object of $\RigSpc^{\proper}/S$,
to the $\infty$-category $\RigSH^{(\hyp)}_{\tau}(X;\Lambda)$.

\item It sends an arrow $(g,f):(S,X) \to (T,Y)$, consisting of morphisms
$g:T \to S$ and $f:Y \to T\times_S X$, to the functor
$f^*\circ g'^*$ with $g':T\times_S X \to X$ the base change of $g$.

\end{itemize}
Said differently, 
\eqref{eq-thm:exist-functo-exceptional-image-3}
is the composition of 
$$\int_{\RigSpc^{\op}}(\RigSpc^{\proper})^{\op}
\to \RigSpc^{\op} 
\xrightarrow{\RigSH^{(\hyp)}_{\tau}(-;\,\Lambda)^*}
\Prl$$
where the first functor is given by $(S,X)\mapsto X$.
We will apply to the functor
\eqref{eq-thm:exist-functo-exceptional-image-3}
the following general construction.

\begin{cons}
\label{cons:adjointable-cocart-fibration}
Let $B$ be a simplicial set, $p:\mathcal{E} \to B$ a 
coCartesian fibration and 
$\mathfrak{D}:\mathcal{E} \to \CAT_{\infty}$ a functor.
We assume the following condition.
\begin{itemize}

\item[($\star$)] For every commutative square
$$\xymatrix{X \ar[r]^f \ar[d]^-g & Y\ar[d]^-{g'}\\
X' \ar[r]^-{f'} & Y'}$$
in $\mathcal{E}$, such that $g$ and $g'$ are 
$p$-coCartesian, and $p(f)$ and $p(f')$ are identity morphisms,
the associated square
$$\xymatrix{\mathfrak{D}(X) \ar[r] \ar[d] & \mathfrak{D}(Y) \ar[d]\\
\mathfrak{D}(X') \ar[r] & \mathfrak{D}(Y')}$$
is right adjointable.

\end{itemize}
Let $p':\mathcal{E}'\to B$ be a coCartesian fibration 
which is opposite to $\mathcal{E}$, i.e., if $p$ is classified 
by a diagram $\mathcal{C}:B\to \CAT_{\infty}$, then 
$p'$ is classified by the diagram $\mathcal{C}^{\op}:B\to 
\CAT_{\infty}$ obtained by composing $\mathcal{C}$ 
with the autoequivalence $(-)^{\op}$ of $\CAT_{\infty}$. 
In particular, for $b\in B$, the fiber
$\mathcal{E}'_b$ of $p'$ at $b$ is equivalent to the opposite
of the fiber $\mathcal{E}_b$ of $p$ at $b$. Similarly, 
given a $p$-coCartesian edge $A \to B$ in $\mathcal{E}$, 
there is an associated $p'$-coCartesian edge $A'\to B'$ in 
$\mathcal{E}'$ such that $A'$ and $B'$ are the images of 
$A$ and $B$ by the equivalences between the fibers of 
$p$ and the opposite of the fibers of $p'$.

Then, there exists a diagram 
$\mathfrak{D}':\mathcal{E}'\to \CAT_{\infty}$ 
which admits the following informal description.
\begin{enumerate}

\item[(1)] For $b\in B$, the functor
$\mathfrak{D}'|_{\mathcal{E}'_b}:\mathcal{E}'_b\to 
\CAT_{\infty}$ lands in $\CAT_{\infty}^{\Rder}$ and it 
is deduced from the functor $\mathfrak{D}|_{\mathcal{E}_b}:
\mathcal{E}_b\to \CAT_{\infty}^{\Lder}$ using the equivalences 
$\mathcal{E}'_b\simeq (\mathcal{E}_b)^{\op}$ and 
$\CAT_{\infty}^{\Rder}\simeq (\CAT_{\infty}^{\Lder})^{\op}$.

\item[(2)] Given a $p$-coCartesian edge $A \to B$ in $\mathcal{E}$
with corresponding $p'$-coCartesian edge $A'\to B'$, 
the associated functor $\mathfrak{D}(A) \to \mathfrak{D}(B)$
is equivalent to the functor 
$\mathfrak{D}'(A') \to \mathfrak{D}'(B')$.

\end{enumerate}
The diagram $\mathfrak{D}'$ is constructed as follows.
Consider the coCartesian fibration 
$q:\mathcal{F} \to \mathcal{E}$ classified by $\mathfrak{D}$. 
By \cite[Proposition 2.4.2.3(3)]{lurie}, 
$p\circ q:\mathcal{F} \to B$ is a coCartesian fibration
and $q$ sends a $p\circ q$-coCartesian edge to a $p$-coCartesian edge.
Applying straightening to $p\circ q$ and $p$, we obtain 
a morphism $\phi:\mathfrak{N}\to \mathcal{C}$ in $\Fun(B,\CAT_{\infty})$
between the diagrams $\mathfrak{N}:B\to \CAT_{\infty}$ 
and $\mathcal{C}:B\to \CAT_{\infty}$ classifying 
$p\circ q$ and $q$ respectively. Note that for $b\in B$, 
the functor $\phi(b):\mathfrak{N}(b) \to \mathcal{C}(b)$ 
is equivalent to the functor $q_b:\mathcal{F}_b\to \mathcal{E}_b$ 
induced on the fibers of 
$p\circ q$ and $p$. Hence, $\phi(b)$ is a coCartesian fibration.
Condition ($\star$) is equivalent to the following one.
\begin{itemize}

\item[($\star'$)] For every $b\in B$, the coCartesian fibration 
$\phi(b):\mathfrak{N}(b) \to \mathcal{C}(b)$ 
is also a Cartesian fibration
and, for every edge $b_0 \to b_1$ in $B$, 
the associated commutative square 
$$\xymatrix{\mathfrak{N}(b_0) \ar[r] \ar[d]^-{\phi(b_0)} & 
\mathfrak{N}(b_1) \ar[d]^-{\phi(b_1)}\\
\mathcal{C}(b_0) \ar[r] & \mathcal{C}(b_1)}$$
is such that the functor $\mathfrak{N}(b_0)\to\mathfrak{N}(b_1)$
takes a $\phi(b_0)$-Cartesian edge to a $\phi(b_1)$-Cartesian edge.

\end{itemize}
Passing to the opposite $\infty$-categories, 
condition ($\star'$) says that the natural transformation 
$\phi^{\op}:\mathfrak{N}^{\op}\to \mathcal{C}^{\op}$
sends a vertex $b\in B$ to a coCartesian fibration and
an edge of $B$ to a functor preserving coCartesian edges.
Applying unstraightening to $\phi^{\op}$, we obtain 
a commutative triangle 
$$\xymatrix{\mathcal{F}' \ar[rr]^-{q'} \ar[dr]_-{p'\circ q'} 
& & \mathcal{E}'\ar[dl]^-{p'}\\
& B &}$$
where $p'$ and $p'\circ q'$ are the 
coCartesian fibrations classified by $\mathcal{C}^{\op}$ and 
$\mathfrak{N}^{\op}$. We may assume that 
$q'$ is a fibration for the coCartesian model structure on 
$(\Set^+_{\Delta})_{/B}$
(see \cite[Proposition 3.1.3.7]{lurie}) which insures that 
$q'$ is an inner fibration (by using \cite[Remark 3.1.3.4]{lurie}). 
In this case, 
$q'$ is also a coCartesian 
fibration. To prove this, we argue as in the proof of Lemma 
\ref{lem:cocart-fibration-M-Xi}. More precisely, by 
\cite[Proposition 2.4.2.11]{lurie}, we know that 
$q'$ is a locally coCartesian fibration and, by 
\cite[Proposition 2.4.2.8]{lurie}, it remains to check that 
locally $q'$-coCartesian edges can be composed. 
This follows from the characterisation of locally 
$q'$-coCartesian edges given in 
\cite[Proposition 2.4.2.11]{lurie}
and condition ($\star'$).
That said, the announced diagram $\mathfrak{D}': 
\mathcal{E}'\to \CAT_{\infty}$ is the one obtained from $q'$
by straightening and composing with the autoequivalence 
$(-)^{\op}$ of $\CAT_{\infty}$.
\end{cons}

\begin{rmk}
\label{rmk:adjointable-cocart-fibration-cocart-section}
Continuing with the notation and assumptions of 
Construction \ref{cons:adjointable-cocart-fibration}, 
let $s:B\to \mathcal{E}$ be a coCartesian section.
This corresponds, by straightening, to a natural transformation
from the constant diagram $\{*\}:B \to \CAT_{\infty}$
to $\mathcal{C}$. Passing to opposite functors and unstraightening, 
we obtain another coCartesian section $s':B\to\mathcal{E}'$.
It follows from the construction that the two composites
$$B\xrightarrow{s}\mathcal{E}\xrightarrow{\mathfrak{D}}
\CAT_{\infty} \qquad\text{and}\qquad 
B\xrightarrow{s'}\mathcal{E}'\xrightarrow{\mathfrak{D}'}
\CAT_{\infty}$$
are the same.
\end{rmk}

\begin{lemma}
\label{lem:adjointable-cocart-fibration-applies}
The condition ($\star$) in Construction 
\ref{cons:adjointable-cocart-fibration} is satisfied for $p$ the 
coCartesian fibration
$\int_{\RigSpc^{\op}}(\RigSpc^{\proper})^{\op}\to \RigSpc^{\op}$, 
given by $(S,X)\mapsto S$, and $\mathfrak{D}$ the functor 
\eqref{eq-thm:exist-functo-exceptional-image-3} composed with 
the inclusion $\Prl\to\CAT_{\infty}$.
\end{lemma}

\begin{proof}
A commutative square as in condition ($\star$) corresponds to 
a square of the form
$$\xymatrix@C=3pc{(S,X) \ar[r]^-{(\id_S,\,f)} 
\ar[d]_-{(g,\,\id_{X'})} & 
(S,Y)\ar[d]^-{(g,\,\id_{Y'})}\\
(T,X') \ar[r]^{(\id_T,\,f')} & (T,Y'),}$$
where $g:T\to S$ is a morphism of rigid analytic spaces, 
$f:Y\to X$ a morphism in $\RigSpc^{\proper}/S$, 
and $f':Y'\to X'$ 
the base change of $f$ along $g$. 
Letting $g':X'\to X$ and $g'':Y'\to Y$ 
be the base changes of $g$, the functor 
\eqref{eq-thm:exist-functo-exceptional-image-3}
takes the above square to the commutative square of 
$\infty$-categories
$$\xymatrix{\RigSH^{(\hyp)}_{\tau}(X;\Lambda) \ar[r]^-{f^*} 
\ar[d]^-{g'^*} & 
\RigSH^{(\hyp)}_{\tau}(Y;\Lambda) \ar[d]^-{g''^*}\\
\RigSH^{(\hyp)}_{\tau}(X';\Lambda) \ar[r]^-{f'^*} & 
\RigSH^{(\hyp)}_{\tau}(Y';\Lambda).\!}$$
The morphism $f$, 
being a morphism of proper rigid analytic $S$-spaces,
is proper.
Thus, the right adjointability of the above square follows from 
Theorem \ref{thm:prop-base}(1).
\end{proof}

By Lemma \ref{lem:adjointable-cocart-fibration-applies},
we may use Construction 
\ref{cons:adjointable-cocart-fibration}
to obtain a functor
\begin{equation}
\label{eq-thm:exist-functo-exceptional-image-5}
\RigSH^{(\hyp)}_{\tau}(-;\Lambda)^*_*:
\int_{\RigSpc^{\op}}\RigSpc^{\proper} \to \Prl.
\end{equation}
More precisely, the composition of 
\eqref{eq-thm:exist-functo-exceptional-image-5}
with $\Prl\to \CAT_{\infty}$
is the functor $\mathfrak{D}'$ when we take for 
$\mathfrak{D}$ the composition of
\eqref{eq-thm:exist-functo-exceptional-image-3}
with $\Prl\to\CAT_{\infty}$; that the 
resulting functor $\mathfrak{D}'$ lands in $\Prl$ follows from 
Corollary \ref{cor:f-lower-star-commute-colimits-rigan}.
The functor 
\eqref{eq-thm:exist-functo-exceptional-image-5} 
admits the following informal description.
\begin{itemize}

\item It sends a pair $(S,X)$, with $S$ a
rigid analytic space and $X$ an object of $\RigSpc^{\proper}/S$,
to the $\infty$-category $\RigSH^{(\hyp)}_{\tau}(X;\Lambda)$.

\item It sends an arrow $(g,f):(S,Y) \to (T,X)$, consisting of morphisms 
$g:T \to S$ and $f:T\times_S Y \to X$, to the functor
$f_*\circ g'^*$ with $g':T\times_S Y \to Y$ the base change of $g$.

\end{itemize}
Integrating the functors $\mathfrak{w}_S$ 
from Notation \ref{nota:rigspc-wc-}, 
we obtain a functor 
\begin{equation}
\label{eq-thm:exist-functo-exceptional-image-frak-w}
\mathfrak{w}:\int_{\RigSpc^{\op}}\WComp
\to \int_{\RigSpc^{\op}}\RigSpc^{\proper}.
\end{equation}
Composing with 
\eqref{eq-thm:exist-functo-exceptional-image-5},
we obtain a functor 
\begin{equation}
\label{eq-thm:exist-functo-exceptional-image-7}
\RigSH^{(\hyp)}_{\tau}(\mathfrak{w}(-);\Lambda)^*_*:
\int_{\RigSpc^{\op}}\WComp \to \Prl.
\end{equation}

\begin{nota}
\label{nota:subcat-rigsh-star-!}
Given $(S,(X,W))\in \int_{\RigSpc^{\op}}
\WComp$, we denote by 
$\RigSH^{(\hyp)}_{\tau}((X,W);\Lambda)^*_!$ the full 
sub-$\infty$-category 
of $\RigSH^{(\hyp)}_{\tau}(W;\Lambda)^*_*$ 
introduced in Notation 
\ref{nota:subcat-rigsh-!-wcomp}, i.e., 
the essential image of the fully faithful embedding
\eqref{eq-nota:subcat-rigsh-!-wcomp-1}.
\symn{$\RigSH(-)^*_{\shriek}$}
\end{nota}

The next statement is a strengthening of Proposition
\ref{prop:h-star-preserves-!-}(1).

\begin{prop}
\label{prop:subcat-rigsh-star-!}
Given an arrow $(g,(f,h)):
(S,(Y,Q))\to (T,(X,P))$
in $\int_{\RigSpc^{\op}} \WComp$, 
the associated functor
\begin{equation}
\label{eq-prop:subcat-rigsh-star-!-1}
\RigSH^{(\hyp)}_{\tau}(Q;\Lambda)^*_*
\to 
\RigSH^{(\hyp)}_{\tau}(P;\Lambda)^*_*
\end{equation}
takes $\RigSH^{(\hyp)}_{\tau}((Y,Q);\Lambda)^*_!$
into
$\RigSH^{(\hyp)}_{\tau}((X,P);\Lambda)^*_!$
and induces a functor 
\begin{equation}
\label{eq-prop:subcat-rigsh-star-!-2}
\RigSH^{(\hyp)}_{\tau}((Y,Q);\Lambda)^*_!
\to 
\RigSH^{(\hyp)}_{\tau}((X,P);\Lambda)^*_!.
\end{equation}
\end{prop}

\begin{proof}
Using Proposition
\ref{prop:h-star-preserves-!-}(1), 
we only need to treat the case of a morphism of the form
$$(g,\id,\id):(S,(Y,Q))\to (T,T\times_S Y,T\times_S Q).$$
In this case, we need to show that the functor 
$$g'^*:\RigSH^{(\hyp)}_{\tau}(Q;\Lambda)
\to \RigSH^{(\hyp)}_{\tau}(T\times_SQ;\Lambda),$$
with $g':T\times_S Q\to Q$ the base change of $g$, 
sends the essential image of $i_!$, with $i:Y\to Q$ the given 
immersion, to the essential image of $i'_!$, with 
$i':T\times_S Y \to T\times_S Q$ the base change of $i$. 
This follows immediately from 
Proposition \ref{prop:base-change-for-i-!-locally closed}(1).
\end{proof}

Combining Proposition
\ref{prop:subcat-rigsh-star-!} with Lemma
\ref{lem:fully-faith-sub-functor}, 
we deduce a functor
\begin{equation}
\label{eq-thm:exist-functo-exceptional-image-11}
\RigSH^{(\hyp)}_{\tau}((-,-);\Lambda)^*_!:
\int_{\RigSpc^{\op}}\WComp \to \Prl,
\end{equation}
and this functor restricts to 
\eqref{eq-prop:composing-f-!-127}
on $\WComp/S$ for every rigid analytic space $S$.
Integrating the functors $\mathfrak{d}_S$ 
from Notation \ref{nota:rigspc-wc-}, 
we obtain a functor  
\begin{equation}
\label{eq-thm:exist-functo-exceptional-image-frak-d}
\mathfrak{d}:\int_{\RigSpc^{\op}}\WComp
\to \int_{\RigSpc^{\op}}\RigSpc^{\wc}
\end{equation}
given by $(S,(X,W))\mapsto (S,X)$. 
By left Kan extension along the functor 
\eqref{eq-thm:exist-functo-exceptional-image-frak-d},
we obtain from  
\eqref{eq-thm:exist-functo-exceptional-image-11} 
a functor 
\begin{equation}
\label{eq-thm:exist-functo-exceptional-image-13}
\RigSH^{(\hyp)}_{\tau}(-;\Lambda)^*_!:
\int_{\RigSpc^{\op}}\RigSpc^{\wc} \to \Prl.
\end{equation}
We gather a few properties satisfied by this functor in
the following lemma.

\begin{lemma}
\label{lem:left-kan-extension-properties}
$\empty$

\begin{enumerate}

\item[(1)] The obvious natural transformation 
$$\RigSH^{(\hyp)}_{\tau}((-,-);\Lambda)^*_!\to 
\RigSH^{(\hyp)}_{\tau}(-;\Lambda)^*_!\circ \mathfrak{d}$$ 
is an equivalence.

\item[(2)] Composing 
\eqref{eq-thm:exist-functo-exceptional-image-13}
with the diagonal functor 
$$\RigSpc^{\op}\to 
\int_{\RigSpc^{\op}}\RigSpc^{\wc}$$ 
yields the ordinary functor 
$\RigSH^{(\hyp)}_{\tau}(-;\Lambda)^*:\RigSpc^{\op}
\to \Prl$.

\item[(3)] For a rigid analytic space $S$, the restriction of 
\eqref{eq-thm:exist-functo-exceptional-image-13}
to $\RigSpc^{\wc}/S$ is equivalent to the functor 
$\RigSH^{(\hyp)}_{\tau}(-;\Lambda)_!:\RigSpc^{\wc}/S
\to \Prl$
of Proposition \ref{prop:composing-f-!}.

\end{enumerate}
\end{lemma}

\begin{proof}
The third assertion follows from 
\cite[Proposition 4.3.3.10]{lurie}. 
Using this and Lemma \ref{lem:left-kan-extension-properties-first},
we deduce the first assertion. 
For the second assertion we argue as follows.
By the first assertion, it suffices to describe the composition of 
\eqref{eq-thm:exist-functo-exceptional-image-11}
with the diagonal functor
$$\RigSpc^{\op}\to 
\int_{\RigSpc^{\op}}\WComp$$
given by $S\mapsto (S,(S,S))$.
In this composition, we may replace
\eqref{eq-thm:exist-functo-exceptional-image-11} 
by
\eqref{eq-thm:exist-functo-exceptional-image-7} 
without changing the result.
In other words, our functor is the composition of
$$\RigSpc^{\op}\xrightarrow{\Delta}
\int_{\RigSpc^{\op}}\RigSpc^{\proper}\xrightarrow{\RigSH^{(\hyp)}_{\tau}(-;\,\Lambda)^*_*}\Prl$$
where $\Delta$ is the diagonal functor given by $S\mapsto (S,S)$.
Since $\Delta$ is a coCartesian section, the result follows from
Remark \ref{rmk:adjointable-cocart-fibration-cocart-section}.
\end{proof}

For later use, we also record the following fact.

\begin{lemma}
\label{lem:rigsh-!-*-inverse-image}
Let $S$ be a rigid analytic space, and let $X\in \RigSpc^{\wc}/S$.
Then, the composition of 
\eqref{eq-thm:exist-functo-exceptional-image-13}
with the functor 
$$(\RigSpc/S)^{\op} \to 
\int_{\RigSpc^{\op}}\RigSpc^{\wc},$$
given by $T\mapsto (T,T\times_S X)$,
is equivalent to the functor
$$\RigSH^{(\hyp)}_{\tau}(-\times_S X;\Lambda)^*:
(\RigSpc/S)^{\op}\to \Prl.$$
\end{lemma}

\begin{proof}
We first reduce to the case where the rigid analytic 
$S$-space $X$ is proper. To do so, we fix a weak 
compactification $W$ of $X$, and consider the functors 
$$\Delta_X:(\RigSpc/S)^{\op} \to 
\int_{\RigSpc^{\op}}\RigSpc^{\wc}
\quad \text{and} \quad
\Delta_W:(\RigSpc/S)^{\op} \to 
\int_{\RigSpc^{\op}}\RigSpc^{\wc}$$
given by $T\mapsto (T,T\times_SX)$ and 
$T\mapsto (T,T\times_SW)$ respectively. 
The given immersion $i:X\to W$ 
induces a natural transformation 
$\mathfrak{i}:\Delta_X \to \Delta_W$. 
Applying \eqref{eq-thm:exist-functo-exceptional-image-13},
we obtain a natural transformation 
$$\mathfrak{i}_!:
\RigSH^{(\hyp)}_{\tau}(-;\Lambda)^*_!\circ \Delta_X \to 
\RigSH^{(\hyp)}_{\tau}(-;\Lambda)^*_!\circ \Delta_W.$$
On $T\in \RigSpc/S$, the natural transformation 
$\mathfrak{i}_!$ 
is given by the fully faithful embedding 
$(T\times_S i)_!$. It follows that 
$\RigSH^{(\hyp)}_{\tau}(-;\Lambda)^*_!\circ \Delta_X$
can be obtained from 
$\RigSH^{(\hyp)}_{\tau}(-;\Lambda)^*_!\circ \Delta_W$
by applying Lemma \ref{lem:fully-faith-sub-functor}
to the essential images of the functors 
$(T\times_S i)_!$, for $T\in \RigSpc/S$. 
Using Proposition \ref{prop:base-change-for-i-!-locally closed}(1), 
we see that it is enough to prove that 
$\RigSH^{(\hyp)}_{\tau}(-;\Lambda)^*_!\circ \Delta_W$
is given by $\RigSH^{(\hyp)}_{\tau}(-\times_S W;\Lambda)^*$.
Said differently, we may assume that $X$ is proper over $S$.

We now prove the lemma assuming that 
$X$ is proper over $S$. (The argument is the same 
as the one used for the proof of Lemma
\ref{lem:left-kan-extension-properties}(2).)
By Lemma \ref{lem:left-kan-extension-properties}(1),
it is enough to prove the same conclusion for the composition of 
\eqref{eq-thm:exist-functo-exceptional-image-11}
with the functor 
$$\Delta'_X:(\RigSpc/S)^{\op} \to 
\int_{\RigSpc^{\op}}\WComp,$$
given by $T\mapsto (T,(T\times_S X,T\times_S X))$. 
In this composition, we may replace
\eqref{eq-thm:exist-functo-exceptional-image-11} 
by
\eqref{eq-thm:exist-functo-exceptional-image-7} 
without changing the result.
Since $\Delta'_X$ is a coCartesian section, the result follows from
Remark \ref{rmk:adjointable-cocart-fibration-cocart-section}.
\end{proof}

By Lemmas 
\ref{lem:left-kan-extension-properties}
and 
\ref{lem:rigsh-!-*-inverse-image},
the functor \eqref{eq-thm:exist-functo-exceptional-image-13}
admits the following informal description.
\begin{itemize}

\item It sends an object $(S,X)$, with $S$ a 
rigid analytic space and $X$ an object of $\RigSpc^{\wc}/S$,
to the $\infty$-category $\RigSH^{(\hyp)}_{\tau}(X;\Lambda)$.

\item It sends an arrow $(S,Y) \to (T,X)$, consisting of morphisms 
$g:T \to S$ and $f:T\times_S Y \to X$, to the functor
$f_!\circ g'^*$ with $g':T\times_S Y \to Y$ the base change of $g$.

\end{itemize}
Finally, we define the functor 
\begin{equation}
\label{eq-thm:exist-functo-exceptional-image-37}
\RigSH^{(\hyp)}_{\tau}(-;\Lambda)^*_!:
\int_{\RigSpc^{\op}}\RigSpc^{\lft} \to \Prl
\end{equation}
to be the left Kan extension of 
\eqref{eq-thm:exist-functo-exceptional-image-13}
along the fully faithful inclusion 
\begin{equation}
\label{eq-thm:exist-functo-exceptional-image-iota}
\iota:\int_{\RigSpc^{\op}}\RigSpc^{\wc}
\to 
\int_{\RigSpc^{\op}}\RigSpc^{\lft}.
\end{equation}
Note that the functor 
\eqref{eq-thm:exist-functo-exceptional-image-37}
is an extension of 
\eqref{eq-thm:exist-functo-exceptional-image-13}
in the usual sense, i.e., the restriction of 
\eqref{eq-thm:exist-functo-exceptional-image-37} 
along $\iota$ is indeed the functor 
\eqref{eq-thm:exist-functo-exceptional-image-13}.

\begin{prop}
\label{prop:restriction-rigsh-lft-!-star}
For a rigid analytic 
space $S$, the restriction of 
\eqref{eq-thm:exist-functo-exceptional-image-37}
to $\RigSpc^{\lft}/S$ is equivalent to the functor 
$$\RigSH^{(\hyp)}_{\tau}(-;\Lambda)_!:\RigSpc^{\lft}/S \to \Prl$$ 
of Corollary \ref{cor:ext-f-!-all-ft-mor}.
\end{prop}

\begin{proof}
By \cite[Proposition 4.3.3.10]{lurie}, 
it is enough to show that the functor 
$\RigSH^{(\hyp)}_{\tau}(-;\Lambda)_!$
in Corollary \ref{cor:ext-f-!-all-ft-mor}
is a left Kan extension of the same-named functor in Proposition 
\ref{prop:composing-f-!}. Using the equivalence
$(\Prl)^{\op}\simeq \Prr$, it is equivalent to show that 
the functor 
$\RigSH^{(\hyp)}_{\tau}(-;\Lambda)^!$
in Corollary \ref{cor:ext-f-!-all-ft-mor}
is the right Kan extension of the same-named functor 
in Remark \ref{rmk:highlighting-dependence-S-!}.
Since the former was defined as the unique 
$\Prr$-valued sheaf for the analytic topology 
extending the latter, the 
result follows from Lemma
\ref{lemma:sheaf-extension-vs-right-Kan-exte} 
below.
\end{proof}

\begin{lemma}
\label{lemma:sheaf-extension-vs-right-Kan-exte}
Let $(\mathcal{C}',\tau')$ be a site with $\mathcal{C}'$ 
an ordinary category admitting finite limits. Let 
$\mathcal{C}\subset \mathcal{C}'$ be a full subcategory 
closed under finite limits and let $\tau$ be the induced 
topology on $\mathcal{C}$. Assume that the morphism of 
sites $(\mathcal{C}',\tau')\to (\mathcal{C},\tau)$ 
induces an equivalence between the associated ordinary topoi. 
(Equivalently, every object of $\mathcal{C}'$ admits 
a cover by objects in $\mathcal{C}$.)
Let $\mathcal{D}$ be an $\infty$-category admitting limits 
and let $F:\mathcal{C}^{\op}\to \mathcal{D}$ be a 
$\mathcal{D}$-valued $\tau$-sheaf on $\mathcal{C}$. 
Then, the right Kan extension $F':\mathcal{C}'^{\op}\to 
\mathcal{D}$ of $F$ along the inclusion 
$\mathcal{C}^{\op}\to \mathcal{C}'^{\op}$
is a $\tau'$-sheaf. More precisely, 
$F'$ is the image of $F$ by the equivalence of 
$\infty$-categories
$\Shv_{\tau}(\mathcal{C};\mathcal{D})\xrightarrow{\sim} 
\Shv_{\tau'}(\mathcal{C}';\mathcal{D})$.
\end{lemma}

\begin{proof}
By Lemma \ref{lem:equi-of-sites-infty-topoi}, 
we have an equivalence of 
$\infty$-topoi $\Shv_{\tau'}(\mathcal{C}')
\simeq \Shv_{\tau}(\mathcal{C})$.
Since $\Shv_{\tau}(\mathcal{C};\mathcal{D})$ can be identified with the 
$\infty$-category of limit-preserving functors from 
$\Shv_{\tau}(\mathcal{C})$ to $\mathcal{D}$, and similarly for 
$\mathcal{C}'$, we deduce an equivalence of $\infty$-categories
$\Shv_{\tau'}(\mathcal{C}';\mathcal{D})\simeq 
\Shv_{\tau}(\mathcal{C};\mathcal{D})$.
This equivalence is given by the restriction functor. 
Since the restriction of $F'$ to $\mathcal{C}$ is equivalent to 
$F$, we only need to prove that $F'$ is a $\tau'$-sheaf.
For $d\in \mathcal{D}$, denote by $\yon(d):\mathcal{D}\to \mathcal{S}$
the copresheaf corepresented by $d$. The functors 
$\yon(d)$, for $d\in \mathcal{D}$, 
form a conservative family of limit-preserving 
functors. Thus, it is enough to show that 
$\yon(d)(F')$ is a $\tau'$-sheaf for every $d\in \mathcal{D}$. 
Since 
$\yon(d)(F')$ is the right Kan extension of $\yon(d)(F)$, 
we are reduced to prove the lemma with $\mathcal{D}$
the $\infty$-category of spaces $\mathcal{S}$.

Recall that we need to show that $F'$ is a sheaf. Since $\mathcal{D}=\mathcal{S}$, we have at our disposal the sheafification 
functors, and these commute with restriction along the 
inclusion $\mathcal{C}\to \mathcal{C}'$. 
Let $F''$ be the $\tau'$-sheaf associated to $F'$. 
Since $F'|_{\mathcal{C}}\simeq F$ is already a $\tau$-sheaf, 
it follows that $F'\to F''$ induces an equivalence 
after restriction to $\mathcal{C}$.
By the universal property of the right Kan extension, there must be 
a map $F'' \to F'$ such that $F'\to F'' \to F'$ is homotopic 
to the identity of $F'$. 
Thus, $F'$ is a retract of the $\tau'$-sheaf $F''$. This proves
that $F'$ is also a $\tau'$-sheaf (and that $F'\simeq F''$).
\end{proof}

\begin{rmk}
\label{rmk:analytic-topology-on-int-rigspc}
The category 
$$\mathfrak{Q}=\left(\int_{\RigSpc^{\op}}\RigSpc^{\lft}\right)^{\op}$$
admits a natural topology, called the analytic topology and denoted by 
``$\an$''. It is induced by a pretopology 
${\rm Cov}_{\an}$ in the sense of 
\cite[Expos\'e II, D\'efinition 1.3]{SGAIV1},
which is given as follows.
For $(S,X)\in \mathfrak{Q}$, a family 
$((S_i,X_i)\to (S,X))_i$ belongs to 
${\rm Cov}_{\an}(S,X)$ 
if $(S_i\to S)_i$ is an open cover of $S$
and the morphisms $S_i\times_S X \to X_i$ are isomorphisms.
\end{rmk}

\begin{prop}
\label{prop:sheaf-property-for-rigsh-star-!-star-direct}
The functor 
\eqref{eq-thm:exist-functo-exceptional-image-37}
is a sheaf for the analytic topology on 
$\mathfrak{Q}$. 
\end{prop}

\begin{proof}
Fix an object $(S_{-1},X)$ in $\mathfrak{Q}$ 
and let $S_{\bullet}$ be a truncated
hypercover of $S_{-1}$ in the analytic topology.
We assume that the $S_n$'s are coproducts of  
open subspaces of $S_{-1}$.
For $n\in \N$, we set $X_n=S_n\times_{S_{-1}}X$
and similarly for every rigid analytic $S_{-1}$-space.
We need to show that 
\begin{equation}
\label{eq-prop:sheaf--for-rigsh-star-!-star-direct-1}
\RigSH^{(\hyp)}_{\tau}((S_{-1},X);\Lambda)^*_!
\to \lim_{[n]\in \mathbf{\Delta}}
\RigSH^{(\hyp)}_{\tau}((S_n,X_n);\Lambda)^*_!
\end{equation}
is an equivalence. 
By Lemma \ref{lemma:sheaf-extension-vs-right-Kan-exte}, the functor
$$\RigSH^{(\hyp)}_{\tau}(S_n\times_{S_{-1}}-;\Lambda)_!:
\Op/X \to \Prl$$
is the left Kan extension of its restriction to the subcategory
$\Op^{\wc}/X\subset \Op/X$ spanned by those open subspaces of 
$X$ which are weakly compactifiable over $S_{-1}$.
Using Proposition \ref{prop:restriction-rigsh-lft-!-star},
we deduce that 
$$\RigSH^{(\hyp)}_{\tau}((S_n,X_n);\Lambda)^*_!
\simeq \underset{U\in \Op^{\wc}/X}{\colim}\;
\RigSH^{(\hyp)}_{\tau}((S_n,U_n);\Lambda)^*_!$$
where the colimit is taken in $\Prl$.
Thus, we are reduced to showing that 
\begin{equation}
\label{eq-prop:sheaf--for-rigsh-star-!-star-direct-3}
\underset{U\in \Op^{\wc}/X}{\colim}\; 
\RigSH^{(\hyp)}_{\tau}((S_{-1},U);\Lambda)^*_!
\to \lim_{[n]\in \mathbf{\Delta}} 
\underset{U\in \Op^{\wc}/X}{\colim}\;
\RigSH^{(\hyp)}_{\tau}((S_n,U_n);\Lambda)^*_!
\end{equation}
is an equivalence. We want to apply 
\cite[Proposition 4.7.4.19]{lurie:higher-algebra}
for commuting the limit with the colimit in the right-hand side
of \eqref{eq-prop:sheaf--for-rigsh-star-!-star-direct-3}.
For this, we need to show that for every $[n']\to [n]$ 
in $\mathbf{\Delta}$ and every inclusion $U\to U'$ in 
$\Op^{\wc}/X$, 
the associated square
$$\xymatrix{\RigSH^{(\hyp)}_{\tau}((S_n,U_n);\Lambda)_!^*
\ar[r] \ar[d] & 
\RigSH^{(\hyp)}_{\tau}((S_n,U'_n);\Lambda)_!^* \ar[d] \\
\RigSH^{(\hyp)}_{\tau}((S_{n'},U_{n'});\Lambda)_!^*
\ar[r] & \RigSH^{(\hyp)}_{\tau}((S_{n'},U'_{n'});\Lambda)_!^*}$$
is right adjointable. Let $g:S_{n'}\to S_n$ 
be the morphism induced by $[n']\to [n]$, and let
$g':U_{n'}\to U_n$ and $g'':U'_{n'}\to U'_n$
be the morphisms obtained by base change. 
Let $u:U\to U'$ be the obvious inclusion, and let 
$u_n:U_n\to U'_n$ and $u_{n'}:U_{n'}\to U'_{n'}$ 
be the morphisms obtained by base change.
Then, using Lemma \ref{lem:left-kan-extension-properties}, and 
looking back at the construction of 
\eqref{eq-thm:exist-functo-exceptional-image-11}, 
we see that the above square is equivalent to 
$$\xymatrix{\RigSH^{(\hyp)}_{\tau}(U_n;\Lambda)
\ar[r]^-{u_{n,\,\sharp}} \ar[d]^-{g'^*} & 
\RigSH^{(\hyp)}_{\tau}(U'_n;\Lambda) \ar[d]^-{g''^*} \\
\RigSH^{(\hyp)}_{\tau}(U_{n'};\Lambda)
\ar[r]^-{u_{n',\,\sharp}} & 
\RigSH^{(\hyp)}_{\tau}(U'_{n'};\Lambda)}$$
which is clearly right adjointable.
Thus, \cite[Proposition 4.7.4.19]{lurie:higher-algebra}
applies, and we are left to showing that 
\begin{equation}
\label{eq-prop:sheaf--for-rigsh-star-!-star-direct-7}
\RigSH^{(\hyp)}_{\tau}((S_{-1},U);\Lambda)^*_!
\to \lim_{[n]\in \mathbf{\Delta}} 
\RigSH^{(\hyp)}_{\tau}((S_n,U_n);\Lambda)^*_!
\end{equation}
is an equivalence for every $U\in \Op^{\wc}/X$.
Said differently, we may assume that 
$X$ is weakly compactifiable. In this case, we may use Lemma
\ref{lem:rigsh-!-*-inverse-image} to rewrite 
\eqref{eq-prop:sheaf--for-rigsh-star-!-star-direct-1}
as follows:
\begin{equation}
\label{eq-prop:sheaf--for-rigsh-star-!-star-direct-19}
\RigSH^{(\hyp)}_{\tau}(X;\Lambda)^*
\to \lim_{[n]\in \mathbf{\Delta}}
\RigSH^{(\hyp)}_{\tau}(X_n;\Lambda)^*
\end{equation}
which is indeed an equivalence by Theorem
\ref{thm:hyperdesc}.
\end{proof}

At this stage, Theorem 
\ref{thm:exist-functo-exceptional-image}
is proven, except for the assertion that 
the functors in
\eqref{eq-thm:exist-functo-exceptional-image-1}
take an object $(S,X)$ to 
$\RigSH^{(\hyp)}_{\tau}(X;\Lambda)$. 
We do know this when $X$ weakly compactifiable
over $S$. In order to establish this
in general, we will need a few more results about the functors in  
\eqref{eq-thm:exist-functo-exceptional-image-1.81}. 
We first introduce a notation
which is useful in discussing these results.

\begin{nota}
\label{mot:momentary-notation-!-S}
The functors in 
\eqref{eq-thm:exist-functo-exceptional-image-1.81}
depend on $S$. To highlight this dependency,
we use ``$!_S$'' in subscript and superscript instead of
``$!$''. More explicitly, we denote  
by $\RigSH^{(\hyp)}_{\tau}(-;\Lambda)_{!_S}$ and 
$\RigSH^{(\hyp)}_{\tau}(-;\Lambda)^{!_S}$ these functors.
Also, given a morphism $f:Y \to X$ in $\RigSpc^{\lft}/S$, 
we sometimes denote by $f_{!_S}$ and $f^{!_S}$ the images 
of $f$ by these functors.
\symn{$\shriek_S$}
\end{nota}

\begin{lemma}
\label{lem:right-adjointability-in-RigSH-star-!}
Let $S$ be a rigid analytic space and $f:Y \to X$ a
morphism in $\RigSpc^{\lft}/S$.
Let $g:S' \to S$ be a morphism of rigid analytic spaces, 
and consider the Cartesian square
$$\xymatrix{Y' \ar[r]^-{g''} \ar[d]^-{f'} & Y \ar[d]^-f\\
X' \ar[r]^-{g'} & X}$$
where $f'$ is the base change of $f$ by $g$.
Consider the commutative square
$$\xymatrix{\RigSH^{(\hyp)}_{\tau}(X';\Lambda)^{!_{S'}} 
\ar[r]^-{g'_*} \ar[d]^-{f'^{!_{S'}}} & 
\RigSH^{(\hyp)}_{\tau}(X;\Lambda)^{!_S} 
\ar[d]^-{f^{!_S}}\\
\RigSH^{(\hyp)}_{\tau}(Y';\Lambda)^{!_{S'}} 
\ar[r]^-{g''_*} & 
\RigSH^{(\hyp)}_{\tau}(Y;\Lambda)^{!_S},\!}$$
where $g'_*$ is obtained by applying the second functor in 
\eqref{eq-thm:exist-functo-exceptional-image-1}
to the arrow $(g,\id_{X'}):(S',X')\to (S,X)$ and similarly for $g''_*$.
This square is left adjointable if $f$ or $g$ is an open immersion.
\end{lemma}

\begin{proof}
We may consider the commutative square in the statement 
as a morphism $(f^{!_{S'}},f^{!_S})$ in $\Fun({\Delta^1},\CAT_{\infty})$ between 
the functors $g'_*$ and $g''_*$, and our goal is 
to show that this morphism belongs to the sub-$\infty$-category
$\Fun^{\rm LAd}({\Delta^1},\CAT_{\infty})$
introduced in \cite[Definition 4.7.4.16]{lurie:higher-algebra}.
By \cite[Corollary 4.7.4.18]{lurie:higher-algebra}, 
it would be enough to show that the morphism $(f^{!_{S'}},f^{!_S})$
is the limit of an inverse system of morphisms in  
$\Fun^{\rm LAd}({\Delta^1},\CAT_{\infty})$.
By Proposition 
\ref{prop:restriction-rigsh-lft-!-star},
the morphism $(f^{!_{S'}},f^{!_S})$ is the limit of morphisms 
in $\Fun({\Delta^1},\CAT_{\infty})$
given by the following commutative squares 
$$\xymatrix{\RigSH^{(\hyp)}_{\tau}(S'\times_S U;\Lambda)^{!_{S'}} 
\ar[r] \ar[d] & 
\RigSH^{(\hyp)}_{\tau}(U;\Lambda)^{!_S} 
\ar[d]\\
\RigSH^{(\hyp)}_{\tau}(S'\times_S V;\Lambda)^{!_{S'}} 
\ar[r] & 
\RigSH^{(\hyp)}_{\tau}(V;\Lambda)^{!_S},\!}$$
where $U\subset X$ and $V \subset Y\times_X U$ are open subspaces 
which are weakly compactifiable over $S$.
Moreover, the transition maps in this inverse
system are given by commutative squares of the same type. 
Therefore, it is enough the show that these squares
are left adjointable, and thus we may assume that 
$X$ and $Y$ are weakly compactifiable over $S$.
In this case, we may use the explicit construction in 
Definition \ref{dfn:f-!-up-down} and 
Theorem \ref{thm:prop-base}(2) to conclude.
\end{proof}

\begin{lemma}
\label{lem:inverse-image-star-vs-!-open-immersion}
Let $S$ be a rigid analytic space and $j:S'\to S$
an open immersion. Let $Y\in \RigSpc^{\lft}/S$ such that 
the structure morphism $Y \to S$ factors through $S'$. 
Then, there exists an equivalence of $\infty$-categories
\begin{equation}
\label{eq-cor:inverse-image-star-vs-!-open-immersion-1}
\RigSH^{(\hyp)}_{\tau}(Y;\Lambda)^{!_S}\simeq 
\RigSH^{(\hyp)}_{\tau}(Y;\Lambda)^{!_{S'}}
\end{equation}
such that the following condition is satisfied.
For every morphism $f:Y \to X$ in 
$\RigSpc^{\lft}/S$, the functor 
$f^{!_S}$ is equivalent, modulo 
\eqref{eq-cor:inverse-image-star-vs-!-open-immersion-1}, 
to the composition of 
$$\RigSH^{(\hyp)}_{\tau}(X;\Lambda)^{!_S}
\xrightarrow{j'^*} 
\RigSH^{(\hyp)}_{\tau}(X';\Lambda)^{!_{S'}}
\xrightarrow{f'^{!_{S'}}} 
\RigSH^{(\hyp)}_{\tau}(Y;\Lambda)^{!_{S'}},$$
where $X'=S'\times_S X$, and $j':X'\to X$ 
and $f':Y\to X'$ are the obvious morphisms.
\end{lemma}

\begin{proof}
The image of the arrow
$(j,\id):(S,Y) \to (S',Y)$
by the first functor in 
\eqref{eq-thm:exist-functo-exceptional-image-1}
is a functor 
\begin{equation}
\label{eq-cor:inverse-image-star-vs-!-open-immersion-3}
\RigSH^{(\hyp)}_{\tau}(Y;\Lambda)^{!_S}\to 
\RigSH^{(\hyp)}_{\tau}(Y;\Lambda)^{!_{S'}}
\end{equation}
such that, for every $f:Y \to X$ as in the statement, the square 
$$\xymatrix{\RigSH^{(\hyp)}_{\tau}(X;\Lambda)^{!_S} 
\ar[d]^-{j'^*}\ar[r]^-{f^{!_S}} &
\RigSH^{(\hyp)}_{\tau}(Y;\Lambda)^{!_S} \ar[d]^-{
\eqref{eq-cor:inverse-image-star-vs-!-open-immersion-3}}\\
\RigSH^{(\hyp)}_{\tau}(X';\Lambda)^{!_{S'}} \ar[r]^-{f'^{!_{S'}}} &
\RigSH^{(\hyp)}_{\tau}(Y;\Lambda)^{!_{S'}}}$$
is commutative by Lemma
\ref{lem:right-adjointability-in-RigSH-star-!}.
Thus, to finish the proof, it is enough to show that 
\eqref{eq-cor:inverse-image-star-vs-!-open-immersion-3}
is an equivalence of $\infty$-categories. By Proposition
\ref{prop:restriction-rigsh-lft-!-star}, the 
question is local on $Y$. (Indeed, we may as well prove that
the right adjoint of
\eqref{eq-cor:inverse-image-star-vs-!-open-immersion-3}
is an equivalence of $\infty$-categories.) 
Thus, we may assume that
$Y$ is weakly compactifiable over $S$. 
In this case, we may use the explicit construction in 
Definition \ref{dfn:f-!-up-down} to conclude.
\end{proof}

\begin{lemma}
\label{lem:lower-!-belongs-to-Prl-omega}
Let $S$ be a rigid analytic space and let
$j:U \to X$ an open immersion in $\RigSpc^{\lft}/S$.
Then the functor 
$$j^{!_S}:\RigSH^{(\hyp)}_{\tau}(X;\Lambda)^{!_S}
\to \RigSH^{(\hyp)}_{\tau}(U;\Lambda)^{!_S}$$
belongs to $\Prl$ and hence admits a right adjoint, which we denote by 
$j_{?_S}$.
\symn{$?_S$}
\end{lemma}

\begin{proof}
Indeed, by Proposition 
\ref{prop:restriction-rigsh-lft-!-star},
$j^{!_S}$ is a limit in $\CAT_{\infty}$ of functors of the form 
$$\RigSH^{(\hyp)}_{\tau}(V;\Lambda)^{!_S}
\to \RigSH^{(\hyp)}_{\tau}(U\cap V;\Lambda)^{!_S}$$
for open subspaces $V\subset X$ which are compactifiable over 
$S$. By \cite[Proposition 5.5.3.13]{lurie}, 
it is thus enough to prove that $j^{!_S}$ is in $\Prl$ 
when $j$ is an open immersion between weakly 
compactifiable rigid analytic $S$-spaces.
In this case, we know that $j^{!_S}$ is equivalent to $j^*$, 
and the result follows.
\end{proof}

\begin{lemma}
\label{lem:ex-right-adjoint-j-?-adjointability}
Let $S$ be a rigid analytic space, and consider a Cartesian square
in $\RigSpc^{\lft}/S$
$$\xymatrix{V \ar[r]^-v \ar[d]^-g & Y \ar[d]^-f\\
U \ar[r]^-u & X,\!}$$
with $u$ an open immersion (resp. a closed immersion). 
Then, the commutative square
$$\xymatrix{\RigSH^{(\hyp)}_{\tau}(X;\Lambda)^{!_S} \ar[r]^-{u^{!_S}} 
\ar[d]^-{f^{!_S}} & 
\RigSH^{(\hyp)}_{\tau}(U;\Lambda)^{!_S} \ar[d]^-{g^{!_S}}\\
\RigSH^{(\hyp)}_{\tau}(Y;\Lambda)^{!_S} \ar[r]^-{v^{!_S}} & 
\RigSH^{(\hyp)}_{\tau}(V;\Lambda)^{!_S},}$$
is right adjointable (resp. left adjointable).
\end{lemma}

\begin{proof}
We only consider the case of open immersions; the case of closed 
immersions is similar.
Using Proposition 
\ref{prop:restriction-rigsh-lft-!-star}, 
\cite[Corollary 4.7.4.18]{lurie:higher-algebra}
and arguing as in the proof of Lemma 
\ref{lem:right-adjointability-in-RigSH-star-!}, 
we reduce to show the lemma when $X$ and $Y$ are 
weakly compactifiable over $S$. In this case, 
the commutative square of the statement coincides with the 
one deduced by adjunction from 
$$\xymatrix{\RigSH^{(\hyp)}_{\tau}(V;\Lambda)\ar[r]^-{v_{\sharp}}
\ar[d]^-{g_!} & \RigSH^{(\hyp)}_{\tau}(Y;\Lambda) \ar[d]^-{f_!}\\
\RigSH^{(\hyp)}_{\tau}(U;\Lambda)\ar[r]^-{u_{\sharp}}
\ar[r] & \RigSH^{(\hyp)}_{\tau}(X;\Lambda).}$$
The right adjointability of this square is clear:
it follows from the construction of the exceptional direct 
image functors given in Definition
\ref{dfn:f-!-up-down} and Proposition \ref{prop:6f1}(3).
\end{proof}

\begin{cons}
\label{cons:functor-i-?-for-a-locally-closed-immersion}
Let $S$ be a rigid analytic space and let
$i:Z \to X$ be a locally closed immersion in 
$\RigSpc^{\lft}/S$. 
We define a functor 
$$i_{?_S}:\RigSH^{(\hyp)}_{\tau}(Z;\Lambda)_{!_S}\to 
\RigSH^{(\hyp)}_{\tau}(X;\Lambda)_{!_S}$$
as follows. Choose an open subspace $U\subset X$ 
containing $Z$ as a closed subspace, and let 
$s:Z \to U$ and $j:U \to X$ be the obvious immersions. 
Define $i_{?_S}$ to be the composite functor 
$j_{?_S}\circ s_{!_S}$.
\symn{$?_S$}
\end{cons}

\begin{lemma}
\label{lem:i-?-indep-choice}
Keep the notations of Construction 
\ref{cons:functor-i-?-for-a-locally-closed-immersion}.
The functor $i_{?_S}$ is independent of the choice of 
the open neighbourhood $U$.
\end{lemma}

\begin{proof}
Let $U'\subset U$ be an open neighbourhood of $Z$ 
contained in $U$. Let $s':Z\to U'$ and $u:U'\to U$ 
be the obvious immersions. We need to show that 
$u_{?_S}\circ s'_{!_S}\simeq s_{!_S}$.
We have a Cartesian square
$$\xymatrix{Z \ar@{=}[r] \ar[d]^-{s'} 
& Z\ar[d]^-s\\
U' \ar[r]^-u & U}$$
which induces an equivalence 
$s'^{!_S}\simeq s^{!_S}\circ u_{?_S}$ by Lemma
\ref{lem:ex-right-adjoint-j-?-adjointability}. From this equivalence, 
we deduce a natural transformation 
$s_{!_S}\to u_{?_S}\circ s'_{!_S}$. This natural transformation
is an equivalence. Indeed, it is enough to check this after 
applying $u^{!_S}$ and $v^{!_S}$, with 
$v:U\smallsetminus Z \to U$ the obvious inclusion,
and this is easily seen to be true using Lemma
\ref{lem:ex-right-adjoint-j-?-adjointability} again.
\end{proof}

\begin{lemma}
\label{lem:right-adjointability-in-RigSH-star-?}
Let $S$ be a rigid analytic space and $i:Z \to X$ a locally closed 
immersion in $\RigSpc^{\lft}/S$.
Let $g:S' \to S$ be a morphism of rigid analytic spaces, 
and consider the Cartesian square
$$\xymatrix{Z' \ar[r]^-{g''} \ar[d]^-{i'} & Z \ar[d]^-i\\
X' \ar[r]^-{g'} & X}$$
where $i'$ is the base change of $i$ by $g$.
Then, there is a commutative square of $\infty$-categories 
$$\xymatrix{\RigSH^{(\hyp)}_{\tau}(Z';\Lambda)^{!_{S'}} 
\ar[r]^-{g''_*} \ar[d]^-{i'_{?_{S'}}} & 
\RigSH^{(\hyp)}_{\tau}(Z;\Lambda)^{!_S} \ar[d]^-{i_{?_S}}\\
\RigSH^{(\hyp)}_{\tau}(X';\Lambda)^{!_{S'}} \ar[r]^-{g'_*} & 
\RigSH^{(\hyp)}_{\tau}(X;\Lambda)^{!_S}}$$
(In the above square, 
$g'_*$ is obtained by applying the second functor in 
\eqref{eq-thm:exist-functo-exceptional-image-1}
to the arrow $(g,\id_{X'}):(S',X')\to (S,X)$ and similarly for $g''_*$.)
\end{lemma}

\begin{proof}
When $i$ is an open immersion, this follows from 
Lemma \ref{lem:right-adjointability-in-RigSH-star-!}.
Thus, we may assume that $i$ is a closed immersion, and we 
need to prove the analogous statement for the functors
$i_{!_S}$ and $i'_{!_{S'}}$. Arguing as in the proof of 
Lemma \ref{lem:right-adjointability-in-RigSH-star-!},
we reduce to the case where $X$ is weakly compactifiable. 
In this case, the functors $i_{!_{S}}$ and $i'_{!_{S'}}$
coincide with $i_*$ and $i'_*$, and the result follows.
\end{proof}

\begin{thm}
\label{thm:exist-functo-exceptional-image-contnd}
Let $S$ be a rigid analytic space and let $T\in \RigSpc^{\lft}/S$.
There is a commutative triangle
$$\xymatrix{(\RigSpc^{\lft}/T)^{\op} \ar[rr] 
\ar@/_/[dr]|-{\RigSH^{(\hyp)}_{\tau}(-;\,\Lambda)^{!_T}\quad}
& & (\RigSpc^{\lft}/S)^{\op}
\ar@/^/[dl]|-{\quad\RigSH^{(\hyp)}_{\tau}(-;\,\Lambda)^{!_S}}\\
& \Prr, &}$$
where the horizontal arrow is the forgetful functor.
For $X\in \RigSpc^{\lft}/T$, the induced equivalence 
of $\infty$-categories
\begin{equation}
\label{eq-thm:exist-functo-exceptional-image-contnd-1}
\RigSH^{(\hyp)}_{\tau}(X;\Lambda)^{!_T} \xrightarrow{\sim}
\RigSH^{(\hyp)}_{\tau}(X;\Lambda)^{!_S}
\end{equation}
is obtained as follows.
Consider the commutative diagram with a Cartesian square
$$\xymatrix{X \ar[r]^-{\delta_X} \ar[dr] & T\times_S X \ar[r]^-{\pr_X} \ar[d] & X \ar[d]\\
& T \ar[r]^-g & S.\!}$$
Then, the equivalence 
\eqref{eq-thm:exist-functo-exceptional-image-contnd-1}
is the composition of
\begin{equation}
\label{eq-thm:exist-functo-exceptional-image-contnd-13}
\RigSH^{(\hyp)}_{\tau}(X;\Lambda)^{!_T} \xrightarrow{(\delta_X)_{?_T}} 
\RigSH^{(\hyp)}_{\tau}(T\times_S X;\Lambda)^{!_T}
\xrightarrow{(\pr_X)_*} \RigSH^{(\hyp)}_{\tau}(X;\Lambda)^{!_S}.
\end{equation}
(Here, we denote by 
$(\pr_X)_*$ the image by the functor 
$\RigSH^{(\hyp)}_{\tau}(-;\Lambda)_*^!$
of the arrow
$(g,\id_{T\times_S X}):(S,X) \to (T,T\times_S X)$.)
\end{thm}

\begin{proof}
By Proposition
\ref{prop:restriction-rigsh-lft-!-star}
and Lemma \ref{lem:ex-right-adjoint-j-?-adjointability},
the composite functors 
\eqref{eq-thm:exist-functo-exceptional-image-contnd-13}
are part of a morphism of $\Prr$-valued 
sheaves on $\RigSpc^{\lft}/T$
$$\RigSH^{(\hyp)}_{\tau}(-;\Lambda)^{!_T}
\to \RigSH^{(\hyp)}_{\tau}(-;\Lambda)^{!_S}|_{\RigSpc^{\lft}/T}.$$
Thus, it is enough to prove that the composite functor
\eqref{eq-thm:exist-functo-exceptional-image-contnd-13}
is an equivalence under the following assumptions:
\begin{itemize}

\item $X$ is weakly compactifiable over $S$;

\item $X\to T$ factors by an open
subspace $T'\subset T$ which is weakly compactifiable over $S$.

\end{itemize}
The morphism $\delta_{X}:X\to T\times_S X$ is the composition of the
open immersion $j:T'\times_S X \to T\times_S X$ and the 
morphism $\delta'_X:X\to T'\times_S X$.
We deduce that the composition of 
\eqref{eq-thm:exist-functo-exceptional-image-contnd-13} 
is equivalent to the composition of 
$$\xymatrix{\RigSH^{(\hyp)}_{\tau}(X;\Lambda)^{!_T} 
\ar[r]^-{(\delta'_X)_{?_T}} &  
\RigSH^{(\hyp)}_{\tau}(T'\times_S X;\Lambda)^{!_T}
\ar[d]^-{j_{?_T}} & \\
& \RigSH^{(\hyp)}_{\tau}(T\times_S X;\Lambda)^{!_T}
\ar[r]^-{(\pr_X)_*} & \RigSH^{(\hyp)}_{\tau}(X;\Lambda)^{!_S}}$$
By Lemma
\ref{lem:inverse-image-star-vs-!-open-immersion}, the functor
$j^{!_T}$ is equivalent to the composition of 
$$\RigSH^{(\hyp)}_{\tau}(T\times_S X;\Lambda)^{!_T}
\xrightarrow{j^*}
\RigSH^{(\hyp)}_{\tau}(T'\times_S X;\Lambda)^{!_{T'}}\simeq 
\RigSH^{(\hyp)}_{\tau}(T'\times_S X;\Lambda)^{!_T}.$$ 
It follows that the functor $j_{?_T}$ is 
equivalent to the composition of 
$$\RigSH^{(\hyp)}_{\tau}(T'\times_S X;\Lambda)^{!_T}\simeq 
\RigSH^{(\hyp)}_{\tau}(T'\times_S X;\Lambda)^{!_{T'}}
\xrightarrow{j_*} \RigSH^{(\hyp)}_{\tau}(T\times_S X;\Lambda)^{!_T}.$$
Thus, modulo the equivalence 
$\RigSH^{(\hyp)}_{\tau}(X;\Lambda)^{!_T}\simeq 
\RigSH^{(\hyp)}_{\tau}(X;\Lambda)^{!_{T'}}$, 
the composition of
\eqref{eq-thm:exist-functo-exceptional-image-contnd-13} 
is equivalent to the composition of 
$$\RigSH^{(\hyp)}_{\tau}(X;\Lambda)^{!_{T'}} 
\xrightarrow{(\delta'_X)_{?_{T'}}} 
\RigSH^{(\hyp)}_{\tau}(T'\times_S X;\Lambda)^{!_{T'}}
\xrightarrow{(\pr'_X)_*} \RigSH^{(\hyp)}_{\tau}(X;\Lambda)^{!_S}$$
where $\pr_X'=\pr_X\circ j$.
Therefore, it is enough to prove the theorem with $T$
replaced by $T'$. Said differently, we may assume that 
$X$ and $T$ are weakly compactifiable over $S$. 
In this case, the diagram
\eqref{eq-thm:exist-functo-exceptional-image-contnd-13}
can be identified with
$$\RigSH^{(\hyp)}_{\tau}(X;\Lambda)^* \xrightarrow{(\delta_X)_*} 
\RigSH^{(\hyp)}_{\tau}(T\times_S X;\Lambda)^*
\xrightarrow{(\pr_X)_*} \RigSH^{(\hyp)}_{\tau}(X;\Lambda)^*$$
whose composition is clearly an equivalence.
\end{proof}

\begin{cor}
\label{cor:rigsh-!-equiv-rigsh-star}
For every rigid analytic space $S$ and every $X\in \RigSpc^{\lft}/S$, 
there is an equivalence of $\infty$-categories 
\begin{equation}
\label{eq-cor:rigsh-!-equiv-rigsh-star}
\RigSH^{(\hyp)}_{\tau}(X;\Lambda)^* \xrightarrow{\sim}
\RigSH^{(\hyp)}_{\tau}(X;\Lambda)^{!_S}.
\end{equation}
Moreover, these equivalences satisfy the following 
properties.
\begin{enumerate}

\item[(1)] Given a Cartesian square of rigid analytic spaces
$$\xymatrix{X'\ar[r]^-{g'} \ar[d]^-{f'} & X \ar[d]^-f\\
S'\ar[r]^-g & S}$$
with $f$ locally of finite type, 
there is a commutative square of $\infty$-categories
$$\xymatrix{\RigSH^{(\hyp)}_{\tau}(X';\Lambda)^* 
\ar[d]^-{\sim}  \ar[r]^-{g'_*}
& \RigSH^{(\hyp)}_{\tau}(X;\Lambda)^* \ar[d]^-{\sim} \\
\RigSH^{(\hyp)}_{\tau}(X';\Lambda)^{!_{S'}} \ar[r]^-{g'_*}
& \RigSH^{(\hyp)}_{\tau}(X;\Lambda)^{!_{S}}.\!}$$

\item[(2)] Given a rigid analytic space $S$ 
and an open immersion $j:X'\to X$ in $\RigSpc^{\lft}/S$, 
we have a commutative square
$$\xymatrix{\RigSH^{(\hyp)}_{\tau}(X;\Lambda)^* 
\ar[d]^-{\sim} \ar[r]^-{j^*} & 
\RigSH^{(\hyp)}_{\tau}(X';\Lambda)^* \ar[d]^-{\sim}\\
\RigSH^{(\hyp)}_{\tau}(X;\Lambda)^{!_S} \ar[r]^-{j^{!_S}} & 
\RigSH^{(\hyp)}_{\tau}(X';\Lambda)^{!_S}.\!}$$

\end{enumerate}
\end{cor}

\begin{proof}
The equivalence
\eqref{eq-cor:rigsh-!-equiv-rigsh-star} 
is the equivalence 
\eqref{eq-thm:exist-functo-exceptional-image-contnd-1}
when $X=T$. Property (1) follows easily from the construction
of the equivalence
\eqref{eq-cor:rigsh-!-equiv-rigsh-star}
and Lemma 
\ref{lem:right-adjointability-in-RigSH-star-?}.
Property (2) follows from Theorem
\ref{thm:exist-functo-exceptional-image-contnd}
combined with Corollary
\ref{cor:restriction-Open-rigsh-!}.
\end{proof}

Theorem 
\ref{thm:exist-functo-exceptional-image-contnd}
shows that the exceptional functors are independent of the base, i.e., 
the functors $f_{!_S}$ and $f^{!_S}$ are independent of 
$S$ up to equivalence. Note also that, by Proposition 
\ref{prop:restriction-rigsh-lft-!-star} and 
Corollary \ref{cor:ext-f-!-all-ft-mor}, these functors extend 
the ones of Definition \ref{dfn:f-!-up-down}.
This justifies the following definition.

\begin{dfn}
\label{dfn:exceptional-functors-in-general}
\ncn{exceptional functors}
Let $f:Y \to X$ be a morphism of rigid analytic spaces which
is locally of finite type. The functors in adjunction
$$f_!:\RigSH^{(\hyp)}_{\tau}(Y;\Lambda)
\rightleftarrows \RigSH^{(\hyp)}_{\tau}(X;\Lambda):f^!$$
are defined to be the images of the arrow 
$(\id_X,f):(X,Y) \to (X,X)$ by the functors in 
\eqref{eq-thm:exist-functo-exceptional-image-1}
modulo the equivalence 
$\RigSH^{(\hyp)}_{\tau}(Y;\Lambda)\simeq 
\RigSH^{(\hyp)}_{\tau}(Y;\Lambda)^{!_X}$
given by Corollary 
\ref{cor:rigsh-!-equiv-rigsh-star}.
The functors $f_!$ and $f^!$ are called the exceptional 
direct and inverse image functors.
\symn{$(-)_{\shriek}, (-)^{\shriek}$}
\end{dfn}

\begin{rmk}
\label{rmk:composition-of-exceptional-images}
Given two morphisms $f:Y\to X$ and $g:Z\to Y$
which are locally of finite type, we have 
equivalences $f_!\circ g_!\simeq (f\circ g)_!$ and
$g^!\circ f^!\simeq (f\circ g)^!$. (This follows from 
the construction and the equivalences 
$f_{!_X}\circ g_{!_X}\simeq (f\circ g)_{!_X}$ and
$g^{!_X}\circ f^{!_X}\simeq (f\circ g)^{!_X}$.)
Therefore, one expects to have functors, 
from the wide subcategory of $\RigSpc$ spanned by 
locally of finite type morphisms, to $\Prl$ and $(\Prr)^{\op}$,
sending a morphism $f$ to the functors $f_!$ and $f^!$.
Our method does not give readily such a functor, 
but techniques from \cite[Part III]{Gait-Rozen-I}
might do. (See Theorem \ref{thm:corresp-6-funct}
and Remark \ref{rmk:corresp-monoidal-struct}
below.)
\end{rmk}

\begin{prop}
\label{prop:exchange-base-change-in-general-1}
Consider a Cartesian square of rigid analytic spaces
$$\xymatrix{Y' \ar[r]^-{g'} \ar[d]^-{f'} & Y \ar[d]^-f\\
X' \ar[r]^-g & X}$$
with $f$ locally of finite type.  
Then, there is a commutative square of $\infty$-categories
$$\xymatrix{\RigSH^{(\hyp)}_{\tau}(Y;\Lambda)
\ar[r]^-{g'^*} \ar[d]^-{f_!} & 
\RigSH^{(\hyp)}_{\tau}(Y';\Lambda)
\ar[d]^-{f'_!}\\
\RigSH^{(\hyp)}_{\tau}(X;\Lambda) \ar[r]^-{g^*}
& \RigSH^{(\hyp)}_{\tau}(X';\Lambda).\!}$$
\end{prop}

\begin{proof}
Applying the first functor in 
\eqref{eq-thm:exist-functo-exceptional-image-1}
to the commutative square
$$\xymatrix{(X,Y) \ar[r]\ar[d] & (X',Y') \ar[d]\\
(X,X) \ar[r] & (X',X'),}$$
we get a commutative square of $\infty$-categories
$$\xymatrix{\RigSH^{(\hyp)}_{\tau}(Y;\Lambda)^{!_X} 
\ar[r]^-{g'^*} \ar[d]^-{f_{!_X}} & 
\RigSH^{(\hyp)}_{\tau}(Y';\Lambda)^{!_{X'}}
\ar[d]^-{f'_{!_{X'}}}\\
\RigSH^{(\hyp)}_{\tau}(X;\Lambda)^{!_X} \ar[r]^-{g^*}
& \RigSH^{(\hyp)}_{\tau}(X';\Lambda)^{!_{X'}}.\!}$$
The result follows then from Corollary 
\ref{cor:rigsh-!-equiv-rigsh-star}.
\end{proof}

\begin{prop}
\label{prop:f-proper-f-star-f-!}
The composition of the first functor in 
\eqref{eq-thm:exist-functo-exceptional-image-1}
with the obvious inclusion 
$$\int_{\RigSpc^{\op}} \RigSpc^{\proper} 
\to \int_{\RigSpc^{\op}} \RigSpc^{\lft}$$
is equivalent to the functor  
\eqref{eq-thm:exist-functo-exceptional-image-5}.
In particular, if $f:Y \to X$ is a proper 
morphism of rigid analytic spaces,
there is an equivalence $f_!\simeq f_*$.
\end{prop}

\begin{proof}
This is a direct consequence of the construction.
\end{proof}

\begin{cor}
\label{cor:f!-nice-f}
Let $f:Y \to X$ be a morphism of rigid analytic spaces.
Assume that $f$ admits a factorization $f=p\circ j$ 
where $j$ is an open immersion and $p$ is a proper morphism. 
Then, there is an equivalence $f_!\simeq p_*\circ j_{\sharp}$.
\end{cor}

\begin{proof}
This follows from 
Corollary \ref{cor:rigsh-!-equiv-rigsh-star}(2),
Remark \ref{rmk:composition-of-exceptional-images}
and Proposition \ref{prop:f-proper-f-star-f-!}.
\end{proof}

\begin{thm}[Ambidexterity]
\label{thm:f-!-for-f-smooth}
\ncn{ambidexterity}
Let $f:Y \to X$ be a smooth morphism between rigid analytic spaces.
There are equivalences $f_!\simeq f_{\sharp}\circ \Th^{-1}(\Omega_f)$ 
and $f^!\simeq \Th(\Omega_f)\circ f^*$.
\end{thm}

\begin{proof}
We first construct a natural transformation
$\alpha_f:f_{\sharp} \to f_!\circ \Th(\Omega_f)$.
Consider the commutative diagram with a Cartesian square
$$\xymatrix{Y \ar[dr]|-{\Delta_f} \ar@{=}@/^1.2pc/[drr] 
\ar@{=}@/_1.2pc/[ddr] & & \\
& Y\times_X Y \ar[r]^-{p_2} \ar[d]_-{p_1} & Y \ar[d]^-f\\
& Y \ar[r]^-f & X.}$$
By Proposition
\ref{prop:exchange-base-change-in-general-1}, 
we have an equivalence $p_{1,\, !}\circ p_2^*\simeq f^*\circ f_!$.
Using the adjunctions $(f_{\sharp},f^*)$ and 
$(p_{2,\,\sharp},p_2^*)$, we deduce a natural transformation 
$f_{\sharp} \circ p_{1,\,!}
\to f_! \circ p_{2,\,\sharp}$.
Applying the latter to $\Delta_{f,\,!}$
and using the equivalences $p_{1,\,!}\circ \Delta_{f,\,!}\simeq \id$
and $p_{2,\,\sharp} \circ \Delta_{f,\,!}\simeq \Th(\Omega_f)$, 
we get $\alpha_f$.

We next show that $\alpha_f$ is an equivalence. It is easy to see
that $\alpha_f$ is compatible with composition, i.e., that the 
analogue of \cite[Proposition 1.7.3]{ayoub-th1} is satisfied.
Moreover, if $j$ is an open immersion, 
$\alpha_j$ is the equivalence $j_{\sharp}\simeq j_!$.
Thus, to show that $\alpha_f$ is invertible, we may argue locally 
on $Y$ for the analytic topology.
Thus, we may assume that $Y$ is weakly compactifiable over $X$. 
Choose a weak compactification $i:Y \to W$ and let 
$g:W\to X$ be the structural morphism.
To prove that $\alpha_f$ is invertible, it is enough 
to show that the natural transformation 
$f_{\sharp}\circ p_{1,\,!}\to f_!\circ p_{2,\,\sharp}$
is invertible. 
Unwinding the definitions, we see that it is enough
to prove that the natural transformation 
$f_{\sharp}\circ \overline{q}_*
\to \overline{f}_*\circ q_{\sharp}$
associated to the Cartesian square
$$\xymatrix{Y\times_X W \ar[r]^-{f'}
\ar[d]^-{g'} & W
\ar[d]^-g \\
Y \ar[r]^-f & X}$$
is an equivalence. This is indeed true by 
Theorem \ref{thm:prop-base}(2).
\end{proof}

There is another way to encapsulate 
much of the six-functor formalism using $(\infty,2)$-categories of 
correspondences (aka., spans). This gives an alternative 
approach to the constructions of this subsection 
which is more elegant and more powerful.
The technology needed to carry out this approach 
is developed in \cite[Part III]{Gait-Rozen-I}
but relies, unfortunately, on yet unproven hypotheses 
in the theory of $(\infty,2)$-categories; 
see \cite[Chapter 10, \S 0.4]{Gait-Rozen-I}.
It is for this reason that we decided to develop 
a more self-contained approach. 
However, for the reader who is willing to 
accept the unproven hypotheses in loc.~cit., 
we briefly explain how this is supposed to work.
For a similar discussion in the context of 
equivariant motives, see \cite[\S 6.2]{hoyois-6op}.

\begin{rmk}
\label{rmk:infty-cat-of-correspondences}
Given an $\infty$-category $\mathcal{C}$ with 
finite limits, there is 
an associated $(\infty,2)$-category $\Corr(\mathcal{C})$
having the same objects as $\mathcal{C}$, and where $1$-morphisms
between $X$ and $Y$ are given by spans
$$\xymatrix@C=1pc{& \ar[dr]^-f Z \ar[dl]_-g & \\
X & & Y,}$$
i.e., maps $(f,g):Z \to X\times Y$. 
Given a second span $(f',g'):Z' \to X\times Y$, a $2$-morphism $(f',g')\Rightarrow (f,g)$ is a morphism $h:Z'\to Z$ such that $g'=gh$ and 
$f'=fh$. If $P_1$, $P_2$ and $P_3$ are properties of morphisms in 
$\mathcal{C}$, we denote by 
${\rm Corr}(\mathcal{C})_{P_1,P_2}^{P_3}$ the subcategory 
obtained by imposing $P_1$, $P_2$ and $P_3$ on the morphisms 
$f$, $g$ and $h$ above. For details, on the $(\infty,2)$-category 
$\Corr(\mathcal{C})$, we refer the reader to 
\cite[Chapter 7, \S 1.2]{Gait-Rozen-I}. Below, we will be 
interested in the $(\infty,2)$-category 
$\Corr(\RigSpc)^{\rm proper}_{{\rm all},\, \wc}$, 
where $2$-morphisms are 
given by proper maps, and right legs of spans are requested to 
be weakly compactifiable while no condition is imposed on left legs.
\symn{$\Corr$}
\end{rmk}

\begin{thm}
\label{thm:corresp-6-funct}
There is a $2$-functor 
\begin{equation}
\label{eq-thm:corresp-6-funct}
\RigSH^{(\hyp)}_{\tau}(-;\Lambda):({\rm Corr}(\RigSpc)^{\rm proper}_{{\rm all},\,\wc})^{2\text{-}\op} \to \Prl
\end{equation}
sending a span of the form 
$X \xleftarrow{f} Y \xrightarrow{\id} Y$ to $f^*$ 
and a span of the form
$Y\xleftarrow{\id} Y \xrightarrow{f} X$ to $f_!$.
(Above, $\Prl$ is considered as an $(\infty,2)$-category in the natural 
way, i.e., where $2$-morphisms are given by natural transformations.)
\end{thm}

\begin{proof}
We denote by ``$\proper$'' (resp. ``${\rm iso}$'', ``${\rm open}$'', 
``${\rm closed}$'', ``${\rm imm}$'') 
the class of proper morphisms (resp. isomorphisms, 
open immersions, closed immersions, locally closed immersions) 
in $\RigSpc$.
By \cite[Chapter 7, Theorem 3.2.2]{Gait-Rozen-I} and
Theorem \ref{thm:prop-base}(1), 
there exists a unique $2$-functor 
\begin{equation}
\label{eq-thm:corresp-6-funct-3}
\RigSH^{(\hyp)}_{\tau}(-;\Lambda):
({\rm Corr}(\RigSpc)^{\proper}_{{\rm all}, \, 
\proper})^{2\text{-}\op}
\to \CAT_{\infty}
\end{equation}
extending the functor 
$\RigSH^{(\hyp)}_{\tau}(-;\Lambda)^*:
\RigSpc^{\op}\to \Prl$.
Also, by the same theorem of loc.~cit., there exists a unique $2$-functor 
\begin{equation}
\label{eq-thm:corresp-6-funct-7}
\RigSH^{(\hyp)}_{\tau}(-;\Lambda):
({\rm Corr}(\RigSpc)^{{\rm iso}}_{{\rm all}, \, 
{\rm open}})^{2\text{-}\op}
\to \CAT_{\infty}
\end{equation}
extending the same functor. In particular, 
these two extensions coincide on $(\RigSpc^{\rm qcqs})^{\op}$.
By \cite[Chapter 7, Theorem 5.2.4]{Gait-Rozen-I} and 
Proposition \ref{prop:loc1}, we may glue uniquely 
\eqref{eq-thm:corresp-6-funct-7}
with the restriction of 
\eqref{eq-thm:corresp-6-funct-3}
to $({\rm Corr}(\RigSpc)^{\rm iso}_{{\rm all}, \, 
{\rm closed}})^{2\text{-}\op}$ and get a $2$-functor 
\begin{equation}
\label{eq-thm:corresp-6-funct-17}
\RigSH^{(\hyp)}_{\tau}(-;\Lambda):
({\rm Corr}(\RigSpc)^{{\rm iso}}_{{\rm all}, \, 
{\rm imm}})^{2\text{-}\op}
\to \CAT_{\infty}
\end{equation}
By a second application of 
\cite[Chapter 7, Theorem 5.2.4]{Gait-Rozen-I} and 
using Proposition \ref{prop:base-change-for-i-!-locally closed}, 
we can glue uniquely
\eqref{eq-thm:corresp-6-funct-3} and
\eqref{eq-thm:corresp-6-funct-17} to get the $2$-functor 
\eqref{eq-thm:corresp-6-funct}
in the statement.
\end{proof}

\begin{rmk}
\label{rmk:corresp-monoidal-struct}
We denote by ``$\lft$'' the class of 
morphisms which are locally of finite type.
It is conceivable that the 
$2$-functor \eqref{eq-thm:corresp-6-funct}
can be extended to a $2$-functor 
\begin{equation}
\label{eq-thm:corresp-6-funct-lft}
\RigSH^{(\hyp)}_{\tau}(-;\Lambda):({\rm Corr}(\RigSpc)^{\rm proper}_{{\rm all},\,\lft})^{2\text{-}\op} \to \Prl
\end{equation}
sending a span of the form 
$X \xleftarrow{f} Y \xrightarrow{\id} Y$ to $f^*$ 
and a span of the form
$Y\xleftarrow{\id} Y \xrightarrow{f} X$ to the functor $f_!$
of Definition 
\ref{dfn:exceptional-functors-in-general}.
We do not pursue this here.
\end{rmk}

\subsection{Projection formula}

$\empty$

\smallskip

\label{subsect:exceptional-complem}

In this subsection, we explain 
how to incorporate the projection formula for the exceptional 
direct image functors into the functor 
$\RigSH^{(\hyp)}_{\tau}(-;\Lambda)^*_!$
of Theorem \ref{thm:exist-functo-exceptional-image}.

\begin{thm}
\label{thm:projection-formula-!-star-rigsh}
The functor $\RigSH^{(\hyp)}_{\tau}(-;\Lambda)^*_!$ from
Theorem \ref{thm:exist-functo-exceptional-image}
admits a structure of a module over the composite functor
\begin{equation}
\label{eq-thm:projection-formula-!-star-rigsh-1}
\int_{\RigSpc^{\op}}\RigSpc^{\lft}
\to \RigSpc^{\op} \xrightarrow{\RigSH^{(\hyp)}_{\tau}
(-;\,\Lambda)^{\otimes}} \CAlg(\Prl),
\end{equation}
considered as a commutative algebra in the 
$\infty$-category of functors from 
$\int_{\RigSpc^{\op}}\RigSpc^{\lft}$ to $\Prl$. 
(The first functor in 
\eqref{eq-thm:projection-formula-!-star-rigsh-1} 
is the one given by $(S,X)\mapsto S$.)
Said differently, there is a functor 
\begin{equation}
\label{eq-thm:projection-formula-!-star-rigsh-17}
\RigSH^{(\hyp)}_{\tau}(-;\,\Lambda)^{\otimes}_!:
\int_{\RigSpc^{\op}}\RigSpc^{\lft}
\to \Mod(\Prl)
\end{equation}
which is a lifting of the
functor $\RigSH^{(\hyp)}_{\tau}(-;\Lambda)^*_!$ 
and which is part of a commutative square
$$\xymatrix{
\int_{\RigSpc^{\op}}\RigSpc^{\lft}
\ar[rrr]^-{\RigSH^{(\hyp)}_{\tau}(-;\,\Lambda)^{\otimes}_!} 
\ar[d] &&& 
\Mod(\Prl) \ar[d] \\
\RigSpc^{\op} 
\ar[rrr]^-{\RigSH^{(\hyp)}_{\tau}(-;\,\Lambda)^{\otimes}}
&&& \CAlg(\Prl).\!}$$
\end{thm}

\begin{proof}
We only sketch the argument, leaving some details to the 
reader. The proof consists in revisiting the 
construction of the functor $\RigSH^{(\hyp)}_{\tau}(-;\Lambda)^*_!$ 
of Theorem \ref{thm:exist-functo-exceptional-image},
exhibiting step by step a natural module structure over
a suitable variant of the algebra 
\eqref{eq-thm:projection-formula-!-star-rigsh-1}.
We start by remarking that the functor
\eqref{eq-thm:exist-functo-exceptional-image-3}
lifts to a functor 
$$\RigSH^{(\hyp)}_{\tau}(-;\Lambda)^{\otimes,\,\otimes}:
\int_{\RigSpc^{\op}} (\RigSpc^{\proper})^{\op}\to 
\CAlg(\Prl)$$
admitting a natural transformation from the composite functor 
$$\int_{\RigSpc^{\op}}(\RigSpc^{\proper})^{\op}
\to \RigSpc^{\op}
\xrightarrow{\RigSH^{(\hyp)}_{\tau}(-;\,\Lambda)^{\otimes}}
\CAlg(\Prl).$$
(The first functor in the composition above is given by 
$(S,X)\mapsto S$.)
Retaining merely the induced module structure on 
\eqref{eq-thm:exist-functo-exceptional-image-3},
we obtain a commutative square
$$\xymatrix{\int_{\RigSpc^{\op}}(\RigSpc^{\proper})^{\op}
\ar[rrr]^-{\RigSH^{(\hyp)}_{\tau}(-;\,\Lambda)^{\otimes,\,*}} \ar[d] 
& & & \Mod(\Prl) \ar[d]\\
\RigSpc^{\op} \ar[rrr]^-{\RigSH^{(\hyp)}_{\tau}
(-;\,\Lambda)^{\otimes}} & & & \CAlg(\Prl).\!}$$
With ${\rm K}$ as in Construction 
\ref{cons:Mod-C-otimes-}, 
we set ${\rm K}_1=\langle 1 \rangle \times_{\Fin,\,e_0}{\rm K}$.
We may view the upper horizontal arrow in the previous square 
as a functor 
\begin{equation}
\label{eq-thm:projection-formula-!-star-rigsh-3}
\RigSH^{(\hyp)}_{\tau}
(-;\,\Lambda)^{\otimes,\,*}:
\left(\int_{\RigSpc^{\op}}(\RigSpc^{\proper})^{\op}\right)
\times {\rm K}_1 \to \Prl.
\end{equation}
Informally, this functor takes a pair of objects 
$((S,X), r:\langle 1\rangle \to \langle m \rangle)$
to the tensor product in $\Prlmon$ of copies of 
$\RigSH^{(\hyp)}_{\tau}(S;\Lambda)$, one for each 
$i\in \{1,\ldots, m\}$
different from $r(1)$, and a copy of 
$\RigSH^{(\hyp)}_{\tau}(X;\Lambda)$, only when $r(1)\in \{1,\ldots, m\}$.
Moreover, an arrow of the form 
$((\id_S,\id_X),s:\langle m \rangle \to \langle n \rangle)$
is sent to a functor induced by the tensor product on
$\RigSH^{(\hyp)}_{\tau}(S;\Lambda)$,
and the tensor product of an object of 
$\RigSH^{(\hyp)}_{\tau}(S;\Lambda)$ with an object of 
$\RigSH^{(\hyp)}_{\tau}(X;\Lambda)$, i.e., the functor 
$$\RigSH^{(\hyp)}_{\tau}(S;\Lambda)
\otimes \RigSH^{(\hyp)}_{\tau}(X;\Lambda)\to 
\RigSH^{(\hyp)}_{\tau}(X;\Lambda),$$
given by $(M,N)\mapsto f^*(M)\otimes N$
where $f:X \to S$ is the structural morphism.
Using this description, it follows from Theorem
\ref{thm:prop-base}(1) 
and Proposition 
\ref{prop:proper-projective-formula-weakly-propro} 
that the condition ($\star$) in Construction 
\ref{cons:adjointable-cocart-fibration}
is satisfied for the functor
\eqref{eq-thm:projection-formula-!-star-rigsh-3}. 
(What plays the role of the simplicial set ``$S$''
in that construction is the category 
$\RigSpc^{\op}\times {\rm K}_1$.)
Applying Construction 
\ref{cons:adjointable-cocart-fibration}, 
we obtain a functor 
\begin{equation}
\label{eq-thm:projection-formula-!-star-rigsh-7}
\RigSH^{(\hyp)}_{\tau}
(-;\,\Lambda)^{\otimes}_*:
\left(\int_{\RigSpc^{\op}}(\RigSpc^{\proper})\right)
\times {\rm K}_1 \to \Prl.
\end{equation}
This functor 
is easily seen to correspond to a $\Mod(\Prl)$-valued functor 
$\RigSH^{(\hyp)}_{\tau}(-;\,\Lambda)^{\otimes}_*$
which is a lift of
\eqref{eq-thm:exist-functo-exceptional-image-5}
and which is part of a commutative square
$$\xymatrix{\int_{\RigSpc^{\op}}(\RigSpc^{\proper})
\ar[rrr]^-{\RigSH^{(\hyp)}_{\tau}(-;\,\Lambda)^{\otimes}_*} \ar[d] 
& & & \Mod(\Prl) \ar[d]\\
\RigSpc^{\op} \ar[rrr]^-{\RigSH^{(\hyp)}_{\tau}
(-;\,\Lambda)^{\otimes}} & & & \CAlg(\Prl).\!}$$
Given $(S,(X,W))\in \int_{\RigSpc^{\op}}
\WComp$, the sub-$\infty$-category
$$\RigSH^{(\hyp)}((X,W);\Lambda)^*_!
\subset \RigSH^{(\hyp)}_{\tau}(W;\Lambda)^*_*$$
(see Notations \ref{nota:subcat-rigsh-!-wcomp} and 
\ref{nota:subcat-rigsh-star-!})
is stable by tensoring with any object of 
$\RigSH^{(\hyp)}_{\tau}(W;\Lambda)^*_*$
and, in particular, by the inverse image of any object of 
$\RigSH^{(\hyp)}_{\tau}(S;\Lambda)^*$.
(This is an immediate consequence of 
Proposition \ref{prop:6f1}(2).)
Applying Lemma 
\ref{lem:fully-faith-sub-functor}
to the restriction of the functor 
\eqref{eq-thm:projection-formula-!-star-rigsh-7}
to the category $\int_{\RigSpc^{\op}}
\WComp$, we obtain a functor 
\begin{equation}
\label{eq-thm:projection-formula-!-star-rigsh-43}
\RigSH^{(\hyp)}_{\tau}(-;\Lambda)^{\otimes}_!:
\int_{\RigSpc^{\op}}
\WComp \to \Mod(\Prl)
\end{equation}
which is a lift of 
\eqref{eq-thm:exist-functo-exceptional-image-11}
and which is part of a commutative square as above.
The remainder of the construction follows closely 
the construction of the functor 
$\RigSH^{(\hyp)}_{\tau}(-;\Lambda)^*_!$ of 
Theorem \ref{thm:exist-functo-exceptional-image}.
Namely, we take a left Kan extension of
\eqref{eq-thm:projection-formula-!-star-rigsh-43}
along the functor 
\eqref{eq-thm:exist-functo-exceptional-image-frak-d}, 
and then a second left Kan extension along the 
fully faithful embedding
\eqref{eq-thm:exist-functo-exceptional-image-iota}. 
That the resulting functor 
\begin{equation}
\label{eq-thm:projection-formula-!-star-rigsh-91}
\RigSH^{(\hyp)}_{\tau}(-;\Lambda)^{\otimes}_!:
\int_{\RigSpc^{\op}}
\RigSpc^{\lft} \to \Mod(\Prl)
\end{equation}
is a lift of \eqref{eq-thm:exist-functo-exceptional-image-13}
follows from \cite[Proposition 4.3.3.10]{lurie}
and \cite[Corollary 3.4.4.6(2)]{lurie:higher-algebra}.
\end{proof}

\begin{prop}
\label{prop:equiv-RigSH-!-star-compat-module}
Let $S$ be a rigid analytic space and 
$X\in \RigSpc^{\lft}/S$. 
There exists an equivalence of 
$\RigSH^{(\hyp)}_{\tau}(S;\Lambda)^{\otimes}$-modules
\begin{equation}
\label{eq-prop:equiv-RigSH-!-star-compat-module-1}
\RigSH^{(\hyp)}_{\tau}((S,X);\Lambda)^{\otimes}_!
\simeq \RigSH^{(\hyp)}_{\tau}(X;\Lambda)^{\otimes}
\end{equation}
which is a lift of the equivalence of $\infty$-categories
provided by Corollary
\ref{cor:rigsh-!-equiv-rigsh-star}.
\end{prop}

\begin{proof}
We want to show that the inverse of the equivalence
\eqref{eq-cor:rigsh-!-equiv-rigsh-star}
can be naturally lifted to a morphism of 
$\RigSH^{(\hyp)}_{\tau}(S;\Lambda)^{\otimes}$-modules.
This equivalence is given by 
the composition of
$$\RigSH^{(\hyp)}_{\tau}((S,X);\Lambda)^*_!
\xrightarrow{(\pr_2)^*} 
\RigSH^{(\hyp)}_{\tau}((X,X\times_SX);\Lambda)^*_!
\xrightarrow{(\delta_X)^?} 
\RigSH^{(\hyp)}_{\tau}((X,X);\Lambda)^*_!$$
where: 
\begin{itemize}

\item $\pr_2:X\times_S X \to X$ is the 
projection to the second factor and 
$\delta_X:X \to X\times_S X$ is the diagonal embedding;

\item $(\delta_X)^?$ is the left adjoint of the functor
$(\delta_X)_?$ as in Construction 
\ref{cons:functor-i-?-for-a-locally-closed-immersion}.

\end{itemize}
The existence of $(\delta_X)^?$ follows from Proposition
\ref{prop:f-proper-f-star-f-!}
which insures that the functor $i_{!_X}$, for $i$ a closed immersion
of rigid analytic $X$-spaces, admits a left adjoint.
The functor $(\pr_2)^*$ admits a natural lift to a 
morphism of 
$\RigSH^{(\hyp)}_{\tau}(S;\Lambda)^{\otimes}$-modules.
So, we are left to prove the same for 
$(\delta_X)^?$. More generally, it is enough 
to prove the following assertions (with $T$ a rigid analytic space).
\begin{enumerate}

\item[(1)] If $j:V \to Y$ is an open immersion in 
$\RigSpc^{\lft}/T$, the functor 
$$j^!:\RigSH^{(\hyp)}_{\tau}((T,Y);\Lambda)^*_!
\to \RigSH^{(\hyp)}_{\tau}((T,V);\Lambda)^*_!$$
lifts to a morphism of 
$\RigSH^{(\hyp)}_{\tau}(Y;\Lambda)^{\otimes}$-modules.

\item[(2)] If $i:Z \to Y$ is a closed immersion in $\RigSpc^{\lft}/T$, 
the functor 
$$i^?:\RigSH^{(\hyp)}_{\tau}((T,Y);\Lambda)^*_!
\to \RigSH^{(\hyp)}_{\tau}((T,Z);\Lambda)^*_!$$
lifts to a morphism of $\RigSH^{(\hyp)}_{\tau}(Y;\Lambda)^{\otimes}$-modules.

\end{enumerate}
For the first assertion, starting with the morphism of 
$\RigSH^{(\hyp)}_{\tau}(Y;\Lambda)^{\otimes}$-modules
$j_!$, we need to show that the morphism
\begin{equation}
\label{eq-prop:equiv-RigSH-!-star-compat-module-3}
j^!(A)\otimes B \to j^!(A\otimes B)
\end{equation}
is an equivalence for $A\in \RigSH^{(\hyp)}_{\tau}((T,Y);\Lambda)^*_!$
and $B\in \RigSH^{(\hyp)}_{\tau}(Y;\Lambda)$.
This can be checked locally on $Y$, and thus we may 
assume that $Y$ is weakly compactifiable over $T$.
In this case, the morphism
\eqref{eq-prop:equiv-RigSH-!-star-compat-module-3} 
can be identified with the equivalence
$j^*(A)\otimes j^*(B) \simeq j^*(A\otimes B)$.
Similarly, for the second assertion,  
starting with the morphism of 
$\RigSH^{(\hyp)}_{\tau}(Y;\Lambda)^{\otimes}$-modules
$i_!$, we need to show that the morphism
\begin{equation}
\label{eq-prop:equiv-RigSH-!-star-compat-module-7}
i^?(A\otimes B) \to i^?(A)\otimes B
\end{equation}
is an equivalence for $A\in \RigSH^{(\hyp)}_{\tau}((T,Y);\Lambda)^*_!$
and $B\in \RigSH^{(\hyp)}_{\tau}(Y;\Lambda)$.
This can checked locally on $Y$, and thus we may 
assume that $Y$ is weakly compactifiable.
In this case, the morphism
\eqref{eq-prop:equiv-RigSH-!-star-compat-module-7}
can be identified with the equivalence
$i^*(A\otimes B)\simeq i^*(A)\otimes i^*(B)$.
\end{proof}

\begin{cor}[Projection formula]
\label{cor:proj-form-f-!}
\ncn{projection formula}
Let $f:Y \to X$ be a morphism of rigid analytic 
spaces which is locally of finite type. 
Then, the functor 
$$f_!:\RigSH^{(\hyp)}_{\tau}(Y;\Lambda)
\to \RigSH^{(\hyp)}_{\tau}(X;\Lambda),$$ 
as in Definition 
\ref{dfn:exceptional-functors-in-general}, admits a 
lift to a morphism of 
$\RigSH^{(\hyp)}_{\tau}(X;\Lambda)^{\otimes}$-modules.
In particular, there is an equivalence 
$$M\otimes f_!N\simeq f_!(f^*M\otimes N)$$
for every $M\in \RigSH^{(\hyp)}_{\tau}(X;\Lambda)$ and 
$N\in \RigSH^{(\hyp)}_{\tau}(Y;\Lambda)$.
\end{cor}

\begin{proof}
This is an immediate consequence of Theorem
\ref{thm:projection-formula-!-star-rigsh}
and Proposition
\ref{prop:equiv-RigSH-!-star-compat-module}.
\end{proof}

\begin{cor}
\label{cor:weak--dual-}
Let $f:Y \to X$ be a morphism of 
rigid analytic spaces which is locally of finite type. 
Then there are equivalences 
$$f^!\underline{\Hom}(M,M')\simeq \underline{\Hom}(f^*M,f^!M')
\quad \text{and} \quad
\underline{\Hom}(f_!N,M)\simeq f_*\underline{\Hom}(N,f^!M)$$
for $M,M'\in \RigSH^{(\hyp)}_{\tau}(X;\Lambda)$ and 
$N\in \RigSH^{(\hyp)}_{\tau}(Y;\Lambda)$.
\end{cor}

\begin{proof}
These are obtained by adjunction from the equivalences
$$(M\otimes -)\circ f_! \simeq f_!\circ (f^*M\otimes -)
\qquad \text{and} \qquad 
(-\otimes f_!N)\simeq f_!\circ (-\otimes N) \circ f^*$$
which are provided by Corollary
\ref{cor:proj-form-f-!}.
\end{proof}

\subsection{Compatibility with the analytification functor}

$\empty$

\smallskip

\label{subsect:compat-with-analytif}

In this last subsection, we prove the 
compatibility of the exceptional functors 
with the analytification functor
\eqref{eqn:analytific-funct-ran-star}. 
We first start with the algebraic analogue of Theorem 
\ref{thm:exist-functo-exceptional-image}.
(Below, for a scheme $S$, we denote by $\Sch^{\lft}/S$
the category of locally of finite type $S$-schemes.)

\begin{thm}
\label{thm:exist-functo-exceptional-image-algebraic}
There are functors 
\begin{equation}
\label{eq-thm:exist-functo-exceptional-image-alg-1}
\begin{array}{rcl}
\SH^{(\hyp)}_{\tau}(-;\Lambda)_!^* & : & 
{\displaystyle \int_{\Sch^{\op}}\Sch^{\lft} \to \Prl}\\
& \vspace{-.3cm} & \\
\SH^{(\hyp)}_{\tau}(-;\Lambda)_*^! & : &
{\displaystyle
\left(\int_{\Sch^{\op}}\Sch^{\lft}\right)^{\op} \to \Prr}
\end{array}
\end{equation}
which are exchanged by the equivalence 
$(\Prl)^{\op}\simeq \Prr$ and which admit
the following informal description.
\symn{$\SH(-)_{\shriek}^*$}
\symn{$\SH(-)^{\shriek}_*$}
\begin{itemize}

\item These functors send an object $(S,X)$, with $S$ a
scheme and $X$ an object of $\Sch^{\lft}/S$,
to the $\infty$-category $\SH^{(\hyp)}_{\tau}(X;\Lambda)$.

\item These functors send an arrow $(g,f):(S,Y) \to (T,X)$, 
consisting of morphisms 
$g:T \to S$ and $f:T\times_S Y \to X$, to the functors
$f_!\circ g'^*$ and $g'_*\circ f^!$ respectively,
with $g':T\times_S Y \to Y$ the base change of $g$.

\end{itemize}
Moreover, the functors in
\eqref{eq-thm:exist-functo-exceptional-image-alg-1}
satisfy the following properties.
\begin{enumerate}

\item[(1)] The ordinary functors 
\begin{equation}
\label{eq-thm:exist-functo-exceptional-image-alg-1.6}
\begin{array}{rcl}
\SH^{(\hyp)}_{\tau}(-;\Lambda)^* & : &
\Sch^{\op}\to \Prl\\
& \vspace{-.3cm} &\\
\SH^{(\hyp)}_{\tau}(-;\Lambda)_* & : &
\Sch^{\op}\to \Prr
\end{array}
\end{equation}
are obtained from the functors in 
\eqref{eq-thm:exist-functo-exceptional-image-alg-1}
by composition with the functor $\Sch^{\op}\to 
\int_{\Sch^{\op}}\Sch^{\lft}$, given by 
$S\mapsto (S,S)$.

\item[(2)] For a scheme $S$, consider the functors
\begin{equation}
\label{eq-thm:exist-functo-exceptional-image-alg-1.81}
\begin{array}{rcl}
\SH^{(\hyp)}_{\tau}(-;\Lambda)_! & : & 
\Sch^{\lft}/S \to \Prl \\
& \vspace{-.3cm} & \\
\SH^{(\hyp)}_{\tau}(-;\Lambda)^! & : & 
\Sch^{\lft}/S \to \Prr
\end{array}
\end{equation}
obtained from the functors in
\eqref{eq-thm:exist-functo-exceptional-image-alg-1}
by restriction to $\Sch^{\lft}/S$. For a morphism
$f:Y \to X$ in $\Sch^{\lft}/S$, denote by 
$f_!$ and $f^!$ the images of $f$ by these functors respectively.
If $f$ is proper there is an equivalence 
$f_!\simeq f_*$ and if $f$ is smooth
there is an equivalence $f^!\simeq \Th(\Omega_f)\circ f^*$.

\item[(3)] The functor $\SH^{(\hyp)}_{\tau}(-;\Lambda)^*_!$ 
can be lifted to a functor 
$$\SH^{(\hyp)}_{\tau}(-;\Lambda)_!^{\otimes}:  
\int_{\Sch^{\op}}\Sch^{\lft} \to \Mod(\Prl)$$
which is part of a commutative square 
$$\xymatrix{\int_{\Sch^{\op}}\Sch^{\lft}
\ar[rr]^-{\SH^{(\hyp)}_{\tau}(-;\,\Lambda)^{\otimes}_!} \ar[d] && 
\Mod(\Prl) \ar[d]\\
\Sch^{\op} \ar[rr]^{\SH^{(\hyp)}_{\tau}(-;\,\Lambda)^{\otimes}} && 
\CAlg(\Prl).\!}$$

\end{enumerate}
\end{thm}

\begin{proof}
This is the algebraic analogue of the combination of 
Theorems 
\ref{thm:exist-functo-exceptional-image}
and
\ref{thm:projection-formula-!-star-rigsh}.
The proof in the algebraic setting is totally similar
to the proof in the rigid analytic setting. 
However, we spend some lines discussing
the construction of the functors in
\eqref{eq-thm:exist-functo-exceptional-image-alg-1} 
in order to introduce some notation which will be useful 
for the proof of Theorem
\ref{thm:compatibility-of-exceptional-with-Ran-star} 
below.

Given a scheme $S$, we denote by $\Sch^{\proper}/S$
the category of proper $S$-schemes.
We also denote by $\Sch^{\cp}/S$ the category of 
compactifiable $S$-schemes, i.e., those $S$-schemes admitting 
an open immersion into a proper $S$-scheme. 
We have an inclusion $\Sch^{\cp}/S \subset \Sch^{\sft}/S$
which is an equality when $S$ is quasi-compact and quasi-separated
by Nagata's compactification theorem
(see \cite[Theorem 4.1]{nagata-deligne}).
We denote by $\Comp/S$ the category whose objects are pairs
$(X,\overline{X})$ where $X$ is an $S$-scheme and 
$\overline{X}$ is a compactification of $X$ over $S$.
We have a functor 
$\mathfrak{d}_S:\Comp/S \to \Sch^{\cp}/S$,
given by $(X,\overline{X}) \mapsto X$.

The construction of the functors in 
\eqref{eq-thm:exist-functo-exceptional-image-alg-1}
starts with the functor
\begin{equation}
\label{eq-thm:exist-funct-exceptional-image-alg-17}
\SH^{(\hyp)}_{\tau}(-;\Lambda)^{*,\,*}:\int_{\Sch^{\op}}
(\Sch^{\proper})^{\op}\to \Prl
\end{equation}
obtained from $\SH^{(\hyp)}_{\tau}(-;\Lambda)^*$ 
by composition with the functor 
$\int_{\Sch^{\op}} (\Sch^{\proper})^{\op} \to
\Sch^{\op}$, given by $(S,X)\mapsto X$.
The condition ($\star$) in Construction 
\ref{cons:adjointable-cocart-fibration}
is satisfied for \eqref{eq-thm:exist-funct-exceptional-image-alg-17} 
by the proper base change theorem (see Proposition
\ref{prop:prop-base-proformal}(1)). 
Using this construction, we obtain a functor  
\begin{equation}
\label{eq-thm:exist-funct-exceptional-image-alg-31}
\SH^{(\hyp)}_{\tau}(-;\Lambda)^*_*:
\int_{\Sch^{\op}}
\Sch^{\proper}\to \Prl
\end{equation}
sending an arrow $(g,f):(S,Y)\to (T,X)$, consisting of morphisms 
$g:T \to S$ and $f:T\times_S Y \to X$, to the composite functor 
$f_*\circ g'^*:\SH^{(\hyp)}_{\tau}(Y;\Lambda) \to 
\SH^{(\hyp)}_{\tau}(X;\Lambda)$, with $g':T\times_SY\to Y$
the base change of $g$. 
Let $S$ be a scheme. For $(X,\overline{X})$ in $\Comp/S$, 
we denote by 
$\SH^{(\hyp)}_{\tau}((X,\overline{X});\Lambda)^*_!$
the essential image of the fully faithful embedding 
$$v_{\sharp}:\SH^{(\hyp)}_{\tau}(X;\Lambda)
\to \SH^{(\hyp)}_{\tau}(\overline{X};\Lambda)$$
where $v:X \to \overline{X}$ is the given open immersion.
By Proposition
\ref{prop:prop-base-proformal}(2), 
the analogue of Proposition
\ref{prop:subcat-rigsh-star-!}
holds true for the functor
\eqref{eq-thm:exist-funct-exceptional-image-alg-31}.
Thus, we may apply Lemma \ref{lem:fully-faith-sub-functor}
to obtain a functor 
\begin{equation}
\label{eq-thm:exist-funct-exceptional-image-alg-53}
\SH^{(\hyp)}_{\tau}((-,-);\Lambda)^*_!:
\int_{\Sch^{\op}} \Comp\to \Prl.
\end{equation}
By left Kan extension along the functor 
$\mathfrak{d}:\int_{\Sch^{\op}}
\Comp \to \int_{\Sch^{\op}} \Sch^{\cp}$,
we deduce from 
\eqref{eq-thm:exist-funct-exceptional-image-alg-53}
the functor 
\begin{equation}
\label{eq-thm:exist-funct-exceptional-image-alg-61}
\SH^{(\hyp)}_{\tau}(-;\Lambda)^*_!:
\int_{\Sch^{\op}} \Sch^{\cp}\to \Prl.
\end{equation}
The analogue of Lemma 
\ref{lem:left-kan-extension-properties}
is also valid here. 
Finally, the first functor in 
\eqref{eq-thm:exist-functo-exceptional-image-alg-1}
is obtained by left Kan extension along 
$\int_{\Sch^{\op}}\Sch^{\cp}
\to \int_{\Sch^{\op}}\Sch^{\lft}$
from \eqref{eq-thm:exist-funct-exceptional-image-alg-61}.
\end{proof}

\begin{rmk}
\label{rmk:star-!-general-homotopical-functor}
Theorem \ref{thm:exist-functo-exceptional-image-algebraic}
holds true with the same proof for any 
stable homotopical functor in the sense of 
\cite[D\'efinition 1.4.1]{ayoub-th1}.
More precisely, given a functor 
$\mathsf{H}^*:\Sch^{\op}\to \Prl$, $f\mapsto f^*$ satisfying 
the $\infty$-categorical versions of the
properties (1)--(6) listed in \cite[\S 1.4.1]{ayoub-th1}, 
there are functors 
\begin{equation}
\label{eq-rmk:star-!-general-homotopical-functor}
\begin{array}{rcl}
\mathsf{H}(-)_!^* & : & 
{\displaystyle \int_{\Sch^{\op}}\Sch^{\lft} \to \Prl}\\
& \vspace{-.3cm} & \\
\mathsf{H}(-)_*^! & : &
{\displaystyle
\left(\int_{\Sch^{\op}}\Sch^{\lft}\right)^{\op} \to \Prr}
\end{array}
\end{equation}
satisfying the properties (1) and (2) of Theorem 
\ref{thm:exist-functo-exceptional-image-algebraic}.
Moreover, if $\mathsf{H}$ admits a lift to a functor 
$\mathsf{H}^{\otimes}:\Sch^{\op}\to \CAlg(\Prl)$
such that the projection formula holds, then 
property (3) of Theorem 
\ref{thm:exist-functo-exceptional-image-algebraic}
is also satisfied.
\end{rmk}

\begin{thm}
\label{thm:compatibility-of-exceptional-with-Ran-star}
Let $A$ be an adic ring. Set $S=\Spf(A)^{\rig}$ and 
$U=\Spec(A)\smallsetminus \Spec(A/I)$
where $I\subset A$ is an ideal of definition.
There is a commutative cube of $\infty$-categories
$$\xymatrix{\int_{(\Sch^{\lft}/U)^{\op}}\Sch^{\lft} 
\ar[rr]^-{(-)^{\an}} 
\ar[dd] \ar[dr]|-{\SH^{(\hyp)}_{\tau}(-;\,\Lambda)^{\otimes}_!} & &
\int_{(\RigSpc^{\lft}/S)^{\op}}\RigSpc^{\lft} \ar'[d][dd] 
\ar[dr]|-{\RigSH^{(\hyp)}_{\tau}(-;\,\Lambda)^{\otimes}_!} &\\
& \Mod(\Prl)  \ar@{=}[rr] \ar[dd] & ^{\;} & \Mod(\Prl) \ar[dd] \\
(\Sch^{\lft}/U)^{\op} \ar'[r][rr]^-{(-)^{\an}} 
\ar[dr]|-{\SH^{(\hyp)}_{\tau}(-;\,\Lambda)^{\otimes}} & & (\RigSpc^{\lft}/S)^{\op} \ar[dr]|-{\RigSH^{(\hyp)}_{\tau}(-;\,\Lambda)^{\otimes}} & \\
& \CAlg(\Prl)  \ar@{=}[rr] & & \CAlg(\Prl).\!}$$
In particular, there is a natural transformation 
\begin{equation}
\label{eq-thm:compatibility-of-exceptional-with-Ran-star-1}
\An^*:\SH^{(\hyp)}_{\tau}(-;\Lambda)^*_! \to 
\RigSH^{(\hyp)}_{\tau}((-)^{\an};\Lambda)^*_!
\end{equation}
between functors from 
$\int_{(\Sch^{\lft}/U)^{\op}}\Sch^{\lft}$
to $\Prl$ which extends the morphism of $\Prl$-valued 
presheaves $\An^*$ underlying \eqref{eq-prop:analytif-presheaf-1} 
in Proposition \ref{prop:analytif-presheaf}.
\end{thm}

\begin{proof}
For simplicity, we only construct the natural 
transformation 
\eqref{eq-thm:compatibility-of-exceptional-with-Ran-star-1}.
It will be clear from the construction how to lift 
this natural transformation into a commutative square 
which is part of a commutative cube as in the statement.

We use the notation introduced in the proof of Theorem
\ref{thm:exist-functo-exceptional-image-algebraic}.
By construction, the functor 
$$\SH^{(\hyp)}_{\tau}(-;\Lambda)^*_!:
\int_{(\Sch^{\lft}/U)^{\op}}\Sch^{\lft}\to \Prl$$ 
is a left Kan extension along the functor 
$$\mathfrak{d}':
\int_{(\Sch^{\lft}/U)^{\op}} \Comp \to 
\int_{(\Sch^{\lft}/U)^{\op}}\Sch^{\lft},$$
given by $(S,(X,\overline{X}))\mapsto (S,X)$,
of the functor 
$$\SH^{(\hyp)}_{\tau}((-,-);\Lambda)^*_!:
\int_{(\Sch^{\lft}/U)^{\op}} \Comp
\to \Prl$$
obtained from
\eqref{eq-thm:exist-funct-exceptional-image-alg-53}
by restriction.
(Here, we are combining the two left Kan extensions
from the proof of Theorem 
\ref{thm:exist-functo-exceptional-image-algebraic}.)
By the universal property of left Kan extensions, 
it is thus enough to construct a natural transformation 
$$\An^*:\SH^{(\hyp)}_{\tau}((-,-);\Lambda)^*_!
\to \RigSH^{(\hyp)}_{\tau}((-)^{\an};\Lambda)^*_!\circ 
\mathfrak{d}'$$
between functors from 
$\int_{(\Sch^{\lft}/U)^{\op}}\Comp$
to $\Prl$. Now, consider the functors 
$$\mathfrak{w}:\int_{(\Sch^{\lft}/U)^{\op}}\Comp
\to \int_{(\Sch^{\lft}/U)^{\op}}\Sch^{\proper}
\quad \text{and} \quad 
\mathfrak{w}':\int_{(\Sch^{\lft}/U)^{\op}}\Comp
\to \int_{(\Sch^{\lft}/U)^{\op}}\Sch^{\lft}$$
given by $(S,(X,\overline{X}))\mapsto (S,\overline{X})$.
The obvious natural transformation 
$v:\mathfrak{d}'\to \mathfrak{w}'$
induces a natural transformation
$$v^{\an}_!:\RigSH^{(\hyp)}_{\tau}((-)^{\an};\Lambda)^*_!
\circ \mathfrak{d}'
\to \RigSH^{(\hyp)}_{\tau}((-)^{\an};\Lambda)^*_!\circ \mathfrak{w}'$$
which is objectwise a fully faithful embedding. Thus, we may 
obtain $\RigSH^{(\hyp)}_{\tau}((-)^{\an};\Lambda)^*_!\circ \mathfrak{d}'$
from $\RigSH^{(\hyp)}_{\tau}((-)^{\an};\Lambda)^*_!\circ \mathfrak{w}'$
by applying Lemma \ref{lem:fully-faith-sub-functor} 
to the essential images of the fully faithful embeddings 
$$v_!^{\an}: \RigSH^{(\hyp)}_{\tau}(X^{\an};\Lambda) \to
\RigSH^{(\hyp)}_{\tau}(\overline{X}{}^{\an};\Lambda)$$
for the objects $(S,(X,\overline{X}))$.
Since $\SH^{(\hyp)}_{\tau}((-,-);\Lambda)^*_!$ 
is constructed from 
$\SH^{(\hyp)}_{\tau}(-;\Lambda)^*_*\circ 
\mathfrak{w}$ in the same way, we are left to 
construct a natural transformation 
$$\SH^{(\hyp)}_{\tau}(-;\Lambda)^*_*\circ 
\mathfrak{w}\to 
\RigSH^{(\hyp)}_{\tau}((-)^{\an};\Lambda)^*_!\circ \mathfrak{w}'.$$
The functor $\mathfrak{w}'$ factors through 
$\int_{(\Sch^{\lft}/U)^{\op}}\Sch^{\proper}$.
Thus, by Proposition 
\ref{prop:f-proper-f-star-f-!}, 
it is enough to construct a natural transformation 
$$\SH^{(\hyp)}_{\tau}(-;\Lambda)^*_* \to 
\RigSH^{(\hyp)}_{\tau}((-)^{\an};\Lambda)^*_*$$
between functors from
$\int_{(\Sch^{\lft}/U)^{\op}}\Sch^{\proper}$ to $\Prl$.
Equivalently, we need to construct a functor 
$$\int_{(\Sch^{\lft}/U)^{\op}\times \Delta^1}
\Sch^{\proper}\to \Prl,$$
which restricts to 
$\SH^{(\hyp)}_{\tau}(-;\Lambda)^*_*$ over $\{0\}\subset \Delta^1$ 
and to $\RigSH^{(\hyp)}_{\tau}((-)^{\an};\Lambda)^*_*$ over
$\{1\}\subset \Delta^1$.
For this, we apply Construction
\ref{cons:adjointable-cocart-fibration}
to the composite functor 
$$\int_{(\Sch^{\lft}/U)^{\op}\times \Delta^1}
(\Sch^{\proper})^{\op}\to \Delta^1\times (\Sch^{\lft}/U)^{\op}
\to \Prl$$
where the first functor is given by 
$((S,\epsilon),X) \mapsto (\epsilon,X)$
and the second one classifies the natural transformation 
$\An^*:\SH^{(\hyp)}_{\tau}(-;\Lambda)
\to \RigSH^{(\hyp)}_{\tau}((-)^{\an};\Lambda)$
underlying \eqref{eq-prop:analytif-presheaf-1} 
in Proposition \ref{prop:analytif-presheaf}.
That condition ($\star$) in Construction 
\ref{cons:adjointable-cocart-fibration} 
is satisfied, follows from Propositions 
\ref{prop:ran-star-com-f-lower-star-prelim} and
\ref{prop:prop-base-proformal}(1),
and Theorem \ref{thm:prop-base}(1).
\end{proof}

\clearpage
\printnomenclature

\clearpage
\printindex
\end{document}